\documentclass[a4paper,11pt]{article}%
 \usepackage[margin=1in,footskip=0.5in]{geometry}
 
\usepackage{setspace}
 \usepackage[symbol]{footmisc}
 \usepackage{caption}
\usepackage{subcaption}

\usepackage{amsthm}
\usepackage{amsmath}
\usepackage{amssymb,amsmath,amsthm,amsfonts,mathrsfs,xcolor}
\usepackage{natbib}
\usepackage[colorlinks,citecolor=blue,urlcolor=blue,filecolor=blue,backref=page]{hyperref}
\usepackage{graphicx}

 \usepackage{tcolorbox}
\theoremstyle{plain}
\newtheorem{theorem}{Theorem}
 \newtheorem{remark}{Remark}
  \newtheorem{claim}{Claim}
 \newtheorem{corollary}{Corollary}
 \newtheorem{lemma}{Lemma}
 \newtheorem{definition}{Definition}
 \newtheorem*{definition*}{Definition}
 
\usepackage{times}

\usepackage{amssymb,amsmath,amsthm,amsfonts,bbm,bm,mathrsfs,xcolor}
 \usepackage{float}  
 \usepackage{xr}
\usepackage{multirow}

\usepackage{bigints}

\def\m{\mathcal}
\def\mb{\mathbb}

\def\ms{\mathscr}
\def\sd{{[d]}}

\def\wt{\widetilde}
\def\wh{\widehat}

\def\dd{{\rm d}}

\def\sk{{[k]}}
\def\mf{\mathfrak}

\def\ov{\overline}

\newcommand{\mnorm}[1]{{\vert\kern-0.25ex\vert\kern-0.25ex\vert #1 
    \vert\kern-0.25ex\vert\kern-0.25ex\vert}}
\newcommand{\bmnorm}[1]{{\big\vert\kern-0.25ex\big\vert\kern-0.25ex\big\vert #1 
    \big\vert\kern-0.25ex\big\vert\kern-0.25ex\big\vert}}
\newcommand{\Bmnorm}[1]{{\Big\vert\kern-0.25ex\Big\vert\kern-0.25ex\Big\vert #1 
    \Big\vert\kern-0.25ex\Big\vert\kern-0.25ex\Big\vert}}
\newcommand{\bbmnorm}[1]{{\bigg\vert\kern-0.25ex\bigg\vert\kern-0.25ex\bigg\vert #1 
    \bigg\vert\kern-0.25ex\bigg\vert\kern-0.25ex\bigg\vert}}
\newcommand{\BBmnorm}[1]{{\Bigg\vert\kern-0.25ex\Bigg\vert\kern-0.25ex\Bigg\vert #1 
    \Bigg\vert\kern-0.25ex\Bigg\vert\kern-0.25ex\Bigg\vert}}

\title{Minimax Optimal Rates for Regression on Manifolds and Distributions}

 \author{Rong Tang\textsuperscript{$\ast$} and Yun Yang\textsuperscript{$\dagger$}}

\date{\textsuperscript{$\ast$}Department of Mathematics, Hong Kong University of Science and Technology\\
\textsuperscript{$\dagger$}Department of Mathematics, University of Maryland, College Park}
 
\begin{document}

\maketitle

\begin{abstract}
Distribution regression seeks to estimate the conditional distribution of a multivariate response given a continuous covariate. This approach offers a more complete characterization of dependence than traditional regression methods. Classical nonparametric techniques often assume that the conditional distribution has a well-defined density, an assumption that fails in many real-world settings. These include cases where data contain discrete elements or lie on complex low-dimensional structures within high-dimensional spaces. In this work, we establish minimax convergence rates for distribution regression under nonparametric assumptions, focusing on scenarios where both covariates and responses lie on low-dimensional manifolds. We derive lower bounds that capture the inherent difficulty of the problem and propose a new hybrid estimator that combines adversarial learning with simultaneous least squares to attain matching upper bounds. Our results reveal how the smoothness of the conditional distribution and the geometry of the underlying manifolds together determine the estimation accuracy.
\medskip

\noindent {\bf Keywords}: Conditional distribution estimation; Manifold learning; Distribution regression; Minimax rate; Conditional generative models;  Adversarial learning.
\end{abstract}
 
 \tableofcontents

\section{Introduction}
Distribution regression (or more precisely, distribution-on-vector regression), where the goal is to estimate the conditional distribution $\mu^*_{Y|X}$ of a random response vector $Y\in\mathbb{R}^{D_Y}$ given a continuous covariate $X\in\mathbb{R}^{D_X}$, is a fundamental problem in statistics and machine learning. Unlike traditional regression~\citep{christensen2002plane,hardle1990applied,koenker2005quantile} or classification~\citep{bishop2006pattern}, which typically involves a univariate response (i.e., $D_Y = 1$) and predicts scalar or categorical outcomes, distribution regression aims to recover the full conditional distribution of a potentially multivariate response, providing a more comprehensive characterization of the dependence between $X$ and $Y$~\citep{dinardo2001nonparametric}, which may represent complex objects encoded or embedded as numerical vectors, including images, texts, or other structured data. In particular, distribution regression allows for capturing how $\mu_{Y|X=x}$ evolves as the covariate $x$ varies, enabling a richer understanding of conditional variability, skewness, uncertainty and multiple-modality~\citep{rodriguez2025density}. This framework is especially important in applications where characterizing the entire distribution, rather than just its mean or quantiles, is crucial, such as in biomedical sciences~\citep{krishnaswamy2014conditional}, climate modeling~\citep{guinness2018compression} and econometrics~\citep{li2007nonparametric}. 

There is a vast literature on nonparametric density regression (conditional density estimation), where the conditional distribution $\mu^*_{Y|X}$ is assumed to have a density function with respect to the Lebesgue measure on $\mathbb{R}^{D_Y}$. However, many existing methods, particularly classical nonparametric estimators based on kernel smoothing~\citep{bashtannyk2001bandwidth, izbicki2016nonparametric, li2022minimax}, have several notable limitations. A primary drawback of these approaches is their reliance on the existence of a conditional density function, an assumption that often fails when the response variable $Y$ contains discrete components or is embedded in a high-dimensional ambient space with low-dimensional singular structures, as is common in structured data environments~\citep{wang2020systematic,bellet2013survey}. As a consequence, these methods are primarily effective in low-dimensional settings but struggle as dimensionality increases, ultimately suffering from the curse of dimensionality~\citep{pope2021intrinsic,latorre2021effect}. Furthermore, classical density regression methods generally lack adaptability to the intrinsic geometric structure of data, such as underlying manifold structures that are common in modern high-dimensional datasets~\citep{gong2019intrinsic,aghajanyan2020intrinsic}. This inability to exploit low-dimensional representations limits their effectiveness in capturing complex dependencies and accurately modeling conditional distributions in modern data environments, which often involve high-dimensional, complex data such as images in computer vision~\citep{parker2010algorithms}, medical imaging~\citep{suetens2017fundamentals}, and signal processing~\citep{francca2021overview}, as well as text in data mining~\citep{zhai2016text}, natural language processing~\citep{kao2007natural}, and public health~\citep{yang2022large}.

These limitations of classical density regression motivate us to study the statistical properties of distribution regression, which can accommodate general data types and singular distributions. In particular, the recent surge in conditional generative models---such as conditional generative adversarial networks~\citep{mirza2014conditional}, conditional diffusion models~\citep{songscore,zhang2023adding}, and conditional normalizing flows~\citep{abdelhamed2019noise,winkler2019learning}---demonstrates their effectiveness and efficiency in generating new data given a covariate (or control variable) in complex environments. These models approximate complex conditional distributions by learning the underlying data-generating processes, making them powerful tools for tasks such as image-to-image translation~\citep{isola2017image}, medical image synthesis~\citep{dar2019image}, and super-resolution imaging~\citep{zhao2019simultaneous}.
Consequently, conditional generative models can be regarded as implicit distribution regression methods, as they do not explicitly estimate the conditional density or cumulative distribution function but instead generate samples that follow the underlying conditional distribution.
However, despite their empirical success, the theoretical understanding of their statistical properties remains limited. In particular, it is unclear how well these models approximate the true conditional distribution and under what conditions they achieve optimal performance. This gap highlights the need to establish rigorous theoretical guarantees for distribution regression, particularly in terms of minimax rates, to provide a deeper understanding of the fundamental limits of learning conditional distributions.

In this work, we investigate the minimax convergence rates for distribution regression under nonparametric settings, where both the response variable $Y$ and the covariate $X$ may be high-dimensional but possess an underlying low-dimensional manifold structure. This setting is particularly relevant for modern conditional generative models using deep neural networks~\citep{sohn2015learning,salakhutdinov2015learning}, as many complex data types---such as images, text, and other structured objects---reside on low-dimensional manifolds despite being represented in high-dimensional ambient spaces. Moreover, deep neural networks are naturally suitable for learning low-dimensional nonlinear features, making them inherently adaptive to such data structures~\citep{schmidt2019deep,kohler2021rate,10.1214/19-AOS1875}. Unlike the unconditional distribution estimation setting (e.g.,~\cite{10.1214/23-AOS2291}), where the data is supported on a single manifold, the conditional distribution setting is more intricate. Both the covariate $X$ and the response $Y$ can reside on distinct manifolds, and more importantly, the manifold supporting $Y$ may depend on $X$. This dependence transforms the problem of recovering the support of $Y$ into a manifold regression problem, which is highly nontrivial and remains largely unexplored in the existing literature.

Concretely, we consider a random design distribution regression setting where the covariate $X$ follows a marginal distribution $\mu^*_X$ supported on a $d_X$-dimensional submanifold $\mathcal{M}_X$ within the ambient covariate space $\mathbb{R}^{D_X}$. Furthermore, our target of interest, the conditional distribution $\mu^*_{Y|X=x}$ of $Y$ given $X = x$, is supported on a $d_Y$-dimensional, $\beta_Y$-smooth submanifold $\mathcal{M}_{Y|x}$ (c.f.~Definition~\ref{def:manifold1}) within the ambient response space $\mathbb{R}^{D_Y}$, which may or may not vary with $x$. In cases where $\mathcal{M}_{Y|x}$ depends on $x$, we assume that its dependence on $x$ is $\beta_X$-smooth (c.f.~Definition~\ref{def:manifold}). We also assume that $\mu^*_{Y|X=x}$ admits a density function with respect to the volume measure of $\mathcal{M}_{Y|x}$ that is $\alpha_Y$-smooth in $y$ and $\alpha_X$-smooth in $x$ (c.f.~Definition~\ref{def:marginalsmooth}). 
The data $\big\{(X_i, Y_i)\big\}_{i=1}^n$ consists of $n$ i.i.d.~copies of $(X, Y)$ sampled from the joint distribution $\mu^*_{XY} = \mu^*_X \mu^*_{Y|X}$.To quantify the estimation error, we follow \cite{tang2023minimax} and use the integral probability metric (IPM) \citep{muller1997integral}, also known as the adversarial loss in the machine learning literature (e.g.,~\cite{singh2018nonparametric}), to measure the closeness between two probability measures, which may be mutually singular in our setting. In particular, we consider the $\gamma$-H\"{o}lder IPM, denoted as $d_\gamma$, which is indexed by a smoothness parameter $\gamma \geq 0$ that regulates the strength of the metric and balances the trade-off between distribution support mismatch and relative density differences over the support (c.f.~equation\eqref{def:dgamma} and the subsequent discussion). Notably, the $d_\gamma$ metric includes classical total variation distance ($\gamma = 0$) and $1$-Wasserstein distance $W_1$ ($\gamma = 1$) as special cases.

\subsection{Related Work}
There is a vast literature on nonparametric density regression in both the statistics and machine learning communities, where proposed estimators range from classical nonparametric methods based on kernel smoothing~\citep{key254987m, fan2004crossvalidation, kernel-smoothing, bashtannyk2001bandwidth,izbicki2016nonparametric, li2022minimax} to Bayesian nonparametric approaches \citep{norets2017adaptive} and more recent methods leveraging deep neural networks \citep{rothfuss2019conditional}. These works mostly consider the classical setting where the dimension $D_X$ of the covariate $X$ is low and the conditional distribution $\mu^*_{Y|X=\cdot}$ admits a density function, corresponding to our setting with $d_X=D_X$, $d_Y=D_Y$, $\m M_X=\mb R^{D_X}$, $\m M_{Y|x} = \mb R^{D_Y}$ for each $x\in \mb R^{D_X}$, and $\beta_X=\beta_Y=\infty$. 
In this classical setting, \cite{li2022minimax} establish the minimax rate of conditional density estimation as $n^{-1/(2+D_X/\alpha_X+D_Y/\alpha_Y)}$ under the total variation metric (corresponding to the $d_0$ IPM) when $\alpha_X \in [0,1]$. \cite{10.1214/23-AOS2270} studies the minimax rate for conditional density estimation under the Kullback-Leibler (KL) risk, providing both upper and lower bounds expressed in terms of empirical Hellinger entropy.

Recent work has also addressed statistical inference problems involving manifolds. Several studies \citep{tang2023minimax, ozakin2009submanifold, berenfeld2024estimating, berenfeld2021density, cholaquidis2022level, divol2022measure} focus on the problem of (unconditional) distribution estimation on an unknown manifold. Notably, \cite{tang2023minimax} derives the minimax rate $n^{-1/2} + n^{-(\alpha_Y+\gamma)/(2\alpha_Y+d_Y)} + n^{-\gamma\beta_Y/d_Y}$ for estimating an $\alpha_Y$-smooth distribution supported on a $d_Y$-dimensional, $\beta_Y$-smooth ($\beta_Y \geq \alpha_Y + 1$) manifold $\mathcal{M}_Y \subset \mathbb{R}^{D_Y}$ with respect to the $d_\gamma$ metric for all $\gamma \geq 0$. Under a similar setting, \cite{divol2022measure} shows that the minimax rate under the $p$-Wasserstein distance is $n^{-1/2} + n^{-(\alpha_Y+\gamma)/(2\alpha_Y+d_Y)}$ for any $p > 0$. Some other studies \citep{genovese2012manifold, genovese2012minimax, divol2021minimax, 10.1214/18-AOS1685} consider the problem of manifold estimation, which corresponds to support recovery for singular distributions. For instance, \cite{10.1214/18-AOS1685} establishes that the minimax rate for estimating a boundaryless, $\beta_Y$-smooth ($\beta_Y \geq 2$), $d_Y$-dimensional submanifold under the Hausdorff distance is $n^{-\beta_Y/d_Y}$.

There has been a recent line of work leveraging generative models, such as generative adversarial networks (GAN) and diffusion models, for implicit (conditional) distribution estimation via sampling, such as~\cite{oko2023diffusion,chen2023score,wang2024diffusion,li2024adapting,de2021diffusion,lee2022convergence,chen2022sampling,lee2023convergence,chen2023restoration,tang2024adapt,li2024sharp,huang2022sampling,liu2021wasserstein,chen2024overview,fu2024unveil,li2024towards,azangulov2024convergence,tang2024conditional}. To name a few most relevant to our problem, in the unconditional distribution estimation case of estimating the distribution $\mu^*_Y$ without covariate $X$, 
\cite{oko2023diffusion} show that diffusion models can achieve the respective minimax rate $n^{-1/(2+D_Y/\alpha_Y)}$ under the total variation metric and $n^{-(1+1/\alpha_Y)/(2+D_Y/\alpha_Y)}$  under the $W_1$ distance. Furthermore, \cite{tang2024adapt} extend the results of \cite{oko2023diffusion} to the manifold setting and derive a convergence rate $n^{-1/2} + n^{-(\alpha_Y+\gamma)/(2\alpha_Y+d_Y)} + n^{-\gamma\beta_Y/(2\alpha_Y+d_Y)}$ under the $d_\gamma$ distance, showing that diffusion models are minimax-optimal at least when $\gamma\in[0,1]$, covering the total variation distance and the $W_1$ distance. In the conditional generative model setting, \cite{huang2022sampling} propose a conditional GAN-based approach~\citep{mirza2014conditional} and establish the consistency of the resulting conditional density estimator, though no convergence rates or error bounds are provided. Meanwhile, \cite{liu2021wasserstein} adopt a Wasserstein generative approach for conditional distribution estimation and derive a convergence rate of $n^{-1/(1+D_X+D_Y)}$ under the $W_1$ distance (corresponding to the $d_1$ IPM). For conditional diffusion models, \cite{chen2024overview} provides a recent survey of related theoretical investigations on conditional score estimation and the resulting sample complexity.
Among the most relevant works to ours, \cite{fu2024unveil} explore the theoretical properties of conditional diffusion models under the classical setting without a manifold structure and derive a convergence rate of $n^{-\alpha_Y/(2\alpha_Y+D_X+D_Y)}$ relative to the total variation distance, under the special case where the conditional distribution has the same smoothness level in $X$ and $Y$, i.e., $\alpha_X=\alpha_Y$. More recently, \cite{tang2024conditional} consider the manifold setting and derive a convergence rate of $n^{-\alpha_X/(2\alpha_X +d_X)} + n^{-(\alpha_Y+1)/(2\alpha_Y+d_Y+d_X\alpha_Y/\alpha_X)}$ relative to the $W_1$ distance for conditional diffusion models when all manifolds are sufficiently smooth.

\subsection{Main Contribution}
In this work, we investigate the minimax convergence rates for distribution regression under nonparametric settings. We derive lower bounds that characterize the fundamental difficulty of the problem and provide matching upper bounds achieved by a new hybrid estimator combining adversarial learning and simultaneous least squares estimation. Our results reveal how the smoothness of the conditional distribution and the geometric properties of the underlying manifolds influence estimation accuracy. Moreover, we extend our analysis to the case where both the response variable and the covariate are high-dimensional but admit an underlying low-dimensional manifold structure. This setting is highly relevant, as many complex data types---such as images, text, and other structured objects---lie on low-dimensional manifolds despite being represented in high-dimensional ambient spaces. Since conditional generative models are particularly effective in modeling such complex data distributions, understanding the minimax rates in this setting provides valuable insights into the theoretical foundations of generative modeling. 

By developing a rigorous theoretical framework, our results precisely characterize the statistical complexity of the problem and establish benchmarks for evaluating the performance of modern conditional generative models. This is particularly relevant in high-dimensional settings, where leveraging low-dimensional structures enhances statistical efficiency. Specifically, our main results on the minimax rate across different regimes are summarized 
below. For all regimes considered, we assume that the covariant space $\m M_X$ exhibits a low-dimensional structure with an intrinsic/effective dimension of $d_X$ (c.f. Definition~\ref{def:Mink}). The minimax rates are presented excluding logarithmic factors.
\begin{description}
    \item[\emph{Regime 1. Classic density regression.}]
    In this setting, the conditional distribution $\mu^*_{Y|X}$ is assumed to admit a density with respect to the Lebesgue measure on $\mathbb{R}^{D_Y}$. The minimax convergence rate takes the form $n^{-{\alpha_X}/{(2\alpha_X+d_X)}} + n^{-{(\alpha_Y+\gamma)}/{(2\alpha_Y+D_Y+\frac{\alpha_Y}{\alpha_X}d_X)}}$. The first term corresponds to the classical minimax rate for estimating an $\alpha_X$-smooth regression function under the $L^2$ loss, as established by~\cite{stone1982optimal}. The second term captures the inherent difficulty of nonparametric conditional density estimation. Notably, when $\alpha_X \in [0,1]$ and $\gamma = 0$, the rate coincides with the minimax rate for conditional density estimation under the total variation metric derived in~\cite{li2022minimax}.

 \item [\emph{Regime 2.  Distribution regression with covariate-independent response space.}]
 In this regime, the support $\m M_{Y|X}$ of the conditional distribution $\mu^*_{Y|X}$ is assumed to be independent of $X$, with a common support $\m M_{Y|X} = \m M_Y$ that is an unknown $\beta_Y$-smooth submanifold of intrinsic dimension $d_Y$. The minimax rate for this setting is
$n^{-{\alpha_X}/{(2\alpha_X+d_X)}} + n^{-{(\alpha_Y+\gamma)}/{(2\alpha_Y+d_Y+\frac{\alpha_Y}{\alpha_X}d_X)}} + n^{-\gamma\beta_Y/d_Y}$.
The first two terms are analogous to those in Regime 1, with the ambient dimension $D_Y$ replaced by the intrinsic dimension $d_Y$, reflecting the lower complexity of the support. The third term accounts for the intrinsic difficulty of estimating the unknown submanifold $\m M_Y$. When $\gamma = 1$, this term matches the minimax rate for estimating a $\beta_Y$-smooth submanifold under the Hausdorff distance, as established in~\cite{10.1214/18-AOS1685}. Additionally, when $d_X = 0$, the minimax rate reduces to that of unconditional distribution estimation on unknown submanifolds, as shown in~\cite{10.1214/23-AOS2291}.

 \item [\emph{Regime 3.  Distribution regression with covariate-dependent response space.} ]
In this regime, the support $\m M_{Y|x}$ of the conditional distribution $\mu^*_{Y|x}$ varies with the covariate $x$, where the collection of conditional response supports $\{\mathcal{M}_{Y|x} :\, x \in \mathcal{M}_X\}$ forms an unknown family of submanifolds that is $(\beta_Y, \beta_X)$-smooth (see Definition~\ref{def:manifold} for details). The minimax rate for this setting is
$n^{-{\alpha_X}/{(2\alpha_X+d_X)}} + n^{-{(\alpha_Y+\gamma)}/{(2\alpha_Y+d_Y+\frac{\alpha_Y}{\alpha_X}d_X)}} + n^{-{\gamma}/{(\frac{d_Y}{\beta_Y}+\frac{d_X}{\beta_X})}}$.
The first two terms are analogous to those in Regime 2, with the key difference arising in the third term, which accounts for the complexity of estimating the submanifold family $\{\mathcal{M}_{Y|x} :, x \in \mathcal{M}_X\}$. As established in Theorem~\ref{th:mainifold}, this term corresponds to the minimax optimal rate for manifold regression over a $(\beta_Y, \beta_X)$-smooth family.
\end{description}


\noindent The remainder of the paper is organized as follows. Section 2 reviews and formalizes key concepts, including multivariate functions with separate smoothness, smooth manifolds, and covariate-dependent manifolds. Sections 3 and 4 present our main theoretical results on the minimax rates for distribution regression under covariate-independent and covariate-dependent response supports, respectively. In Section 5, we introduce conditional distribution estimators that attain the minimax upper bounds across the different regimes. Finally, some concluding discussion are offered in Section 6.


\section{Background and preliminary results}
In this section, we begin by introducing notation. We then present a formal definition of functions with separate smoothness, which will be used to characterize the conditional distribution functions $\mu^*_{Y|X=x}$ of $Y$ given $X=x$ and their supporting manifolds $\mathcal{M}_{Y|x}$. Finally, we provide a brief review of submanifolds, with more detailed background material included in Appendix~\ref{sec:ssmanifold} of the supplementary material. We also formally define the covariate-dependent manifold $\mathcal{M}_{Y|x}$ and characterize its joint smoothness in $Y$ and $x$.

\subsection{Notation}
Recall from the introduction section that $X\in \mb R^{D_X}$ denotes the covariate, distributed as $\mu^*_X$, and $Y\in \mb R^{D_Y}$ denotes the response variable, distributed as $\mu^*_Y$, with the superscript $*$ indicating the ground truth. The conditional distribution of $Y$ given $X=x$ is denoted by $\mu^*_{Y|X=x}$, leading to the joint distribution of $(X,Y)$ as $\mu^*_{XY}=\mu^*_X\mu^*_{Y|X}$, where $\mu^*_X\mu^*_{Y|X}$ represents the generation process of first generating $X\sim \mu^*_X$ and then $[Y|X=x]\sim \mu^*_{Y|X=x}$. When no ambiguity arises, we also use the shorthand $\mu^*$ to denote $\mu^*_{XY}$.
We use $\mu^{*,\otimes n}$ to denote the $n$-fold product of $\mu^*$, and $\big\{(X_i,Y_i)\big\}_{i=1}^n$ to denote a sample of size $n$ drawn from $\mu^*$. The support of $\mu^*_X$ is denoted by $\m M_X$, and the support of $\mu^*_{Y|x}$ is denoted by $\m M_{Y|x}$. We write $\m M=\{(x,y):\,x\in \m M_X, y\in \m M_{Y|x}\}$ as the joint space of $(X,Y)$ and $\m M_Y=\bigcup_{x\in \m M_X} \m M_{Y|x}$ as the marginal space of $Y$.

We use $\|x\|$ to denote the Euclidean norm of a vector $x\in \mb R^d$, and $\mathbf{0}_d$ to represent the $d$-dimensional zero vector. For a set $U \subseteq \mathbb{R}^d$, we denote by $\mathbb{B}_{U}(x, r) = \{ y \in U :\, \|y - x\| < r \}$ the ball of radius $r$ centered at $x$ and contained in $U$. For a measure $\mu$ on $\mb R^d$, we write $\mu|_U$ as the restriction of $\mu$ to $U$, i.e., $\mu|_{U}(A)=\mu(A\cap U)$ for any measurable set $A\subseteq \mb R^d$.  The floor and ceiling functions for $\alpha\in \mathbb{R}$ are denoted by $\lfloor \alpha\rfloor$ and $\lceil \alpha\rceil$, respectively, which round $\alpha$ to the nearest smaller and larger integers. For two real numbers $a,b$, we write $a\vee b$ and $a\wedge b$ as the maximal and minimal value between $a$ and $b$ respectively.    
For any sequence $\{a_n:\,n\geq 1\}$, we use the notation ${\m O}(a_n)$ to mean of order $a_n$ up to multiplicative constant, and use $\wt{\m O}(a_n)$ to mean of order $a_n$ up to multiplicative constant and logarithmic terms of $n$. 

For a positive integer $m$, we use the shorthand $[m]:=\{1,\cdots,m\}$.  We denote by $\mb N$  the set of non-negative integers, $\mb N_{+}$ the set of positive integers, and write $\mathbb{N}_0^d=\{(j_1,\cdots,j_d)\,|\, j_i\in \mathbb{N}, \, \forall i\in [d]\}$ as the set of all multi-indices with $d$ components. For a multi-index $j=(j_1,\cdots,j_d)\in\mathbb{N}_0^d$, we use $|j| = \sum_{i=1}^d j_i$ to mean its size and $j!=\prod_{i=1}^dj_i !$ as the multi-index factorial.  
For a multivariate function $f:\,\mb R^{d_1}\times \mb R^{d_2}\to\mb R$ and two multi-indices $j\in \mb N_0^{d_1}$ and $j'\in \mb N_0^{d_2}$, we denote by $f^{(j,j')}(x,y)$ the mixed partial derivative $\frac{\partial^{|j|+|j'|}f(x,y)}{\partial x^{j_{1}}\cdots \partial x^{j_{d_1}}\partial y^{j'_{1}}\cdots\partial y^{j'_{d_2}}}$ evaluated at $(x,y)$. Moreover,  for a vector-valued function $(x,y)\mapsto f(x,y)=(f_1(x,y),f_2(x,y),\cdots, f_d(x,y))\in\mb R^d$, the notation $f^{(j,j')}(x,y)$ represents the vector of mixed partial derivatives $(f_1^{(j,j')}(x,y),f_2^{(j,j')}(x,y),\cdots, f_d^{(j,j')}(x,y))$ evaluated at $(x,y)$. For a vector $x\in\mb R^d$, we use $x_i$ to denote its $i$-th element.   For $x,y \in \mathbb{R}^d$ and $j\in \mb N_0^d$, we use the shorthand $(x-y)^{j}$ to represent $\prod_{i=1}^d (x_i-y_i)^{j_i}$.

\subsection{Functions with separate smoothness}
In order to allow conditional distributions $\mu^*_{Y|X=x}$ and their supporting manifolds $\m M_{Y|x}$ to have different smoothness levels in $x$ and $y$, we consider two classes of functions with separate smoothness: a weaker class that requires differentiability along each coordinate separately and a stronger class that requires joint differentiability.

Before that, recall the classical definition of the $\alpha$-smooth H\"{o}lder function class $\m H^{\alpha}_r(\mb R^{d})$ with radius $r>0$ over $\mathbb{R}^{d}$, which assumes a uniform smoothness level across all its components, that is,
 \begin{equation*}
    \begin{aligned}
        &\m H^{\alpha}_r(\mb R^{d}):=\Big\{f:\, \mb R^{d}\rightarrow \mathbb{R}:\,\|f\|_{\m H^{\alpha}(\mb R^{d})}=\sum_{j\in \mb N_0^d, |j|<\alpha}\sup_{x\in \mb R^{d}}| f^{(j)}(x)|\\
        &\qquad\qquad\qquad\qquad+\sum_{j\in \mb N_0^d,  \alpha-1\leq|j|<\alpha}\sup_{x, y\in\mb R^{d},\,x\neq y} \big|f^{(j)}(x)-f^{(j)}(y)\big| /\|x-y\|^{\alpha-\lfloor\alpha\rfloor} \leq r \Big\};
    \end{aligned}
\end{equation*}
when $\alpha>0$ and $\m H^{\alpha}_r(\mb R^{d})=\big\{f:\, \mb R^{d}\rightarrow \mathbb{R}:\, \sup_{x\in \mb R^{d}}|f(x)|\leq r \big\}$ when $\alpha=0$.  Additionally,  for any subset $U \in \mb R^{d}$ and a function $f:U \to \mb R$, we say $f\in  {\m H}_r^{\alpha}(U)$ if there exists an extension $\ov f\in   {\m H}^{\alpha}_r(\mb R^{d})$ of $f$ from $U$ to $\mb R^d$, that is, $\ov f|_{U}=f$. For any integer $D>1$, we use ${\m H}_{r,D}^{\alpha}(U)=\big\{f=(f_1,f_2,\cdots,f_D):\, U \to \mb R^D:\, \forall i\in [D], f_i\in   {\m H}_r^{\alpha}(U)\big\}$ to denote the corresponding vector-valued function space. 

There are multiple ways to define a multivariate function with separate smoothness levels across its components. We first introduce a class of smooth multivariate functions, denoted as $\ov{\m H}^{\alpha_1,\alpha_2}_r(\mb R^{d_1},\mb R^{d_2})$, which includes functions that exhibit different marginal smoothness across components. This definition corresponds to the so-called anisotropic function class in the literature~\citep{barron1999risk,nicolas2005adaptive,bhattacharya2014anisotropic}, which we use to characterize our conditional distribution function class, as the marginal smoothness constraint is sufficient for controlling the complexity of the function class through the covering entropy. 
 \begin{definition}\label{def:marginalsmooth}
A function $f:\, \mb R^{d_1}\times \mb R^{d_2}\to \mb R$ belongs to the family $  \ov{\m H}^{\alpha_1,\alpha_2}_r(\mb R^{d_1},\mb R^{d_2})$ if for any $y\in \mb R^{d_2}$, $f(\cdot,y)\in \m H^{\alpha_1}_r(\mb R^{d_1})$ and for any $x\in \mb R^{d_1}$, $f(x,\cdot)\in \m H^{\alpha_2}_r(\mb R^{d_2})$.  
  \end{definition}

Next, we introduce a second, stronger definition of multivariate functions with separate smoothness, denoted as ${\m H}^{\alpha_1, \alpha_2}$, which not only requires marginal smoothness but also imposes constraints on the boundedness of certain mixed partial derivatives for both components. This definition will be used to characterize the covariate-dependent supporting manifold $\m M_{Y|x}$, as it is necessary to ensure that the smoothness definition of the manifold is intrinsic---that is, compatible across different parameterizations; see Remark~\ref{rek:smooth_mani1}  for further details.
 \begin{definition}\label{defsmooth2}
 The class ${\m H}^{\alpha_1,\alpha_2}_r(\mb R^{d_1},\mb R^{d_2})$ consists of all   functions $f:\mb R^{d_1}\times \mb R^{d_2}\to \mb R$ so that
\begin{equation*}
    \begin{aligned}
    &\sum_{(j_1,j_2)\in \m J^{d_1,d_2}_{\alpha_1,\alpha_2}} \sup_{(x,y)\in \mb R^{d_1}\times \mb R^{d_2}} |f^{(j_1,j_2)}(x,y)|  +\sum_{(j_1,j_2)\in \m J^{d_1,d_2}_{\alpha_1,\alpha_2}\atop\frac{|j_1|+1}{\alpha_1}+\frac{|j_2|}{\alpha_2}\geq 1}\underset{x,x_0\in \mb R^{d_1}, y\in \mb R^{d_2}\atop x\neq x_0}{\sup}\frac{| f^{(j_1,j_2)}(x,y)- f^{(j_1,j_2)}(x_0,y)|}{\|x-x_0\|^{\alpha_1-|j_1|-\frac{\alpha_1}{\alpha_2}|j_2|}}\\
  &\qquad+\sum_{(j_1,j_2)\in \m J^{d_1,d_2}_{\alpha_1,\alpha_2}\atop\frac{|j_1|}{\alpha_1}+\frac{|j_2|+1}{\alpha_2}\geq 1}\underset{x\in \mb R^{d_1}, y,y_0\in \mb R^{d_2}\atop y\neq y_0}{\sup}\frac{| f^{(j_1,j_2)}(x,y)- f^{(j_1,j_2)}(x,y_0)|}{\|y-y_0\|^{\alpha_2-|j_2|-\frac{\alpha_2}{\alpha_1}|j_1|}} \leq r,
    \end{aligned}
\end{equation*}
where $\m J^{d_1,d_2}_{\alpha_1,\alpha_2}=\{j_1\in \mb N_0^{d_1}, j_2\in \mb N_0^{d_2}:\,\frac{|j_1|}{\alpha_1}+\frac{|j_2|}{\alpha_2}<1 \}$.
 \end{definition}

\noindent Specifically, when $\alpha_1 = \alpha_2 = \alpha$, the class ${\m H}^{\alpha_1,\alpha_2}_r(\mb R^{d_1},\mb R^{d_2})$ reduces to the classical $\alpha$-smooth H\"{o}lder function class $\m H^{\alpha}_r(\mb R^{d})$ on the joint space $\mb R^{d_1}\times\mb R^{d_2}$. For this reason, we call functions in $\ov{\m H}^{\alpha_1,\alpha_2}_r(\mb R^{d_1},\mb R^{d_2})$ marginally smooth functions, while functions in ${\m H}^{\alpha_1, \alpha_2}$ will be referred to as jointly smooth functions.

The stronger smoothness criterion in Definition~\ref{defsmooth2} requires the existence of mixed derivatives of $f(x,y)$ up to a certain order and enables a local polynomial approximation of $f$ up to certain degree, which is crucial for controlling approximation error when building local polynomial approximations of smooth manifold charts during our estimator construction. Specifically, the following lemma shows that in the vicinity of any point $(x_0,y_0)$, the function $f(x,y)$ can be approximated by a polynomial function with an error of $\m O(\|x-x_0\|^{\alpha_1}+\|y-y_0\|^{\alpha_2})$.

\begin{lemma}[Local polynomial approximation for $ {\m H}^{\alpha_1,\alpha_2}$-smooth functions]\label{le:poly}
    Suppose $f\in  {\m H}^{\alpha_1,\alpha_2}_r(\mb R^{d_1},\mb R^{d_2})$, then there exists a constant $r_1$ so that
     for any $(x_0,y_0)\in \mb R^{d_1}\times\mb R^{d_2}$, it holds for any $(x,y)\in \mb R^{d_1}\times\mb R^{d_2}$ that,
     \begin{equation*}
       \Big| f(x,y)-\sum_{(j_1,j_2)\in \m J^{d_1,d_2}_{\alpha_1,\alpha_2}}  \frac{f^{(j_1,j_2)}(x_0,y_0)}{j_1!j_2!}(x-x_0)^{j_1}(y-y_0)^{j_1}\Big|\leq r_1(\|x-x_0\|^{\alpha_1}+\|y-y_0\|^{\alpha_2}).
     \end{equation*}
\end{lemma}

\noindent The function class ${\m H}^{\alpha_1,\alpha_2}_r(\mb R^{d_1},\mb R^{d_2})$ is closely related to the class $\ov {\m H}^{\alpha_1,\alpha_2}_r(\mb R^{d_1},\mb R^{d_2})$ defined in Definition~\ref{def:marginalsmooth}. On one hand, we have ${\m H}^{\alpha_1,\alpha_2}_r(\mb R^{d_1},\mb R^{d_2})\subset \ov{\m H}^{\alpha_1,\alpha_2}_r(\mb R^{d_1},\mb R^{d_2})$, and this inclusion is strict since marginal differentiability does not imply joint differentiability. On the other hand, the following lemma shows that over a fixed compact set, each function in $\ov {\m H}^{\alpha_1, \alpha_2}_r(\mb R^{d_1},\mb R^{d_2})$ can be approximated by a function in ${\m H}^{\alpha_1, \alpha_2}_{cr(\log \varepsilon)^2}(\mb R^{d_1},\mb R^{d_2})$ for any given error tolerance $\varepsilon>0$, where $c$ is a constant independent of $\varepsilon$.

\begin{lemma}[Relationship between ${\m H}^{\alpha_1,\alpha_2}_r$ and  $ \ov{\m H}^{\alpha_1,\alpha_2}_r$]\label{le:appro}
      Consider an arbitrary function $\bar f\in  \ov{\m H}^{\alpha_1,\alpha_2}_r(\mb R^{d_1},\mb R^{d_2})$, and two compact sets $U_1\subset\mb R^{d_1}$ and $U_2\subset\mb R^{d_2}$, then there exists a constant $c$ so that for any $0<\varepsilon\leq e^{-1}$,  there exists  a function $ f\in  {\m H}^{\alpha_1,\alpha_2}_{cr(\log \varepsilon)^2}(\mb R^{d_1},\mb R^{d_2})$ so that 
        \begin{equation*}
            \underset{x\in U_1,y\in U_2}{\sup} |f(x,y)-\bar f(x,y)|\leq \varepsilon \quad\text{and}  \quad\underset{x\in U_1,y\in U_2\atop (l_1,l_2)\in \m J^{d_1,d_2}_{\alpha_1,\alpha_2}}{\sup} |f^{(l_1,l_2)}(x,y)|\leq cr.  
        \end{equation*}
\end{lemma}

\noindent The approximation property of this lemma allows, in many cases, the two smoothness criteria to be used interchangeably up to a logarithmic term. However, the stronger smoothness condition in Definition~\ref{defsmooth2} is necessary to rigorously define the smoothness of the covariate-dependent supporting manifold $\mathcal{M}_{Y|x}$ through its local charts; see Remarks~\ref{rek:smooth_mani1} in the following subsection for further details.

\subsection{Smooth manifolds and covariate-dependent manifolds}
We focus on distribution regression in settings where both the covariate and the response may exhibit low-dimensional structure. A natural way to describe such structure mathematically is through the \emph{manifold hypothesis}. In its simplest form, this hypothesis asserts that high-dimensional data of interest (including both $X$ and $Y$ in our context) often lie on an unknown $d$-dimensional submanifold $\mathcal{M}$ of $\mathbb{R}^D$, where $d < D$. To formally study distribution regression under the manifold hypothesis, we introduce several key concepts and definitions related to submanifolds in this subsection, which will be used throughout the paper. In particular, we define a class of covariate-dependent manifolds to characterize the support $\mathcal{M}_{Y|x}$ of the response variable $Y$, which may vary with the covariate value $X = x$.

We follow~\cite{10.1214/23-AOS2291,10.1214/18-AOS1685,divol2022measure} in defining the class of regular manifolds. A key quantity that determines the regularity of a manifold, first introduced in~\cite{federer1959curvature}, is the reach $r_{\mathcal{M}}$, defined as
\begin{equation*}
\begin{aligned} 
r_{\mathcal{M}} & :=\sup \Big\{\varepsilon\, \Big|\, \forall x \in  \m M^{\varepsilon}, \text{ there exists unique }y \in \m M, \text { so that } \operatorname{dist}(x, \m M)=\|x-y\|\Big\}, \\ 
&\text { where } \ \ \operatorname{dist}(x, M)=\underset{y\in \m M}{\inf}\|x-y\|, \ \mbox{ and }\  \m M^{\varepsilon}  =\left\{x \in \mathbb{R}^D: \operatorname{dist}(x, \m M)<\varepsilon\right\}.
\end{aligned}
\end{equation*}
The reach $r_{\mathcal{M}}$ quantifies the largest radius of a neighborhood around $\mathcal{M}$ within which every point has a unique projection onto the manifold. A lower bound on the reach (i.e., $r_{\mathcal{M}} \geq \tau > 0$) is crucial, as it prevents the manifold from becoming nearly self-intersecting and ensures a uniform upper bound on its curvature, given by $r_{\mathcal{M}}^{-1} \leq \tau^{-1}$. For a more detailed discussion on the importance of this assumption, we refer the reader to~\cite{10.1214/18-AOS1685}.

Following standard differential geometry texts such as \cite{do2016differential}, the smoothness of a submanifold $\mathcal{M}$ of $\mb R^D$---a manifold embedded in $\mb R^D$---is defined by the smoothness of its local charts. Specifically, for every point $y_0 \in \mathcal{M}$, the manifold $\m M$ can locally be represented as the graph of an $\mathcal{H}^\beta$-smooth, one-to-one function $\phi_{y_0}: V_{y_0} \to \mathbb{R}^D$, where $V_{y_0}$ is an open subset of $\mathbb{R}^d$ containing the origin $\mathbf{0}_d$, and $\phi_{y_0}(\mathbf{0}_d) = y_0$~\citep{10.1214/23-AOS2291}. The pair $( \phi_{y_0}(V_{y_0}),\phi_{y_0}^{-1})$ is referred to as a $\mathcal{H}^\beta$-smooth local chart on $\mathcal{M}$.
In \cite{divol2022measure}, the function $\phi_{y_0}$ is alternatively defined as the inverse of the orthogonal projection ${\rm Proj}_{T{y_0} \mathcal{M}}$ of a local neighborhood of $y_0$ in $\m M$ onto the tangent space $T_{y_0} \mathcal{M}$. Here,  the tangent space $T_{y_0} \mathcal{M}$ is identified with a $d$-dimensional subspace of $\mb R^D$ that pass through the origin, and consists of all vectors tangential to $\mathcal{M}$ at $y_0$.   For precise definitions and additional background on submanifolds and tangent spaces, please refer to Appendix~\ref{sec:ssmanifold}.  These two definitions for smooth submanifolds are shown to be equivalent in Lemma~\ref{le:defmanifold} in Appendix~\ref{app: regularitymanifold}. For clarity and consistency, we adopt the latter definition of the class of $\beta$-smooth submanifolds as described in \cite{divol2022measure} throughout this paper, which is stated as follows.
 
\begin{definition}[$\beta$-Smooth submanifold]\label{def:manifold1}
  A $d$-dimensional submanifold $\m M$ in $\mb R^D$ is said to belong to the manifold class $\ms{M}_{\tau,\tau_1,L}^{\beta}(d,D)$ if: 1.~$\m M$ is closed; 2.~it has reach larger than $\tau$; and 3.~for all $y_0 \in \m M$, there exists a neighborhood $U_{y_0}$ of $y_0$ on $\m M$ so that the projection 
  $\wt\pi_{y_0}:\m M \rightarrow T_{y_0} \mathcal{M}$ defined by $\wt\pi_{y_0}(y)=\operatorname{Proj}_{T_{y_0} \mathcal{M}}(y-y_0)$,  when  {restricted to $U_{y_0}$}, is a   {diffeomorphism}, with  {inverse function} $\phi_{y_0}$ defined on $\mb B_{T_{y_0}{\m M}}(0,\tau_1)$, and  {$\phi_{y_0}\in \m H^{\beta}_{L,D}(\mb B_{T_{y_0}{\m M}}(0,\tau_1))$} (recall that ${\m H}_{L,D}^{\beta}$ denotes the $\beta$-smooth H\"{o}lder class of $\mb R^D$-valued functions with H\"{o}lder norm bounded by $L$).
\end{definition}

Next, we formally define a family of manifolds $\big\{\mathcal{M}_{Y|x}:\,x \in \mathcal{M}_X\big\}$ that is indexed by $x$ on its own support $\m 
M_X$ in the covariate space $\mb R^{D_X}$ and varies smoothly with respect to $x \in \m M_X$. The notion of (joint) smoothness in $(x, y)$ for the family $\big\{\mathcal{M}_{Y|x}\, : x \in \mathcal{M}_X\big\}$ is based on characterizing the joint smoothness of the local charts (which now also depends on $x$) introduced in Definition~\ref{def:manifold1}.
Specifically, for each $(x_0, y_0) \in \mathcal{M}$, we consider the orthogonal projection ${\rm Proj}_{T{y_0} \mathcal{M}_{Y|x_0}}(\cdot - y_0)$. When restricted to a local neighborhood of $y_0$ on $\mathcal{M}_{Y|x}$, this projection should be invertible for each $x$ near $x_0$, provided that the tangent spaces of $\mathcal{M}_{Y|x}$ at points near $y_0$ remain sufficiently aligned with $T_{y_0} \mathcal{M}_{Y|x_0}$. The (joint) smoothness of the manifold family is then defined through the (joint) smoothness of the inverse of this projection in a neighborhood of $(x_0, y_0)$. The precise definition is given below.

\begin{definition}[$(\beta_Y,\beta_X)$-Smooth submanifold family]\label{def:manifold}
A submanifold family $\big\{\mathcal{M}_{Y|x}:\,x \in \mathcal{M}_X\big\}$ is said to belong to $\ms M^{\beta_Y\,\beta_X}_{\tau,\tau_1,L}(d,D,\m M_X)$, if 
for any $x\in \m M_X$: 1.~the manifold $\m M_{Y|x}$ is a closed $d$-dimensional submanifold in  $\mb R^D$; 2.~it has reach larger that $\tau$; and 3.~if, for any $w_0=(x_0,y_0)\in \m M$, there exists a neighborhood $U_{\omega_0}$ of $y_0$ on $\m M_Y$, so that for any $x\in \mb B_{\m M_X}(x_0,\tau)$, the function $\wt \pi_{w_0}: \m M_Y  \rightarrow T_{y_0} \mathcal{M}_{Y|x_0}$ defined by  $\wt \pi_{w_0}(y)=\operatorname{Proj}_{T_{y_0} \mathcal{M}_{Y|x_0}}\left(y-y_0\right)$, when restricted to $U_{\omega_0}\cap \m M_{Y|x}$, is a  diffeomorphism with  inverse function $\phi_{\omega_0,x}(\cdot)$ defined on $\mb B_{T_{y_0}{\m M_{Y|x_0}}}(0,\tau_1)$. Moreover, the function $\Phi_{\omega_0}: \mb B_{T_{y_0}{\m M_{Y|x_0}}}(0,\tau_1)\times  \mb B_{\m M_X}(x_0,\tau)\to \mb R^{D_Y}$ defined as {$\Phi_{\omega_0}(z,x)=\phi_{\omega_0,x}(z)$ belongs to $ {\m H}^{\beta_Y,
 \beta_X}_{L,D_Y}(\mb B_{T_{y_0}{\m M_{Y|x_0}}}(0,\tau_1), \mb B_{\m M_X}(x_0,\tau))$}.
\end{definition}

\begin{remark}\label{rek:smooth_mani1}
When $\beta_Y\geq 2$, and $\beta_Y\geq \beta_X$, assuming the manifold family $\big\{\mathcal{M}_{Y|x}:\,x \in \mathcal{M}_X\big\}$ to be $(\beta_Y, \beta_X)$-smooth is  {equivalent to} assuming the  {existence} of   {$x$-dependent} and  $ {\m H}^{\beta_Y,
 \beta_X}$-smooth local charts to characterize the manifold family. Specifically, this means that for any  point $\omega_0=(x_0,y_0)\in \m M$ and any  $x$ near $x_0$, the manifold $\m M_{Y|x}$ can be locally represented as the graph of a injective  function $\wt g_{\omega_0,x}:  \mb B_{\mb R^{d}}(\mathbf{0},\wt\tau_1)\to \mb R^D$ indexed by $x$, for some positive constant $\wt \tau_1$; in addition, this function changes smoothly in both $x$ and $y$, i.e., the multivariate function $\wt G_{\omega_0}$ defined by $\wt G_{\omega_0}(z,x)=\wt g_{\omega_0,x}(z)$ is ${\m H}^{\beta_Y,
 \beta_X}$-smooth. It's also equivalent to the assumption that locally, the manifold family can be described as set of {solution manifolds} indexed by $x$, with the function $F(y,x)$ that define the equation system being $ {\m H}^{\beta_Y,\beta_X}$-smooth. See Lemma~\ref{le:defmanifold} in Appendix~\ref{app: regularitymanifold} for details.
  
\end{remark}

 \section{Minimax Rate for Distribution Regression with Covariate-independent Response Space}

In this section, we establish the minimax rate of convergence for distribution regression with $n$ i.i.d.~samples $\big\{(X_i,Y_i)\big\}_{i=1}^n$ drawn from $\mu^*_{XY} = \mu^*_X \mu^*_{Y|X}$, under a relatively simpler setting where the support $\mathcal{M}_{Y|x}$ of $\mu^*_{Y|x}$ is \emph{independent} of $x$. Specifically, we assume $\mathcal{M}_{Y|x} = \mathcal{M}_Y$ for all $x \in \mathcal{M}_X$. This setting includes the classical case of density regression when $Y$ is supported on the ambient space $\mathbb{R}^{D_Y}$. We will study the more general, covariate-dependent case in the next section.

We analyze the minimax rate relative to the integral probability metric \citep[IPM,][]{muller1997integral}, which is also called the adversarial loss in the machine learning literature~\citep{singh2018nonparametric,10.1214/23-AOS2291,JMLR:v22:20-911}. Specifically, we consider the following IPM, induced by a H\"{o}lder test function class indexed by a smoothness parameter $\gamma\geq 0$, referred to as the ($\gamma$-)H\"{o}lder IPM,
\begin{equation}\label{def:dgamma}
  d_{\gamma}(\mu,\nu)=\underset{f\in \m H^{\gamma}_1(\mb R^{D_Y})}{\sup} \bigg|\int_{\mb R^{D_Y}} f(y)\,\dd \mu-\int_{\mb R^{D_Y}} f(y)\,\dd \nu\,\bigg|,
\end{equation}
for any two distributions $\mu$ and $\nu$ over $\mb R^{D_Y}$. This metric quantifies the maximum discrepancy in expected test function values between the two distributions $\mu$ and $\nu$, evaluated over test functions from the H\"{o}lder space $\m H^{\gamma}_1(\mb R^{D_Y})$.
The smoothness parameter $\gamma$ controls the strength of the metric. Larger values of $\gamma$ correspond to smoother test functions, which average out local distortions. This makes $d_{\gamma}$ less sensitive to fine details, such as differences in the supports of the distributions, and more responsive to significant global differences in the allocation of probability mass. In contrast, smaller values of $\gamma$ make $d_{\gamma}$ more sensitive to structural changes in the distributions, allowing it to detect subtle variations in shape, such as support misalignment and small bumps in density. Many common probability metrics are special cases of the H\"{o}lder IPM. For example, the 1-Wasserstein distance $W_1$ corresponds to the choice $\gamma = 1$, while the total variation distance $d_{\rm TV}$ corresponds to choosing $\gamma=0$.

To further compare two conditional distributions, such as when evaluating the quality of a conditional distribution estimator $\widehat{\mu}_{Y|X}$ for approximating $\mu^*_{Y|X}$, we adopt the expected H\"{o}lder IPM, i.e., $\mathbb{E}_{\mu^*_X}\big[d_{\gamma}\big(\widehat{\mu}_{Y|X}, \mu^*_{Y|X}\big)\big]$, which takes the expectation with respect to the marginal distribution $\mu^*_X$ over the covariate $X$.
More concretely, we consider two regimes for analyzing the minimax rate of conditional distribution estimation under the expected $d_{\gamma}$ metric. The first regime, referred to as Regime 1, assumes that $\mathcal{M}_Y = \mathbb{R}^{D_Y}$ and that $\mu^\ast_{Y|x}$ is absolutely continuous with respect to the Lebesgue measure on the ambient space. In this case, the response variable $Y$ does not exhibit any low-dimensional manifold structure. The second regime, referred to as Regime 2, assumes that $\mathcal{M}_Y$ is an \emph{unknown}, $\beta_Y$-smooth, $d_Y$-dimensional submanifold with $d_Y < D_Y$, and that $\mu^\ast_{Y|x}$ admits a density with respect to the volume measure on $\mathcal{M}_Y$ (see Appendix~\ref{sec:ssmanifold} for the precise definition).
In both regimes, we allow $\mathcal{M}_X$ to exhibit low-dimensional structure by imposing conditions on its Minkowski dimension, defined below. Recall that for any $\varepsilon > 0$, a set $P \subseteq S$ is called an $\varepsilon$-packing of $S$ if $\|x - x'\| > \varepsilon$ for every pair of distinct points $x, x' \in P$.
\begin{definition}\label{def:Mink} (Covariate space Minkowski dimension)
   We say that a topological space $\m M_X\subset\mb R^{D_X}$ has Minkowski dimension at most $d_X$, or write $\m M_X \in \ms M_X(D_X,d_X,L)$ for some $L>0$, if $\m M_X\in \mb B_{\mb R^{D_X}}(\mathbf{0},L)$, and for any $ 0<\varepsilon \leq 1$,  the maximal cardinality of an $\varepsilon$-packing of $\mathcal{M}_X$ is at most $L \varepsilon^{-d_X}$.
\end{definition}
\noindent This assumption is less restrictive than the manifold assumption in Definition~\ref{def:manifold1}, as it does not impose any conditions on the smoothness or reach of the manifold. In particular, any compact $d_X$-dimensional submanifold of $\mathbb{R}^{D_X}$ has Minkowski dimension (at most) $d_X$.

\subsection{Density regression in Euclidean spaces}
In this subsection, we analyze Regime 1, which corresponds to classical density regression, where the conditional distribution $\mu^*_{Y|x}$ is characterized by a conditional density function $u^*(y\,|\,x)$ with respect to the Lebesgue measure on $\mathbb{R}^{D_Y}$. We further assume that $u^*(y\,|\,x)$ is $\alpha_Y$-smooth in $y$ (marginally) and $\alpha_X$-smooth in $x$, which defines the class of conditional density functions considered below.

\medskip
 \noindent
\emph{\textbf{Regime 1 (Euclidean response space).} For dimensions $D_Y,D_X\in \mb N_{+}$, $d_X\in \mb N\cap [0,D_X]$, smoothness parameters $\alpha_Y,\alpha_X\in(0,\infty)$, and a constant $L>0$, we define the distribution family $\m P^*_1=\m P^*_1(D_Y,D_X,d_X,\alpha_Y,\alpha_X,L)$
that consists of all joint distributions $\mu^*=\mu^*_X\mu^*_{Y|X}=\mu^*_Y\mu^*_{X|Y}$ so that 
\begin{enumerate}
    \item The support $\m M_X$ of $\mu^*_X$ belongs to the family $\ms M_X(D_X,d_X,L)$ and the support $\m M_Y$ of $\mu^*_Y$ is a subset of $\mb B_{\mb R^{D_Y}}(\mathbf{0},L)$. 
    \item  For any $x\in \m M_X$, $\mu^*_{Y|x}$ has a density function $u^*(\,\cdot\,|\,x)$ with respect to Lebesgue measure on $\mb R^{D_Y}$, and $u^*(y\,|\,x)\in  \ov{\m H}^{\alpha_Y,\alpha_X}_L(\mb R^{D_Y},\m M_X)$.
\end{enumerate}}
 
\medskip
\noindent We also allow $d_X = 0$ in the above definition, which corresponds either to unconditional distribution estimation or to settings where the covariate $X$ is discrete and takes finitely many values. The assumption that $\mu^*_Y$ is compactly supported is made primarily for technical convenience. However, the analysis can be extended to cases with non-compact support, provided that $u^*(\,\cdot \,|\, x)$ exhibits sufficiently light tails (e.g., exponential decay) for every $x \in \m M_X$. In such cases, it is sufficient to restrict the analysis to a compact region with radius on the order of $\m O(\log n)$. Additionally, we assume that the density $\mu^*(y\,|\,x)$ decays $\alpha_Y$-smoothly to zero near the boundary of ${\rm supp}(\mu^*_{Y|x})$, although the exact boundary need not be known and is allowed to vary with $x$. The following theorem summarizes our result on the minimax rate for estimating the family of conditional distributions $\{\mu^*_{Y|X=x}:\,x\in \m M_X\}$ under this regime. A proof of the theorem is provided in Appendix~\ref{APP:proofEuclidean}.

\begin{theorem}[Minimax rate under Regime 1]\label{th1:lower} For each $\gamma\geq 0$, there exist a constant $L_0$ so that when $L,n\geq L_0$, it holds that
\begin{equation*}
    \begin{aligned}
        C\,\bigg( n^{-\frac{\alpha_X}{2\alpha_X+d_X}}\ +\ & n^{-\frac{\alpha_Y+\gamma}{2\alpha_Y+D_Y+\frac{\alpha_Y}{\alpha_X}d_X}}\,  \bigg)\leq \underset{\wh \mu_{Y|X}}{\inf}\,\underset{\mu^*=\mu^*_X \mu^*_{Y|X} \in \m P^*_1}{\sup}\mb{E}_{\mu^{*,\otimes n}}\Big[\mb{E}_{\mu^*_X}\big[d_{\gamma}\big(\mu^*_{Y|X},\,\wh \mu_{Y|X}\big)\big]\Big]\\
&\qquad\leq C_1\,\bigg(  \sqrt{\log n}\cdot \Big(\frac{n}{\log n}\Big)^{-\frac{\alpha_X}{2\alpha_X+d_X}}\ +\ \Big(\frac{n}{\log n}\Big)^{-\frac{\alpha_Y+\gamma}{2\alpha_Y+D_Y+\frac{\alpha_Y}{\alpha_X}d_X}}\,\bigg),       \end{aligned}
    \end{equation*}
 where $(C,C_1)$ are constants independent of $n$, and the infimum is taken over all conditional distribution estimators $\wh \mu_{Y|X}$ based on data $\big\{(X_i,Y_i)\big\}_{i=1}^n$ sampled from $\mu^{*,\otimes n}$. The shorthand $\m P^*_1$ stands for $\m P^*_1(D_Y,D_X,d_X,\alpha_Y,\alpha_X,L)$.
\end{theorem}
\noindent Here the assumption $L\geq L_0$ is used for deriving the minimax lower bound. The proof involves constructing distributions that are difficult to distinguish and applying reduction techniques to transform the estimation problem into a multiple testing problem. The constant $L_0$ serves as a threshold ensuring that the constructed distributions satisfy the assumptions of Regime 1. We observe a phase transition in the minimax convergence rate as the parameter $\gamma$ varies. When $\gamma \geq \frac{d_Y \alpha_X}{2\alpha_X + d_X}$, the dominant term in the rate is $\wt{\m O}\big(n^{- \alpha_X/(2\alpha_X + d_X)}\big)$, which matches the classical minimax rate for estimating an $\alpha_X$-smooth regression function under $L_2$ loss~\citep{stone1982optimal}. This is because smoother test functions average out local fluctuations in the conditional density, making the metric $d_\gamma$ primarily responsive to the overall dependence trend, such as the conditional mean of $Y$ given $X$. In this regime, the complexity is governed solely by the smoothness $\alpha_X$ and intrinsic dimension $d_X$ of the covariate $X$.
In contrast, when $\gamma \leq \frac{d_Y \alpha_X}{2\alpha_X + d_X}$, the metric $d_\gamma$ becomes more sensitive to local features of the density, and the dominant term in the minimax rate becomes $\wt{\m O}\big(n^{- (\alpha_Y + \gamma)/(2\alpha_Y + D_Y + \frac{\alpha_Y}{\alpha_X} d_X)}\big)$. This rate reflects the intrinsic difficulty of nonparametric conditional density estimation, and improves as either the smoothness of the conditional density increases or the intrinsic dimensions $d_X$ and $D_Y$ decrease. The rate also improves with larger values of $\gamma$, as the metric gradually shifts its sensitivity from local irregularities toward global structural differences.

A related work by \cite{10.1214/23-AOS2270} studies the minimax rate for conditional density estimation under the Kullback-Leibler (KL) risk, providing both upper and lower bounds expressed in terms of empirical Hellinger entropy. Under the assumption that $\m M_X$ and $\m M_Y$ are unit cubes in $\mb R^{D_X}$ and $\mb R^{D_Y}$, respectively, and that the partial derivatives $(u^*)^{(j_1, j_2)}(y \mid x)$ exist and are bounded for all multi-indices $j_1 \in \mb N_0^{D_Y}$ and $j_2 \in \mb N_0^{D_X}$ with $|j_1| \leq \alpha_Y$ and $|j_2| \leq \alpha_X$, they derive an upper bound of $\wt{\m O}\big(n^{-{\alpha_Y}/(\alpha_Y + D_Y + \frac{\alpha_Y}{\alpha_X}D_X)}\big)$ for the KL risk when $\alpha_X, \alpha_Y \in \mb N_{+}$. This result further implies an upper bound of $\wt{\m O}\big(n^{-{\alpha_Y}/(2\alpha_Y + 2D_Y + 2\frac{\alpha_Y}{\alpha_X}D_X)}\big)$ for the expected total variation distance, via Pinsker's inequality.
In contrast, by setting $\gamma = 0$ in Theorem~\ref{th1:lower}, we obtain a sharper upper bound of $\wt{\m O}\big(n^{-{\alpha_Y}/(2\alpha_Y + D_Y + \frac{\alpha_Y}{\alpha_X} d_X)}\big)$ for the expected total variation distance, along with a matching lower bound. Our result further accommodates low-dimensional structure in the covariate space $\m M_X$ and relies on a weaker smoothness assumption that does not require the existence of mixed partial derivatives of order up to $(\alpha_X + \alpha_Y)$.

In another line of work, \cite{10.1214/22-EJS2037} show that for $\alpha_X \in [0,1]$, a properly designed kernel estimator achieves the minimax rate under the expected total variation distance. Our result extends beyond this setting, allowing for general $\alpha_X > 0$ and covering a broader class of metrics $\{d_\gamma:\,\gamma\geq 0\}$.
Finally, \cite{tang2024conditionaldiffusionmodelsminimaxoptimal} study the convergence rate of conditional diffusion models~\citep{song2020score, batzolis2021conditional, NEURIPS2021_cfe8504b} under the expected 1-Wasserstein distance, which corresponds to $d_\gamma$ with $\gamma = 1$. Their derived upper bound, up to logarithmic factors, matches ours in Theorem~\ref{th1:lower} for $\gamma = 1$, although they do not provide a matching lower bound. When combined with our minimax lower bound, their result implies that conditional diffusion models are minimax optimal under the expected $1$-Wasserstein metric.

\subsection{Distribution regression with low-dimensional manifold structures}
In this subsection, we consider the regime where the response space $\m M_Y$ is an \emph{unknown} $\beta_Y$-smooth submanifold of intrinsic dimension $d_Y$, embedded in the ambient space $\mb R^{D_Y}$. The conditional distribution $\mu^*_{Y|x}$ is characterized by a density function $u^*(y\mid x)$ defined with respect to the \emph{volume measure} on $\m M_Y$. We assume that $u^*$ exhibits marginal smoothness of order $\alpha_Y$ in the $y$-component and $\alpha_X$ in the $x$-component (c.f.~Definition~\ref{def:marginalsmooth}). We refer to this setting as ``distribution regression" rather than ``density regression", since $\mu^*_{Y|x}$ is not absolutely continuous with respect to the Lebesgue measure on the ambient space $\mathbb{R}^{D_Y}$, nor with respect to any known base measure, due to the supporting manifold $\mathcal{M}_Y$ being unknown. A formal definition of this regime is given below.

\medskip
\noindent\emph{
\textbf{Regime 2 (Covariate-independent manifold response space).} For dimensions $D_Y,d_Y,D_X\in \mb N_{+}$, $d_X\in\mb N_0$, smoothness parameters $\beta_Y,\alpha_Y,\alpha_X>0$,  a function $g:\mb R^{+}\to \mb R^{+}$, and absolute constants $\tau,\tau_1,L>0$, we define the following distribution family 
\begin{equation*}
\begin{aligned}
& \m P^*_2=\m P^*_2(D_Y,D_X,d_Y,d_X,\beta_Y,\alpha_Y,\alpha_X,\tau,\tau_1,g,L),
 \end{aligned}
\end{equation*}
 which consists of all $\mu^*=\mu^*_X\mu^*_{Y|X}$ so that 
 \begin{enumerate}
     \item The support $\m M_X$ of $\mu^*_X$ belongs to $\ms M_X(D_X,d_X,L)$.
     \item For any $x\in \m M_X$, the conditional distribution $\mu^*_{Y|x}$ is supported on a submanifold $\m M_{Y}$ independent of $x$, and  has a density function $u^*(\,\cdot\,|\,x)$ with respect to the volume measure of $\m M_{Y}$, where $\m M_{Y}\in \ms M^{\beta_Y}_{\tau,\tau_1,L}(d_Y,D_Y)$ and $u^*\in  \ov{\m H}^{\alpha_Y,\alpha_X}_{L}(\m M_Y,\m M_X)$.
     \item For any $x_0\in \m M_X$, $y_0\in \m M_{Y|x_0}$ and $0<r\leq 1$,  the measure $\mu^*_X$ on the ball $\mb B_{\m M_X}(x_0,r)$ is bounded below by $g(r)r^{d_X}/L$, and the measure $\mu^*_{Y|x_0}$ of $Y$ given $X=x_0$,  on the ball $\mb B_{\m M_{Y|x}}(y_0,r)$, is bounded below by $ g(r)r^{d_Y}/L$.
      \end{enumerate}
}

 \medskip
\noindent The function $g(\cdot)$ in Item 3 of Regime 2 is introduced for technical purposes, serving to control over the constant term that captures the decay behavior of $u^*(\,\cdot \,|\, x)$ when taking the supremum over all measures in $\mathcal{P}_2^*$.  Setting $g(r) \equiv 1$ corresponds to the case where $u^*(y\,|\,x)$ is uniformly bounded away from zero for any $y\in \m M_{Y|x}$ and $x\in \m M_X$. However, our framework accommodates greater generality by requiring only that  $g(r) > 0$ for all $r > 0$. As an illustrative example, consider the distribution $\mu^*(y\,|\,x)=\mu^*(y)=G_{\#}\nu$  supported on the unit sphere $\mb S_1$, where $G(\theta)=(\sin(\pi\theta),\cos(\pi\theta))$, and $\nu$  is a probability measure with density $v(\theta)\propto\theta^2(1-\theta)^2\mathbf{1}(0<\theta<1)+\theta^2(\theta+1)^2\mathbf{1}(-1<\theta<0)$. It can be shown that the density  of $\mu^*$ with respect to the volume measure on  $\mb S_1$ is given by $u^*(y_1,y_2)\propto \arccos(y_1)^2(\pi-\arccos(y_1))^2$, which is uniformly Lipschitz continuous but not bounded away from zero. Nonetheless, by choosing $g(r)=r^2$, the inequality $\mu^*(\mb B_{\mb S_1}(y,r))>\frac{0.15}{\pi^3} r g(r)$  holds for every $y\in\mb S_1$ and $0<r\leq 1$. A similar argument applies to $\mu^*_X$, where we likewise do not require the measure to admit a density uniformly bounded away from zero. We are now prepared to present our result on the minimax rate of convergence under Regime 2.

\begin{theorem}[Minimax rate under Regime 2]\label{th2:lower1} For each $\gamma>0$, if $\beta_Y\geq 2 \vee (\alpha_Y+1)$, then there exists a constant $L_0$ so that when $L,\tau,\tau_1,n\geq L_0$, it holds that
    \begin{equation*}
    \begin{aligned}
       C\,\bigg(n^{-\frac{\alpha_X}{2\alpha_X+d_X}}\ +\ & n^{-\frac{\alpha_Y+\gamma}{2\alpha_Y+d_Y+\frac{\alpha_Y}{\alpha_X}d_X}} \ + \  n^{-\frac{\gamma}{\frac{d_Y}{\beta_Y}}}\,\bigg) \\
       &\leq \underset{\wh \mu_{Y|X}}{\inf}\,\underset{\mu^*=\mu^*_X\mu^*_{Y|X}\in \m P^*_2}{\sup}\mb{E}_{\mu^{*,\otimes n}}\Big[\mb{E}_{\mu^*_X}\big[d_{\gamma}\big(\mu^*_{Y|X},\, \wh \mu_{Y|X}\big)\big]\Big]\\
&\qquad \leq C_1\bigg((\log n)^3\cdot n^{-\frac{\alpha_X}{2\alpha_X+d_X}}\ + \ \Big(\frac{n}{\log n}\Big)^{-\frac{\alpha_Y+\gamma}{2\alpha_Y+d_Y+\frac{\alpha_Y}{\alpha_X}d_X}}\ +\ n^{-\frac{\gamma}{\frac{d_Y}{\beta_Y}}}\,\bigg),
    \end{aligned}
    \end{equation*}
    where $(C,C_1)$ are constants independent of $n$, and the infimum is taken over all conditional distribution estimators $\wh \mu_{Y|X}$ based on data $\big\{(X_i,Y_i)\big\}_{i=1}^n$ sampled from $\mu^{*,\otimes n}$. The shorthand $\m P^*_2$ stands for $\m P^*_2(D_Y,D_X,d_Y,d_X,\beta_Y,\alpha_Y,\alpha_X,\tau,\tau_1,g,L)$. 
\end{theorem}

Given that $\m M_Y$ is unknown and we only observe $n$ i.i.d.~samples, the estimator $\wh \mu_{Y|x}$ and the true conditional distribution $\mu^*_{Y|x}$ are almost surely mutually singular, as they are supported on different submanifolds. Consequently, the total variation distance between them is identically $1$ and fails to meaningfully reflect distributional closeness. To address this issue, we restrict attention to metrics $d_\gamma$ with $\gamma > 0$ in this regime. The condition $\beta_Y \geq 2$ ensures that the submanifold $\m M_Y$ has bounded curvature, while the assumption $\beta_Y \geq \alpha_Y + 1$ guarantees that the smoothness parameter $\alpha_Y$ is compatible and invariant to the choice of the local charts of the manifold $\m M_Y$. For further discussion, see \cite{10.1214/23-AOS2291}. Compared to Theorem~\ref{th1:lower}, the minimax rate in Theorem~\ref{th2:lower1} contains an additional term $n^{-\gamma\beta_Y/d_Y}$, which reflects the intrinsic difficulty of estimating the unknown submanifold $\m M_Y$ from i.i.d.~samples $\{Y_i\}_{i=1}^n$ drawn on the manifold. Moreover, in settings where $\mu_X^*$ is discrete (i.e., $d_X = 0$) or where $Y$ is independent of $X$ (corresponding to the limiting case of $\alpha_X \to \infty$), the minimax rate simplifies to $n^{-1/2} + n^{-(\alpha_Y+\gamma)/(2\alpha_Y + d_Y)} + n^{-\gamma\beta_Y/d_Y}$, which recovers the rate for unconditional distribution estimation on unknown submanifolds obtained in \cite{10.1214/23-AOS2291}.

\begin{figure}[H]
     \centering
     \begin{subfigure}{0.45\textwidth}
         \centering
         \includegraphics[width=\textwidth]{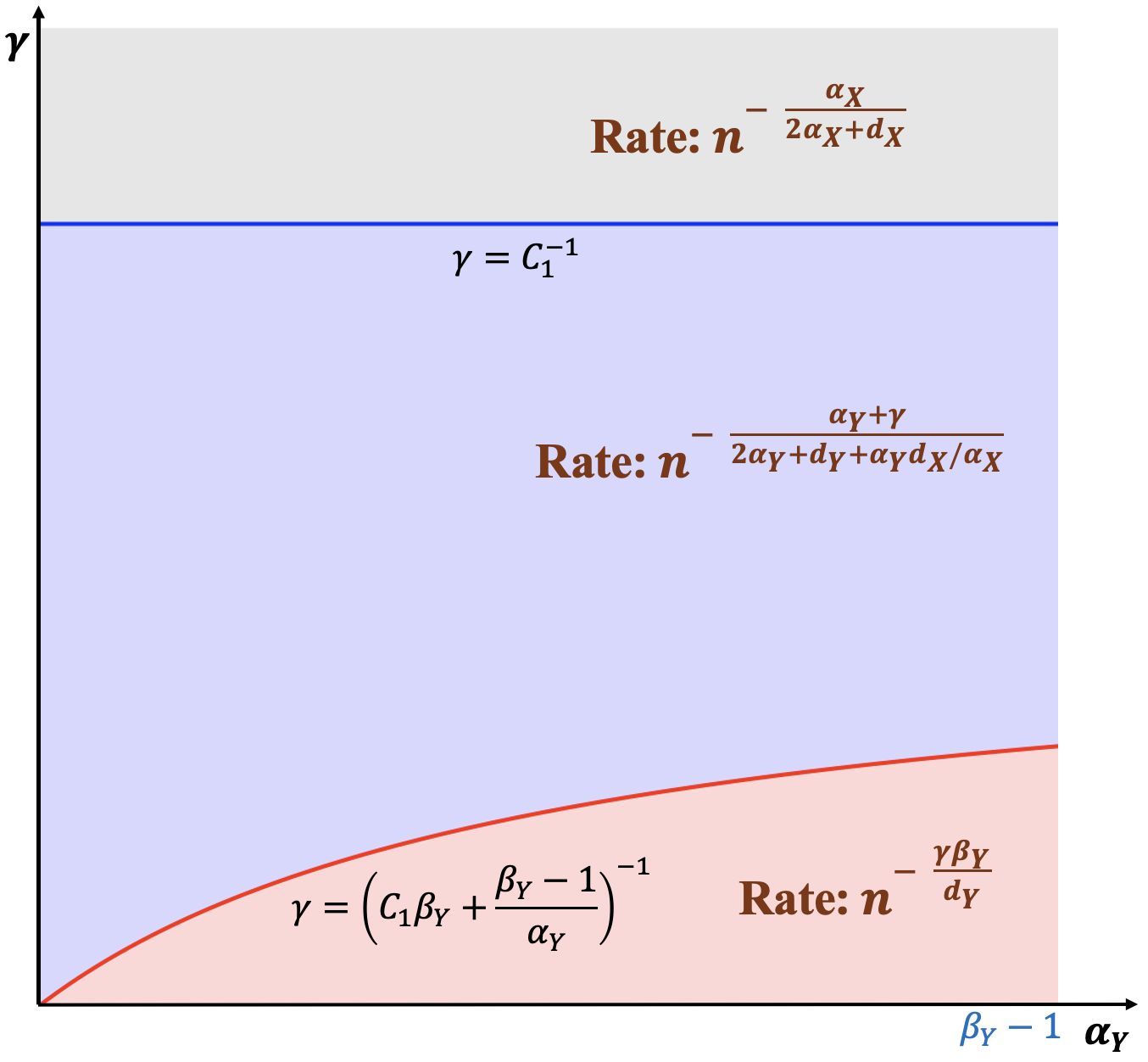}
         \caption{Dominant term in the minimax rate for varying $\alpha_Y$ and $\gamma$, where $C_1=\frac{2\alpha_X+d_X}{d_Y\alpha_X}$. }
     \end{subfigure}
     \hfill
     \begin{subfigure}{0.45\textwidth}
         \centering
         \includegraphics[width=\textwidth]{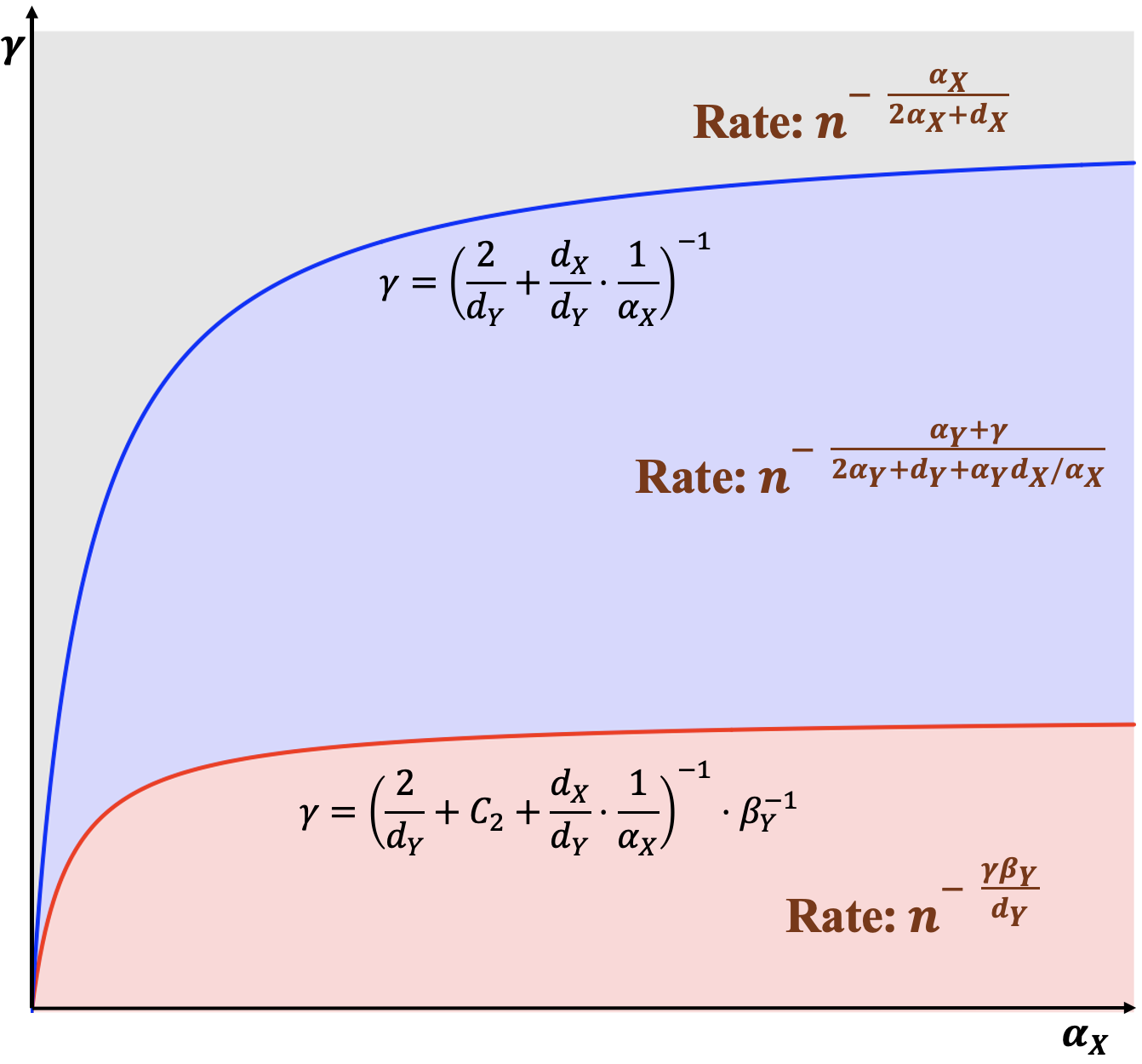}
         \caption{Dominant term in the minimax rate for varying $\alpha_X$ and $\gamma$, where $C_2=\frac{\beta_Y-1}{\beta_Y\alpha_Y}$.}
         
     \end{subfigure}
        \caption{Diagram for the minimax rate under Regime 2 for fixed $d_X\in \mb N$, $d_Y\in \mb N^{+}$ and $\beta_Y\geq 2$.}
        \label{fig:1}
\end{figure}

Figure~\ref{fig:1} illustrates the three regimes of problem characteristics identified in Theorem~\ref{th2:lower1}, based on varying values of $\alpha_X$, $\alpha_Y$, and $\gamma$. Each regime is determined by which of the three terms in the minimax rate dominates. The diagram reveals two critical transition points for $\gamma$: $\gamma =\frac{d_Y \alpha_Y}{2\alpha_Y\beta_Y + d_Y(\beta_Y - 1) + d_X \beta_Y \alpha_Y/\alpha_X}$ and $\gamma=\frac{d_y\alpha_X}{2\alpha_X+d_X}$. The first transition occurs between two dominant error regimes. For smaller values of $\gamma$, the error is governed by support estimation with rate $n^{-{\gamma\beta_Y}/d_Y}$. As $\gamma$ increases, the dominant term shifts to that of nonparametric conditional density estimation, with error rate $n^{-(\alpha_Y+\gamma)/(2\alpha_Y + d_Y + \frac{\alpha_Y}{\alpha_X}d_X)}$. This transition reflects a key sensitivity of the $d_\gamma$ metric: for small $\gamma$, it is more responsive to support misalignment than to discrepancies in mass allocation across the support. Consequently, support estimation dominates when $\gamma$ is small.  Moreover, this transition point increases with larger $\alpha_X$ and $\alpha_Y$, indicating that higher smoothness in the covariate or response reduces the complexity of density recovery. As a result, a larger $\gamma$ is needed to render support estimation errors negligible in comparison to those in density estimation. The second transition point marks the shift from conditional density estimation for smaller $\gamma$ values to global dependence recovery, characterized by the rate $n^{-\frac{\alpha_X}{2\alpha_X+d_X}}$. This threshold depends on $\alpha_X$ but not on $\alpha_Y$, and it increases with larger $\alpha_X$. A higher $\alpha_X$ reduces the difficulty of capturing the dependence between $X$ and $Y$, thereby requiring a smoother test function (i.e., larger $\gamma$) to adequately smooth out local variations in the conditional distribution.

\section{Minimax Rate for Distribution Regression with Covariate-dependent Response Space}
In this section, we investigate a more complex setting where the support $\m M_{Y|x}$ of the conditional distribution $\mu^*_{Y|x}$ depends on the covariate $x$. This additional flexibility requires estimating not a single submanifold $\m M_Y$, but a family of submanifolds $\big\{\m M_{Y|x} :\, x \in \m M_X\big\}$ indexed by $x$. We refer to this task as \emph{manifold regression}, where the goal is to use the data $\big\{(X_i, Y_i)\big\}_{i=1}^n$ i.i.d.~drawn from $\mu^\ast = \mu^\ast_X \mu^\ast_{Y|X}$, to estimate or predict the submanifold $\m M_{Y|x}$, which serves as the support of $\mu^\ast_{Y|X = x}$, for any given $x$. We begin by introducing the formal setup and deriving the minimax rate for manifold regression. We then extend the analysis to obtain the minimax rate for distribution regression in this more general setting with covariate-dependent supporting manifold.

\subsection{Manifold regression} \label{sec:manifold_reg}

\noindent Recall that we observe i.i.d.~data $\big\{(X_i, Y_i)\big\}_{i=1}^n$ drawn from a joint distribution $\mu^* = \mu^*_X \mu^*_{Y|X}$, where the conditional distribution $\mu^*_{Y|x}$ has support $\m M_{Y|x}$.  In this subsection, our goal is to analyze the minimax rate for estimating the family of submanifolds $\big\{\mathcal{M}_{Y|x} : x \in \mathcal{M}_X\big\}$ based on the observed data, under the assumption that this family is $(\beta_Y, \beta_X)$-smooth (c.f.~Definition~\ref{def:manifold}). This problem  is highly relevant to various real-world applications. For example,  consider the face image data $Y$  conditioned on specific attributes $X$ such as age and gender~\cite{antipov2017face,lu2018attribute,ding2021ccgan}.  For a given  value of $X$, it is reasonable to assume that the image dataset lies in (or close to) a submanifold~\citep{wang2008manifold}, while different $X$ values  may correspond to distinct manifolds. For instance, the face image dataset for age $18$ might be quite different 
 from that for age $80$, making it reasonable to model these as two distinct submanifolds.  It is worth noting that when $\m M_{Y|x}=\m M_Y$ for any $x\in \m M_X$, the problem reduces to (single) manifold estimation,  a topic previously explored in various literature~\citep{10.1214/18-AOS1685,genovese2012minimax, JMLR:v23:21-0338}. Therefore, our framework can be viewed as an extension of these prior works to the conditional setting under a noiseless model. We measure the estimation error using the maximal Hausdorff distance evaluated over the covariate space $\m M_X$, defined as ${\sup}_{x \in \m M_X} \mb H\big(\m M_{Y|x}, \wh{\m M}_{Y|x}\big)$, where the Hausdorff distance $\mb H(\m M_1, \m M_2)$ between two sets $\m M_1$ and $\m M_2$ is defined as $\sup_{x\in \m M_1}\inf_{y\in \m M_2}\|x-y\|+\sup_{x\in \m M_2}\inf_{y\in \m M_1}\|x-y\|$. The Hausdorff distance $\mb H$ is commonly used to evaluate errors in manifold estimation~\citep{10.1214/18-AOS1685,genovese2012minimax}. Our analysis will be carried out over a class $\m P^*$ of distributions $\mu^*$ defined as follows.


\medskip
\noindent\emph{\textbf{Regime 3a (Manifold regression).} For dimensions $d_X,D_X,d_Y,D_Y\in \mb N_{+}$, smoothness parameters $\beta_Y$, $\beta_X>0$, and absolute constants $\tau,\tau_1,L>0$,  we define the following distribution family 
$\m P^*=\m P^*(D_Y,D_X,d_Y,d_X,\beta_Y,\beta_X,\tau,\tau_1,L)$, 
which consists of all $\mu^*=\mu^*_X\mu^*_{Y|X}$ so that 
\begin{enumerate}
    \item  $\mu^*_X$ has a support $\m M_X\in \m M^{\beta_X \vee 2}_{\tau,\tau_1,L}(d_X,D_X)$ and has a density $u_X$ function with respect to the volume measure on $\m M_X$ such that $1/L\leq u^*_X(x)\leq L$ for any $x\in \m  M_X$. 
    \item For any $x\in \m M_X$, the conditional distribution $\mu^*_{Y|x}$ is supported on a manifold $\m M_{Y|x}$, and admits a density function $u^*(\cdot\,|\,x)$  with respect to the volume measure on $\m M_{Y|x}$ so that $1/L\leq u^*(y|x)\leq L$ for any $y\in \m  M_{Y|x}$,  and $\{\m M_{Y|x}\,:\, x\in \m M_X\}\in\ms M^{\beta_Y\,\beta_X}_{\tau,\tau_1,L}(d_Y,D_Y,\m M_X)$.
  \end{enumerate}
}

\noindent Different from Regimes 1 and 2, here we no longer impose any smoothness conditions on the conditional density function $u^*(y\,|\,x)$, since the goal is to recover the support of $\mu^\ast_{Y|X}$. However, we require the covariate space $\m M_X$ to be a smooth submanifold, as the regularity of $\m M_X$ facilitates the control of the ``worst-case'' sense error in terms of Hausdorff distance through a localized mean squared error, simplifying the problem to controlling 
 an ``average'' sense error.  Note that in the subsequent subsection, where we focus on estimating the conditional distribution and the error metric is directly defined in an ``average'' sense (rather than a worst-case one), this stronger assumption on $\m M_X$ can be relaxed to requiring only that $\m M_X$ has bounded upper Minkowski dimension, as specified in Definition~\ref{def:Mink}.   Moreover, here we requires the density function $u^*_X$ to be bounded away from zero to ensure that there are sufficiently many samples around each $x\in \m M_X$, which is crucial for  controlling the maximal Hausdorff distance. We conjecture that this condition could be relaxed by considering an average Hausdorff distance, for example, $\mb E_{\mu^*_X}\big[\mb H\big(\m M_{Y|x},\,\wh {\m M}_{Y|x}\big)\big]$. With these assumptions in place, we are now ready to present our main result on the minimax rate of convergence for manifold regression. The proof is provided in Appendix~\ref{prooftheorem1}.

\begin{theorem}[Minimax rate for manifold regression]\label{th:mainifold} Suppose $\beta_Y\geq 2$ and $\beta_Y\geq \beta_X$, then there exists a constant $L_0$ so that when $L,\tau,\tau_1,n\geq L_0$, it holds that
    \begin{equation*}
    \begin{aligned}
       &C\, n^{-\frac{1}{\frac{d_Y}{\beta_Y}+\frac{d_X}{\beta_X}}}   \  \leq \underset{\wh {\m M}_{Y|x},\,x\in\m M_X}{\inf}\,\underset{\mu^*\in \m P^*}{\sup}\mb{E}_{\mu^{*,\otimes n}}\Big[\underset{x\in \m M_X}{\sup}\mb H\big(\m M_{Y|x},\,\wh {\m M}_{Y|x}\big)]\Big]\leq \  C_1 \Big(\frac{n}{\log n}\Big)^{-\frac{1}{\frac{d_Y}{\beta_Y}+\frac{d_X}{\beta_X}}},
    \end{aligned}
    \end{equation*}
     where $(C,C_1)$ are constants independent of $n$,  and the infimum is taken over all estimators $\{\wh {\m M}_{Y|x}:\,x\in\m M_X\}$ based on data $\big\{(X_i,Y_i)\big\}_{i=1}^n$ sampled from $\mu^{*,\otimes n}$. Here, the shorthand $\m P^*$ stands for $\m P^*(D_Y,D_X,d_Y,d_X,\beta_Y,\beta_X,\tau,\tau_1,L)$. 
\end{theorem}

Compared to the minimax rate for estimating a $\beta_Y$-smooth, $d_Y$-dimensional submanifold~\citep{10.1214/18-AOS1685}, which is $n^{-1/(d_Y/\beta_Y)}$, our rate includes an additional term $d_X/\beta_X$ in the denominator of the exponent. This reflects the increased statistical complexity of estimating an entire family of submanifolds $\big\{\mathcal{M}_{Y|x}:\,x \in \mathcal{M}_X\big\}$, rather than a single submanifold $\m M_{Y|x}$. Our results indicate that higher smoothness $\beta_X$ with respect to the covariate $X$ makes the manifold regression problem easier, as it facilitates information sharing across different covariate values, which in turn leads to faster convergence of the minimax rate.  Our estimator generalizes the local polynomial estimator from~\citep{10.1214/18-AOS1685} by incorporating the covariate $X$. Specifically,  for each data point $w_k=(X_k, Y_k)$,  we select nearby data samples $(X,Y)$  such that $\|Y-Y_k\|\leq h_1$ and $\|X-X_k\|\leq h_2$, where $h_1\asymp ({\log n}/{n})^{\frac{\beta_X}{d_Y\beta_X+{d_X\beta_Y}}}$ and $h_2\asymp({\log n}/{n})^{\frac{\beta_Y}{ {d_Y\beta_X}+d_X\beta_Y}}$. We then learn a local polynomial estimator by minimizing the average reconstruction loss between $Y$ and $G(Q(Y),X)$, where $Q(\cdot)=V^T(\cdot-Y_k)$ with $V$ being  $D_y\times d_Y$ orthonormal matrices targeting  one of the orthonormal basis $V^*_{k}$ of the tangent space $T_{Y_k}\m M_{Y|X_k}$. The function $G(\cdot,X)$ consists of  polynomial functions designed to approximate $\Phi_{w_k}(V^*_{k}z, x)$, where $\Phi_{w_k}(\cdot,X)$ is the inverse of ${\rm Proj}_{T_{Y_k}\m M_{Y|X_k}}(\cdot-Y_k)$ when restricted to the neighborhood of $w_k$ on $\m M_{Y|X}$, as defined in Definition~\ref{def:manifold}.  The assumption $\beta_Y\geq \beta_X$ ensures that $h_2\leq h_1$,  allowing us to establish the equivalence between the distance $\|Y-Y_k\|$  and the distance of the projections $\|V^*_k(Y-Y_k)\|$, up to multiplicative constants. This equivalence enables the analysis to be carried out in the low-dimensional coordinates $V^*_k{}^T(Y-Y_k)$ by employing polynomial approximations of the $\m H^{\beta_Y,\beta_X}$ smooth functions $\Phi_{w_k}(z, x)$ within the regions $\|z\|\leq h_1$ and $\|x\|\leq h_2$. Similar to~\citep{10.1214/18-AOS1685}, the final estimator  is then  constructed by assembling a union of polynomial patches. Further details of the estimator are provided in Appendix~\ref{prooftheorem1}.



\subsection{Distribution regression with covariate-dependent manifolds}
In this subsection, we study the problem of distribution regression under the setting where the conditional response supports $\big\{\mathcal{M}_{Y|x} :\, x \in \mathcal{M}_X\big\}$ form an unknown family of submanifolds that is $(\beta_Y, \beta_X)$-smooth (c.f.~Definition~\ref{def:manifold}). We still use $u^*(y \,|\, x)$ to denote conditional density function of the conditional distribution $\mu^*_{Y|x}$ with respect to the volume measure on its (covariate-dependent) supporting manifold $\mathcal{M}_{Y|x}$. Due to the variability in the response space and its associated volume measure across different values of $x$, we employ the  stronger smoothness criteria $\m H^{\alpha_Y,\alpha_X}$ defined in Definition~\ref{defsmooth2} to quantify the smoothness of $u^*$, and we will discuss its implications later in Remark~\ref{Resmoothness}.  A formal definition of this regime is presented below.


\medskip
\noindent\emph{\textbf{Regime 3b (Covariate-dependent manifold response space).} For dimensions $D_Y,d_Y,D_X,d_X\in \mb N_{+}$,  smoothness parameters $\beta_Y,\beta_X,\alpha_Y,\alpha_X>0$,  a function $g:\mb R^{+}\to \mb R^{+}$, and absolute constants $\tau,\tau_1,L>0$, we define the following distribution family 
\begin{equation*}
\begin{aligned}
& \m P^*_3=\m P^*_3(D_Y,D_X,d_Y,d_X,\beta_Y,\beta_X,\alpha_Y,\alpha_X,\tau,\tau_1,g,L),
 \end{aligned}
\end{equation*}
 which consists of all $\mu^*=\mu^*_X\mu^*_{Y|X}$ so that 
 \begin{enumerate}
     \item The supporting manifold $\m M_X$ of $\mu^*_X$ belongs to $\ms M_X(D_X,d_X,L)$.
\item  For any $x\in \m M_X$, the conditional distribution $\mu^*_{Y|x}$ supported on a submanifold $\m M_{Y|x}$ and has a density function $u^*(\,\cdot\,|\,x)$ with respect to the volume measure of $\m M_{Y|x}$ so that $\big\{\m M_{Y|x}:\,x\in \m M_X\big\}\in \ms M^{\beta_Y\,\beta_X}_{\tau,\tau_1,L}(d_Y,D_Y,\m M_X)$ and  there exists a function $\ov u^*\in  {\m H}^{\alpha_Y,\alpha_X}_{L}(\mb R^{D_Y},\mb R^{D_X})$  so that $u^*(y|x)=\ov u^*(y,x)$ for any $(x,y)\in \m M$.
     \item For any $x_0\in \m M_X$, $y_0\in \m M_{Y|x_0}$ and any $0<r\leq 1$,  it holds that $\mu^*_X(\mb B_{\m M_X}(x_0,r))\geq   g(r)\,r^{d_X}/L$ and $\mu^*_{Y|x_0}(\mb B_{\m M_{Y|x}}(y_0,r))\geq  g(r)\,r^{d_Y}/L$.
\end{enumerate}
}

\medskip
\noindent 
Compared to Regime 2, Regime 3b introduces an additional parameter $\beta_X$ that characterizes the smoothness of the manifold family $\big\{\mathcal{M}_{Y|x} : \, x \in \mathcal{M}_X\big\}$ with respect to the index variable $x$. Unlike Regime 3a, Regime 3b imposes weaker conditions on the covariate distribution, whose support is not necessarily a smooth submanifold and does not require a density that is bounded away from zero.
 On the other hand, since Regime 3b focuses on distribution estimation, it  requires a smoothness condition on the conditional density function $u^*$. Here, the conditional density function $u^*(y\,|\,x)$ operates on the joint space $\m M_{YX}=\{(y,x)\,:\, x\in 
 \m M_X,y\in \m M_{Y|x}\}$,  which cannot be decomposed into a product form $U_1\times U_2$ for the spaces of $y$ and $x$, due to the dependency of $\m M_{Y|x}$ on $x$. To quantify the smoothness of $u^*$ with respect to $y$ and $x$, we assume that $u^*$ can be expressed as the restriction of a function that is $\m H^{\alpha_Y,\alpha_X}$-smooth over the entire space $\mb R^{D_Y}\times \mb R^{D_X}$.   
 We are now ready to present our result on the minimax rate of convergence for distribution regression under Regime 3b.

\begin{theorem}\label{th2:lower}(Minimax rate for distribution regression under Regime 3b)  For each $\gamma>0$, if $\beta_Y\geq 2 \vee (\alpha_Y+1) \vee \beta_X$, $\beta_X\geq  \alpha_X+\frac{\alpha_X}{\alpha_Y}$ and $\alpha_Y\geq \alpha_X$, then there exits a constant $L_0$ so that when $L,\tau,\tau_1,n\geq L_0$, it holds that
    \begin{equation*}
    \begin{aligned}
       &C\,\bigg(n^{-\frac{\alpha_X}{2\alpha_X+d_X}}\ +\ n^{-\frac{\alpha_Y+\gamma}{2\alpha_Y+d_Y+\frac{\alpha_Y}{\alpha_X}d_X}}\ +\ n^{-\frac{\gamma}{\frac{d_Y}{\beta_Y}+\frac{d_X}{\beta_X}}} \ \bigg)\\
       &\qquad\leq \underset{\wh \mu_{Y|X}}{\inf}\,\underset{\mu^*=\mu^*_X\mu^*_{Y|X}\in \m P^*_3}{\sup}\mb{E}_{\mu^{*,\otimes n}}\Big[\mb{E}_{\mu^*_X}\big[d_{\gamma}\big(\mu^*_{Y|X},\,\wh \mu_{Y|X}\big)\big]\Big]\\
&\qquad\quad\leq C_1\bigg((\log n)^3\cdot n^{-\frac{\alpha_X}{2\alpha_X+d_X}}\ +\ (\log n)\cdot\Big(\frac{n}{\log n}\Big)^{-\frac{\alpha_Y+\gamma}{2\alpha_Y+d_Y+\frac{\alpha_Y}{\alpha_X}d_X}}\ +\ (\log n)\cdot n^{-\frac{\gamma}{\frac{d_Y}{\beta_Y}+\frac{d_X}{\beta_X}}}\ \bigg),        
    \end{aligned}
    \end{equation*}
    where $(C,C_1)$ are constants independent of $n$, and the infimum is taken over all conditional distribution estimators $\wh \mu_{Y|X}$ based on data $\big\{(X_i,Y_i)\big\}_{i=1}^n$ sampled from $\mu^{*,\otimes n}$. Here, the shorthand $\m P^*_3$ stands for $\m P^*_3(D_Y,D_X,d_Y,d_X,\beta_Y,\beta_X,\alpha_Y,\alpha_X,\tau,\tau_1,g,L)$
\end{theorem}

\noindent 
\begin{remark}\label{Resmoothness}
This theorem assumes that $\beta_Y \geq 2 \vee (\alpha_Y + 1) \vee \beta_X$, $\beta_X \geq \alpha_X + \frac{\alpha_X}{\alpha_Y}$, $\alpha_Y \geq \alpha_X$, and requires the stronger smoothness criteria $\m H^{\alpha_Y,\alpha_X}$ on the conditional density function. These conditions  enable a suitable decomposition of the distribution regression problem into two main tasks: manifold regression and density regression.    Specifically, for each fixed point $w_0=(x_0,y_0)$ in $\m M$ and for any $x$ near $x_0$,  we can perform localized analysis by restricting the measure $\mu^*_{Y|x}$ to $U_{w_0}\cap \m M_{Y|x}$, where $U_{w_0}$ is a defined neighborhood of $y_0$ on $\m M_Y$ (see Definition~\ref{def:manifold}). We then map the high-dimensional data points into a lower-dimensional latent space by projecting them onto a fixed tangent space $T_{y_0}\m M_{Y|x_0}$, that is, $\wt \pi_{w_0}(y)={\rm Proj}_{T_{y_0}\m M_{Y|x_0}}(y-y_0)$, and noting that each tangent vector can be uniquely represented  by a $d_y$-dimensional coordinate. The resulting push forward measure $[\wt \pi_{w_0}]_{\#}[\mu^*_{Y|x}|_{U_{w_0}}]$  admits a  density  function $v_{w_0}(z|x)$ with respect the volume measure on $T_{y_0}\m M_{Y|x_0}$,  given by
\begin{equation*}
    v_{w_0}(z|x)= u^*(\Phi_{w_0}(z,x)|x)\cdot \sqrt{|J_{\Phi_{\omega_0}(\cdot,x)}(z)^TJ_{ \Phi_{\omega_0}(\cdot,x)}(z)|_{+}}, \quad z\in \mb B_{T_{y_0}{\m M_{Y|x_0}}}(0,\tau_1).\footnote{Here we use $|\cdot|_{+}$ to denote the pseduo-determinant, which is the product of all non-zero eigenvalues of a square matrix. We use $\mathbf{J}_f(x)$ to denote the Jacobian matrix of $f$ evaluate at $x$, so that the $(i,j)$ element of  $\mathbf{J}_f(x)$  is $\frac{\partial f_i(x)}{\partial x_j}$.}
\end{equation*}
The $\m H^{\alpha_Y,\alpha_X}$-smoothness of $u^*(\cdot,\cdot)$ and the $\m H^{\beta_Y,\beta_X}$-smoothness of $\Phi_{w_0}(\cdot,\cdot)$, together with the conditions $\beta_Y\geq \alpha_Y+1$, $\beta_X\geq \alpha_X+\frac{\alpha_X}{\alpha_Y}$ and $\alpha_Y\geq \alpha_X$, then ensure that $ v_{w_0}(\cdot,\cdot)$ is $\ov {\m H}^{\alpha_Y,\alpha_X}$-smooth. Therefore,  if the tangent space at $w_0$ can be exactly recovered,  learning the  local conditional distribution near $w_0$ can be divided into: (1) learning  the  ${\m H}^{\beta_Y,\beta_X}$-smooth map $\Phi_{w_0}(\cdot,\cdot)$ (manifold regression); and (2) learning the $\ov{\m H}^{\alpha_Y,\alpha_X}$-smooth conditional density function $v_{w_0}$ (density regression). However, it  is generally impossible  to exactly recover the tangent space with only a finite number of samples around $w_0$. Nevertheless, it is possible to  approximate  a hyperplane $\wh T$ close to $T_{y_0}\m M_{Y|x_0}$. By adding the condition $\beta_Y\geq 2\vee \beta_X$,  it is ensured that for $x\approx x_0$, the function ${\rm Proj}_{\wh T}(y-y_0)$, which operates on $y\in \m M_{Y|x}\cap U_{w_0}$, is invertible. Moreover, the inverse function is $\m H^{\beta_Y,\beta_X}$-smooth when treating $x$ as an input (c.f. Lemma~\ref{le:projection} of Appendix~\ref{app: regularitymanifold}). Furthermore,  the push forward measure $[{\rm Proj}_{\wh T}(\cdot-y_0)]_{\#}[\mu^*_{Y|x}|_{U_{w_0}}]$  also admits an $\ov{\m H}^{\alpha_Y,\alpha_X}$-smooth conditional density  function (c.f.,  Lemma~\ref{le:defdistribution} of Appendix~\ref{app: regularitymanifold}). This allows for a similar decomposition of the problem, even if the tangent space cannot be precisely recovered.
    \end{remark}


When comparing the minimax rate in Theorem~\ref{th2:lower} with that from Theorem~\ref{th2:lower1}, the key difference lies in the last term related to supporting manifold estimation. Specifically, by setting $\gamma=1$, the term $n^{-{\gamma}/{(\frac{d_Y}{\beta_Y}+\frac{d_X}{\beta_X})}}$ in Theorem~\ref{th2:lower}, up to logarithmic terms, matches the minimax optimal rate for manifold regression on the $(\beta_Y,\beta_X)$-smooth manifold family, as obtained in Theorem~\ref{th:mainifold}. To simplify notation, we use $\upsilon_1 = {d_Y}/{\beta_Y}$ and $\upsilon_2 = {d_X}/{\beta_X}$ to denote the complexity indices characterizing the supporting manifolds associated with the response variable $Y$ and the covariate variable $X$, respectively—defined as the intrinsic dimensions scaled by manifold smoothness; and $\upsilon_3 = {d_Y}/{\alpha_Y}$, $\upsilon_4 = {d_X}/{\alpha_X}$ to denote the complexity indices characterizing the conditional distribution class with respect to inputs $y$ and $x$, defined as the input intrinsic dimensions scaled by the corresponding density marginal smoothness. The minimax rates in Theorem~\ref{th2:lower} can then be expressed as $\wt{\m O}\big(n^{-\frac{1}{2+\upsilon_4}}+n^{-\frac{1+\gamma/\alpha_Y}{2+\upsilon_3+\upsilon_4}}+n^{-\frac{\gamma}{\upsilon_1+\upsilon_2}}\big)$, which depends on the magnitude of the intrinsic dimensions relative to the smoothness levels, the value of $\gamma$, and its proportion relative to $\alpha_Y$. Similar to Regime 1 and 2, the dominant term in the overall rate varies with different values of $\gamma$. When $\gamma$ is sufficiently small, specifically, $\gamma\leq \frac{\upsilon_1+\upsilon_2}{2+\upsilon_3+\upsilon_4-\frac{1}{\alpha_Y}\cdot ((\upsilon_1+\upsilon_2)\wedge (\upsilon_3\alpha_Y))}$, the manifold regression hardness becomes the bottleneck, and the dominant term in the minimax rate is $n^{-\frac{\gamma}{\upsilon_1+\upsilon_2}}$. When $\frac{\upsilon_1+\upsilon_2}{2+\upsilon_3+\upsilon_4-\frac{1}{\alpha_Y}\cdot ((\upsilon_1+\upsilon_2)\wedge (\upsilon_3\alpha_Y))}\leq \gamma\leq \frac{(\upsilon_1+\upsilon_2)\vee (\upsilon_3\alpha_Y)}{2+\upsilon_4}$, the term $n^{-\frac{1+\gamma/\alpha_Y}{2+\upsilon_3+\upsilon_4}}$ related to nonparametric conditional density estimation becomes dominant. If $\gamma$ increases beyond $\frac{(\upsilon_1+\upsilon_2)\vee (\upsilon_3\alpha_Y)}{2+\upsilon_4}$,  the dominant term becomes the nonparametric mean regression risk $n^{-\frac{1}{2+\upsilon_4}}$, reflecting the overall dependence trend of $Y$ on $X$ (see the discussion after Theorem~\ref{th1:lower}).   

\begin{figure}
    \centering
    \includegraphics[width=0.5\linewidth]{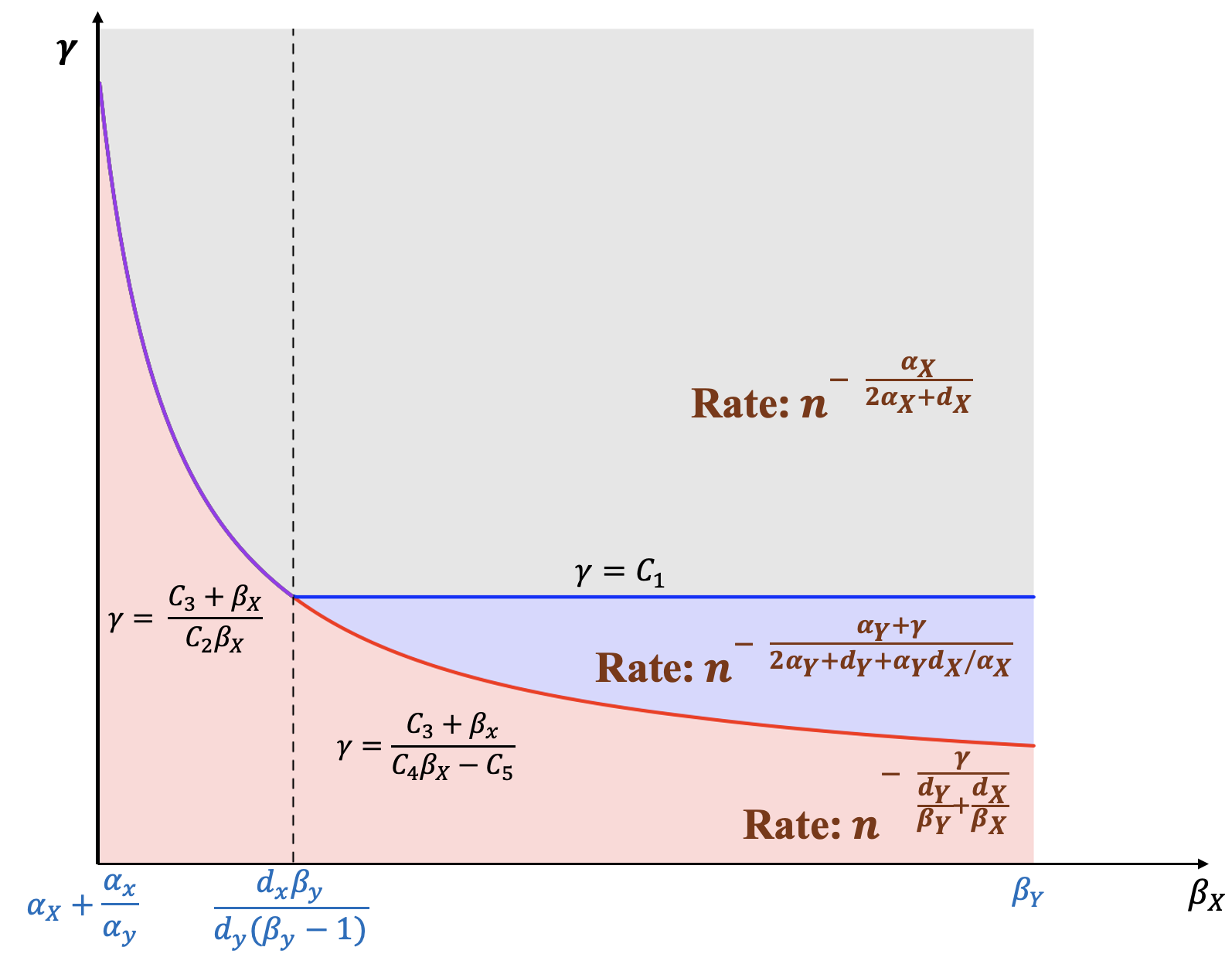}
    \caption{Diagram for the minimax rate under Regime 3 for varying $\gamma$ and $\beta_X$, where $C_1=\frac{d_Y\alpha_X}{2\alpha_X+d_X}$, $C_2=(2+\frac{d_X}{\alpha_X})\frac{\beta_Y}{d_Y}$, $C_3=\frac{d_X\beta_Y}{d_Y}$, $C_4=(2+\frac{d_X}{\alpha_X}+\frac{d_Y}{\alpha_Y})\frac{\beta_Y}{d_Y}-\frac{1}{\alpha_Y}$ and $C_5=\frac{d_X\beta_Y}{d_Y\alpha_Y}$.}
    \label{fig:2}
\end{figure}

Figure~\ref{fig:2} illustrates these three regimes  with varying $\beta_X$ and $\gamma$.  When $\beta_X$  falls within the interval $ (\alpha_X+\frac{\alpha_X}{\alpha_Y}, \frac{d_X\beta_Y}{d_Y(\beta_Y-1)})$---assuming this interval is non-empty---there are only two regimes, where the rate for nonparametric conditional density estimation is either dominated by that for manifold regression or by that for  nonparametric mean regression. The transition point in terms of $\gamma$ decreases with increasing $\beta_X$, as a larger $\beta_X$ reduces the challenges for manifold regression,  allowing for less smooth test functions to be effective in averaging out minor irregularities in the support, thereby focusing more on the global dependence of $Y$ on $X$. When $\beta_X\notin (\alpha_X+\alpha_X/\alpha_Y, \frac{d_X\beta_Y}{d_Y(\beta_Y-1)})$, all three regimes become possible. For the first transitions  (from manifold regression to  nonparametric conditional density regression), the transition point in terms of $\gamma$ decreases with increasing $\beta_X$, as large $\beta_X$ ease the manifold regression, prompting an earlier shift in challenges of nonparametric conditional density estimation. While the second transitions point (from nonparametric conditional density regression to mean regression) remains constant relative to $\beta_X$, as the rates for these tasks are independent of $\beta_X$.

A natural extension beyond our current setting is the noisy case, corresponding to a singular measure deconvolution problem in which the observed data are contaminated by additive noise. Specifically, we observe $n$ i.i.d.~samples $\big\{(X_i, Z_i)\big\}_{i=1}^n$ generated according to the model $X_i \sim \mu^*_X$, $Y_i \sim \mu^*_{Y|X_i}$, and $Z_i = Y_i + \varepsilon_i$, where $\{\varepsilon_i\}_{i=1}^n$ are i.i.d.~zero-mean errors independent of $\big\{(X_i, Y_i)\big\}_{i=1}^n$. The goal is to recover $\mu^*_{Y|x}$ and the underlying manifolds $\m M_{Y|x}$ for each ``noiseless" covariate value $x$, based on noisy measurements. However, even the support recovery problem of estimating a single manifold from noisy observations is intrinsically difficult: for instance, \cite{10.1214/12-AOS994} show that when the noise is Gaussian, the minimax rate of manifold estimation under the Hausdorff distance is lower bounded by $C(\log n)^{-1}$. One way to mitigate this slow convergence is to assume that the noise variance $\sigma^2$ decreases with the sample size. For clarity and simplicity, this paper focuses on the noiseless scenario and defers the detailed exploration of the deconvolution problem, including the analysis of the minimax rate in terms of both $n$ and $\sigma^2$, to future work.

\section{Minimax Optimal Estimators for Distribution Regression}
In this section, we introduce our conditional distribution  estimators designed to achieve minimax upper bounds across different regimes.  We will start with the simpler estimator for Regime 1,  where the response space is Euclidean. Following this, we will proceed to describe the more complex estimators for Regimes 2 and 3,  where the response variable lies in a low-dimensional manifold.

A key component of our approach is the use of multi-scale function decomposition via wavelet, which provides a robust framework for analyzing functions by separating them into components at different levels of detail. This methodology is particularly effective for characterizing H\"{o}lder regularity, as it captures both local and global smoothness properties through a hierarchical analysis of the structure of the function. The core of this decomposition is the concept of a wavelet, defined as a rapidly decaying and localized oscillating function. Commonly used constructions include the Haar basis~\citep{triebel2010bases}, Meyer basis~\citep{Triebel2006,meyer1992wavelets}, and Daubechies basis~\citep{daubechies1988orthonormal}, among others. A fundamental aspect of wavelet analysis is the concept of scaling, which involves stretching or shrinking the wavelet to adapt to different features of the target function. By stretching the analyzing function, one can capture slowly varying and global trends, whereas shrinking it allows the detection of abrupt changes and fine details. 

To put things more formally, consider the space $\m L^2(\mb R^d)$ of square-integrable functions on $\mb R^d$. Within this space, one can construct a complete orthonormal basis $\bigcup_{j\geq 0}\ov\Psi_j^{d}$ formed by localized oscillatory functions. The level-zero basis $\ov\Psi_0^d$ is generated by shifting a compactly supported scaling function, while the higher-level bases $\ov\Psi_j^d$ are formed by shifting and scaling a compactly supported oscillatory function by a factor of $2^{-(j-1)}$. As a result, any function $p\in \m L^2(\mb R^d)$ admits a unique expansion of the form 
\begin{equation*} p(x)=\sum_{j=0}^{\infty} \sum_{\psi\in \ov\Psi_j^d} p_{\psi} \,\psi(x)\quad\mbox{with}\quad p_{\psi}=\int_{\mb R^d} \psi(x)p(x)\,\dd x. \end{equation*} 
The coefficients $p_{\psi}$ reveal different aspects of the function $p(x)$: at lower levels (i.e., $\psi\in \Psi_j^d$ with small $j$), they capture broad, slowly varying trends, whereas at higher levels, they are sensitive to fine details and abrupt variations. Consequently, for a smooth function $p(\cdot)$ without significant local oscillations, the wavelet coefficients $p_\psi$ tend to be small in absolute value for higher levels. In particular, if $p(\cdot)$ belongs to the H\"{o}lder space $\m H^{\alpha}$ with bounded norm, then for any $j\in \mb N$ and $\psi\in \Psi_j^d$, the coefficients satisfy the bound $|p_{\psi}|\leq C\, 2^{-\frac{dj}{2}-j\alpha}$ for some constant $C$ independent of $j$. Further details on wavelet theory are provided in Appendix~\ref{app:asswav}.

Throughout the following, for any dimension $d$, we use $\bigcup_{j\geq 0}\ov\Psi_j^{d}$ to denote an orthonormal wavelet basis satisfying appropriate smoothness conditions, as specified in Lemma~\ref{le:wavelet} of Appendix~\ref{app:asswav} (for example, the Daubechies basis~\citep{daubechies1988orthonormal}).  The precise smoothness requirements for different regimes are detailed in Appendix~\ref{app:estimatordetail}.

\subsection{Minimax optimal estimator for Euclidean response space}\label{sec:estregime1}
In this subsection, we focus on Regime 1, where the conditional distribution $\mu^*_{Y|X}$ is characterized by a $\ov{\m H}^{\alpha_Y,\alpha_X}$-smooth conditional density function $u^*(\,\cdot\,|\,\cdot\,)$. Our goal is to construct an estimator for this conditional density. For any $x\in \m M_X$, given that $u^*(\,\cdot\,|\,x)$ is assumed to be compactly supported  within $\mb B_{\mb R^{D_Y}}(\mathbf{0},L)$,  we define for any $j\in \mb N$, 
\begin{equation}\label{defPsi_j1}
      \Psi_j^{D_Y}=\big\{\psi\in \ov \Psi_j^{D_Y}:\, \text{supp}(\psi)\cap \mb B_{\mb R^{D_Y}}(\mathbf{0},L)\neq \emptyset\big\}.
\end{equation}
Then the function $u^*(\,\cdot\,|\,x)$ has a wavelet expansion as 
\begin{equation*}
    u^*(y\,|\,x)=\sum_{j=0}^{\infty} \sum_{\psi\in \Psi_j^{D_Y}} u^*_{\psi}(x)  \,\psi(y)\quad\mbox{with}\quad u^*_{\psi}(x)=\mb{E}_{\mu^*_{Y|x}}[\psi(y)]=\int_{\mb R^{D_Y}} \psi(y)\,u^*(y\,|\,x)\,\dd y.
\end{equation*}
Since $u^*(\,\cdot\,|\,x)$ belongs to $\m H^{\alpha_Y}_{L}(\mb R^{D_Y})$, we truncate its wavelet expansion at a finite level $J$ to eliminate high-frequency fluctuations. The value of $J$ will be carefully chosen later to balance the bias-variance trade-off. Consequently, the problem of jointly estimating $u^*(y\,|\,x)$ over $x\in \m M_X$ reduces to the joint estimation of the wavelet coefficients $u^*_{\psi}(x)$ for $j\in \mb N$, $\psi\in \Psi_j^{D_Y}$, and $x\in \m M_X$. Observing that each coefficient $u^*_{\psi}(x)$ can be expressed as the conditional mean $\mb{E}{\mu^*_{Y|x}}[\psi(y)]$, the estimation of $u^*_{\psi}(x)$ for different $\psi$ can be formulated as a collection of regression problems, where the response variables are $\big\{\psi(Y_i)\big\}_{i=1}^n$ with covariates $\{X_i\}_{i=1}^n$.

For each level $j\in {0}\cup [J]$, we consider an approximation family $\ms S_j$ consisting of functions mapping $\mb R^{D_X}$ to $\mb R$. For each $\psi\in \Psi_j^{D_Y}$, we minimize the mean squared error to obtain
\begin{equation*} 
\begin{aligned}
\wh u_{\psi}(\cdot)=\underset{ u\in  \ms S_j}{\arg\min} \ \frac{1}{n}\sum_{i=1}^n  \big(\psi(Y_i)-u(X_i)\big)^2.
\end{aligned}
   \end{equation*}
Note that this estimation procedure uses the same approximation family $\ms S_j$ for coefficients of $\psi$  at each specific level  $j$, while $\ms S_j$ varies across different levels $j$.

To construct the approximation family $\ms S_j$, we leverage the fact that, for each $\psi\in \Psi_j^{D_Y}$, the conditional mean $\mb{E}_{\mu^*_{Y|x}}[\psi(Y)] = \int_{\mb R^{D_Y}} \psi(y) \,u^*(y\,|\,x)\,\dd y$ is a $\m H^{\alpha_X}$-smooth function of $x$, with its H\"{o}lder norm bounded by $\m O\big(2^{-\frac{D_Yj}{2}}\big)$. This property motivate us to define the following approximating family by utilizing local polynomial approximations for H\"{o}lder-smooth functions,
 \begin{equation}\label{defmsS_j}
\begin{aligned}
     \ms S_j&=\Bigg\{u(x)= \frac{\sum_{i=1}^{W_j} \sum_{k\in \mb N_0^{D_X}, |k|<\alpha_X}a_{ik} (x-b_i)^k\, \rho\Big(\frac{\|x-b_i\|}{\varepsilon^x_j}\Big)}{\sum_{i=1}^{W_j} \rho\Big(\frac{\|x-b_i\|}{\varepsilon^x_j}\Big)+\frac{1}{n}}:\, b_i\in \mb B_{\mb R^{D_X}}(\mathbf{0},L), \\
    & \qquad\qquad\qquad   a_{ik}\in \Big[-\frac{C}{2^{D_Yj/2}},\frac{C}{2^{D_Yj/2}}\Big], \ \  \text{ for any } i\in[W_j]\  \mbox{and multi-index}  \ k\ \Bigg\},
\end{aligned}
\end{equation}
where $\varepsilon^x_j=2^{jD_Y/({2\alpha_X+d_X)}}\big(\frac{n}{\log n}\big)^{-{1}/{(2\alpha_X+d_X)}}$, $W_j=C_1\,(\varepsilon^x_j)^{-d_X}$, and $(C,C_1)$ are some sufficiently large constants. Here, $\rho$ is a smooth transition function satisfying $\rho(t)=1$ for $t\in [0,1]$ and $\rho(t)=0$ for $t\geq 2$.
The function $\rho\big(\|x-b_i\|/\varepsilon^x_j\big)$ effectively partitions the covariate space $\m M_X$ into local neighborhoods, where the radius (bandwidth) and number of neighborhoods scale with the effective dimension $d_X$ of $\m M_X$. Within each neighborhood, the conditional mean $\mb{E}_{\mu^*_{Y|x}}[\psi(Y)]$ can be well approximated by a low-degree polynomial in $x$.

By substituting the estimator $\wh u_{\psi}(x)$ into the truncated wavelet expansion of $u^*(y\,|\,x)$, we can derive a conditional density estimator as
\begin{equation}\label{est:u_c}
   \wh u(y\,|\,x)=\sum_{j=0}^J \sum_{\psi\in \Psi_j^{D_Y}}  \wh u_{\psi}(x)\,\psi(y),\quad x\in \m M_X.
\end{equation}
The following theorem shows that the conditional distribution estimator $\wh \mu_{Y|X}$, whose density function is $\wh u(y\,|\,x)$, can achieve the minimax upper bound stated in Theorem~\ref{th1:lower} \emph{simultaneously} for all $\gamma\geq 0$.

\begin{theorem}[Convergence rate for density regression estimator in Regime 1]\label{th1}
Let $\m P_1^*$ be the target distribution class defined in Theorem~\ref{th1:lower}. Suppose $\big\{(X_i,Y_i)\big\}_{i=1}^n$ are $n$ i.i.d.~samples from $\mu^*$, and set $J=\lceil {\frac{1}{2\alpha_Y+D_Y+d_X\frac{\alpha_Y}{\alpha_X}}}\cdot \log_2 (\frac{n}{\log n})\rceil$.  For any $\mu^*=\mu^*_X\mu^*_{Y|X}\in \m P^*_1$, the following holds with probability at least $1-n^{-1}$: for any $\gamma\geq 0$, the conditional density estimator $\wh u$ defined in~\eqref{est:u_c} satisfies
    \begin{equation*}
    \begin{aligned}
\mb{E}_{\mu^*_X}\big[d_{\gamma}\big(\mu^*_{Y|X},\,\wh \mu_{Y|X}\big)\big]\leq C_{\gamma}\bigg(\sqrt{\log n}\cdot \Big(\frac{n}{\log n}\Big)^{-\frac{\alpha_X}{2\alpha_X+d_X}}\ + \ \Big(\frac{n}{\log n}\Big)^{-\frac{\alpha_Y+\gamma}{2\alpha_Y+D_Y+\frac{\alpha_Y}{\alpha_X}d_X}}\, \bigg),
    \end{aligned}
    \end{equation*}
    where $C_{\gamma}$ is a constant independent $n$.
\end{theorem}
%
\noindent A complete proof of Theorem~\ref{th1} is provided in Appendix~\ref{proof:th1}, and further details on the estimator construction are summarized in Appendix~\ref{R1gamma>0}. A key observation is that the bandwidth $\varepsilon_j^x$ increases with the level $j$, in contrast to the bandwidth $2^{-j}$ in $Y$, which decreases as $j$ increases. This asymmetric scaling is crucial for balancing the overall complexity of jointly estimating $\mb{E}_{\mu^*_{Y|x}}[\psi(Y)]$ across different levels $j$. Intuitively, as $j$ increases, the resolution in $Y$ becomes finer because the bandwidth decreases, allowing the model to capture more detailed variations in $Y$. At the same time, the resolution in $X$ becomes coarser because the bandwidth increases, meaning that the model mainly captures broad, global patterns in $X$ while finer structures in $Y$ are being learned.
     
The multiresolution analysis underlying wavelet decompositions shares a close connection with score-based forward backward diffusion models~\citep{song2020score} for implicit distribution estimation. For example, in the backward diffusion model, the data generation process gradually builds structure by transforming white noise into realistic data, following a progression from coarse to fine details. This process parallels how multiresolution analysis decomposes a function, first capturing global trends and then progressively refining finer structures. In particular, when comparing the conditional diffusion model with the wavelet-based conditional distribution estimator, both approaches can be viewed as solving multiple mean regression problems across different resolution levels. More specifically, the index $j\in \mb N$ in the preceding wavelet estimator and the time variable $t\in \mb R^{+}$ in the backward diffusion model and both represent levels of resolution, controlling the scale of analysis from coarse to fine details.

 \subsection{Minimax optimal estimator for manifold response space}
In this subsection, we focus on Regimes 2 and 3, where, given $X=x$ for $x\in \m M_X$, the conditional distribution $\mu^*_{Y|x}$ is supported on a $d_Y$-dimensional submanifold $\m M_{Y|x}$. 
Since the conditional density with respect to the Lebesgue measure does not exist in these regimes, we reformulate the conditional distribution estimation problem as one that involves simultaneously estimating the conditional expectations $\m J^*(f,x) := \mb{E}_{\mu^*_{Y|x}}[f(Y)]$ for a class of test functions $f\in \m H_1^{\gamma}(\mb R^{D_Y})$ and covariate values $x\in \m M_X$, where $\gamma \geq 0$ corresponds to the same smoothness index used in defining the H\"{o}lder IPM $d_\gamma$. In other words, we construct an explicit estimator for the conditional expectation functional, denoted by $\m{\wh J}:\, \m H_1^{\gamma}(\mb R^{D_Y}) \times\mb R^{D_X} \to \mb R$, and evaluate its performance using the \emph{simultaneous estimation risk}
\begin{align}\label{eq:ser}
    \mb{E}_{\mu^*_X}\bigg[\sup_{f\in \m H_1^{\gamma}(\mb R^{D_Y})} \Big|\,\m{\wh J}(f,x)-\m J^*(f,x)\,\Big|\bigg].
\end{align}
There exists a one-to-one correspondence between the conditional distribution $\mu^*_{Y|x}$ and the conditional expectation functional $\m J^*(\,\cdot\,, x)$ evaluated over any rich enough class of test functions that is dense in $\m L_2(\mb R^{D_Y})$, such as $\m H^{\gamma}(\mb R^{D_Y})$. As a result, estimating the conditional distribution is equivalent to estimating its associated conditional expectation functional.  Specifically,  as discussed in~\citep{tang2023minimax}, for any fixed $\gamma>0$ and $x$, one can employ adversarial training with $\gamma$-smooth test functions to obtain a conditional distribution estimator
\begin{equation*}
 \wh \mu^{\gamma}_{Y|x}= \underset{\mu \in \m P^*_{Y}}{\arg\min}  \underset{f \in \m H^{\gamma}_1(\mb R^{D_Y})}{\sup} |\mb{E}_{\mu}[f(y)] - \wh{\m J}(f, x)|
\end{equation*}
 where for a suitable $\m P_Y^*$, the estimation error of $\wh \mu^{\gamma}_{Y|x}$ under the $d_{\gamma}$ metric can be bounded from above by twice the maximal deviation between  $\wh {\m J}(f,x)$ and $\m J^*(f,x)$ over $f\in\m H_1^{\gamma}(\mb R^D)$.  Furthermore, given a suitable set $\Gamma$ of $\gamma$ values, consider the estimator:
\begin{equation*}
    \wh{\mu}_{Y|x} = \underset{\mu \in \m P^*_{Y}}{\arg\min} \sum_{\gamma \in \Gamma} \frac{1}{\delta_{n,\gamma}} \cdot \underset{f \in \m H^{\gamma}_1(\mb R^{D_Y})}{\sup}  \left[\mb{E}_{\mu}[f(y)] - \wh{\m J}(f, x)\right],
\end{equation*}
with appropriate choices for $\m P^*_Y$ and $\delta_{n,\gamma}$. This estimator is simultaneously minimax optimal up to logarithmic factors for all $\gamma > 0$ (cf. Corollary~\ref{co1} and Corollary~\ref{co2} in Appendix~\ref{app:estimatordetail}). This optimality is attained by incorporating a jointly optimal $\wh{\m J}(\,\cdot\,,\, \cdot\,)$, the construction of which will be detailed below.


To construct an estimator $\wh{\m J}(f,x)$ for $\m J^*(f,x)$, we first observe that, since $\mu^*_{Y|x}$ is compactly supported for any $x\in \m M_X$, it suffices to restrict our analysis to test functions $f\in \m L^2(\mb R^{D_Y}) \cap \m H_1^{\gamma}(\mb R^{D_Y})$. Each such function admits a wavelet expansion
\begin{equation*}
    f(y)=\sum_{j=0}^{\infty}\sum_{\psi\in \ov \Psi_j^{D_Y}} f_{\psi} \psi(y) \ \ \ \text{ with } \ f_{\psi}=\int_{\mb R^{D_Y}} f(y)\,\psi(y)\,\dd y.
\end{equation*}
We fix a finite truncation level $J$ (to be specified later) and consider the wavelet thresholding approximation $f_J$ of $f$:
 \begin{equation*}
     f(y) \ \ = \ \ \underbrace{\sum_{j=0}^{J}\sum_{\psi\in \ov \Psi_j^{D_Y}} f_{\psi} \psi(y)}_{f_J(y)} \ \ + \ \ \underbrace{\sum_{j=J+1}^{\infty}\sum_{\psi\in \ov \Psi_j^{D_Y}} f_{\psi} \psi(y)}_{f_J^{\perp}(y)}\,,
 \end{equation*}
with $f_J^{\perp}$ denoting the corresponding remainder term.
The thresholding approximation $f_J(\cdot)$ primarily captures the slowly varying and global structure of the function $f$, while the remainder term $f_J^{\perp}(\cdot)$ accounts for the more abrupt, localized variations and oscillations. By decomposing the conditional expectation as $\m J^*(f,x)=\mb{E}_{\mu^*_{Y|x}}[f_J(Y)]+\mb{E}_{\mu^*_{Y|x}}[f_J^\perp(Y)]$, we then estimate the two components $\mb{E}_{\mu^*_{Y|x}}[f_J(Y)]$ and $\mb{E}_{\mu^*_{Y|x}}[f_J^\perp(Y)]$ using different strategies.

\subsubsection{Estimator for coarse-scale component $\mb{E}_{\mu^*_{Y|x}}[f_J(Y)]$}   Given the inherent smoothing effect of the truncation operation in $f_J$, minor irregularities in the conditional distribution $\mu^*_{Y|x}$ have limited impact and are effectively averaged out. Based on this observation, we construct a estimator for $\mb{E}_{\mu^*_{Y|x}}[f_J(Y)]$ by treating $\mu^*_{Y|x}$ as if it admits a density with respect to the Lebesgue measure on $\mb R^{D_Y}$. Specifically, we estimate the coarse-scale component $\mb{E}_{\mu^*_{Y|x}}[f_J(y)]$ in $\m J^*(f,x)$ by $\int_{\mb R^{D_Y}} f_J(y)\,\wh u(y\,|\,x)\,\dd y$, where the ``conditional density" estimator $\wh u(y\,|\,x)$ is constructed solely to define this integral. The construction follows a strategy similar to that introduced in Section~\ref{sec:estregime1} for Regime 1, as detailed below.

To construct the conditional density estimator $\wh u$, we begin by simultaneously estimating the conditional means $\mb{E}_{\mu^*_{Y|x}}\big[2^{j(d_Y - D_Y)/2}\psi(y)\big]$ for all $j\in \{0\}\cup [J]$, $\psi\in \Psi_j^{D_Y}$ and $x\in \m M_X$. The scaling factor $2^{j(d_Y - D_Y)/2}$ is introduced to account for the intrinsic dimension $d_Y$ of the support of $\mu^*_{Y|x}$, ensuring that the second moment $\mb{E}_{\mu^*_{Y|x}}\big[(2^{j(d_Y - D_Y)/2}\psi(y))^2\big]$ remains bounded, i.e., of order $\m O(1)$. 
In contrast to the method used in Section~\ref{sec:estregime1}, where each conditional expectation was estimated independently through separate mean regression problems, we adopt a joint estimation strategy that better reflects the low-dimensional structure of the support $\m M_{Y|x}$. Estimating each $\mb{E}_{\mu^*_{Y|x}}[\psi(y)]$ separately may ignore geometric dependencies and lead to inefficient use of data.
Instead, we treat the wavelet function $\psi$ as an additional input, alongside $x$, and formulate a \emph{joint mean regression} problem over the product space $\Psi_j^{D_Y} \times \mb R^{D_X}$. This leads us to define an estimator $\wh S_j^{\dagger}$ satisfying $\wh S_j^{\dagger}(\psi,x)\approx \mb{E}_{\mu^*_{Y|x}}\big[2^{j(d_Y - D_Y)/2}\psi(y)\big]$. To this end, for each $j \in {0} \cup [J]$, we introduce a function class $\m S_j^{\dagger}$ consisting of mappings $S:\, \Psi_j^{D_Y} \times \mb R^{D_X} \to \mb R$, and formulate the following joint mean regression problem by minimizing the aggregated squared loss over all $\psi \in \Psi_j^{D_Y}$:
\begin{equation}\label{eq:cde}
    \wh S_j^{\dagger}={\arg\min}_{S\in \m S_j^{\dagger}}\ \ \frac{1}{n}\sum_{i=1}^n \sum_{\psi\in \Psi_j^{D_Y}}\Big(2^{\frac{j(d_Y-D_Y)}{2}}\psi(y)-S(\psi,X_i)\Big)^2.
\end{equation}
Note that the separate mean regressions described in Section~\ref{sec:estregime1} can be viewed as a special case of the joint mean regression framework introduced above. In that setting, the approximation family for $S$ is separable in $\psi$ and takes the form $\m S_j=\big\{S(\psi,x)=\sum_{\psi'\in \Psi_j^{D_Y}}s_{\psi'}(x)\cdot\mathbf{1}(\psi'=\psi),\ \text{such that } s_{\psi'}\in\ms S_j \text{ for each }\psi'\in \Psi_j^{D_Y}\big\}$.   However, this separable approximation family does not allow the sharing of information across different $\psi$. Specifically, due to the manifold structure of the response space, only a subset of the functions $\psi(\cdot)$ have non-zero conditional means. This inherent sparsity is not fully utilized in separate mean regression. In contrast, by choosing $\m S_j$ in a non-separable form, one can more effectively exploit this structure. Further details on these constructions are provided in Appendix~\ref{R2gamma} (for Regime 2) and Appendix~\ref{R3gamma>0} (for Regime 3b). The conditional density estimator $\wh u(y\,|\,x)$ is then defined as 
\begin{equation*}
    \wh u(y\,|\,x)= \sum_{j=0}^{J}\sum_{\psi\in \Psi_j^{D_Y}}  2^{\frac{j(D_Y-d_Y)}{2}} \,\wh S_j^{\dagger}(\psi,x) \,\psi(y), 
\end{equation*}
and the associated estimator for $\mb{E}_{\mu^*_{Y|x}}[f_J(Y)] $ is given by 
\begin{equation*}
    \int_{\mb R^{D_Y}} f_J(y)\, \wh u(y\,|\,x)\,\dd y=  \sum_{j=0}^{J}\sum_{\psi\in \Psi_j^{D_Y}}  2^{\frac{j(D_Y-d_Y)}{2}}   \, f_{\psi}\,\wh S_j^{\dagger}(\psi,x). 
\end{equation*}

\subsubsection{Estimator for fine-scale component $\mb{E}_{\mu^*_{Y|x}}[f_J^\perp(Y)]$}
This term is more sensitive to fine-scale structure and to potential misalignment in the support of the distributions resulting from manifold estimation. To address this, we incorporate an explicit manifold estimation step by learning $x$-dependent local charts of the submanifold $\m M_{Y|x}$. Specifically, for each local patch of the joint space $\m M$, we learn an \emph{encoder} $Q:\, \mb R^{D_Y} \to \mb R^{d_Y}$ and a \emph{conditional decoder} $G: \,\mb R^{d_Y} \times \mb R^{D_X} \to \mb R^{D_Y}$ such that the reconstruction relation $y \approx G(Q(y), x)$ holds for each $(x, y)$ in the patch. 

These estimated charts allow us to map the data into a low-dimensional latent space $\mb R^{d_Y}$, where subsequent analysis becomes more tractable.
In the second step, we will perform density regression in the latent space using the transformed samples $\big\{\big(X_i, Q(Y_i)\big)\big\}_{i=1}^n$ to estimate conditional density functions associated with the latent distributions. This encoder–decoder framework, which shifts the analysis from the ambient to a lower-dimensional latent space, is widely used in practice, including in methods such as latent diffusion models~\cite{Rombach_2022_CVPR}, variational autoencoders~\cite{kingma2013auto}, and Wasserstein autoencoders~\cite{tolstikhin2017wasserstein}, among others.

The final estimator is formulated as a mixture of conditional generative models, given by 
\begin{align}\label{eq:mgm}
    \sum_{k\in \wh {\m K}}\big[\wh G_\sk (\,\cdot\,,x)\big]_{\#}\wh\nu_\sk(\,\cdot\,|\,x),
\end{align}
where $\wh{\m K}$ is a data-dependent index set, $\wh G_\sk: \mb R^{d_Y} \times \mb R^{D_X} \to \mb R^{D_Y}$ is a learned decoding map from the latent space $\mb R^{d_Y}$ to the data ambient space $\mb R^{D_Y}$, and $\wh \nu_\sk(\,\cdot\,|\,x)$ is an estimated conditional distribution of the latent variable on $\mb R^{d_Y}$. This pushforward measure serves as a surrogate for $\mu^*_{Y|x}$ in the estimation of the fine-scale component $\mb{E}_{\mu^*_{Y|x}}[f_J^\perp(Y)]$.

For technical simplicity, we adopt a data-splitting strategy to divide the dataset into two disjoint subsets: $I_1 = [\lfloor n/2 \rfloor]$ and $I_2 = [n] \setminus I_1$. The two-step estimation procedure described above can be summarized in the following concrete algorithm.


\medskip
\noindent{\bf \underline{Manifold estimation}:} Let $\{\omega_k = (x_k, y_k)\}_{k=1}^{K}$ be a $\tau_2$-covering set of $\mb B_{\mb R^{D_X}}(\mathbf{0}, L) \times \mb B_{\mb R^{D_Y}}(\mathbf{0}, L)$, where $\tau_2$ is a sufficiently small absolute constant. Define
 \begin{equation*}
     \wh{\m K}=\big\{k\in [K]:\, \exists i\in I_1, \|(X_i,Y_i)-\omega_k\|\leq \sqrt{2}\,\tau_2\big\}.
 \end{equation*}
Let $\m G$ be a family of functions $G: \mb R^{d_Y}\times \mb R^{D_X}\to \mb R^{D_Y}$. For each $k\in \wh{\m K}$, we define the estimator
\begin{equation*}
\begin{aligned}
    & (\wh G_{[k]}, \wh V_{[k]})\\
    &=\underset{G\in  \m G\atop V\in \mb O(D_Y,d_Y)}{\arg\min} \frac{1}{|I_1|}\sum_{i\in I_1} \big\|Y_i-G(V^T(Y_i-y_k),X_i)\big\|^2\cdot\mathbf{1}\big(Y_i\in \mb B_{\mb R^{D_Y}}(y_k,2\tau_2)\big)\cdot \mathbf{1}\big(X_i\in \mb B_{\mb R^{D_X}}(x_k,2\tau_2)\big),
\end{aligned}
\end{equation*}
where $\mb O(D_Y,d_Y)=\{A\in \mb R^{D_Y\times d_Y}:\, A^TA=I_{d_Y}\}$.

\medskip
\noindent{\bf \underline{Density regression on the latent space}:} Denote $\wh Q_{[k]}(y)=\wh V_{[k]}^T(y-y_k)$. For any $j\in \mb N$,  we define $ \Psi_j^{d_Y}=\big\{\psi\in \ov \Psi_j^{d_Y}:\, \text{supp}(\psi)\cap \mb B_{\mb R^{d_Y}}(\mathbf{0},2\tau_2)\neq \emptyset\big\}$. Let $\ms S_j$ denote a class of functions $v:\mb R^{D_X}\to \mb R$.  For each $k\in \wh {\m K}$, $j\in \{0\}\cup [J]$ and $\psi\in \Psi_j^{d_Y}$,  we define the estimator
\begin{equation}\label{eqn:whv}
\begin{aligned}
    & \wh v_{k\psi}(\cdot)={\arg\min}_{v\in  \ms S_j}\ \frac{1}{|I_2|}\sum_{i\in I_2} \Big[\,\psi\big(\wh Q_{[k]}(Y_i)\big)\cdot \rho_\sk(X_i,Y_i)\,-\,v(X_i)\,\Big]^2,
\end{aligned}
   \end{equation}
where $\rho_\sk(x,y)=\frac{\rho( {\|(x,y)-(x_k,y_k)\|^2}/{\tau_2^2})}{\sum_{k=1}^K \rho( {\|(x,y)-(x_k,y_k)\|^2}/{\tau_2^2})}$ with $\rho$ being a smooth transition function taking value $1$ on $[0,1]$ and zero on $[2,\infty)$. Here, the functions $\{\rho_\sk\}_{k=1}^K$ serve as a partition of unity, allowing the local constructions around each $\omega_k$ to be smoothly combined into a global estimator.

\medskip
\noindent{\bf \underline{Final estimator for $\mb{E}_{\mu^*_{Y|x}}[f_J^\perp(Y)]$}:} Denote $\wh \nu_\sk(\cdot|x)$ as the measure that has a density function $$\sum_{j=0}^J\sum_{\psi\in \Psi_j^{d_Y}}\psi(\cdot)\,\wh v_{k\psi}(x)$$ with respect to the Lebesgue measure on $\mb R^{d_Y}$. By using $\sum_{k\in \wh{\m K}} \big[\wh G_{[k]}(\,\cdot\,,x)\big]_{\#}\wh\nu_\sk(\,\cdot\,|\,x)$ as an estimator for $\mu^*_{Y|x}$, we define the plug-in estimator for $\mb{E}_{\mu^*_{Y|x}}[f_J^\perp(Y)]$ as
 \begin{equation*}
 \sum_{k\in \wh{\m K}} \int_{\mb R^{d_Y}} f_J^{\perp}\big(\wh G_{[k]}(z,x)\big) \bigg\{\sum_{j=0}^J\sum_{\psi\in \Psi_j^{d_Y}}\psi(z)\,\wh v_{k\psi}(x)\bigg\}\,\dd z.
\end{equation*}

\subsubsection{Convergence rate of the estimator for $\mb{E}_{\mu^*_{Y|x}}[f(Y)]$}

   For any $\m L^2$ integrable function $f: \, \mb R^{D_Y} \to \mb R$ and any $x \in \mb R^{D_X}$, our estimator $\wh{\m J}(f, x)$ for $\mb{E}_{\mu^*_{Y|x}}[f(Y)]$ is constructed by combining the estimators for the coarse-scale and fine-scale components,
\begin{equation}\label{whJgamma>0}
\begin{aligned}
       & \wh{\m J}(f,x)=\underbrace{\sum_{j=0}^J\sum_{\psi\in \Psi_j^{D_Y}} 2^{\frac{j(D_Y-d_Y)}{2}} f_{\psi}  \, \wh S_j(\psi,x)}_{\text{estimator for $\mb{E}_{\mu^*_{Y|x}}[f_J(Y)]$}} \\
       &\quad+\ \underbrace{\sum_{k\in \wh{\m K}} \int_{\mb R^{d_Y}} f_J^{\perp}\big(\wh G_{[k]}(z,x)\big) \bigg\{\sum_{j=0}^J\sum_{\psi\in \Psi_j^{d_Y}}\psi(z)\,\wh v_{k\psi}(x)\bigg\}\,\dd z}_{\text{estimator for  $\mb{E}_{\mu^*_{Y|x}}[f^{\perp}_J(y)]$}}\,,\\
       &\text{where}\ \ f_{\psi}=\int_{\mb R^{D_Y}} f(y)\,\psi(y)\,\dd y \ \ \mbox{ and }\ \ f^{\perp}_J(y)=f(y)-\sum_{j=0}^J \sum_{\psi\in  \ov\Psi_j^{D_Y}} f_{\psi}\, \psi(y).
\end{aligned}
\end{equation}
Suppose $\big\{(X_i, Y_i)\big\}_{i=1}^n$ are $n$ i.i.d.~samples from $\mu^*$, and let $\m P_2^*$ and $\m P_3^*$ denote the target distribution classes defined in Theorem~\ref{th2:lower1} and Theorem~\ref{th2:lower}, respectively. 
The following theorem shows that, by setting $J=\big\lceil {\frac{1}{2\alpha_Y+d_Y+d_X \alpha_Y/\alpha_X}}\cdot \log_2 (\frac{n}{\log n})\big\rceil$, there exist suitable choices of $\m G$, $\m S_j^\dagger$, and $\ms S_j$ for Regime 2 ($\mu^* \in \m P_2^*$) and Regime 3b ($\mu^* \in \m P_3^*$) such that the estimator $\wh{\m J}$ \emph{simultaneously} achieves the minimax upper bound for all $\gamma > 0$.

\begin{theorem}[Convergence rates for distribution regression estimators in Regimes 2 and 3b] \label{th:combined}

For Regimes 2 and 3b, there exist distinct families $\m G$ and $\{\m S_j^\dagger\}_{j=0}^J$ tailored for each regime, alongside families $\{\ms S_j\}_{j=0}^J$ that are applicable to both regimes, so that for any $\mu^* = \mu^*_X \mu^*_{Y|X} \in \m P^*_i$ (where $i = 2$ for Regime 2, and $i = 3$ for Regime 3b), the following holds with probability at least $1 - n^{-1}$: for any $\gamma > 0$, the conditional expectation functional estimator $\wh{\m J}$ defined in~\eqref{whJgamma>0} satisfies
\begin{equation*}
    \begin{aligned}
    & \mb{E}_{\mu^*_X} \bigg[ \underset{f\in \m H_1^{\gamma}(\mb R^{D_Y})}{\sup} \Big|\,\wh{\m J}(f,x) \,-\, \mb{E}_{\mu^*_{Y|X}} f(y)\,\Big|\bigg] \\
  &\leq C_{\gamma}  \left\{
    \begin{array}{cc}
     \bigg(  (\log n)^3\cdot n^{-\frac{\alpha_X}{2\alpha_X+d_X}}\, +\,\Big(\frac{n}{\log n}\Big)^{-\frac{\alpha_Y+\gamma}{2\alpha_Y+d_Y+\frac{\alpha_Y}{\alpha_X}d_X}} \, +\, n^{-\frac{\gamma}{\frac{d_Y}{\beta_Y}}}\,\bigg),    & \text{for } i=2, \\
    \bigg(  (\log n)^3\cdot n^{-\frac{\alpha_X}{2\alpha_X+d_X}}\, +\,\log n \cdot \Big(\frac{n}{\log n}\Big)^{-\frac{\alpha_Y+\gamma}{2\alpha_Y+d_Y+\frac{\alpha_Y}{\alpha_X}d_X}} \, +\, \log n\cdot n^{-\frac{\gamma}{\frac{d_Y}{\beta_Y}+\frac{d_X}{\beta_X}}}\,\bigg),      & \text{for } i=3,
    \end{array}
    \right.
    \end{aligned}
\end{equation*}
 for some constant  $C_{\gamma}$   independent of $n$.  
\end{theorem}

\noindent The proof of Theorem~\ref{th:combined} is provided in Appendix~\ref{proofth:thgamma>0} (for Regime 2) and Appendix~\ref{proofth3} (for Regime 3b).  The estimator $\wh{\m J}(f, x)$ leverages the strengths of density regression performed in both the ambient space and the latent space.  By considering the wavelet expansion of $\m L^2$-integrable functions, the task of estimating $\mb{E}_{\mu^*_{Y|x}}[f(Y)]$ for $f \in \m L^2(\mb R^{D_Y})$ becomes equivalent to jointly estimating the coefficients $\mb{E}_{\mu^*_{Y|x}}[\psi(y)]$ over $\psi \in \bigcup_{j \geq 0} \ov \Psi_j^{D_Y}$. Moreover, when $f \in \m H^{\gamma}_1(\mb R^{D_Y})$, the collection of wavelet coefficients $\big\{\mb{E}_{\mu^*_{Y|x}}[\psi(y)]:\, \psi \in \ov \Psi_j^{D_Y}\big\}$ contribute to $\mb{E}_{\mu^*_{Y|x}}[f(Y)]$ with different levels of importance, depending on the resolution level $j$ and the smoothness parameter $\gamma$.
Notably, the difficulty of jointly estimating the coefficients over $\Psi_j^{D_Y}$ decreases as $j$ becomes smaller, due to the lower complexity of the basis functions at coarse scales. This property can be exploited in density regression over the ambient space by using a joint mean regression strategy, with a function class $\m S_j^{\dagger}$ of reduced complexity selected for lower levels $j$. Accordingly, the conditional density estimator $\wh u(y\,|\,x)$ defined in~\eqref{eq:cde} is particularly effective for estimating $\mb{E}_{\mu^*_{Y|x}}[f(Y)]$ when $f$ is smooth (i.e., for large $\gamma$). Specifically, by choosing $J=\big\lceil \frac{1}{d_Y}\log_2(\frac{n}{\log n})\big\rceil$, there exists an appropriate choice of the function families $\big\{\m S_j^{\dagger}\big\}_{j \in {0} \cup [J]}$ such that the estimator $\wh u(y\,|\,x)$ achieves the minimax upper bound for all $\gamma \geq \frac{d_Y \alpha_X}{2\alpha + d_X}$ under Regime 2. This result is detailed in Theorem~\ref{th:thgamma>1} in Appendix~\ref{R2gamma>1}.

However, without explicitly estimating the manifold, this approach integrates manifold estimation and conditional density estimation into a single process of joint mean regression. While efficient, it may overlook finer local details of the supporting manifolds, and can fail to achieve minimax optimality for small $\gamma$, where the loss $d_\gamma$ becomes more sensitive to such fine-scale structures and misalignments arising from manifold estimation. In contrast, density regression in the latent space---augmented by an explicit manifold estimation step---can achieve the minimax rate when $\gamma$ is small. Specifically,  by setting  $J=\big\lceil {\frac{1}{2\alpha_Y+d_Y+d_X \alpha_Y/{\alpha_X}}}\cdot \log_2 (\frac{n}{\log n})\big\rceil$, the mixture of generative models given by~\eqref{eq:mgm} serving as an estimator for the conditional distribution $\mu^*_{Y|x}$, can simultaneously achieve minimax optimality for all $\gamma \leq 1$ under Regime 2, up to logarithmic factors, as detailed in Theorem~\ref{thgamma<1} of Appendix~\ref{R2gamma<1}.

On the other hand, for large values of $\gamma$, this encoder–decoder–based manifold estimation approach fails to fully exploit the higher-order smoothness of the test functions. In such cases, better convergence rates are achievable through alternative strategies.
The estimator defined in~\eqref{whJgamma>0} addresses this trade-off by combining the strengths of both approaches: it uses density regression in the ambient space to estimate the coarse-scale component, while employing density regression in the latent space to recover finer-scale details.

\section{Discussion}
In this paper, we explored the minimax rate of distribution regression under a non-parametric setting, where both the response variable and the covariate may exhibit low-dimensional structures. Our analysis extended to settings in which the conditional response space varies with the covariate, thereby generalizing the classical manifold estimation and support recovery problems into a manifold regression framework.
The minimax rates derived for manifold regression rely on regularity assumptions in the covariate space, including the condition that the covariate density is bounded away from zero. An important direction for future work is to explore the possibility of relaxing or eliminating these assumptions, either through more refined analytical techniques or by adopting weaker evaluation metrics.  Additionally, the rate-optimal procedure for distribution regression developed in this work is primarily theoretical, designed to establish the minimax upper bound. Developing a computationally efficient algorithm that achieves similar statistical guarantees remains an open challenge. For example, our procedure employs density regression in the ambient space to capture global structure and in the latent space to resolve fine-scale details. Given the conceptual similarities between our multiscale approach and techniques used in forward-backward diffusion models \cite{song2020score,ho2020denoising}, it would be worthwhile to investigate whether ideas from our estimator could enhance score-based generative models~\cite{chen2022sampling,oko2023diffusion,tang2024conditionaldiffusionmodelsminimaxoptimal}. Specifically, one could envision a new class of diffusion-based models that estimate global structure in the conditional distribution using diffusion processes in the ambient space~\cite{song2020score}, while capturing fine-scale features via latent diffusion methods~\cite{Rombach_2022_CVPR}.

 \bibliographystyle{plainnat} 
\bibliography{ref.bib,refb.bib}


\newpage
\appendix
\begin{center}
    \Large \textbf{Supplementary Materials to ``Minimax Optimal Rates for Regression on Manifolds and Distributions"}
\end{center}

\medskip
\noindent\textbf{Notation:} We adopt the notations in the manuscript, and further introduce the following additional notations for the technical proofs.
For a set $U\subset \mb R^d$, we write $\mathbf{1}_U(x)$ the indicator function of  $x\in U$.  For two vectors $a,b\in \mb R^d$, we use $\|a-b\|=\sqrt{\sum_{j=1}^d (a_i-b_i)^2}$ to denote the Euclidean distance between them.  For two sequences $\{a_n\}$ and $\{b_n\}$,  the notations $a_n \lesssim b_n$ and $a_n \gtrsim b_n$ imply $a_n \leq Cb_n$ and $a_n \geq C b_n$, respectively, for some constant $C>0$ independent of $n$. Additionally, $a_n \asymp b_n$ indicates that both $a_n \lesssim b_n$ and $a_n\gtrsim b_n$ are hold.  For sequences  $\{a_n\}$, $\{b_n\}$, $\{c_n\}$ We write $a_n=b_n+\m O(c_n)$ if $\|a_n-b_n\|\lesssim c_n$.  For a function $f:\mb R^{d_1}\times \mb R^{d_2}$, we use $\mathbf{J}_f(x)$ to denote the Jacobian matrix of $f$ evaluate at $x$, so that the $(i,j)$ element of  $\mathbf{J}_f(x)$  is $\frac{\partial f_i(x)}{\partial x_j}$.  We denote the $d$-dimensional zero vector as $\mathbf{0}_d$ and may omit the subscript $d$ when it does not lead to ambiguity. 
For a function $f: U\to \mb R$, we use ${\rm supp}(f)=\{x\in U\,:\, f(x)\neq 0\}$ to denote the support of $f$.
\section{Omitted Definitions and Results in Main Text}

\subsection{Smooth Submanifold}\label{sec:ssmanifold}
This subsection provides an introduction to Riemannian submanifolds, Intuitively speaking, a manifold is a topological space that locally resembles the Euclidean space. A submanifold in the ambient space $\mb R^D$ can be viewed as a nonlinear ``subspace'' and is formally defined as follows.

\begin{definition}[Submanifold]
 A subset $\m M$ of $\mathbb{R}^D$ is a $d$-dimensional Riemannian submanifold if for every point $x$ in $\m M$, there exists a neighbourhood $V$ of $x$ on $\m M$ and an open set $U \subseteq \mathbb{R}^d$, such that that there exists a homeomorphism  $\xi$ that maps $U$ to $V$, that is, $\xi: \, U\rightarrow V$ is bijective and both $\xi$ and $\xi^{-1}$ are continuous maps. Moreover,  the differential $D_y \xi$ of $\xi(\cdot)$ at $y$ exists and be injective for every $y\in U$.\footnote{Here, the differential of $\xi(\cdot)$ at $y$, denoted as $D_y \xi$, is a linear map defined by $D_y \xi [v] = \lim_{t \to 0} \frac{\xi(y + tv) - \xi(y)}{t} = J_{\xi}(y)v$ for $v \in \mb R^d$. The injectiveness of $D_y \xi$ is equivalent to the Jacobian matrix $J_{\xi}(y)$ having full rank.} 
We call $(V,\xi)$ a local coordinate chart of $\m M$ near $x$,  and $\xi$ a coordinate map around $x$. We refer to $D$ as the ambient dimension and $d$ as the intrinsic dimension of $\m M$. 
\end{definition}

\begin{definition}[Atlas]
     A collection of $d$-dimensional charts $\ms A = \{(U_\lambda, \varphi_\lambda)\}_{\lambda\in \Lambda}$ is called an atlas on $\m M$ if 1. $\m M = \bigcup_{\lambda\in\Lambda} U_\lambda$. 2. Each chart $(U_{\lambda}, \varphi_{\lambda})$ in atlas $\ms A$ consists of a homeomorphism $\varphi_{\lambda}:\, U_{\lambda} \to \widetilde U_{\lambda}$,  from an open set $U_{\lambda} \subset \m M$ to an open set $\widetilde U_{\lambda} \subset \mb R^d$. 3. Any two charts $(U, \varphi)$ and $(V, \psi)$ in atlas $\mathscr{A}$ are compatible, meaning that the transition map $\varphi \circ \psi^{-1}: \psi(U \cap V) \rightarrow \varphi(U \cap V)$ is a diffeomorphism.
\end{definition}

The tangent space $T_{\theta}\m M$ is the linearization of $\m M$ at $\theta$. When $\m M$ is an embedded submanifold of a Euclidean space $\mb R^D$, the tangent spaces of $\m M$ are linear subspaces of $\mb R^{D}$ that pass through the origin and have dimensions that equal the intrinsic dimension $d$ of $\m M$. The formal definition is as follows.

\begin{definition}[Tangent space]  For a submanifold $\m M$ embedded in $\mb R^D$, we denote the tangent space of $\m M$ at $\theta$ as $T_{\theta}\m M=\{c'(0)\,|\, c: I\to \m M \text{ is }C^{1}\text{-smooth around } 0 \text{ and } c(0)=\theta\}$, where $I$ is any open interval containing $t=0$.  That is, $v$ is in $T_{\theta}\m M$ if and only if there exists a smooth curve on $\m M$ passing through $x$ with velocity $v$. Vectors in $T_{\theta}\m M$ are called tangent vectors to $\m M$ at $\theta$.  The collection $T \m M=\{(\theta,v)\,:\, \theta\in \m M, v\in T_{\theta} \m M\}$ is called the tangent bundle of $\m M$.
\end{definition}
 
To manage multiple local coordinate charts in the underlying data manifold representation, we will use the mathematical technique of \emph{partition of unity} as defined below.

\begin{definition}[partition of unity]
 A partition of unity subordinate to altas $\ms A=\{(U_{\lambda},\varphi_{\lambda})\}_{\lambda \in \Lambda}$ is a collection of smooth  functions $\{\rho_\lambda\}_{\lambda\in \Lambda}$ on $\mathcal M$ so that \begin{enumerate} 
     \item $0\leq \rho_\lambda\leq 1$ for all $\lambda\in \Lambda$, and $\sum_{\lambda\in \Lambda}\rho_\lambda(x)=1$ for all $x\in \mathcal{M}$.
     \item ${\rm supp}(\rho_\lambda)\subset U_{\lambda}$ for any $\lambda\in\Lambda$.
     \item Each point $x\in \mathcal{M}$ has a neighborhood which intersects ${\rm supp}(\rho_\lambda)$ for only finitely many $\lambda\in \Lambda$.
 \end{enumerate}
 \end{definition}

  Using the partition of unity, one can glue constructions in the local charts to form a global construction on the manifold. Such a global construction usually does not rely on the choice of the partition of unity. Conversely, the partition of unity enables the decomposition of a global estimation problem into local ones, which resembles the data localization in local (polynomial) regression~\citep{loader2006local,bickel2007local}. 
\begin{definition}[Riemannian volume measure of submanifold]
Suppose $\ms A=\{(U_{\lambda},\varphi_{\lambda})\}_{\lambda \in \Lambda}$ is an atlas on  a submanifold $\m M$ and  $\{\rho_\lambda\}_{\lambda\in \Lambda}$ is a  partition of unity subordinate to altas $\ms A$. Then 
the Riemannian volume measure $\mu_{\m M}$ can be written as
 \begin{equation*}
     \,\dd \mu_{\m M}=\sum_{\lambda\in \Lambda} \rho_{\lambda}(\varphi_\lambda^{-1}(z))\sqrt{{\rm det}(J_{\varphi_\lambda^{-1}}(z)^TJ_{\varphi_\lambda^{-1}}(z))}\,\dd z,
 \end{equation*}
 where $\dd z$ is the Lebesgue measure on $\mb R^{d}$.
 A measure $\mu$  on $\m M$ is said to have a density $f$ (with respect to the volume measure $\mu_{\m M}$) if for any measurable subset $A \subset \m M$, 
 \begin{equation*}
     \mu(A)=\int_{A} f \,\dd  \mu_{\m M}=\sum_{\lambda\in \Lambda}\int_{\varphi_\lambda(U_{\lambda}\cap A)} \rho_{\lambda}(\varphi_\lambda^{-1}(z))\cdot f(\varphi_\lambda^{-1}(z))\sqrt{{\rm det}(J_{\varphi_\lambda^{-1}}(z)^TJ_{\varphi_\lambda^{-1}}(z))}\,\dd z.
 \end{equation*}
\end{definition}
Note that the Riemannian volume measure and the density function with respect to it are independent of the choice of atlas and partition of unity.
\begin{definition}[Reach]
     The reach of a closed subset $\m M \subset \mathbb{R}^D$ is defined as
\begin{equation*}
 r_M=\sup \Big\{\varepsilon\, \Big|\, \forall x \in  \m M^{\varepsilon}, \text{ there exists unique }y \in \m M, \text { so that } \operatorname{dist}(x, \m M)=\|x-y\|\Big\},
\end{equation*}
where  ${\rm dist}(z,\m M)=\inf_{p\in \m M}\|p-z\|$ denotes the distance function to $\m M$, and $\m M^{\varepsilon}  =\left\{x \in \mathbb{R}^D: \operatorname{dist}(x, \m M)<\varepsilon\right\}$ is the $\varepsilon$-offset of $\m M$. 
\end{definition}
A lower bound on the reach  prevents the manifold from becoming nearly self-intersecting and ensures a uniform upper bound on its curvature.  We also restate the definition of a  $\beta$-smooth submanifold  as described in Definition~\ref{def:manifold1} of the main text for completeness.

\begin{definition*}[$\beta$-Smooth submanifold]
 A $d$-dimensional submanifold $\m M$ in $\mb R^D$ is said to belong to the manifold class $\ms{M}_{\tau,\tau_1,L}^{\beta}(d,D)$ if: 1.~$\m M$ is closed; 2.~it has reach larger than $\tau$; and 3.~for all $y_0 \in \m M$, there exists a neighborhood $U_{y_0}$ of $y_0$ on $\m M$ so that the projection 
  $\wt\pi_{y_0}:\m M \rightarrow T_{y_0} \mathcal{M}$ defined by $\wt\pi_{y_0}(y)=\operatorname{Proj}_{T_{y_0} \mathcal{M}}(y-y_0)$,  when  {restricted to $U_{y_0}$}, is a  diffeomorphism, with  inverse function $\phi_{y_0}$ defined on $\mb B_{T_{y_0}{\m M}}(0,\tau_1)$, and  {$\phi_{y_0}\in \m H^{\beta}_{L,D}(\mb B_{T_{y_0}{\m M}}(0,\tau_1))$}.
\end{definition*}

\noindent\textbf{Geometric Properties of $\beta$-smooth submanifolds with positive reach:} (see for example, Lemma 20 of~\cite{divol2022measure})Suppose $\m M\in \ms{M}_{\tau,\tau_1,L}^{\beta}(d,D)$ with $\beta\geq 2$. Then 
\begin{enumerate}
    \item If $h\leq \frac{\tau}{4}$, then there exist some constants $(c,C)$ so that for any $x\in \m M$, 
    \begin{equation*}
        c\, h^d\leq {\rm vol}_{\m M}(\mb B_{\m M}(x,h))\leq C\, h^d,
    \end{equation*}
    where ${\rm vol}_{\m M}$ denotes the volume measure of $\m M$.
    \item For any $h\leq r_0=\tau_1\wedge ((\tau\wedge L)/4)$ and $x\in \m M$, $\mb B_{\m M}(x,h)\subset \phi_x\big(\mb B_{T_x \m M}(\mathbf 0,h)\big)\subset \mb B_{\m M}(x,8h/7)$.
\item If ${\rm Proj}_{\m M}(z)=x$ for some $z$ satisfying ${\rm dist}(z,\m M)<\tau$, then $z-x\in T_x\m M^{\perp}$.
\end{enumerate}

\subsection{Smooth submanifold family and smooth conditional distributions}\label{app: regularitymanifold}
Firstly we recall the definition of  $(\beta_Y,\beta_X)$-smooth manifold family defined in Definition~\ref{def:manifold} of the main text.

\begin{definition*}[$(\beta_Y,\beta_X)$-smooth submanifold family] 
A submanifold family $\big\{\mathcal{M}_{Y|x}:\,x \in \mathcal{M}_X\big\}$ is said to belong to $\ms M^{\beta_Y\,\beta_X}_{\tau,\tau_1,L}(d,D,\m M_X)$, if 
for any $x\in \m M_X$: 1.~the manifold $\m M_{Y|x}$ is a closed $d$-dimensional submanifold in  $\mb R^D$; 2.~it has reach larger that $\tau$; and 3.~if, for any $w_0=(x_0,y_0)\in \m M$, there exists a neighborhood $U_{\omega_0}$ of $y_0$ on $\m M_Y$, so that for any $x\in \mb B_{\m M_X}(x_0,\tau)$, the function $\wt \pi_{w_0}: \m M_Y  \rightarrow T_{y_0} \mathcal{M}_{Y|x_0}$ defined by  $\wt \pi_{w_0}(y)=\operatorname{Proj}_{T_{y_0} \mathcal{M}_{Y|x_0}}\left(y-y_0\right)$, when restricted to $U_{\omega_0}\cap \m M_{Y|x}$, is a  diffeomorphism with  inverse function $\phi_{\omega_0,x}(\cdot)$ defined on $\mb B_{T_{y_0}{\m M_{Y|x_0}}}(0,\tau_1)$. Moreover, the function $\Phi_{\omega_0}: \mb B_{T_{y_0}{\m M_{Y|x_0}}}(0,\tau_1)\times  \mb B_{\m M_X}(x_0,\tau)\to \mb R^{D_Y}$ defined as {$\Phi_{\omega_0}(z,x)=\phi_{\omega_0,x}(z)$ belongs to $ {\m H}^{\beta_Y,
 \beta_X}_{L,D_Y}(\mb B_{T_{y_0}{\m M_{Y|x_0}}}(0,\tau_1), \mb B_{\m M_X}(x_0,\tau))$}. 
    
\end{definition*} 

 \medskip
\noindent We have the following lemma which provides an equivalent definition of $(\beta_Y,\beta_X)$-smooth manifold family, whose proof is given in Appendix~\ref{proof:lemmamanifoldpro}.

\begin{lemma}\label{le:defmanifold}(Properties of Smooth submanifold family)
  Suppose $\beta_Y\geq 2$ and $\beta_Y\geq \beta_X$. Consider a submanifold faimly $\big\{\mathcal{M}_{Y|x}:\,x \in \mathcal{M}_X\big\}$, the following statements are equivalent:
    \begin{enumerate}
        \item  There exist constants $(\tau,\tau_1, L)$ so that $\{{\m M}_{Y|x}\}_{x\in \m M_X}\in  {\ms M}^{\beta_Y\,\beta_X}_{\tau,\tau_1,L}(d_Y,D_Y,\m M_X)$.
        \item (Existence of $x$-dependent $\m H^{\beta_Y,\beta_X}$-smooth 
  local charts)  There exist constants $(\wt\tau,\wt \tau_1,\wt L)$ so that for any $w_0=(x_0,y_0)\in \m M$, there exists a neighborhood $\wt U_{y_0}$ of $y_0$ on $\m M_{Y}$ such that  for any $x\in \mb B_{\m M_X}(x_0,\wt\tau)$,  it holds that $\mb B_{\m M_{Y|x}}(y_0,\wt\tau)\subset \wt U_{y_0}\cap \m M_{Y|x}\subset \mb R^{D_Y}$ and there exists a uniformly $\wt L$-Lipschitz diffeomorphism $\wt Q_{\omega_0}(\cdot,x)$ that maps $ \wt U_{y_0}\cap \m M_{Y|x}$ to $\mb B_{\mb R^{d_Y}}(\mathbf{0},\wt\tau_1)$ with inverse denoted as $\wt g_{\omega_0,x}(\cdot)$, so that $\wt Q_{\omega_0}(y_0,x_0)=\mathbf{0}$ and the function $\wt  G_{\omega_0}: \mb B_{\mb R^{d_Y}}(\mathbf{0},\wt\tau_1)\times \mb B_{\m M_X}(x_0,\wt\tau)\to \mb R^{D_Y}$ defined as $ \wt G_{\omega_0}(z,x)=\wt g_{\omega_0,x}(z)$ satisfies that $ \wt G_{\omega_0}\in  {\m H}^{\beta_Y,
 \beta_X}_{\wt L,D_Y}(\mb B_{\mb R^{d_Y}}(\mathbf{0},\wt\tau_1), \mb B_{\m M_X}(x_0,\wt\tau))$. 
        \item (Solution manifold with $\m H^{\beta_Y,\beta_X}$-smooth defining functions) There exist constants $(\ov \tau,\ov \tau_1,\ov L)$ so that $\m M_Y\subset \mb B_{\mb R^{D_Y}}(\mathbf{0},\ov L)$ and for any $\omega_0=(x_0,y_0)\in \m M$, there exists a function $F_{\omega_0}\in  {\m H}^{\beta_Y,
 \beta_X}_{\ov L,D_Y-d_Y}(\mb B_{\mb R^{D_Y}}(y_0,\ov \tau),\mb B_{\m M_X}(x_0,\ov \tau))$ so that for any $x\in  \mb B_{\m M_X}(x_0,\ov \tau)$, it holds that $\mb B_{\m M_{Y|x}}(y_0,\ov\tau)=\{y\in\mb B_{\mb R^{D_Y}}(y_0,\ov \tau):\, F_{\omega_0}(y,x)=\mathbf 0\}$,  and for any $(x,y)\in \mb B_{\m M_X}(x_0,\ov\tau)\times \mb B_{\mb R^{D_Y}}(y_0,\ov \tau)$, it holds that $J_{F_{\omega_0}(\cdot,x)}(y)J_{F_{\omega_0}(\cdot,x)}(y)^T\succeq  \ov\tau_1 I_{D_Y-d_Y}$.

    \end{enumerate}
\end{lemma}

As a crucial intermediate result for proving Lemma~\ref{le:defmanifold}, the following lemma states that if $x$-dependent $\mathcal{H}^{\beta_Y,\beta_X}$-smooth local charts exist, then for an appropriate choice of $V$, the function $V^T(\cdot-y_0)$, when restricted to $\mathcal{M}_{Y|x}$, will be locally invertible around $y_0$.

\begin{lemma}\label{le:projection}

Suppose the family of submanifolds \(\{\mathcal{M}_{Y|x} : x \in \mathcal{M}_X\}\) meets the conditions specified in Point 2 of Lemma~\ref{le:defmanifold}, with \(\beta_Y \geq \max(2, \beta_X)\). For any \(\omega_0 = (x_0, y_0) \in \mathcal{M}\), consider \(P_{\omega_0}\) as the projection matrix onto \(T_{\mathcal{M}_{Y|x_0}} y_0\) and let \(V_{\omega_0} \in \mathbb{R}^{D_Y \times d_Y}\) be an arbitrary orthonormal matrix such that \(V_{\omega_0}^T P_{\omega_0} V_{\omega_0} \succeq \tau_0 I_{d_Y}\) for some positive constant \(\tau_0\). 
Then, there exist constants \((\tau, \tau_1, L)\) such that for any \(\omega_0 = (x_0, y_0) \in \mathcal{M}\), there is a subset \(U_{\omega_0}\) of \(\mathcal{M}_Y\) satisfying the following conditions:
1. For any \(x \in \mathbb{B}_{\mathcal{M}_X}(x_0, \tau)\), \(\mathbb{B}_{\mathcal{M}_{Y|x}}(y_0, \tau) \subset U_{\omega_0} \cap \mathcal{M}_{Y|x}\).
2. The function \(V_{\omega_0}^T(\cdot - y_0)\), when restricted to domain \(U_{\omega_0} \cap \mathcal{M}_{Y|x}\), is a diffeomorphism onto its image, with the inverse function denoted by \(g_{\omega_0, x}\), defined on \(\mathbb{B}_{\mathbb{R}^{d_Y}}(\mathbf{0}, \tau_1)\).
3. The function \(G_{\omega_0}: \mathbb{B}_{\mathbb{R}^{d_Y}}(\mathbf{0}, \tau_1) \times \mathbb{B}_{\mathcal{M}_X}(x_0, \tau) \to \mathbb{R}^{D_Y}\), defined by \(G_{\omega_0}(z, x) = g_{\omega_0, x}(z)\),belongs to \(\mathcal{H}^{\beta_Y, \beta_X}_{L, D_Y}(\mathbb{B}_{\mathbb{R}^{d_Y}}(\mathbf{0}, \tau_1), \mathbb{B}_{\mathcal{M}_X}(x_0, \tau))\).
\end{lemma}
The proof of Lemma~\ref{le:projection} is given in Appendix~\ref{proof:le:projection}.

 \medskip
 \noindent For ease of notation, we make the following definition to the smooth conditional distributions on submanifolds.

\begin{definition} (Smooth conditional distributions)
The conditional distribution $\{\mu^*_{Y|x}\}_{x\in \m M_X}$ is said to be inside ${\ms C}^{\beta_Y,\beta_X,\alpha_Y,\alpha_X}_{\tau,\tau_1,L}(d_Y,D_Y,\m M_X)$ if for any $x\in \m M_X$, $\mu^*_{Y|x}$ is supported on a submanifold $\m M_{Y|x}$ and has a density function $u^*(\,\cdot\,|\,x)$ with respect to the volume measure of $\m M_{Y|x}$ so that $\{\m M_{Y|x}\}_{x\in \m M_X}\in {\ms   M}^{\beta_Y\,\beta_X}_{\tau,\tau_1,L}(d_Y,D_Y,\m M_X)$ and  there exists a function $\ov u^*\in  {\m H}^{\alpha_Y,\alpha_X}_{L}(\mb R^{D_Y},\mb R^{D_X})$  so that $u^*(y|x)=\ov u^*(y,x)$ for any $(x,y)\in \m M$. 
  \end{definition}

\noindent The following lemma, whose proof is given in Appendix~\ref{proof:lemmadensitypro}, shows that the smoothness of the density function of $\mu^*_{Y|x}$ w.r.t. the volume measure of $\m M_{Y|x}$ is equivalent to the smoothness of the latent distributions defined through the $x$-dependent local charts of the submanifolds. 
  
\begin{lemma}\label{le:defdistribution}
(Equivalence between smoothness of density function  and smoothness of latent distribution) Consider the conditional distribution $\{\mu^*_{Y|x}\}_{x\in \m M_X}$ supported on $\{\m M_{Y|x}\}_{x\in \m M_X}\in {\ms M}^{\beta_Y\,\beta_X}_{\tau,\tau_1,L}(d_Y,D_Y,\m M_X)$ with $\beta_Y\geq 2$ and $\beta_Y\geq \beta_X$, then for any $\alpha_Y,\alpha_X> 0$ satisfying $\alpha_Y\geq \alpha_X$, $\beta_Y\geq  \alpha_Y+1 $ and $\beta_X\geq \alpha_X+\frac{\alpha_X}{\alpha_Y}$, we have 
\begin{enumerate}
    \item If for any $\omega_0=(x_0,y_0)\in \m M$, the push-forward measure $[\operatorname{Proj}_{T_{y_0} \mathcal{M}_{Y|x_0}}\left(\cdot-y_0\right)]_{\#}(\mu^*_{Y|x}|_{U_{\omega_0}\cap \m M_{Y|x}})$\footnote{Here we have adopted the notation introduced in the definition  of $(\beta_Y,\beta_X)$-smooth submanifold family.} exists with a density function with respect to the volume measure of  
    $T_{\m M_{Y|x_0}} y_0$, denoted as $\nu_{\omega_0}(\cdot|x)$, and it satisfies that $\nu_{\omega_0}(z,|,x)\in  {\m H}^{\alpha_Y,\alpha_X}_{L}(\mb B_{T_{\m M_{Y|x_0}} y_0}(\mathbf{0},\tau_1),\mb B_{\m M_X}(x_0,\tau) )$. Then there exists $L'$ so that  $\{\mu^*_{Y|x}\}_{x\in \m M_X}\in {\ms C}^{\beta_Y,\beta_X,\alpha_Y,\alpha_X}_{\tau,\tau_1,L'}(d_Y,D_Y,\m M_X)$.

 \item  If  $\{\mu^*_{Y|x}\}_{x\in \m M_X}\in {\ms C}^{\beta_Y,\beta_X,\alpha_Y,\alpha_X}_{\tau,\tau_1,L}(d_Y,D_Y,\m M_X)$. Then there exists a constant $L'
$ so that for any $\omega_0=(x_0,y_0)\in \m M$ and any $\wt Q_{\omega_0}$ that satisfies the conditions specified in Point 2 of Lemma~\ref{le:defmanifold}, 
the density of the push forward measure $[\wt Q_{\omega_0}(\cdot,x)]_{\#}(\mu^*_{Y|x}|_{\wt U^{\omega_0}_{Y|x}})$ with respect to the Lebesgue measure on $\mb R^{d_Y}$, denoted as $\wt\nu_{\omega_0}(\cdot|x)$, exists and satisfies  $\wt\nu_{\omega_0}(z,|,x)\in   {\m H}^{\alpha_Y,\alpha_X}_{L'}(\mb B_{\mb R^{d_Y}}(\mathbf{0},\tau_1),\mb B_{\m M_X}(x_0,\tau) )$.
\end{enumerate}
\end{lemma}

 \noindent The following result shows that smooth conditional distributions on submanifold can be expressed as mixture of  conditional generative models.

 \begin{lemma}{(Expressing $\mu^*_{Y|X}$ as mixture of  conditional generative models)}\label{legenerative}
Suppose $\{\mu^*_{Y|x}\}_{x\in \m M_X}\in {\ms C}^{\beta_Y,\beta_X,\alpha_Y,\alpha_X}_{\tau,\tau_1,L}(d_Y,D_Y,\m M_X)$.
 For any $\tau_2$ with $0<\tau_2\leq (\tau\wedge \tau_1)/4$, let $\{(x_k^*,y_k^*)\}_{k=1}^{K^*}\subset \m M$ be a $\tau_2$-covering set of $\m M$.  Then  for any $k\in [K^*]$, there exist functions  $G^*_{[k]}\in  {\m H}^{\beta_Y,
 \beta_X}_{L_1,D_Y}(\mb R^{d_Y},\mb R^{D_X})$, $v^*_{[k]}(z,x)\in  {\m H}^{\alpha_Y,\alpha_X}_{L_1}(\mb R^{D_Y},\mb R^{D_X}) $ with some constant $L_1$, such that for any $z\in \mb B_{\mb R^{d_Y}}(\mathbf{0},\tau_1)$ and $x\in \m M_X$,  $v^*_{[k]}(z,x)=0$ if either $\|x-x_k^*\|\geq  \sqrt{2}\tau_2$  or $\|G^*_\sk(z,x)-y_k^*\|\geq  \sqrt{2}\tau_2$. Moreover, for any  $x\in \m M_X$ and any continuous function $g:\m M_{Y|x}\to \mb R$, it holds that
         \begin{equation*}
             \mb{E}_{y\sim\mu^*_{Y|x}}[g(y)]=\sum_{k=1}^{K^*}\int_{\mb B_{\mb R^{d_Y}}(\mathbf{0},\tau_1)} g(G^*_{[k]}(z,x)) v^*_{[k]}(z|x)\,\dd z.
         \end{equation*}

 \end{lemma}
\noindent The proof of Lemma~\ref{legenerative} is given in Appendix~\ref{prooflegenerative}.

\subsection{Wavelet}\label{app:asswav}
In this section, we give a brief introduction to the wavelet. Let $\phi_{\mf M}$ and $\phi_{\mf F}$ be a compactly supported wavelet and scaling function, respectively, for example Daubechies wavelets~\citep{ daubechies1992ten,meyer1992wavelets}. This implies that 
\begin{equation*}
 \left\{
\begin{array}{ll}
   \phi_{\mf F} (x-k) & j=0, k\in \mb Z,  \\
    2^{(j-1)/2} \phi_{\mf M}(2^{j-1}x-k), &  j\in \mb N^{+}, k\in \mb Z,
\end{array}
\right.
\end{equation*}
is an orthonormal basis of $\m L^2(\mb R)$, where we use $\m L^2$ to denote the set of square integrable functions. To obtain a basis of $\m L^2(\mb R^d)$ for an integer $d>1$, set 
\begin{equation*}
     \mf G = \{\mf F,\,\mf M\}^d\setminus \{(\mf F,\ldots,\mf F)\}.
\end{equation*}
Then for any multi-index $k\in \mb Z^d$, the level zero basis $\phi_k^{[d]}$ is obtained by translating the $d$-fold tensor product $\phi_{\mf F}^{\otimes d}$ by $k$ as $\phi_{k}^{[d]}(x) = \prod_{i=1}^d \phi_{\mf F}(x_i-k_i)$ for $x=(x_1,\ldots,x_d)\in\mb R^d$, and for any $j\geq 1$, the level $j$ basis $\big\{\psi_{ljk}^{[d]}:\, l\in[2^d-1]\big\}$ with translation $k$ is any ordering of the following $2^d-1$ functions,
\begin{align*}
    \psi_{gjk}^{[d]}(x)=2^{\frac{d(j-1)}{2}} \,\prod_{i=1}^d \phi_{g_i}^{[d]}\big(2^{j-1}x_i - k_i\big), \quad \forall g\in \mf G.
\end{align*}
This gives the orthornormal basis  
\begin{equation*}
\left\{
\begin{array}{ll}
   \phi_k^{[d]}(x), & j=0,l=0, k\in \mb Z^d,  \\
    \psi_{ljk}^{[d]}(x), &  j\in \mb N^{+},l\in [2^d-1], k\in \mb Z^d.
\end{array}
\right.
\end{equation*}
Denote $\ov\Psi^{d}_0=\{\phi_{k}^\sd(\cdot)\,:\, k\in \mb Z^d\}$ as the set of level zero basis and $\ov\Psi^{d}_j=\{\psi_{ljk}^\sd(\cdot)\,:\, l\in [2^d-1], k\in \mb Z^d\}$ as the set of level $j$ basis for $j\in \mb N^{+}$.   We can define the Besov space $B^s_{p,q}(\mb R^d)$ consists of functions $f$ that admits the wavelet expansion
\begin{equation*}
  f(x)
 = \sum_{j=0}^{\infty}\sum_{\psi\in \ov\Psi_j^d}f_{\psi}\psi(x),
\end{equation*}
where  $f_{\psi}:=\int f(x) \psi(x)\,\dd x$, and is equipped with the norm 
 \begin{equation*}
     \|f\|_{B^s_{p,q}}:= \Big\|2^{js}2^{dj(\frac{1}{2}-\frac{1}{p})}\|f_j\|_{p}\Big\|_{q},
 \end{equation*}
 with  $f_j=\{f_{\psi}\}_{\psi\in \ov\Psi_j^d}$. The following Theorem collects the relationship between the Besov space and H\"{o}lder space.
\begin{theorem}
(Theorem 1.122 of~\cite{Triebel2006} and Proposition 4.3.30 of~\cite{gine_nickl_2015}) Let $\alpha>0$, if $\alpha$ is not integer, then
\begin{equation*}
    \m H^{\alpha}(\mb R^d)=B^{\alpha}_{\infty,\infty}(\mb R^d);
\end{equation*}
if $\alpha$ is integer, then
\begin{equation*}
    B^{\alpha}_{1,\infty}(\mb R^d)\subset \m H^{\alpha}(\mb R^d)\subset B^{\alpha}_{\infty,\infty}(\mb R^d).
\end{equation*}
\end{theorem}
\noindent Focusing on the H\"{o}lder space, we can find a wavelet basis that satisfies the following property.
\begin{lemma}\label{le:wavelet}
  For any positive integer $\alpha$,  there exists an orthonormal basis $\bigcup_{j\geq 0}\ov \Psi_j^d$ for $\m L^2(\mb R^d)$, so that there exist constants $C_R,C_L,C_L',C_L^{\dagger},C_L^{\ddagger},C_W,C_I$ such that for any integer $j\geq 0$,
    \begin{enumerate}
    \item (Regularity) $\sup_{x\in \mb R^d}|\psi^{(l)}(x)|\leq C_R 2^{j|l|+\frac{dj}{2}}$ holds for any $l\in \mb N_0^d$ with $|l|\leq \alpha$ and $\psi\in \ov\Psi_j^{d}$.
        \item (Locality) for any $\psi\in \ov\Psi_j^{d}$, there exists a rectangle $I_{\psi}$ such that 
        \begin{enumerate}
            \item   for any $l\in \mb N_0^d$ with $|l|\leq \alpha$, ${\rm supp}(\psi^{(l)})\subset I_{\psi}$ and  the diameter of $I_{\psi}$ is smaller than $C_L 2^{-j}$
            \item $\sup_{x\in \mb R^d}\sum_{\psi\in \ov \Psi_j^d} \mathbf{1}(x\in I_{\psi})\leq C_L'$
            \item for any $R\geq 1$, $\big|\{\psi\in \ov\Psi_j^{d}\,: \,I_{\psi}\cap \mb B_{\mb R^d}(0,R)\neq \emptyset\}\big|\leq C_L^{\dagger} R \,2^{jd}$.
            \item  for any $j\geq 1$ and $x\in \mb R^d$, $\big|\{\psi\in \ov\Psi_j^{d}\,: \,I_{\psi}\cap \mb B_{\mb R^d}(x,2^{-(j-1)})\neq \emptyset\}\big|\leq C_L^{\ddagger}.$
        \end{enumerate}
        \item (Wavelet coefficients of smooth function) for any $\alpha_1\leq \alpha$, $r>0$ and $f\in \m H^{\alpha_1}_{r}(\mb R^d) $,  it holds for any $\psi\in \ov\Psi_j^{d}$ that the wavelet coefficient $f_{\psi}=\int_{\mb R^{d}} f(x)\psi(x)\,\dd x$ is bounded by $C_W r 2^{-\frac{dj}{2}-j\alpha_1}$ in absolute value.
        \item (Index of Wavelet basis) for any $R'>0$, let $\Psi_j^{d}=\{\psi\in \ov \Psi_j^d:\, {\rm supp}(\psi)\cap \mb B_{\mb R^{d}}(\mathbf{0},R')\neq \emptyset\}$, then $\Psi_j^{d}$ can be written as an index set 
        \begin{equation*}
   \Psi _j^{d}=\{\psi_{j\iota}(\cdot):\, \iota\in \ms I_j\subset [0,1]^{d+1}\},
\end{equation*}
where $\ms I_j$ is $C_I2^{-j}/(R'+C_L)$-separated. 
\end{enumerate}

\end{lemma}
\noindent  The proof of Lemma~\ref{le:wavelet} is provided in Appendix~\ref{proof:le:wavelet}. 
The following lemma presents the wavelet truncation approximation for marginal smooth functions, the proof of which is given in Appendix~\ref{proof:le:approwaveletsmooth}.
 \begin{lemma}\label{le:approwaveletsmooth}
Suppose $f\in \ov{\m H}^{\alpha_1,\alpha_2}_L(\mb R^{d_1},\mb R^{d_2})$. Consider two wavelet basis $\{\ov\Psi_j^{d_1}\}_{j\geq 0}$ and $\{\ov\Psi_j^{d_2}\}_{j\geq 0}$ that both satisfy the properties in Lemma~\ref{le:wavelet} with smoothness $\alpha=\lceil\alpha_1\vee\alpha_2\rceil$ and constants $C_R,C_L,C_L',C_L^{\dagger},C_L^{\ddagger},C_W, C_I$. It holds  for any $x\in \mb R^{d_1}$ and $y\in \mb R^{d_2}$ that
\begin{equation*}
    \begin{aligned}
        \Big|f(x,y)-\sum_{j_1=0}^{J_1}\sum_{j_2=0}^{J_2} \sum_{\psi\in \ov\Psi_{j_1}^{d_1}}\sum_{\psi\in \ov\Psi_{j_2}^{d_2}}   f_{\psi,\phi} \psi(x)\phi(y)\Big|\leq  C_RC_L'C_W L\, 2^{-J_1\alpha_1}+2^{d_1}C_R^3 C_L'{}^2C_W C_L^{d_1} L J_1 2^{-J_2\alpha_2},
        \end{aligned}
 \end{equation*}
 where $f_{\psi,\phi}=\int_{\mb R^{d_2}}\int_{\mb R^{d_1}} f(x,y)\psi(x)\phi(y)\,\dd x\dd y.$
\end{lemma}
\subsection{Matching error for Joint Mean Regression}
In this subsection, we present a general result for bounding the matching error in joint mean regression. This result will be frequently applied in the proofs of the main results that follow. Let $\Lambda$ be a countable set and consider a function class $\{\psi_{\lambda}(\cdot)\}_{\lambda\in \Lambda}$ on $\mb R^{D_Y}$, the joint mean regression aim to find a  $\wh S(\lambda,x)$ that solves
 \begin{equation}\label{jointre}
    \underset{S\in \m S}{\arg\min}\frac{1}{n} \sum_{i=1}^n\sum_{\lambda\in \Lambda } (S(\lambda,X_i)-\psi_{\lambda}(Y_i))^2 ,
 \end{equation}
 where $\m S$ is a suitable approximation family for $S$.  This can be think of using the function $\wh S(\lambda,X)$ that depend both on the index $\lambda$ and the covariate $X$ to form a global estimator to the conditional expectation of $\mb{E}[\psi_{\lambda}(Y)|X]$ over $\lambda\in \Lambda,x\in \m M_X$.  
 We derive the following theorem to study the matching error of the joint mean regression, the proof of which is given in Appendix~\ref{proof:theoremjointregression}.
 \begin{theorem}\label{theoremjointregression}
 Suppose $\{(X_i,Y_i)\}_{i=1}^n$  are $n$ i.i.d samples from $\mu^*=\mu^*_X\mu^*_{Y|X}$ supported on $\m M$.  Consider the estimator $\wh S(\cdot,\cdot)$ defined in~\eqref{jointre}. Assume that there exist positive constants $C,C_1$ so that the following assumptions are satisfied:
   \begin{enumerate}
       \item It holds for any $S\in \m S$ that $\underset{(x,y)\in \m M}{\sup}\sum_{\lambda\in \Lambda }  S(\lambda,x)^2+|\psi_{\lambda}(y)S(\lambda,x)|\leq C$.
       \item Denote $   \ell(x,y,S)=\sum_{\lambda\in \Lambda }   \|S(\lambda,x)\|^2-2\psi_{\lambda}(y)^T S(\lambda,x)$, then for any $S,S'\in \m S$, it holds that 
       \begin{equation*}
           \mb{E}_{\mu^*}\Big[\big(  \ell(X,Y,S)-  \ell(X,Y,S')\big)^2\Big]\leq C\, \mb{E}_{\mu^*_X}\Big[\sum_{\lambda\in \Lambda }  \big(S(\lambda,X)-S'(\lambda,X)\big)^2 \Big].
       \end{equation*}
       \item Define the distance $d_n$ as $d_n(S,S')=\sqrt{\frac{1}{n}\sum_{i=1}^n(\ell(X_i,Y_i,S)-\ell(X_i,Y_i,S'))^2}$ and $\mathbf N(\m S,d_n,\varepsilon)$ be the $\varepsilon$-covering number of $\m S$ with respect to $d_n$, Then, for some terms $W_n,T_n> 1$ that may depend on $n$, it holds for any $0<\varepsilon\leq\sup_{S,S'\in \m S}d_n(S,S')$ that
       \begin{equation*}
         \mathbf N(\m S,d_n,\varepsilon)\leq (\frac{T_n}{\varepsilon})^{W_n}.
       \end{equation*}
   \end{enumerate}
   Then for any constant $c>0$, there exists a constant $C_1$ so that it holds with probability at least $1-n^{-c}$ that 
   \begin{equation*}
   \begin{aligned}
       & \mb{E}_{\mu^*_X}\Big[\sum_{\lambda\in \Lambda }  \Big(\wh S(\lambda,X)-\mb{E}_{\mu^*_{Y|X}}[\psi_{\lambda}(y)] \Big)^2 \Big]\\
        &\leq C_1\frac{W_n(\log n+\log T_n)}{n}+C_1\,\underset{S\in \m S}{\min}\,  \mb{E}_{\mu^*_X}\Big[\sum_{\lambda\in \Lambda }  \Big(S(\lambda,X)-\mb{E}_{\mu^*_{Y|X}}[\psi_{\lambda}(y)] \Big)^2 \Big].
   \end{aligned}
   \end{equation*}
\end{theorem}

\section{Details of Miniax Optimal Estimators}\label{app:estimatordetail}

\subsection{Minimax Optimal Estimator for Regime 1}\label{R1gamma>0}

Consider a wavelet basis $\bigcup_{j\geq 0}\ov\Psi_j^{D_Y}$ that satisfies the properties stated in Lemma~\ref{le:wavelet}, where the parameter $\alpha$ is greater than $\lceil\alpha_Y\rceil\vee\lceil\alpha_X\rceil$. For any $j\in \{0\}\cup [J]$ with $J=\lceil {\frac{1}{2\alpha_Y+D_Y+d_X\frac{\alpha_Y}{\alpha_X}}}\cdot \log_2 (\frac{n}{\log n})\rceil$, define $\Psi_j^{D_Y}$ as the subset of the wavelet basis   $\bigcup_{j\geq 0}\ov\Psi_j^{D_Y}$ for which
\begin{equation*}
      \Psi_j^{D_Y}=\{\psi\in \ov \Psi_j^{D_Y}:\, \text{supp}(\psi)\cap \mb B_{\mb R^{D_Y}}(\mathbf{0},L)\neq \emptyset\}.
\end{equation*}
 Consider a smooth transition function $\rho:\mb R\to [0,1]$ defined by
 \begin{equation}\label{def:transition}
      \rho(t)=\left\{
      \begin{array}{cc}
        0   & |t|\geq 2 \\
        1   & |t|\leq 1\\
       \frac{1}{1+ \exp(\frac{3-2t}{(t-1)(t-2)})}& 1<t<2\\
        \frac{1}{1+ \exp(\frac{2t+3}{(t+1)(2+t)})}& -2<t<-1.\\
      \end{array}
      \right.
  \end{equation}
This function ensures $\rho(t) = 1$ for $t \in [0,1]$ and $\rho(t) = 0$ for $t \in [2,\infty)$. For any $j\in [J]$, define a class of functions $\ms S_j$  on $ \mb R^{D_X}$  as
\begin{equation}\label{msSjR1}
\begin{aligned}
     \ms S_j&=\Bigg\{S(x)=\frac{\sum_{i=1}^{W_j} \sum_{k\in \mb N_0^{D_X}, |k|<\alpha_X}a_{ik}  (x-b_i)^k \rho(\frac{\|x-b_i\|}{\varepsilon^x_j})}{\sum_{i=1}^{W_j} \rho(\frac{\|x-b_i\|}{\varepsilon^x_j})+\frac{1}{n}}:\,
     b_i\in \mb B_{\mb R^{D_X}}(\mathbf{0},L),  a_{ik} \in [-\frac{C}{2^{D_Yj/2}},\frac{C}{2^{D_Yj/2}}]\\
     &\qquad  \text{for any } i\in [W_j] \text{ and } k\in \mb N_0^{D_X} \text{ with } |k|<\alpha_X \Bigg\},
\end{aligned}
\end{equation}
where $\varepsilon^x_j=2^{\frac{jD_Y}{2\alpha_X+d_X}}(\frac{n}{\log n})^{-\frac{1}{2\alpha_X+d_X}}$, $W_j=C_1\,(\varepsilon^x_j)^{-d_X} $ and $C,C_1$ are large enough constants. Consider the  estimator
\begin{equation*}
    \wh u_{\psi}(\cdot)=\underset{S\in  \ms S_j}{\arg\min} \frac{1}{n}\sum_{i=1}^n  (\psi(Y_i)-S(X_i))^2,\quad j\in \{0\}\cup [J],\psi\in \Psi_j^{D_Y}.
\end{equation*}
Finally, define a conditional density estimator for $\mu^*_{Y|x}$ as
\begin{equation*}
  \wh u(\cdot|x)=\sum_{j=0}^J \sum_{\psi\in \Psi_j^{D_Y}} \psi(\cdot) \wh u_{\psi}(x).
 \end{equation*}

\subsection{Minimax Optimal Estimator for Regime 2}\label{R2gamma}

\subsubsection{Density regression in the ambient space}\label{R2gamma>1}

For any $j\in \mb N$, define
\begin{equation*}
    \Psi_j^{D_Y}=\{\psi\in\ov\Psi _j^{D_Y}:\,{\rm supp}(\psi)\cap \mb B_{\mb R^{D_Y}}(\mathbf{0},L)\neq \emptyset\},
\end{equation*}
where  $\bigcup_{j\geq 0}\ov\Psi_j^{D_Y}$ is a wavelet basis that satisfies the properties stated in Lemma~\ref{le:wavelet} with the parameter $\alpha$ being greater than $\lceil\alpha_Y\rceil\vee\lceil\alpha_X\rceil\vee \lceil \frac{d_Y\alpha_Y}{2\alpha_X+d_X}\rceil\vee \lceil\beta_Y\rceil$. For any $j\in \mb N$, consider the estimator
 \begin{equation}\label{eqn:whSjapp}
 \begin{aligned}
         &\wh S_j^{\dagger}(\,\cdot\,,\,\cdot\,)={\arg\min}_{S\in  {\m S}_j^\dagger} \frac{1}{n}\sum_{i=1}^n \sum_{\psi\in \Psi_j^{D_Y}} (2^{\frac{j(d_Y-D_Y)}{2}}\psi(Y_i)-S(\psi,X_i))^2.
 \end{aligned}
\end{equation}
To construct the families $\m S_j^\dagger$, we leverage the fact that, for any $\psi \in \Psi_j^{D_Y}$,  the term 
\begin{align*}
    \mb{E}_{\mu^*_{Y|x}}\big[\,2^{jd_Y-\frac{jD_Y}{2}}\psi(y)\big]=\int_{\m M_Y} 2^{jd_Y-\frac{jD_Y}{2}}\,\psi(y)\,u^*(y\,|\,x)\, {\rm vol}_{\m M_Y}(\dd y),
\end{align*}
where $\mathrm{vol}_{\m M_Y}(\dd y)$ denotes the volume measure on the manifold $\m M_Y$, is $\m H^{\alpha_X}$-smooth as a function of $x$ and has a bounded H\"{o}lder norm. As a result, each conditional expectation $\mb{E}_{\mu^*_{Y|x}}[\psi(y)]$ can be effectively approximated using local polynomial approximation techniques.  Furthermore, since the response space $\m M_Y$ lies on a low-dimensional submanifold, only $\mathcal{O}(2^{d_Y j})$ of the functions $\psi(\cdot)$ will have non-zero conditional means. This observation allows us to construct parametric families $\m S_j^\dagger$ whose complexity depends only on the level $j$, the intrinsic dimensions $d_x,d_Y$ and the smoothness level $\alpha_X$.
According to Lemma~\ref{le:wavelet}, for any $j\in \mb N$, we can express $\Psi_j^{D_Y}$ using an index set as follows:
\begin{equation*}
     \Psi _j^{D_Y}=\big\{\psi_{j\iota}(\cdot):\, \iota\in \ms I_j\subset [0,1]^{D_Y+1}\big\},
\end{equation*}
where $\ms I_j$ is a $c\,2^{-j}$-separated set for some constant $c > 0$. We denote the index of $\psi \in \Psi_j^{D_Y}$ by $\m I_j(\psi)$; that is, for $\psi = \psi_{j\iota}$, we write $\m I_j(\psi) = \iota$. Then we define $\m S^{\dagger}_j$ as
\begin{equation}\label{defSdaggerr2}
\begin{aligned}
     \m S_j^\dagger&=\Bigg\{S(\psi,x)=\frac{\sum_{i_1=1}^{W_j}\sum_{i_2=1}^{W_j'} \sum_{k\in \mb N_0^{D_X}, |k|<\alpha_X}a_{i_1i_2k}(x-b_{i_2})^k \rho\Big(\frac{\|x-b_{i_2}\|}{\varepsilon^x_j}\Big)\rho\Big(\frac{\|\m I_j(\psi)-e_{i_1}\|}{\varepsilon_j^y}\Big)}{\sum_{i_1=1}^{W_j}\sum_{i_2=1}^{W_j'} \rho\Big(\frac{\|x-b_{i_2}\|}{\varepsilon^x_j}\Big)\rho\Big(\frac{\|\m I_j(\psi)-e_{i_1}\|}{\varepsilon_j^y}\Big)+\frac{1}{n}}:\, \\
     & \qquad\qquad
     b_{i_2}\in \mb B_{\mb R^{D_X}}(\mathbf{0},L),\,  a_{i_1i_2k}\in [-\frac{C}{2^{d_Yj/2}},\frac{C}{2^{d_Yj/2}}], \,e_{i_1}\in [0,1]^{D_Y+1} \text{ for any } i_1,i_2,k \Bigg\},
\end{aligned}
\end{equation}
where $\rho$ is a smooth transition function defined in~\eqref{def:transition}; $\varepsilon_j^y = \frac{2^{-j}}{C_1}$ and $\varepsilon_j^x = 2^{{j d_Y}/{(2\alpha_X + d_X)}} \big(\frac{n}{\log n}\big)^{-1/{(2\alpha_X + d_X)}}$ are the bandwidth parameters in the $y$ and $x$ directions, respectively. The quantities $W_j = C_3\, (\varepsilon_j^y)^{-d_Y}$ and $W_j' = C_2\, (\varepsilon_j^x)^{-d_X}$ represent the number of local neighborhoods in $y$ and $x$, respectively, over which the partition of unity is defined. The numbers $C$, $C_1$, $C_2$, and $C_3$ are sufficiently large constants. Now we define
\begin{equation*}
    \wh u(y\,|\,x)=\sum_{j=0}^J\sum_{\psi\in \Psi_j^{D_Y}} 2^{\frac{j(D_Y-d_Y)}{2}} \wh S_j^{\dagger}(\psi,x)\psi(y).
\end{equation*}
The measure $\wh u(y\,|\,x)\,\dd y$, when utilized directly as an estimator for the conditional distribution $\mu^*_{Y|x}$, can achieve minimax optimality under the condition when $\gamma>\frac{d_Y\alpha_X}{2\alpha_X+d_X}$. This is formally established in the theorem presented below.



\begin{theorem}\label{th:thgamma>1}
 Let $J=\lceil {\frac{1}{d_Y}}\cdot \log_2 (\frac{n}{\log n})\rceil$. With the choice of $\m S_j^{\dagger}$ defined in~\eqref{defSdaggerr2},  consider any  distribution $\mu^*=\mu^*_X\mu^*_{Y|X}\in \m P_2$,  it holds with probability at least $1-\frac{1}{n}$ that, for any $\gamma>\frac{d_Y\alpha_X}{2\alpha_X+d_X}$,  
    \begin{equation*}
            \mb{E}_{\mu^*_X} \Big[ \underset{f\in \m H_1^{\gamma}(\mb R^{D_Y})}{\sup} \big|\mb{E}_{\mu^*_{Y|x}} f(y)- \int_{\mb R^{D_Y}} f(y) \wh u(y\,|\,x)\,\dd y\big|\Big]\lesssim (\log n) \cdot  n^{-\frac{\alpha_X}{2\alpha_X+d_X}}.
    \end{equation*}
\end{theorem}
\noindent The proof of Theorem~\ref{th:thgamma>1} is given in Appendix~\ref{proofth:thgamma>1}.

\subsubsection{Density regression in the latent space }\label{R2gamma<1}
 We split the data into two subsets by considering $I_1=[\lfloor n/2\rfloor]$ and $I_2=[n]\setminus I_1$.
  Let $\{\omega_k=(x_k,y_k)\}_{k=1}^{K}$ be a $\tau_2$-covering set of $\mb B_{\mb R^{D_X}}(\mathbf{0},L)\times \mb B_{\mb R^{D_Y}}(\mathbf{0},L)$, where $\tau_2$ is a sufficiently small absolute constant. Define
 \begin{equation}\label{setmK}
     \wh{\m K}=\{k\in [K]:\, \exists i\in I_1, \|(x_i,y_i)-\omega_k\|\leq \sqrt{2}\tau_2\}.
 \end{equation}
Consider a wavelet basis $\bigcup_{j\geq 0}\ov\Psi _j^{d_Y}$  that satisfies the properties of Lemma~\ref{le:wavelet}, where the parameter $\alpha$ is greater than $\lceil\alpha_Y\rceil\vee\lceil\alpha_X\rceil\vee\lceil\beta_Y\rceil\vee \lceil\beta_X\rceil$. Then for any $j\in \mb N$, we denote
 \begin{equation*}
      \Psi_j^{d_Y}=\{\psi\in \ov \Psi_j^{d_Y}:\, \text{supp}(\psi)\cap \mb B_{\mb R^{d_Y}}(\mathbf{0},2\tau_2)\neq \emptyset\}.
\end{equation*}
For any $k\in \wh{\m K}$, we consider the estimator
\begin{equation}\label{est:GV}
    (\wh G_{[k]}, \wh V_{[k]})=\underset{G\in  {\m G}\atop V\in \mb O(D_Y,d_Y)}{\arg\min} \frac{1}{|I_1|}\sum_{i\in I_1} \|Y_i-G(V^T(Y_i-y_k),X_i)\|^2\mathbf{1}\big(X_i\in \mb B_{\mb R^{D_X}}(x_k,2\tau_2)\big)\mathbf{1}\big(Y_i\in \mb B_{\mb R^{D_Y}}(y_k,2\tau_2)\big),
\end{equation}
where $\mb O(D_Y,d_Y)=\{A\in \mb R^{D_Y\times d_Y}:\, A^TA=I_{d_Y}\}$.  To clarify the selection of $\m G$, we note that the choice depends on whether the submanifold $\m M_{Y|x}$ varies with $x$. Here in Regime 2, since $\m M_{Y|x}$ remains invariant across $x$, we define $\m G$ as a function class that operates solely on the latent space $\mb R^{d_Y}$ and does not depend on the covariate. Given that the global manifold $\m M_Y$ is $\beta$-smooth, we construct each function $G:\,\mb R^{d_Y}\times \mb R^{D_X} \to\mb R^{D_Y}$ in $\m G$ by truncating the wavelet expansion of $\m H^{\beta}$-smooth functions at a finite resolution level. Specifically, the function class $\m G$ is defined as
\begin{equation}\label{defGr2}
\begin{aligned}
       {\m G}&=\bigg\{G(z,x)=G(z)=\sum_{j_1=0}^{J_1} \sum_{\psi_1\in  \Psi_{j_1}^{d_Y} } g_{\psi_1 }\,\psi_1(z) \,:  g_{\psi_1}\in [-L_1\, \delta_{j_1},L_1\, \delta_{j_1}]^{D_Y} \text{ for each }\psi_1\bigg\},
\end{aligned}
\end{equation}
where $J_1= \lceil\log_2 (n^{-1/{d_Y}})\rceil$, $\delta_{j_1}=2^{-{d_Yj_1}/{2}-(j_1\beta_Y)}$ and $L_1$ is a sufficiently large constant. Then we denote $\wh Q_{[k]}(\cdot)=\wh V_{[k]}^T(\cdot-y_k)$.
For any $k\in \wh {\m K}$, $j\in \{0\}\cup [J]$ and $\psi\in \Psi_j^{d_Y}$,  we consider the estimator
\begin{equation}\label{eqn:whS1}
\begin{aligned}
    & \wh v_{k\psi}(\cdot)={\arg\min}_{S\in  \ms S_j} \frac{1}{|I_2|}\sum_{i\in I_2}\sum_{\psi\in \Psi_j^{d_Y}} (\psi(\wh Q_{[k]}(Y_i)) \rho_\sk(X_i,Y_i)-S(X_i))^2,
\end{aligned}
   \end{equation}
where $\rho_\sk(x,y)=\frac{\rho( {\|(x,y)-(x_k,y_k)\|^2}/{\tau_2^2})}{\sum_{k=1}^K \rho( {\|(x,y)-(x_k,y_k)\|^2}/{\tau_2^2})}$ with $\rho$ being defined in~\eqref{def:transition}. 
Note that the construction of $\ms S_j$ in both Regime 2 and the later Regime 3b is based on the construction in Equation~\eqref{msSjR1}, originally developed for density regression in Euclidean space (Regime 1). The key modification is the substitution of the ambient dimension $D_Y$ with the intrinsic dimension $d_Y$. Specifically, 
for any  $j\in \{0\}\cup [J]$ with $J=\lceil {\frac{1}{2\alpha_Y+d_Y+d_X\frac{\alpha_Y}{\alpha_X}}}\cdot \log_2 (\frac{n}{\log n})\rceil$, define 
 \begin{equation}\label{defMSddaggerr2}
\begin{aligned}
     \ms S_j&=\Bigg\{S(x)=\frac{\sum_{i=1}^{W_j} \sum_{k\in \mb N_0^{D_X}, |k|<\alpha_X}a_{ik} (x-b_i)^k \rho(\frac{\|x-b_i\|}{\varepsilon^x_j})}{\sum_{i=1}^{W_j} \rho(\frac{\|x-b_i\|}{\varepsilon^x_j})+\frac{1}{n}}:\,\\
     &\qquad\qquad\qquad b_i\in \mb B_{\mb R^{D_X}}(\mathbf{0},L),  a_{ik}\in [-\frac{C}{2^{d_Yj/2}},\frac{C}{2^{d_Yj/2}}], \text{for any } i,k \Bigg\},
\end{aligned}
\end{equation}
where $W_j=C_1\, (\varepsilon_j^x)^{-d_X}$, $\varepsilon_j^x= 2^{\frac{jd_Y}{2\alpha_X+d_X}}(\frac{n}{\log n})^{-\frac{1}{2\alpha_X+d_X}}$, $C_1,C$ are large enough constants  and $\rho$ is a smooth transition function defined in~\eqref{def:transition}. 
Then denote $\wh \nu_\sk(\cdot|x)$ as the measure  that has a density function $\sum_{j=0}^J\sum_{\psi\in \Psi_j^{d_Y}}\psi(\cdot)\wh v_{k\psi}(x)$ with respect to the Lebesgue measure on $\mb R^{d_Y}$. We can define  a mixture of conditional generative models $ \sum_{k\in \wh{\m K}} \wh G_{[k]}(\cdot,x)_{\#}\wh\nu_\sk(\cdot|x)$, which, as an estimator of the conditional distribution, can achieve minimax optimality when $\gamma\leq 1$, as detailed in the following theorem.

\begin{theorem}\label{thgamma<1}
 Let $J=\lceil {\frac{1}{2\alpha_Y+d_Y+d_X\frac{\alpha_Y}{\alpha_X}}}\cdot \log_2 (\frac{n}{\log n})\rceil$. With the choices of $\m G$ and $\ms S_j$  defined in~\eqref{defGr2} and ~\eqref{defMSddaggerr2} respectively. Consider any distribution $\mu^*=\mu^*_X\mu^*_{Y|X}\in \m P_2^*$, it holds with probability at least $1-\frac{1}{n}$ that for any $\gamma\leq 1$,
    \begin{equation*}
    \begin{aligned}
           &\mb{E}_{\mu^*_X}\Big[\underset{f\in \m H^{\gamma}_1(\mb R^{D_Y})}{\sup} \big|\mb{E}_{\mu^*_{Y|x}}[f(Y)]-\sum_{k\in \wh{\m K}} \int_{\mb R^{d_Y}} f(\wh G_\sk(z,x)) \sum_{j=0}^J\sum_{\psi\in \Psi_j^{d_Y}}\psi(z)\wh v_{k\psi}(x)\,\dd z\big|\Big]\\
           &\lesssim  (\log n)^2\cdot n^{-\frac{\alpha_X}{2\alpha_X+d_X}}+(\frac{n}{\log n})^{-\frac{\alpha_Y+\gamma}{2\alpha_Y+d_Y+\frac{\alpha_Y}{\alpha_X}d_X}}+n^{-\frac{\gamma}{\frac{d_Y}{\beta_Y}}}.
    \end{aligned}
    \end{equation*}
\end{theorem}
 \noindent The proof of Theorem~\ref{thgamma<1} is given in Appendix~\ref{proofthgamma<1}.
 \subsubsection{Simultaneous minimax optimal estimator for $\gamma>0$}\label{R2gamma>0}

 Choose $J=\lceil {\frac{1}{2\alpha_Y+d_Y+d_X\frac{\alpha_Y}{\alpha_X}}}\cdot \log_2 (\frac{n}{\log n})\rceil$, then we define an operator $\m J(f,x)$ such that for any continuous function $f:\mb R^{D_Y}\to \mb R$ and $x\in \mb R^{D_X}$,
\begin{equation*}
\begin{aligned}
       & \wh{\m J}(f,x)=\sum_{j=0}^J\sum_{\psi\in  \Psi_j^{D_Y}} f_{\psi}2^{\frac{j(D_Y-d_Y)}{2}} \wh S_j^\dagger(\psi,x)+\sum_{k\in \wh{\m K}} \int_{\mb R^{d_Y}} f^{\perp}_J(\wh G_\sk (z)) \sum_{j=0}^J \sum_{\psi\in \Psi_j^{d_Y}}\psi(z)  \wh v_{k\psi}(x)\,\dd z,\\
       &\qquad f_{\psi}=\int_{\mb R^{D_Y}} f(y)\psi(y)\,\dd y,\quad f^{\perp}_J(y)=f(y)-\sum_{j=0}^J \sum_{\psi\in  \ov\Psi_j^{D_Y}} f_{\psi} \psi(y),
\end{aligned}
\end{equation*}
where all notations are adopted from subsections~\ref{R2gamma>1} and~\ref{R2gamma<1}. The estimator  $\wh{\m J}(f,x)$ can achieve the upper bound specified in Theorem~\ref{th:combined} (Regime 2). By utilizing $\wh{\m J}(f,x)$, we can also derive a valid conditional distribution estimator that is simultaneous optimal for all $\gamma>0$ using the steps described below.

\medskip
\noindent Consider the set $\Gamma=\{\frac{1}{\log n}, \frac{2}{\log n}, \cdots, \frac{s}{\log n}\}$ with $s=\lceil \frac{d_Y\alpha_X}{2\alpha_X+d_X}\log n\rceil$,  and define 
\begin{equation*}
\delta_{n,\gamma}= C_{\gamma}\Big((\log n)^3\cdot n^{-\frac{\alpha_X}{2\alpha_X+d_X}}+(\frac{n}{\log n})^{-\frac{\alpha_Y+\gamma}{2\alpha_Y+d_Y+\frac{\alpha_Y}{\alpha_X}d_X}}+  n^{-\frac{\gamma}{\frac{d_Y}{\beta_Y}}}\Big). 
\end{equation*}
For any $x\in \m M_X$,   consider the estimator
\begin{equation*}
    \wh \mu_{Y|x}=\underset{\mu\in \m P^*_{Y}}{\arg\min}  \sum_{\gamma\in \Gamma} \frac{1}{\delta_{n,\gamma}}\cdot \underset{f\in \m H^{\gamma}_1(\mb R^{D_Y})}{\sup}\Big[\mb{E}_{\mu}[f(y)]- \wh{\m J}(f,x)\Big],
    \end{equation*}
where $\m P_Y^*$ includes all probability measures of $\mu$ that are supported on a submanifold $\m M_{Y}$ and have a density function $u(\cdot)$ with respect to the volume measure of $\m M_{Y}$ such  that $\m M_{Y}\in \ms M^{\beta_Y}_{\tau,\tau_1,L}(d_Y,D_Y)$, $\mu\in \m H^{\alpha_Y}_L(\m M_Y)$.

\begin{corollary}\label{co1}
   With the choices of $\m S_j^\dagger, \ms S_j,\m G$  defined in~\eqref{defSdaggerr2},~\eqref{defMSddaggerr2} and~\eqref{defGr2} respectively, alongside $J=\lceil {\frac{1}{2\alpha_Y+d_Y+d_X\frac{\alpha_Y}{\alpha_X}}}\cdot \log_2 (\frac{n}{\log n})\rceil$. For any $\mu^*=\mu^*_X\mu^*_{Y|X}\in \m P^*_2$, it holds with probability at least $1-\frac{1}{n}$ that for any $\gamma> 0$,  
 \begin{equation*}
     \mb{E}_{\mu^*_X}[d_{\gamma}(\mu^*_{Y|X},\wh \mu_{Y|X})]\lesssim  (\log n)^4\cdot n^{-\frac{\alpha_X}{2\alpha_X+d_X}}+\log n\cdot (\frac{n}{\log n})^{-\frac{\alpha_Y+\gamma}{2\alpha_Y+d_Y+\frac{\alpha_Y}{\alpha_X}d_X}}+ \log n\cdot  n^{-\frac{\gamma}{\frac{d_Y}{\beta_Y}}}.
 \end{equation*}
\end{corollary}
\noindent The proof of Corollary~\ref{co1} is given in Appendix~\ref{proofco1}.

\subsection{Minimax Optimal Estimator for Regime 3b}\label{R3gamma>0}
The estimator is formulated similarly to that for Regime 2.  Choose $J=\lceil {\frac{1}{2\alpha_Y+d_Y+d_X\frac{\alpha_Y}{\alpha_X}}}\cdot \log_2 (\frac{n}{\log n})\rceil$,  we define an operator $\wh{\m J}(f,x)$ so that for any continuous function $f:\mb R^{D_Y}\to \mb R$ and $x\in \mb B_{\mb R^{D_X}}(\mathbf{0},L)$,
\begin{equation*}
\begin{aligned}
       & \wh{\m J}(f,x)= {\sum_{j=0}^J\sum_{\psi\in  \Psi_j^{D_Y}} f_{\psi}2^{\frac{j(D_Y-d_Y)}{2}} \wh S_j^\dagger(\psi,x)}+ {\sum_{k\in \wh{\m K}} \int_{\mb R^{d_Y}} f^{\perp}_J(\wh G_\sk (z,x)) \sum_{j=0}^J \sum_{\psi\in \Psi_j^{d_Y}}\psi(z) \wh v_{k\psi}(x)\,\dd z},\\
  &\qquad  \Psi_j^{D_Y}=\{\psi\in\ov\Psi _j^{D_Y}:\,{\rm supp}(\psi)\cap \mb B_{\mb R^{D_Y}}(\mathbf{0},L)\neq \emptyset\},\\
       &\qquad f_{\psi}=\int_{\mb R^{D_Y}} f(y)\psi(y)\,\dd y,\quad f^{\perp}_J(y)=f(y)-\sum_{j=0}^J \sum_{\psi\in \ov \Psi_j^{D_Y}} f_{\psi} \psi(y), 
       \end{aligned}
\end{equation*}
where
\begin{equation*} 
    (\wh G_{[k]}, \wh V_{[k]})=\underset{G\in  {\m G}\atop V\in \mb O(D_Y,d_Y)}{\arg\min} \frac{1}{|I_1|}\sum_{i\in I_1} \|Y_i-G(V^T(Y_i-y_k),X_i)\|^2\mathbf{1}\big(X_i\in \mb B_{\mb R^{D_X}}(x_k,2\tau_2)\big)\mathbf{1}\big(Y_i\in \mb B_{\mb R^{D_Y}}(y_k,2\tau_2)\big),
\end{equation*}
and $\wh S_j^\dagger(\,\cdot\,,\,\cdot\,)$,  $\wh v_{k\psi}(\cdot)$  are the estimators defined in~\eqref{eqn:whSjapp} and~\eqref{eqn:whS1} respectively. For the approximation families, $\ms S_j$ is defined as in~\eqref{defMSddaggerr2}. For the family $\m G$, 
unlike Regime 2,  Regime 3b involves  scenarios where the submanifold $\mathcal{M}_{Y|x}$ varies with $x$. In this context, we construct $G$ using the tensor-product expansion of ${\m H}^{\beta_Y,\beta_X}$-smooth functions. Specifically, we use the basis functions $\big\{\psi_1(z)\cdot\psi_2(x) : \,\psi_1 \in \bigcup_{j=0}^{\infty} \Psi_j^{d_Y},\ \psi_2 \in \bigcup_{j=0}^{\infty} \ov\Psi_j^{D_X}\big\}$ and truncate the expansion at finite resolution levels. Accordingly, the function class $\m G$ is defined as:
\begin{equation}\label{defGr3}
\begin{aligned}
      {\m G}&=\{G(z,x)=\sum_{j_1=0}^{J_1}\sum_{j_2=0}^{J_2}\sum_{\psi_1\in  \Psi_{j_1}^{d_Y} }\sum_{\psi_2\in \ov {\Psi}_{j_2}^{D_X}}g_{\psi_1\psi_2}\psi_1(z)\psi_2(x)\,: \\
   &\qquad g_{\psi_1\psi_2}\in [-L_1\, \delta_{j_1j_2},L_1\, \delta_{j_1j_2}]^{D_Y}, \text{ for  each }\psi_1, \psi_2\},
\end{aligned}
\end{equation}
where $J_1= \lceil\log_2 (n^{-\frac{1}{d_Y+d_X\frac{\beta_Y}{\beta_X}}})\rceil$, $J_2= \lceil\log_2 (n^{-\frac{1}{d_X+d_Y\frac{\beta_X}{\beta_Y}}})\rceil$, $\delta_{j_1j_2}=2^{-\frac{d_Yj_1+D_Xj_2}{2}-((j_1\beta_Y)\vee (j_2\beta_X))}$ and $\Psi_j^{d_Y}=\{\psi\in \ov \Psi_j^{d_Y}:\, \text{supp}(\psi)\cap \mb B_{\mb R^{d_Y}}(\mathbf{0},2\tau_2)\neq \emptyset\}$.

Now, let's define the class ${\m S}_j^\dagger$.  Compared with Regime 2,
the construction of $\m S_j^{\dagger}$ becomes more challenging in Regime 3b , where the conditional response space $\m M_{Y|x}$ varies with $x$. In this setting, the conditional distribution $\mu^*_{Y|x}$ can be expressed as a mixture of conditional generative models, $\mu^*_{Y|x}=\sum_{k=1}^{K^*} G_\sk^*(\,\cdot\,,x)_{\#}\nu^*_\sk(\,\cdot\,|x)$, where the generators $G_k^*$ are $x$-dependent and ${\m H}^{\beta_Y, \beta_X}$-smooth (see Lemma~\ref{legenerative} in Appendix~\ref{app: regularitymanifold}). The conditional mean
\begin{equation}\label{cond:gen}
    \mb{E}_{\mu^*_{Y|x}}\big[\,2^{jd_Y-\frac{jD_Y}{2}}\,\psi(y)\,\big]=\sum_{k=1}^{K^*} \int_{\mb R^{d_Y}} 2^{jd_Y-\frac{jD_Y}{2}}\,\psi\big(G_\sk^*(z,x)\big)\,\nu^*_\sk(z\,|\,x)\,\dd z,
\end{equation}
 may not be uniformly $\m H^{\alpha_X}$-smooth in $x$ because the gradients of $\psi \in \Psi_j^{D_Y}$ grow rapidly with $j$. To address this challenge, we propose a hybrid strategy for constructing approximation families for $\mb{E}_{\mu^*_{Y|x}}[\psi(y)]$, applied over all $\psi \in \Psi_j^{D_Y}$ and $x \in \m M_X$. The first component involves building parametric approximation families for the generators ${G^*_{\sk}}$ and the latent distributions ${\nu_\sk^*}$, enabling direct approximation of the integral in~\eqref{cond:gen}. While effective for high-resolution levels (large $j$), this approach does not fully exploit the smoothness properties of $\psi$ when $j$ is small.
 The second component adopts a strategy similar to that used in Regime 2: for each $\psi \in \Psi_j^{D_Y}$, we treat $\mb{E}_{\mu^*_{Y|x}}[\psi(y)]$ as an $\m H^{\alpha_X}$-smooth function of $x$, and use local polynomial approximations. However, as noted earlier, this strategy becomes less effective at large $j$ due to the growing instability of the wavelet basis functions.
 
To combine these two strategies in a resolution-adaptive manner, we first define parametric function classes for approximating $G_\sk^*$ and $\nu_\sk^*$, where the number of parameters increases with $j$, allowing the approximation accuracy to improve as resolution increases. We then use local polynomial approximations to model the residual, capturing smooth variation in $x$.
Specifically, letting $\m T_a(x) = \max\big(-a, \min(a, x)\big)$ 
be a truncation operator and $\wt \beta_X = \alpha_X + \frac{\alpha_X}{\alpha_Y}$,  and
recall that for any $j\in \mb N$, $\Psi _j^{D_Y}$ can be written as an index set
\begin{equation*}
     \Psi _j^{D_Y}=\{\psi_{j\iota}(\cdot):\, \iota\in \ms I_j\subset [0,1]^{D_Y+1}\},
\end{equation*}
where $\ms I_j$ is $c 2^{-j}$-separated and we use $\m I_j(\psi)$ to denote the index of $\psi\in\Psi _j^{D_Y}$ (i.e., $\m I_j(\psi_{j\iota})=\iota$). 
We define $\m S_j^{\dagger}$ for Regime 3b as the class of mappings $S:\, \Psi_j^{D_Y} \times \mb R^{D_X} \to \mb R$ structured as follows:
\begin{equation}\label{defSdaggerr3} 
    \begin{aligned}
    &S(\psi,x)=\m T_{C_12^{-\frac{d_Yj}{2}}}\bigg(\ \sum_{i_1=1}^{W_j}\sum_{i_2=1}^{W_j'} \bigg[\, \rho\Big(\frac{\|x-b_{i_2}\|}{\varepsilon_j^x}\Big)\, \rho\Big(\frac{|\m I_j(\psi)-e_{i_1i_2}|}{\varepsilon_j^y}\Big)\\
 &\quad\cdot\Big\{\sum_{k=1}^{K^*} \int_{\mb B_{\mb R^{d_Y}}(\mathbf{0},\tau_1)} 2^{\frac{j(d_Y-D_Y)}{2}}\,\psi\big(G_{k,i_2}(z,x)\big)\, \nu_{k,i_2}(z,x)\,\dd z \ \ +\sum_{|l|\leq\lfloor\wt\beta_X\rfloor^2 +\atop \lfloor \alpha_X\rfloor,\, l\in\mb N_{0}^{D_X}} a_{i_1i_2l}\,  (x-b_{i_2})^l \Big)\bigg]\\
 &\qquad\qquad\qquad\qquad\qquad\qquad\cdot \frac{1}{\sum_{i_1=1}^{W_j}\sum_{i_2=1}^{W_j'}\rho\Big(\frac{\|x-b_{i_2}\|}{\varepsilon_j^x}\Big)\,\rho\Big(\frac{|\m I_j(\psi)-e_{i_1i_2}|}{\varepsilon_j^y}\Big)+\frac{1}{n^2}}\ \bigg), \\
  &\text{where} \quad  G_{k,i_2}(z,x)=\sum_{s=0}^{j} \sum_{\psi\in \wt\Psi_s^{d_Y}}\sum_{l\in \mb N_0^{D_X}\atop |l|<\beta_X}  g_{k,i_2,s,\psi,l} \,(x-b_{i_2})^{l} \, \psi(z), \\
  &\qquad\qquad\qquad \text{ and}\quad \nu_{k,i_2}(z,x)=\sum_{s=0}^{j} \sum_{\psi\in \wt\Psi_s^{d_Y}}\sum_{l\in \mb N_0^{D_X}\atop |l|<\alpha_X}  v_{k,i_2,s,\psi,l}\, (x-b_{i_2})^{l} \, \psi(z). 
   \end{aligned}
\end{equation}
 Here,  $\rho$ is a smooth transition function defined in~\eqref{def:transition}; $\wt\Psi_j^{d_Y}=\{\psi\in \ov \Psi_j^{d_Y}:\, \text{supp}(\psi)\cap \mb B_{\mb R^{d_Y}}(\mathbf{0},\tau_1)\neq \emptyset\}$; $K^*$ and $C_1$ are sufficiently large constants.  $\varepsilon_j^y = \frac{2^{-j}}{C_1}$ and $\varepsilon_j^x = 2^{{j d_Y}/{(2\alpha_X + d_X)}} \big(\frac{n}{\log n}\big)^{-1/{(2\alpha_X + d_X)}}$ are the bandwidth parameters in the $y$ and $x$ directions, respectively. The quantities $W_j = C_3\, (\varepsilon_j^y)^{-d_Y}$ and $W_j' = C_2\, (\varepsilon_j^x)^{-d_X}$ for large enough constants $C_2$ and $C_3$.  The parameters are constrained as follows: $g_{k,i_2,s,\psi,l} \in [-C_1, C_1]^{D_Y}$, $v_{k,i_2,s,\psi,l} \in [-C_1, C_1]$, and $a_{i_1i_2l} \in [-C_1 n, C_1 n]$. The indices $e_{i_1i_2}$ lie in $[0,2]^{D_Y + 1}$. The centers $\{b_1, b_2, \dots, b_{W_j'}\}$ are $\varepsilon_j^x$-separated, meaning that $\|b_i - b_k\| \geq \varepsilon_j^x$ for any $i \neq k$ in $[W_j']$, and all lie within the ball $\mb B_{\mb R^{D_X}}(\mathbf{0}, L_1)$ for a large enough constant $L_1$.


 \medskip
\noindent Similar to Regime 2, by utilizing $\wh{\m J}(f,x)$, we can also develop a conditional distribution estimator by considering the set $\Gamma=\{\frac{1}{\log n}, \frac{2}{\log n}, \cdots, \frac{s}{\log n}\}$ with $s=\lceil (d_Y\vee (\frac{d_Y}{\beta_Y}+\frac{d_X}{\beta_X}))\frac{\alpha_X}{2\alpha_X+d_X}\log n\rceil$,  and define 
\begin{equation*}
\delta_{n,\gamma}=C_{\gamma}\cdot\Big((\log n)^3\cdot n^{-\frac{\alpha_X}{2\alpha_X+d_X}}+(\log n)\cdot(\frac{n}{\log n})^{-\frac{\alpha_Y+\gamma}{2\alpha_Y+d_Y+\frac{\alpha_Y}{\alpha_X}d_X}}+(\log n) \cdot n^{-\frac{\gamma}{\frac{d_Y}{\beta_Y}+\frac{d_X}{\beta_X}}}\Big). 
\end{equation*}
For any $x\in \m M_X$,   consider the estimator
\begin{equation*}
    \wh \mu_{Y|x}=\underset{\mu\in \m P^*_{Y}}{\arg\min}  \sum_{\gamma\in \Gamma} \frac{1}{\delta_{n,\gamma}} \cdot\underset{f\in \m H^{\gamma}_1(\mb R^{D_Y})}{\sup} \Big[\mb{E}_{\mu}[f(y)]- \wh{\m J}(f,x)\Big],
    \end{equation*}
 where $\m P_Y^*$ is defined as in Appendix~\ref{R2gamma<1}.
 \begin{corollary}\label{co2}
     With the choice of $\m S_j^\dagger, \ms S_j,\m G$  defined in~\eqref{defSdaggerr3},~\eqref{defMSddaggerr2} and~\eqref{defGr3} respectively, alongside $J=\lceil {\frac{1}{2\alpha_Y+d_Y+d_X\frac{\alpha_Y}{\alpha_X}}}\cdot \log_2 (\frac{n}{\log n})\rceil$.  For any $\mu^*=\mu^*_X\mu^*_{Y|X}\in \m P^*_3$, it holds with probability at least $1-\frac{1}{n}$ that for any $\gamma> 0$,  
 \begin{equation*}
     \mb{E}_{\mu^*_X}[d_{\gamma}(\mu^*_{Y|X},\wh \mu_{Y|X})]\lesssim (\log n)^4\cdot n^{-\frac{\alpha_X}{2\alpha_X+d_X}}+(\log n)^2\cdot(\frac{n}{\log n})^{-\frac{\alpha_Y+\gamma}{2\alpha_Y+d_Y+\frac{\alpha_Y}{\alpha_X}d_X}}+(\log n)^{2}\cdot n^{-\frac{\gamma}{\frac{d_Y}{\beta_Y}+\frac{d_X}{\beta_X}}}.
 \end{equation*}
 \end{corollary}

\section{Proof for Distribution Regression with Euclidean Response}\label{APP:proofEuclidean}
 \subsection{Proof of Theorem~\ref{th1} (minimax upper bound for Regime 1)}\label{proof:th1}

For any $j\in \{0\}\cup [J]$ with $J=\lceil {\frac{1}{2\alpha_Y+D_Y+d_X\frac{\alpha_Y}{\alpha_X}}}\cdot \log_2 (\frac{n}{\log n})\rceil$, we define a class of mappings $\m S_j$  on $\Psi_j^{D_Y}\times \mb B_{\mb R^{D_X}}(\mathbf{0},L)$  as
\begin{equation}\label{defSr1}
\begin{aligned}
     \m S_j&=\Bigg\{S(\psi,x)=\sum_{\wt\psi\in \Psi_J^{D_Y}}\frac{\sum_{i=1}^{W_j} \sum_{k\in \mb N_0^{D_X}, |k|<\alpha_X}a_{ik}^{\wt\psi} (x-b_i)^k \rho(\frac{\|x-b_i\|}{\varepsilon^x_j})}{\sum_{i=1}^{W_j} \rho(\frac{\|x-b_i\|}{\varepsilon^x_j})+\frac{1}{n}}\cdot \mathbf{1}(\wt\psi=\psi):\,
     b_i\in \mb B_{\mb R^{D_X}}(\mathbf{0},L), \\ &a_{ik}^{\wt\psi}\in [-\frac{C}{2^{D_Yj/2}},\frac{C}{2^{D_Yj/2}}],  \text{ for any } i\in [W_j],\wt\psi\in \Psi_j^{D_Y} \text{ and } k\in \mb N_0^{D_X} \text{ with } |k|<\alpha_X \Bigg\},
\end{aligned}
\end{equation}
where $\varepsilon^x_j=2^{\frac{jD_Y}{2\alpha_X+d_X}}(\frac{n}{\log n})^{-\frac{1}{2\alpha_X+d_X}}$, $W_j=C_1\,(\varepsilon^x_j)^{-d_X} $ and $C,C_1$ are large enough constants. Then consider the estimator 
\begin{equation}\label{Mestimation}
    \wh S_j={\arg\min}_{S\in  \m S_j} \frac{1}{n}\sum_{i=1}^n \sum_{\psi\in \Psi_j^{D_Y}} (\psi(Y_i)-S(\psi,X_i))^2.
\end{equation}
It is straightforward to verify that $\wh S_j(\psi,x)=\wh u_{\psi}(x)$ for any $\psi\in \Psi_j^{D_Y}$ and $x\in \m M_X$, and we can express
 \begin{equation*}
  \wh u(\cdot|x)=\sum_{j=0}^J \sum_{\psi\in \Psi_j^{D_Y}} \psi(\cdot) \wh S_j(\psi,x).
 \end{equation*}
We then present the following lemma to bound the mean squared error between $\wh S_j(\psi,x)$ and $u_{\psi}^*(x)=\int_{\mb R^{D_Y}} u^*(y\,|\,x)\psi(y)\,\dd y$, where $u^*(y\,|\,x)$ is the density function of $\mu^*_{Y|x}$ with respect to the Lebesgue measure on $\mb R^{D_Y}$.

\begin{lemma}\label{lemma1.1}
   Suppose $\mu^*\in \m P^*_1$ and  with the choices of $\m S_j$ defined in~\eqref{defSr1}, there exists a constant $C$ so that  it holds with probability at least $1-\frac{1}{n}$ that for any $j\in [J]$,
      \begin{equation*}
          \mb{E}_{\mu^*_X} \bigg[\sum_{\psi\in \Psi_j^{D_Y}} (\wh S_j(\psi,X)-u^*_{\psi}(X))^2\bigg]\leq C\, 2^{\frac{2j\alpha_XD_Y}{2\alpha_X+d_X}} (\frac{n}{\log n})^{-\frac{2\alpha_X}{2\alpha_X+d_X}}.
      \end{equation*}
\end{lemma}
\noindent The proof of Lemma~\ref{lemma1.1} is given in Appendix~\ref{proof:lemma1.1}. For ease of notation, we define $\wh S_j(\psi,x)=0$ for any $j>J$.  Then, the estimator $\wh u(y\,|\,x)$ can be rewritten as 
\begin{equation*}
    \wh u(y\,|\,x)=\sum_{j=0}^{\infty}\sum_{\psi\in \Psi_j^{D_Y}}\psi(y)\wh S_j(\psi,x).
\end{equation*}
For any $\gamma\geq 0$, we can obtain the following bound:
  \begin{equation*}
  \begin{aligned}
        & \mb{E}_{\mu^*_X} \Big[\underset{f\in \m H_1^{\gamma}(\mb R^{D_Y})}{\sup}\int_{\mb R^{D_Y}} f(y)u^*(y\,|\,X)\,\dd y-\int_{\mb R^{D_Y}} f(y)\wh u(y\,|\,X)\,\dd y \Big]\\
       &=\mb{E}_{\mu^*_X} \Big[\underset{f\in \m H^{\gamma}_1(\mb R^{D_Y})}{\sup}\sum_{j=0}^{\infty} \sum_{\psi\in \Psi_j^{D_Y}}f_{\psi} (u^*_{\psi}(X)-\wh S(\psi,X))\Big] \\
&\leq \mb{E}_{\mu^*_X} \Big[ \underset{f\in \m H_1^{\gamma}(\mb R^{D_Y})}{\sup}\sum_{j=0}^{J} \sum_{\psi\in \Psi_j^{D_Y}} f_{\psi}\cdot\big(u^*_{\psi}(X)-\wh S(\psi,X)\big)\Big]+\mb{E}_{\mu^*_X} \Big[ \underset{f\in \m H_1^{\gamma}(\mb R^{D_Y})}{\sup}\sum_{j=J+1}^{\infty} \sum_{\psi\in \Psi_j^{D_Y}} f_{\psi}\cdot u^*_{\psi}(X)\Big]\\
           &\overset{(i)}{\leq} C\, \sum_{j=0}^{J}\sum_{\psi\in \Psi_j^{D_Y}}2^{-j\gamma-\frac{jD_Y}{2}}\sqrt{\mb{E}_{\mu^*_X}\Big[ (u^*_{\psi}(X)-\wh S(\psi,X))^2\Big]}+ C\,\sum_{j=J+1}^{\infty} \sum_{\psi\in \Psi_j^{D_Y}} 2^{-j(\gamma+\alpha_Y)-jD_Y}\\
     &\leq C\,  \sum_{j=0}^{J}\sqrt{\sum_{\psi\in \Psi_j^{D_Y}}2^{-2j\gamma}\mb{E}_{\mu^*_X}\Big[ \ (u^*_{\psi}(X)-\wh S(\psi,X))^2\Big]}+C\,\sum_{j=J+1}^{\infty} \sum_{\psi\in \Psi_j^{D_Y}} 2^{-j(\gamma+\alpha_Y)-jD_Y},
  \end{aligned}
  \end{equation*}
  where $(i)$ uses $|f_{\psi}|\lesssim 2^{-j\gamma-jD_Y/2}$ for $\psi\in \Psi_j^{D_Y}$, and $u^*(\,\cdot\,|\,x)\in \m H^{\alpha_Y}_L(\mb R^{D_Y})$, implying that for $\psi\in \Psi_j^{D_Y}$, $|u^*_{\psi}(x)|\lesssim 2^{-\frac{D_Yj}{2}-j\alpha_Y}$, alongside the Jensen's inequality; the last inequality is derived using Cauchy-Schwarz inequality and $|\Psi_j^{D_Y}|\lesssim 2^{D_Yj}$.
  Finally, using Lemma~\ref{lemma1.1}, we can get it holds with probability at least $1-\frac{1}{n}$ that for any $\gamma\geq 0$,
    \begin{equation*}
  \begin{aligned}
        & \mb{E}_{\mu^*_X} \Big[\underset{f\in \m H_1^{\gamma}(\mb R^{D_Y})}{\sup}\int_{\mb R^{D_Y}} f(y)u^*(y\,|\,X)\,\dd y-\int_{\mb R^{D_Y}} f(y)\wh u(y\,|\,X)\,\dd y \Big]\\
        &\leq C\, \sum_{j=0}^J 2^{-j\gamma}  2^{\frac{j\alpha_XD_Y}{2\alpha_X+d_X}} (\frac{n}{\log n})^{-\frac{\alpha_X}{2\alpha_X+d_X}} +C\, 2^{-J(\gamma+\alpha_Y)}\\
       &\leq C_1\, (\log n)\cdot (\frac{n}{\log n})^{-\frac{\alpha_X}{2\alpha_X+d_X}}+C_1\,(\frac{n}{\log n})^{-\frac{\alpha_Y+\gamma}{2\alpha_Y+D_Y+\frac{\alpha_Y}{\alpha_X}d_X}}.
  \end{aligned}
  \end{equation*}
  This completes the proof.

\subsection{Proof of Theorem~\ref{th1:lower} (minimax lower bound for Regime 1)}\label{proofth1:lower}
The upper bound can be directly derived from Theorem~\ref{th1}, so here we focus solely on establishing the lower bound. Notice that the lower bound for $d_X=0$ follows directly follows from the minimax rate for the unconditional case (see for example, Theorem 4 of~\cite{JMLR:v22:20-911}). Therefore, we will assume $d_X\in \mb N_{+}$ in the following.
\subsubsection{Proof for the lower bound of $n^{-\frac{\alpha_Y+\gamma}{2\alpha_Y+D_Y+\frac{\alpha_Y}{\alpha_X}d_X}}$}
Define the covariate space $\m M_X=[0,1]^{d_X}\times \mathbf 0_{D_X-d_X}$ and  let $\mu^*_X$ be the uniform distribution over $\m M_X$. Then let $\wt m_1 = \lceil\, b\,n^{\frac{1}{2\alpha_Y+D_Y+\frac{\alpha_Y}{\alpha_X}d_X}}\rceil$ and $\wt m_2 = \lceil b\, n^{\frac{1}{2\alpha_X+d_X+\frac{\alpha_X}{\alpha_Y}D_Y}}\rceil$ where $b$ is a large enough positive constant. Consider the following bump function 
\begin{equation}\label{eqnk.2}
 \wt k(t)=\left\{\begin{array}{l}
(1-t)^{\alpha_Y\vee \alpha_X\vee\gamma+1} t^{\alpha_Y\vee\alpha_X\vee\gamma+1}(t-\frac{1}{2}), \quad t \in(0,1) \\
0, \quad \text { o.w. }
\end{array}\right.
\end{equation}
so that $\int_{-\infty}^\infty \wt k(t) \,\dd t=0$, and the corresponding localized bump function over $\mb R^{D_Y}\times \mb R^{D_X}$,

\begin{equation}\label{eqnk.2.2}
    \wt \psi_{\xi_1,\xi_2}(y,x)=\prod_{i=1}^{D_Y} \wt k\Big(\wt m_1\sqrt{\frac{D_Y}{2}}y_i+\frac{\wt m_1}{2}-\xi_{1i}\Big)\prod_{i=1}^{d_X} \wt k\Big(\wt m_2\sqrt{2d_X}x_i-\xi_{2i}\Big),\quad\forall y\in \mb B_{\mb R^{D_Y}}(\mathbf{0},1),
\end{equation}
indexed by the $D_Y$-dimensional grid point $\xi_1=(\xi_{11},\ldots,\xi_{1D_Y})\in[\wt m_1]^{D_Y}$ and the $d_X$-dimensional grid $\xi_2=(\xi_{21},\ldots,\xi_{2d_X})\in[\wt m_2]^{d_X}$, where we have used the notation $[m]^d=\{(i_1,i_2,\cdots,i_d):\, i_k\in \{1,2,\cdots, m\},\,  \forall k\in \{1,2,\cdots,d\}\}$. Define the baseline density function
\begin{equation}\label{eqnk.2.3}
    \nu_0(y)=\left\{\begin{array}{cc}
\frac{     \prod_{i=1}^{D_Y} (1-y_i)^{\alpha_Y\vee \gamma+1} (y_i+1)^{\alpha_Y\vee\gamma+1}}{(\int_0^1  (1-t)^{\alpha_Y\vee \gamma+1} (t+1)^{\alpha_Y\vee\gamma+1}\,\dd t)^{D_Y}}  & y\in [-1,1]^{D_Y} \\
       0  & o.w.
    \end{array}
    \right.
\end{equation}
and two function sets 
\begin{equation}\label{defpsialpha}
\begin{aligned}
\Psi_{\alpha_Y,\alpha_X}&=\Big\{ \nu_{\omega}(y,x)=\nu_0(y)+\Big(\frac{1}{\wt m_1}\Big)^{\alpha_Y}\sum_{\xi_1\in[\wt m_1]^{D_Y}}\sum_{\xi_2\in[\wt m_2]^{d_X}}  \omega_{\xi_1,\xi_2} \,\wt \psi_{\xi_1,\xi_2}(y,x)\\
&\qquad:\  \omega=\{\omega_{\xi_1,\xi_2}\}_{\xi_1\in[\wt m_1]^{D_Y},\xi_2\in[\wt m_2]^{d_X}}\in \{0,1\}^{\wt m_1^{D_Y}\times \wt m_2^{d_X}}\Big\},\\
\Lambda_{\gamma}&=\Big\{ f_{v}(y,x)=\Big(\frac{1}{\wt m_1}\Big)^{\gamma}\sum_{\xi_1\in[\wt m_1]^{D_Y}}\sum_{\xi_2\in[\wt m_2]^{d_X}}  v_{\xi_1,\xi_2} \,\wt \psi_{\xi_1,\xi_2}(y,x)\\
&\qquad:\  v=\{v_{\xi_1,\xi_2}\}_{\xi_1\in[\wt m_1]^{D_Y},\xi_2\in[\wt m_2]^{d_X}}\in \{-1,1\}^{\wt m_1^{D_Y}\times \wt m_2^{d_X}}\Big\},
\end{aligned}
\end{equation}
 Here, $\Psi_{\alpha_Y,\alpha_X}$ consists of all perturbed conditional densities around $\nu_0(\cdot)$ and $\Lambda_\gamma$ serves as set of discriminators for discriminating the conditional densities in $\Psi_{\alpha_Y,\alpha_X}$. Moreover, $\wt{\psi}_{\xi_1,\xi_2}$'s with distinct indices $(\xi_1,\xi_2)$'s have disjoint supports and when $b$ is sufficiently large, we have for each $\nu\in \Psi_{\alpha_Y,\alpha_X}$: $\nu(y,x) =\nu_0(y,x)$ for all $(y,x)\notin \mb B_{\mb R^{D_Y}}(\mathbf{0},3/4)\times \mb B_{\m M_X}(\mathbf{0},3/4)$; and $\nu(y,x) \geq \inf_{y\in \mb B_{\mb R^{D_Y}}(\mathbf{0},3/4)}\nu_0(y)- b^{-\alpha_Y} \sup_{t\in (0,1)}|\wt k(t)|^{D_Y+d_X} >0$ for all $y\in \mb R^{D_Y}(\mathbf{0},3/4)$ and $x\in \m M_X$, which makes $\nu$ non-negative. In addition, since $\int_{-\infty}^\infty \wt k(t) \,\dd t=0$, we have $\int_{\mb R^{D_Y}} \nu(y,x)\,\dd y = \int_{\mb R^{D_Y}} \nu_0(y)\,\dd y =1$. Therefore, all functions in $\Psi_{\alpha_Y,\alpha_X}$ are valid conditional probability density functions.  Furthermore, we state the following lemma that verifies the smoothness of functions in $\Psi_{\alpha_Y,\alpha_X}$ and $\Lambda_{\gamma}$, the proof of which is given in Appendix~\ref{proof:lemmalowerboundsmooth}.
 
 \begin{lemma}\label{lemmalowerboundsmooth}
   Let $\phi_1\in \m H^{\lceil\alpha_1\rceil}_{L}(\mb R^{d_1})$, $\phi_2\in \m H^{\lceil\alpha_2\rceil}_{L}(\mb R^{d_2})$ be two compactly supported functions.  Consider the function 
     \begin{equation*}
         f(x,y)=(\frac{1}{m_1})^{\alpha_1}  \sum_{\xi_1\in [m_1]^{d_1}}\sum_{\xi_2\in [m_2]^{d_2}}  \omega_{\xi_1,\xi_2}\phi_1(m_1x-\xi_1)\phi_2(m_2y-\xi_2).
     \end{equation*}
 For any positive constants $C,C_1,C_2$, there exists a constant $L_1$ so that for any $m_1,m_2\in \mb N_{+}$ with $C_1m_2^{\alpha_2}\leq m_1^{\alpha_1}\leq C_2m_2^{\alpha_2}$,  and any $\omega_{\xi_1,\xi_2}\in [-C,C]$, it holds that $f\in  {\m H}^{\alpha_1,\alpha_2}_{L_1}(\mb R^{d_1},\mb R^{d_2})$.
 \end{lemma}
\noindent Therefore, there exist constants $(L_1,L_2)$ such that $\Psi_{\alpha_Y,\alpha_X}\subset  {\m H}^{\alpha_Y,\alpha_X}_{L_1}(\mb R^{D_Y},\m M_X)\subset \ov{\m H}^{\alpha_Y,\alpha_X}_{L_2}(\mb R^{D_Y},\m M_X)$ and for any  $f\in \Lambda_\gamma$ and $x\in \m M_X$, it holds that $f(\cdot,x)\in \m H^{\gamma}_{L_2}(\mathbb{R}^{D_Y})$.  Then for each $\omega \in \{0,1\}^{\wt m_1^{D_Y}\times \wt m_2^{d_X}}$, we define the conditional distribution $\mu^{\omega}_{Y|X}$ of $Y|X$ as $\mu^{\omega}_{Y|X}=\nu_{\omega}(y,X) \dd y$ and the joint distribution of $(X,Y)$ as $ \mu^{\omega}=\mu^*_X\mu^{\omega}_{Y|X}$. Then there exists a constant $L$ so that $\mu^{\omega}\in \m P^*_1(D_Y,D_X,d_X,\alpha_Y,\alpha_X,L)$. Next,  by the Varshamov-Gilbert lemma~\citep{Tsybakov2009},  there exists a set $\{\omega^{(0)},\cdots , \omega^{(H)}\}\subset \{0,1\}^{\wt m_1^{D_Y}\times \wt m_2^{d_X}}$ such that $\log H \geq \frac{\wt m_1^{D_Y}\wt m_2^{d_X}}{8}\log 2$ and the Hamming distance $\|\omega^{(j)}-\omega^{(k)}\|_{\rm H}\geq \frac{\wt m_1^{D_Y}\wt m_2^{d_X}}{8}$ for any distinct pair $j,k\in [H]$. Therefore, for any distinct $j, k\in[H]$, we have by our construction of $\mu^{\omega}$'s that 
 \begin{equation*}
\begin{aligned}
&\mb{E}_{\mu^*_X}[d_{\gamma}(\mu^{\omega^{(j)}}_{Y|X}, \mu^{\omega^{(k)}}_{Y|X})]=\mb{E}_{\mu^*_X}\Big[\underset{f\in \m H^{\gamma}_1(\mathbb{R}^{D_Y})}{\sup}\int_{\mb R^{D_Y}} f(y)\cdot (\nu_{\omega^{(j)}}(y,X)-\nu_{\omega^{(k)}}(y,X)) \,\dd y\Big]\\
&\geq \frac{1}{L_2}  \mb{E}_{\mu^*_X}\Big[\underset{f\in \Lambda_{\gamma}}{\sup}\int_{\mb R^{D_Y}} f(y,X)\cdot (\nu_{\omega^{(j)}}(y,X)-\nu_{\omega^{(k)}}(y,X)) \,\dd y\Big]\\
&\geq \frac{1}{L_2}  \underset{f\in \Lambda_{\gamma}}{\sup} \, \mb{E}_{\mu^*_X}\Big[\int_{\mb R^{D_Y}} f(y,X)\cdot (\nu_{\omega^{(j)}}(y,X)-\nu_{\omega^{(k)}}(y,X)) \,\dd y\Big]\\
&=\frac{1}{L_2}\underset{v\in \{-1,1\}^{\wt m_1^{D_Y}\times \wt m_2^{d_X}}}{\sup} (\frac{1}{\wt m})^{\alpha_Y+\gamma} \int_{[0,1]^{d_X}}\int_{\mb R^{D_Y}} \sum_{\xi_1\in[\wt m_1]^{D_Y}}\sum_{\xi_2\in[\wt m_2]^{d_X}}  v_{\xi_1,\xi_2}(\omega_{\xi_1,\xi_2}^{(j)}-\omega_{\xi_1,\xi_2}^{(k)})    \wt \psi^2_{\xi_1,\xi_2}(y,(x, \mathbf{0}_{D_X-d_X}))\,\dd y\dd x\\
&\gtrsim (\frac{1}{\wt m_1})^{\alpha_Y+\gamma}  \cdot (\frac{1}{\wt m_1})^{D_Y}(\frac{1}{\wt m_2})^{d_X}\cdot\|\omega^{(j)}-\omega^{(k)}\|_{\rm H}\\
&\gtrsim (\frac{1}{\wt m_1})^{\alpha_Y+\gamma}.
 \end{aligned}
\end{equation*}
  Moreover, we have
\begin{equation*}
\begin{aligned}
&D_{\rm KL}(    \mu^{\omega^{(j)}},\mu^{\omega^{(k)}})\\
&=\mb{E}_{\mu^*_X}\bigg[\int_{\mb R^{D_Y}} -\log\bigg(\underbrace{
\frac{\nu_0(y)+\Big(\frac{1}{\wt m_1}\Big)^{\alpha_Y}\sum_{\xi_1\in[\wt m_1]^{D_Y}}\sum_{\xi_2\in[\wt m_2]^{d_X}}  \omega^{(k)}_{\xi_1,\xi_2} \wt \psi_{\xi_1,\xi_2}(y,x)}{\nu_0(y)+\Big(\frac{1}{\wt m_1}\Big)^{\alpha_Y}\sum_{\xi_1\in[\wt m_1]^{D_Y}}\sum_{\xi_2\in[\wt m_2]^{d_X}}  \omega^{(j)}_{\xi_1,\xi_2} \wt \psi_{\xi_1,\xi_2}(y,x)}}_{:\,=1+u(y,x)}
\bigg) \nu_{\omega^{(j)}}(y,x) \, \dd y\bigg]
 \end{aligned}
\end{equation*}
For sufficiently large $b$, we have $|u(y,x)| \leq 1/4$ so that $-\log(1+u(y,x)) \leq u^2(y,x) - u(y,x)$. 
This leads to
\begin{equation}\label{eqnboundKL}
\begin{aligned}
&D_{\rm KL}(    \mu_{\omega^{(j)}},\mu_{\omega^{(k)}})\
\leq C \Big(\frac{1}{\wt m_1}\Big)^{2\alpha_Y} \\
&+ \Big(\frac{1}{\wt m_1}\Big)^{\alpha_Y}\int_{[0,1]^{d_X}}\int_{\mb R^{D_Y}}
 \bigg\{\sum_{\xi_1\in[\wt m_1]^{D_Y}}\sum_{\xi_2\in[\wt m_2]^{d_X}} ( \omega_{\xi}^{(j)}-\omega_{\xi}^{(k)})\cdot \psi_{\xi_1,\xi_2}(y,(x,\mathbf 0_{D_X-d_X}))  \bigg\}\,  \dd y\,\dd x =C\, (\frac{1}{\wt m})^{2\alpha_Y},
\end{aligned}
\end{equation}
where we used the fact that $\int_{\mb R^{D_Y}} \psi_{\xi_1,\xi_2}(y,(x,\mathbf 0_{D_X-d_X}))\,\dd y=0$. Then we can apply Fano's lemma (proposition 15.12 of~\cite{wainwright2019high}) to obtain
\begin{equation}\label{fano}
    \begin{aligned}
    & \underset{\wh{\mu}_{Y|X}}{\inf}\underset{\mu\in \mathcal{P}_1^{\ast}} {\sup} \,\mb{E}_{\mu^{\otimes n} }\mathbb{E}_{\mu_X} \big[d_{\gamma}(\wh{\mu}_{Y|X}, \mu_{Y|X})\big]\geq     \underset{\wh{\mu}_{Y|X}}{\inf}\underset{j\in [H]} {\sup} \,\mb{E}_{\mu^{\omega^{j}\otimes n} }\mathbb{E}_{\mu_X^*} \big[d_{\gamma}(\wh{\mu}_{Y|X}, \mu^{\omega^{(j)}}_{Y|X})\big]\\
    &  \geq \frac{1}{2}\underset{\wh{\mu}_{Y|X}\in\{ \mu^{\omega^{(j)}}_{Y|X}:\, j\in [H]\} }{\inf}\,\underset{j\in [H]}{\sup} \,\mb{E}_{\mu^{\omega^{j}\otimes n} }\mathbb{E}_{\mu_X^*} \big[d_{\gamma}(\wh{\mu}_{Y|X}, \mu^{\omega^{(j)}}_{Y|X})\big]\\
    & \geq   \,  \frac{1}{2}\,\underset{j,k\in [H]\atop j\neq k} {\inf} \mb{E}_{\mu^*_X}[d_{\gamma}(\mu^{\omega^{(j)}}_{Y|X},\mu^{\omega^{(k)}}_{Y|X})] \cdot\bigg(1-\frac{\log 2+\frac{n}{H^2} \sum_{j,k =1}^{H}D_{\rm KL}(    \mu^{\omega^{(j)}},\mu^{\omega^{(k)}})}{\log H}\bigg)\\
    &\gtrsim  n^{-\frac{\alpha_Y+\gamma}{2\alpha_Y + D_Y+d_X\frac{\alpha_Y}{\alpha_X}}}.
\end{aligned}
\end{equation}

\subsubsection{Proof for the lower bound of $n^{-\frac{\alpha_X}{2\alpha_X+d_X}}$}
Consider the same covariate space $\m M_X=[0,1]^{d_X}\times \mathbf 0_{D_X-d_X}$ and  uniform distribution $\mu^*_X$ over $\m M_X$. Define $\wt m = \lceil b n^{\frac{1}{2\alpha_X+d_X}}\rceil$, where $b$ is a large enough positive constant. Consider $\wt k(\cdot)$ as defined in~\eqref{eqnk.2} and the  localized bump function over $\mb R^{D_X}$,
\begin{equation} 
    \wt \psi_{\xi}(x)=\prod_{i=1}^{d_X} \wt k\Big(\wt m\sqrt{2d_X}x_i-\xi_{i}\Big)
\end{equation}
indexed by the $d_X$-dimensional grid $\xi=(\xi_{1},\ldots,\xi_{d_X})\in[\wt m]^{d_X}$. Then define two function sets 
\begin{equation}\label{defpsialpha1}
\begin{aligned}
\Psi_{\alpha_X}&=\Big\{ \nu_{\omega}(y,x)=\nu_0(y)+\Big(\frac{1}{\wt m}\Big)^{\alpha_X}\sum_{\xi\in[\wt m]^{d_X}}  \omega_{\xi} \,\wt \psi_{\xi}(x) \prod_{i=1}^{D_Y} \wt k(y_i):\, \omega=\{\omega_{\xi}\}_{\xi\in[\wt m]^{d_X}}\in \{0,1\}^{\wt m^{d_X}}\Big\},\\
\Lambda_{\gamma}&=\Big\{ f_{v}(y,x)=\sum_{\xi\in[\wt m]^{d_X}} v_{\xi} \,\wt \psi_{\xi}(x) \prod_{i=1}^{D_Y} \wt k(y_i):\,  v=\{v_{\xi}\}_{\xi\in[\wt m]^{d_X}}\in \{-1,1\}^{\wt m^{d_X}}\Big\},
\end{aligned}
\end{equation}
where  $\nu_0$ is defined in~\eqref{eqnk.2.3}. Then  it is straightforward to verify that there exist constants $(L_1,L_2)$ such that $\Psi_{\alpha_X}\subset {\m H}^{\alpha_Y,\alpha_X}_{L_1}(\mb R^{D_Y},\m M_X)$ and for any $f\in \Lambda_\gamma$ and $x\in \m M_X$, $f(\cdot,x)\in \m H^{\gamma}_{L_2}(\mathbb{R}^{D_Y})$.  Moreover,  $\nu_{\omega}$'s in $\Psi_{\alpha_X}$ are valid probability density functions. Then for each $\omega \in \{0,1\}^{ \wt m^{d_X}}$,  we define the conditional distribution $\mu^{\omega}_{Y|X}$ of $Y|X$ as $\mu^{\omega}_{Y|X}=\nu_{\omega}(y,X) \dd y$ and the joint distribution of $(X,Y)$ as $ \mu^{\omega}=\mu^*_X\mu^{\omega}_{Y|X}$. Then there exists a constant $L$ so that $\mu^{\omega}\in \m P^*_1(D_Y,D_X,d_X,\alpha_Y,\alpha_X,L)$.  Next,  by the Varshamov-Gilbert lemma~\citep{Tsybakov2009},  there exists a set $\{\omega^{(0)},\cdots , \omega^{(H')}\}\subset \{0,1\}^{\wt m^{d_X}}$ such that $\log H' \geq \frac{\wt m^{d_X}}{8}\log 2$ and the Hamming distance $\|\omega^{(j)}-\omega^{(k)})\|_{\rm H}\geq \frac{\wt m^{d_X}}{8}$ for any distinct pair $j,k\in [H']$. Therefore, for any distinct $j, k\in[H']$, we have by our construction of $\mu^{\omega}$'s that 
 \begin{equation*}
\begin{aligned}
&\mb{E}_{\mu^*_X}[d_{\gamma}(\mu^{\omega^{(j)}}_{Y|X}, \mu^{\omega^{(k)}}_{Y|X})]=\mb{E}_{\mu^*_X}\Big[\underset{f\in \m H^{\gamma}_1(\mathbb{R}^{D_Y})}{\sup}\int_{\mb R^{D_Y}} f(y)\cdot (\nu_{\omega^{(j)}}(y,X)-\nu_{\omega^{(k)}}(y,X)) \,\dd y\Big]\\
&\geq \frac{1}{L_2}  \mb{E}_{\mu^*_X}\Big[\underset{f\in \Lambda_{\gamma}}{\sup}\int_{\mb R^{D_Y}} f(y,X)\cdot (\nu_{\omega^{(j)}}(y,X)-\nu_{\omega^{(k)}}(y,X)) \,\dd y\Big]\\
&=\frac{1}{L_2} \underset{f\in \Lambda_{\gamma}}{\sup} \, \mb{E}_{\mu^*_X}\Big[\int_{\mb R^{D_Y}} f(y,X)\cdot (\nu_{\omega^{(j)}}(y,X)-\nu_{\omega^{(k)}}(y,X)) \,\dd y\Big]\\
&=\frac{1}{L_2} \underset{v\in \{-1,1\}^{\wt m^{d_X}}}{\sup} (\frac{1}{\wt m})^{\alpha_X} \int_{[0,1]^{d_X}}\int_{\mb R^{D_Y}} \sum_{\xi\in[\wt m]^{d_X}}  v_{\xi}(\omega_{\xi}^{(j)}-\omega_{\xi}^{(k)})    \wt \psi^2_{\xi}(x,\mathbf{0}_{D_X-d_X}) \prod_{i=1}^{D_Y} \wt k(y_i)^2\,\dd y\dd x\\
&\gtrsim (\frac{1}{\wt m})^{\alpha_X}  \cdot (\frac{1}{\wt m})^{d_X}\cdot\|\omega^{(j)}-\omega^{(k)}\|_{\rm H}\\
&\gtrsim (\frac{1}{\wt m})^{\alpha_X}.
 \end{aligned}
\end{equation*}
  Moreover, similar to~\eqref{eqnboundKL}, we can derive
  \begin{equation*}
\begin{aligned}
&D_{\rm KL}(    \mu^{\omega^{(j)}},\mu^{\omega^{(k)}})\\
&=\mb{E}_{\mu^*_X}\int_{\big[0,1\big]^{D_Y}} -\log\bigg(
\frac{\nu_0(y)+\Big(\frac{1}{\wt m}\Big)^{\alpha_X}\sum_{\xi\in[\wt m]^{d_X}}  \omega_{\xi}^{(j)} \,\wt \psi_{\xi}(x) \prod_{i=1}^{D_Y} \wt k(y_i) }{\nu_0(y)+\Big(\frac{1}{\wt m}\Big)^{\alpha_X}\sum_{\xi\in[\wt m]^{d_X}}  \omega_{\xi}^{(k)} \,\wt \psi_{\xi}(x) \prod_{i=1}^{D_Y} \wt k(y_i)} 
\bigg) \nu_{\omega^{(j)}}(y,x) \, \dd y  
\lesssim (\frac{1}{\wt m})^{2\alpha_X},
\end{aligned}
\end{equation*}
where we used the fact that $\int_{\mb R^{D_Y}} \prod_{i=1}^{D_Y} \wt k(y_i)\,\dd y=0$. Then we can apply Fano's lemma to obtain
\begin{equation*}
    \begin{aligned}
   & \underset{\wh{\mu}_{Y|X}}{\inf}\underset{\mu\in \mathcal{P}_1^{\ast}} {\sup} \,\mb{E}_{\mu^{\otimes n} }\mathbb{E}_{\mu_X} \big[d_{\gamma}(\wh{\mu}_{Y|X}, \mu_{Y|X})\big]\geq     \underset{\wh{\mu}_{Y|X}}{\inf}\underset{j\in [H]} {\sup} \,\mb{E}_{\mu^{\omega^{j}\otimes n} }\mathbb{E}_{\mu_X^*} \big[d_{\gamma}(\wh{\mu}_{Y|X}, \mu^{\omega^{(j)}}_{Y|X})\big]\\
    &  \geq \frac{1}{2}\underset{\wh{\mu}_{Y|X}\in\{ \mu^{\omega^{(j)}}_{Y|X}:\, j\in [H]\} }{\inf}\,\underset{j\in [H]}{\sup} \,\mb{E}_{\mu^{\omega^{j}\otimes n} }\mathbb{E}_{\mu_X^*} \big[d_{\gamma}(\wh{\mu}_{Y|X}, \mu^{\omega^{(j)}}_{Y|X})\big]\\
    & \geq   \,  \frac{1}{2}\,\underset{h,\ell\in [H']\atop h\neq \ell} {\inf} \mb{E}_{\mu^*_X}[d_{\gamma}(\mu^{\omega^{(h)}}_{Y|X},\mu^{\omega^{(\ell)}}_{Y|X})] \cdot\bigg(1-\frac{\log 2+\frac{n}{H'{}^2} \sum_{h,\ell =1}^{H'}D_{\rm KL}(    \mu^{\omega^{(h)}},\mu^{\omega^{(\ell)}})}{\log H'}\bigg)\\
    &\gtrsim  n^{-\frac{\alpha_X}{2\alpha_X +  d_X}}.
\end{aligned}
\end{equation*}
\subsection{Proof of Lemma~\ref{lemma1.1}}\label{proof:lemma1.1}
We first derive an oracle inequality  in the following lemma,
\begin{lemma}\label{lestimation}
Suppose $\mu^*\in \m P^*_1$ and with the choices of $\m S_j$ defined in~\eqref{defSr1},  it holds with probability larger than  $1-\frac{1}{n}$ that for any $j \in [J]$,
\begin{equation*}
\begin{aligned}
  \mb{E}_{\mu^*_X}  \bigg[\sum_{\psi\in \Psi_j^{D_Y}} (\wh S_j(\psi,x)-u^*_{\psi}(x))^2  \bigg]
   \lesssim  \frac{ 2^{D_Yj} W_j \log n}{n}
 +\underset{S\in \m S_j}{\min} \, \mb{E}_{\mu^*_X}  \bigg[\sum_{\psi\in \Psi_j^{D_Y}} ( S(\psi,x)-u^*_{\psi}(x))^2 \bigg].\\
\end{aligned}
\end{equation*}
\end{lemma}
\noindent  The proof of Lemma~\ref{lestimation} is provided in Appendix~\ref{proof:lestimation}. Then we provide an upper bound for the approximation error given by $\underset{S\in \m S_j}{\min} \, \mb{E}_{\mu^*_X}\big[ \sum_{\psi\in \Psi_j^{D_Y}} ( S(\psi,x)-u^*_{\psi}(x))^2\big]$. Fix an arbitrary $j\in [J]$ and considering $u^*\in  {\m H}^{\alpha_Y,\alpha_X}_{L}(\mb R^{D_Y},\m M_X)$,  there exists $\ov u^*\in  {\m H}^{\alpha_Y,\alpha_X}_{L}(\mb R^{D_Y},\mb R^{D_X})$ so that $\ov u^*|_{\mb R^{D_Y}\times \m M_X}=u^*$. Consequently, there exists a constant $L_1$ so that for any $\psi\in \Psi_j^{D_Y}$,
\begin{equation*}
    2^{\frac{D_Yj}{2}}u^*_{\psi}(x)=2^{\frac{D_Yj}{2}}\int_{\mb R^{D_Y}} \psi(y) \ov u^*_{\psi}(y|x)\, \dd y\in \m H^{\alpha_X}_{L_1}(\mb R^{D_X}),
\end{equation*}
where we have used the fact that the support of $\psi(y)$ has  a volume of $\m O(2^{-jD_Y})$ and $|\psi(y)|=\m O(2^{\frac{D_Y j}{2}})$. Let ${\m N}^x_{\varepsilon^x_j}$ denote the largest $\varepsilon^x_j$-packing set of $\m M_X$, then   for large enough constant $C_1$, we have
$|{\m N}^x_{\varepsilon^x_j}|\leq W_j=C_1(\varepsilon_j^x)^{-d_X}$. Then we define a set $ \ov{\m N}^x_{\varepsilon^x_j}= {\m N}^x_{\varepsilon^x_j}\cup \m X$, where $\m X$ is an arbitrary subset of $\m M_X\setminus {\m N}^x_{\varepsilon^x_j}$ with $|\m X|=W_j-| {\m N}^x_{\varepsilon^x_j}|$.  For any $\psi\in \Psi_j^{D_Y}$,  we define
\begin{equation*}
 \wt u_{\psi}(x)=\frac{\sum_{\wt x\in \ov{\m N}^x_{\varepsilon^x_j}} \sum_{k\in \mb N_0^{D_X}, |k|<\alpha_X}u^*_{\psi}{}^{(k)}(\wt x) (x-\wt x)^k \rho(\frac{\|x-\wt x\|}{\varepsilon^x_j})}{\sum_{\wt x\in \ov{\m N}^x_{\varepsilon^x_j}}  \rho(\frac{\|x-\wt x\|}{\varepsilon^x_j})}
\end{equation*}
and for any $x\in \m M_X$,
\begin{equation*}
\begin{aligned}
       S^*_j(\psi,x)&=\frac{\sum_{\wt x\in \ov{\m N}^x_{\varepsilon^x_j}} \sum_{k\in \mb N_0^{D_X}, |k|<\alpha_X}u^*_{\psi}{}^{(k)}(\wt x) (x-\wt x)^k \rho(\frac{\|x-\wt x\|}{\varepsilon^x_j})}{\sum_{\wt x\in \ov{\m N}^x_{\varepsilon^x_j}}  \rho(\frac{\|x-\wt x\|}{\varepsilon^x_j})+\frac{1}{n}}.
\end{aligned}
\end{equation*}
It holds that $S^*_j(\psi,x)\in \m S_j$ and for any $x\in \m M_X$, $\psi\in \Psi_j^{D_Y}$,
\begin{equation*}
    \begin{aligned}
       | \wt u_{\psi}(x)- S^*_j(\psi,x)|&=  \frac{|\sum_{\wt x\in \ov{\m N}^x_{\varepsilon^x_j}} \sum_{k\in \mb N_0^{D_X}, |k|<\alpha_X}u^*_{\psi}{}^{(k)}(\wt x) (x-\wt x)^k \rho(\frac{\|x-\wt x\|}{\varepsilon^x_j})|}{n \cdot(\sum_{\wt x\in \ov{\m N}^x_{\varepsilon^x_j}}  \rho(\frac{\|x-\wt x\|}{\varepsilon^x_j})+\frac{1}{n})(\sum_{\wt x\in \ov{\m N}^x_{\varepsilon^x_j}}  \rho(\frac{\|x-\wt x\|}{\varepsilon^x_j}))}\\
       &\leq \frac{1}{n}\frac{\sum_{\wt x\in \ov{\m N}^x_{\varepsilon^x_j}} \sum_{k\in \mb N_0^{D_X}, |k|<\alpha_X}|u^*_{\psi}{}^{(k)}(\wt x) (x-\wt x)^k |\rho(\frac{\|x-\wt x\|}{\varepsilon^x_j})}{\sum_{\wt x\in \ov{\m N}^x_{\varepsilon^x_j}}  \rho(\frac{\|x-\wt x\|}{\varepsilon^x_j})}\\
       &\leq \frac{1}{n}\cdot\underset{\wt x\in \ov{\m N}^x_{\varepsilon^x_j}, x\in \m M_X}{\sup}\sum_{k\in \mb N_0^{D_X}, |k|<\alpha_X}|u^*_{\psi}{}^{(k)}(\wt x) (x-\wt x)^k |\\
       &\lesssim    2^{-\frac{D_Yj}{2}}n^{-1};
    \end{aligned}
\end{equation*}
\begin{equation*}
    \begin{aligned}
          | \wt u_{\psi}(x)- u^*_{\psi}(x)|&=\frac{\big|\sum_{\wt x\in \ov{\m N}^x_{\varepsilon^x_j}} \big(\sum_{k\in \mb N_0^{D_X}, |k|<\alpha_X}u^*_{\psi}{}^{(k)}(\wt x) (x-\wt x)^k -u^*_{\psi}(x)\big)\rho(\frac{\|x-\wt x\|}{\varepsilon^x_j})\big|}{\sum_{\wt x\in \ov{\m N}^x_{\varepsilon^x_j}}  \rho(\frac{\|x-\wt x\|}{\varepsilon^x_j})}\\
          &\leq \underset{\wt x\in \ov{\m N}^x_{\varepsilon^x_j}, \,x\in \mb{B}_{\m M_X}(\wt x,2\varepsilon_j^x)}{\sup}\big|\sum_{k\in \mb N_0^{D_X}, |k|<\alpha_X}u^*_{\psi}{}^{(k)}(\wt x) (x-\wt x)^k -u^*_{\psi}(x)\big|\\
          &\lesssim  2^{-\frac{D_Yj}{2}} (\varepsilon^x_j)^{\alpha_X}.
    \end{aligned}
\end{equation*}
We can then get 
\begin{equation*}
\begin{aligned}
&\underset{S\in \m S_j}{\min} \, \mb{E}_{\mu^*_X} \Big[\sum_{\psi\in \Psi_j^{D_Y}} ( S(\psi,x)-u^*_{\psi}(x))^2\Big]\\
   &\leq  \mb{E}_{\mu^*_X} \Big[\sum_{\psi\in \Psi_j^{D_Y}}( S^*_j(\psi,x)- u^*_{\psi}(x))^2\Big]\\
   &\lesssim \sum_{\psi\in \Psi_j^{D_Y}}  2^{-{D_Yj}} ((\varepsilon^x_j)^{-\alpha_X}+\frac{1}{n})^2\\
   &\lesssim (\varepsilon^x_j)^{2\alpha_X}+\frac{1}{n^2}.
   \end{aligned}
\end{equation*}
Finally, by substituting $\varepsilon^x_j=2^{\frac{jD_Y}{2\alpha_X+d_X}}(\frac{n}{\log n})^{-\frac{1}{2\alpha_X+d_X}}$ and $W_j\asymp (\varepsilon^x_j)^{-d_X} $,   the desired result follows  directly from lemma~\ref{lestimation}.
\subsection{Proof of Lemma~\ref{lestimation}}\label{proof:lestimation}
To show the desired result, we will apply Theorem~\ref{theoremjointregression} with $\{\psi_{\lambda}(\cdot)\}_{\lambda\in \Lambda}=\Psi_j^{D_Y}$.  We will then proceed by verifying the three assumptions in Theorem~\ref{theoremjointregression}. For the first assumption, note that  for any  $S(\psi,x)=\frac{\sum_{i=1}^{W_j} \sum_{k\in \mb N_0^{D_X}, |k|<\alpha_X}a_{ik}^\psi (x-b_i)^k \rho(\frac{\|x-b_i\|}{\varepsilon^x_j})}{\sum_{i=1}^{W_j} \rho(\frac{\|x-b_i\|}{\varepsilon^x_j})+\frac{1}{n}}\in \m S_j$,
\begin{equation*}
\begin{aligned}
   &\sup_{x\in \m M_X}\sup_{\psi\in \Psi_j^{D_Y}} \frac{\sum_{i=1}^{W_j} \sum_{k\in \mb N_0^{D_X}, |k|<\alpha_X}a_{ik}^\psi (x-b_i)^k \rho(\frac{\|x-b_i\|}{\varepsilon^x_j})}{\sum_{i=1}^{W_j} \rho(\frac{\|x-b_i\|}{\varepsilon^x_j})+\frac{1}{n}}\\
   &\leq \underset{i\in [W_j]}{\sup}\underset{x\in \m M_X}{\sup}\underset{\psi\in \Psi_j^{D_Y}}{\sup } \sum_{k\in \mb N_0^{D_X}, |k|<\alpha_X}a_{ik}^\psi (x-b_i)^k \\
   &\lesssim 2^{-\frac{D_Yj}{2}}.
\end{aligned}
\end{equation*}
Moreover, for any $y$, there exists only a constant-order number of $\psi\in \Psi_j^{D_Y}$ so that $\psi(y)\neq 0$. Therefore, it holds that
\begin{equation*}
\begin{aligned}
    &\sup_{(x,y)\in \m M}\sup_{S\in \m S_j}\sum_{\psi\in \Psi_j^{D_Y}}S^2(\psi,x)+|\psi(y)S(\psi,x)|\\
     &\lesssim \sup_{(x,y)\in \m M}\sum_{\psi\in \Psi_j^{D_Y}} 2^{-D_Y j}+ \sum_{\psi\in \Psi_j^{D_Y}}|\psi (y)|\cdot 2^{-\frac{D_Yj}{2}}\\
     &=\m O(1),
\end{aligned}
\end{equation*}
which verifies the first assumption.
 For the second assumption, let 
 \begin{equation*}
     \ell(x,y,S)=\sum_{\psi\in \Psi_j^{D_Y}} S^2(\psi,x)-2\psi(y)S(\psi,x).
 \end{equation*}
  It holds that
\begin{equation*}
\begin{aligned}
&=\mb{E}_{\mu^*}[(\ell(X,Y,S)-\ell(X,Y,S'))^2]\\
&= \mb{E}_{\mu^*}\Big[\Big(\sum_{\psi\in \Psi_j^{D_Y}}\big(S^2(\psi,X)-S'{}^2(\psi,X)\big)-2\psi(Y) \big(S(\psi,X)-S'(\psi,X)\big)\Big)^2\Big]\\
&=\mb{E}_{\mu^*}\Big[ \Big(\sum_{\psi\in \Psi_j^{D_Y}}\big(S(\psi,X)+S'(\psi,X)-2\psi(Y)\big)\cdot \big(S(\psi,X)-S'(\psi,X)\big)\Big)^2\Big]
\\
&\leq 8\,  \mb{E}_{\mu^*}\Big[ \Big(\sum_{\psi\in \Psi_j^{D_Y}}\psi(Y)\cdot \big(S(\psi,X)-S'(\psi,X)\big)\Big)^2\Big]\\
&\qquad\qquad+2\,\mb{E}_{\mu^*_X}\Big[ \Big(\sum_{\psi\in \Psi_j^{D_Y}}\big(S(\psi,X)+S'(\psi,X)\big)\cdot \big(S(\psi,X)-S'(\psi,X)\big)\Big)^2\Big]\\
&\lesssim  \mb{E}_{\mu^*}\Big[ \Big(\sum_{\psi\in \Psi_j^{D_Y}}\psi(Y)\cdot \big(S(\psi,X)-S'(\psi,X)\big)\Big)^2\Big]+\mb{E}_{\mu^*_X}\Big[ \Big(\sum_{\psi\in \Psi_j^{D_Y}}2^{-\frac{D_Yj}{2}}\cdot \big|S(\psi,X)-S'(\psi,X)\big|\Big)^2\Big].
\end{aligned}
\end{equation*}
Then notice that 
\begin{equation*}
    \begin{aligned}
         &\mb{E}_{\mu^*}\Big[ \Big(\sum_{\psi\in \Psi_j^{D_Y}}\psi(Y)\cdot \big(S(\psi,X)-S'(\psi,X)\big)\Big)^2\Big]\\
         &=\mb{E}_{\mu^*}\Big[ \sum_{\psi_1\in \Psi_j^{D_Y}}\sum_{\psi_2\in \Psi_j^{D_Y}}\psi_1(Y)\psi_2(Y)\cdot \big(S(\psi_1,X)-S'(\psi_1,X)\big)\big(S(\psi_2,X)-S'(\psi_2,X)\big) \Big]\\
         &=\mb{E}_{\mu^*}\Big[ \sum_{\psi_1,\psi_2\in \Psi_j^{D_Y}\atop {\rm supp}(\psi_1)\cap {\rm supp}(\psi_2)\neq \emptyset}\psi_1(Y)\psi_2(Y)\cdot \big(S(\psi_1,X)-S'(\psi_1,X)\big)\big(S(\psi_2,X)-S'(\psi_2,X)\big) \Big]\\
         &=\mb{E}_{\mu^*_X}\Big[  \sum_{\psi_1,\psi_2\in \Psi_j^{D_Y}\atop {\rm supp}(\psi_1)\cap {\rm supp}(\psi_2)\neq \emptyset}\mb{E}_{\mu^*_{Y|X}}[\psi_1(y)\psi_2(y)]\cdot \big(S(\psi_1,X)-S'(\psi_1,X)\big)\big(S(\psi_2,X)-S'(\psi_2,X)\big) \Big]\\
         &\lesssim\mb{E}_{\mu^*_X}\Big[  \sum_{\psi_1,\psi_2\in \Psi_j^{D_Y}\atop {\rm supp}(\psi_1)\cap {\rm supp}(\psi_2)\neq \emptyset}\big(S(\psi_1,X)-S'(\psi_1,X)\big)^2+\big(S(\psi_2,X)-S'(\psi_2,X)\big)^2 \Big]\\
         &\lesssim \mb{E}_{\mu^*_X}\Big[  \sum_{\psi_1\in \Psi_j^{D_Y}}\sum_{\psi_2\in \Psi_j^{D_Y}\atop {\rm supp}(\psi_1)\cap {\rm supp}(\psi_2)\neq \emptyset}\big(S(\psi_1,X)-S'(\psi_1,X)\big)^2  \Big]\\
         &\lesssim \mb{E}_{\mu^*_X}\Big[  \sum_{\psi_1\in \Psi_j^{D_Y}}\big(S(\psi_1,X)-S'(\psi_1,X)\big)^2  \Big],
    \end{aligned}
\end{equation*}
where the last inequality uses the fact that for any $\psi_1\in \Psi_j^{D_Y}$, there are only constant number of $\psi_2\in \Psi_j^{D_Y}$ so that ${\rm supp}(\psi_1)\cap {\rm supp}(\psi_2)\neq \emptyset$. Moreover,
\begin{equation*}
\begin{aligned}
   & \mb{E}_{\mu^*_X}\Big[ \Big(\sum_{\psi\in \Psi_j^{D_Y}}2^{-\frac{D_Yj}{2}}\cdot \big|S(\psi,X)-S'(\psi,X)\big|\Big)^2\Big]\\
 &\leq  \mb{E}_{\mu^*_X}\Big[ \sum_{\psi\in \Psi_j^{D_Y}}2^{-D_Yj}\cdot \sum_{\psi\in \Psi_j^{D_Y}} \big(S(\psi,X)-S'(\psi,X)\big)^2\Big]\\
   &\lesssim  \mb{E}_{\mu^*_X}\Big[  \sum_{\psi\in \Psi_j^{D_Y}} \big(S(\psi,X)-S'(\psi,X)\big)^2\Big].
\end{aligned}
\end{equation*}
Therefore, it holds for some constant $C$ that
\begin{equation*}
\begin{aligned}
     &\mb{E}_{\mu^*}\Big[\big(\ell(X,Y,S)-\ell(X,Y,S')\big)^2\Big]\leq C\, \mb{E}_{\mu^*_X}\Big[  \sum_{\psi\in \Psi_j^{D_Y}} \big(S(\psi,x)-S'(\psi,x)\big)^2\Big],\\
\end{aligned}
\end{equation*}
which  verifies the second assumption. Now we verify the last assumption. Note that for any $S,S'  \in \m S_j$, it holds that
\begin{equation*}
 \begin{aligned}
&d_n( S , S')\\
&=\sqrt{\frac{1}{n}\sum_{i=1}^n \Big(\sum_{\psi\in \Psi_j^{D_Y}}\big(S^2(\psi,X_i)-S'{}^2(\psi,X_i)\big)-2\psi(Y_i) \big(S(\psi,X_i)-S'(\psi,X_i)\big)\Big)^2}\\
&= \sqrt{\frac{1}{n}\sum_{i=1}^n \Big(\sum_{\psi\in \Psi_j^{D_Y}}\big(S(\psi,X_i)+S'(\psi,X_i)-2\psi(Y_i)\big)\cdot \big(S(\psi,X_i)-S'(\psi,X_i)\big)\Big)^2}\\
&\leq  \sqrt{\frac{1}{n}\sum_{i=1}^n  \sum_{\psi\in \Psi_j^{D_Y}}\big(S(\psi,X_i)+S'(\psi,X_i)-2\psi(Y_i)\big)^2\cdot  \sum_{\psi\in \Psi_j^{D_Y}}\big(S(\psi,X_i)-S'(\psi,X_i)\big)^2}\\
&\lesssim 2^{\frac{D_Y j}{2}}\sqrt{\frac{1}{n}\sum_{i=1}^n     \sum_{\psi\in \Psi_j^{D_Y}}\big(S(\psi,X_i)-S'(\psi,X_i)\big)^2},
 \end{aligned}
 \end{equation*}
where the last inequality uses that for any $(x,y)\in \m M$ and $S,S'\in \m S_j$, 
\begin{equation*}
    \begin{aligned}
      \sum_{\psi\in \Psi_j^{D_Y}}\big(S(\psi,x)+S'(\psi,x)-2\psi(y)\big)^2
       \lesssim | \Psi_j^{D_Y}|\cdot 2^{-D_Yj}+\sum_{\psi\in \Psi_j^{D_Y}} \psi(y)^2
        \lesssim 2^{D_Yj}.
    \end{aligned}
\end{equation*}
Furthermore, for any $\psi\in \Psi_j^{D_Y}$, $x\in \m M_X$, $S(\psi,x)=\frac{\sum_{i=1}^{W_j} \sum_{k\in \mb N_0^{D_X}, |k|<\alpha_X}a_{ik}^\psi (x-b_i)^k \rho(\frac{\|x-b_i\|}{\varepsilon^x_j})}{\sum_{i=1}^{W_j} \rho(\frac{\|x-b_i\|}{\varepsilon^x_j})+\frac{1}{n}}$ and $S'(\psi,x)=\frac{\sum_{i=1}^{W_j} \sum_{k\in \mb N_0^{D_X}, |k|<\alpha_X}a_{ik}^\psi{}' (x-b_i')^k \rho(\frac{\|x-b'_i\|}{\varepsilon^x_j})}{\sum_{i=1}^{W_j} \rho(\frac{\|x-b'_i\|}{\varepsilon^x_j})+\frac{1}{n}}$, it holds that
\begin{equation*}
    \begin{aligned}
        &|S(\psi,x)-S'(\psi,x)|\\
        &\leq \Big|\frac{\sum_{i=1}^{W_j} \sum_{k\in \mb N_0^{D_X}, |k|<\alpha_X}a_{ik}^\psi (x-b_i)^k \rho(\frac{\|x-b_i\|}{\varepsilon^x_j})}{\sum_{i=1}^{W_j} \rho(\frac{\|x-b_i\|}{\varepsilon^x_j})+\frac{1}{n}}-\frac{\sum_{i=1}^{W_j} \sum_{k\in \mb N_0^{D_X}, |k|<\alpha_X}a_{ik}^\psi{}' (x-b_i)^k \rho(\frac{\|x-b_i\|}{\varepsilon^x_j})}{\sum_{i=1}^{W_j} \rho(\frac{\|x-b_i\|}{\varepsilon^x_j})+\frac{1}{n}}\Big|\\
        &+\Big|\frac{\sum_{i=1}^{W_j} \sum_{k\in \mb N_0^{D_X}, |k|<\alpha_X}a_{ik}^\psi{}'(x-b_i)^k \rho(\frac{\|x-b_i\|}{\varepsilon^x_j})}{\sum_{i=1}^{W_j} \rho(\frac{\|x-b_i\|}{\varepsilon^x_j})+\frac{1}{n}}-\frac{\sum_{i=1}^{W_j} \sum_{k\in \mb N_0^{D_X}, |k|<\alpha_X}a_{ik}^\psi{}' (x-b_i')^k \rho(\frac{\|x-b_i'\|}{\varepsilon^x_j})}{\sum_{i=1}^{W_j} \rho(\frac{\|x-b_i\|}{\varepsilon^x_j})+\frac{1}{n}}\Big|\\
        &+\Big|\frac{\sum_{i=1}^{W_j} \sum_{k\in \mb N_0^{D_X}, |k|<\alpha_X}a_{ik}^\psi{}'(x-b_i')^k \rho(\frac{\|x-b_i'\|}{\varepsilon^x_j})}{\sum_{i=1}^{W_j} \rho(\frac{\|x-b_i\|}{\varepsilon^x_j})+\frac{1}{n}}-\frac{\sum_{i=1}^{W_j} \sum_{k\in \mb N_0^{D_X}, |k|<\alpha_X}a_{ik}^\psi{}' (x-b_i')^k \rho(\frac{\|x-b_i'\|}{\varepsilon^x_j})}{\sum_{i=1}^{W_j} \rho(\frac{\|x-b_i'\|}{\varepsilon^x_j})+\frac{1}{n}}\Big|\\
        &\lesssim \underset{i\in [W_j]}{\max}\sum_{k\in \mb N_0^{D_X}, |k|<\alpha_X}|a_{ik}^\psi-a_{ik}^\psi{}'|+2^{-\frac{D_Yj}{2}}\frac{n}{\varepsilon_j^x} \sum_{i=1}^{W_j}\|b_i-b_i'\|.
    \end{aligned}
\end{equation*}
Therefore, we have 
\begin{equation*}
    \begin{aligned}
        &\sqrt{\frac{1}{n}\sum_{i=1}^n     \sum_{\psi\in \Psi_j^{D_Y}}\big(S(\psi,X_i)-S'(\psi,X_i)\big)^2}\\
        &\lesssim \sqrt{\sum_{\psi\in \Psi_j^{D_Y}} \Big[\underset{i\in [W_j]}{\max}\big(\sum_{k\in \mb N_0^{D_X}, |k|<\alpha_X}|a_{ik}^\psi-a_{ik}^\psi{}'|\big)^2+2^{-D_Yj}\frac{n^2}{(\varepsilon_j^x)^2} \big(\sum_{i=1}^{W_j}\|b_i-b_i'\|\big)^2\Big]}\\
        &\lesssim \sqrt{\sum_{\psi\in \Psi_j^{D_Y}} \sum_{i\in [W_j]}  \sum_{k\in \mb N_0^{D_X}, |k|<\alpha_X}(a_{ik}^\psi-a_{ik}^\psi{}')^2+\frac{n^2W_j}{(\varepsilon_j^x)^2}  \sum_{i=1}^{W_j}\|b_i-b_i'\|^2}.
    \end{aligned}
\end{equation*}
Using the fact that the $\varepsilon$-covering number of a $d$-dimensional ball with radius $R$ is bounded by $(\frac{3R}{\varepsilon})^d$,  there exists a constant $C$ so that   for any $0<\varepsilon\leq \sup_{S,S'\in \m S_j} d_n(S,S')$,
 \begin{equation*}
    \log \mathbf{N}(\m S,d_n,\varepsilon)\leq C\, W_j 2^{jD_Y} \log \frac{n}{\varepsilon}.
 \end{equation*}
 which verifies the third assumption. The desired result is obtained by setting $W_n=C\,W_j 2^{jD_Y}$ and $T_n=n$ in Theorem~\ref{theoremjointregression},  and applying a union bound over $j\in [J]$.
 \section{Proof for Distribution Regression with Manifold Responses}\label{APP:B}
In the forthcoming analysis, let  $\m M_X$ denote the support of $\mu^*_X$,  and let $\m M_{Y|x}$ denote the support of $\mu^*_{Y|X}$.  We define $\m M=\{(x,y):\, x\in \m M_X, y\in \m M_{Y|x} \}$ as the support of the joint distribution $\mu^*=\mu^*_X\mu^*_{Y|X}$. Let $u^*(\,\cdot\,|\,x)$ represent the density function of $\mu^*_{Y|x}$ with respect to the volume measure of $\m M_{Y|x}$. Moreover, $\m M_Y=\bigcup_{x\in \m M_X} \m M_{Y|x}$ is the support of the marginal distribution of $Y$. 

 \medskip
 We will also  refer to  the notations from  the definition of the $(\beta_Y,\beta_X)$-smooth submanifold family as outlined in Definition~\ref{def:manifold} in the main text, and provide a recapitulation here: for any $w_0=(x_0,y_0)\in \m M$, there exists a neighborhood $U_{\omega_0}$ of $y_0$ on $\m M_Y$, so that for any $x\in \mb B_{\m M_X}(x_0,\tau)$, the function $\operatorname{Proj}_{T_{y_0} \mathcal{M}_{Y|x_0}}\left(y-y_0\right): \mathcal{M}_Y \rightarrow T_{y_0} \mathcal{M}_{Y|x_0}$, when {restricted to $U_{\omega_0}\cap \m M_{Y|x}$}, is a   {diffeomorphism} with  inverse function $\phi_{\omega_0,x}(\cdot)$ defined on $\mb B_{T_{\m M_{Y|x_0}}y_0}(0,\tau_1)$. Moreover, the function $\Phi_{\omega_0}: \mb B_{T_{\m M_{Y|x_0}}y_0}(0,\tau_1)\times  \mb B_{\m M_X}(x_0,\tau)\to \mb R^{D_Y}$ define as {$\Phi_{\omega_0}(z,x)=\phi_{\omega_0,x}(z)$ belongs to ${\m H}^{\beta_Y,
 \beta_X}_{L,D_Y}(\mb B_{T_{\m M_{Y|x_0}}y_0}(0,\tau_1), \mb B_{\m M_X}(x_0,\tau))$}.

  \medskip
 For any point $w_0\in \m M$, the terms $U_{w_0}$, $\Phi_{w_0}$ will be used to denote the neighborhood and function described above respectively. In the scenario where the response space remains invariant across different covariates (referred to as Regime 2), we have $\m  M_{Y|x}=\m M_Y$ for all $x\in \m M_X$. Consequently,  $\Phi_{\omega_0}(z,x)$ is independent of $x$, allowing us to simplify the notation to $\Phi_{\omega_0}(z)=\Phi_{\omega_0}(z,x)$, and we have $\Phi_{\omega_0}(z)\in \m H_{L,D_Y}^{\beta_Y}(\mb B_{T_{\m M_{Y|x_0}}y_0}(0,\tau_1))$.

\subsection{Proof of Theorem~\ref{th:thgamma>1}}\label{proofth:thgamma>1}
 We consider the estimator defined in Appendix~\ref{R2gamma>1}. For any $j\in \{0\}\cup [J]$ with $J=\lceil \frac{1}{d_Y}\cdot \log_2 (\frac{n}{\log n})\rceil$, the following lemma  provides a bound for the mean squared error between $\wh S^{\dagger}_j(\psi,x)$ and $u_{\psi}^*(x)=\mb{E}_{\mu^*_{Y|x}}[2^{\frac{j(d_Y-D_Y)}{2}}\psi(y)]$.
\begin{lemma}\label{lemmaregime2MSE}
     Suppose $\mu^*\in \m P_2$ and  with the choices of $\m S_j^\dagger$ defined in~\eqref{defSdaggerr2}, it holds with probability at least $1-\frac{1}{n^2}$ that for any $j\in [J]$,
      \begin{equation*}
          \mb{E}_{\mu^*_X} \bigg[\sum_{\psi\in \Psi_j^{D_Y}} (\wh S_j^\dagger(\psi,x)-u^*_{\psi}(x))^2\bigg]\lesssim 2^{\frac{2j\alpha_Xd_Y}{2\alpha_X+d_X}} (\frac{n}{\log n})^{-\frac{2\alpha_X}{2\alpha_X+d_X}}.
      \end{equation*}
\end{lemma}

\noindent The proof of Lemma~\ref{lemmaregime2MSE} is given in Appendix~\ref{proof:lemmaregime2MSE}. Then let $\eta=\frac{d_Y\alpha_X}{2\alpha_X+d_X}$.    Utilizing the property that for any function $f\in \mb H^{\eta}_1(\mb R^{D_Y})$ and $\psi\in \Psi_j^{D_Y}$, it holds that $|f_{\psi}|=|\int_{\mb R^{D_Y}} f(y)\psi(y)\,\dd y|\lesssim 2^{-j\eta-jD_y}$,  we can deduce that with probability at least $1-\frac{1}{n^2}$,
           \begin{equation*}
     \begin{aligned}
          & \mb{E}_{\mu^*_X}  \Big[\underset{f\in \m H_1^{\eta}(\mb R^{D_Y})}{\sup} \big|\mb{E}_{\mu^*_{Y|X}} f(y)- \sum_{j=0}^J\sum_{\psi\in\Psi_j^{D_Y}} f_{\psi}2^{\frac{j(D_Y-d_Y)}{2}} \wh S_j^\dagger(\psi,X)\big|\Big]\\
          &\leq \mb{E}_{\mu^*_X}\Big[\underset{f\in  \m H_1^{\eta}(\mb R^{D_Y})}{\sup}\big|\sum_{j=0}^J\sum_{\psi\in \Psi_j^{D_Y}} f_{\psi}\mb{E}_{\mu^*_{Y|X}}[\psi(y)]-\sum_{j=0}^J\sum_{\psi\in\Psi_j^{D_Y}} f_{\psi}2^{\frac{j(D_Y-d_Y)}{2}} \wh S_j^\dagger(\psi,X)\big|\Big]\\
          &+ \underset{f\in  \m H_1^{\eta}(\mb R^{D_y})}{\sup}  \underset{y\in \m M_Y}{\sup} \Big|\sum_{j=J+1}^{\infty}\sum_{\psi\in \Psi_j^{D_Y}} f_{\psi} \psi(y)\Big|\\
            &\lesssim \sum_{j=0}^{J}\sqrt{\sum_{\psi\in \Psi_j^{D_Y}}2^{-2j\eta}\mb{E}_{\mu^*_X}\Big[ \big(\mb{E}_{\mu^*_{Y|x}}[ 2^{\frac{j(d_Y-D_Y)}{2}} \psi(y)]-\wh S_j^\dagger(\psi,x)\big)^2\Big]}+2^{-J\eta}\\
     &\lesssim (\log n)\cdot (\frac{n}{\log n})^{-\frac{\alpha_X}{2\alpha_X+d_X}}.
           \end{aligned}
 \end{equation*}
 So for any $\gamma\geq \eta=\frac{d_Y\alpha_X}{2\alpha_X+d_X}$, it holds that 
 \begin{equation*}
     \begin{aligned}
     & \mb{E}_{\mu^*_X}  \Big[\underset{f\in \m H_1^{\gamma}(\mb R^{D_Y})}{\sup} \big|\mb{E}_{\mu^*_{Y|X}} f(y)- \sum_{j=0}^J\sum_{\psi\in\Psi_j^{D_Y}} f_{\psi}2^{\frac{j(D_Y-d_Y)}{2}} \wh S_j^\dagger(\psi,X)\big|\Big]\\
                & \leq \mb{E}_{\mu^*_X}  \Big[\underset{f\in \m H_1^{\eta}(\mb R^{D_Y})}{\sup} \big|\mb{E}_{\mu^*_{Y|X}} f(y)- \sum_{j=0}^J\sum_{\psi\in\Psi_j^{D_Y}} f_{\psi}2^{\frac{j(D_Y-d_Y)}{2}} \wh S_j^\dagger(\psi,X)\big|\Big]\\
        &\lesssim (\log n)\cdot (\frac{n}{\log n})^{-\frac{\alpha_X}{2\alpha_X+d_X}}.
     \end{aligned}
 \end{equation*}

\subsection{Proof of Theorem~\ref{thgamma<1}}\label{proofthgamma<1}

We first derive the following results concerning the population-level reconstruction error for the first step of manifold recovery, the proof of which is given in Appendix~\ref{proof:lemma:manifoldlearningerrorr2}.
 
\begin{lemma}\label{lemma:manifoldlearningerrorr2}
Suppose $\mu^*\in \m P^*_2$ and  with the choices of $\m G$ defined in~\eqref{defGr2}, there exist positive constants $C,C_1$ so that it holds with probability at least $1-\frac{1}{n^2}$ that 
 \begin{enumerate}
     \item  
 For any $k\in \wh{\m K}$ and $\gamma_1\in (0,1]$,
    \begin{equation*}
    \begin{aligned}
           &\mb{E}_{\mu^*_X}\mb{E}_{\mu^*_{Y|X}}[\|Y-\wh G_{[k]}(\wh Q_{[k]}(Y))\|^{\gamma_1}\cdot\mathbf{1}(X\in \mb B_{\mb R^{D_X}}(x_k,2\tau_2))\mathbf{1}(Y\in \mb B_{\mb R^{D_Y}}(y_k,2\tau_2))]\\
          & \lesssim  \left\{
          \begin{array}{cc}
           C \, \frac{(\log n)^{1+\gamma_1}}{\sqrt{n}}   & \frac{d_Y}{\beta_Y} \leq 2\gamma_1, \\
C \,(\log n \wedge \frac{1}{d_Y-2\gamma_1\beta_Y)} )^{1+\gamma_1}\cdot n^{-\frac{\gamma_1}{ \frac{d_Y}{\beta_Y}}}   & \frac{d_Y}{\beta_Y}>2\gamma_1.
          \end{array}
        \right.
        \end{aligned}
        \end{equation*}
\item  For any $k\in \wh {\m K}$, there exists $(x^*_k,y^*_k)\in \mb B_{\m M}((x_k,y_k),\sqrt{2}\tau_2)$  such that 
\begin{equation*}
    \wh V_\sk^T P^*_\sk\wh V_\sk\gtrsim C_1I_{d_Y},
\end{equation*}
where $\m P^*_\sk$ is the projection matrix of $T_{\m M_Y} y^*_k$.
\end{enumerate}

\end{lemma}

Given the assumption that $\m M_{Y|x}=\m M_Y$ for any $x\in \m M_x$,  and note that if a function $f(y,x)$ is $\m H^{\beta_Y,\beta_X}$-smooth for some $\beta_X>0$,  and if  $f(y,x)$ is independent of $x$, then  $f(y)=f(y,x)$ must inherently be $\m H^{\beta_Y}$-smooth. Conversely, if $f(y)$ is $\m H^{\beta_Y}$-smooth, defining $f(y,x)=f(y)$ will result in a function being $\m H^{\beta_Y,\beta_X}$-smooth for any $\beta_X>0$. Consequently, we can use Lemma~\ref{le:projection} from Appendix~\ref{app: regularitymanifold} to obtain the invertibility of $\wh V_\sk^T(\cdot-y_k)$. Specifically,  when $\tau_2$ is small enough,   given the second statement in Lemma~\ref{lemma:manifoldlearningerrorr2}, for any $k\in \wh {\m K}$, there exists a subset $\wh U^\sk_{Y}$ so that $\mb B_{\m M_{Y}}(y_k^*,3\tau_2)\subset \wh U^\sk_{Y}\subset \m M_{Y}$, and the function $\wh Q_\sk(\cdot)=\wh V_\sk^T(\cdot-y_k)$, when restricted to domain $\wh U^\sk_{Y}$, is a diffeomorphism  that maps $\wh U^\sk_{Y}$ to $\mb B_{\mb R^{d_Y}}(\wh V_\sk^T(y_k^*-y_k),3\tau_2)$ with inverse denoted as $[\wh Q_\sk(\cdot)]^{-1}$. The function $\wh G^{\dagger}_\sk: \mb B_{\mb R^{d_Y}}(\wh V_\sk^T(y_k^*-y_k),3\tau_2)\to \mb R^{D_Y}$ defined as  $\wh G^{\dagger}_\sk(z)=[\wh Q_\sk(\cdot)]^{-1}(z)$ belongs to $ \m H^{\beta_Y}_{L_1,D_Y}(\mb B_{\mb R^{d_Y}}(\wh V_\sk^T(y_k^*-y_k),3\tau_2))$ for some constant $L_1$. Based on this fact, the push forward measure $\wh Q_{\sk\#}(\mu^*_{Y|x}|_{\wh U^\sk_{Y}})$, has a density 
$\wh \nu_{\sk}(z|x)=u^*(\wh G^\dagger_\sk(z)|x)\cdot \sqrt{{\rm det}(J_{\wh G^\dagger_\sk}(z)^TJ_{\wh G^\dagger_\sk}(z))}$ for $z\in \mb B_{\mb R^{d_Y}}(\wh V_\sk^T(y_k^*-y_k),3\tau_2)$,  where $u^*(\,\cdot\,|\,x)$ is the density of $\mu^*_{Y|x}$ with respect to the volume measure of $\m M_Y$. Since $\beta_Y\geq \alpha_Y+1$,  there exists a constant $L_2$ so that $\wh \nu_{\sk}(z,|,x)\in  \ov{\m H}^{\alpha_Y,\alpha_X}_{L_2}(\mb B_{\mb R^{d_Y}}(\wh V_\sk^T(y_k^*-y_k),3\tau_2),\m M_X)$.
Furthermore, for any $j\in\mb N$ and $\psi\in \Psi_j^{d_Y}$, we have 
\begin{equation*}
\begin{aligned}
     & \mb{E}_{\mu^*_{Y|x}}[\psi(\wh Q_{\sk}(Y))\rho_\sk(x,Y)]\\
     &=  \mb{E}_{\mu^*_{Y|x}}[\psi(\wh Q_{\sk}(Y))\rho_\sk(x,Y)\mathbf{1}(Y\in \wh U^\sk_{Y})\mathbf{1}(x\in \mb B_{\m M_x}(x_k^*,2\tau_2))]\\
    & =  \mb{E}_{\mu^*_{Y|x}}[\psi(\wh Q_{\sk}(Y))\rho_\sk(x,\wh G_\sk^{\dagger}(\wh Q_{\sk}(Y)))\mathbf{1}(Y\in \wh U^\sk_{Y})]\\
    &=\int_{\mb B_{\mb R^{d_Y}}(\wh V_\sk^T(y_k^*-y_k),3\tau_2)} \psi(z) \rho_\sk(x, \wh G_\sk^{\dagger}(z))\wh v_\sk (z|x)\,\dd z.
\end{aligned}
\end{equation*}
Let $\ov \nu_{\sk}(z,|,x)\in  \ov{\m H}^{\alpha_Y,\alpha_X}_{L_2}(\mb R^{d_Y},\mb R^{D_X})$ be a smooth extension of $\wh \nu_{\sk}(z,|,x)$ to $\mb R^{d_Y}\times \mb R^{D_X}$. Define 
\begin{equation*}
    \wt v_\sk (z,x)=\left\{
    \begin{array}{cc}
     \rho_\sk(x, \wh G_\sk^{\dagger}(z))\ov \nu_\sk (z|x),    &  \text{ if }z\in \mb B_{\mb R^{d_Y}}(\wh V_\sk^T(y_k^*-y_k),3\tau_2), x\in \mb B_{\mb R^{D_X}}(x_k^*,3\tau_2) \\
      0   & \text{otherwise}.
    \end{array}
    \right.
\end{equation*}
We can verify that $ \wt v_\sk (z,x)\in  \ov{\m H}^{\alpha_Y,\alpha_X}_{L_3}(\mb R^{d_Y},\mb R^{D_X})$ with a constant $L_3$. Therefore, for any $j\in \mb N$, $\psi\in \Psi_j^{d_Y}$ and $x\in \m M_X$, it holds that
\begin{equation*}
     2^{\frac{d_Yj}{2}} \mb{E}_{\mu^*_{Y|x}}[\psi(\wh Q_{\sk}(Y))\rho_\sk(x,Y)]=2^{\frac{d_Yj}{2}}\int_{\mb R^{d_Y}} \psi(z) \wt v_\sk (z|x)\,\dd z,
\end{equation*}
and $2^{\frac{d_Yj}{2}}\int_{\mb R^{d_Y}} \psi(z) \wt v_\sk (z|\cdot)\,\dd z \in \m H^{\alpha_X}_{L_4}(\mb R^{D_X})$ for some constant $L_4$. Moreover, for any $x\in \m M_X$, given that  $\wt v_\sk (\cdot|x)\in \m H^{\alpha_Y}_{L_3}(\mb R^{d_Y})$, it follows  that for any $x\in \m M_X$,
\begin{equation*}
      \Big|\mb{E}_{\mu^*_{Y|x}}[\psi(\wh Q_{\sk}(Y))\rho_\sk(x,Y)]\Big|=\Big|\int_{\mb R^{d_Y}} \psi(z) \wt v_\sk (z|x)\,\dd z\Big|\lesssim 2^{-\frac{d_Yj}{2}-j\alpha_Y}.
\end{equation*}
Let $J=\lceil {\frac{1}{2\alpha_Y+d_Y+d_X\frac{\alpha_Y}{\alpha_X}}}\cdot \log_2 (\frac{n}{\log n})\rceil$. For $j\in \{0\}\cup [J]$, denote
\begin{equation*}
    \m S_j^{\ddagger}=\{S: \Psi_j^{d_Y}\times \mb R^{D_X}\to \mb R:\, S(\psi,x)=\sum_{\psi_1\in \Psi_j^{d_Y}} s_{\psi_1}(x),\text{ where } s_{\psi_1}\in \ms S_j \text{ for each }\psi_1\in \Psi_j^{d_Y}\},
\end{equation*}
where $\ms S_j$ is defined in~\eqref{defMSddaggerr2}. Using the independence of $\{X_i\}_{i\in I_1}$ and $\{X_i\}_{i\in I_2}$, and mirroring the analysis from the proof of Lemma~\ref{lemma1.1}---where we replace $D_Y$ with $d_Y$, and modify $\psi(Y)$ to $\psi(\wh Q_{\sk}(Y))\rho_\sk(X,Y)$. To apply Theorem~\ref{theoremjointregression}, we set $\{\psi_{\lambda}((X,Y))\}_{\lambda\in \Lambda}=\{\psi(\wh Q_{\sk}(Y))\rho_\sk(X,Y):\, \psi\in \Psi_{j}^{d_Y}\}$, where the response variable $Y$ is redefined as the joint vector of $(X,Y)$, alongside $\m S= \m S_j^{\ddagger}$ --- we can show that, by applying a union argument over $j\in [J]$ and $k\in\wh {\m K}$, it holds with probability at least $1-\frac{1}{n^2}$ that for any $j\in [J]$ and $k\in\wh {\m K}$,
      \begin{equation}\label{eqnregime2gamma<1}
          \mb{E}_{\mu^*_X} \Big[\sum_{\psi\in \Psi_j^{d_Y}} (\wh v_{k\psi}(X)-  \mb{E}_{\mu^*_{Y|X}}[\psi(\wh Q_{\sk}(Y))\rho_\sk(X,Y)])^2\Big]\lesssim 2^{\frac{2j\alpha_Xd_Y}{2\alpha_X+d_X}} (\frac{n}{\log n})^{-\frac{2\alpha_X}{2\alpha_X+d_X}}.
      \end{equation}
Furthermore, recall that 
     \begin{equation*}
     \wh{\m K}=\{k\in [K]:\, \exists i\in I_1, \|(X_i,Y_i)-(x_k,y_k)\|\leq \sqrt{2}\tau_2\}.
 \end{equation*}
 So for any $k\in [K]\setminus\wh {\m K}$, it holds that 
 \begin{equation*}
     \frac{1}{n}\sum_{i\in I_1}\rho_\sk(X_i,Y_i)\leq   \frac{1}{n}\sum_{i\in I_1}\mathbf{1}(\|(X_i,Y_i)-(x_k,y_k)\|\leq \sqrt{2}\tau_2)=0,
 \end{equation*}
and by Bernstein's inequality, it holds with probability at least $1-\frac{1}{n^2}$ that for any $k\in [K]\setminus\wh {\m K}$,
\begin{equation*}
    \mb{E}_{\mu^*}[\rho_\sk(X,Y)]\lesssim\sqrt{\frac{\log n}{n}}.
\end{equation*}
Denote $\wh \mu_{Y|x}= \sum_{k\in \wh{\m K}} \wh G_{[k]}(\cdot,x)_{\#}\wh\nu_\sk(\cdot|x)$, it holds with probability at least $1-\frac{1}{n^2}$ that 
\begin{equation*}
    \begin{aligned}
            & \mb{E}_{\mu^*_X} \Big[\underset{f\in \m H_1^{\gamma}(\mb R^{D_Y})}{\sup}\big|\int  f(y)\dd \mu^*_{Y|X} -\int  f(y)\dd \wh\mu_{Y|X}\big| \Big]\\
         & =\mb{E}_{\mu^*_X} \Big[\underset{f\in \m H_1^{\gamma}(\mb R^{D_Y})}{\sup}\int f(y)\dd \mu^*_{Y|X} -\int f(y)\dd \wh\mu_{Y|X} \Big]\\
        &= \mb{E}_{\mu^*_X} \Big[\underset{f\in \m H_1^{\gamma}(\mb R^{D_Y})}{\sup}\int \sum_{k=1}^K f(y)\rho_\sk(X,y)\dd \mu^*_{Y|X} -\int f(y)\dd \wh\mu_{Y|X}  \Big]\\
          &\leq\mb{E}_{\mu^*_X} \Big[\underset{f\in \m H_1^{\gamma}(\mb R^{D_Y})}{\sup}\int \sum_{k\in \wh{\m K}}f(y)\rho_\sk(X,y)\dd \mu^*_{Y|X}-\int f(y)\dd \wh\mu_{Y|X}\Big]\\
          &\qquad+\mb{E}_{\mu^*_X} \Big[\underset{f\in \m H_1^{\gamma}(\mb R^{D_Y})}{\sup}\int \sum_{k\in [K]\setminus\wh{\m K}}f(y) \rho_\sk(X,y)\dd \mu^*_{Y|X}\Big]\\
            &\leq\mb{E}_{\mu^*_X} \Big[\underset{f\in \m H_1^{\gamma}(\mb R^{D_Y})}{\sup}\int \sum_{k\in \wh{\m K}}f(y) \rho_\sk(X,y)\dd \mu^*_{Y|X}-\int f(y)\dd \wh\mu_{Y|X}  \Big]  +\mb{E}_{\mu^*} \Big[  \sum_{k\in [K]\setminus\wh{\m K}}\rho_\sk(X,Y) \Big]\\
              &\leq\mb{E}_{\mu^*_X} \Big[\underset{f\in \m H_1^{\gamma}(\mb R^{D_Y})}{\sup}\int \sum_{k\in \wh{\m K}}f(y) \rho_\sk(X,y)\dd \mu^*_{Y|X}-\int f(y)\dd \wh\mu_{Y|X}  \Big]  + C\,\sqrt{\frac{\log n}{n}}.\\
    \end{aligned}
\end{equation*}
To simplify the notation,  for any $j>J$, $k\in \wh {\m K}$ and $\psi\in \Psi_j^{d_Y}=\{\psi\in \ov \Psi_j^{d_Y}:\, \text{supp}(\psi)\cap \mb B_{\mb R^{d_Y}}(\mathbf{0},2\tau_2)\neq \emptyset\}$, we set $\wh v_{k\psi}(\cdot)\equiv 0$. Then denote $f_{\psi}=\int f(y)\psi(y)\,\dd y$,  it holds with probability at least $1-\frac{3}{n^2}$ that for any $\gamma\in (0,1]$,
 
  \begin{equation*}
  \begin{aligned}
     &\mb{E}_{\mu^*_X} \Big[\underset{f\in \m H_1^{\gamma}(\mb R^{D_Y})}{\sup}\int \sum_{k\in \wh{\m K}}f(y) \rho_\sk(X,y)\dd\mu^*_{Y|X}-\int f(y)\dd \wh\mu_{Y|X}  \Big]\\
    &\leq\mb{E}_{\mu^*_X} \Big[\underset{f\in \m H_1^{\gamma}(\mb R^{D_Y})}{\sup} \sum_{k\in \wh{\m K}} \int \left(f(y)-f\big(\wh G_\sk(\wh Q_\sk(y))\big)\right)\rho_\sk(X,y)\dd\mu^*_{Y|X}\Big]\\
&+\mb{E}_{\mu^*_X} \Big[\underset{f\in \m H_1^{\gamma}(\mb R^{D_Y})}{\sup} \sum_{k\in \wh{\m K}} \int f(\wh G_\sk(\wh Q_\sk(y)))\rho_\sk(X,y)\dd\mu^*_{Y|X}-\sum_{k\in \wh{\m K}} \int f(\wh G_{\sk}(z)) \sum_{j=0}^{\infty} \sum_{\psi\in \Psi_j^{d_Y}}\psi(z)\wh v_{k\psi}(X)\,\dd z\Big]\\
&\lesssim \sum_{k\in \wh{\m K}} \mb{E}_{\mu^*}\big[\|Y-\wh G_\sk(\wh Q_\sk(Y))\|^{\gamma}\rho_\sk(X,Y)\big]\\
&\qquad+\mb{E}_{\mu^*_X} \Big[ \sum_{k\in \wh{\m K}}  \underset{f\in \m H_1^{\gamma}(\mb R^{d_Y})}{\sup}\int f(\wh Q_\sk(y))\rho_\sk(X,y)\dd\mu^*_{Y|X}-\int f(z) \sum_{j=0}^{\infty} \sum_{\psi\in \Psi_j^{d_Y}}\psi(z) \wh v_{k\psi}(X)\,\dd z\Big]\\
&= \sum_{k\in \wh{\m K}} \mb{E}_{\mu^*}\big[\|Y-\wh G_\sk(\wh Q_\sk(Y))\|^{\gamma}\rho_\sk(X,Y)\big]\\
&\qquad+\mb{E}_{\mu^*_X} \Big[ \sum_{k\in \wh{\m K}}  \underset{f\in \m H_1^{\gamma}(\mb R^{d_Y})}{\sup}\int \sum_{j=0}^{\infty}\sum_{\psi\in \Psi_{j}^{d_Y}} f_{\psi}\cdot
\psi(\wh Q_\sk(y))\rho_\sk(X,y)\dd\mu^*_{Y|X}- \sum_{j=0}^{\infty}\sum_{\psi\in \Psi_{j}^{d_Y}}  f_{\psi}\wh v_{k\psi}(X)\Big]\\
&\leq \underbrace{\sum_{k\in \wh{\m K}}  \mb{E}_{\mu^*}\big[\|Y-\wh G_{[k]}(\wh Q_{[k]}(Y))\|^{\gamma}\cdot\mathbf{1}(X\in \mb B_{\mb R^{D_X}}(x_k,2\tau_2))\mathbf{1}(Y\in \mb B_{\mb R^{D_Y}}(y_k,2\tau_2))\big]}_{ (E_A)} \\
&\qquad+\underbrace{\mb{E}_{\mu^*_X} \Big[ \sum_{k\in \wh{\m K}}  \underset{f\in \m H_1^{\gamma}(\mb R^{d_Y})}{\sup}\sum_{j=0}^{J} \sum_{\psi\in \Psi_j^{d_Y}} f_{\psi}\cdot\Big( \mb{E}_{\mu^*_{Y|X}}[\psi(\wh Q_\sk(Y))\rho_\sk(X,Y)]- \wh v_{k\psi}(X) \Big)\Big]}_{ (E_B)}\\
&\qquad+\underbrace{\mb{E}_{\mu^*_X} \Big[ \sum_{k\in \wh{\m K}}  \underset{f\in \m H_1^{\gamma}(\mb R^{d_Y})}{\sup}\sum_{j=J+1}^{\infty} \sum_{\psi\in \Psi_j^{d_Y}} f_{\psi}\cdot \mb{E}_{\mu^*_{Y|X}}[\psi(\wh Q_\sk(Y))\rho_\sk(X,Y)] \Big]
}_{(E_C)}.
  \end{aligned}
  \end{equation*}
  By Lemma~\ref{lemma:manifoldlearningerrorr2}, we have 
  \begin{equation*}
       (E_A)\lesssim  (\log n \wedge \frac{\mathbf{1}(d_Y/\beta_Y>2\gamma)}{\beta_Yd_Y-2\gamma\beta_Y)} )^2\cdot n^{-\frac{\gamma}{\frac{d_Y}{\beta_Y}}}+\frac{(\log n)^2}{\sqrt{n}}.
  \end{equation*}
  Moreover,  since $|f_{\psi}|\lesssim 2^{-j\gamma-jd_Y/2}$ for $\psi\in \Psi_j^{d_Y}$, we have
  \begin{equation*}
      \begin{aligned}
           (E_B)&\lesssim \sum_{k\in \wh{\m K}}\sum_{j=0}^{J}\sum_{\psi\in \Psi_j^{d_Y}}2^{-j\gamma-\frac{jd_Y}{2}}\sqrt{\mb{E}_{\mu^*_X}\Big[ \Big(\mb{E}_{\mu^*_{Y|X}}[\psi(\wh Q_\sk(Y))\rho_\sk(X,Y)]- \wh v_{k\psi}(X)\Big)^2\Big]}\\
          &\lesssim \sum_{k\in \wh{\m K}}\sum_{j=0}^{J}\sqrt{\sum_{\psi\in \Psi_j^{d_Y}}2^{-2j\gamma}\mb{E}_{\mu^*_X}\Big[ \Big(\mb{E}_{\mu^*_{Y|X}}[\psi(\wh Q_\sk(Y))\rho_\sk(X,Y)]- \wh v_{k\psi}(X)\Big)^2\Big]}\\
          &\lesssim (\log n)\cdot (\frac{n}{\log n})^{-\frac{\alpha_X}{2\alpha_X+d_X}}+C_1\,(\frac{n}{\log n})^{-\frac{\alpha_Y+\gamma}{2\alpha_Y+d_Y+\frac{\alpha_Y}{\alpha_X}d_X}},
      \end{aligned}
  \end{equation*}
  and 
 \begin{equation*}
     \begin{aligned}
                 (E_C) & \lesssim \sum_{k\in \wh{\m K}}\sum_{j=J+1}^{\infty} \sum_{\psi\in \Psi_j^{d_Y}} 2^{-j(\gamma+\alpha_Y)-jd_Y}  \\
     &  \lesssim(\frac{n}{\log n})^{-\frac{\alpha_Y+\gamma}{2\alpha_Y+d_Y+\frac{\alpha_Y}{\alpha_X}d_X}} \\
     \end{aligned}
 \end{equation*}
 Finally, we have 
 \begin{equation*}
     \begin{aligned}
      & \mb{E}_{\mu^*_X} \Big[\underset{f\in \m H_1^{\gamma}(\mb R^{D_Y})}{\sup}\big|\int f(y)\dd \mu^*_{Y|X} -\int f(y)\dd \wh\mu_{Y|X}\big| \Big]\\
      &\lesssim \sqrt{\frac{\log n}{n}}+ (E_A)+ (E_B)+(E_C)\\
      &  \lesssim (\log n)^2\cdot n^{-\frac{\alpha_X}{2\alpha_X+d_X}}+(\frac{n}{\log n})^{-\frac{\alpha_Y+\gamma}{2\alpha_Y+d_Y+\frac{\alpha_Y}{\alpha_X}d_X}}
       +n^{-\frac{\gamma}{\frac{d_Y}{\beta_Y}}}.
     \end{aligned}
 \end{equation*}

\subsection{Proof of Theorem~\ref{th:combined} (minimax upper bound for Regime 2 and 3b)}
\subsubsection{Proof  for Regime 2}\label{proofth:thgamma>0}
We consider the estimator detailed in Appendix~\ref{R2gamma>0}. For
\begin{equation*}
\begin{aligned}
       & \wh{\m J}(f,x)= {\sum_{j=0}^J\sum_{\psi\in\Psi_j^{D_Y}} f_{\psi}2^{\frac{j(D_Y-d_Y)}{2}} \wh S_j^\dagger(\psi,x) }+ {\sum_{k\in \wh{\m K}} \int_{\mb R^{d_Y}} f^{\perp}_J(\wh G_\sk (z)) \sum_{j=0}^J \sum_{\psi\in \Psi_j^{d_Y}}\psi(z) \wh v_{k\psi}(x)\,\dd z},\\
       &\qquad f_{\psi}=\int_{\mb R^{D_Y}} f(y)\psi(y)\,\dd y,\quad f^{\perp}_J(y)=f(y)-\sum_{j=0}^J \sum_{\psi\in  \ov\Psi_j^{D_Y}} f_{\psi} \psi(y),
\end{aligned}
\end{equation*}
where $J=\lceil {\frac{1}{2\alpha_Y+d_Y+d_X\frac{\alpha_Y}{\alpha_X}}}\cdot \log_2 (\frac{n}{\log n})\rceil$.
 We can get 
  \begin{equation*}
     \begin{aligned}
          & \mb{E}_{\mu^*_X}  \Big[\underset{f\in \m H_1^{\gamma}(\mb R^{D_Y})}{\sup} \big|\mb{E}_{\mu^*_{Y|X}} [f(Y)]- \wh{\m J}(f,X)\big|\Big]\\
           &\leq   \underbrace{ \mb{E}_{\mu^*_X} \Big[ \underset{f\in \m H_1^{\gamma}(\mb R^{D_Y})}{\sup} \big| \sum_{j=0}^J \sum_{\psi\in \Psi_j^{D_Y}} f_{\psi} \big(\mb{E}_{\mu^*_{Y|X}}[\psi(Y)]- 2^{\frac{j(D_Y-d_Y)}{2}} \wh S_j^\dagger(\psi,X)\big)\big|\Big]}_{ (E_A)}  \\
          &+ \underbrace{\mb{E}_{\mu^*_X}  \Big[\underset{f\in \m H_1^{\gamma}(\mb R^{D_Y})}{\sup} \Big|\mb{E}_{\mu^*_{Y|X}}\big[\sum_{k\in [K]\setminus\wh{\m K}}\rho_\sk (X,Y)f_J^{\perp}(Y)+\sum_{k\in \wh{\m K}}\rho_\sk(X,Y)\big(f_J^{\perp}(Y)-f_J^{\perp}(\wh G_\sk(\wh Q_\sk(Y)))\big)\big]\Big|\Big]}_{ (E_B)} \\
          &+ \underbrace{\mb{E}_{\mu^*_X}  \Big[\underset{f\in \m H_1^{\gamma}(\mb R^{D_Y})}{\sup} \Big|\sum_{k\in \wh{\m K}} \mb{E}_{\mu^*_{Y|X}}\big[f_J^{\perp}(\wh G_\sk(\wh Q_\sk(Y)))\rho_\sk(X,Y)\big]-\int_{\mb R^{d_Y}} f^{\perp}_J(\wh G_\sk (z)) \sum_{j=0}^J \sum_{\psi\in \Psi_j^{d_Y}}\psi(z) \wh v_{k\psi}(X)\,\dd z\Big|\Big]}_{(E_C)}.\\
           \end{aligned}
 \end{equation*}
We first bound term $ (E_A)$, notice that 
\begin{equation*}
    \begin{aligned}
         (E_A)\ &\leq \  C\, \sum_{j=0}^{J}\sum_{\psi\in  \Psi_j^{D_Y}}2^{-j\gamma-\frac{jD_Y}{2}}\sqrt{\mb{E}_{\mu^*_X}\Big[ \Big(\mb{E}_{\mu^*_{Y|X}}\big[\, 2^{\frac{j(d_Y-D_Y)}{2}} \psi(Y)\,\big]-\wh S_j^\dagger(\psi,X)\Big)^2\Big]}\\
        & \leq  C\,   \sum_{j=0}^{J} \sqrt{\sum_{\psi\in \Psi_j^{D_Y}} 2^{-jD_Y}}\cdot   2^{-j\gamma}\sqrt{\sum_{\psi\in \Psi_j^{D_Y}}\mb{E}_{\mu^*_X}\Big[ \Big(\mb{E}_{\mu^*_{Y|X}}\big[\, 2^{\frac{j(d_Y-D_Y)}{2}} \psi(Y)\,\big]-\wh S_j^\dagger(\psi,X)\Big)^2\Big]}\\
     & \leq \ C_1\,  \sum_{j=0}^{J}2^{-j\gamma}\sqrt{\sum_{\psi\in \Psi_j^{D_Y}}\mb{E}_{\mu^*_X}\Big[ \Big(\mb{E}_{\mu^*_{Y|X}}\big[\, 2^{\frac{j(d_Y-D_Y)}{2}} \psi(Y)\,\big]-\wh S_j^\dagger(\psi,X)\Big)^2\Big]},
    \end{aligned}
\end{equation*}
where the first inequalities uses $f\in \m H^{\gamma}_1(\mb R^{D_Y})$, which implies $|f_{\psi}|\lesssim 2^{-j\gamma-\frac{jD_Y}{2}}$, alongside Jensen's inequality;  the second inequality is derived using the Cauchy-Schwarz inequality, while the final inequality uses the fact that $|\Psi_j^{D_Y}|=\m O(2^{D_Yj})$.
We then bound the mean squared error $\sum_{\psi\in \Psi_j^{D_Y}}\mb{E}_{\mu^*_X}\Big[ \big(\mb{E}_{\mu^*_{Y|X}}[ 2^{\frac{j(d_Y-D_Y)}{2}} \psi(Y)]-\wh S_j^\dagger(\psi,X)\big)^2\Big]$ for each $j$ by applying Lemma~\ref{lemmaregime2MSE}, which yields
\begin{equation}
    \sum_{\psi\in \Psi_j^{D_Y}}\mb{E}_{\mu^*_X}\Big[ \big(\mb{E}_{\mu^*_{Y|X}}[ 2^{\frac{j(d_Y-D_Y)}{2}} \psi(Y)]-\wh S_j^\dagger(\psi,X)\big)^2\Big]\lesssim 2^{\frac{2j\alpha_Xd_Y}{2\alpha_X+d_X}} \big(\frac{n}{\log n}\big)^{-\frac{2\alpha_X}{2\alpha_X + d_X}}.
\end{equation}
 This further implies
\begin{equation*}
      (E_A)\lesssim \sum_{j=0}^{J}2^{-j\gamma}\big(2^{\frac{j\alpha_Xd_Y}{2\alpha_X+d_X}} \big(\frac{n}{\log n}\big)^{-\frac{\alpha_X}{2\alpha_X + d_X}}\big) \lesssim   (\log n)\cdot (\frac{n}{\log n})^{-\frac{\alpha_X}{2\alpha_X+d_X}}+(\frac{n}{\log n})^{-\frac{\alpha_Y+\gamma}{2\alpha_Y+d_Y+\frac{\alpha_Y}{\alpha_X}d_X}}.
 \end{equation*}
Note that when $\gamma > \frac{d_Y \alpha_X}{2\alpha_X + d_X}$, the dominant term in the summation is at $j=0$, indicating that the bottleneck lies in learning the overall dependence of $Y$ on $X$, reflected by the conditional mean of the wavelets at smaller levels, leading to a term of $n^{-\alpha_X/(2\alpha_X + d_X)}$. Conversely, when $\gamma < \frac{d_Y \alpha_X}{2\alpha_X + d_X}$, the dominant term is at $j=J$, suggesting that the bottleneck is in learning finer irregularities of the conditional distribution, captured by the conditional mean of the wavelets at higher levels, resulting in a term of $n^{-(\alpha_Y + \gamma)/\big(2\alpha_Y + D_Y + \frac{\alpha_Y}{\alpha_X}d_X\big)}$.
 Then for the term $ (E_B)$ and $(E_C)$,  notice that
 \begin{equation*}
     f^{\perp}_J(y)=f(y)-\sum_{j=0}^J \sum_{\psi\in \ov \Psi_j^{D_Y}} f_{\psi} \psi(y)=\sum_{j=J+1}^{\infty} \sum_{\psi\in \ov \Psi_j^{D_Y}} f_{\psi} \psi(y)\lesssim 2^{-J\gamma},
 \end{equation*}
 and there exists a constant $C$ so that for any $y,y'\in \mb R^{D_Y}$, $j\in \mb N$ and $\psi\in \ov\Psi_j^{D_Y}$, 
 \begin{equation*}
     |\psi(y)-\psi(y')|\leq C\, 2^{j+\frac{jD_y}{2}} \|y-y'\|,
 \end{equation*}
 and 
  \begin{equation*}
     |\psi(y)-\psi(y')|\leq  |\psi(y)|+|\psi(y')|\leq C\, 2^{\frac{jD_y}{2}}.
 \end{equation*}
So let $J'=-\log_2(n^{-\frac{1}{\frac{d_Y}{\beta_Y}}}+ n^{-\frac{\alpha_X}{(2\alpha_X+d_X)\gamma}})$, when $1\leq \gamma\leq \frac{d_Y\alpha_Y}{2\alpha_X+d_X}$,
\begin{equation*}
\begin{aligned}
        &\big|f^{\perp}_J(y)-f^{\perp}_J(\wh G_{[k]}(\wh Q_{[k]}(y)))\big|\\
        &\leq \big|\sum_{j=J+1}^{J'} \sum_{\psi\in \ov \Psi_j^{D_Y}} f_{\psi} \psi(y)-\sum_{j=J+1}^{J'} \sum_{\psi\in \ov \Psi_j^{D_Y}} f_{\psi} \psi(\wh G_{[k]}(\wh Q_{[k]}(y)))\big|+n^{-\frac{\gamma}{\frac{d_Y}{\beta_Y}}}+n^{-\frac{\alpha_X}{2\alpha_X+d_X}}\\
        &\lesssim \sum_{j=J+1}^{J'} 2^{-j(\gamma-1)}(2^{-j}\wedge\|y-\wh G_{[k]}(\wh Q_{[k]}(y))\|)+n^{-\frac{\gamma}{\frac{d_Y}{\beta_Y}}}+n^{-\frac{\alpha_X}{2\alpha_X+d_X}};\\
\end{aligned}
\end{equation*}
when $\gamma\leq 1$, let $\gamma_1=(\frac{2\alpha_X+d_X}{2\alpha_X}\gamma)\wedge 1$, then 
\begin{equation*}
\begin{aligned}
        &\big|f^{\perp}_J(y)-f^{\perp}_J(\wh G_{[k]}(\wh Q_{[k]}(y)))\big|\\
        &\leq |\sum_{j=J+1}^{J'} \sum_{\psi\in \ov \Psi_j^{D_Y}} f_{\psi} \psi(y)-\sum_{j=J+1}^{J'} \sum_{\psi\in \ov \Psi_j^{D_Y}} f_{\psi} \psi(\wh G_{[k]}(\wh Q_{[k]}(y)))|+n^{-\frac{\gamma}{\frac{d_Y}{\beta_Y}}}+n^{-\frac{\alpha_X}{2\alpha_X+d_X}}\\
        & \lesssim \sum_{j=J+1}^{J'} 2^{-j(\gamma-1)}(2^{-j}\wedge\|y-\wh G_{[k]}(\wh Q_{[k]}(y))\|)+n^{-\frac{\gamma}{\frac{d_Y}{\beta_Y}}}+n^{-\frac{\alpha_X}{2\alpha_X+d_X}}\\
        &\lesssim  \sum_{j=J+1}^{J'} 2^{j(\gamma_1-\gamma)}\|y-\wh G_{[k]}(\wh Q_{[k]}(y))\|^{\gamma_1}+n^{-\frac{\gamma}{\frac{d_Y}{\beta_Y}}}+n^{-\frac{\alpha_X}{2\alpha_X+d_X}}.\\
\end{aligned}
\end{equation*}
Moreover, as demonstrated in the proof of Theorem~\ref{thgamma<1} in Appendix~\ref{proofthgamma<1},  it holds with probability at least $1-\frac{3}{n^2}$ that
\begin{enumerate}
    \item  for any  $\gamma_1\in(0,1]$ and $k\in \wh{\m K}$,
\begin{equation*}
    \begin{aligned}
   &\mb{E}_{\mu^*_X}\mb{E}_{\mu^*_{Y|X}}[\|Y-\wh G_{[k]}(\wh Q_{[k]}(Y))\|^{\gamma_1}\cdot\mathbf{1}(X\in \mb B_{\mb R^{D_X}}(x_k,2\tau_2))\mathbf{1}(Y\in \mb B_{\mb R^{D_Y}}(y_k,2\tau_2))]\\
          & \lesssim  \left\{
          \begin{array}{cc}
           C \, \frac{(\log n)^{1+\gamma_1}}{\sqrt{n}}   & \frac{d_Y}{\beta_Y} \leq 2\gamma_1, \\
C \,(\log n \wedge \frac{1}{d_Y-2\gamma_1\beta_Y)} )^{1+\gamma_1}\cdot n^{-\frac{\gamma_1}{ \frac{d_Y}{\beta_Y}}}   & \frac{d_Y}{\beta_Y}>2\gamma_1;
          \end{array}
        \right.
        \end{aligned}
        \end{equation*}
        \item for any $j\in \{0\}\cup [J]$ with $J=\lceil {\frac{1}{2\alpha_Y+d_Y+d_X\frac{\alpha_Y}{\alpha_X}}}\cdot \log_2 (\frac{n}{\log n})\rceil$,
           \begin{equation*} 
          \mb{E}_{\mu^*_X} \Big[\sum_{\psi\in \Psi_j^{d_Y}} (\wh v_{k\psi}(X)-  \mb{E}_{\mu^*_{Y|X}}[\psi(\wh Q_{\sk}(Y))\rho_\sk(X,Y)])^2\Big]\lesssim 2^{\frac{2j\alpha_Xd_Y}{2\alpha_X+d_X}} (\frac{n}{\log n})^{-\frac{2\alpha_X}{2\alpha_X+d_X}};
      \end{equation*}
      \item    for any $j\in \mb N$, $\psi\in \Psi_j^{d_Y}$ and $x\in \m M_X$, 
 
    \begin{equation*}
      \Big|\mb{E}_{\mu^*_{Y|x}}[\psi(\wh Q_{\sk}(Y))\rho_\sk(x,Y)]\Big|\lesssim 2^{-\frac{d_Yj}{2}-j\alpha_Y};
\end{equation*}
 
\item for any $k\in [K]\setminus \wh{\m K}$,  $\mb{E}_{\mu^*}[\rho_\sk(X,Y)]\lesssim \sqrt{\frac{\log n}{n}}$.
\end{enumerate}

    \medskip
\noindent So for any $1<\gamma\leq \frac{d_Y\alpha_Y}{2\alpha_X+d_X}$,
\begin{equation*}
     \begin{aligned}
      &   (E_B)\lesssim\mb{E}_{\mu^*}[\sum_{k\in [K]\setminus \wh{\m K}}\rho_\sk(X,Y)]+\mb{E}_{\mu^*}  \Big[\underset{f\in \m H_1^{\gamma}(\mb R^{D_Y})}{\sup} \Big|\sum_{k\in \wh{\m K}}\rho_\sk(X,Y)\big(f_J^{\perp}(Y)-f_J^{\perp}(\wh G_\sk(\wh Q_\sk(Y)))\big)\big]\Big|\Big]\\
  &\lesssim \sqrt{\frac{\log n}{n}}+ (\log n)\cdot 2^{-J(\gamma-1)}\sum_{k\in \wh{\m K}}\mb{E}_{\mu^*} \big[\rho_\sk(X,Y)\|Y-\wh G_{[k]}(\wh Q_{[k]}(Y))\|\big]+n^{-\frac{\gamma}{\frac{d_Y}{\beta_Y}}}+n^{-\frac{\alpha_X}{2\alpha_X+d_X}}\\
  &\lesssim \frac{(\log n)^3}{\sqrt{n}}+(\log n)\cdot (\frac{n}{\log n})^{-\frac{\gamma-1}{2\alpha_Y+d_Y+d_X\frac{\alpha_Y}{\alpha_X}}-\frac{1}{\frac{d_Y}{\beta_Y}}}+n^{-\frac{\gamma}{\frac{d_Y}{\beta_Y}}}+n^{-\frac{\alpha_X}{2\alpha_X+d_X}}\\
  &\lesssim (\frac{n}{\log n})^{-\frac{\alpha_Y+\gamma}{2\alpha_Y+d_Y+\frac{\alpha_Y}{\alpha_X}d_X}}+ (\log n)^3 n^{-\frac{\alpha_X}{2\alpha_X+d_X}},
    \end{aligned}
 \end{equation*}
 where the last inequality uses the fact that $\beta_Y\geq \alpha_Y+1$. Similarly, we can get when $\gamma\leq 1$, 
 \begin{equation*}
     \begin{aligned}
      &   (E_B)\lesssim\mb{E}_{\mu^*}[\sum_{k\in [K]\setminus \wh{\m K}}\rho_\sk(X,Y)]+\mb{E}_{\mu^*}  \Big[\underset{f\in \m H_1^{\gamma}(\mb R^{D_Y})}{\sup} \Big|\sum_{k\in \wh{\m K}}\rho_\sk(X,Y)\big(f_J^{\perp}(Y)-f_J^{\perp}(\wh G_\sk(\wh Q_\sk(Y)))\big)\big]\Big|\Big]\\
  &\lesssim   \sqrt{\frac{\log n}{n}}+n^{-\frac{\gamma}{\frac{d_Y}{\beta_Y}}} +n^{-\frac{\alpha_X}{2\alpha_X+d_X}}+(\log n)\cdot  \sum_{k\in \wh{\m K}}\mb{E}_{\mu^*} \big[\rho_\sk(X,Y)\|Y-\wh G_{[k]}(\wh Q_{[k]}(Y))\|\big]\\
  &+(n^{\frac{\beta_Y(\gamma_1-\gamma)}{d_Y}}\wedge n^{\frac{\alpha_X(\gamma_1/\gamma-1)}{2\alpha_X+d_X}})  \sum_{k\in \wh{\m K}}\mb{E}_{\mu^*} \big[\rho_\sk(X,Y)\|Y-\wh G_{[k]}(\wh Q_{[k]}(Y))\|^{\gamma_1}\big]\\
    &\lesssim  n^{-\frac{\gamma}{\frac{d_Y}{\beta_Y}}} +(\log n)^3\cdot n^{-\frac{\alpha_X}{2\alpha_X+d_X}}+  (\frac{n}{\log n})^{-\frac{\alpha_Y+\gamma}{2\alpha_Y+d_Y+\frac{\alpha_Y}{\alpha_X}d_X}},
    \end{aligned}
 \end{equation*}
where the last inequality uses that for $\gamma_1=(\frac{2\alpha_X+d_X}{2\alpha_X}\gamma)\wedge 1$, it holds that 
\begin{equation*}
    \begin{aligned}
        &(n^{\frac{\beta_Y(\gamma_1-\gamma)}{d_Y}}\wedge n^{\frac{\alpha_X(\gamma_1/\gamma-1)}{2\alpha_X+d_X}})  \sum_{k\in \wh{\m K}}\mb{E}_{\mu^*} \big[\rho_\sk(X,Y)\|Y-\wh G_{[k]}(\wh Q_{[k]}(Y))\|^{\gamma_1}\big]\\
        &\lesssim    (n^{\frac{\beta_Y(\gamma_1-\gamma)}{d_Y}}\wedge n^{\frac{\alpha_X(\gamma_1/\gamma-1)}{2\alpha_X+d_X}})\cdot (n^{-\frac{\gamma_1}{\frac{d_Y}{\beta_Y}}}+\frac{(\log n)^2}{\sqrt{n}})\\
        &\lesssim  n^{-\frac{\gamma}{\frac{d_Y}{\beta_Y}}} + \frac{(\log n)^2}{\sqrt{n}}  n^{\frac{\alpha_X(\frac{2\alpha_X+d_X}{2\alpha_X}-1)}{2\alpha_X+d_X}}\\
        &\lesssim  n^{-\frac{\gamma}{\frac{d_Y}{\beta_Y}}}+(\log n)^2\cdot n^{-\frac{\alpha_X}{2\alpha_X+d_X}}.
    \end{aligned}
\end{equation*}
When $\gamma\geq \eta=\frac{d_Y\alpha_Y}{2\alpha_X+d_X}$,  we have
\begin{equation*}
\begin{aligned}
         (E_B)&\lesssim \sqrt{\frac{\log n}{n}}+\mb{E}_{\mu^*}  \Big[\underset{f\in \m H_1^{\eta}(\mb R^{D_Y})}{\sup} \Big|\sum_{k\in \wh{\m K}}\rho_\sk(X,Y)\big(f_J^{\perp}(Y)-f_J^{\perp}(\wh G_\sk(\wh Q_\sk(Y),X))\big)\big]\Big|\Big]\\
        &\lesssim  (\log n)^3 n^{-\frac{\alpha_X}{2\alpha_X+d_X}}.
\end{aligned}
\end{equation*}
Finally, for term $(E_C)$, it holds that
\begin{equation*}
     \begin{aligned}
      &  (E_C)=\mb{E}_{\mu^*_X}  \Big[\underset{f\in \m H_1^{\gamma}(\mb R^{D_Y})}{\sup} \Big|\sum_{k\in \wh{\m K}} \mb{E}_{\mu^*_{Y|X}}\big[f_J^{\perp}(\wh G_\sk(\wh Q_\sk(Y)))\rho_\sk(X,Y)\big]-\int_{\mb R^{d_Y}} f^{\perp}_J(\wh G_\sk (z)) \sum_{j=0}^J \sum_{\psi\in \Psi_j^{d_Y}}\psi(z) \wh v_{k\psi}(X)\,\dd z\Big|\Big]\\
        &\overset{(i)}\lesssim  (\frac{n}{\log n})^{-\frac{\gamma}{2\alpha_Y+d_Y+d_X\frac{\alpha_Y}{\alpha_X}}}\mb{E}_{\mu^*_X} \Big[\underset{f:\mb R^{d_Y}\to \mb R\atop \int  f^2(y)\,\dd y\leq 1}{\sup}  \sum_{k\in \wh{\m K}} \mb{E}_{\mu^*_{Y|X}} [f(\wh Q_{[k]}(Y))\rho_\sk(X,Y)]\\
        &\qquad\qquad-\sum_{k\in \wh{\m K}}\int_{\mb R^{d_Y}} f(z)\sum_{j=0}^J \sum_{\psi\in \Psi_j^{d_Y}}\psi(z) \wh v_{k\psi}(X)\,\dd z \Big]\\
    &\lesssim   (\frac{n}{\log n})^{-\frac{\gamma}{2\alpha_Y+d_Y+d_X\frac{\alpha_Y}{\alpha_X}}}\mb{E}_{\mu^*_X} \Big[\sum_{k\in \wh{\m K}}\underset{f:\mb R^{d_Y}\to \mb R\atop \int f^2(y)\,\dd y\leq 1}{\sup}\sqrt{\sum_{j=0}^{\infty} \sum_{\psi\in \Psi_j^{d_Y}} f_{\psi}^2 }\\
    &\qquad \cdot\sqrt{\sum_{j=0}^{J}\sum_{\psi\in \Psi_j^{d_Y}}\Big( \mb{E}_{\mu^*_{Y|X}}[\psi(\wh Q_\sk(Y))\rho_\sk(X,Y)]- \wh v_{k\psi}(X) \Big)^2+\sum_{j=J+1}^{\infty}\sum_{\psi\in \Psi_j^{d_Y}}\Big( \mb{E}_{\mu^*_{Y|X}}[\psi(\wh Q_\sk(Y))\rho_\sk(X,Y)] \Big)^2}\Big] \\
  &\lesssim    (\frac{n}{\log n})^{-\frac{\alpha_Y+\gamma}{2\alpha_Y+d_Y+\frac{\alpha_Y}{\alpha_X}d_X}} \\
  &\qquad+  (\frac{n}{\log n})^{-\frac{\gamma}{2\alpha_Y+d_Y+d_X\frac{\alpha_Y}{\alpha_X}}}\sum_{k\in \wh{\m K}} \sqrt{\mb{E}_{\mu^*_X}\Big[\sum_{j=0}^{J}\sum_{\psi\in \Psi_j^{d_Y}}\Big( \mb{E}_{\mu^*_{Y|X}}[\psi(\wh Q_\sk(Y))\rho_\sk(X,Y)]-\wh v_{k\psi}(X)\Big)^2\Big]}\\
  &\lesssim (\frac{n}{\log n})^{-\frac{\alpha_Y+\gamma}{2\alpha_Y+d_Y+\frac{\alpha_Y}{\alpha_X}d_X}},
    \end{aligned}
 \end{equation*}
 where to derive $(i)$, we utilize the property that there exists a positive constant $c$ so that for any $k\in \wh {\m K}$, $x\in \m M_X$ and $y\in \mb R^{D_Y}$:
 \begin{equation*}
 f_J^{\perp}(\wh G_\sk(\wh Q_\sk(y)))\rho_\sk(x,y)=f_J^{\perp}(\wh G_\sk(\wh Q_\sk(y)))\rho_\sk(x,y) \cdot \rho(\|\wh Q_\sk(y)\|^2/c^2),
 \end{equation*}
 and for any $z\in \mb R^{d_Y}$,
 \begin{equation*}
     f^{\perp}_J(\wh G_\sk (z)) \sum_{j=0}^J \sum_{\psi\in \Psi_j^{d_Y}}\psi(z)\wh v_{k\psi}(x)=f^{\perp}_J(\wh G_\sk (z)) \sum_{j=0}^J \sum_{\psi\in \Psi_j^{d_Y}}\psi(z)\wh v_{k\psi}(x) \cdot \rho(\|z\|^2/c^2),
 \end{equation*}
 where $\rho$ is the smooth transition function defined in~\eqref{def:transition}.
Furthermore, there exists a constant $C$ so that for any $f\in \m H_1^{\gamma}(\mb R^{D_Y})$, the function $2^{J\gamma}f_J^{\perp}(\wh G_\sk(z))\cdot \rho(\|z\|^2/c^2)$ satisfies
\begin{equation*}
\begin{aligned}
   & \int_{\mb R^{d_Y}} \Big(2^{J\gamma}f_J^{\perp}(\wh G_\sk(z))\cdot \rho(\|z\|^2/c^2)\Big)^2\,\dd z\\
&\leq \int_{\mb B_{\mb R^{d_Y}}(\mathbf{0},\sqrt{2}c)} \Big(2^{J\gamma}f_J^{\perp}(\wh G_\sk(z)) \Big)^2\,\dd z\\
&\leq  \int_{\mb B_{\mb R^{d_Y}}(\mathbf{0},\sqrt{2}c)} \dd z\cdot \sup_{z\in  \mb B_{\mb R^{d_Y}}(\mathbf{0},\sqrt{2}c)} \Big(2^{J\gamma}f_J^{\perp}(\wh G_\sk(z)) \Big)^2 \\
&\leq C.
\end{aligned}
 \end{equation*}
Therefore, it holds for any $x\in \m M_X$ that
\begin{equation*}
    \begin{aligned}
 &\underset{f\in \m H_1^{\gamma}(\mb R^{D_Y})}{\sup} \Big|\sum_{k\in \wh{\m K}} \mb{E}_{\mu^*_{Y|X}}\big[f_J^{\perp}(\wh G_\sk(\wh Q_\sk(Y)))\rho_\sk(X,Y)\big]-\int_{\mb R^{d_Y}} f^{\perp}_J(\wh G_\sk (z)) \sum_{j=0}^J \sum_{\psi\in \Psi_j^{d_Y}}\psi(z) \wh v_{k\psi}(X)\,\dd z\Big|\\
  &=\underset{f\in \m H_1^{\gamma}(\mb R^{D_Y})}{\sup} \Big|\sum_{k\in \wh{\m K}} \mb{E}_{\mu^*_{Y|X}}\big[f_J^{\perp}(\wh G_\sk(\wh Q_\sk(Y))) \rho(\|\wh Q_\sk(Y)\|^2/c^2)\rho_\sk(X,Y)\big]\\
  &\qquad\qquad\quad-\int_{\mb R^{d_Y}} f^{\perp}_J(\wh G_\sk (z))\rho(\|z\|^2/c^2) \sum_{j=0}^J \sum_{\psi\in \Psi_j^{d_Y}}\psi(z)  \wh v_{k\psi}(X)\,\dd z\Big|\\
        & =2^{-J\gamma} \underset{ \wt f(z)=2^{J\gamma}f_J^{\perp}(\wh G_\sk(z))\cdot \rho(\|z\|^2/c^2)\atop
           f\in \m H_1^{\gamma}(\mb R^{D_Y})}{\sup} \Big|\sum_{k\in \wh{\m K}} \mb{E}_{\mu^*_{Y|X}}\big[\wt f(\wh Q_\sk(Y))\rho_\sk(X,Y)\big]-\int_{\mb R^{d_Y}} \wt f(z) \sum_{j=0}^J \sum_{\psi\in \Psi_j^{d_Y}}\psi(z) \wh v_{k\psi}(X)\,\dd z \Big|\\
        &\leq   \sqrt{C}\cdot2^{-J\gamma}\underset{f:\mb R^{d_Y}\to \mb R\atop \int f^2(y)\,\dd y\leq 1}{\sup}  \Big|\sum_{k\in \wh{\m K}} \mb{E}_{\mu^*_{Y|X}} [f(\wh Q_{[k]}(Y))\rho_\sk(X,Y)]-\sum_{k\in \wh{\m K}}\int_{\mb R^{d_Y}} f(z)\sum_{j=0}^J \sum_{\psi\in \Psi_j^{d_Y}}\psi(z) \wh v_{k\psi}(X)\,\dd z \Big|,
    \end{aligned}
\end{equation*}
 which further substantiates inequality $(i)$. Finally, by combining the bounds for term $ (E_A)$, $ (E_B)$, $(E_C)$, we can then get the desired results.

\subsubsection{Proof for Regime 3b}\label{proofth3}
The overall structure  of the  proof  mirrors that for Regime 2, as detailed in Appendix~\ref{proofth:thgamma>0}.  We consider the estimator $ \wh{\m J}(f,x)$ defined in Appendix~\ref{R3gamma>0}:
\begin{equation*}
\begin{aligned}
       & \wh{\m J}(f,x)= {\sum_{j=0}^J\sum_{\psi\in\Psi_j^{D_Y}} f_{\psi}2^{\frac{j(D_Y-d_Y)}{2}} \wh S_j^\dagger(\psi,x) }+ {\sum_{k\in \wh{\m K}} \int_{\mb R^{d_Y}} f^{\perp}_J(\wh G_\sk (z,x)) \sum_{j=0}^J \sum_{\psi\in \Psi_j^{d_Y}}\psi(z) \wh v_{k\psi}(x)\,\dd z},\\
       &\qquad f_{\psi}=\int_{\mb R^{D_Y}} f(y)\psi(y)\,\dd y,\quad f^{\perp}_J(y)=f(y)-\sum_{j=0}^J \sum_{\psi\in  \ov\Psi_j^{D_Y}} f_{\psi} \psi(y),
\end{aligned}
\end{equation*}
where $J=\lceil {\frac{1}{2\alpha_Y+d_Y+d_X\frac{\alpha_Y}{\alpha_X}}}\cdot \log_2 (\frac{n}{\log n})\rceil$.
 We can get 
  \begin{equation*}
     \begin{aligned}
          & \mb{E}_{\mu^*_X}  \Big[\underset{f\in \m H_1^{\gamma}(\mb R^{D_Y})}{\sup} \big|\mb{E}_{\mu^*_{Y|X}} [f(Y)]- \wh{\m J}(f,X)\big|\Big]\\
           &\leq   \underbrace{ \mb{E}_{\mu^*_X} \Big[ \underset{f\in \m H_1^{\gamma}(\mb R^{D_Y})}{\sup} \big| \sum_{j=0}^J \sum_{\psi\in \Psi_j^{D_Y}} f_{\psi} \big(\mb{E}_{\mu^*_{Y|X}}[\psi(Y)]- 2^{\frac{j(D_Y-d_Y)}{2}} \wh S_j^\dagger(\psi,X)\big)\big|\Big]}_{ (E_A)}  \\
          &+ \underbrace{\mb{E}_{\mu^*_X}  \Big[\underset{f\in \m H_1^{\gamma}(\mb R^{D_Y})}{\sup} \Big|\mb{E}_{\mu^*_{Y|X}}\big[\sum_{k\in [K]\setminus\wh{\m K}}\rho_\sk (X,Y)f_J^{\perp}(Y)+\sum_{k\in \wh{\m K}}\rho_\sk(X,Y)\big(f_J^{\perp}(y)-f_J^{\perp}(\wh G_\sk(\wh Q_\sk(Y),X))\big)\big]\Big|\Big]}_{ (E_B)} \\
          &+ \mb{E}_{\mu^*_X}  \Big[\underset{f\in \m H_1^{\gamma}(\mb R^{D_Y})}{\sup} \Big|\sum_{k\in \wh{\m K}} \mb{E}_{\mu^*_{Y|X}}\big[f_J^{\perp}(\wh G_\sk(\wh Q_\sk(Y),X))\rho_\sk(X,Y)\big]\\
          &\underbrace{\qquad\qquad\qquad-\int_{\mb R^{d_Y}} f^{\perp}_J(\wh G_\sk (z,X)) \sum_{j=0}^J \sum_{\psi\in \Psi_j^{d_Y}}\psi(z) \wh v_{k\psi}(X)\,\dd z\Big|\Big]}_{(E_C)}.\\
           \end{aligned}
 \end{equation*}
To bound term $ (E_A)$, notice that similarly to Regime 2, we have
\begin{equation*}
    \begin{aligned}
        &\mb{E}_{\mu^*_X} \Big[ \underset{f\in \m H_1^{\gamma}(\mb R^{D_Y})}{\sup} \big| \sum_{j=0}^J \sum_{\psi\in   \Psi_j^{D_Y}} f_{\psi} \big(\mb{E}_{\mu^*_{Y|X}}[\psi(Y)]- 2^{\frac{j(D_Y-d_Y)}{2}} \wh S_j^\dagger(\psi,X)\big)\big|\Big]\\
 &\lesssim \sum_{j=0}^{J}\sum_{\psi\in  \Psi_j^{D_Y}}2^{-j\gamma-\frac{jD_Y}{2}}\sqrt{\mb{E}_{\mu^*_X}\Big[ \big(\mb{E}_{\mu^*_{Y|X}}[ 2^{\frac{j(d_Y-D_Y)}{2}} \psi(Y)]-\wh S_j^\dagger(\psi,X)\big)^2\Big]}\\
     &\lesssim  \sum_{j=0}^{J}\sqrt{\sum_{\psi\in  \Psi_j^{D_Y}}2^{-2j\gamma}\mb{E}_{\mu^*_X}\Big[ \big(\mb{E}_{\mu^*_{Y|X}}[ 2^{\frac{j(d_Y-D_Y)}{2}} \psi(Y)]-\wh S_j^\dagger(\psi,X)\big)^2\Big]}.
    \end{aligned}
\end{equation*}
Then we bound $\sum_{\psi\in  \Psi_j^{D_Y}}\mb{E}_{\mu^*_X}\Big[ \big(\mb{E}_{\mu^*_{Y|X}}[ 2^{\frac{j(d_Y-D_Y)}{2}} \psi(Y)]-\wh S_j^\dagger(\psi,X)\big)^2\Big]$, where recall that
\begin{equation*}
 \begin{aligned}
         &\wh S_j^{\dagger}={\arg\min}_{S\in  {\m S}_j^\dagger} \frac{1}{|I_1|}\sum_{i\in I_1} \sum_{\psi\in  \Psi_j^{D_Y}} (2^{\frac{j(d_Y-D_Y)}{2}}\psi(Y_i)-S(\psi,X_i))^2.\\
 \end{aligned}
\end{equation*}
\begin{lemma}\label{le2:app}
 Suppose $\mu^*\in \m P^*_3$ and with the choices of $\m S_j^\dagger$ defined in~\eqref{defSdaggerr3},   it holds with probability larger than  $1-\frac{1}{n^2}$ that for any $j \in [J]$,
\begin{equation*}
\begin{aligned}
  &\mb{E}_{\mu^*_X} \Big[\sum_{\psi\in   \Psi_j^{D_Y}} \big(\mb{E}_{\mu^*_{Y|X}}[2^{\frac{j(d_Y-D_Y)}{2}}\psi(Y)]-\wh S_j^\dagger(\psi,X)\big)^2\Big]
   \lesssim   \frac{\log n}{n} 2^{jd_Y} (\varepsilon_j^x)^{-d_X} +(\log n)^2\cdot (\varepsilon_j^x)^{2\alpha_X},
\end{aligned}
\end{equation*}
where recall $\varepsilon_j^x= 2^{\frac{jd_Y}{2\alpha_X+d_X}}(\frac{n}{\log n})^{-\frac{1}{2\alpha_X+d_X}}$ and $J=\lceil {\frac{1}{2\alpha_Y+d_Y+d_X\frac{\alpha_Y}{\alpha_X}}}\cdot \log_2 (\frac{n}{\log n})\rceil$.
\end{lemma}
 
\noindent The proof of Lemma~\ref{le2:app} is provided in Appendix~\ref{proofle2:app}. So it holds with probability at least $1-\frac{1}{n^2}$ that for any $\gamma\in (0,1]$,
 \begin{equation*}
 \begin{aligned}
        (E_A)&\lesssim \sum_{j=0}^J 2^{-j\gamma} \Big(\sqrt{\frac{\log n}{n}} 2^{jd_Y/2} (\varepsilon_j^x)^{-d_X/2} +(\log n)\cdot (\varepsilon_j^x)^{\alpha_X}\Big)\\
      & \lesssim  (\log n)^2\cdot (\frac{n}{\log n})^{-\frac{\alpha_X}{2\alpha_X+d_X}}+ (\log n)\cdot (\frac{n}{\log n})^{-\frac{\alpha_Y+\gamma}{2\alpha_Y+D_Y+\frac{\alpha_Y}{\alpha_X}d_X}}.
 \end{aligned}
 \end{equation*}
Now we bound term $ (E_B)$. Follow the same procedure as in the proof for Regime 2, let $J'=-\log_2(n^{-\frac{1}{\frac{d_Y}{\beta_Y}+\frac{d_X}{\beta_X}}}+ n^{-\frac{\alpha_X}{(2\alpha_X+d_X)\gamma}})$,  we can get, when $1\leq \gamma\leq \frac{\big(d_Y\vee (\frac{d_Y}{\beta_Y}+\frac{d_X}{\beta_X})\big)\alpha_Y}{2\alpha_X+d_X}$,
\begin{equation*}
\begin{aligned}
        &\big|f^{\perp}_J(y)-f^{\perp}_J(\wh G_{[k]}(\wh Q_{[k]}(y),x))\big|\lesssim \sum_{j=J+1}^{J'} 2^{-j(\gamma-1)}(2^{-j}\wedge\|y-\wh G_{[k]}(\wh Q_{[k]}(y),x)\|)+n^{-\frac{\gamma}{\frac{d_Y}{\beta_Y}+\frac{d_X}{\beta_X}}}+n^{-\frac{\alpha_X}{2\alpha_X+d_X}};\\
\end{aligned}
\end{equation*}
and when $\gamma\leq 1$, let $\gamma_1=(\frac{2\alpha_X+d_X}{2\alpha_X}\gamma)\wedge 1$, it holds that 
\begin{equation*}
\begin{aligned}
        &\big|f^{\perp}_J(y)-f^{\perp}_J(\wh G_{[k]}(\wh Q_{[k]}(y),x))\big|\lesssim  \sum_{j=J+1}^{J'} 2^{j(\gamma_1-\gamma)}\|y-\wh G_{[k]}(\wh Q_{[k]}(y),x)\|^{\gamma_1}+n^{-\frac{\gamma}{\frac{d_Y}{\beta_Y}+\frac{d_X}{\beta_X}}}+n^{-\frac{\alpha_X}{2\alpha_X+d_X}}.\\
\end{aligned}
\end{equation*}
Then we establish a bound on the population-level reconstruction error in the following lemma, the proof of which is given in Appendix~\ref{proof:lemma:manifoldlearningr3}.
 
\begin{lemma}\label{lemma:manifoldlearningr3}
Suppose $\mu^*\in \m P_3$ and  with the choices of $\m G$ defined in~\eqref{defGr3},  for any $0<\gamma_1\leq 1$,  it holds with probability at least $1-\frac{1}{n^2}$ that
  \begin{enumerate}
 \item  For any $k\in \wh{\m K}$ and $\gamma_1\in(0,1]$,
    \begin{equation*}
    \begin{aligned}
           &\mb{E}_{\mu^*_X}\mb{E}_{\mu^*_{Y|X}}[\|Y-\wh G_{[k]}(\wh Q_{[k]}(Y),X)\|^{\gamma_1}\cdot\mathbf{1}(X\in \mb B_{\mb R^{D_X}}(x_k,2\tau_2))\mathbf{1}(Y\in \mb B_{\mb R^{D_Y}}(y_k,2\tau_2))]\\
          & \lesssim \left\{
          \begin{array}{cc}
            \frac{(\log n)^{1+\gamma_1}}{\sqrt{n}}   & \frac{d_Y}{\beta_Y}+\frac{d_X}{\beta_X}\leq 2\gamma_1, \\
\Big((\log n \wedge \frac{1}{\beta_Y(d_Y/\beta_Y+d_X/\beta_X-2\gamma_1)} )^{1+\gamma_1}+(\log n)^{\gamma_1})\cdot n^{-\frac{\gamma_1}{\frac{d_X}{\beta_X}+\frac{d_Y}{\beta_Y}}}   & \frac{d_Y}{\beta_Y}+\frac{d_X}{\beta_X}>2\gamma_1.
          \end{array}
        \right.
        \end{aligned}
        \end{equation*}
        \item  For any $k\in \wh {\m K}$, there exists $(x^*_k,y^*_k)\in \mb B_{\m M}((x_k,y_k),\sqrt{2}\tau_2)$  such that 
\begin{equation*}
    \wh V_\sk^T P^*_\sk\wh V_\sk\gtrsim C_1I_{d_Y},
\end{equation*}
where $\m P^*_\sk$ is the projection matrix of $T_{\m M_{Y|x^*_k}} y^*_k$.
\end{enumerate}
\end{lemma}
\noindent Moreover, since for any $k\in [K]\setminus\wh {\m K}$, it holds that 
 \begin{equation*}
     \frac{1}{n}\sum_{i\in I_1}\rho_\sk(X_i,Y_i)\leq   \frac{1}{n}\sum_{i\in I_1}\mathbf{1}(\|(X_i,Y_i)-(x_k,y_k)\|\leq \sqrt{2}\tau_2)=0.
 \end{equation*}
By Bernstein's inequality, it holds with probability at least $1-\frac{1}{n^2}$ that for any $k\in [K]\setminus\wh {\m K}$,
\begin{equation*}
    \mb{E}_{\mu^*}[\rho_\sk(X,Y)]\lesssim\sqrt{\frac{\log n}{n}}.
\end{equation*}
Therefore it holds with probability at least $1-\frac{2}{n^2}$ that for any $1\leq \gamma\leq \frac{\big(d_Y\vee (\frac{d_Y}{\beta_Y}+\frac{d_X}{\beta_X})\big)\alpha_Y}{2\alpha_X+d_X}$,
\begin{equation*}
     \begin{aligned}
      &   (E_B)\lesssim \sqrt{\frac{\log n}{n}}+\mb{E}_{\mu^*}  \Big[\underset{f\in \m H_1^{\gamma}(\mb R^{D_Y})}{\sup} \Big|\sum_{k\in \wh{\m K}}\rho_\sk(X,Y)\big(f_J^{\perp}(Y)-f_J^{\perp}(\wh G_\sk(\wh Q_\sk(Y),X))\big)\big]\Big|\Big]\\
  &\lesssim \sqrt{\frac{\log n}{n}}+ (\log n)\cdot 2^{-J(\gamma-1)}\sum_{k\in \wh{\m K}}\mb{E}_{\mu^*} \big[\rho_\sk(X,Y)\|Y-\wh G_{[k]}(\wh Q_{[k]}(Y),X)\|\big]+n^{-\frac{\gamma}{\frac{d_Y}{\beta_Y}+\frac{d_X}{\beta_X}}}+n^{-\frac{\alpha_X}{2\alpha_X+d_X}}\\
  &\lesssim (\log n)^{2}\cdot (\frac{n}{\log n})^{-\frac{\gamma-1}{2\alpha_Y+d_Y+d_X\frac{\alpha_Y}{\alpha_X}}-\frac{1}{\frac{d_Y}{\beta_Y}+\frac{d_X}{\beta_X}}}+\frac{(\log n)^3}{\sqrt{n}}+n^{-\frac{\gamma}{\frac{d_Y}{\beta_Y}+\frac{d_X}{\beta_X}}}+n^{-\frac{\alpha_X}{2\alpha_X+d_X}}\\
  &\lesssim (\log n)^3\cdot n^{-\frac{\alpha_X}{2\alpha_X+d_X}}+ (\log n)\cdot (\frac{n}{\log n})^{-\frac{\alpha_Y+\gamma}{2\alpha_Y+D_Y+\frac{\alpha_Y}{\alpha_X}d_X}},
    \end{aligned}
 \end{equation*}
 where the last inequality uses $\beta_Y\geq \alpha_Y+1$ and $\beta_X\geq \alpha_X+\frac{\alpha_X}{\alpha_Y}$.
 Similarly, we can get when $\gamma\leq 1$, 
 \begin{equation*}
     \begin{aligned}
      &   (E_B)\lesssim \sqrt{\frac{\log n}{n}}+\mb{E}_{\mu^*_X}  \Big[\underset{f\in \m H_1^{\gamma}(\mb R^{D_Y})}{\sup} \Big|\sum_{k\in \wh{\m K}}\rho_\sk(x,y)\big(f_J^{\perp}(y)-f_J^{\perp}(\wh G_\sk(\wh Q_\sk(y),x))\big)\big]\Big|\Big]\\
  &\lesssim  n^{-\frac{\gamma}{\frac{d_Y}{\beta_Y}+\frac{d_X}{\beta_X}}} +n^{-\frac{\alpha_X}{2\alpha_X+d_X}}+(\log n)\cdot  \sum_{k\in \wh{\m K}}\mb{E}_{\mu^*} \big[\rho_\sk(X,Y)\|Y-\wh G_{[k]}(\wh Q_{[k]}(Y),X)\|\big]\\
  &+(n^{\frac{(\gamma_1-\gamma)}{\frac{d_Y}{\beta_Y}+\frac{d_X}{\beta_X}}}\wedge n^{\frac{\alpha_X(\gamma_1/\gamma-1)}{2\alpha_X+d_X}})  \sum_{k\in \wh{\m K}}\mb{E}_{\mu^*} \big[\rho_\sk(X,Y)\|Y-\wh G_{[k]}(\wh Q_{[k]}(Y),X)\|^{\gamma_1}\big]\\
    &\lesssim (\log n)\cdot n^{-\frac{\gamma}{\frac{d_Y}{\beta_Y}+\frac{d_X}{\beta_X}}} +(\log n)^2\cdot n^{-\frac{\alpha_X}{2\alpha_X+d_X}}+\frac{(\log n)^3}{\sqrt{n}},
    \end{aligned}
 \end{equation*}
 where the last inequality uses that for any $\gamma_1=(\frac{2\alpha_X+d_X}{2\alpha_X}\gamma)\wedge 1$,
\begin{equation*}
    \begin{aligned}
        &(n^{\frac{(\gamma_1-\gamma)}{\frac{d_Y}{\beta_Y}+\frac{d_X}{\beta_X}}}\wedge n^{\frac{\alpha_X(\gamma_1/\gamma-1)}{2\alpha_X+d_X}})  \sum_{k\in \wh{\m K}}\mb{E}_{\mu^*} \big[\rho_\sk(X,Y)\|Y-\wh G_{[k]}(\wh Q_{[k]}(Y),X)\|^{\gamma_1}\big]\\
        &\lesssim    (n^{\frac{(\gamma_1-\gamma)}{\frac{d_Y}{\beta_Y}+\frac{d_X}{\beta_X}}}\wedge n^{\frac{\alpha_X(\gamma_1/\gamma-1)}{2\alpha_X+d_X}})\cdot (n^{-\frac{\gamma_1}{\frac{d_Y}{\beta_Y}+\frac{d_X}{\beta_X}}}+\frac{(\log n)^2}{\sqrt{n}})\\
        &\lesssim  n^{-\frac{\gamma}{\frac{d_Y}{\beta_Y}+\frac{d_X}{\beta_X}}} + \frac{(\log n)^2}{\sqrt{n}}  n^{\frac{\alpha_X(\frac{2\alpha_X+d_X}{2\alpha_X}-1)}{2\alpha_X+d_X}}\\
        &\lesssim  n^{-\frac{\gamma}{\frac{d_Y}{\beta_Y}+\frac{d_X}{\beta_X}}}+(\log n)^2\cdot n^{-\frac{\alpha_X}{2\alpha_X+d_X}}.
    \end{aligned}
\end{equation*}
When $\gamma\geq \eta=\frac{\big(d_Y\vee (\frac{d_Y}{\beta_Y}+\frac{d_X}{\beta_X})\big)\alpha_Y}{2\alpha_X+d_X}$, we have 
\begin{equation*}
\begin{aligned}
         (E_B)&\lesssim \sqrt{\frac{\log n}{n}}+\mb{E}_{\mu^*_X}  \Big[\underset{f\in \m H_1^{\eta}(\mb R^{D_Y})}{\sup} \Big|\sum_{k\in \wh{\m K}}\rho_\sk(x,y)\big(f_J^{\perp}(y)-f_J^{\perp}(\wh G_\sk(\wh Q_\sk(y),x))\big)\big]\Big|\Big]\\
        &\lesssim  (\log n)^3 n^{-\frac{\alpha_X}{2\alpha_X+d_X}}.
\end{aligned}
\end{equation*}
 Finally, we bound term $(E_C)$.  Given the second statement in Lemma~\ref{lemma:manifoldlearningr3}, we can use Lemma~\ref{le:projection} in Appendix~\ref{app: regularitymanifold} to obtain the invertibility of $\wh V_\sk^T(\cdot-y_k)$. Specifically, when $\tau_2$ is small enough,  for any $k\in \wh{\m K}$ and $x\in \mb B_{\m M_x}(x_k^*,3\tau_2)$, there exists a subset $\wh U^\sk_{Y|x}$ so that $\mb B_{\m M_{Y|x}}(y_k^*,3\tau_2)\subset \wh U^\sk_{Y|x}\subset \m M_{Y|x}$, and the function $\wh Q_\sk(\cdot)=\wh V_\sk^T(\cdot-y_k)$, when restricted to domain $\wh U^\sk_{Y|x}$, is a diffeomorphism  that maps $\wh U^\sk_{Y|x}$ to $\mb B_{\mb R^{d_Y}}(\wh V_\sk^T(y_k^*-y_k),3\tau_2)$ with inverse denoted as $[\wh Q_\sk(\cdot,x)]^{-1}$. The function $\wh G^{\dagger}_\sk: \mb B_{\mb R^{d_Y}}(\wh V_\sk^T(y_k^*-y_k),3\tau_2)\times  \mb B_{\m M_X}(x_k^*,3\tau_2)\to \mb R^{D_Y}$ defined as  $\wh G^{\dagger}_\sk(z,x)=[\wh Q_\sk(\cdot,x)]^{-1}(z)$ belongs to $ {\m H}^{\beta_Y,
 \beta_X}_{L,D_Y}(\mb B_{\mb R^{d_Y}}(\wh V_\sk^T(y_k^*-y_k),3\tau_2),\mb B_{\m M_X}(x_k^*,3\tau_2))$. Then by Lemma~\ref{le:defdistribution},  the push forward measure $\wh Q_{\sk\#}(\mu^*_{Y|x}|_{\wh U^\sk_{Y|x}})$ has a density $\wh \nu_{\sk}(z,|,x)\in  {\m H}^{\alpha_Y,\alpha_X}_{L_1}(\mb B_{\mb R^{d_Y}}(\wh V_\sk^T(y_k^*-y_k),3\tau_2),\mb B_{\m M_X}(x_k^*,3\tau_2))$ for some constant $L_1$.
Furthermore, for any $j\in \mb N$ and $\psi\in \Psi_j^{d_Y}$, we have 
\begin{equation*}
\begin{aligned}
     & \mb{E}_{\mu^*_{Y|x}}[\psi(\wh Q_{\sk}(Y))\rho_\sk(x,Y)]\\
     &=  \mb{E}_{\mu^*_{Y|x}}[\psi(\wh Q_{\sk}(Y))\rho_\sk(x,Y)\mathbf{1}(Y\in \wh U^\sk_{Y|x})\mathbf{1}(x\in \mb B_{\m M_x}(x_k^*,2\tau_2))]\\
    & =  \mb{E}_{\mu^*_{Y|x}}[\psi(\wh Q_{\sk}(Y))\rho_\sk(x,\wh G_\sk^{\dagger}(\wh Q_{\sk}(Y),x))\mathbf{1}(Y\in \wh U^\sk_{Y|x})]\\
    &=\int_{\mb B_{\mb R^{d_Y}}(\wh V_\sk^T(y_k^*-y_k),3\tau_2)} \psi(z) \rho_\sk(x, \wh G_\sk^{\dagger}(z,x))\wh v_\sk (z|x)\,\dd z.
\end{aligned}
\end{equation*}
Let $\ov \nu_{\sk}(z,|,x)\in  {\m H}^{\alpha_Y,\alpha_X}_{L_1}(\mb R^{d_Y},\mb R^{D_X})$ be a smooth extension of $\wh \nu_{\sk}(z,|,x)$ to $\mb R^{d_Y}\times \mb R^{D_X}$. Then we define a function $ \wt v_\sk:\mb R^{d_Y}\times \mb R^{D_X}\to \mb R$:
\begin{equation*}
    \wt v_\sk (z,x)=\left\{
    \begin{array}{cc}
     \rho_\sk(x, \wh G_\sk^{\dagger}(z,x))\ov \nu_\sk (z|x),    & \text{ if } z\in \mb B_{\mb R^{d_Y}}(\wh V_\sk^T(y_k^*-y_k),3\tau_2), x\in \mb B_{\mb R^{D_X}}(x_k^*,3\tau_2), \\
      0   & \text{otherwise.}
    \end{array}
    \right.
\end{equation*}
We can verify that $ \wt v_\sk (z,x)\in  \ov{\m H}^{\alpha_Y,\alpha_X}_{L_2}(\mb R^{d_Y},\mb R^{D_X})$ for some constant $L_2$. So for any $j\in \mb N$, $\psi\in \Psi_j^{d_Y}$ and $x\in \m M_X$,
\begin{equation*}
     2^{\frac{d_Yj}{2}} \mb{E}_{\mu^*_{Y|x}}[\psi(\wh Q_{\sk}(Y))\rho_\sk(x,Y)]=2^{\frac{d_Yj}{2}}\int_{\mb R^{d_Y}} \psi(z) \wt v_\sk (z|x)\,\dd z,
\end{equation*}
and $2^{\frac{d_Yj}{2}}\int_{\mb R^{d_Y}} \psi(z) \wt v_\sk (z|\cdot)\,\dd z\in \m H^{\alpha_X}_{L_3}(\mb R^{D_X})$ for some constant $L_3$. Moreover, for any $x\in \m M_X$, since $\wt v_\sk (\cdot|x)\in \m H^{\alpha_Y}_{L_4}(\mb R^{d_Y})$, it holds that 
\begin{equation*}
      \Big|\mb{E}_{\mu^*_{Y|x}}[\psi(\wh Q_{\sk}(Y))\rho_\sk(x,Y)]\Big|=\Big|\int_{\mb R^{d_Y}} \psi(z) \wt v_\sk (z|x)\,\dd z\Big|\lesssim 2^{-\frac{d_Yj}{2}-j\alpha_Y}.
\end{equation*}
  Let $J=\lceil {\frac{1}{2\alpha_Y+d_Y+d_X\frac{\alpha_Y}{\alpha_X}}}\cdot \log_2 (\frac{n}{\log n})\rceil$. For $j\in \{0\}\cup [J]$, denote
\begin{equation*}
    \m S_j^{\ddagger}=\{S: \Psi_j^{d_Y}\times \mb R^{D_X}\to \mb R:\, S(\psi,x)=\sum_{\psi_1\in \Psi_j^{d_Y}} s_{\psi_1}(x),\text{ where } s_{\psi_1}\in \ms S_j \text{ for each }\psi_1\in \Psi_j^{d_Y}\},
\end{equation*}
where $\ms S_j$ is defined in~\eqref{defMSddaggerr2}. Using the independence of $\{X_i\}_{i\in I_1}$ and $\{X_i\}_{i\in I_2}$, and mirroring the analysis from the proof of Lemma~\ref{lemma1.1}---where we replace $D_Y$ with $d_Y$, and modify $\psi(Y)$ to $\psi(\wh Q_{\sk}(Y))\rho_\sk(X,Y)$. To apply Theorem~\ref{theoremjointregression}, we set $\{\psi_{\lambda}((X,Y))\}_{\lambda\in \Lambda}=\{\psi(\wh Q_{\sk}(Y))\rho_\sk(X,Y):\, \psi\in \Psi_{j}^{d_Y}\}$, where the response variable $Y$ is redefined as the joint vector of $(X,Y)$, alongside $\m S= \m S_j^{\ddagger}$ --- we can show that, by applying a union argument over $j\in [J]$ and $k\in\wh {\m K}$, it holds with probability at least $1-\frac{1}{n^2}$ that for any $j\in [J]$ and $k\in\wh {\m K}$,
      \begin{equation*} 
          \mb{E}_{\mu^*_X} \Big[\sum_{\psi\in \Psi_j^{d_Y}} (\wh v_{k\psi}(X)-  \mb{E}_{\mu^*_{Y|X}}[\psi(\wh Q_{\sk}(Y))\rho_\sk(X,Y)])^2\Big]\lesssim 2^{\frac{2j\alpha_Xd_Y}{2\alpha_X+d_X}} (\frac{n}{\log n})^{-\frac{2\alpha_X}{2\alpha_X+d_X}}.
      \end{equation*}
Thus, employing a similar strategy to that used in Regime 2, we can demonstrate that it holds with probability at least $1-\frac{1}{n^2}$ that for any $\gamma>0$,

\begin{equation*}
     \begin{aligned}
      &  (E_C)=\mb{E}_{\mu^*_X}  \Big[\underset{f\in \m H_1^{\gamma}(\mb R^{D_Y})}{\sup} \Big|\sum_{k\in \wh{\m K}} \mb{E}_{\mu^*_{Y|X}}\big[f_J^{\perp}(\wh G_\sk(\wh Q_\sk(Y),X))\rho_\sk(X,Y)\big]\\
          &\qquad\qquad\qquad-\int_{\mb R^{d_Y}} f^{\perp}_J(\wh G_\sk (z,X)) \sum_{j=0}^J \sum_{\psi\in \Psi_j^{d_Y}}\psi(z) \wh v_{k\psi}(X)\,\dd z\Big|\Big]\\
        &\lesssim  (\frac{n}{\log n})^{-\frac{\gamma}{2\alpha_Y+d_Y+d_X\frac{\alpha_Y}{\alpha_X}}}\mb{E}_{\mu^*_X} \Big[\underset{f:\mb R^{d_Y}\to \mb R\atop \int f^2(y)\,\dd y\leq 1}{\sup}  \sum_{k\in \wh{\m K}} \mb{E}_{\mu^*_{Y|X}} [f(\wh Q_{[k]}(Y))\rho_\sk(X,Y)]\\
        &\qquad\qquad-\sum_{k\in \wh{\m K}}\int_{\mb R^{d_Y}} f(z)\sum_{j=0}^J \sum_{\psi\in \Psi_j^{d_Y}}\psi(z) \wh v_{k\psi}(X)\,\dd z \Big]\\
    &\lesssim   (\frac{n}{\log n})^{-\frac{\gamma}{2\alpha_Y+d_Y+d_X\frac{\alpha_Y}{\alpha_X}}}\mb{E}_{\mu^*_X} \Big[\sum_{k\in \wh{\m K}}\underset{f:\mb R^{d_Y}\to \mb R\atop \int f^2(y)\,\dd y\leq 1}{\sup}\sqrt{\sum_{j=0}^{\infty} \sum_{\psi\in \Psi_j^{d_Y}} f_{\psi}^2 }\\&\qquad \cdot\sqrt{\sum_{j=0}^{J}\sum_{\psi\in \Psi_j^{d_Y}}\Big( \mb{E}_{\mu^*_{Y|X}}[\psi(\wh Q_\sk(Y))\rho_\sk(X,Y)]- \wh v_{k\psi}(X) \Big)^2+\sum_{j=J+1}^{\infty}\sum_{\psi\in \Psi_j^{d_Y}}\Big( \mb{E}_{\mu^*_{Y|X}}[\psi(\wh Q_\sk(Y))\rho_\sk(X,Y)] \Big)^2}\Big] \\
  &\lesssim    (\frac{n}{\log n})^{-\frac{\alpha_Y+\gamma}{2\alpha_Y+d_Y+\frac{\alpha_Y}{\alpha_X}d_X}} \\
  &\qquad\qquad+  (\frac{n}{\log n})^{-\frac{\gamma}{2\alpha_Y+d_Y+d_X\frac{\alpha_Y}{\alpha_X}}}\sum_{k\in \wh{\m K}} \sqrt{\mb{E}_{\mu^*_X}\Big[\sum_{j=0}^{J}\sum_{\psi\in \Psi_j^{d_Y}}\Big( \mb{E}_{\mu^*_{Y|X}}[\psi(\wh Q_\sk(Y))\rho_\sk(X,Y)]-\wh v_{k\psi}(X)\Big)^2\Big]}\\
  &\lesssim (\frac{n}{\log n})^{-\frac{\alpha_Y+\gamma}{2\alpha_Y+d_Y+\frac{\alpha_Y}{\alpha_X}d_X}}.
    \end{aligned}
 \end{equation*}
 By combining the bounds for term $ (E_A)$, $ (E_B)$, $(E_C)$, we can then get the desired results.

\subsection{Proof of Theorem~\ref{th2:lower1} (minimax lower bound for Regime 2)}\label{proofth2:lower1}
The upper bound is established by Theorem~\ref{th:combined} and Corollary~\ref{co1}; hence, our focus here is on establishing the lower bound.  The term $n^{-\frac{\beta_Y\gamma}{d_Y}}$ in the lower bound is directly derived from the minimax lower bound for the unconditional case as specified in Theorem 3.1 of~\cite{10.1214/23-AOS2291}.  Moreover, the lower bound for $d_X=0$ also  follows directly from the minimax rate in the unconditional case. Consequently, our analysis will concentrate on the terms $n^{-\frac{\alpha_X}{2 \alpha_X+d_X}}$ and $n^{-\frac{\alpha_Y+\gamma}{2 \alpha_Y+d_Y+\frac{\alpha_Y}{\alpha_X} d_X}}$ for $d_X\in \mb N_{+}$.  Define the covariate space $\m M_X=[-1,1]^{d_X}\times \mathbf 0_{D_X-d_X}$, with $\mu^*_X$ representing the uniform distribution over $\m M_X$. Let $\m M_0=\mb S_2^{d_Y}\times \mathbf{0}_{D_Y-d_Y-1}=\{y\in \mb R^{D_Y}:\, \|y_{1:d_Y+1}\|^2=2,\,y_{d_Y+2:D_Y}=\mathbf 0_{D_Y-d_Y-1}\}$ denote the $d_Y$-dimensional sphere embedded in $\mb R^{D_Y}$ and let $\wt{\m{M}}_0=\{y\in \mb{R}^{D_Y}:\, y_{1:d_Y}\in \mb B_{\mb R^{d_Y}}(\mathbf{0},1),\, y_{d+1}=\sqrt{2-\|y_{1:d_Y}\|^2},\, y_{d_Y+2:D_Y}=\mathbf 0_{D_Y-d_Y-1}\}$ denote the middle area of $\m M_0$.   Then $\wt{\m{M}}_0$
admits a global parametrization $G_0:\, \mb B_{\mb R^{d_Y}}(\mathbf{0},1)\to \wt{\m{M}}_0$ defined as $G_0(z)=(z,\,\sqrt{2-\|z\|_2^2}\,,\,\mathbf{0}_{D_Y-d_Y-1})$ for $z\in \mb B_{\mb R^{d_Y}}(\mathbf{0},1)$. So we can define $\nu_0$ as the density function on $ \mb B_{\mb R^{d_Y}}(\mathbf{0},1)$ so that $[G_0]_\# \nu_0$ is the normalized restriction of $\mu_0$ on $\wt{\m{M}_0}$, or
\begin{align*}
    \nu_0(z)=\frac{1}{\wt{C}}\sqrt{{\rm det}(\mathbf{J}_{G_0}(z)^T\mathbf{J}_{G_0}(z))}, \quad \forall z\in  \mb B_{\mb R^{d_Y}}(\mathbf{0},1),
\end{align*}
where  $\mathbf{J}_{G_0}$ denotes the Jacobian matrix of $G_0$ and $\wt{C}=\int_{ \mb B_{\mb R^{d_Y}}(\mathbf{0},1)}\sqrt{{\rm det}(\mathbf{J}_{G_0}(z)^T\mathbf{J}_{G_0}(z))}\,\dd z$ is the normalizing constant.
Let $\Psi_{\alpha_Y,\alpha_X}$ be a conditional density function class of $z|x$ indexed by a parameter $\omega$,  so that for any $\nu_{\omega}(z,x)\in \Psi_{\alpha_Y,\alpha_X}$ and $x\in \m M_X$, $\nu_{\omega}(z,x)=\nu_0(z)$ if  $z\notin \mb B_{\mb R^{d_Y}}(\mathbf{0},3/4)$. Then for any $\nu_{\omega}(z,x)\in \Psi_{\alpha_Y,\alpha_X}$ and $x\in \m M_X$, we define the following distribution over $\m M_0$ as
\begin{equation*}
    \mu^{\omega}_{Y|x}=\Big(1-\frac{\wt{C}}{C}\Big)\cdot\mu_1+\frac{\wt{C}}{C}\cdot [G_0]_{\#}[\nu_{\omega}(z,x)\dd z],
\end{equation*}
where $\mu_1$ represents the uniform distribution over $\wt{\m M}_1=\m M_0\setminus \wt {\m M}_0$. and $C$ is the surface area of $\mb{S}_2^{d_Y}$.
Then $\mu^{\omega}_{Y|x}$ has the following conditional density function with respect to the volume measure of $\m{M}_0$,
\begin{equation*}
\begin{aligned}
       \mu^{\omega}(y|x)&=\frac{1}{C}\, \mathbf{1}(y\in \wt{\m M}_1)+\frac{\wt{C}}{C}\cdot\frac{\nu_{\omega}(y_{1:d_Y},x)}{\sqrt{{\rm det}(\mathbf{J}_{G_0}(y_{1:d_Y})^T\mathbf{J}_{G_0}(y_{1:d_Y}))}}\cdot \mathbf{1}(y\in {\wt{\m M}}_0), \, \forall y\in \m M_0\\
       &=\left\{\begin{array}{cc}
           \frac{\wt C}{C} \frac{\nu_{\omega}(y_{1:d_Y},x)}{\sqrt{{\rm det}(\mathbf{J}_{G_0}(y_{1:d_Y})^T\mathbf{J}_{G_0}(y_{1:d_Y}))}}, &  \quad y\in \wt{\m M}_0=\{ y=(z,\sqrt{2-\|z\|^2_2},\mathbf{0}_{D_Y-d_Y-1})\,:\,\|z\|\leq 1\} \\
           \frac{1}{C}, & \quad y\in \wt{\m M}_1= {\m M}_0\setminus \wt{\m M}_0
       \end{array}
       \right.\\
        &=\left\{\begin{array}{cc}
           \frac{\wt C}{C} \frac{\nu_{\omega}(y_{1:d_Y},x)}{\sqrt{{\rm det}(\mathbf{J}_{G_0}(y_{1:d_Y})^T\mathbf{J}_{G_0}(y_{1:d_Y}))}}, &  \quad y\in \{ y=(z,\sqrt{2-\|z\|^2_2},\mathbf{0}_{D_Y-d_Y-1})\,:\,\|z\|\leq 3/4\} \\
           \frac{1}{C}, & \quad y\in \m M_0\setminus \{ y=(z,\sqrt{2-\|z\|^2_2},\mathbf{0}_{D_Y-d_Y-1})\,:\,\|z\|\leq 3/4\}.
       \end{array}
       \right.
\end{aligned}
\end{equation*}
 Moreover, we have
\begin{equation*}
    \begin{aligned}
       & d_{\gamma}(\mu^{\omega}_{Y|x},\mu^{\omega'}_{Y|x})=\frac{\wt C}{C}\underset{f\in \m H^{\gamma}_1(\mb R^{D_Y})}{\sup} \int_{\mb B_{\mb R^{d_Y}}(\mathbf{0},1)} f(G_0(z)) (\nu_{\omega}(z,x)-\nu_{\omega'}(z,x))\, \dd z\\
        &\geq \frac{\wt C}{C} \underset{f\in \m H^{\gamma}_1(\mb R^{d_Y})}{\sup} \int_{\mb B_{\mb R^{d_Y}}(\mathbf{0},1)} f(z) (\nu_{\omega}(z,x)-\nu_{\omega'}(z,x))\, \dd z,\\
    \end{aligned}
\end{equation*}
and 
\begin{equation*}
    \begin{aligned}
D_{\rm KL}(\mu^{\omega}_{Y|x},\mu^{\omega'}_{Y|x})=\frac{\wt C}{C}\int_{\mb B_{\mb R^{d_Y}}(\mathbf{0},1)} -\log \frac{\nu_{\omega'}(z,x)}{\nu_{\omega}(z,x)} \nu_{\omega}(z,x)\,\dd z.
    \end{aligned}
\end{equation*}
Therefore, selecting 
\begin{equation*} 
\begin{aligned}
\Psi_{\alpha_Y,\alpha_X}&=\Big\{ \nu_{\omega}(z,x)=\nu_0(z)+\Big(\frac{1}{\wt m_1}\Big)^{\alpha_Y}\sum_{\xi_1\in[\wt m_1]^{D_Y}}\sum_{\xi_2\in[\wt m_2]^{d_X}}  \omega_{\xi_1,\xi_2} \,\wt \psi_{\xi_1,\xi_2}(z,x)\\
&\qquad:\  \omega=\{\omega_{\xi_1,\xi_2}\}_{\xi_1\in[\wt m_1]^{d_Y},\xi_2\in[\wt m_2]^{d_X}}\in \{0,1\}^{\wt m_1^{d_Y}\times \wt m_2^{d_X}}\Big\},\\
&  \wt m_1 = \lceil\, b\,n^{\frac{1}{2\alpha_Y+d_Y+\frac{\alpha_Y}{\alpha_X}d_X}}\rceil, \quad \wt m_2 = \lceil b\, n^{\frac{1}{2\alpha_X+d_X+\frac{\alpha_X}{\alpha_Y}d_Y}}\rceil,
\end{aligned}
 \end{equation*}
 where 
 \begin{equation*}
     \wt \psi_{\xi_1,\xi_2}(y,x)=\prod_{i=1}^{d_Y} \wt k\Big(\wt m_1\sqrt{\frac{d_Y}{2}}y_i+\frac{\wt m_1}{2}-\xi_{1i}\Big)\prod_{i=1}^{d_X} \wt k\Big(\wt m_2\sqrt{2d_X}x_i-\xi_{2i}\Big),\quad\forall y\in \mb B_{\mb R^{d_Y}}(\mathbf{0},1),
 \end{equation*}
 and
\begin{equation*} 
 \wt k(t)=\left\{\begin{array}{l}
(1-t)^{\alpha_Y\vee \alpha_X\vee\gamma+1} t^{\alpha_Y\vee\alpha_X\vee\gamma+1}(t-\frac{1}{2}), \quad t \in(0,1) \\
0, \quad \text { o.w. }
\end{array}\right.
\end{equation*}
We can verify that there exists a constant $L$ so that for any $\nu_{\omega}\in \Psi_{\alpha_Y,\alpha_X}$,  the function $\ov\mu^{\omega}:\mb R^{D_Y}\times \mb R^{D_X}\to \mb R$ defined by 
\begin{equation*}
\begin{aligned}
       \ov\mu^{\omega}(y,x)& =\left\{\begin{array}{cc}
           \frac{\wt C}{C} \frac{\nu_{\omega}(y_{1:d_Y},x)}{\sqrt{{\rm det}(\mathbf{J}_{G_0}(y_{1:d_Y})^T\mathbf{J}_{G_0}(y_{1:d_Y}))}}, &  \quad \|y_{1:d_Y}\|\leq 1 \\
           \frac{1}{C}, & \quad \|y_{1:d_Y}\|> 1
       \end{array}
       \right.\\
        &=\left\{\begin{array}{cc}
           \frac{\wt C}{C} \frac{\nu_{\omega}(y_{1:d_Y},x)}{\sqrt{{\rm det}(\mathbf{J}_{G_0}(y_{1:d_Y})^T\mathbf{J}_{G_0}(y_{1:d_Y}))}}, & \quad \|y_{1:d_Y}\|\leq 3/4 \\
           \frac{1}{C}, & \quad \quad \|y_{1:d_Y}\|> 3/4.
       \end{array}
       \right.
\end{aligned}
\end{equation*}
satisfies that $\ov \mu^{\omega}\in \m H_L^{\alpha_Y,\alpha_X}(\mb R^{D_Y},\mb R^{D_X})$, and $\ov\mu^{\omega}(y,x)=\mu^{\omega}(y|x)$ holds for any $y\in \m M_0$ and $x\in \m M_X$. Therefore, let $g(\cdot)\equiv 1$, for any $\beta_Y>0$, there exist constants $\tau,\tau_1,L$ so that
\begin{equation*}
    \begin{aligned}
        &\Big\{\mu=\mu^*_X\mu_{Y|x}\,:\,\mu_{Y|x}=\Big(1-\frac{\wt{C}}{C}\Big)\cdot\mu_1+\frac{\wt{C}}{C}\cdot [G_0]_{\#}[\nu( z,x)\dd z],\quad \nu \in \Psi_{\alpha_Y,\alpha_X}\Big\}\\
        &\subset \m P^*_2(D_Y,D_X,d_Y,d_X,\beta_Y,\alpha_Y,\alpha_X,\tau,\tau_1,g,L).
    \end{aligned}
\end{equation*}
Following the same procedure as in the proof of Theorem~\ref{th1:lower} (see Appendix~\ref{proofth1:lower}), we can then get the desired lower bound  of $n^{-\frac{\alpha_Y+\gamma}{2\alpha_Y+d_Y+\frac{\alpha_Y}{\alpha_X}d_X}}$.  Similarly, to attain the desired lower bound of  $n^{-\frac{\alpha_X}{2\alpha_X+d_X}}$, we can follow the same step as in the proof of Theorem~\ref{th1:lower}, but this time opting for $\Psi_{\alpha_Y,\alpha_X}$ as
\begin{equation*} 
\begin{aligned}
\Psi_{\alpha_Y,\alpha_X}&=\Psi_{\alpha_X}=\Big\{ \nu_{\omega}(y,x)=\nu_0(y)+\Big(\frac{1}{\wt m}\Big)^{\alpha_X}\sum_{\xi\in[\wt m]^{d_X}}  \omega_{\xi} \,\wt \psi_{\xi}(x) \prod_{i=1}^{d_Y} \wt k(y_i):\, \omega=\{\omega_{\xi}\}_{\xi\in[\wt m]^{d_X}}\in \{0,1\}^{\wt m^{d_X}}\Big\},\\
&\wt m = \lceil b n^{\frac{1}{2\alpha_X+d_X}}\rceil,  \quad \wt \psi_{\xi}(x)=\prod_{i=1}^{d_X} \wt k\Big(\wt m\sqrt{2d_X}x_i-\xi_{i}\Big).
\end{aligned}
\end{equation*}

\subsection{Proof of Theorem~\ref{th2:lower} (minimax lower bound for Regime 3b)}\label{proofth2lower}
 The upper bound is established by Theorem~\ref{th:combined} and Corollary~\ref{co2}, so our focus here is solely on the lower bound. The lower bound of $n^{-\frac{\alpha_X}{2\alpha_X+d_X}}+n^{-\frac{\alpha_Y+\gamma}{2\alpha_Y+d_Y+\frac{\alpha_Y}{\alpha_X}d_X}}$ can be directly derived  from the proof of the lower bound in Theorem~\ref{th2:lower1} (see Appendix~\ref{proofth2:lower1}). So the remaining task is to  show the lower bound of $n^{-\frac{\gamma}{\frac{d_Y}{\beta_Y}+\frac{d_X}{\beta_X}}}$. Notice that when $\gamma>1$,  we can observe that
    \begin{equation*}
        \begin{aligned}
            &\frac{\gamma}{\frac{d_Y}{\beta_Y}+\frac{d_X}{\beta_X}}=\frac{\gamma(\alpha_Y+1)}{\frac{d_Y(\alpha_Y+1)}{\beta_Y}+\frac{d_X(\alpha_Y+1)}{\beta_X}}\\
            &> \frac{\alpha_Y+\gamma}{\frac{d_Y(\alpha_Y+1)}{\beta_Y}+\frac{d_X(\alpha_Y+1)}{\beta_X}}\qquad (\gamma>1)\\
            &\geq \frac{\alpha_Y+\gamma}{d_Y+\frac{d_X\alpha_Y}{\alpha_X}}\qquad (\beta_Y\geq \alpha_Y+1\quad \beta_X\geq \alpha_X+\frac{\alpha_X}{\alpha_Y})\\
            &> \frac{\alpha_Y+\gamma}{2\alpha_Y+d_Y+\frac{d_X\alpha_Y}{\alpha_X}}.
        \end{aligned}
    \end{equation*}
Hence, the term  $n^{-\frac{\gamma}{\frac{d_Y}{\beta_Y}+\frac{d_X}{\beta_X}}}$ will be dominated by  $n^{-\frac{\alpha_X}{2\alpha_X+d_X}}+n^{-\frac{\alpha_Y+\gamma}{2\alpha_Y+d_Y+\frac{\alpha_Y}{\alpha_X}d_X}}$. So here we only focus on the scenario where $\gamma\leq 1$. Define the covariate space $\m M_X=[-1,1]^{d_X}\times 0_{D_X-d_X}$ and let $\mu^*_X$ be the uniform distribution over $\m M_X$. Then any $x\in \m M_X$ can be expressed as a $d_X$-dimensional vector by removing the last $D_X-d_X$ element. So in the following, we write $x=(x_1,x_2,\cdots,x_d)$ when no ambiguity may arise. Let $\m M_0=\mb S_2^{d_Y}\times \mathbf{0}_{D_Y-d_Y-1}=\{y\in \mb R^{D_Y}:\, \|y_{1:d_Y+1}\|^2=2,\,y_{d_Y+2:D_Y}=0_{D_Y-d_Y-1}\}$ denote the $d_Y$-dimensional sphere embedded in $\mb R^{D_Y}$, with $\mu_0$ representing the uniform distribution over $\m M_0$.  Let $\wt{\m{M}}_0=\{y\in \mb{R}^{D_Y}:\, y_{1:d_Y}\in \mb B_{\mb R^{d_Y}}(\mathbf{0},1),\, y_{d+1}=\sqrt{2-\|y_{1:d_Y}\|^2},\, y_{d_Y+2:D_Y}=0_{D_Y-d_Y-1}\}$ denote the middle area of $\m M_0$.  Then $\wt{\m{M}}_0$
admits a global parametrization $G_0:\, \mb B_{\mb R^{d_Y}}(\mathbf{0},1)\to \wt{\m{M}}_0$ defined as $G_0(z)=(z,\,\sqrt{2-\|z\|_2^2}\,,\,\mathbf{0}_{D_Y-d_Y-1})$ for $z\in \mb B_{\mb R^{d_Y}}(\mathbf{0},1)$. So we can define $\nu_0$ as the density function on $ \mb B_{\mb R^{d_Y}}(\mathbf{0},1)$ so that $[G_0]_\# \nu_0$ is the normalized restriction of $\mu_0$ on $\wt{\m{M}_0}$, or
\begin{align*}
    \nu_0(z)=\frac{1}{\wt{C}}\sqrt{{\rm det}(\mathbf{J}_{G_0}(z)^T\mathbf{J}_{G_0}(z))}, \quad \forall z\in  \mb B_{\mb R^{d_Y}}(\mathbf{0},1),
\end{align*}
where  $\mathbf{J}_{G_0}$ denotes the Jacobian matrix of $G_0$ and $\wt{C}=\int_{ \mb B_{\mb R^{d_Y}}(\mathbf{0},1)}\sqrt{{\rm det}(\mathbf{J}_{G_0}(z)^T\mathbf{J}_{G_0}(z))}\,\dd z$ is the normalizing constant. Moreover,
there exist positive constants $c_1,c_2$ so that for any $z\in  \mb B_{\mb R^{d_Y}}(\mathbf{0},1)$, $c_1I_d\preccurlyeq \mathbf{J}_{G_0}(z)^T\mathbf{J}_{G_0}(z)\preccurlyeq c_2I_d$.  Next we will add small bumps to function $G_0$ to construct perturbations of $\wt{\m{M}}_0$, whose unions with the spherical cap $\wt{\m M}_1:\,=\m M_0\setminus \wt{\m{M}}_0$ form our constructed perturbed $x$-dependent manifolds.

 Let $m_1=\lceil b\,n^{\frac{1}{d_Y+d_X\frac{\beta_Y}{\beta_X}}}\rceil$ and $m_2=\lceil b\,n^{\frac{1}{d_X+d_Y\frac{\beta_X}{\beta_Y}}}\rceil$, where $b$ is a large enough constant.
Then consider  a bump function 
\begin{equation}\label{eqnk.1}
 k(t)=\left\{\begin{array}{l}
(1-t)^{\beta_Y+1} t^{\beta_Y+1}, \quad t \in(0,1), \\
0, \quad \text { o.w. }
\end{array}\right.
 \end{equation}
and for $\xi_1=(\xi_{11},\xi_{12},\cdots,\xi_{1d_Y})\in[\wt m_1]^{d_Y}$, $\xi_2=(\xi_{21},\xi_{22},\cdots,\xi_{2d_X})\in[\wt m_2]^{d_X}$,
\begin{align*}
    \psi_{\xi_1,\xi_2}(z,x)=\prod_{i=1}^{d_Y} k\Big(m_1\,\sqrt{\frac{d_Y}{2}}\,z_i+\frac{m_1}{2}-\xi_{1i}\Big)\prod_{i=1}^{d_X} k\Big(m_2\,\sqrt{\frac{d_X}{2}}\,x_i+\frac{m_2}{2}-\xi_{2i}\Big).
\end{align*}
For any $\omega=(\omega_{\xi_1,\xi_2})_{\{\xi_1\in [m_1]^{d_Y},\xi_2\in[m_2]^{d_X}\}} \in \{0,1\}^{m_1^{d_Y}\times m_2^{d_X}}$, we define the multi-bump function 
\begin{align*}
    g_{\omega}(z,x)=\sum_{\xi_1\in [m_1]^{d_Y},\xi_2\in[m_2]^{d_X}}\frac{1}{m_1^{\beta_Y}} \,\omega_{\xi_1,\xi_2} \, \psi_{\xi_1,\xi_2}(z,x),
\end{align*}
whose bumps correspond to the non-zero components of $\omega$. Finally, we define $G_\omega(z,x) = G_0(z)+(\mathbf{0}_{d_Y},\,g_{\omega}(z,x),\,\mathbf{0}_{D_Y-d_Y-1})$ as the perturbed $x$-dependent generative map parametrized by  the binary tensor $\omega$. By Lemma~\ref{lemmalowerboundsmooth}, it is straightforward to verify that there exists a constant $L$ so that  $G_\omega$ belongs to $ {\m H}^{\beta_Y,\beta_X}_{L,D_Y}(\mb B_{\mb R^{d_Y}}(\mathbf{0},1),\m M_X)$. Furthermore, by Lemma F.3 in~\cite{10.1214/23-AOS2291},  which is a two-sided version of the Varshamov-Gilbert lemma~\citep{Tsybakov2009}, there exists a subset $\{\omega^{(1)},\cdots , \omega^{(H_0)}\}\subset \{0,1\}^{m_1^{d_Y}\times m_2^{d_X}}$ such that:
\begin{enumerate} 
\item $\log H_0 \geq \frac{m_1^{d_Y}m_2^{d_X}}{8}-\log 2$;
\item for any $j,k\in [H_0]$ with $j\neq k$, the Hamming distance $\|\omega^{(j)}-\omega^{(k)}\|_{\rm H}$ between $\omega^{(j)}$ and $\omega^{(k)}$ satisfies $\frac{m_1^{d_Y}m_2^{d_X}}{4}\leq \|\omega^{(j)}-\omega^{(k)}\|_{\rm H}\leq \frac{3m_1^{d_Y}m_2^{d_X}}{4}$.
\end{enumerate}
For each $\omega\in \{0,1\}^{m_1^{d_Y}\times m_2^{d_X}}$, define $\bar{\omega}=1-\omega$ in the element-wise manner. We may expand the above $H_0$ tensors into $H=2H_0$ ones, ordered as
\begin{equation*}
\{{\omega}^{(1)},\cdots, {\omega}^{(H)}\}= \{\omega^{(1)},\cdots, \omega^{(H_0)},
\bar{\omega}^{(1)},\cdots,\bar{\omega}^{(H_0)}\}.
\end{equation*}
Then $\log H \geq \frac{m_1^{d_Y}m_2^{d_X}}{8}$ and for any $i,j\in [H]$ with $i\neq j$, it holds that $\|{\omega}^{(i)}- {\omega}^{(j)}\|_{\rm H} \geq \frac{m_1^{d_Y}m_2^{d_X}}{4}$. 

Next, for each $i\in[H]$ and $x\in \m M_X$, let ${\m{M}}^{{\omega}^{(i)}}_{Y|x}=G_{{\omega}^{(i)}}(\mb B_{\mb R^{d_Y}}(\mathbf{0},1),x)=\{G_{{\omega}^{(i)}}(z,x):\, z\in \mb B_{\mb R^{d_Y}}(\mathbf{0},1)\}$ denote the perturbed manifold from $G_{{\omega}^{(i)}}(\cdot,x)$. We define a perturbation to $\mu_0$ by smoothly gluing together the uniform distribution $\mu_1$ over $\wt{\m M}_1$ and $\mu^{\omega^{(i)}}_{Y|x}:\,= [G_{\omega^{(i)}}(\cdot,x)]_\# \nu_0$ over $\m M^{\omega^{(i)}}_{Y|x}$ as
\begin{equation*}
    \mu^{i}_{Y|x}=\Big(1-\frac{\wt{C}}{C}\Big)\cdot\mu_1
    +\frac{\wt{C}}{C}\cdot{\mu}^{{\omega}^{(i)}}_{Y|x},
\end{equation*}
where $C$ is the volume of $\m M_0$ so that $C^{-1}$ is the density function of the uniform distribution over $\m M_0$ and $ C>\wt C$. Then 
$ \mu^{i}_{Y|x}$ is supported over the manifold $\m M_{Y|x}^i:\,= \wt{\m M}_1\cup \m M^{\omega^{(i)}}_{Y|x}$. Given that $G_{\omega^{(i)}}\in   {\m H}^{\beta_Y,\beta_X}_{L,D_Y}(\mb B_{\mb R^{d_Y}}(\mathbf{0},1),\m M_X)$ and, by construction, $G_{\omega^{(i)}}(z,x)=G_0(z)$ for any $z\in \mb B_{\mb R^{d_Y}}(\mathbf{0},1)\setminus \mb B_{\mb R^{d_Y}}(\mathbf{0},\frac{3}{4})$, it follows from Lemma~\ref{le:defmanifold} that $\{\m M^{\omega^{(i)}}_{Y|x}\}_{x\in \m M_X}\in \ms M^{\beta_Y,\beta_X}_{\tau,\tau_1,L}(d_Y,D_Y,\m M_X)$ for some small enough $\tau,\tau_1$ and large enough $L$. Furthermore, the density function of distribution $   \mu^{i}_{Y|x}$ with respect to the volume measure of $\m M^{i}_{Y|x}$ is given by
\begin{equation*}
    u_i(y|x)=\frac{1}{C}\, \mathbf{1}(y\in \wt{\m M_1})+\frac{1}{C}\, \frac{\sqrt{{\rm det}(\mathbf{J}_{G_0}(y_{1:{d_Y}})^T\mathbf{J}_{G_0}(y_{1:d_Y}))}}{\sqrt{{\rm det}(\mathbf{J}_{G_{{\omega}^{(i)}}(\cdot,x)}(y_{1:d_Y})^T\mathbf{J}_{G_{{\omega}^{(i)}}(\cdot,x)}(y_{1:d_Y}))}}\, \mathbf{1}(y\in {\m{M}}^{{\omega}^{(i)}}_{Y|x}).
    \end{equation*}
Then consider the smooth transition function
 \begin{equation}\label{eqnrho1}
      \rho_a(t)=\left\{
      \begin{array}{cc}
        0   & |t|\geq a \\
        1   & |t|\leq 1\\
       \frac{1}{1+ \exp(\frac{(a+1)-2t}{(t-1)(t-a)})}& 1<t<a\\
        \frac{1}{1+ \exp(\frac{(a+1)+2t}{(t+1)(a+t)})}& -a<t<-1,\\
      \end{array}
      \right.
  \end{equation}
and define
\begin{equation*}
    \ov u_i(y,x)=\frac{1}{C}+\frac{1}{C}\Big(\frac{\sqrt{{\rm det}(\mathbf{J}_{G_0}(y_{1:{d_Y}})^T\mathbf{J}_{G_0}(y_{1:d_Y}))}}{\sqrt{{\rm det}(\mathbf{J}_{G_{{\omega}^{(i)}}(\cdot,x)}(y_{1:d_Y})^T\mathbf{J}_{G_{{\omega}^{(i)}}(\cdot,x)}(y_{1:d_Y}))}} -1\Big)\rho_{\frac{16}{9}}(\frac{\|y_{1:d_Y}\|^2}{\frac{9}{16}}),
\end{equation*}
    Note that function $\frac{\sqrt{{\rm det}(\mathbf{J}_{G_0}(y_{1:{d_Y}})^T\mathbf{J}_{G_0}(y_{1:d_Y}))}}{\sqrt{{\rm det}(\mathbf{J}_{G_{{\omega}^{(i)}}(\cdot,x)}(y_{1:d_Y})^T\mathbf{J}_{G_{{\omega}^{(i)}}(\cdot,x)}(y_{1:d_Y}))}} =1$ for $y_{1:d_Y}\in \mb B_{\mb R^{d_Y}}(\mathbf{0},1)\setminus \mb B_{\mb R^{d_Y}}(\mathbf{0},\frac{3}{4})$.
Consequently,  for any $x\in \m M_X$ and $y\in \m M^{i}_{Y|x}$, it holds that  $\ov \mu_i(y,x)=\mu_i(y|x)$, and there exists a constant $L$ such that $\ov\mu_i(y,x)\in  {\m H}^{\beta_Y-1,\beta_X-\frac{\beta_X}{\beta_Y}}_{L}(\mb R^{D_Y}, \mb R^{D_X})\subset  {\m H}^{\alpha_Y,\alpha_X}_{L}(\mb R^{D_Y}, \mb R^{D_X})$. Therefore, for any $i\in [H]$, it holds that $\mu^*_X\mu^{i}_{Y|X}\in\mathcal{P}_3^{\ast}(D_Y,D_X,d_Y,d_X,\beta_X,\beta_Y,\alpha_Y,\alpha_X,\tau,\tau_1,g,L)$, where $g(\cdot)\equiv 1$.  Then let $\bar{\mu}=\frac{1}{H} \sum_{i=1}^{H} \mu^*_X\mu^{i}_{Y|X}$ be the averaged distribution. Since for any fixed index $\xi\in [m_1]^{d_Y}\otimes [m_2]^{d_X}$, there are equal numbers of $0$'s and $1$'s in the sequence $\big({\omega}_{\xi}^{(1)},\cdots,{\omega}_{\xi}^{(H)}\big)$, we have 
\begin{equation*}
    \begin{aligned}
        &D_{\rm KL}(\mu^*_X\mu^{i}_{Y|X},\bar{\mu})=\mb{E}_{\mu^*_X}[D_{\rm KL}(\mu^{i}_{Y|x},\bar{\mu}_{Y|x})]\leq \log 2.
    \end{aligned}
\end{equation*}
Moreover, for  any pair of $j,k\in [H]$ with $j\neq k$, by construction we have $\|\omega^{(j)}- \omega^{(k)}\|_{\rm H} \geq \frac{m_1^{d_Y}m_2^{d_X}}{4}$. Define
\begin{equation*}
\widetilde{f}(z,x)=\sum_{\xi_1\in[m_1]^{d_Y}}\sum_{\xi_2\in[m_2]^{d_X}}\left(\frac{1}{m_1}\right)^{\gamma} v_{\xi_1,\xi_2} \psi_{\xi_1\xi_2}(z,x),
\end{equation*}
where
\begin{equation*}
v_{\xi_1,\xi_2}=\left\{\begin{array}{l}
1, \quad  \omega_{\xi_1,\xi_2}^{(j)}=1\text{ and }\omega_{\xi_1,\xi_2}^{(k)}=0;\text{ or }\omega_{\xi_1,\xi_2}^{(j)}=\omega_{\xi_1,\xi_2}^{(k)},\\
-1, \quad   \omega_{\xi_1,\xi_2}^{(j)}=0\text{ and }\omega_{\xi_1,\xi_2}^{(k)}=1.
\end{array}\right.
\end{equation*}
By the definition of $g_{\omega}(z,x)$, there exists a constant $c$ such that  for any $j \in[H]$ and $x\in \m M_X$,  it holds that ${\rm supp}(\mu^{\omega^{(j)}}_{Y|x})\subset  \mathbb{R}^{d_Y}\times\{y_{d+1}: |y_{d+1}-\sqrt{2-\|y_{1:d}\|^2}|\leq \frac{c}{m_1^{\beta_Y}}\}\times \{(y_{d_Y+2},\cdots,x_{D_Y})^T=\mathbf{0}_{D_Y-d_Y-1}\}$.  Define function $h: \mathbb{R}\to \mathbb{R}$ by $h(x)=\max(- \frac{c}{m_1^{\beta_Y}},\,\min(\frac{c}{m_1^{\beta_Y}},x))$, then $h$ is a 1-Lipschitz function over $\mathbb{R}$. Consider function $\chi:\mathbb{R}\to \mathbb{R}$ defined by $\chi(t)=e^{-1/t}$ for $t>0$ and $\chi(t)=0$ for $t\leq 0$. For $z\in \mathbb{R}^{d_Y}$, we define 
\begin{equation*}
    q(z)=\left\{
    \begin{array}{cc}
      \sqrt{2-\|z\|^2}\cdot\frac{\chi(5/4-\|z\|_2)}{\chi(5/4-\|z\|_2)+\chi(\|z\|_2-1)}   & \|z\|\leq \frac{5}{4} \\
       0  &  \|z\|>\frac{5}{4}.
    \end{array}
    \right.
\end{equation*}
 Note that when $z\in \mb B_{\mb R^{d_Y}}(\mathbf{0},1)$, $q(z)=\sqrt{2-\|z\|^2}$ and we multiply $\sqrt{2-\|z\|^2}$ by $\frac{\chi(5/4-\|z\|_2)}{\chi(5/4-\|z\|_2)+\chi(\|z\|_2-1)}$ to smoothly extend $\sqrt{2-\|z\|^2}$ from $\mb B_{\mb R^{d_Y}}(\mathbf{0},1)$ to the entire space. Now define 
\begin{equation*}
f(y,x)=\widetilde{f}(y_{1:d},x)\, h\big(y_{d+1}-q(y_{1:d})\big)\, m_1^{\gamma-\gamma\beta_Y+\beta_Y}.
\end{equation*}
We then prove that $f(\cdot,x)$ is $\gamma$-smooth with bounded H\"{o}lder norm. Since for any $y, y' \in \mathbb{R}^{D_Y}$, it holds that $|h\big(y_{d_Y+1}-q(y_{1:d_Y})\big)|\leq \frac{c}{m_1^{\beta_Y}}$ and $|h\big(y'_{d_Y+1}-q(y'_{1:d_Y})\big)|\leq \frac{c}{m_1^{\beta_Y}}$. Therefore, we have
\begin{equation*}
\begin{aligned}
    |h\big(y_{d_Y+1}-q(y_{1:d_Y})\big)-h\big(y'_{d_Y+1}-q(y'_{1:d_Y})\big)|&\leq \frac{(2c)^{1-\gamma} }{m^{\beta_Y(1-\gamma)}} |h\big(y_{d_Y+1}-q(y_{1:d_Y})\big)-h\big(y'_{d_Y+1}-q(y'_{1:d_Y})\big)|^{\gamma}\\
& \lesssim \frac{1}{m^{\beta_Y(1-\gamma)}}\|y-y'\|^{\gamma}.
\end{aligned}
\end{equation*}
Moreover, for any $z,z'\in \mathbb{R}^{d_Y}$, there exists a constant $c_1$ such that 
\begin{equation*}
|\widetilde{f}(z,x)-\widetilde{f}(z',x)|\leq c_1 \frac{1}{m_1^{\gamma-1}}\|z-z'\|.
\end{equation*}
Therefore, in the case $\|z-z'\|\leq \frac{1}{m_1}$, we have $\|z-z'\|\leq \frac{1}{m_1^{1-\gamma}}\|z-z'\|_2^\gamma$, and thus $|\widetilde{f}(z,x)-\widetilde{f}(z,x)|\leq c_1 \|z-z'\|_2^\gamma$; in the case $\|z-z'\|_2> \frac{1}{m_1}$, since there exists a constant $c_2$ such that $\sup_{z\in \mathbb{R}^d}|\widetilde{f}(z,x)|\leq \frac{c_2}{m_1^{\gamma}}$, it holds that $|\widetilde{f}(z,x)-\widetilde{f}(z,x)|\leq 2c_2 \|z-z'\|_2^{\gamma}$. Putting pieces together, we have that for any $y, y' \in \mathbb{R}^{D_Y}$ and $x\in \mb R^{D_X}$, there exist constants $c_3,c_4$ such that 
\begin{equation*}
\begin{aligned}
|f(y,x)-f(y',x)|&\leq m_1^{\gamma-\gamma\beta_Y+\beta_Y}\bigg(\Big|\widetilde{f}(y_{1:d_Y},x) \cdot \Big(h\big(y_{d_Y+1}-q(y_{1:d_Y})\big)-h\big(y'_{d_Y+1}-q(y'_{1:d_Y})\big)\Big)\Big|\\
&\qquad\qquad+\Big|h\big(y'_{d_Y+1}-q(y'_{1:d_Y})\big) \cdot (\widetilde{f}(y_{1:d_Y},x)-\widetilde{f}(y_{1:d_Y}',x))\Big|\bigg)\\
&\leq c_3\big( \|y-y'\|^{\gamma}+ m_1^{\gamma(1-\beta_Y)}\|y-y'\|^{\gamma}\big)\\
&\leq c_4 \|y-y'\|^{\gamma},
\end{aligned}
\end{equation*}
where the last inequality is due to $\beta_Y>1$. Consequently, we have for any $x\in \m M_X$, $\frac{1}{c_4}f(\cdot,x)\in \m H^{\gamma}_1(\mathbb{R}^D_Y)$ (recall that we only consider $\gamma<1$). Then
\begin{equation*}
\begin{aligned}
&\mb{E}_{\mu^*_X}[d_{\gamma}(\mu^{j}_{Y|x}, \mu^{k}_{Y|x})]\\
&\geq \frac{1}{c_4}\cdot\mb{E}_{\mu^*_X}\Big[\int f(y,x) \,\dd\mu^{j}_{Y|x}-\int f(y,x) \,\dd\mu^{k}_{Y|x}\Big]\\
&=\frac{\wt{C}}{c_1\cdot C\, }m_1^{\gamma-\gamma\beta_Y+\beta_Y}\mb{E}_{\mu^*_X}\Big[ \int_{\mb R^{d_Y}} \widetilde{f}(z,x)\cdot \big(g_{\omega^{(j)}}(z,x)-g_{\omega^{(k)}}(z,x)\big) \,\nu_0(z)\,\dd z\Big]\\
&=\frac{\wt{C}}{c_1\cdot C}\int_{[0,1]^{d_X}} m_1^{-\gamma\beta_Y} \int_{\mb B_{\mb R^{d_Y}}(\mathbf{0},1)}\sum_{\xi_1\in[m_1]^{d_Y}}\sum_{\xi_2\in[m_2]^{d_X}}v_{\xi_1,\xi_2} \psi_{\xi_1,\xi_2}(z,x) \\
&\qquad\qquad\cdot\sum_{\xi_1\in[m_1]^{d_Y}}\sum_{\xi_2\in[m_2]^{d_X}} \big(\omega_{\xi_1,\xi_2}^{(j)}-\omega_{\xi_1,\xi_2}^{(k)}\big)\,\psi_{\xi_1\xi_2}(z,x)\, \nu_0(z)\, \dd z\dd x \\
&=\frac{\wt{C}}{c_1\cdot C}\int_{[0,1]^{d_X}} m_1^{-\gamma\beta_Y} \int_{\mb B_{\mb R^{d_Y}}(\mathbf{0},1)}\sum_{\xi_1\in[m_1]^{d_Y}}\sum_{\xi_2\in[m_2]^{d_X}}v_{\xi_1,\xi_2} \big(\omega_{\xi_1,\xi_2}^{(j)}-\omega_{\xi_1,\xi_2}^{(k)}\big)\psi^2_{\xi_1,\xi_2}(z,x) \nu_0(z)  \dd z\dd x \\
&\gtrsim  m_1^{-\gamma\beta_Y}.
 \end{aligned}
\end{equation*}
Then similarly to the proof of Theorem~\ref{th1:lower}, we can apply Fano's lemma to obtain
\begin{equation*}
    \begin{aligned}
   & \underset{\wh{\mu}_{Y|X}}{\inf}\underset{\mu\in \mathcal{P}_2^{\ast}} {\sup} \,\mb{E}_{\mu^{\otimes} n}\mathbb{E}_{\mu_X} \big[d_{\gamma}(\wh{\mu}_{Y|X}, \mu_{Y|X})\big]\\
    & \geq   \,  \frac{1}{2}\,\underset{j,k\in [H]\atop j\neq k} {\inf} \mb{E}_{\mu^*_X}[d_{\gamma}(\mu^{j}_{Y|X},\mu^{k}_{Y|X})] \cdot\bigg(1-\frac{\log 2+\frac{n}{H^2} \sum_{j=1}^{H}D_{\rm KL}(  \mu^*_X\mu^j_{Y|X} ,\bar\mu)}{\log H}\bigg)\\
   &\gtrsim  n^{-\frac{\gamma}{d_Y/\beta_Y+d_X/\beta_X}}.
\end{aligned}
\end{equation*}

 \subsection{Proof of Corollary~\ref{co1} and Corollary~\ref{co2}}\label{proofco1}
 We will show  Corollary~\ref{co1} here. The proof of Corollary~\ref{co2} follows the same approach. Note that for any $x\in \m M_X$, since $\mu^*_{Y|x}\in \m P_Y^*$, it holds that
  \begin{equation*}
        \begin{aligned}
 \sum_{\gamma\in \Gamma} \underset{f\in \m H^{\gamma}_1(\mb R^{D_Y})}{\sup}\frac{1}{\delta_{n,\gamma}} \Big[\mb{E}_{\wh \mu_{Y|x}}[f(Y)]- \wh{\m J}(f,x)\Big]\leq   \sum_{\gamma\in \Gamma}  \underset{f\in \m H^{\gamma}_1(\mb R^{D_Y})}{\sup} \frac{1}{\delta_{n,\gamma}}\Big[\mb{E}_{\mu^*_{Y|x}}[f(Y)]- \wh{\m J}(f,x)\Big].
        \end{aligned}
    \end{equation*}
Therefore,  we have 
\begin{equation*}
    \begin{aligned}
        &\sup_{\gamma\in \Gamma}  \mb{E}_{\mu^*_X}\Big[\underset{f\in \m H^{\gamma}_1(\mb R^{D_Y})}{\sup} \frac{1}{\delta_{n,\gamma}}\Big( \mb{E}_{\wh \mu_{Y|X}}[f(Y)]- \wh{\m J}(f,X)\Big) \Big]\\
         &\leq \mb{E}_{\mu^*_X}\Big[\sum_{\gamma\in \Gamma} \underset{f\in \m H^{\gamma}_1(\mb R^{D_Y})}{\sup}\frac{1}{\delta_{n,\gamma}}\Big( \mb{E}_{\wh \mu_{Y|X}}[f(Y)]- \wh{\m J}(f,X)\Big) \Big]\\
          &\leq \mb{E}_{\mu^*_X}\Big[\sum_{\gamma\in \Gamma} \underset{f\in \m H^{\gamma}_1(\mb R^{D_Y})}{\sup}\frac{1}{\delta_{n,\gamma}}\Big( \mb{E}_{ \mu^*_{Y|X}}[f(Y)]- \wh{\m J}(f,X)\Big) \Big]\\
          &= \sum_{\gamma\in \Gamma} \mb{E}_{\mu^*_X}\Big[ \underset{f\in \m H^{\gamma}_1(\mb R^{D_Y})}{\sup}\frac{1}{\delta_{n,\gamma}}\Big( \mb{E}_{ \mu^*_{Y|X}}[f(Y)]- \wh{\m J}(f,X)\Big) \Big]\\
          &\lesssim (\log n)\cdot \sup_{\gamma\in \Gamma} \mb{E}_{\mu^*_X}\Big[ \underset{f\in \m H^{\gamma}_1(\mb R^{D_Y})}{\sup}\frac{1}{\delta_{n,\gamma}}\Big( \mb{E}_{ \mu^*_{Y|X}}[f(Y)]- \wh{\m J}(f,X)\Big) \Big].
    \end{aligned}
\end{equation*}
Furthermore,  by Theorem~\ref{th:combined},  it holds with probability at least $1-\frac{1}{n}$ that 
\begin{equation*}
    \begin{aligned}
   \sup_{\gamma\in \Gamma} \frac{1}{\delta_{n,\gamma}} \mb{E}_{\mu^*_X}\Big[ \underset{f\in \m H^{\gamma}_1(\mb R^{D_Y})}{\sup} \Big(\mb{E}_{\mu^*_{Y|x}}[f(Y)]- \wh{\m J}(f,x)\Big)\Big] \leq 1.
        \end{aligned}
    \end{equation*}
Therefore,
    \begin{equation*}
        \begin{aligned}
 &\sup_{\gamma\in \Gamma}  \mb{E}_{\mu^*_X}\Big[\underset{f\in \m H^{\gamma}_1(\mb R^{D_Y})}{\sup} \frac{1}{\delta_{n,\gamma}}\Big( \mb{E}_{\wh \mu_{Y|X}}[f(Y)]- \mb{E}_{ \mu^*_{Y|X}}[f(Y)]\Big) \Big]\\
 &\leq   \sup_{\gamma\in \Gamma}  \mb{E}_{\mu^*_X}\Big[\underset{f\in \m H^{\gamma}_1(\mb R^{D_Y})}{\sup} \frac{1}{\delta_{n,\gamma}}\Big( \mb{E}_{\wh \mu_{Y|X}}[f(y)]- \wh {\m J}(f,X)\Big) \Big]\\
 &\qquad+ \sup_{\gamma\in \Gamma}  \mb{E}_{\mu^*_X}\Big[\underset{f\in \m H^{\gamma}_1(\mb R^{D_Y})}{\sup} \frac{1}{\delta_{n,\gamma}}\Big( \mb{E}_{ \mu^*_{Y|X}}[f(Y)]- \wh {\m J}(f,X)\Big) \Big]
\lesssim \log n.
        \end{aligned}
    \end{equation*}
    Then for any $\gamma>0$, if $\gamma<\frac{1}{\log n}$, then $$\mb{E}_{\mu^*_X}\Big[ \underset{f\in \m H^{\gamma}_1(\mb R^{D_Y})}{\sup} \Big(\mb{E}_{\wh \mu_{Y|x}}[f(Y)]- \mb{E}_{\mu^*_{Y|x}}[f(Y)]\Big)\Big]\leq 2= 2\exp(\frac{\beta_Y}{d_Y}) n^{-\frac{\frac{1}{\log n}}{\frac{d_Y}{\beta_Y}}}\leq 2\exp(\frac{\beta_Y}{d_Y}) n^{-\frac{\gamma}{\frac{d_Y}{\beta_Y}}}.$$
    If $\frac{1}{\log n}\leq\gamma\leq \frac{d_Y\alpha_X}{2\alpha_X+d_X}$, then there exists $k\in [s]$, so that $\frac{k}{\log n}\leq\gamma\leq \frac{k+1}{\log n}$, thus 
    \begin{equation*}
        \begin{aligned}
            \mb{E}_{\mu^*_X}\Big[ \underset{f\in \m H^{\gamma}_1(\mb R^{D_Y})}{\sup} \Big(\mb{E}_{\wh \mu_{Y|x}}[f(y)]- \mb{E}_{\mu^*_{Y|x}}[f(Y)]\Big)\Big]&\leq        \mb{E}_{\mu^*_X}\Big[ \underset{f\in \m H^{\frac{k}{\log n}}_1(\mb R^{D_Y})}{\sup} \Big(\mb{E}_{\wh \mu_{Y|x}}[f(y)]- \mb{E}_{\mu^*_{Y|x}}[f(Y)]\Big)\Big]\\
            &\lesssim \log n\cdot\delta_{n,\frac{k}{\log n}}\asymp \log n\cdot\delta_{n,\gamma}.
        \end{aligned}
    \end{equation*}
    If $\gamma>\frac{d_Y\alpha_X}{2\alpha_X+d_X}$, then 
     \begin{equation*}
        \begin{aligned}
           \mb{E}_{\mu^*_X}\Big[  \underset{f\in \m H^{\gamma}_1(\mb R^{D_Y})}{\sup} \Big(\mb{E}_{\wh \mu_{Y|x}}[f(y)]- \mb{E}_{\mu^*_{Y|x}}[f(Y)]\Big)\Big]&\leq        \mb{E}_{\mu^*_X}\Big[ \underset{f\in \m H^{\frac{d_Y\alpha_X}{2\alpha_X+d_X}}_1(\mb R^{D_Y})}{\sup} \Big(\mb{E}_{\wh \mu_{Y|x}}[f(y)]- \mb{E}_{\mu^*_{Y|x}}[f(Y)]\Big)\Big]\\
            &\lesssim \log n \cdot\delta_{n,\frac{d_Y\alpha_X}{2\alpha_X+d_X}}\asymp  (\log n)^4\cdot n^{-\frac{\alpha_X}{2\alpha_X+d_X}}.
        \end{aligned}
    \end{equation*}
    Proof is completed.

 \subsection{Proof of Lemma~\ref{lemmaregime2MSE}}\label{proof:lemmaregime2MSE}
We will first show an oracle inequality for the estimator 
 \begin{equation*}
 \begin{aligned}
         \wh S_j^{\dagger}&={\arg\min}_{S\in  {\m S}_j^\dagger} \frac{1}{n}\sum_{i=1}^n \sum_{\psi\in \Psi_j^{D_Y}} (2^{\frac{j(d_Y-D_Y)}{2}}\psi(Y_i)-S(\psi,X_i))^2 ,\quad  j\in \{0\}\cup [J] \text{ with }J=\lceil \frac{1}{d_Y}\cdot \log_2 (\frac{n}{\log n})\rceil,
 \end{aligned}
\end{equation*}
with a general choice of ${\m S}_j^\dagger$. For $S,S'\in {\m S}_j^\dagger$, we denote 
\begin{equation*}
   d_S(S,S')=\underset{x\in \m M_X}{\sup} \sqrt{\sum_{\psi\in \Psi_j^{D_Y}}(S(\psi,x)-S'(\psi,x))^2},
\end{equation*}
and let $\mathbf{N}({\m S}_j^\dagger,d_S,\varepsilon)$ denote the $\varepsilon$-covering number of ${\m S}_j^\dagger$ under the pseudo-distance  $d_S$.

\begin{lemma}\label{lestimation:2}
Suppose  $\big\{(X_i,Y_i)\big\}_{i=1}^n$ are $n$ i.i.d data from $\mu^*=\mu^*_X\mu^*_{Y|X}$, and the following assumptions are satisfied: (1) for any $x\in \m M_X={\rm supp}(\mu^*_X)$, $\mu^*_{Y|x}$ supported on a submanifold, denoted as $\m M_{Y|x}$, and  has a density function $u^*(\,\cdot\,|\,x)$ with respect to the volume measure of $\m M_{Y|x}$, and there exist constants $\beta_Y\geq 2, \beta_X,\alpha_Y,\alpha_X>0$ and a function $\ov u^*\in    \ov{\m H}^{\alpha_Y,\alpha_X}_{L}(\mb R^{D_Y},\mb R^{D_X})$ so that $\{\m M_{Y|x}\}_{x\in \m M_X}\in \ms M^{\beta_Y\,\beta_X}_{\tau,\tau_1,L}(d_Y,D_Y,\m M_X)$ and $\ov u^*(y,x)=u(y|x)$ for any $(x,y)\in \m M$;  (2) there exists a constant $C$ so that for any $x\in \m M_X$,  $j\in \{0\}\cup [J]$ and $S\in {\m S}_j^\dagger$,
\begin{equation*}
 \underset{\psi\in \Psi_j^{D_Y}}{\sup} |S(\psi,x)|\leq C\, 2^{-\frac{d_Yj}{2}}\text{ and }  \sum_{\psi\in \Psi_j^{D_Y}} \mathbf{1}(S(\psi,x)\neq 0)\leq C \,2^{d_Y j},
\end{equation*}
 and  $\log \mathbf{N}({\m S}_j^\dagger,d_S,\varepsilon)\leq \m W_j\log(\frac{n}{\varepsilon})$ for any $\varepsilon<\sup_{S,S'\in \m S_j^\dagger}d_S(S,S')$. Then it holds with probability at least $1-\frac{1}{n^2}$ that for any $j \in [J]$,
\begin{equation*}
\begin{aligned}
  &\mb{E}_{\mu^*_X} \Big[\sum_{\psi\in \Psi_j^{D_Y}} \big(\mb{E}_{\mu^*_{Y|X}}[2^{\frac{j(d_Y-D_Y)}{2}}\psi(Y)]-\wh S_j^\dagger(\psi,X)\big)^2\Big]
   \lesssim   \frac{\log n}{n} \m W_j \\
 &\qquad\qquad+\underset{S\in \m S_j^\dagger}{\min} \, \mb{E}_{\mu^*_X}\Big[ \sum_{\psi\in \Psi_j^{D_Y}} ( \mb{E}_{\mu^*_{Y|X}}[2^{\frac{j(d_Y-D_Y)}{2}}\psi(Y)]-S(\psi,X))^2\Big].\\
\end{aligned}
\end{equation*}
\end{lemma}

\noindent The proof of Lemma~\ref{lestimation:2} is provided in Appendix~\ref{proof:lestimation:2}. Then for the family $\m S^{\dagger}_j$ defined as
\begin{equation*}
\begin{aligned}
     \m S_j^\dagger&=\Bigg\{S(\psi,x)=\frac{\sum_{i_1=1}^{W_j}\sum_{i_2=1}^{W_j'} \sum_{k\in \mb N_0^{D_X}, |k|<\alpha_X}a_{i_1i_2k}(x-b_{i_2})^k \rho\big(\frac{\|x-b_{i_2}\|}{\varepsilon^x_j}\big)\rho\big(\frac{\|\m I_j(\psi)-e_{i_1}\|}{\varepsilon_j^y}\big)}{\sum_{i_1=1}^{W_j}\sum_{i_2=1}^{W_j'} \rho\big(\frac{\|x-b_{i_2}\|}{\varepsilon^x_j}\big)\rho\big(\frac{\|\m I_j(\psi)-e_{i_1}\|}{\varepsilon_j^y}\big)+\frac{1}{n}}:\, \\
     & \qquad\qquad\text{for any } i_1\in [W_j],i_2\in [W_j'], \text{ and } k\in \mb N_0^{D_X} \text{ with } |k|<\alpha_X,\\
     &\qquad\qquad\quad\text{ it holds that }
     b_{i_2}\in \mb B_{\mb R^{D_X}}(\mathbf{0},L),  a_{i_1i_2k}\in [-\frac{C}{2^{d_Yj/2}},\frac{C}{2^{d_Yj/2}}], e_{i_1}\in [0,1]^{D_Y+1} \Bigg\},
\end{aligned}
\end{equation*}
where $\varepsilon_j^x= 2^{\frac{jd_Y}{2\alpha_X+d_X}}(\frac{n}{\log n})^{-\frac{1}{2\alpha_X+d_X}}$,  $ \varepsilon_j^y=\frac{2^{-j}}{C_1}$, $W_j'=C_2\, (\varepsilon_j^x)^{-d_X}$,  $W_j=C_3\, (\varepsilon_j^y)^{-d_Y}$. It holds for any $S\in S_j^\dagger$ that 
\begin{equation*}
    \underset{\psi\in \Psi_j^{D_Y}}{\sup}\underset{x\in \m M_X}{\sup} |S(\psi,x)|\leq \underset{x\in \m M_X}{\sup} \sum_{k\in \mb N_0^{D_X}, |k|<\alpha_X}|a_{i_1i_2k}|\cdot\|x-b_{i_2}\|^k \lesssim 2^{-\frac{d_Yj}{2}}.
\end{equation*}
Moreover, since for any $\psi,\psi'\in \Psi_j^{D_Y}$ with $\psi\neq \psi'$, it holds that $\|\m I_j(\psi)-\m I_j(\psi')\|>c\, 2^{-j}$.  If $\varepsilon_j^y\leq \frac{c}{4}2^{-j}$, then for any $e\in [0,1]^{D_Y+1}$, there are at least one $\psi\in \Psi_j^{D_Y}$ so that $\rho\big(\frac{\|\m I_j(\psi)-e\|}{\varepsilon_j^y}\big)\neq 0$. Therefore, there are at least $W_j=\m O(2^{jd_Y})$ number of $\psi\in \Psi_j^{D_Y}$ so that $S(\psi,x)\neq 0$.  So
\begin{equation*}
     \sum_{\psi\in \Psi_j^{D_Y}} \mathbf{1}(S(\psi,x)\neq 0)\leq C \,2^{d_Y j}.
\end{equation*}
Furthermore, consider
 \begin{equation*}
     S(\psi,x)=\frac{\sum_{i_1=1}^{W_j}\sum_{i_2=1}^{W_j'} \sum_{k\in \mb N_0^{D_X}, |k|<\alpha_X}a_{i_1i_2k}(x-b_{i_2})^k \rho\big(\frac{\|x-b_{i_2}\|}{\varepsilon^x_j}\big)\rho\big(\frac{\|\m I_j(\psi)-e_{i_1}\|}{\varepsilon_j^y}\big)}{\sum_{i_1=1}^{W_j}\sum_{i_2=1}^{W_j'} \rho\big(\frac{\|x-b_{i_2}\|}{\varepsilon^x_j}\big)\rho\big(\frac{\|\m I_j(\psi)-e_{i_1}\|}{\varepsilon_j^y}\big)+\frac{1}{n}}\in S_j^\dagger
 \end{equation*}
and
 \begin{equation*}
     S'(\psi,x)=\frac{\sum_{i_1=1}^{W_j}\sum_{i_2=1}^{W_j'} \sum_{k\in \mb N_0^{D_X}, |k|<\alpha_X}a'_{i_1i_2k}(x-b'_{i_2})^k \rho\big(\frac{\|x-b'_{i_2}\|}{\varepsilon^x_j}\big)\rho\big(\frac{\|\m I_j(\psi)-e'_{i_1}\|}{\varepsilon_j^y}\big)}{\sum_{i_1=1}^{W_j}\sum_{i_2=1}^{W_j'} \rho\big(\frac{\|x-b'
     _{i_2}\|}{\varepsilon^x_j}\big)\rho\big(\frac{\|\m I_j(\psi)-e'
     _{i_1}\|}{\varepsilon_j^y}\big)+\frac{1}{n}}\in S_j^\dagger.
 \end{equation*}
 It holds that for any $\psi\in \Psi_j^{D_Y}$ and $x\in \m M_X$,  
 \begin{equation*}
    \begin{aligned}
        &|S(\psi,x)-S'(\psi,x)|\\
        &\lesssim    \underset{i_1\in [W_j]}{\sup}\underset{i_2\in [W_j']}{\sup}\sum_{k\in \mb N_0^{D_X}, |k|<\alpha_X} |a_{i_1i_2k}-a_{i_1i_2k}'|+\frac{W_jn}{\varepsilon^x_j}\sum_{i_2=1}^{W_j'}\|b_{i_2}-b_{i_2}'\|+\frac{W_j'n}{\varepsilon^y_j}\sum_{i_1=1}^{W_j}\|e_{i_1}-e_{i_1}'\|.
    \end{aligned}
\end{equation*}
So 
\begin{equation*}
\begin{aligned}
      d_S(S,S')&=\underset{x\in \m M_X}{\sup} \sqrt{\sum_{\psi\in \Psi_j^{D_Y}}(S(\psi,x)-S'(\psi,x))^2}\\
      & \leq \underset{x\in \m M_X}{\sup} \sqrt{\sum_{\psi\in \Psi_j^{D_Y}}(S(\psi,x)-S'(\psi,x))^2 (\mathbf{1}(S(\psi,x)\neq 0)+\mathbf{1}(S'(\psi,x)\neq 0))}\\
      &\lesssim 2^{\frac{jd_Y}{2}}\cdot\Big(\underset{i_1\in [W_j]}{\sup}\underset{i_2\in [W_j']}{\sup}\sum_{k\in \mb N_0^{D_X}, |k|<\alpha_X} |a_{i_1i_2k}-a_{i_1i_2k}'|+\frac{W_jn}{\varepsilon^x_j}\sum_{i_2=1}^{W_j'}\|b_{i_2}-b_{i_2}'\|+\frac{W_j'n}{\varepsilon^y_j}\sum_{i_1=1}^{W_j}\|e_{i_1}-e_{i_1}'\|\Big).
\end{aligned}
\end{equation*}
Then, using the fact that the $\varepsilon$-covering number of a $d$-dimensional ball with radius $R$ is being bounded by $(\frac{3R}{\varepsilon})^d$, we have 
 \begin{equation*}
    \log \mathbf{N}(\m S_j^\dagger,d_S,\varepsilon)\lesssim  W_jW_j' \log \frac{n}{\varepsilon} \lesssim 2^{jd_Y} (\varepsilon_j^x)^{-d_X}\log \frac{n}{\varepsilon}.
 \end{equation*}
Now we bound the approximation error.  Let ${\rm vol}_{\m M}$ denote the volume measure of $\m M$ and let $\ov u^*\in \m H_L^{\alpha_Y,\alpha_X}(\mb R^{D_Y},\mb R^{D_X})$ be a smooth extension of $u^*$. We have
\begin{equation*}
    u^*_{\psi}(x)=\int 2^{\frac{j(d_Y-D_Y)}{2}}\psi(y) \ov u^*(y\,|\,x)\,\dd {\rm vol}_{\m M_Y}\in \m H^{\alpha_X}_{L_12^{- {d_Yj}/{2}}}(\mb R^{D_X}), 
    \end{equation*}
where we have used the fact that
\begin{equation*}
    \int 2^{\frac{j(d_Y-D_Y)}{2}}|\psi(y)|\,\dd {\rm vol}_{\m M_Y}\lesssim 2^{\frac{jd_Y}{2}}\int \mathbf{1}(\psi(y)\neq 0)  \,\dd {\rm vol}_{\m M_Y}\lesssim 2^{-\frac{d_Yj}{2}}.
\end{equation*}
   Let $\m N^x_{\varepsilon^x_j}$ denote the largest $\varepsilon^x_j$-packing set of $\m M_X$, then its cardinality satisfies $|\m N^x_{\varepsilon^x_j}|\leq C_2(\varepsilon^x_j)^{-d_X}=W_j'$ when $C_2$ is large enough.  Then we define a set $ \ov{\m N}^x_{\varepsilon^x_j}= {\m N}^x_{\varepsilon^x_j}\cup \m X$, where $\m X$ is an arbitrary subset of $\m M_X\setminus {\m N}^x_{\varepsilon^x_j}$ with $|\m X|=W_j'-| {\m N}^x_{\varepsilon^x_j}|$. Denote
   \begin{equation*}
     \Psi_j^*:=\{\psi\in  \ov \Psi_j^{D_Y}:\, {\rm supp}(\psi)\cap \m M_Y\neq \emptyset\},
\end{equation*}
it holds that $\Psi_j^*\subset  \Psi_j^{D_Y}$ and $|\Psi_j^*|\leq C_3(\varepsilon_j^y)^{-d_Y}=W_j$ when $C_3$ is large enough. Moreover, define $\ov\Psi_j^*=\Psi_j^*\cup \Phi_j$,  where $\Phi_j$  is an arbitrary subset of $\Psi_j^{D_Y}\setminus\Psi_j^*$ with $|\Phi_j|=W_j-|\Psi_j^*|$. For any $\psi\in \Psi_j^{D_Y}$, we define
\begin{equation*}
\begin{aligned}
     \wt u_{\psi}(x)=\left\{
     \begin{array}{cc}
       \frac{\sum_{\wt x\in \ov{\m N}^x_{\varepsilon^x_j}} \sum_{k\in \mb N_0^{D_X}, |k|<\alpha_X}u^*_{\psi}{}^{(k)}(\wt x) (x-\wt x)^k \rho(\frac{\|x-\wt x\|}{\varepsilon^x_j})}{\sum_{\wt x\in \ov{\m N}^x_{\varepsilon^x_j}}  \rho(\frac{\|x-\wt x\|}{\varepsilon^x_j})},    &  \psi\in \Psi_j^* \\
        0  &  o.w.
     \end{array}
     \right.
\end{aligned}
\end{equation*}
and 
\begin{equation*}
\begin{aligned}
  S^*_j(\psi,x)&=\frac{\sum_{\wt x\in \ov{\m N}^x_{\varepsilon^x_j}} \sum_{k\in \mb N_0^{D_X}, |k|<\alpha_X}u^*_{\psi}{}^{(k)}(\wt x) (x-\wt x)^k \rho(\frac{\|x-\wt x\|}{\varepsilon^x_j})}{\sum_{\wt x\in \ov{\m N}^x_{\varepsilon^x_j}}  \rho(\frac{\|x-\wt x\|}{\varepsilon^x_j})+\frac{1}{n}}\\
  &=\frac{\sum_{\psi_1\in \ov\Psi_j^*}\sum_{\wt x\in \ov{\m N}^x_{\varepsilon^x_j}} \sum_{k\in \mb N_0^{D_X}, |k|<\alpha_X}u^*_{\psi_1}{}^{(k)}(\wt x) (x-\wt x)^k \rho(\frac{\|x-\wt x\|}{\varepsilon^x_j})\mathbf{1}(\psi=\psi_1)}{\sum_{\wt x\in \ov{\m N}^x_{\varepsilon^x_j}}  \rho(\frac{\|x-\wt x\|}{\varepsilon^x_j})+\frac{1}{n}}\\
    &=\frac{\sum_{\psi_1\in \ov\Psi_j^*}\sum_{\wt x\in \ov{\m N}^x_{\varepsilon^x_j}} \sum_{k\in \mb N_0^{D_X}, |k|<\alpha_X}u^*_{\psi_1}{}^{(k)}(\wt x) (x-\wt x)^k \rho(\frac{\|x-\wt x\|}{\varepsilon^x_j})\mathbf{1}(\psi=\psi_1)}{\sum_{\psi_1\in \ov\Psi_j^*}\sum_{\wt x\in \ov{\m N}^x_{\varepsilon^x_j}}  \rho(\frac{\|x-\wt x\|}{\varepsilon^x_j})\mathbf{1}(\psi=\psi_1)+\frac{1}{n}}\\
   & =\frac{\sum_{\psi_1\in \ov\Psi_j^*}\sum_{\wt x\in \ov{\m N}^x_{\varepsilon^x_j}} \sum_{k\in \mb N_0^{D_X}, |k|<\alpha_X}u^*_{\psi_1}{}^{(k)}(\wt x) (x-\wt x)^k \rho(\frac{\|x-\wt x\|}{\varepsilon^x_j})\rho(\frac{\|\m I_j(\psi)-\m I_j(\psi_1)\|}{\varepsilon_j^y})}{\sum_{\psi_1\in \ov\Psi_j^*}\sum_{\wt x\in \ov{\m N}^x_{\varepsilon^x_j}}  \rho(\frac{\|x-\wt x\|}{\varepsilon^x_j})\rho(\frac{\|\m I_j(\psi)-\m I_j(\psi_1)\|}{\varepsilon_j^y})+\frac{1}{n}}.
\end{aligned}
\end{equation*}
It holds that $S^*_j(\psi,x)\in \m S_j^\dagger$,
Moreover, for any   $\psi\in \Psi_j^{D_Y}\setminus\Psi_j^{*}$, it holds that $ u^*_{\psi}(\cdot)\equiv 0$, and therefore $S_j^*(\psi,\cdot)=\mu^*_{\psi}(\cdot)\equiv 0$. Moreover, for any $x\in \m M_X$ and $\psi\in \Psi_j^*$, we have 
\begin{equation*}
    \begin{aligned}
       | \wt u_{\psi}(x)- S^*_j(\psi,x)|&=  \frac{|\sum_{\wt x\in \ov{\m N}^x_{\varepsilon^x_j}} \sum_{k\in \mb N_0^{D_X}, |k|<\alpha_X}u^*_{\psi}{}^{(k)}(\wt x) (x-\wt x)^k \rho(\frac{\|x-\wt x\|}{\varepsilon^x_j})|}{n \cdot(\sum_{\wt x\in \ov{\m N}^x_{\varepsilon^x_j}}  \rho(\frac{\|x-\wt x\|}{\varepsilon^x_j})+\frac{1}{n})(\sum_{\wt x\in \ov{\m N}^x_{\varepsilon^x_j}}  \rho(\frac{\|x-\wt x\|}{\varepsilon^x_j}))}\\
       &\leq \frac{1}{n}\frac{\sum_{\wt x\in \ov{\m N}^x_{\varepsilon^x_j}} \sum_{k\in \mb N_0^{D_X}, |k|<\alpha_X}|u^*_{\psi}{}^{(k)}(\wt x) (x-\wt x)^k |\rho(\frac{\|x-\wt x\|}{\varepsilon^x_j})}{\sum_{\wt x\in \ov{\m N}^x_{\varepsilon^x_j}}  \rho(\frac{\|x-\wt x\|}{\varepsilon^x_j})}\\
       &\leq \frac{1}{n}\cdot\underset{\wt x\in \ov{\m N}^x_{\varepsilon^x_j}, x\in \m M_X}{\sup}\sum_{k\in \mb N_0^{D_X}, |k|<\alpha_X}|u^*_{\psi}{}^{(k)}(\wt x) (x-\wt x)^k |\\
       &\lesssim    2^{-\frac{d_Yj}{2}}n^{-1};
    \end{aligned}
\end{equation*}
\begin{equation*}
    \begin{aligned}
          | \wt u_{\psi}(x)- u^*_{\psi}(x)|&=\frac{\big|\sum_{\wt x\in \ov{\m N}^x_{\varepsilon^x_j}} \big(\sum_{k\in \mb N_0^{D_X}, |k|<\alpha_X}u^*_{\psi}{}^{(k)}(\wt x) (x-\wt x)^k -u^*_{\psi}(x)\big)\rho(\frac{\|x-\wt x\|}{\varepsilon^x_j})\big|}{\sum_{\wt x\in \ov{\m N}^x_{\varepsilon^x_j}}  \rho(\frac{\|x-\wt x\|}{\varepsilon^x_j})}\\
          &\leq \underset{\wt x\in \ov{\m N}^x_{\varepsilon^x_j}, \,x\in \mb{B}_{\m M_X}(\wt x,2\varepsilon_j^x)}{\sup}\big|\sum_{k\in \mb N_0^{D_X}, |k|<\alpha_X}u^*_{\psi}{}^{(k)}(\wt x) (x-\wt x)^k -u^*_{\psi}(x)\big|\\
          &\lesssim  2^{-\frac{d_Yj}{2}} (\varepsilon^x_j)^{\alpha_X}.
    \end{aligned}
\end{equation*}
 We can get 
\begin{equation*}
\begin{aligned}
&\underset{S\in \m S_j}{\min} \, \mb{E}_{\mu^*_X} \Big[\sum_{\psi\in \Psi_j^{D_Y}} ( S(\psi,X)-u^*_{\psi}(X))^2\Big]\\
   &\leq  \mb{E}_{\mu^*_X} \Big[\sum_{\psi\in \Psi_j^{D_Y}}( S^*_j(\psi,X)- u^*_{\psi}(X))^2\Big]\\
   &= \mb{E}_{\mu^*_X} \Big[\sum_{\psi\in \Psi_j^*}( S^*_j(\psi,X)- u^*_{\psi}(X))^2\Big]\\
   &\lesssim \sum_{\psi\in \Psi_j^*}  2^{-{d_Yj}} ((\varepsilon^x_j)^{-\alpha_X}+\frac{1}{n})^2\\
   &\lesssim (\varepsilon^x_j)^{2\alpha_X}+\frac{1}{n^2}.
   \end{aligned}
\end{equation*}
The desired result then follows by substituting $\varepsilon_j^x= 2^{\frac{jd_Y}{2\alpha_X+d_X}}(\frac{n}{\log n})^{-\frac{1}{2\alpha_X+d_X}}$.

 \subsection{Proof of Lemma~\ref{lemma:manifoldlearningerrorr2}}\label{proof:lemma:manifoldlearningerrorr2}
 We begin by establishing a general lemma to bound the population-level reconstruction error. Consider arbitrary points $x_0\in \mb B_{\mb R^{D_X}}(\mathbf{0},L)$ and $y_0\in \mb B_{\mb R^{D_Y}}(\mathbf{0},L)$, and  consider the estimator
\begin{equation*}
    (\wh G, \wh V)=\underset{G\in  {\m G}\atop V\in \mb O(D_Y,d_Y)}{\arg\min} \frac{1}{n}\sum_{i=1}^n \|Y_i-G(V^T(Y_i-y_0),X_i)\|^2\mathbf{1}(X_i\in \mb B_{\mb R^{D_X}}(x_0,2\tau_2))\mathbf{1}(Y_i\in \mb B_{\mb R^{D_Y}}(y_0,2\tau_2)),
\end{equation*}
where $\big\{(X_i,Y_i)\big\}_{i=1}^n$ are i.i.d. samples from $\mu^*$ and $\m G$ represents an arbitrary class of functions  $G:\mb R^{d_Y}\times \mb R^{D_X}\to \mb R^{D_Y}$.

\begin{lemma}\label{lemma:reconstruct}
Suppose  $\big\{(X_i,Y_i)\big\}_{i=1}^n$ are $n$ i.i.d data  from $\mu^*=\mu^*_X\mu^*_{Y|X}$, and the following assumptions are satisfied: (1) for any $x\in \m M_X={\rm supp}(\mu^*_X)$, the conditional distribution of $Y$ given $X=x$, denoted as $\mu^*_{Y|x}$, is supported on a submanifold $\m M_{Y|x}$, and has a density function $u^*(\,\cdot\,|\,x)$ with respect to the volume measure of $\m M_{Y|x}$. There exist constants $\beta_Y\geq 2, \beta_X,\alpha_Y,\alpha_X,L>0$ and a function $\ov u^*\in    \ov{\m H}^{\alpha_Y,\alpha_X}_{L}(\mb R^{D_Y},\mb R^{D_X})$ such that $\{\m M_{Y|x}\}_{x\in \m M_X}\in \ms M^{\beta_Y\,\beta_X}_{\tau,\tau_1,L}(d_Y,D_Y,\m M_X)$ and $\ov u^*|_{\m M}=u$;  (2) there exists a function $g:\mb R^{+}\to \mb R^{+}$ such that for any $x_0\in \m M_X$, $y_0\in \m M_{Y|x_0}$ and for all $0<r\leq 1$, it holds that $\mu^*_X(\mb B_{\m M_X}(x_0,r))\geq  g(r)$ and $\mu^*_{Y|x_0}(\mb B_{\m M_{Y|x}}(y_0,r))\geq g(r)$; (3) there exist constants $L>0$ and $\beta>1$ such that that for any $G(z,x)\in \m G$, it holds for any $x\in \m M_X$ that $G(\cdot,x)\in \m H^{\beta}_{L,D_Y}(\mb R^{d_Y})$. Then 
\begin{enumerate}
    \item If there exists $G \in \m G$ and $V \in \mb O(D_Y,d_Y)$ such that for any $(x,y)\in \m M=\{(x,y):\, x\in \m M_X,y\in \m M_{Y|x}\}$ with $x\in \mb B_{\mb R^{D_X}}(x_0,2\tau_2)$ and $y\in \mb B_{\mb R^{D_Y}}(y_0,2\tau_2)$,  it holds that $\|y-G(V^T(y-y_0),x)\|\leq \varepsilon^*$. Consider any $\gamma_1\in (0,1]$ and denote  $\mathbf N(\m G,d_{\infty}^{\gamma_1},\varepsilon)$
  as the $\varepsilon$-covering number of $\m G$ with respect to the $d_{\infty}^{\gamma_1}$ distance, where $d_{\infty}^{\gamma_1}(G_1,G_2)=\underset{z\in \mb R^{d_Y},x\in \mb R^{D_X}}{\sup}\|G_1(z,x)-G_2(z,x)\|^{\gamma_1}$. There exists a constant $C$ so that, with probability at least $1-\frac{1}{n^3}$,
  \begin{equation*}
 \begin{aligned}
 &\mb{E}_{\mu^*}[\|Y-\wh G(\wh V^T(Y-y_0),X)\|^{\gamma_1}\cdot\mathbf{1}(X\in \mb B_{\mb R^{D_X}}(x_0,2\tau_2))\mathbf{1}(Y\in \mb B_{\mb R^{D_Y}}(y_0,2\tau_2))]\\
 &\leq C\,  \Big(\frac{1}{\sqrt{n}}\int_0^{\infty} \sqrt{ \log \mathbf{N}(\m G,d_{\infty}^{\gamma_1},\varepsilon/2)}\dd \varepsilon+\sqrt{\frac{\log n}{n}}+(\varepsilon^*)^{\gamma_1}\Big).
      \end{aligned}
\end{equation*} 
    \item If  there exists $(x^*,y^*)\in \mb B_{\m M}((x_0,y_0),\sqrt{2}\tau_2)$, and $\tau_2<\frac{\tau_1\wedge \tau}{2}$. Then let $ P^*$ be the projection matrix of $T_{\m M_{Y|x^*}}y^*$,  there exist positive constants $c,c_1$ so that if  $\mb{E}_{\mu^*}[\|Y-\wh G(\wh V^T(Y-y_0),X)\|\cdot\mathbf{1}(X\in \mb B_{\mb R^{D_X}}(x_0,2\tau_2))\mathbf{1}(Y\in \mb B_{\mb R^{D_Y}}(y_0,2\tau_2))]\leq c$, then $\wh V^T P^*\wh V^T\geq c_1 I_{d_Y}$.
\end{enumerate} 
\end{lemma}
\noindent The proof of Lemma~\ref{lemma:reconstruct} can be found in Appendix~\ref{proof:lemma:reconstruct}. Given Lemma~\ref{lemma:reconstruct}, it suffices to demonstrate the first statement of Lemma~\ref{lemma:manifoldlearningerrorr2}. The second statement of Lemma~\ref{lemma:manifoldlearningerrorr2} naturally follows from the second statement of Lemma~\ref{lemma:reconstruct}. Consider the family
\begin{equation}\label{defmG}
\begin{aligned}
       {\m G}&=\{G(z)=\sum_{j_1=0}^{J_1} \sum_{\psi_1\in  \Psi_{j_1}^{d_Y} } g_{\psi_1 }\psi_1(z) \,:  g_{\psi_1}\in [-L_1\, \delta_{j_1},L_1\, \delta_{j_1}]^{D_Y}, \text{ for }\psi_1\in \Psi_{j_1}^{d_Y} \},
\end{aligned}
\end{equation}
where $J_1= \lceil\log_2 (n^{-\frac{1}{d_Y}})\rceil$ and $\delta_{j_1}=2^{-\frac{d_Yj_1}{2}-(j_1\beta_Y)}$. It is straightforward to verify that for any $\beta<\beta_Y$, there exists a constant $L$ so that $\m G\subset \m H^{\beta}_{L,D_Y}(\mb R^{d_Y})$. Moreover, we can derive the following lemma that control the covering number of $\m G$, the proof of which is given in Appendix~\ref{proof:lemma:coveringGregime2}.
\begin{lemma}\label{lemma:coveringGregime2}
With the choice of $\m G$ in~\eqref{defmG}, there exists a constant $C_1$ so that  for any $\gamma_1\in (0,1]$,  the $\varepsilon$-covering number $\mathbf N(\m G,d_{\infty}^{\gamma_1},\varepsilon)$  of $\m G$ with respect to the $d_{\infty}^{\gamma_1}$ distance, satisfies that 
 \begin{equation*}
 \begin{aligned}
    &   \log \mathbf N(\m G,d_{\infty}^{\gamma_1},\varepsilon) \\
       &\leq \left\{\begin{array}{cc}
   C_1\sum_{j_1=0}^{J_1}  2^{d_Yj_1}\log\left(\frac{C_1 J_12^{-\frac{d_Yj_1}{4\gamma_1}-\frac{j_1\beta_Y}{2}}}{\varepsilon^{\frac{1}{\gamma_1}}}\vee 1\right)    &  \frac{d_Y}{\beta_1} \leq 2\gamma_1, \\
      C_1\sum_{j_1=0}^{J_1}  2^{d_Yj_1}\log \left(\frac{C_1(J_1\wedge c(\beta_Y,d_Y,d_X,\gamma_1))2^{-j_1\beta_Y}}{\varepsilon^{\frac{1}{\gamma_1}}  s_{j_1}}\vee 1\right)    &  \frac{d_Y}{\beta_Y} > 2\gamma_1,
  \end{array}
  \right.
 \end{aligned}
 \end{equation*}
 where $c(\beta_Y,d_Y,d_X,\gamma_1)= \frac{2^{\frac{(d_Y-2\beta_Y\gamma_1)}{4\gamma_1}}}{2^{\frac{(d_Y-2\beta_Y\gamma_1)}{4\gamma_1}}-1}$  and $ s_{j_1}=\sqrt{\frac{2^{\frac{d_Yj_1}{2\gamma_1}-j_1\beta_Y}}{2^{\frac{d_YJ_1}{2\gamma_1}-J_1\beta_Y}}}$.
\end{lemma}
\noindent  Notice that $|I_1|=\lfloor n/2\rfloor \asymp n$,  we will  then bound the integral $\frac{1}{\sqrt{n}}\int_0^{\infty} \sqrt{ \log \mathbf N(\m G, d_{\infty}^{\gamma_1},\varepsilon)}\dd \varepsilon$. When $\frac{d_Y}{\beta_Y}\leq 2\gamma_1$, we have 
\begin{equation*}
    \begin{aligned}
         &\frac{1}{\sqrt{n}}\int_0^{\infty} \sqrt{ \log \mathbf N(\m G, d_{\infty}^{\gamma_1},\varepsilon)}\dd \varepsilon\\
         &\lesssim \frac{1}{\sqrt{n}}\int_0^{\infty} \sqrt{\sum_{j_1=0}^{J_1} \log \left(\frac{C_1 J_12^{-\frac{d_Yj_1}{4\gamma_1}-\frac{j_1\beta_Y}{2}}}{\varepsilon^{\frac{1}{\gamma_1}}}\vee 1\right) 2^{d_Yj_1}}\,\dd \varepsilon\\
      &\lesssim \frac{1}{\sqrt{n}} \sum_{j_1=0}^{J_1}  \int_0^{\infty}\sqrt{\log \left(\frac{C_1 J_12^{-\frac{d_Yj_1}{4\gamma_1}-\frac{j_1\beta_Y}{2}}}{\varepsilon^{\frac{1}{\gamma_1}}}\vee 1\right) 2^{d_Yj_1 }}\,\dd \varepsilon\\
 &\lesssim \frac{1}{\sqrt{n}} \sum_{j_1=0}^{J_1}  \int_0^{\big(C_1 J_12^{-\frac{d_Yj_1}{4\gamma_1}-\frac{j_1\beta_Y)}{2}}\big)^{\gamma_1}}\sqrt{\log \left(\frac{C_1 J_12^{-\frac{d_Yj_1}{4\gamma_1}-\frac{j_1\beta_Y}{2}}}{\varepsilon^{1/\gamma_1}}\right) 2^{d_Yj_1}}\,\dd \varepsilon\\
&\lesssim \frac{1}{\sqrt{n}} \sum_{j_1=0}^{J_1} J_1^{\gamma_1}2^{\frac{d_Yj_1}{4}-\gamma_1\frac{j_1\beta_Y}{2}} \lesssim \frac{(\log n)^{1+\gamma_1}}{\sqrt{n}}.
    \end{aligned}
\end{equation*}
When $\frac{d_Y}{\beta_Y}>2\gamma_1$, we have 
\begin{equation*}
    \begin{aligned}
         &\frac{1}{\sqrt{n}}\int_0^{\infty} \sqrt{ \log \mathbf N(\m G, d_n^{\gamma_1},\varepsilon)}\dd \varepsilon\\
         &\lesssim \frac{1}{\sqrt{n}}\int_0^{\infty} \sqrt{\sum_{j_1=0}^{J_1}  \log \left(\frac{C_1(J_1\wedge c(\beta_Y,d_Y,d_X,\gamma_1))2^{-j_1\beta_Y}}{\varepsilon^{1/\gamma_1} s_{j_1}}\vee 1\right) 2^{d_Yj_1}}\,\dd \varepsilon\\
      &\lesssim \frac{1}{\sqrt{n}} \sum_{j_1=0}^{J_1}\int_0^{\infty}\sqrt{\log \left(\frac{C_1(J_1\wedge c(\beta_Y,d_Y,d_X,\gamma_1))2^{-j_1\beta_Y}}{\varepsilon^{1/\gamma_1} s_{j_1}}\vee 1\right) 2^{d_Yj_1}}\,\dd \varepsilon\\
&\lesssim \frac{1}{\sqrt{n}} \sum_{j_1=0}^{J_1} \frac{1}{s_{j_1}^{\gamma_1}}(J_1\wedge c(\beta_Y,d_Y,d_X,\gamma_1))^{\gamma_1}2^{-\gamma_1j_1\beta_Y} 2^{\frac{d_Yj_1}{2}}\\
&\lesssim  (J_1\wedge c(\beta_Y,d_Y,d_X,\gamma_1))^{1+\gamma_1}   \frac{1}{\sqrt{n}} 2^{-\gamma_1J_1\beta_Y} 2^{\frac{d_YJ_1}{2}}\\
&\lesssim (\log n \wedge \frac{1}{ d_Y-2\gamma_1\beta_Y} )^{1+\gamma_1}\cdot n^{-\frac{\gamma_1}{\frac{d_Y}{\beta_Y}}}.
    \end{aligned}
\end{equation*}
 Then it remains to bound the term  $\varepsilon^*$.  Fix an arbitrary $k\in [K]$. If $\mb B_{\mb R^{D_X+D_Y}}((x_k,y_k),\sqrt{2}\tau_2)\cap \m M=\emptyset$, then $k\notin \wh{\m K}$. Otherwise, there exists $(x_k^*,y_k^*)\in \mb B_{\m M}((x_k,y_k),\sqrt{2}\tau_2)$. Let $V_{\sk}^*$ be an arbitrary  orthonormal basis of $T_{\m M_{Y}}y_k^*$, and denote $Q_\sk^*(y)=(V_\sk^*)^T(y-y_k)$ and $G_\sk^*(z)= \Phi_{y_k^*}(V_\sk^*(z+(V_\sk^*)^T(y_k-y_k^*)))$. Then  $G_\sk^*\in \m H^{\beta_Y}_{L,D_Y}\big(\mb B_{\mb R^{d_Y}}(V_\sk^*{}^T(y_k^*-y_k),\tau_1)
 \big)$ and for any $y\in \m M_Y$ with  $\|y-y_k^*\|<\tau_1$,  we have $y= G_\sk^*(Q_\sk^*(y))$.   Moreover, by leveraging the decay of wavelet coefficients for $ \m H^{\beta_Y}$-smooth functions as stated in Lemma~\ref{le:wavelet}, when $J_1= \lceil\log_2 (n^{-\frac{1}{d_Y}})\rceil$ and $\tau_2<\frac{\tau_1\wedge \tau}{4}$, it holds that for any $z\in \mb B_{\mb R^{d_Y}}(\mathbf{0},2\tau_2)\subset \mb B_{\mb R^{d_Y}}(V_\sk^*{}^T(y_k^*-y_k),\tau_1)$ that,
 \begin{equation*}
    \Big\|G_\sk^*(z)- \sum_{j=1}^{J_1} \sum_{\psi\in \Psi_j^{d_Y}} \int_{\mb R^{d_Y}} G_{\sk}^*(z)\psi(z)\,\dd z\cdot \psi(z)\Big\|\leq C\,  n^{-\frac{1}{\frac{d_Y}{\beta_Y}}}.
 \end{equation*}
Moreover, we have $G^{\dagger}_\sk(z)=\sum_{j=1}^{J_1} \sum_{\psi\in \Psi_j^{d_Y}} \int_{\mb R^{d_Y}} G_{\sk}^*(z)\psi(z)\,\dd z\cdot \psi(z) \in \m G$, and for any $y\in \m M_Y$ with $\|y-y_k\|\leq 2\tau_2$, 
\begin{equation*}
\begin{aligned}
     \|y-G^{\dagger}_\sk((V_\sk^*)^T(y-y_0))\|\leq     \|y-G^*_\sk(Q^*_\sk(y))\|+C\,n^{-\frac{1}{\frac{d_Y}{\beta_Y}}}&=C\,n^{-\frac{1}{\frac{d_Y}{\beta_Y}}}.
\end{aligned}
\end{equation*}
Therefore,  let $\wh Q_\sk(\cdot)=\wh V_\sk^T(\cdot-y_k)$, using Lemma~\ref{lemma:coveringGregime2}, we can conclude that for any $\gamma_1\in (0,1]$, there exists a constant $C_{\gamma_1}$ so that it holds with probability at least $1-\frac{c}{n^3}$ that for any $k\in \wh {\m K}$, 

\begin{equation*}
    \begin{aligned}
       &\mb{E}_{\mu^*}[\|Y-\wh G_\sk(\wh Q_\sk(Y))\|^{\gamma_1}\cdot\mathbf{1}(X\in \mb B_{\mb R^{D_X}}(x_0,2\tau_2))\mathbf{1}(Y\in \mb B_{\mb R^{D_Y}}(y_0,2\tau_2))]\\
       &\leq  \left\{
          \begin{array}{cc}
           C_{\gamma_1}\, \frac{(\log n)^{1+\gamma_1}}{\sqrt{n}}   & \frac{d_Y}{\beta_Y} \leq 2\gamma_1, \\
C_{\gamma_1}\,(\log n \wedge \frac{1}{d_Y-2\gamma_1\beta_Y)} )^{1+\gamma_1}\cdot n^{-\frac{\gamma_1}{ \frac{d_Y}{\beta_Y}}}   & \frac{d_Y}{\beta_Y}>2\gamma_1.
          \end{array}
        \right.
    \end{aligned}
\end{equation*}
Then if $\frac{d_Y}{2\beta_Y}>1$, set $\gamma_1=1$,  it holds with probability at least $1-\frac{c}{n^3}$ that for any $k\in  \wh {\m K}$, 
\begin{equation*}
    \begin{aligned}
       &\mb{E}_{\mu^*}[\|Y-\wh G_\sk(\wh Q_\sk(Y))\|\cdot\mathbf{1}(X\in \mb B_{\mb R^{D_X}}(x_0,2\tau_2))\mathbf{1}(Y\in \mb B_{\mb R^{D_Y}}(y_0,2\tau_2))]\lesssim  
 n^{-\frac{1}{ \frac{d_Y}{\beta_Y}}}.
     \end{aligned}
\end{equation*}
Therefore, with probability at least $1-\frac{c}{n^3}$, it holds for any $k\in  \wh {\m K}$ and any $\gamma_1\in (0,1]$ that
\begin{equation*}
    \begin{aligned}
       &\mb{E}_{\mu^*}[\|Y-\wh G_\sk(\wh Q_\sk(Y))\|^{\gamma_1}\cdot\mathbf{1}(X\in \mb B_{\mb R^{D_X}}(x_0,2\tau_2))\mathbf{1}(Y\in \mb B_{\mb R^{D_Y}}(y_0,2\tau_2))]\\
       &\leq \Big(\mb{E}_{\mu^*}[\|Y-\wh G_\sk(\wh Q_\sk(Y))\|\cdot\mathbf{1}(X\in \mb B_{\mb R^{D_X}}(x_0,2\tau_2))\mathbf{1}(Y\in \mb B_{\mb R^{D_Y}}(y_0,2\tau_2))]\Big)^{\gamma_1}\lesssim n^{-\frac{\gamma_1}{ \frac{d_Y}{\beta_Y}}}.
     \end{aligned}
\end{equation*}
If $\frac{d_Y}{2\beta_Y}\leq 1$, let $\delta_n=\frac{1-\frac{d_Y}{4\beta_Y}}{\lceil\log n\rceil}$ and  consider the set $\Gamma=\{\frac{d_Y}{4\beta_Y}, \frac{d_Y}{4\beta_Y}+\delta_n,\cdots, \frac{d_Y}{4\beta_Y}+ \delta_n\cdot \lceil\log n\rceil\}$. Then by a union argument,  it holds that with probability at least $1-\frac{c\,\log n}{n^3}$ that for any $k\in  \wh {\m K}$ and any $\gamma_1\in\Gamma$ that
\begin{equation*}
    \begin{aligned}
       &\mb{E}_{\mu^*}[\|Y-\wh G_\sk(\wh Q_\sk(Y))\|^{\gamma_1}\cdot\mathbf{1}(X\in \mb B_{\mb R^{D_X}}(x_0,2\tau_2))\mathbf{1}(Y\in \mb B_{\mb R^{D_Y}}(y_0,2\tau_2))]\\
       &\leq  \left\{
          \begin{array}{cc}
           C \, \frac{(\log n)^{1+\gamma_1}}{\sqrt{n}}   & \frac{d_Y}{\beta_Y} \leq 2\gamma_1, \\
C \,(\log n \wedge \frac{1}{d_Y-2\gamma_1\beta_Y)} )^{1+\gamma_1}\cdot n^{-\frac{\gamma_1}{ \frac{d_Y}{\beta_Y}}}   & \frac{d_Y}{\beta_Y}>2\gamma_1.
          \end{array}
        \right.
    \end{aligned}
\end{equation*}
Under the above event, for any $\gamma_2\in (0,\frac{d_Y}{4\beta_Y})$, by setting $\gamma_1=\frac{d_Y}{4\beta_Y}$, it holds that 

\begin{equation*}
    \begin{aligned}
       &\mb{E}_{\mu^*}[\|Y-\wh G_\sk(\wh Q_\sk(Y))\|^{\gamma_2}\cdot\mathbf{1}(X\in \mb B_{\mb R^{D_X}}(x_0,2\tau_2))\mathbf{1}(Y\in \mb B_{\mb R^{D_Y}}(y_0,2\tau_2))]\\
       &\leq \Big(\mb{E}_{\mu^*}[\|Y-\wh G_\sk(\wh Q_\sk(Y))\|^{\gamma_1}\cdot\mathbf{1}(X\in \mb B_{\mb R^{D_X}}(x_0,2\tau_2))\mathbf{1}(Y\in \mb B_{\mb R^{D_Y}}(y_0,2\tau_2))]\Big)^{\gamma_2/\gamma_1}\\
       &\lesssim n^{-\frac{\gamma_2}{4\gamma_1}}=n^{-\frac{\gamma_2}{d_Y/\beta_Y}}.
     \end{aligned}
\end{equation*}
Moreover, for any $\gamma_2\in [\frac{d_Y}{4\beta_Y},1]$, there exists $\gamma_1\in \Gamma$ so that $\gamma_1\leq\gamma_2\leq\gamma_1+\delta_n$, so
\begin{equation*}
    \begin{aligned}
       &\mb{E}_{\mu^*}[\|Y-\wh G_\sk(\wh Q_\sk(Y))\|^{\gamma_2}\cdot\mathbf{1}(X\in \mb B_{\mb R^{D_X}}(x_0,2\tau_2))\mathbf{1}(Y\in \mb B_{\mb R^{D_Y}}(y_0,2\tau_2))]\\
       &\leq 2L\, \mb{E}_{\mu^*}[\|Y-\wh G_\sk(\wh Q_\sk(Y))\|^{\gamma_1}\cdot\mathbf{1}(X\in \mb B_{\mb R^{D_X}}(x_0,2\tau_2))\mathbf{1}(Y\in \mb B_{\mb R^{D_Y}}(y_0,2\tau_2))] \\
       &\leq  \left\{
          \begin{array}{cc}
       C_1   \frac{(\log n)^{1+\gamma_1}}{\sqrt{n}}   & \frac{d_Y}{\beta_Y} \leq 2\gamma_1, \\
 C_1(\log n \wedge \frac{1}{d_Y-2\gamma_1\beta_Y)} )^{1+\gamma_1}\cdot n^{-\frac{\gamma_1}{ \frac{d_Y}{\beta_Y}}}   & \frac{d_Y}{\beta_Y}>2\gamma_1\\
          \end{array}
        \right.\\
        &\leq \left\{
          \begin{array}{cc}
          C_2\frac{(\log n)^{1+\gamma_2}}{\sqrt{n}}   & \frac{d_Y}{\beta_Y} \leq 2\gamma_2, \\
C_2(\log n \wedge \frac{1}{d_Y-2\gamma_1\beta_Y)} )^{1+\gamma_2}\cdot n^{-\frac{\gamma_2}{ \frac{d_Y}{\beta_Y}}}   & \frac{d_Y}{\beta_Y}>2\gamma_2.
          \end{array}
        \right.
     \end{aligned}
\end{equation*}
This completes the proof  of Lemma~\ref{lemma:manifoldlearningerrorr2}.

\subsection{Proof of Lemma~\ref{le2:app}}\label{proofle2:app}
We will show the desired result using Lemma~\ref{lestimation:2}. For the family  ${\m S}_j^\dagger$  that consists of 
\begin{equation*} 
    \begin{aligned}
    &S(\psi,x)=\m T_{C_1\,2^{-\frac{d_Yj}{2}}}\Big(\\
   & \frac{\sum_{i_1=1}^{W_j}\sum_{i_2=1}^{W_j'} \sum_{k=1}^{K^*} \int_{\mb B_{\mb R^{d_Y}}(\mathbf{0},\tau_1)} 2^{\frac{j(d_Y-D_Y)}{2}}\psi(G_{k,i_2}(z,x))\nu_{k,i_2}(z,x)\,\dd z \rho(\frac{\|x-b_{i_2}\|}{\varepsilon_j^x})\rho(\frac{|\m I_j(\psi)-e_{i_1i_2}|}{\varepsilon_j^y})}{\sum_{i_1=1}^{W_j}\sum_{i_2=1}^{W_j'}\rho(\frac{\|x-b_{i_2}\|}{\varepsilon_j^x})\rho(\frac{|\m I_j(\psi)-e_{i_1i_2}|}{\varepsilon_j^y})+\frac{1}{n^2}}\\
   &+\frac{\sum_{i_1=1}^{W_j}\sum_{i_2=1}^{W_j'}\sum_{l\in\mb N_{0}^{D_X}\atop |l|\leq\lfloor\wt\beta_X\rfloor^2+\lfloor \alpha_X\rfloor} a_{i_1i_2l} (x-b_{i_2})^l \rho(\frac{\|x-b_{i_2}\|}{\varepsilon_j^x})\rho(\frac{|\m I_j(\psi)-e_{i_1i_2}|}{\varepsilon_j^y})}{\sum_{i_1=1}^{W_j}\sum_{i_2=1}^{W_j'}\rho(\frac{\|x-b_{i_2}\|}{\varepsilon_j^x})\rho(\frac{|\m I_j(\psi)-e_{i_1i_2}|}{\varepsilon_j^y})+\frac{1}{n^2}}\Big),\\
  &\text{where } \wt\beta_X=\alpha_X+\frac{\alpha_X}{\alpha_Y},\\
 &G_{k,i_2}(z,x)=\sum_{s=0}^{j} \sum_{\psi\in \wt\Psi_s^{d_Y}}\sum_{l\in \mb N_0^{D_X}\atop |l|<\beta_X}  g_{k,i_2,s,\psi,l} (x-b_{i_2})^{l} \cdot \psi(z) \\
      &\text{and }\nu_{k,i_2}(z,x)=\sum_{s=0}^{j} \sum_{\psi\in \wt\Psi_s^{d_Y}}\sum_{l\in \mb N_0^{D_X}\atop |l|<\alpha_X}  v_{k,i_2,s,\psi,l} (x-b_{i_2})^{l} \cdot \psi(z), \\
   \end{aligned}
\end{equation*}
where the parameters satisfy that $ g_{k,i_2,s,\psi,l}\in[-C_1,C_1]^{D_Y}$, $v_{k,i_2,s,\psi,l}\in [-C_1,C_1]$, $a_{i_1i_2l}\in [-C_1\,n, C_1\,n]$, $e_{i_1i_2}\in [0,2]^{D_Y+1}$, and $\{b_1,b_2,\cdots,b_{W_j'}\}\subset \mb B_{\mb R^{D_X}}(\mathbf{0},L_1)$    are $\varepsilon_j^x-$ separated.  It holds for any $S\in S_j^\dagger$ that 
\begin{equation*}
    \underset{\psi\in \Psi_j^{D_Y}}{\sup}\underset{x\in \m M_X}{\sup} |S(\psi,x)|\leq C_1 2^{-\frac{d_Yj}{2}}.
\end{equation*}
Moreover, for any $\psi\in \Psi_j^{D_Y}$  and $x\in \m M_X$,  $S(\psi,x)$ will be non-zero only if there exist $i_1\in [W_j]$ and $i_2\in [W_j']$ so that $\|x-b_{i_2}\|<2\varepsilon_j^x$ and $|\m I_j(\psi)-e_{i_1i_2}|<2\varepsilon_j^y$.  Given that 
the set $\{b_1,b_2,\cdots,b_{W_j'}\}$ are $\varepsilon_j^x$-separated,  for any $x\in \m M_X$, there are $\m O(1)$ number of $i_2\in [W_j']$ so that $\|x-b_{i_2}\|<2\varepsilon_j^x$.  Moreover, for any  $i_2\in [W_j']$ and $i_1\in [W_j]$, there are at most constant number of $\psi\in \Psi_j^{D_Y}$ so that $|\m I_j(\psi)-e_{i_1i_2}|<2\varepsilon_j^y$. Therefore,  for any $x\in \m M_X$,  there are   $\m O(W_j)=\m O(2^{d_Yj})$  number of $\psi\in \Psi_j^{D_Y}$ so that $S(\psi,x)\neq 0$, and thus
\begin{equation*}
     \sum_{\psi\in \Psi_j^{D_Y}} \mathbf{1}(S(\psi,x)\neq 0)\leq C \,2^{d_Y j}.
\end{equation*}
Furthermore, consider $S,S'\in \m S_j^{\dagger}$ with
\begin{equation*}
    \begin{aligned}
 &S(\psi,x)=\m T_{C\,2^{-\frac{d_Yj}{2}}}\Big(\\
   & \frac{\sum_{i_1=1}^{W_j}\sum_{i_2=1}^{W_j'} \sum_{k=1}^{K^*} \int_{\mb B_{\mb R^{d_Y}}(\mathbf{0},\tau_1)} 2^{\frac{j(d_Y-D_Y)}{2}}\psi(G_{k,i_2}(z,x))v_{k,i_2}(z,x)\,\dd z \rho(\frac{\|x-b_{i_2}\|}{\varepsilon_j^x})\rho(\frac{|\m I_j(\psi)-e_{i_1i_2}|}{\varepsilon_j^y})}{\sum_{i_1=1}^{W_j}\sum_{i_2=1}^{W_j'}\rho(\frac{\|x-b_{i_2}\|}{\varepsilon_j^x})\rho(\frac{|\m I_j(\psi)-e_{i_1i_2}|}{\varepsilon_j^y})+\frac{1}{n^2}} \\
   &+\frac{\sum_{i_1=1}^{W_j}\sum_{i_2=1}^{W_j'}\sum_{l\in\mb N_{0}^{D_X}\atop |l|\leq\lfloor\wt\beta_X\rfloor^2+\lfloor \alpha_X\rfloor} a_{i_1i_2l} (x-b_{i_2})^l \rho(\frac{\|x-b_{i_2}\|}{\varepsilon_j^x})\rho(\frac{|\m I_j(\psi)-e_{i_1i_2}|}{\varepsilon_j^y})}{\sum_{i_1=1}^{W_j}\sum_{i_2=1}^{W_j'}\rho(\frac{\|x-b_{i_2}\|}{\varepsilon_j^x})\rho(\frac{|\m I_j(\psi)-e_{i_1i_2}|}{\varepsilon_j^y})+\frac{1}{n^2}}\Big),\\
  &\text{where } G_{k,i_2}(z,x)=\sum_{s=0}^{j} \sum_{\psi\in \wt\Psi_s^{d_Y}}\sum_{l\in \mb N_0^{D_X}\atop |l|<\beta_X}  g_{k,i_2,s,\psi,l} (x-b_{i_2})^{l} \cdot \psi(z) \\
      &\text{and }v_{k,i_2}(z,x)=\sum_{s=0}^{j} \sum_{\psi\in \wt\Psi_s^{d_Y}}\sum_{l\in \mb N_0^{D_X}\atop |l|<\alpha_X}  v_{k,i_2,s,\psi,l} (x-b_{i_2})^{l} \cdot \psi(z) , \\
    \end{aligned}
\end{equation*}
and
\begin{equation*}
    \begin{aligned}
     &S'(\psi,x)=\m T_{C\,2^{-\frac{d_Yj}{2}}}\Big(\\
   & \frac{\sum_{i_1=1}^{W_j}\sum_{i_2=1}^{W_j'} \sum_{k=1}^{K^*}  \int_{\mb B_{\mb R^{d_Y}}(\mathbf{0},\tau_1)} 2^{\frac{j(d_Y-D_Y)}{2}}\psi(G'_{k,i_2}(z,x))v'_{k,i_2}(z,x)\,\dd z \rho(\frac{\|x-b'_{i_2}\|}{\varepsilon_j^x})\rho(\frac{|\m I_j(\psi)-e'_{i_1i_2}|}{\varepsilon_j^y})}{\sum_{i_1=1}^{W_j}\sum_{i_2=1}^{W_j'}\rho(\frac{\|x-b'_{i_2}\|}{\varepsilon_j^x})\rho(\frac{|\m I_j(\psi)-e'_{i_1i_2}|}{\varepsilon_j^y})+\frac{1}{n^2}} \\
   &+\frac{\sum_{i_1=1}^{W_j}\sum_{i_2=1}^{W_j'}\sum_{l\in\mb N_{0}^{D_X}\atop |l|\leq\lfloor\wt\beta_X\rfloor^2+\lfloor \alpha_X\rfloor} a'_{i_1i_2l} (x-b'_{i_2})^l \rho(\frac{\|x-b'_{i_2}\|}{\varepsilon_j^x})\rho(\frac{|\m I_j(\psi)-e'_{i_1i_2}|}{\varepsilon_j^y})}{\sum_{i_1=1}^{W_j}\sum_{i_2=1}^{W_j'}\rho(\frac{\|x-b'_{i_2}\|}{\varepsilon_j^x})\rho(\frac{|\m I_j(\psi)-e'_{i_1i_2}|}{\varepsilon_j^y})+\frac{1}{n^2}}\Big),\\
  &\text{where }G'_{k,i_2}(z,x)=\sum_{s=0}^{j} \sum_{\psi\in \wt\Psi_s^{d_Y}}\sum_{l\in \mb N_0^{D_X}\atop |l|<\beta_X}  g'_{k,i_2,s,\psi,l} (x-b_{i_2}')^{l} \cdot \psi(z) \\
      &\text{and }v'_{k,i_2}(z,x)=\sum_{s=0}^{j} \sum_{\psi\in \wt\Psi_s^{d_Y}}\sum_{l\in \mb N_0^{D_X}\atop |l|<\alpha_X}  v'_{k,i_2,s,\psi,l} (x-b_{i_2}')^{l} \cdot \psi(z). \\
    \end{aligned}
\end{equation*}
It holds for any $\psi\in \Psi_j^{D_Y}$ and $x\in \m M_X$ that,  
 \begin{equation*}
    \begin{aligned}
        &|S(\psi,x)-S'(\psi,x)|\\
        &\lesssim    j 2^{\frac{jd_Y}{2}+j}\underset{i_2\in [W_j']}{\sup}\underset{k\in [K^*]}{\sup}\underset{s\in\{0,\cdots,j\}\atop \wt\psi\in \wt\Psi_s^{d_Y}}{\sup}\sum_{l\in \mb N_0^{D_X}\atop |l|<\beta_X} \|g'_{k,i_2,s,\wt\psi,l}-g_{k,i_2,s,\wt\psi,l}\|\\
        &+ j 2^{\frac{jd_Y}{2}+j}\underset{i_2\in [W_j']}{\sup}\underset{k\in [K^*]}{\sup}\underset{s\in\{0,\cdots,j\}\atop \wt\psi\in \wt\Psi_s^{d_Y}}{\sup}\sum_{l\in \mb N_0^{D_X}\atop |l|<\beta_X} \|v'_{k,i_2,s,\wt\psi,l}-v_{k,i_2,s,\wt\psi,l}\|\\
         &+\max_{i_1\in [W_j],i_2\in [W_j']}\sum_{l\in\mb N_{0}^{D_X}\atop |l|<\alpha_X} |a_{i_1i_2l}-a'_{i_1i_2l}|+\frac{W_jn^{4}}{\varepsilon^x_j}\sum_{i_2=1}^{W_j'}\|b_{i2}-b_{i2'}\|+\frac{n^{4}}{\varepsilon^y_j}\sum_{i_1=1}^{W_j}\sum_{i_2=1}^{W_j'}\|e_{i_1i2}-e_{i_1i2}'\|.\\
    \end{aligned}
\end{equation*}
Therefore, we have 
\begin{equation*}
    \begin{aligned}
      d_S(S,S')&=\underset{x\in \m M_X}{\sup} \sqrt{\sum_{\psi\in \Psi_j^{D_Y}}(S(\psi,x)-S'(\psi,x))^2}\\
      & \leq \underset{x\in \m M_X}{\sup} \sqrt{\sum_{\psi\in \Psi_j^{D_Y}}(S(\psi,x)-S'(\psi,x))^2 (\mathbf{1}(S(\psi,x)\neq 0)+\mathbf{1}(S'(\psi,x)\neq 0))}\\
     &\lesssim 2^{\frac{jd_Y}{2}}\Bigg(  j 2^{\frac{jd_Y}{2}+j}\underset{i_2\in [W_j']}{\sup}\underset{k\in [K^*]}{\sup}\underset{s\in\{0,\cdots,j\}\atop \wt\psi\in \wt\Psi_s^{d_Y}}{\sup}\sum_{l\in \mb N_0^{D_X}\atop |l|<\beta_X} \|g'_{k,i_2,s,\wt\psi,l}-g_{k,i_2,s,\wt\psi,l}\|\\
        &+  j 2^{\frac{jd_Y}{2}+j}\underset{i_2\in [W_j']}{\sup}\underset{k\in [K^*]}{\sup}\underset{s\in\{0,\cdots,j\}\atop \wt\psi\in \wt\Psi_s^{d_Y}}{\sup}\sum_{l\in \mb N_0^{D_X}\atop |l|<\beta_X} \|v'_{k,i_2,s,\wt\psi,l}-v_{k,i_2,s,\wt\psi,l}\|\\
         &+\max_{i_1\in [W_j],i_2\in [W_j']}\sum_{l\in\mb N_{0}^{D_X}\atop |l|<\alpha_X} |a_{i_1i_2l}-a'_{i_1i_2l}|+\frac{W_jn^{4}}{\varepsilon^x_j}\sum_{i_2=1}^{W_j'}\|b_{i2}-b_{i2'}\|+\frac{n^{4}}{\varepsilon^y_j}\sum_{i_1=1}^{W_j}\sum_{i_2=1}^{W_j'}\|e_{i_1i2}-e_{i_1i2}'\|\Bigg).
    \end{aligned}
\end{equation*}
Then, using the fact that the $\varepsilon$-covering number of a $d$-dimensional ball with radius $R$ is being bounded by $(\frac{3R}{\varepsilon})^d$, we have for any $0<\varepsilon\leq \sup_{S,S'\in S_j^\dagger} d_S(S,S')$,
 \begin{equation*}
    \log \mathbf{N}(S_j^\dagger,d_S,\varepsilon)\lesssim 2^{jd_Y} (\varepsilon_j^x)^{-d_X}\log \frac{n}{\varepsilon}.
 \end{equation*}
Then by Lemma~\ref{lestimation:2}, it holds with probability larger than  $1-\frac{1}{n^2}$ that for any $j \in [J]$,
\begin{equation*}
\begin{aligned}
  &\mb{E}_{\mu^*_X} \Big[\sum_{\psi\in   \Psi_j^{D_Y}} \big(\mb{E}_{\mu^*_{Y|X}}[2^{\frac{j(d_Y-D_Y)}{2}}\psi(Y)]-\wh S_j^\dagger(\psi,X)\big)^2\Big]
   \lesssim   \frac{\log n}{n} 2^{jd_Y} (\varepsilon_j^x)^{-d_X} \\
 &\qquad\qquad+\underset{S\in \m S_j^\dagger}{\min} \, \mb{E}_{\mu^*_X}\Big[ \sum_{\psi\in \Psi_j^{D_Y}} ( \mb{E}_{\mu^*_{Y|X}}[2^{\frac{j(d_Y-D_Y)}{2}}\psi(Y)]-S(\psi,X))^2\Big].\\
\end{aligned}
\end{equation*}

 
\noindent Next, we bound the approximation error ${\min}_{S\in \m S_j^\dagger} \, \mb{E}_{\mu^*_X}\Big[ \sum_{\psi\in \Psi_j^{D_Y}} ( \mb{E}_{\mu^*_{Y|X}}[2^{\frac{j(d_Y-D_Y)}{2}}\psi(Y)]-S(\psi,X))^2\Big]$. Consider a  $\tau_2$-covering set $\{(x_k^*,y_k^*)\}_{k=1}^{K^*}\subset \m M$ of $\m M$,  by Lemma~\ref{legenerative}, we can write 
 \begin{equation*}
     \begin{aligned}
         \mb{E}_{\mu^*_{Y|x}}[2^{\frac{j(d_Y-D_Y)}{2}}\psi(y)]=\sum_{k=1}^{K^*} \int_{\mb B_{\mb R^{d_Y}}(\mathbf{0},\tau_1)} 2^{\frac{j(d_Y-D_Y)}{2}}\psi(G_{[k]}^*(z,x)) v_\sk^*(z,x)\,\dd z, \, x\in \m M_X,   \psi\in \Psi_j^{D_Y}, j\in \{0\}\cup [J],
     \end{aligned}
 \end{equation*}
 where $G^*_{[k]}\in  {\m H}^{\beta_Y,
 \beta_X}_{L_1,D_Y}(\mb R^{d_Y},\mb R^{D_X})$, $v^*_{[k]}\in  {\m H}^{\alpha_Y,\alpha_X}_{L_1}(\mb R^{d_Y},\mb R^{D_X})$. Moreover, for any $z\in \mb B_{\mb R^{d_Y}}(\mathbf{0},\tau_1)$ and $x\in \m M_X$,   $v^*_{[k]}(z,x)$ is zero if $\|x-x_k^*\|\geq  2\tau_2$  or $\|G^*_\sk(z,x)-y_k^*\|\geq  2\tau_2$.   Fix a $j\in \{0\}\cup [J]$ and let $\mb N_{\varepsilon_j^x}^x$ be the largest $\varepsilon_j^x$-packing set of $\m M_X$. Then for any $k\in [K^*]$ and $x^*\in \mb R^{D_X}$, we define  $ G_{[k],x^*}^{\dagger}(\cdot,\cdot)$ and $ v_{[k],x^*}^{\dagger}(\cdot,\cdot)$ as follows.
 \begin{enumerate}
     \item If  $x^*\in \mb N_{\varepsilon_j^x}^x$, and $\|x^*-x_k^*\|\leq \tau_2+2\varepsilon_j^x$, then considering the following local approximation to $G^*_\sk$ and  $v^*_\sk$:
 \begin{equation*}
 \begin{aligned}
       G_{[k],x^*}^{\dagger}(z,x)= \sum_{s=0}^{j} \sum_{\psi\in \wt\Psi_s^{d_Y}}\sum_{l\in \mb N_0^{D_X}\atop |l|<\beta_X} \int_{\mb R^{d_Y}} \frac{1}{l!}G_{[k]}^{*(\mathbf{0},l)}(t,x^*) (x-x^*)^{l}\psi(t)\,\dd t\cdot \psi(z) 
 \end{aligned}
 \end{equation*}
 and 
 \begin{equation*}
 \begin{aligned}
        v_{[k],x^*}^{\dagger}(z,x)=\sum_{s=0}^{j} \sum_{\psi\in \wt \Psi_s^{d_Y}}\sum_{l\in \mb N_0^{D_X}\atop |l|<\alpha_X} \int_{\mb R^{d_Y}} \frac{1}{l!}v_{[k]}^{*(\mathbf{0},l)}(t,x^*) (x-x^*)^{l}\psi(t)\,\dd t\cdot \psi(z),
 \end{aligned}
 \end{equation*}
 where recall $\wt \Psi_s^{d_Y}=\{\psi\in \ov\Psi_s^{d_Y}:\,{\rm supp}(\psi)\cap \mb B_{\mb R^{d_Y}}(\mathbf{0},\tau_1)\neq \emptyset\}$, and we use $G^{(\mathbf 0,l)}(z,x)$ to denote the partial derivative of $G(z,\cdot)$ of order $l$ evaluated at $x$.
 It holds that 
 \begin{equation}\label{eqn:diffG}
 \begin{aligned}
            &\underset{z\in \mb B_{\mb R^{d_Y}}(\mathbf{0},\tau_1)\atop x\in \mb B_{\m M_X}(x^*,2\varepsilon_j^x)}{\sup} \| G_{[k],x^*}^{\dagger}(z,x)-G^*_{[k]}(z,x)\|\\
             &\leq \underset{z\in \mb B_{\mb R^{d_Y}}(\mathbf{0},\tau_1)\atop x\in \mb B_{\m M_X}(x^*,2\varepsilon_j^x)}{\sup} \|  \sum_{s=0}^{j} \sum_{\psi\in \wt\Psi_s^{d_Y}}\ \int_{\mb R^{d_Y}} G_{[k]}^*(t,x)\psi(t)\,\dd t\cdot \psi(z)-G^*_{[k]}(z,x)\|\\
      &\quad  + \underset{z\in \mb B_{\mb R^{d_Y}}(\mathbf{0},\tau_1)\atop x\in \mb B_{\m M_X}(x^*,2\varepsilon_j^x)}{\sup} \| G_{[k],x^*}^{\dagger}(z,x)- \sum_{s=0}^{j} \sum_{\psi\in \wt\Psi_s^{d_Y}} \int_{\mb R^{d_Y}} G_{[k]}^*(t,x)\psi(t)\,\dd t\cdot \psi(z)\|\\
           & \lesssim 2^{-j\beta_Y}+ (\varepsilon_j^x)^{\beta_X} \cdot  \underset{z\in \mb B_{\mb R^{d_Y}}(\mathbf{0},\tau_1)}{\sup} \sum_{s=0}^j \sum_{\psi\in \wt\Psi_s^{d_Y}}\int_{\mb R^{d_Y}} |\psi(t)|\,\dd t \cdot \psi(z)   \lesssim 2^{-j\beta_Y}+ (\log n)\cdot(\varepsilon_j^x)^{\beta_X},
 \end{aligned}
 \end{equation}
 and similarly,
 \begin{equation}\label{eqn:diffv}
   \underset{z\in \mb B_{\mb R^{d_Y}}(\mathbf{0},\tau_1)\atop x\in \mb B_{\m M_X}(x^*,2\varepsilon_j^x)}{\sup} \| v_{[k],x^*}^{\dagger}(z,x)-v^*_{[k]}(z,x)\|\lesssim 2^{-j\alpha_Y}+ (\log n) \cdot (\varepsilon_j^x)^{\alpha_X}.
 \end{equation}
 \item If $x^*\notin \mb N_{\varepsilon_j^x}^x$, or  $x^*\in \mb N_{\varepsilon_j^x}^x$, but   $\|x^*-x_k^*\|> \tau_2+2\varepsilon_j^x$,  we define  $G_{[k],x^*}^{\dagger}(z,x)\equiv\mathbf 0_{D_Y}$ and $v_{[k],x^*}^{\dagger}(z,x)\equiv 0$.
 \end{enumerate}
  Let $\mb N^z_{c\,2^{-j}}$ be a $c\,2^{-j}$-covering set of $\mb B_{\mb R^{d_Y}}(\mathbf{0},\tau_1)$, contained within $\mb B_{\mb R^{d_Y}}(\mathbf{0},\tau_1)$, where $c$ is a small enough positive constant. For any $x^*\in \m M_X$,  denote
\begin{equation*}
 \begin{aligned}
          \Psi_{j}^{D_Y}(x^*)&=\{\psi\in \Psi_j^{D_Y}:\, \exists z^* \in \mb N^z_{c\,2^{-j}}, k\in [{K^*}], l\in \mb N_{0}^{D_Y}\text{ with }  |l|\leq \lfloor\wt \beta_X\rfloor \\
          &\qquad \text{so that }{\rm supp}(\psi^{(l)})\cap \mb B_{\mb R^{D_Y}}(G^*_{[k]}(z^*,x^*),C2^{-j})\neq \emptyset\},
 \end{aligned}
 \end{equation*}
where $C$ is a large enough constant. Then we present the following lemma that decompose    $\mb{E}_{\mu^*_{Y|x}}[2^{\frac{j(d_Y-D_Y)}{2}}\psi(y)]$ into summation of a term that depend on $G_{[k],x^*}^{\dagger},v_{[k],x^*}^{\dagger}$ and a polynomial term. 
\begin{lemma}\label{claimlemma11}
There exist constants $C_1,C_2$ such that for any $x^*\in \mb N_{\varepsilon_j^x}^x$, $\psi^*\in \Psi_j^{D_Y}$,  there exists coefficients $a^*_{\psi^*,x^*,s}\in (-C_1\, n,C_1\, n)$ indexed by  $s\in \mb N_0^{D_X}$ with $  s\leq \lfloor\wt\beta_X\rfloor^2+\lfloor \alpha_X\rfloor$, satisfying the following conditions:
  \begin{enumerate}
      \item It holds for any $x\in \mb B_{\m M_X}(x^*,2\varepsilon_j^x)$  that 
      \begin{equation*}
     \begin{aligned}
               & \Bigg|\sum_{k=1}^{K^*} \int_{\mb B_{\mb R^{d_Y}}(\mathbf{0},\tau_1)} 2^{\frac{j(d_Y-D_Y)}{2}}\psi^*(G_{[k]}^*(z,x)) v_\sk^*(z,x)\,\dd z\\
               & -\bigg(\sum_{k=1}^{K^*} \int_{\mb B_{\mb R^{d_Y}}(\mathbf{0},\tau_1)} 2^{\frac{j(d_Y-D_Y)}{2}}\psi^*(G_{[k],x^*}^\dagger(z,x)) v_{\sk,x^*}^\dagger(z,x)\,\dd z+\sum_{s\in \mb N_0^{D_X}\atop 0\leq s\leq \lfloor\wt\beta_X\rfloor^2+\lfloor \alpha_X\rfloor} a^{*}_{\psi^*,x^*,s}(x-x^*)^s\bigg)\Bigg|\\
           &\leq C_2\, (\log n)\cdot 2^{-\frac{jd_Y}{2}}(\varepsilon_j^x)^{\alpha_X}.
             \end{aligned}
 \end{equation*} 
    \item    If $ \psi^*\in \Psi_j^{D_Y} \setminus\Psi_{j}^{D_Y}(x^*)$, then it holds for any $x\in \mb B_{\m M_X}(x^*,\varepsilon_j^x)$  and $x'\in \mb B_{\mb N^x_{\varepsilon_j^x}}(x,2\varepsilon_j^x)$ that, 
 $$\sum_{k=1}^{K^*} \int_{\mb B_{\mb R^{d_Y}}(\mathbf{0},\tau_1)} 2^{\frac{j(d_Y-D_Y)}{2}}\psi^*(G_{[k]}^*(z,x)) v_\sk^*(z,x)\,\dd z=0$$
 and
 $$\sum_{k=1}^{K^*} \int_{\mb B_{\mb R^{d_Y}}(\mathbf{0},\tau_1)} 2^{\frac{j(d_Y-D_Y)}{2}}\psi^*(G_{[k],x'}^\dagger(z,x)) v_{\sk,x'}^\dagger(z,x)\,\dd z+\sum_{s\in \mb N_0^{D_X}\atop 0\leq s\leq \lfloor\wt\beta_X\rfloor^2+\lfloor \alpha_X\rfloor} a^{*}_{\psi^*,x',s}(x-x')^s=0.$$
  \end{enumerate}

\end{lemma}
\noindent The proof of Lemma~\ref{claimlemma11} is provided in Appendix~\ref{proof:claimlemma11}.
Then  since $\mathscr{I}_{j}$ is $c2^{-j}$ separated, let $c'=c/2$, for any $\iota,\iota'\in \mathscr{I}_{j}$, $\rho(\frac{|\iota-\iota'|}{c'2^{-j}})\neq 0$ if and only if $\iota=\iota'$. 
Applying Lemma~\ref{claimlemma11}, for any $\psi\in \Psi_j^{D_Y}$  and $x\in \m M_X$, if $\sum_{x^*\in \mb N^x_{\varepsilon_j^x}}\sum_{\psi^*\in  \Psi_{j}^{D_Y}(x^*)}\rho(\frac{\|x-x^*\|}{\varepsilon_j^x})\rho(\frac{|\m I_j(\psi)-\m I_j(\psi^*)|}{c'2^{-j}})\geq 1$, then

\begin{equation}\label{eqn:plemmaD31}
    \begin{aligned}
          &\sum_{k=1}^{K^*} \int_{\mb B_{\mb R^{d_Y}}(\mathbf{0},\tau_1)} 2^{\frac{j(d_Y-D_Y)}{2}}\psi(G_{[k]}^*(z,x)) v_\sk^*(z,x)\,\dd z\\
          &=\frac{\sum_{x^*\in \mb N^x_{\varepsilon_j^x}}\sum_{\psi^*\in  \Psi_{j}^{D_Y}(x^*)}\sum_{k=1}^{K^*} \int_{\mb B_{\mb R^{d_Y}}(\mathbf{0},\tau_1)} 2^{\frac{j(d_Y-D_Y)}{2}}\psi(G_{[k]}^*(z,x)) v_\sk^*(z,x)\,\dd z\rho(\frac{\|x-x^*\|}{\varepsilon_j^x})\rho(\frac{|\m I_j(\psi)-\m I_j(\psi^*)|}{c'2^{-j}})}{\sum_{x^*\in \mb N^x_{\varepsilon_j^x}}\sum_{\psi^*\in  \Psi_{j}^{D_Y}(x^*)}\rho(\frac{\|x-x^*\|}{\varepsilon_j^x})\rho(\frac{|\m I_j(\psi)-\m I_j(\psi^*)|}{c'2^{-j}})}\\
          &\leq\frac{\sum_{x^*\in \mb N^x_{\varepsilon_j^x}}\sum_{\psi^*\in  \Psi_{j}^{D_Y}(x^*)}  \sum_{k=1}^{K^*} \int_{\mb B_{\mb R^{d_Y}}(\mathbf{0},\tau_1)} 2^{\frac{j(d_Y-D_Y)}{2}}\psi(G_{[k],x^*}^\dagger(z,x)) v_{\sk,x^*}^\dagger(z,x) \,\dd z\rho(\frac{\|x-x^*\|}{\varepsilon_j^x})\rho(\frac{|\m I_j(\psi)-\m I_j(\psi^*)|}{c'2^{-j}})}{\sum_{x^*\in \mb N^x_{\varepsilon_j^x}}\sum_{\psi^*\in  \Psi_{j}^{D_Y}(x^*)}\rho(\frac{\|x-x^*\|}{\varepsilon_j^x})\rho(\frac{|\m I_j(\psi)-\m I_j(\psi^*)|}{c'2^{-j}})}\\
          &+\frac{\sum_{x^*\in \mb N^x_{\varepsilon_j^x}}\sum_{\psi^*\in  \Psi_{j}^{D_Y}(x^*)}  \sum_{s\in \mb N_0^{D_X}\atop 0\leq s\leq \lfloor\wt\beta_X\rfloor^2+\lfloor \alpha_X\rfloor} a^*_{\psi^*,x^*,s}(x-x^*)^s\rho(\frac{\|x-x^*\|}{\varepsilon_j^x})\rho(\frac{|\m I_j(\psi)-\m I_j(\psi^*)|}{c'2^{-j}})}{\sum_{x^*\in \mb N^x_{\varepsilon_j^x}}\sum_{\psi^*\in  \Psi_{j}^{D_Y}(x^*)}\rho(\frac{\|x-x^*\|}{\varepsilon_j^x})\rho(\frac{|\m I_j(\psi)-\m I_j(\psi^*)|}{c'2^{-j}})}+C_2(\log n)\cdot(\varepsilon_j^x)^{\alpha_X}2^{-\frac{d_Yj}{2}}\\
           &=\frac{\sum_{x^*\in \mb N^x_{\varepsilon_j^x}}\sum_{\psi^*\in  \Psi_{j}^{D_Y}(x^*)}  \sum_{k=1}^{K^*} \int_{\mb B_{\mb R^{d_Y}}(\mathbf{0},\tau_1)} 2^{\frac{j(d_Y-D_Y)}{2}}\psi(G_{[k],x^*}^\dagger(z,x)) v_{\sk,x^*}^\dagger(z,x) \,\dd z\rho(\frac{\|x-x^*\|}{\varepsilon_j^x})\rho(\frac{|\m I_j(\psi)-\m I_j(\psi^*)|}{c'2^{-j}})}{\sum_{x^*\in \mb N^x_{\varepsilon_j^x}}\sum_{\psi^*\in  \Psi_{j}^{D_Y}(x^*)}\rho(\frac{\|x-x^*\|}{\varepsilon_j^x})\rho(\frac{|\m I_j(\psi)-\m I_j(\psi^*)|}{c'2^{-j}})+\frac{1}{n^2}}\\
          &+\frac{\sum_{x^*\in \mb N^x_{\varepsilon_j^x}}\sum_{\psi^*\in  \Psi_{j}^{D_Y}(x^*)}  \sum_{s\in \mb N_0^{D_X}\atop 0\leq s\leq \lfloor\wt\beta_X\rfloor^2+\lfloor \alpha_X\rfloor} a^*_{\psi^*,x^*,s}(x-x^*)^s\rho(\frac{\|x-x^*\|}{\varepsilon_j^x})\rho(\frac{|\m I_j(\psi)-\m I_j(\psi^*)|}{c'2^{-j}})}{\sum_{x^*\in \mb N^x_{\varepsilon_j^x}}\sum_{\psi^*\in  \Psi_{j}^{D_Y}(x^*)}\rho(\frac{\|x-x^*\|}{\varepsilon_j^x})\rho(\frac{|\m I_j(\psi)-\m I_j(\psi^*)|}{c'2^{-j}})+\frac{1}{n^2}}\\
          &+\underbrace{C_2(\log n)\cdot(\varepsilon_j^x)^{\alpha_X}2^{-\frac{d_Yj}{2}}+\frac{C_3}{n}}_{=\m O\big(\log n\cdot(\varepsilon_j^x)^{\alpha_X} 2^{-\frac{d_Yj}{2}}\big)}.\\
    \end{aligned}
\end{equation}
On the other hand, if $\sum_{x^*\in \mb N^x_{\varepsilon_j^x}}\sum_{\psi^*\in  \Psi_{j}^{D_Y}(x^*)}\rho(\frac{\|x-x^*\|}{\varepsilon_j^x})\rho(\frac{|\m I_j(\psi)-\m I_j(\psi^*)|}{c'2^{-j}})< 1$.    Since $\mb N_{\varepsilon_j^x}^x$ is the largest $\varepsilon_j^x$-packing of $\m M_X$, there exists $x^*\in \mb N_{\varepsilon_j^x}^x$ so that $\|x-x^*\|\leq \varepsilon_j^x$ and $\rho(\frac{\|x-x^*\|}{\varepsilon_j^x})=1$. Moreover, since $$\sum_{\psi^*\in  \Psi_{j}^{D_Y}(x^*)}\rho(\frac{\|x-x^*\|}{\varepsilon_j^x})\rho(\frac{|\m I_j(\psi)-\m I_j(\psi^*)|}{c'2^{-j}})\leq  \sum_{x^*\in \mb N^x_{\varepsilon_j^x}}\sum_{\psi^*\in  \Psi_{j}^{D_Y}(x^*)}\rho(\frac{\|x-x^*\|}{\varepsilon_j^x})\rho(\frac{|\m I_j(\psi)-\m I_j(\psi^*)|}{c'2^{-j}})< 1,$$ it holds that $\psi\in  \Psi_j^{D_Y} \setminus\Psi_{j}^{D_Y}(x^*)$. Therefore,
\begin{equation}\label{eqn:32}
    \sum_{k=1}^{K^*} \int_{\mb B_{\mb R^{d_Y}}(\mathbf{0},\tau_1)} 2^{\frac{j(d_Y-D_Y)}{2}}\psi(G_{[k]}^*(z,x)) v_\sk^*(z,x)\,\dd z=0,
\end{equation}
 and for any  $x'\in \mb B_{\mb N^x_{\varepsilon_j^x}}(x,2\varepsilon_j^x)$
 \begin{equation}\label{eqn:33}
     \begin{aligned}
        \sum_{k=1}^{K^*} \int_{\mb B_{\mb R^{d_Y}}(\mathbf{0},\tau_1)} 2^{\frac{j(d_Y-D_Y)}{2}}\psi(G_{[k],x'}^\dagger(z,x)) v_{\sk,x'}^\dagger(z,x)\,\dd z+\sum_{s\in \mb N_0^{D_X}\atop 0\leq s\leq \lfloor\wt\beta_X\rfloor^2+\lfloor \alpha_X\rfloor} a^{*}_{\psi,x',s}(x-x')^s=0.
     \end{aligned}
 \end{equation}
Hence,
 \begin{equation*}
    \begin{aligned}
      &\sum_{k=1}^{K^*} \int_{\mb B_{\mb R^{d_Y}}(\mathbf{0},\tau_1)} 2^{\frac{j(d_Y-D_Y)}{2}}\psi(G_{[k]}^*(z,x)) v_\sk^*(z,x)\,\dd z=0\\
          &=\frac{\sum_{x^*\in \mb N^x_{\varepsilon_j^x}}\sum_{\psi^*\in  \Psi_{j}^{D_Y}(x^*)}  \sum_{k=1}^{K^*} \int_{\mb B_{\mb R^{d_Y}}(\mathbf{0},\tau_1)} 2^{\frac{j(d_Y-D_Y)}{2}}\psi(G_{[k],x^*}^\dagger(z,x)) v_{\sk,x^*}^\dagger(z,x) \,\dd z\rho(\frac{\|x-x^*\|}{\varepsilon_j^x})\rho(\frac{|\m I_j(\psi)-\m I_j(\psi^*)|}{c'2^{-j}})}{\sum_{x^*\in \mb N^x_{\varepsilon_j^x}}\sum_{\psi^*\in  \Psi_{j}^{D_Y}(x^*)}\rho(\frac{\|x-x^*\|}{\varepsilon_j^x})\rho(\frac{|\m I_j(\psi)-\m I_j(\psi^*)|}{c'2^{-j}})+\frac{1}{n^2}}\\
          &+\frac{\sum_{x^*\in \mb N^x_{\varepsilon_j^x}}\sum_{\psi^*\in  \Psi_{j}^{D_Y}(x^*)}  \sum_{s\in \mb N_0^{D_X}\atop 0\leq s\leq \lfloor\wt\beta_X\rfloor^2+\lfloor \alpha_X\rfloor} a^*_{\psi^*,x^*,s}(x-x^*)^s\rho(\frac{\|x-x^*\|}{\varepsilon_j^x})\rho(\frac{|\m I_j(\psi)-\m I_j(\psi^*)|}{c'2^{-j}})}{\sum_{x^*\in \mb N^x_{\varepsilon_j^x}}\sum_{\psi^*\in  \Psi_{j}^{D_Y}(x^*)}\rho(\frac{\|x-x^*\|}{\varepsilon_j^x})\rho(\frac{|\m I_j(\psi)-\m I_j(\psi^*)|}{c'2^{-j}})+\frac{1}{n^2}}.\\
    \end{aligned}
\end{equation*}
 Let  $  W_j(x^*)=  |\Psi_{j}^{D_Y}(x^*)|$, we have $\max_{x^*\in \mb N_{\varepsilon_j^x}^x}  W_j(x^*)\leq W_j=C_3\, (\varepsilon_j^y)^{-d_Y} $ and  $ |\mb N^x_{\varepsilon_j^x}|\leq W_j'= C_2\, (\varepsilon_j^x)^{-d_X}$ when $C_2,C_3$ are sufficiently large.
Let $\m X_j$ be an arbitrary subset of $\mb B_{\mb R^{D_X}}(\mathbf{0},L_1)\setminus \cup_{x\in \m M_X}\mb B_{\mb R^{D_X}}(x,2\varepsilon_j^x)$ so that the points in $\m  X_j$ are $\varepsilon_j^x$-separated and $|\m X_j|=W_j'-|\mb N^x_{\varepsilon_j^x}|$ (note that such a set $\m X_j$ exist if $L_1$ is sufficiently large). Arrange the points in $\mb N^x_{\varepsilon_j^x}$ and $\m X_j$ as $\mb N^x_{\varepsilon_j^x}=(x_{j1}, x_{j2},\cdots, x_{j|\mb N^x_{\varepsilon_j^x}|})$ and $\m X_j=(z_{j1},z_{j2},\cdots,z_{j|\m X_j|})$, we denote 
 \begin{equation*}
     x_{jl}^*=\left\{\begin{array}{ll}
       x_{jl},   & \text{ for }l\in \{1,2,\cdots, |\mb N^x_{\varepsilon_j^x}|\} \\
        z_{jl_1} \text{ with } l_1=l- |\mb N^x_{\varepsilon_j^x}|, & \text{ for }l\in \{|\mb N^x_{\varepsilon_j^x}|+1,|\mb N^x_{\varepsilon_j^x}|+2,\cdots, W_j'\}.
     \end{array}\right.
 \end{equation*}
 Furthermore, denote   $\Psi_{j}^{D_Y}(x^*)=( \psi_{j1}^{x^*},\psi_{j2}^{x^*},\cdots,\psi_{j W_j(x^*)}^{x^*})$. We define 
  \begin{equation*}
     e^*_{ji_1i_2}=\left\{
     \begin{array}{ll}
      \m I_j(\psi_{ji_1}^{x_{ji_2}^*}),    & \text{ if }i_1\leq W_j(x_{ji_2}^*) \text{ and }i_2\leq |\mb N^x_{\varepsilon_j^x}| \\
      (2,2,\cdots,2),  & \text{ if } W_j(x_{ji_2}^*)<i_1 \leq W_j\text{ or }i_2> |\mb N^x_{\varepsilon_j^x}|, 
     \end{array}
     \right.
 \end{equation*}
 \begin{equation*}
     c^*_{ji_1i_2l}=\left\{
     \begin{array}{ll}
    a^*_{\psi_{ji_1}^{x_{ji_2}^*},x_{ji_2}^*,l},   & \text{ if }i_1\leq W_j(x_{ji_2}^*) \text{ and }i_2\leq |\mb N^x_{\varepsilon_j^x}|\\
     0,  & \text{ if } W_j(x_{ji_2}^*)<i_1 \leq W_j \text{ or }i_2> |\mb N^x_{\varepsilon_j^x}|,
     \end{array}
     \right.
 \end{equation*}
 and
 \begin{equation*} 
    \begin{aligned}
    &S^*_j(\psi,x)=\m T_{C_1\,2^{-\frac{d_Yj}{2}}}\Big(\\
   & \frac{\sum_{i_1=1}^{W_j}\sum_{i_2=1}^{W_j'} \sum_{k=1}^{K^*} \int_{\mb B_{\mb R^{d_Y}}(\mathbf{0},\tau_1)} 2^{\frac{j(d_Y-D_Y)}{2}}\psi(G^\dagger_{\sk,x_{ji_2}^*}(z,x))v^\dagger_{\sk,x_{ji_2}^*}(z,x)\,\dd z \rho(\frac{\|x-x^*_{ji_2}\|}{\varepsilon_j^x})\rho(\frac{|\m I_j(\psi)-e^*_{ji_1i_2}|}{\varepsilon_j^y})}{\sum_{i_1=1}^{W_j}\sum_{i_2=1}^{W_j'}\rho(\frac{\|x-x^*_{ji_2}\|}{\varepsilon_j^x})\rho(\frac{|\m I_j(\psi)-e^*_{ji_1i_2}|}{\varepsilon_j^y})+\frac{1}{n^2}}\\
   &+\frac{\sum_{i_1=1}^{W_j}\sum_{i_2=1}^{W_j'}\sum_{l\in\mb N_{0}^{D_X}\atop |l|\leq\lfloor\wt\beta_X\rfloor^2+\lfloor \alpha_X\rfloor} c^*_{ji_1i_2l} (x-x^*_{ji_2})^l \rho(\frac{\|x-x^*_{ji_2}\|}{\varepsilon_j^x})\rho(\frac{|\m I_j(\psi)-e^*_{ji_1i_2}|}{\varepsilon_j^y})}{\sum_{i_1=1}^{W_j}\sum_{i_2=1}^{W_j'}\rho(\frac{\|x-x^*_{ji_2}\|}{\varepsilon_j^x})\rho(\frac{|\m I_j(\psi)-e^*_{ji_1i_2}|}{\varepsilon_j^y})+\frac{1}{n^2}}\Big)\in \m S_j^\dagger.\\
   \end{aligned}
\end{equation*}
 Then for any $x\in \m M_X$, denote
\begin{equation*}
    \wt\Psi_j^{D_Y}(x)=\{\psi\in \Psi_j^{D_Y}:\, \sum_{x^*\in \mb N^x_{\varepsilon_j^x}}\sum_{\psi^*\in  \Psi_{j}^{D_Y}(x^*)}\rho(\frac{\|x-x^*\|}{\varepsilon_j^x})\rho(\frac{|\m I_j(\psi)-\m I_j(\psi^*)|}{c'2^{-j}}\geq 1\}.
\end{equation*}
 We have $\sup_{x\in \m M_X}|\wt\Psi_j^{D_Y}(x)|=\m O(2^{jd_Y})$ and by~\eqref{eqn:32},~\eqref{eqn:33}, it holds for any $\psi\in  \Psi_j^{D_Y}\setminus  \wt\Psi_j^{D_Y}(x)$ that, 
 \begin{equation*}
 \begin{aligned}
      &\mb{E}_{\mu^*_{Y|x}}[2^{\frac{j(d_Y-D_Y)}{2}}\psi(y)]=\sum_{k=1}^{K^*} \int_{\mb B_{\mb R^{d_Y}}(\mathbf{0},\tau_1)} 2^{\frac{j(d_Y-D_Y)}{2}}\psi(G_{[k]}^*(z,x)) v_\sk^*(z,x)\,\dd z=S^*_j(\psi,x)=0.\\
    \end{aligned}
\end{equation*}
Furthermore, since $\sup_{x\in \m M_X}\sup_{\psi\in \Psi_j^{D_Y}}\mb{E}_{\mu^*_{Y|x}}[2^{\frac{j(d_Y-D_Y)}{2}}\psi(y)]\leq C_1\, 2^{-\frac{d_Yj}{2}}$,  we can get for any $x\in \m M_X$ and $\psi\in \Psi_j^{D_Y}$, 
 \begin{equation*} 
    \begin{aligned}
    &S^*_j(\psi,x)=\\
   & =\frac{\sum_{x^*\in \mb N^x_{\varepsilon_j^x}}\sum_{\psi^*\in  \Psi_{j}^{D_Y}(x^*)}  \sum_{k=1}^{K^*} \int_{\mb B_{\mb R^{d_Y}}(\mathbf{0},\tau_1)} 2^{\frac{j(d_Y-D_Y)}{2}}\psi(G_{[k],x^*}^\dagger(z,x)) v_{\sk,x^*}^\dagger(z,x) \,\dd z\rho(\frac{\|x-x^*\|}{\varepsilon_j^x})\rho(\frac{|\m I_j(\psi)-\m I_j(\psi^*)|}{c'2^{-j}})}{\sum_{x^*\in \mb N^x_{\varepsilon_j^x}}\sum_{\psi^*\in  \Psi_{j}^{D_Y}(x^*)}\rho(\frac{\|x-x^*\|}{\varepsilon_j^x})\rho(\frac{|\m I_j(\psi)-\m I_j(\psi^*)|}{c'2^{-j}})+\frac{1}{n^2}}\\
          &+\frac{\sum_{x^*\in \mb N^x_{\varepsilon_j^x}}\sum_{\psi^*\in  \Psi_{j}^{D_Y}(x^*)}  \sum_{l\in \mb N_0^{D_X}\atop 0\leq s\leq \lfloor\wt\beta_X\rfloor^2+\lfloor \alpha_X\rfloor} a^*_{\psi^*,x^*,l}(x-x^*)^l\rho(\frac{\|x-x^*\|}{\varepsilon_j^x})\rho(\frac{|\m I_j(\psi)-\m I_j(\psi^*)|}{c'2^{-j}})}{\sum_{x^*\in \mb N^x_{\varepsilon_j^x}}\sum_{\psi^*\in  \Psi_{j}^{D_Y}(x^*)}\rho(\frac{\|x-x^*\|}{\varepsilon_j^x})\rho(\frac{|\m I_j(\psi)-\m I_j(\psi^*)|}{c'2^{-j}})+\frac{1}{n^2}}.
   \end{aligned}
\end{equation*}
Therefore, using bound~\eqref{eqn:plemmaD31}, we have
 \begin{equation*}
 \begin{aligned}
      &\mb{E}_{\mu^*_X}\Big[ \sum_{\psi\in \Psi_j^{D_Y}} ( \mb{E}_{\mu^*_{Y|X}}[2^{\frac{j(d_Y-D_Y)}{2}}\psi(y)]-S^*_j(\psi,x))^2\Big]\\
      &= \mb{E}_{\mu^*_X}\Big[ \sum_{\psi\in \Psi_j^{D_Y}} \Big( \sum_{k=1}^{K^*} \int_{\mb B_{\mb R^{d_Y}}(\mathbf{0},\tau_1)} 2^{\frac{j(d_Y-D_Y)}{2}}\psi(G_{[k]}^*(z,x)) v_\sk^*(z,x)\,\dd z-S^*_j(\psi,x)\Big)^2\Big]\\
      &=\mb{E}_{\mu^*_X}\Big[ \sum_{\psi\in  \wt\Psi_j^{D_Y}(x)} \Big( \sum_{k=1}^{K^*} \int_{\mb B_{\mb R^{d_Y}}(\mathbf{0},\tau_1)} 2^{\frac{j(d_Y-D_Y)}{2}}\psi(G_{[k]}^*(z,x)) v_\sk^*(z,x)\,\dd z-S^*_j(\psi,x)\Big)^2\Big]\\
 &\lesssim (\log n)^2\cdot (\varepsilon_j^x)^{2\alpha_X} 2^{-jd_Y} \sup_{x\in \m M_X}|\wt\Psi_j^{D_Y}(x)|\lesssim (\log n)^2\cdot (\varepsilon_j^x)^{2\alpha_X},
 \end{aligned}
\end{equation*}
which completes the proof.

  \subsection{Proof of Lemma~\ref{lemma:manifoldlearningr3}}\label{proof:lemma:manifoldlearningr3}
  We will use Lemma~\ref{lemma:reconstruct}
   to show the desired results. Denote 
    \begin{equation*}
  {\Psi}_j^{D_X}=\{\psi\in \ov\Psi_j^{D_X}:\, {\rm supp}(\psi)\cap \m M_x\neq \emptyset\},
 \end{equation*}
 and consider the family
\begin{equation}\label{def:wtG}
 \begin{aligned}
    \wt {\m G}&=\{G(z,x)=\sum_{j_1=0}^{J_1}\sum_{j_2=0}^{J_2}\sum_{\psi_1\in  \Psi_{j_1}^{d_Y} }\sum_{\psi_2\in {\Psi}_{j_2}^{D_X}}g_{\psi_1\psi_2}\psi_1(z)\psi_2(x)\,: \\
   & |g_{\psi_1\psi_2}|\leq L_1\, 2^{-\frac{d_Yj_1+D_Xj_2}{2}-((j_1\beta_Y)\vee (j_2\beta_X))} \text{ for }\psi_1\in \Psi_{j_1}^{d_Y}, \psi_2\in {\Psi}_{j_2}^{D_X} \},
 \end{aligned}
 \end{equation}
 where $J_1= \lceil\log_2 (n^{-\frac{1}{d_Y+d_X\frac{\beta_Y}{\beta_X}}})\rceil$, $J_2= \lceil\log_2 (n^{-\frac{1}{d_X+d_Y\frac{\beta_X}{\beta_Y}}})\rceil$. Since for any $z\in \mb R^{d_Y}$ and $x\in \m M_X$, it holds that  
 \begin{equation*}
     \begin{aligned}
       & \sum_{j_1=0}^{J_1}\sum_{j_2=0}^{J_2}\sum_{\psi_1\in  \Psi_{j_1}^{d_Y} }\sum_{\psi_2\in \ov {\Psi}_{j_2}^{D_X}}g_{\psi_1\psi_2}\psi_1(z)\psi_2(x) \\
        &= \sum_{j_1=0}^{J_1}\sum_{j_2=0}^{J_2}\sum_{\psi_1\in  \Psi_{j_1}^{d_Y} }\sum_{\psi_2\in  {\Psi}_{j_2}^{D_X}}g_{\psi_1\psi_2}\psi_1(z)\psi_2(x), 
     \end{aligned}
 \end{equation*}
we can obtain  
\begin{equation*} 
\begin{aligned}
     (\wh G_{[k]}, \wh V_{[k]})&=\underset{G\in  {\m G}\atop V\in \mb O(D_Y,d_Y)}{\arg\min} \frac{1}{|I_1|}\sum_{i\in I_1} \|Y_i-G(V^T(Y_i-y_k),X_i)\|^2\mathbf{1}(X_i\in \mb B_{\mb R^{D_X}}(x_k,2\tau_2))\mathbf{1}(Y_i\in \mb B_{\mb R^{D_Y}}(y_k,2\tau_1))\\
     &=\underset{G\in  \wt{\m G}\atop V\in \mb O(D_Y,d_Y)}{\arg\min} \frac{1}{|I_1|}\sum_{i\in I_1} \|Y_i-G(V^T(Y_i-y_k),X_i)\|^2\mathbf{1}(X_i\in \mb B_{\mb R^{D_X}}(x_k,2\tau_2))\mathbf{1}(Y_i\in \mb B_{\mb R^{D_Y}}(y_k,2\tau_1)).
\end{aligned}
\end{equation*}
Furthermore,  since $\beta_Y\geq 2$, we have the following smoothness property for functions in $\wt{\m G}$, the proof of which is given in Appendix~\ref{proof:{lemma: smoothG}}.
\begin{lemma}\label{lemma: smoothG}
 With the choice of $\wt {\m G}$ in~\eqref{def:wtG}, there exists a constant $L_1$ so that for any $G\in \wt{\m G}$, $x\in \mb R^{D_X}$ and $z\in \mb R^{d_Y}$, 
  \begin{equation*}
       \| J_{G(\cdot,x)}(z)\|_{F}\leq L_1.
  \end{equation*}
  Moreover, for any $1<\beta<2$, there exists a constant $L_{\beta}$ so that for any $G\in  \wt{\m G}$, $x\in \mb R^{D_X}$ and $z,z'\in \mb R^{d_Y}$
  \begin{equation*}
       \| J_{G(\cdot,x)}(z)-J_{G(\cdot,x)}(z')\|_{F}\leq L_{\beta}\|z-z'\|^{\beta-1}.
  \end{equation*}
\end{lemma}
\noindent Moreover, we can derive the following lemma that control the covering number of $\wt{\m G}$.
\begin{lemma}\label{lemma:coveringG}
 With the choice of $\wt {\m G}$ in~\eqref{def:wtG}, there exists a constant $C_1$ so that for any $0<\gamma_1\leq 1$, the $\varepsilon$-covering number $\mathbf N(\wt{\m G},d_{\infty}^{\gamma_1},\varepsilon)$  of $\wt{\m G}$ with respect to the $d_{\infty}^{\gamma_1}$ distance, defined as $d^{\gamma_1}_{\infty}(G_1,G_2)=\underset{z\in \mb R^{d_Y},x\in \mb R^{D_X}}{\sup}\|G_1(z,x)-G_2(z,x)\|^{\gamma_1}$, satisfies
 \begin{equation*}
 \begin{aligned}
    &   \log \mathbf N(\wt{\m G},d^{\gamma_1}_{\infty},\varepsilon) \\
       &\leq \left\{\begin{array}{cc}
   C_1\sum_{j_1=0}^{J_1}\sum_{j_2=0}^{J_2} 2^{d_Yj_1+d_Xj_2}\log\left(\frac{C_1 (J_1+J_2)2^{-\frac{d_Yj_1+d_Xj_2}{4\gamma_1}-\frac{(j_1\beta_Y)\vee (j_2\beta_X)}{2}}}{\varepsilon^{\frac{1}{\gamma_1}}}\vee 1\right)    &  \frac{d_Y}{\beta_1}+\frac{d_X}{\beta_2}\leq 2\gamma_1, \\
      C_1\sum_{j_1=0}^{J_1}\sum_{j_2=0}^{J_2} 2^{d_Yj_1+d_Xj_2}\log \left(\frac{C_1((J_1+J_2)\wedge c(\beta_Y,\beta_X,d_Y,d_X,\gamma_1))2^{-((j_1\beta_Y)\vee (j_2\beta_X))}}{\varepsilon^{\frac{1}{\gamma_1}}  s_{j_1j_2}}\vee 1\right)    &  \frac{d_Y}{\beta_1}+\frac{d_X}{\beta_2}> 2\gamma_1,
  \end{array}
  \right.
 \end{aligned}
 \end{equation*}
 where $c(\beta_Y,\beta_X,d_Y,d_X,\gamma_1)=\frac{2^{\frac{d_Y+d_X\frac{\beta_Y}{\beta_X}-2\beta_Y\gamma_1}{4\gamma_1}}}{2^{\frac{d_Y+d_X\frac{\beta_Y}{\beta_X}-2\beta_Y\gamma_1}{4\gamma_1}}-1}+\frac{2^{\frac{d_X+d_Y\beta_X/\beta_Y-2\beta_X\gamma_1}{4\gamma_1}}}{2^{\frac{d_X+d_Y\frac{\beta_X}{\beta_Y}-2\beta_X\gamma_1}{4\gamma_1}}-1}$  and $ s_{j_1j_2}=\sqrt{\frac{2^{\frac{d_Yj_1+d_Xj_2}{2\gamma_1}-(j_1\beta_Y\vee j_2\beta_X)}}{2^{\frac{d_YJ_1+d_XJ_2}{2\gamma_1}-(J_1\beta_Y\vee J_2\beta_X)}}}$.
\end{lemma}

\noindent The proof of Lemma~\ref{lemma:coveringG} is provided in Appendix~\ref{proof:lemma:coveringG}. Then we can bound the integral $\frac{1}{\sqrt{n}}\int_0^{\infty} \sqrt{ \log \mathbf N(\wt{\m G}, d_{\infty}^{\gamma_1},\varepsilon)}\dd \varepsilon$. 
When $\frac{d_Y}{\beta_Y}+\frac{d_X}{\beta_X}\leq 2\gamma_1$, we have 
\begin{equation*}
    \begin{aligned}
         &\frac{1}{\sqrt{n}}\int_0^{\infty} \sqrt{ \log \mathbf N(\wt{\m G}, d_{\infty}^{\gamma_1},\varepsilon)}\dd \varepsilon\\
         &\lesssim \frac{1}{\sqrt{n}}\int_0^{\infty} \sqrt{\sum_{j_1=0}^{J_1}\sum_{j_2=0}^{J_2} \log \left(\frac{C_1 (J_1+J_2)2^{-\frac{d_Yj_1+d_Xj_2}{4\gamma_1}-\frac{(j_1\beta_Y)\vee (j_2\beta_X)}{2}}}{\varepsilon^{1/\gamma_1}}\vee 1\right) 2^{d_Yj_1+d_Xj_2}}\,\dd \varepsilon\\
      &\lesssim \frac{1}{\sqrt{n}} \sum_{j_1=0}^{J_1}\sum_{j_2=0}^{J_2} \int_0^{\infty}\sqrt{\log \left(\frac{C_1 (J_1+J_2)2^{-\frac{d_Yj_1+d_Xj_2}{4\gamma_1}-\frac{(j_1\beta_Y)\vee (j_2\beta_X)}{2}}}{\varepsilon^{1/\gamma_1}}\vee 1\right) 2^{d_Yj_1+d_Xj_2}}\,\dd \varepsilon\\
 &\lesssim \frac{1}{\sqrt{n}} \sum_{j_1=0}^{J_1}\sum_{j_2=0}^{J_2} \\
 &\qquad\int_0^{\big(C_1 (J_1+J_2)2^{-\frac{d_Yj_1+d_Xj_2}{4\gamma_1}-\frac{(j_1\beta_Y)\vee (j_2\beta_X)}{2}}\big)^{\gamma_1}}\sqrt{\log \left(\frac{C_1 (J_1+J_2)2^{-\frac{d_Yj_1+d_Xj_2}{4\gamma_1}-\frac{(j_1\beta_Y)\vee (j_2\beta_X)}{2}}}{\varepsilon^{1/\gamma_1}}\right) 2^{d_Yj_1+d_Xj_2}}\,\dd \varepsilon\\
&\lesssim \frac{1}{\sqrt{n}} \sum_{j_1=0}^{J_1}\sum_{j_2=0}^{J_2} (J_1+J_2)^{\gamma_1}2^{\frac{d_Yj_1+d_Xj_2}{4}-\gamma_1\frac{(j_1\beta_Y)\vee (j_2\beta_X)}{2}} \lesssim \frac{(\log n)^{1+\gamma_1}}{\sqrt{n}}.
    \end{aligned}
\end{equation*}
When $\frac{d_Y}{\beta_Y}+\frac{d_X}{\beta_X}>2\gamma_1$, we have 
\begin{equation*}
    \begin{aligned}
         &\frac{1}{\sqrt{n}}\int_0^{\infty} \sqrt{ \log \mathbf N(\wt{\m G}, d_{\infty}^{\gamma_1},\varepsilon)}\dd \varepsilon\\
         &\lesssim \frac{1}{\sqrt{n}}\int_0^{\infty} \sqrt{\sum_{j_1=0}^{J_1}\sum_{j_2=0}^{J_2} \log \left(\frac{C_1((J_1+J_2)\wedge c(\beta_Y,\beta_X,d_Y,d_X,\gamma_1))2^{-((j_1\beta_Y)\vee (j_2\beta_X))}}{\varepsilon^{1/\gamma_1} s_{j_1j_2}}\vee 1\right) 2^{d_Yj_1+d_Xj_2}}\,\dd \varepsilon\\
      &\lesssim \frac{1}{\sqrt{n}} \sum_{j_1=0}^{J_1}\sum_{j_2=0}^{J_2} \int_0^{\infty}\sqrt{ \log \left(\frac{C_1((J_1+J_2)\wedge c(\beta_Y,\beta_X,d_Y,d_X,\gamma_1))2^{-((j_1\beta_Y)\vee (j_2\beta_X))}}{\varepsilon^{1/\gamma_1} s_{j_1j_2}}\vee 1\right) 2^{d_Yj_1+d_Xj_2}}\,\dd \varepsilon\\
&\lesssim \frac{1}{\sqrt{n}} \sum_{j_1=0}^{J_1}\sum_{j_2=0}^{J_2} \frac{1}{s_{j_1j_2}^{\gamma_1}}((J_1+J_2)\wedge c(\beta_Y,\beta_X,d_Y,d_X,\gamma_1))^{\gamma_1}2^{-\gamma_1((j_1\beta_Y)\vee (j_2\beta_X))} 2^{\frac{d_Yj_1+d_Xj_2}{2}}\\
&\lesssim  ((J_1+J_2)\wedge c(\beta_Y,\beta_X,d_Y,d_X,\gamma_1))^{1+\gamma_1}  \frac{1}{\sqrt{n}} 2^{-\gamma_1((J_1\beta_Y)\vee (J_2\beta_X))} 2^{\frac{d_YJ_1+d_XJ_2}{2}}\\
&\lesssim   ((J_1+J_2)\wedge c(\beta_Y,\beta_X,d_Y,d_X,\gamma_1))^{1+\gamma_1} n^{-\frac{\gamma_1}{\frac{d_X}{\beta_X}+\frac{d_Y}{\beta_Y}}}\lesssim (\log n \wedge \frac{1}{\beta_Y(d_Y/\beta_Y+d_X/\beta_X-2\gamma_1)} )^{1+\gamma_1}\cdot n^{-\frac{\gamma_1}{\frac{d_X}{\beta_X}+\frac{d_Y}{\beta_Y}}}.
    \end{aligned}
\end{equation*}
Then it remains to bound the term  $\varepsilon^*$ in  Lemma~\ref{lemma:reconstruct}.  Fix an arbitrary $k\in [K]$. If $\mb B_{\mb R^{D_X+D_Y}}((x_k,y_k),\sqrt{2}\tau_2)\cap \m M=\emptyset$, then $k\notin \wh{\m K}$. Otherwise, there exists $(x_k^*,y_k^*)\in \mb B_{\m M}((x_k,y_k),\sqrt{2}\tau_2)$. Let $V_{\sk}^*$ be an arbitrary  orthonormal basis of $T_{\m M_{Y|x_{k}^*}}y_k^*$. Denote $Q_\sk^*(y)=(V_\sk^*)^T(y-y_k)$ and $G_\sk^*(z,x)= \Phi_{(x_k^*,y_k^*)}(V_\sk^*(z+(V_\sk^*)^T(y_k-y_k^*)),x)$. Then there exists  $\ov G_\sk^*\in  \ov{\m H}^{\beta_Y,\beta_X}_L(\mb R^{d_Y}, \mb R^{D_X})$  so that $\ov G_\sk^*(z,x)$ and $G_\sk^*(z,x)$ coincide within $\mb B_{\mb R^{d_Y}}((V_\sk^*)^T(y_k^*-y_k),\tau_1)\times \mb B_{\m M_X}(x_k^*,\tau)$. Moreover, for any $(x,y)\in \m M$ with $\|x-x_k^*\|<\tau$ and $\|y-y_k^*\|<\tau_1$, it holds that $y=\ov G_\sk^*(Q_\sk^*(y),x)$.  Then   Let 
 \begin{equation*}
      G_\sk^\dagger(z,x)=\sum_{j_1=0}^{J_1}\sum_{j_2=0}^{J_2} \sum_{\psi_1\in  \Psi_{j_1}^{d_Y} }\sum_{\psi_2\in {\Psi}_{j_2}^{D_X}} g_{\sk,\psi_1,\psi_2}^* \psi_1(z)\psi_2(x), \quad  g_{\sk,\psi_1,\psi_2}^* =\int_{\mb R^{D_X}} \int_{\mb R^{d_Y}} \ov G_\sk^*(z,x) \psi_1(z)\psi_2(x)\,\dd z\dd x.
 \end{equation*}
 It holds that $ G_\sk^\dagger\in \wt{\m G}$. Moreover, by leveraging the wavelet approximation for $  \ov{\m H}^{\beta_Y,\beta_X}$-smooth functions as described in Lemma~\ref{le:approwaveletsmooth}, and setting $J_1= \lceil\log_2 (n^{-\frac{1}{d_Y+d_X\frac{\beta_Y}{\beta_X}}})\rceil$, $J_2= \lceil\log_2 (n^{-\frac{1}{d_X+d_Y\frac{\beta_X}{\beta_Y}}})\rceil$ and $\tau_2<\frac{\tau_1\wedge \tau}{4}$, there exists a constant $C$ such that for any $x\in  \mb B_{\m M_X}(x_k,2\tau_2)$ and $z\in \mb B_{R^{d_Y}}(\mathbf{0},2\tau_2)$,
 \begin{equation*}
    \|G_\sk^*(z,x)-  G_\sk^\dagger(z,x)\|\leq C\, (\log n)\cdot n^{-\frac{1}{\frac{d_Y}{\beta_Y}+\frac{d_X}{\beta_X}}}.
 \end{equation*}
 Therefore, for any $y\in \m M_Y$ with $\|y-y_k\|\leq 2\tau_2$ and  $x\in \mb B_{\m M_X}(x_k,2\tau_2)$,
\begin{equation*}
\begin{aligned}
     \|y-G^{\dagger}_\sk((V_\sk^*)^T(y-y_0),x)\|\leq     \|y-G^*_\sk(Q^*_\sk(y),x)\|+C\,(\log n)\cdot n^{-\frac{1}{\frac{d_Y}{\beta_Y}+\frac{d_X}{\beta_X}}}=C\,(\log n)\cdot n^{-\frac{1}{\frac{d_Y}{\beta_Y}+\frac{d_X}{\beta_X}}}.
\end{aligned}
\end{equation*}
Therefore, by Lemma~\ref{lemma:reconstruct},  we can conclude that for any $\gamma_1\in (0,1]$, there exists a constant $C_{\gamma_1}$ so that it holds with probability at least $1-\frac{c}{n^3}$ that for any $k\in \wh{\m K}$, 
  \begin{equation*}
    \begin{aligned}
           &\mb{E}_{\mu^*_X}\mb{E}_{\mu^*_{Y|X}}[\|Y-\wh G_{[k]}(\wh Q_{[k]}(Y),X)\|^{\gamma_1}\cdot\mathbf{1}(X\in \mb B_{\mb R^{D_X}}(x_k,2\tau_2))\mathbf{1}(Y\in \mb B_{\mb R^{D_Y}}(y_k,2\tau_2))]\\
          & \leq C_{\gamma_1} \left\{
          \begin{array}{cc}
            \frac{(\log n)^{1+\gamma_1}}{\sqrt{n}}   & \frac{d_Y}{\beta_Y}+\frac{d_X}{\beta_X}\leq 2\gamma_1, \\
\Big((\log n \wedge \frac{1}{\beta_Y(d_Y/\beta_Y+d_X/\beta_X-2\gamma_1)} )^{1+\gamma_1}+(\log n)^{\gamma_1})\cdot n^{-\frac{\gamma_1}{\frac{d_X}{\beta_X}+\frac{d_Y}{\beta_Y}}}   & \frac{d_Y}{\beta_Y}+\frac{d_X}{\beta_X}>2\gamma_1.
          \end{array}
        \right.
        \end{aligned}
        \end{equation*}
Then if $\frac{d_Y}{2\beta_Y}+\frac{d_X}{2\beta_X}>1$, set $\gamma_1=1$,  it holds with probability at least $1-\frac{c}{n^3}$ that for any $k\in \wh{\m K}$, 
\begin{equation*}
    \begin{aligned}
       &\mb{E}_{\mu^*}[\|Y-\wh G_\sk(\wh Q_\sk(Y),X)\|\cdot\mathbf{1}(X\in \mb B_{\mb R^{D_X}}(x_0,2\tau_2))\mathbf{1}(Y\in \mb B_{\mb R^{D_Y}}(y_0,2\tau_2))]\lesssim  
 (\log n)\cdot n^{-\frac{1}{ \frac{d_Y}{\beta_Y}}}.
     \end{aligned}
\end{equation*}
Therefore, it holds with probability at least $1-\frac{c}{n^3}$ that for any $k\in \wh{\m K}$ and any $\gamma_1\in (0,1]$ that
\begin{equation*}
    \begin{aligned}
       &\mb{E}_{\mu^*}[\|Y-\wh G_\sk(\wh Q_\sk(Y),X)\|^{\gamma_1}\cdot\mathbf{1}(X\in \mb B_{\mb R^{D_X}}(x_0,2\tau_2))\mathbf{1}(Y\in \mb B_{\mb R^{D_Y}}(y_0,2\tau_2))]\\
       &\leq \Big(\mb{E}_{\mu^*}[\|Y-\wh G_\sk(\wh Q_\sk(Y),X)\|\cdot\mathbf{1}(X\in \mb B_{\mb R^{D_X}}(x_0,2\tau_2))\mathbf{1}(Y\in \mb B_{\mb R^{D_Y}}(y_0,2\tau_2))]\Big)^{\gamma_1}\lesssim (\log n)^{\gamma_1} \cdot n^{-\frac{\gamma_1}{ \frac{d_Y}{\beta_Y}+\frac{d_X}{\beta_X}}}.
     \end{aligned}
\end{equation*}
If $\frac{d_Y}{2\beta_Y}+\frac{d_X}{2\beta_X}<1$, let $\delta_n=\frac{1-\frac{d_Y}{4\beta_Y}-\frac{d_X}{4\beta_X}}{\lceil\log n\rceil}$ and  consider the set $\Gamma=\{\frac{d_Y}{4\beta_Y}+\frac{d_X}{4\beta_X}, \frac{d_Y}{4\beta_Y}+\frac{d_X}{4\beta_X}+\delta_n,\cdots, \frac{d_Y}{4\beta_Y}+\frac{d_X}{4\beta_X}+ \delta_n\cdot \lceil\log n\rceil\}$. Then by a union argument,  it holds that with probability at least $1-\frac{c\,\log n}{n^3}$ that for any $k\in \wh{\m K}$ and any $\gamma_1\in\Gamma$ that
  \begin{equation*}
    \begin{aligned}
           &\mb{E}_{\mu^*_X}\mb{E}_{\mu^*_{Y|X}}[\|Y-\wh G_{[k]}(\wh Q_{[k]}(Y),X)\|^{\gamma_1}\cdot\mathbf{1}(X\in \mb B_{\mb R^{D_X}}(x_k,2\tau_2))\mathbf{1}(Y\in \mb B_{\mb R^{D_Y}}(y_k,2\tau_2))]\\
          & \lesssim \left\{
          \begin{array}{cc}
            \frac{(\log n)^{1+\gamma_1}}{\sqrt{n}}   & \frac{d_Y}{\beta_Y}+\frac{d_X}{\beta_X}\leq 2\gamma_1, \\
\Big((\log n \wedge \frac{1}{\beta_Y(d_Y/\beta_Y+d_X/\beta_X-2\gamma_1)} )^{1+\gamma_1}+(\log n)^{\gamma_1})\cdot n^{-\frac{\gamma_1}{\frac{d_X}{\beta_X}+\frac{d_Y}{\beta_Y}}}   & \frac{d_Y}{\beta_Y}+\frac{d_X}{\beta_X}>2\gamma_1.
          \end{array}
        \right.
        \end{aligned}
        \end{equation*}
Then under the above event, for any $\gamma_2\in (0,1]$ with $\gamma_2<\frac{d_Y}{4\beta_Y}+\frac{d_X}{4\beta_X}$, by setting $\gamma_1=\frac{d_Y}{4\beta_Y}+\frac{d_X}{4\beta_X}$, it holds that 

\begin{equation*}
    \begin{aligned}
       &\mb{E}_{\mu^*}[\|Y-\wh G_\sk(\wh Q_\sk(Y),X)\|^{\gamma_2}\cdot\mathbf{1}(X\in \mb B_{\mb R^{D_X}}(x_0,2\tau_2))\mathbf{1}(Y\in \mb B_{\mb R^{D_Y}}(y_0,2\tau_2))]\\
       &\leq \Big(\mb{E}_{\mu^*}[\|Y-\wh G_\sk(\wh Q_\sk(Y),X)\|^{\gamma_1}\cdot\mathbf{1}(X\in \mb B_{\mb R^{D_X}}(x_0,2\tau_2))\mathbf{1}(Y\in \mb B_{\mb R^{D_Y}}(y_0,2\tau_2))]\Big)^{\gamma_2/\gamma_1}\\
       &\lesssim (\log n)^{\gamma_2}\cdot n^{-\frac{\gamma_2}{4\gamma_1}}=(\log n)^{\gamma_2}\cdot n^{-\frac{\gamma_2}{\frac{d_Y}{\beta_Y}+\frac{d_X}{\beta_X}}}.
     \end{aligned}
\end{equation*}
Moreover, for any $\gamma_2\in [\frac{d_Y}{4\beta_Y}+\frac{d_X}{4\beta_X},1]$, there exists $\gamma_1\in \Gamma$ so that $\gamma_1\leq\gamma_2\leq\gamma_1+\delta_n$, and therefore,
\begin{equation*}
    \begin{aligned}
       &\mb{E}_{\mu^*}[\|Y-\wh G_\sk(\wh Q_\sk(Y),X)\|^{\gamma_2}\cdot\mathbf{1}(X\in \mb B_{\mb R^{D_X}}(x_0,2\tau_2))\mathbf{1}(Y\in \mb B_{\mb R^{D_Y}}(y_0,2\tau_2))]\\
       &\leq 2L\, \mb{E}_{\mu^*}[\|Y-\wh G_\sk(\wh Q_\sk(Y),X)\|^{\gamma_1}\cdot\mathbf{1}(X\in \mb B_{\mb R^{D_X}}(x_0,2\tau_2))\mathbf{1}(Y\in \mb B_{\mb R^{D_Y}}(y_0,2\tau_2))] \\
       &\leq    \left\{
          \begin{array}{cc}
            C_1\,\frac{(\log n)^{1+\gamma_1}}{\sqrt{n}}   & \frac{d_Y}{\beta_Y}+\frac{d_X}{\beta_X}\leq 2\gamma_1, \\
C_1\, \big((\log n \wedge \frac{1}{\beta_Y(d_Y/\beta_Y+d_X/\beta_X-2\gamma_1)} )^{1+\gamma_1}+(\log n)^{\gamma_1}\big)\cdot n^{-\frac{\gamma_1}{\frac{d_X}{\beta_X}+\frac{d_Y}{\beta_Y}}}   & \frac{d_Y}{\beta_Y}+\frac{d_X}{\beta_X}>2\gamma_1.
          \end{array}
        \right.\\
        &\leq \ \left\{
          \begin{array}{cc}
            C_2\,\frac{(\log n)^{1+\gamma_2}}{\sqrt{n}}   & \frac{d_Y}{\beta_Y}+\frac{d_X}{\beta_X}\leq 2\gamma_2, \\
C_2\,\big((\log n \wedge \frac{1}{\beta_Y(d_Y/\beta_Y+d_X/\beta_X-2\gamma_2)} )^{1+\gamma_2}+(\log n)^{\gamma_2}\big)\cdot n^{-\frac{\gamma_2}{\frac{d_X}{\beta_X}+\frac{d_Y}{\beta_Y}}}   & \frac{d_Y}{\beta_Y}+\frac{d_X}{\beta_X}>2\gamma_2.
          \end{array}
        \right.
     \end{aligned}
\end{equation*}
This completes the proof for the first statement of Lemma~\ref{lemma:manifoldlearningr3} by combining all pieces. The second statement of Lemma~\ref{lemma:manifoldlearningr3} then directly follows from the second statement of Lemma~\ref{lemma:reconstruct}.

\subsection{Proof of Lemma~\ref{lestimation:2}}\label{proof:lestimation:2}

The proof follows the pipeline of the proof of Lemma~\ref{lestimation} and is included here for completeness.  To show the result for a fixed $j\in [J]$, we will use Theorem~\ref{theoremjointregression} with $\{\psi_{\lambda}(\cdot)\}_{\lambda\in \Lambda}=\Psi_j^{D_Y}$.  Then we will verify the three assumptions in Theorem~\ref{theoremjointregression}.   For the first assumption, it holds for a constant $C_1$ that 
 
    \begin{equation*}
\begin{aligned}
           & \sup_{(x,y)\in \m M}\sup_{S\in \m S_j^\dagger}\sum_{\psi\in  \Psi_j^{D_Y}}S^2(\psi,x)+|2^{\frac{j(d_Y-D_Y)}{2}}\psi(y) S(\psi,x)|\\
          & \leq \sup_{(x,y)\in \m M}\sup_{S\in \m S_j^\dagger}\sum_{\psi\in \Psi_j^{D_Y}} S^2(\psi,x)+ C\,\sum_{\psi\in \Psi_j^{D_Y}}|2^{\frac{j(d_Y-D_Y)}{2}}\psi (y)|\cdot 2^{-\frac{d_Yj}{2}}\\
  &\leq \sup_{x\in \m M_X} \Big\{\sup_{\psi\in \Psi_j^{D_Y}}|S(\psi,x)|^2\cdot\sum_{\psi\in \Psi_j^{D_Y}} \mathbf{1}\big(S(\psi,x)\neq 0\big)\Big\}+C\, \sum_{\psi\in \Psi_j^{D_Y}}|2^{\frac{-D_Yj}{2}}\psi (y)|\leq C_1.
\end{aligned}
\end{equation*}
 Then for the second assumption, we denote 
\begin{equation*}
    \begin{aligned}
       \ell(x,y,S)=\sum_{\psi\in \Psi_j^{D_Y}}S^2(\psi,x)-2^{\frac{j(d_Y-D_Y)}{2}+1}\psi(y) S(\psi,x).
    \end{aligned}
\end{equation*}
It holds for any $S,S'\in S_j^\dagger$ that
\begin{equation*}
\begin{aligned}
&\mb{E}_{\mu^*}\big[(\ell(X,Y,S)-\ell(X,Y,S'))^2\big]\\
&=\mb{E}_{\mu^*}\Big[ \Big(\sum_{\psi\in \Psi_j^{D_Y}}\big(S(\psi,X)+S'(\psi,X)-2^{\frac{j(d_Y-D_Y)}{2}+1}\psi(Y)\big)\cdot \big(S(\psi,X)-S'(\psi,X)\big)\Big)^2\Big]
\\
&\leq 8\, \mb{E}_{\mu^*}\Big[ \Big(\sum_{\psi\in  \Psi_j^{D_Y}}2^{\frac{j(d_Y-D_Y)}{2}+1}\psi(Y)\cdot \big(S(\psi,X)-S'(\psi,X)\big)\Big)^2\Big]\\
&\qquad\qquad+2\,\mb{E}_{\mu^*}\Big[ \Big(\sum_{\psi\in \Psi_j^{D_Y}}\big( S (\psi,X)+S'(\psi,X)\big)\cdot \big(S(\psi,X)-S'(\psi,X)\big)\Big)^2\Big]\\
&\leq  8\,\mb{E}_{\mu^*}\Big[ \Big(\sum_{\psi\in  \Psi_j^{D_Y}}2^{\frac{j(d_Y-D_Y)}{2}}\psi(Y)\cdot \big(  S (\psi,X)-S'(\psi,X)\big)\Big)^2\Big]+8C^2\mb{E}_{\mu^*_X}\Big[ \Big(\sum_{\psi\in \Psi_j^{D_Y}}2^{-\frac{d_Yj}{2}}\cdot \big|S(\psi,X)-S'(\psi,X)\big|\Big)^2\Big].
\end{aligned}
\end{equation*}
Then notice that 
\begin{equation*}
    \begin{aligned}
         &\mb{E}_{\mu^*}\Big[ \Big(\sum_{\psi\in \Psi_j^{D_Y}}2^{\frac{j(d_Y-D_Y)}{2}}\cdot\psi(Y)\cdot \big(S(\psi,X)-S'(\psi,X)\big)\Big)^2\Big]\\
         &=\mb{E}_{\mu^*_X}\Big[  \sum_{\psi_1,\psi_2\in \Psi_j^{D_Y}\atop {\rm supp}(\psi_1)\cap {\rm supp}(\psi_2)\neq \emptyset} 2^{j(d_Y-D_Y)}\mb{E}_{\mu^*_{Y|X}}[\psi_1(Y)\psi_2(Y)]\cdot \big(S(\psi_1,X)-S'(\psi_1,X)\big)\big(S(\psi_2,X)-S'(\psi_2,X)\big) \Big]\\
         &\lesssim\mb{E}_{\mu^*_X}\Big[  \sum_{\psi_1,\psi_2\in  \Psi_j^{D_Y}\atop {\rm supp}(\psi_1)\cap {\rm supp}(\psi_2)\neq \emptyset}\big(S(\psi_1,X)-S'(\psi_1,X)\big)^2+\big(S(\psi_2,X)-S'(\psi_2,X)\big)^2 \Big]\\
         &\lesssim \mb{E}_{\mu^*_X}\Big[  \sum_{\psi_1\in  \Psi_j^{D_Y}}\sum_{\psi_2\in  \Psi_j^{D_Y}\atop {\rm supp}(\psi_1)\cap {\rm supp}(\psi_2)\neq \emptyset}\big(S(\psi_1,X)-S'(\psi_1,X)\big)^2\Big]\\
         &\lesssim \mb{E}_{\mu^*_X}\Big[  \sum_{\psi_1\in  \Psi_j^{D_Y}}\big(  S (\psi_1,X)-S'(\psi_1,X)\big)^2  \Big],
    \end{aligned}
\end{equation*}
where we have used the fact that for any $x\in \m M_X$,
\begin{equation*}
    \mb{E}_{\mu^*_{Y|x}}[\psi_1(Y)\psi_2(Y)]\lesssim\int_{\m M_{Y|x}} \mathbf{1}(y\in {\rm supp}(\psi_1)\cap {\rm supp}(\psi_2)) 2^{D_Yj}\,u^*(y\,|\,x)\,\dd {\rm vol}_{\m M_{Y|x}}\lesssim 2^{(D_Y-d_Y)j}.    
\end{equation*}
Moreover,
\begin{equation*}
\begin{aligned}
   & \mb{E}_{\mu^*_X}\Big[ \Big(\sum_{\psi\in \Psi_j^{D_Y}}2^{-\frac{d_Yj}{2}}\cdot \big|S(\psi,X)-S'(\psi,X)\big|\Big)^2\Big]\\
   &=\mb{E}_{\mu^*_X}\Big[ \Big(\sum_{\psi\in \Psi_j^{D_Y}}2^{-\frac{d_Yj}{2}}\cdot \big|S(\psi,X)-S'(\psi,X)\big|\cdot \mathbf{1}\big(S(\psi,X)\neq 0\text{ or }S'(\psi,X)\neq 0\big)\Big)^2\Big]\\
   &\leq  \mb{E}_{\mu^*_X}\Big[ \sum_{\psi\in   \Psi_j^{D_Y}}2^{-d_Yj} \mathbf{1}\big(S(\psi,X)\neq 0\text{ or }S'(\psi,X)\neq 0\big)\cdot \sum_{\psi\in \wh\Psi_j^{D_Y}} \big(S(\psi,X)-S'(\psi,X)\big)^2\Big]\\
   &\leq 2C\, \mb{E}_{\mu^*_X}\Big[ \sum_{\psi\in  \Psi_j^{D_Y}} \big(S(\psi,x)-S'(\psi,x)\big)^2\Big].
\end{aligned}
\end{equation*}
 Therefore, it holds for some constant $C_2$ that
\begin{equation*}
\begin{aligned}
     &\left\|\ell(x,y,S)-\ell(x,y,S')\right\|^2_2\leq C_2\,   \mb{E}_{\mu^*_X}\Big[ \sum_{\psi\in \Psi_j^{D_Y}} \big(S(\psi,x)-S'(\psi,x)\big)^2\Big],\\
\end{aligned}
\end{equation*}
which verifies the second assumption.
For the last assumption, note that for any $S,S'\in \m S_j^\dagger$,
\begin{equation*}
 \begin{aligned}
 &d_n( S , S')\\
&=\sqrt{\frac{1}{n}\sum_{i=1}^n \Big(\sum_{\psi\in  \Psi_j^{D_Y}}\big(S^2(\psi,X_i)-S'{}^2(\psi,X_i)\big)-2^{\frac{j(d_Y-D_Y)}{2}+1}\psi(Y_i) \big(S(\psi,X_i)-S'(\psi,X_i)\big)\Big)^2}\\
&\leq  \sqrt{\frac{1}{n}\sum_{i=1}^n  \sum_{\psi\in \Psi_j^{D_Y}}\big(S(\psi,X_i)+S'(\psi,X_i)-2^{\frac{j(d_Y-D_Y)}{2}+1}\psi(Y_i)\big)^2\cdot  \sum_{\psi\in \Psi_j^{D_Y}}\big(S(\psi,X_i)-S'(\psi,X_i)\big)^2}\\
&\leq C_3 2^{\frac{d_Y j}{2}}\sqrt{\frac{1}{n}\sum_{i=1}^n     \sum_{\psi\in \Psi_j^{D_Y}}\big(S(\psi,X_i)-S'(\psi,X_i)\big)^2}\\
&\leq C_3 2^{\frac{d_Y j}{2}}d_S(S,S').
 \end{aligned}
 \end{equation*}
Then, using the fact that $\log \mathbf{N}({\m S}_j^\dagger,d^S,\varepsilon)\leq \m W_j\log(\frac{ n}{\varepsilon})$, there exists a constant we have for any $0<\varepsilon\leq \sup_{S,S'\in S_j^\dagger}d_n(S,S')$,
 \begin{equation*}
    \log \mathbf{N}(\m S_j^\dagger,d_n,\varepsilon)\leq \m W_j\log \frac{C_3\,n \cdot 2^{d_YJ/2}}{\varepsilon}\leq 2\,\m W_j\log \frac{n}{\varepsilon}.
 \end{equation*}
 The desired result is obtained by setting $W_n=2\m W_j$ and $T_n=n$ in Theorem~\ref{theoremjointregression},  and applying a union bound over $j\in [J]$.

\subsection{Proof of Lemma~\ref{lemma:reconstruct}}\label{proof:lemma:reconstruct}
 

Denote  $\wh Q(y)=\wh V^T(y-y_0)$, it holds that 
\begin{equation*}
\begin{aligned}
       &\frac{1}{n}\sum_{i=1}^n \|Y_i-\wh G(\wh Q(Y_i),X_i)\|^2\mathbf{1}(X_i\in \mb B_{\mb R^{D_X}}(x_0,2\tau_2))\mathbf{1}(Y_i\in \mb B_{\mb R^{D_Y}}(y_0,2\tau_2))\\
    &= \underset{G\in \m G\atop V\in \mb O(D_Y,d_Y)}{\min}\frac{1}{n}\sum_{i=1}^n \|Y_i- G( V^T(Y_i-y_0),X_i)\|^2\mathbf{1}(X_i\in \mb B_{\mb R^{D_X}}(x_0,2\tau_2))\mathbf{1}(Y_i\in \mb B_{\mb R^{D_Y}}(y_0,2\tau_2))\\
    &\leq (\varepsilon^*)^2.
\end{aligned}
\end{equation*}
Therefore, 
\begin{equation*}
\begin{aligned}
       &\frac{1}{n}\sum_{i=1}^n \|Y_i-\wh G(\wh Q(Y_i),X_i)\|^{\gamma_1}\mathbf{1}(X_i\in \mb B_{\mb R^{D_X}}(x_0,2\tau_2))\mathbf{1}(Y_i\in \mb B_{\mb R^{D_Y}}(y_0,2\tau_2))\\
    &\leq \Big(\frac{1}{n}\sum_{i=1}^n \|Y_i- \wh G(\wh Q(Y_i),X_i)\|^2\mathbf{1}(X_i\in \mb B_{\mb R^{D_X}}(x_0,2\tau_2))\mathbf{1}(Y_i\in \mb B_{\mb R^{D_Y}}(y_0,2\tau_2))\Big)^{\frac{\gamma_1}{2}}\leq (\varepsilon^*)^{\gamma_1}.
\end{aligned}
\end{equation*}
Define the class 
    \begin{equation*}
 \begin{aligned}
 &{\mathcal{F}}=\{f(x,y)=\|y-G(V^T(y-y_0),x)\|^{\gamma_1}:\, G\in \m G, V\in \mb O(D_Y,d_Y)\}.
     \end{aligned}
\end{equation*} 
Then  we have $\|y-\wh G(\wh Q(y),x)\|^{\gamma_1}\in \m F$.  Moreover, It is straightforward to verify that for any $\beta\in (1,\beta_Y)$, there exists a constant $L$  so that for any $G\in \m G$, and $x\in \m M_X$, it holds that $G(\cdot,x)\in \m H^{\beta}_{L,D_Y}(\mb R^{d_Y})$. Then for any $G_1,G_2\in \m G$, $V_1,V_2\in \mb O(D_Y,d_Y)$, and $(x,y)\in \m M$ where  $x\in \mb B_{\mb R^{D_X}}(x_k,2\tau_2)$ and $y\in \mb B_{\mb R^{D_Y}}(y_k,2\tau_2)$, it holds that
\begin{equation}\label{eqn13}
 \begin{aligned}
 &\Big| \|y-G_1(V_1^T(y-y_k),x)\|^{\gamma_1}  -\|y-G_2(V_2^T(y-y_k),x)\|^{\gamma_1} \Big|\\
 &\leq\|G_1(V_1^T(y-y_k),x)-G_2(V_2^T(y-y_k),x)\|^{\gamma_1} \\
 &\leq  \underset{z\in \mb R^{d_Y},x\in \m M_X}{\sup} \|G_1(z,x)-G_2(z,x)\|^{\gamma_1} + (2L\tau_2)^{\gamma_1}\,\|V_1-V_2\|_{\rm op}^{\gamma_1}
     \end{aligned}
\end{equation} 
Consider the distance 
\begin{equation*}
    \begin{aligned}
        d_n(f,f')=\sqrt{\frac{1}{n}\sum_{i=1}^n (f(X_i,Y_i)-f'(X_i,Y_i))^2\mathbf{1}(X_i\in \mb B_{\mb R^{D_X}}(x_0,2\tau_2))\mathbf{1}(Y_i\in \mb B_{\mb R^{D_Y}}(y_0,2\tau_2))}.
    \end{aligned}
\end{equation*}
Using~\eqref{eqn13}, we can bound the $\varepsilon$-covering number $\mathbf{N}(\m F,d_n,\varepsilon)$ of $\m F$ with respect to $d_n$ by 
\begin{equation*}
\begin{aligned}
        \mathbf{N}(\m F,d_n,\varepsilon)&\leq  \mathbf{N}(\m G,d_{\infty}^{\gamma_1},\frac{\varepsilon}{2})\cdot  \mathbf{N}\Big(\mb O(D_Y,d_Y),\|\cdot\|_{\rm op},\frac{\varepsilon^{\frac{1}{\gamma_1}}}{2^{\frac{1}{\gamma_1}+1}L\tau_2}\Big)\\
    &\leq \mathbf{N}(\m G,d_{\infty}^{\gamma_1},\frac{\varepsilon}{2})  \cdot (\frac{C}{\varepsilon^{\frac{1}{\gamma_1}}})^{D_Yd_Y}.
\end{aligned}
\end{equation*}
Then by standard symmetrization and Dudley’s entropy integral bound (see for example,~\cite{wainwright2019high}), we can get that 
\begin{equation*}
 \begin{aligned}
 &\mathbb{E}\bigg[\underset{f\in \mathcal{F}} {\sup}\Big|\frac{1}{n}\sum_{i=1}^n f(X_i,Y_i)\mathbf{1}(X_i\in \mb B_{\mb R^{D_X}}(x_0,2\tau_2))\mathbf{1}(Y_i\in \mb B_{\mb R^{D_Y}}(y_0,2\tau_2))\\
 &-\mathbb{E}_{\mu^*} \big[f(X,Y)\mathbf{1}(X\in \mb B_{\mb R^{D_X}}(x_0,2\tau_2))\mathbf{1}(Y\in \mb B_{\mb R^{D_Y}}(y_0,2\tau_2))\big]\Big|\bigg]\\
 &\leq \frac{C_1}{\sqrt{n}}\int_0^{\infty} \sqrt{ \log \mathbf{N}(\m G,d^{\gamma_1}_{\infty},\frac{\varepsilon}{2} )}\dd \varepsilon+\frac{C_1}{\gamma_1\sqrt{n}}.
      \end{aligned}
\end{equation*} 
Then by Talagrand concentration inequality (see for example, Theorem 3.27 of~\cite{wainwright2019high}), there exists a constant $C_2$, such that it holds with probability at least $1-n^{-3}$ that 
\begin{equation*}
    \begin{aligned}
        \underset{f\in \mathcal{F}} {\sup}&\Big|\frac{1}{n}\sum_{i=1}^n f(X_i,Y_i)\mathbf{1}(X_i\in \mb B_{\mb R^{D_X}}(x_0,2\tau_2))\mathbf{1}(Y_i\in \mb B_{\mb R^{D_Y}}(y_0,2\tau_2))\\
 &-\mathbb{E}_{\mu^*} \big[f(X,Y)\mathbf{1}(X\in \mb B_{\mb R^{D_X}}(x_0,2\tau_2))\mathbf{1}(Y\in \mb B_{\mb R^{D_Y}}(y_0,2\tau_2))\big]\Big|\\
 &\leq C_2 \frac{1}{\sqrt{n}}\int_0^{\infty} \sqrt{ \log \mathbf{N}(\m G,d^{\gamma_1}_{\infty},\frac{\varepsilon}{2} )}\dd \varepsilon+C_2\sqrt{\frac{\log n}{n}}+\frac{C_2}{\gamma_1\sqrt{n}}.
    \end{aligned}
\end{equation*}
So by combining all pieces, it holds with probability at least $1-n^{-3}$ that 
\begin{equation*}
 \begin{aligned}
 &\mb{E}_{\mu^*}[\|Y-\wh G(\wh V(Y-y_0),X)\|^{\gamma_1} \cdot\mathbf{1}(X\in \mb B_{\mb R^{D_X}}(x_0,2\tau_2))\mathbf{1}(Y\in \mb B_{\mb R^{D_Y}}(y_0,2\tau_2))]\\
&=\mb{E}_{\mu^*}[\|Y-\wh G(\wh Q(Y),X)\|^{\gamma_1} \cdot\mathbf{1}(X\in \mb B_{\mb R^{D_X}}(x_0,2\tau_2))\mathbf{1}(Y\in \mb B_{\mb R^{D_Y}}(y_0,2\tau_2))]\\
 &\leq   \underset{f\in \mathcal{F}} {\sup}\Big|\frac{1}{n}\sum_{i=1}^n f(X_i,Y_i)\mathbf{1}(X_i\in \mb B_{\mb R^{D_X}}(x_0,2\tau_2))\mathbf{1}(Y_i\in \mb B_{\mb R^{D_Y}}(y_0,2\tau_2))\\
 &-\mathbb{E}_{\mu^*} \big[f(X,Y)\mathbf{1}(X\in \mb B_{\mb R^{D_X}}(x_0,2\tau_2))\mathbf{1}(Y\in \mb B_{\mb R^{D_Y}}(y_0,2\tau_2))\big]\Big|+(\varepsilon^*)^{\gamma_1}\\
 &\lesssim \frac{1}{\sqrt{n}}\int_0^{\infty} \sqrt{ \log \mathbf{N}(\m G,d_{\infty}^{\gamma_1}, \varepsilon )}\dd \varepsilon+\sqrt{\frac{\log n}{n}}+\frac{1}{\gamma_1\sqrt{n}}+(\varepsilon^*)^{\gamma_1}.
      \end{aligned}
\end{equation*} 
The proof of the first statement is complete. Then we show the second statement.    Let $V^*$ be a $D_Y\times d_Y$ matrix whose column form an orthonormal basis of $T_{\m M_{Y|x^*}}y^*$. Denote $Q^*(y)=(V^*)^T(y-y^*)$ and $G^*(z,x)=\Phi_{(x^*,y^*)}(V^*z,x)$. Then $G^*\in   {\m H}^{\beta_Y,\beta_X}_{L,D_Y}(\mb B_{\mb R^{d_Y}}(\mathbf{0},\tau_1),\mb B_{\m M_X}(x^*,\tau))$, and for any $(x,y)\in \m M$ with $\|x-x^*\|<\tau$ and $\|y-y^*\|<\tau_1$,  we have $y=G^*(Q^*(y),x)$.  Moreover, define 
 \begin{equation*}
     v^*(z,x)=u^*(G^*(z,x)|x)\cdot\sqrt{{\rm det}(J_{G^*(\cdot,x)}(z)^TJ_{G^*(\cdot,x)}(z))}.
 \end{equation*}
Let $\alpha_1=1\wedge \alpha_Y$ and  $\alpha_2=1\wedge \alpha_X\wedge \alpha_Y\wedge (\alpha_Y\beta_X)\wedge (\beta_X-\frac{\beta_X}{\beta_Y})$.  It holds that $v^*\in  \ov{\m H}^{\alpha_1,\alpha_2}_L(\mb B_{\mb R^{d_Y}}(\mathbf{0},\tau_1), \mb B_{\m M_X}(x^*,\tau))$ with a constant $L$. Therefore, there exists a constant $L_1$ so that for any $x',x\in \mb B_{\m M_X}(x^*,\tau)$ and $z,z'\in \mb B_{\mb R^{d_Y}}(\mathbf{0},\tau_1)$,
\begin{equation*}
  \|  v^*(z,x)-  v^*(z',x')\|\leq L_1(\|z-z'\|^{\alpha_1}+\|x-x'\|^{\alpha_2}).
\end{equation*}
Moreover, there exists a constant $\tau_3<\tau_2$ so that when $\|x-x^*\|\leq \tau_3$ and $\|z\|\leq \tau_3$, 
\begin{equation*}
    \|G^*(z,x)-y_0\|\leq \|G^*(z,x)-G^*(\mathbf{0},x^*)\|+\|y^*-y_0\|<2\tau_2.
\end{equation*}
Furthermore, since $\mu^*_{Y|x^*}(\mb B_{\m M_{Y|x^*}}(y^*,\tau_3/2))\geq g(\tau_3/2)/L$, it holds that
\begin{equation*}
    \begin{aligned}
     &  g(\tau_3/2)/L \leq \mu^*_{Y|x^*}(\mb B_{\m M_{Y|x^*}}(y^*,\tau_3/2))\\
       &=\int_{\{z\in \mb B_{\mb R^{d_Y}}(\mathbf{0},\tau_1):\, \|G^*(z,x^*)-y^*\|<\tau_3/2\}}v^*(z,x^*)\,\dd z\\
       &\leq \int_{\mb B_{\mb R^{d_Y}}(\mathbf{0},\tau_3/2)}v^*(z,x^*)\,\dd z\\
       &\leq \underset{z\in \mb B_{\mb R^{d_Y}}(\mathbf{0},\tau_3/2) }{\max} v^*(z,x^*) \frac{\pi^{d_Y / 2}}{(d_Y/2)!} (\tau_3/2)^{d_Y}.
    \end{aligned}
\end{equation*}
Therefore, there exists $\wt z\in \mb B_{\mb R^{d_Y}}(\mathbf{0},\tau_3/2)$ so that $v^*(\wt z,x^*)\geq \frac{g(\tau_3/2)(d_Y/2)!}{\pi^{d_Y / 2}(\tau_3/2)^{d_Y}L}=\tau_4>0$. Then consider a small enough positive constant $\tau_5$ that will be chosen later.   When $\tau_5<\frac{\tau_3}{2}\wedge (\frac{\tau_4}{4L_1})^{\frac{1}{\alpha_1\wedge \alpha_2}}$,  for any $z\in \mb B_{\mb R^{d_Y}}(\wt z, \tau_5)$ and $x\in \mb B_{\m M_X}(x^*,\tau_5)$, it holds that
\begin{equation*}
  v^*(z,x)\leq L_1\|z-\wt z\|^{\alpha_1}+L_1\|x-x^*\|^{\alpha_2}+  v^*(\wt z,x^*)\leq \frac{\tau_4}{2}+ v^*(\wt z,x^*)\leq \frac{3}{2} v^*(\wt z,x^*)
\end{equation*}
and 
\begin{equation*}
    v^*(z,x)\geq v^*(\wt z,x^*)-\frac{\tau_4}{2}\geq \frac{v^*(\wt z,x^*)}{2}.
\end{equation*}
Moreover, since  $\mb{E}_{\mu^*}[\|Y-\wh G(\wh Q(Y),X)\|\cdot\mathbf{1}(X\in \mb B_{\mb R^{D_X}}(x_0,2\tau_2))\mathbf{1}(Y\in \mb B_{\mb R^{D_Y}}(y_0,2\tau_2))]\leq c$, there exists a constant $C_1$ so that
\begin{equation*}
    \begin{aligned}
             &\mb{E}_{\mu^*_X}\Big[\int_{\mb B_{\mb R^{d_Y}}(\wt z,\tau_5)}\|G^*(z,X)-\wh G(\wh Q(G^*(z,X)),X)\|^2\cdot\mathbf{1}(X\in \mb B_{\m M_X}(x^*,\tau_5)) v^*(z,X)\,\dd z\Big]\\
         &\leq \mb{E}_{\mu^*_X}\mb{E}_{\mu^*_{Y|X}}[\|Y-\wh G(\wh Q(Y),X)\|^2\cdot\mathbf{1}(X\in \mb B_{\mb R^{D_X}}(x_0,2\tau_2))\mathbf{1}(Y\in \mb B_{\mb R^{D_Y}}(y_0,2\tau_2))]\\
         &\leq C_1\, \mb{E}_{\mu^*_X}\mb{E}_{\mu^*_{Y|X}}[\|Y-\wh G(\wh Q(Y),X)\|\cdot\mathbf{1}(X\in \mb B_{\mb R^{D_X}}(x_0,2\tau_2))\mathbf{1}(Y\in \mb B_{\mb R^{D_Y}}(y_0,2\tau_2))]\\
         &\leq C_1c.
           \end{aligned}
\end{equation*}
Define $\wh l(z,x)=\wh G(\wh Q(G^*(z,x)),x)$. Given that for any $x\in \m M_X$, $\wh G(\cdot,x)\in\m H^{\beta}_{L,D_Y}(\mb R^{d_Y})$ with $\beta>1$ and a constant $L$, there exists a constant $L_2$ such that for any $x\in B_{\m M_X}(x,\tau)$ and $z\in\mb B_{\mb R^{d_Y}}(\mathbf{0},\tau_1)$, 
\begin{equation*}
\begin{aligned}
        &\|(G^*(z,x)-\wh l(z,x))-(G^*(\wt z,x)-\wh l(\wt z,x)+\big(J_{G^*(\cdot,x)}(\wt z)-J_{\wh l(\cdot,x)}(\wt z)\big)(z-\wt z))\|\\
        &\qquad\leq L_2\|z-\wt z\|^{\beta\wedge 2}.
\end{aligned}
\end{equation*}
Therefore,
\begin{equation*}
    \begin{aligned}
          &\mb{E}_{\mu^*_X}\Big[\int_{\mb B_{\mb R^{d_Y}}(\wt z,\tau_5)}\|G^*(z,X)-\wh G(\wh Q(G^*(z,X)),X)\|^2\cdot v^*(z,X)\,\dd z\cdot \mathbf{1}(X\in \mb B_{\m M_X}(x^*,\tau_5))\Big]\\
          &=\mb{E}_{\mu^*_X|_{B_{\m M_X}(x^*,\tau_5)}}\Big[\int_{\mb B_{\mb R^{d_Y}}(\wt z,\tau_5)}\|G^*(z,X)-\wh l(z,X)\|^2\cdot v^*(z,X)\,\dd z\Big]\\
              &\geq \frac{1}{4}\mb{E}_{\mu^*_X|_{B_{\m M_X}(x^*,\tau_5)}}\Big[\int_{\mb B_{\mb R^{d_Y}}(\wt z,\tau_5)}\|G^*(\wt z,X)-\wh l(\wt z,X)+\big(J_{G^*(\cdot,X)}(\wt z)-J_{\wh l(\cdot,X)}(\wt z)\big)(z-\wt z)\|^2\,\dd z\cdot  v^*(\wt z,x^*)\Big]\\
              &\qquad\qquad-\frac{3\pi^{d_Y / 2}}{2(d_Y/2)!}L_2\,\tau_5^{2(\beta\wedge 2)}   v^*(\wt z,x^*) (\tau_5)^{d_Y}\mu^*_X(B_{\m M_X}(x,\tau_5)) \\
              &\geq L_3 \, \mb{E}_{\mu^*_X|_{B_{\m M_X}(x^*,\tau_5)}}\Big[\|J_{G^*(\cdot,X)}(\wt z)-J_{\wh l(\cdot,X)}(\wt z)\|_{\rm F}^2\Big]\tau_5^{d_Y+2} v^*(\wt z,x^*)-\\
              &\qquad\qquad\frac{3\pi^{d_Y / 2}}{2(d_Y/2)!}L_2\,\tau_5^{2(\beta\wedge 2)}  v^*(\wt z,x^*) (\tau_5)^{d_Y}\mu^*_X(B_{\m M_X}(x,\tau_5)),
    \end{aligned}
\end{equation*}
where the last inequality uses the fact that  for any $d$-variate polynomial $\mathbf{S}(y)=\sum_{j\in \mb N_0^d,\,|j|\leq k} a_{j} y^{j}$, $y\in\mb R^d$, there exists some positive constant $C(d,k)$ only depending on $(d,k)$ such that 
\begin{equation*}
\int_{\mb B_1^d} \mathbf{S}^2(y) \, \dd y \geq C(d,k) \sum_{j\in \mb N_0^d,\,|j|\leq k} a_{j}^2.
  \end{equation*}
So combined with $\mb{E}_{\mu^*_X}\Big[\int_{\mb B_{\mb R^{d_Y}}(\wt z,\tau_5)}\|G^*(z,X)-\wh G(\wh Q(G^*(z,X)),X)\|^2\cdot\mathbf{1}(X\in \mb B_{\m M_X}(x^*,\tau_5)) v^*(z,X)\,\dd z\Big]\leq C_1c$, we can obtain
\begin{equation*}
    \begin{aligned}
        & \mb{E}_{\mu^*_X|_{B_{\m M_X}(x^*,\tau_5)}}\Big[\|J_{G^*(\cdot,X)}(\wt z)-J_{\wh l(\cdot,X)}(\wt z)\|_{\rm F}^2\Big]\leq \frac{C_1c }{\tau_5^{d_Y+2}\tau_4L_3}\\
         &\qquad +\frac{3\pi^{d_Y / 2}L_2}{2(d_Y/2)!L_3}\,\tau_5^{2(\beta\wedge 2)-2}   \mu^*_X(B_{\m M_X}(x^*,\tau_5)).
    \end{aligned}
\end{equation*}
Therefore there exists $\wt x\in B_{\m M_X}(x^*,\tau_5)$,  so that 
\begin{equation*}
    \begin{aligned}
        \|J_{G^*(\cdot,\wt x)}(\wt z)-J_{\wh l(\cdot,\wt x)}(\wt z)\|_{\rm F}^2\leq\frac{C_1c}{\tau_5^{d_Y+2}\tau_4L_3  \mu^*_X(B_{\m M_X}(x^*,\tau_5))}+\frac{3\pi^{d_Y / 2}L_2}{2(d_Y/2)!L_3}\,\tau_5^{2(\beta\wedge 2)-2}.
    \end{aligned}
\end{equation*}
Then notice that 
\begin{equation*}
    J_{\wh l(\cdot,\wt x)}(\wt z)= J_{\wh G(\cdot,\wt x)}(\wh Q(G^*(\wt z,\wt x)))\wh V^T J_{G^*(\cdot,\wt x)}(\wt z), 
\end{equation*}
 \begin{equation*}
   \| J_{G^*(\cdot,\wt x)}(\wt z)-V^*\|= \| J_{G^*(\cdot,\wt x)}(\wt z)-J_{G^*(\cdot, x^*)}(\mathbf 0)\|\leq L\,(\|\wt z\|+\|\wt x-x^*\|^{(\beta_X-\frac{\beta_X}{\beta_Y})\wedge 1}),
 \end{equation*}
and  there exists a constant $L_4$ so that for any $z\in \mb R^{d_Y}$ and $x\in \mb R^{D_X}$,
 \begin{equation*}
     J_{\wh G(\cdot,x)}(z)^T J_{\wh G(\cdot,x)}(z)\preceq L_4 I_{d_Y}.
 \end{equation*}
When $\tau_3$, $\tau_5$ and $c$ are small enough, it holds that
\begin{equation*}
\begin{aligned}
       &\|V^*-J_{\wh G(\cdot,\wt x)}(\wh Q(G^*(\wt z,\wt x)))\wh V^T V^*\|_{\rm F}\\
       &\leq \|V^*-J_{G^*(\cdot,\wt x)}(\wt z)\|_{\rm F}+\|J_{G^*(\cdot,\wt x)}(\wt z)-J_{\wh l(\cdot,\wt x)}(\wt z)\|_{\rm F}\\
       &\qquad+\|J_{\wh G(\cdot,\wt x)}(\wh Q(G^*(\wt z,\wt x)))\wh V^T J_{G^*(\cdot,\wt x)}(\wt z)-J_{\wh G(\cdot,\wt x)}(\wh Q(G^*(\wt z,\wt x)))\wh V^T V^*\|_{\rm F}\\
       &\leq \|V^*-J_{G^*(\cdot,\wt x)}(\wt z)\|_{\rm F}+\|J_{G^*(\cdot,\wt x)}(\wt z)-J_{\wh l(\cdot,\wt x)}(\wt z)\|_{\rm F}\\
       &\qquad+\|J_{\wh G(\cdot,\wt x)}(\wh Q(G^*(\wt z,\wt x)))\wh V^T\|_{\rm op} \|J_{G^*(\cdot,\wt x)}(\wt z)-V^*\|_{\rm F}\\
       &\leq (1+\sqrt{L_4})L(\tau_3/2+\tau_5^{(\beta_X-\frac{\beta_X}{\beta_Y})\wedge 1})+\sqrt{\frac{C_1c}{\tau_5^{d_Y+2}\tau_4L_3  \mu^*_X(B_{\m M_X}(x^*,\tau_5))}+\frac{3\pi^{d_Y / 2}L_2}{2(d_Y/2)!L_3}\,\tau_5^{2(\beta\wedge 2)-2}} \\
       &\leq \frac{1}{4}(1\wedge \frac{1}{L_4}).
\end{aligned}
\end{equation*}
 Therefore, 
\begin{equation*}
  \begin{aligned}
      &\|I_{d_Y}-(V^*)^T\wh VJ_{\wh G(\cdot,\wt x)}(\wh Q(G^*(\wt z,\wt x)))^TJ_{\wh G(\cdot,\wt x)}(\wh Q(G^*(\wt z,\wt x)))\wh V^T V^*\|_{\rm F}\\
      &=   \|(V^*)^TV^*-(V^*)^T\wh VJ_{\wh G(\cdot,\wt x)}(\wh Q(G^*(\wt z,\wt x)))^TJ_{\wh G(\cdot,\wt x)}(\wh Q(G^*(\wt z,\wt x)))\wh V^T V^*\|_{\rm F}\\
      &\leq \|(V^*)^T(V^*-J_{\wh G(\cdot,\wt x)}(\wh Q(G^*(\wt z,\wt x)))\wh V^T V^*)\|_{\rm F}\\
      &\qquad\qquad+\|(V^*-J_{\wh G(\cdot,\wt x)}(\wh Q(G^*(\wt z,\wt x)))\wh V^T V^*)^TJ_{\wh G(\cdot,\wt x)}(\wh Q(G^*(\wt z,\wt x)))\wh V^T V^*\|_{\rm F}\\
      &\leq \frac{1}{2},
  \end{aligned}
\end{equation*}
which, combined with  $ J_{\wh G(\cdot,\wt x)}(\wh Q(G^*(\wt z,\wt x)))^T J_{\wh G(\cdot,\wt x)}(\wh Q(G^*(\wt z,\wt x)))\preceq L_4 I_{d_Y}$ can imply that 
\begin{equation*}
    (V^*)^T\wh V\wh V^TV^*\succeq \frac{1}{2L_4} I_d,
\end{equation*}
and thus

\begin{equation*}
     \wh V^TP^*\wh V=\wh V^TV^*(V^*)^T\wh V\succeq \frac{1}{2L_4} I_d.
\end{equation*}

\subsection{Proof of Lemma~\ref{lemma:coveringGregime2}}\label{proof:lemma:coveringGregime2}
 
    Consider $$G(z)=\sum_{j_1=0}^{J_1}\sum_{\psi_1\in  \Psi_{j_1}^{d_Y} }g_{\psi_1}\psi_1(z)$$ and $$G'(z)=\sum_{j_1=0}^{J_1}\sum_{\psi_1\in  \Psi_{j_1}^{d_Y} }g'_{\psi_1}\psi_1(z).$$ 
     Then there exists a constant $C$ so that 
    \begin{equation*}
        \begin{aligned}
           & \underset{ z\in \mb R^{d_Y}}{\sup}\|G(z)-G'(z)\|\\
            &=\underset{z\in \mb R^{d_Y}}{\sup}\Big\|\sum_{j_1=0}^{J_1}\sum_{\psi_1\in  \Psi_{j_1}^{d_Y} }\big(g_{\psi_1}-g_{\psi_1}'\big)\psi_1(z)\Big\|\\
            &\leq \sum_{j_1=0}^{J_1}\underset{\psi_1\in  \Psi_{j_1}^{d_Y} }{\max}\big\|g_{\psi_1}-g_{\psi_1}'\big\| \cdot\underset{ z\in \mb R^{d_Y}}{\sup} \sum_{\psi_1\in  \Psi_{j_1}^{d_Y} }|\psi_1(z)|\\
            &\leq C\, \sum_{j_1=0}^{J_1}\underset{\psi_1\in  \Psi_{j_1}^{d_Y} }{\max}\big\|g_{\psi_1}-g_{\psi_1'}\big\| \cdot 2^{\frac{d_Yj_1}{2}}.
        \end{aligned}
    \end{equation*}
When $\frac{d_Y}{\beta_Y} \leq 2\gamma_1$, we have 
    \begin{equation*}
        \begin{aligned}
            &\sum_{j_1=0}^{J_1} 2^{\frac{d_Yj_1 }{4\gamma_1}-\frac{j_1\beta_Y}{2}} \leq  (J_1+1)\leq 2J_1.
                   \end{aligned}
    \end{equation*}
So if for any $j_1\in [J_1]$ and $\psi_1\in  \Psi_{j_1}^{d_Y}$,
    \begin{equation*}
 \big\|g_{\psi_1 }-g_{\psi_1 }'\big\|\leq \frac{\varepsilon^{\frac{1}{\gamma_1}}}{2CJ_1} 2^{\frac{d_Yj_1}{4\gamma_1}-\frac{j_1\beta_Y}{2}-\frac{d_Yj_1}{2}},
    \end{equation*}
    then 
    \begin{equation*}
        \underset{z\in \mb R^{d_Y}}{\sup}\|G(z)-G'(z)\|^{\gamma_1}\leq \Big(C\sum_{j_1=0}^{J_1}\frac{\varepsilon^{\frac{1}{\gamma_1}}}{2CJ_1} 2^{\frac{d_Yj_1}{4\gamma_1}-\frac{j_1\beta_Y}{2}}\Big)^{\gamma_1}\leq \varepsilon.
    \end{equation*}
    Therefore, we can get
    \begin{equation*}
    \begin{aligned}
      \mathbf N(\m G,d_{\infty}^{\gamma_1},\varepsilon)&\leq \prod_{j_1=0}^{J_1}\prod_{\psi_1\in  \Psi_{j_1}^{d_Y}} \mathbf N([-L_1\, 2^{-\frac{d_Yj_1}{2}-j_1\beta_Y},
        L_1\, 2^{-\frac{d_Yj_1}{2}-j_1\beta_Y}]^{D_Y},\frac{\varepsilon^{\frac{1}{\gamma_1}}}{2CJ_1} 2^{\frac{d_Yj_1}{4\gamma_1}-\frac{j_1\beta_Y}{2}-\frac{d_Yj_1}{2}},\|\cdot\|)\\
          &\leq \prod_{j_1=0}^{J_1}\prod_{\psi_1\in  \Psi_{j_1}^{d_Y}}\lceil\Big(\frac{12\sqrt{D_Y}L_1C \, J_12^{-\frac{d_Yj_1}{4\gamma_1}-\frac{j_1\beta_Y}{2}}}{\varepsilon^{\frac{1}{\gamma_1}}}\Big)^{D_Y}\rceil \vee 1\\
           &\leq \prod_{j_1=0}^{J_1}\prod_{\psi_1\in  \Psi_{j_1}^{d_Y}} \Big(\frac{24\sqrt{D_Y}L_1C\, J_12^{-\frac{d_Yj_1}{4\gamma_1}-\frac{j_1\beta_Y}{2}}}{\varepsilon^{\frac{1}{\gamma_1}}}\Big)^{D_Y}\vee 1.
         \end{aligned}
    \end{equation*}
Hence there exist  constants $C_1,C_2$ so that for any $\gamma_1\geq \frac{d_Y}{2\beta_Y}$,
 \begin{equation*}
        \log   \mathbf N(\m G,d_{\infty}^{\gamma_1},\varepsilon)\leq C_1\sum_{j_1=0}^{J_1}2^{d_Yj_1}\log\left(\frac{C_2 J_12^{-\frac{d_Yj_1}{4\gamma_1}-\frac{j_1\beta_Y}{2}}}{\varepsilon^{\frac{1}{\gamma_1}}}\vee 1\right).
    \end{equation*}
    When $\frac{d_Y}{\beta_Y}>2\gamma_1$,  denote
    \begin{equation*}
        s_{j_1 }=\sqrt{\frac{2^{\frac{d_Yj_1 }{2\gamma_1}-j_1\beta_Y}}{2^{\frac{d_YJ_1}{2\gamma_1}-J_1\beta_Y}}}.
    \end{equation*}
It holds that
    \begin{equation*}
        \begin{aligned}
    &S=\sum_{j_1=0}^{J_1}s_{j_1}=\sqrt{\frac{1}{2^{\frac{d_YJ_1}{2\gamma_1}-J_1\beta_Y}}}\cdot \sum_{j_1=0}^{J_1}  2^{\frac{d_Yj_1 }{4\gamma_1}-\frac{j_1\beta_Y}{2}}\\
           &=\sqrt{\frac{1}{2^{\frac{d_YJ_1}{2\gamma_1}-J_1\beta_Y}}}\cdot \frac{2^{\frac{(d_Y-2\beta_Y\gamma_1)(J_1+1)}{4\gamma_1}}-1}{2^{\frac{(d_Y-2\beta_Y\gamma_1)}{4\gamma_1}}-1}\ \\
           &\leq J_1\wedge \frac{2^{\frac{(d_Y-2\beta_Y\gamma_1)}{4\gamma_1}}}{2^{\frac{(d_Y-2\beta_Y\gamma_1)}{4\gamma_1}}-1}.
                   \end{aligned}
    \end{equation*}
So if for any $j_1\in [J_1]$ and $\psi_1\in  \Psi_{j_1}^{d_Y}$,
    \begin{equation*}
 \big\|g_{\psi_1 }-g_{\psi_1 }'\big\|\leq \frac{\varepsilon^{\frac{1}{\gamma_1}} s_{j_1}}{CS} 2^{-\frac{d_Yj_1}{2}},
    \end{equation*}
    then 
    \begin{equation*}
    \begin{aligned}
                \underset{x\in \mb R^{D_X}\atop z\in \mb R^{d_Y}}{\sup}\|G(z,x)-G'(z,x)\|^{\gamma_1}&\leq \Big(C\sum_{j_1=0}^J \frac{\varepsilon^{\frac{1}{\gamma_1}} s_{j_1}}{CS}\Big)^{\gamma_1}= \varepsilon.
    \end{aligned}
    \end{equation*}
     Therefore, there exist constants $C_1,C_2$ so that for any $\gamma_1\geq \frac{d_Y}{2\beta_Y}$, 
    \begin{equation*}
    \begin{aligned}
       \log  \mathbf N(\m G,d_{\infty}^{\gamma_1},\varepsilon)&\leq \sum_{j_1=0}^{J_1}\sum_{\psi_1\in  \Psi_{j_1}^{d_Y}} \log \mathbf N([-L_1\, 2^{-\frac{d_Yj_1}{2}-j_1\beta_Y},
        L_1\, 2^{-\frac{d_Yj_1}{2}-j_1\beta_Y}]^{D_Y}, \frac{\varepsilon^{\frac{1}{\gamma_1}} s_{j_1}}{CS} 2^{-\frac{d_Yj_1}{2}},\|\cdot\|)\\
          &\leq C_1 \sum_{j_1=0}^{J_1}  2^{d_Yj_1} \log\Big(\frac{C_2\, S2^{-j_1\beta_Y}}{\varepsilon^{\frac{1}{\gamma_1}} s_{j_1}}\vee 1\Big),
         \end{aligned}
    \end{equation*}
  which completes the proof.


\subsection{Proof of Lemma~\ref{claimlemma11}}\label{proof:claimlemma11}
Fix an $x^*\in \mb N_{\varepsilon_j^x}^x$, then for any $\psi^*\in \Psi_j^{D_Y}$
and  $x\in \mb B_{\m M_X}(x^*,2\varepsilon_j^x)$, it holds that
 \begin{equation*}
     \begin{aligned}
       &\sum_{k=1}^{K^*} \int_{\mb B_{\mb R^{d_Y}}(\mathbf{0},\tau_1)} 2^{\frac{j(d_Y-D_Y)}{2}}\psi^*(G_{[k]}^*(z,x)) v_\sk^*(z,x)\,\dd z-\sum_{k=1}^{K^*} \int_{\mb B_{\mb R^{d_Y}}(\mathbf{0},\tau_1)} 2^{\frac{j(d_Y-D_Y)}{2}}\psi^*(G_{[k],x^*}^\dagger(z,x)) v_{\sk,x^*}^\dagger(z,x) \,\dd z\\
     &=\sum_{k\in [K^*]\atop \|x^*-x_k^*\|\leq \tau_2+2\varepsilon_j^x} \int_{\mb B_{\mb R^{d_Y}}(\mathbf{0},\tau_1)} 2^{\frac{j(d_Y-D_Y)}{2}}\psi^*(G_{[k]}^*(z,x)) v_\sk^*(z,x)\,\dd z\\
     &-\sum_{k\in [K^*]\atop \|x^*-x_k^*\|\leq \tau_2+2\varepsilon_j^x} \int_{\mb B_{\mb R^{d_Y}}(\mathbf{0},\tau_1)} 2^{\frac{j(d_Y-D_Y)}{2}}\psi^*(G_{[k],x^*}^\dagger(z,x)) v_{\sk,x^*}^\dagger(z,x) \,\dd z\\  
         &=\underbrace{\sum_{k\in [K^*]\atop \|x^*-x_k^*\|\leq \tau_2+2\varepsilon_j^x} \int_{\mb B_{\mb R^{d_Y}}(\mathbf{0},\tau_1)} 2^{\frac{j(d_Y-D_Y)}{2}}\big(\psi^*(G_{[k]}^*(z,x))-\psi^*(G_{[k],x^*}^\dagger(z,x))\big) v_{\sk,x^*}^\dagger(z,x) \,\dd z}_{ (E_A)}\\
         &+\underbrace{\sum_{k\in [K^*]\atop \|x^*-x_k^*\|\leq \tau_2+2\varepsilon_j^x} \int_{\mb B_{\mb R^{d_Y}}(\mathbf{0},\tau_1)} 2^{\frac{j(d_Y-D_Y)}{2}}\psi^*(G_{[k]}^*(z,x)) \big(v_\sk^*(z,x)-v_{\sk,x^*}^\dagger(z,x) \big)\,\dd z}_{ (E_B)}.\\
     \end{aligned}
 \end{equation*}
Let $I_{\psi}$ be a rectangle on which $\psi$ is supported  and $y_{\psi}$ denote the center of $I_{\psi}$. Then for any $\psi^*\in \Psi_{j}^{D_Y}$, $x^*\in \mb N_{\varepsilon_j^x}^x$, $x\in \mb B_{\m M_X}(x^*,2\varepsilon_j^x)$, and $k\in [K^*]$ with $\|x^*-x_k^*\|\leq \tau_2+2\varepsilon_j^x$, we have
\begin{equation*}
\begin{aligned}
     &\{z\in \mb B_{\mb R^{d_Y}}(\mathbf{0},\tau_1):\, \psi^*(G_{\sk}^{*}(z,x))-\psi^*(G_{\sk,x^*}^{\dagger}(z,x))\neq 0\}\\
  &\subset \{z\in \mb B_{\mb R^{d_Y}}(\mathbf{0},\tau_1):\, \psi^*(G_{\sk}^{*}(z,x))\neq 0\}\cup \{z\in \mb B_{\mb R^{d_Y}}(\mathbf{0},\tau_1):\, \psi^*(G_{\sk,x^*}^{\dagger}(z,x))\neq 0\}\\
 &\subset \{z\in \mb B_{\mb R^{d_Y}}(\mathbf{0},\tau_1):\, \|y_{\psi^*}-G_{\sk,x^*}^{\dagger}(z,x)\|< C\,2^{-j}\}\cup \{z\in \mb B_{\mb R^{d_Y}}(\mathbf{0},\tau_1):\, \|y_{\psi^*}-G_{\sk}^{*}(z,x)\|< C\,2^{-j}\}\\
 &\subset \{z\in \mb B_{\mb R^{d_Y}}(\mathbf{0},\tau_1):\, \|y_{\psi^*}-G_{\sk}^{*}(z,x)\|< C_1\,2^{-j}\},
\end{aligned}
\end{equation*}
where we have used the fact that 
\begin{equation*}
    \begin{aligned}
       \|G_{\sk,x^*}^{\dagger}(z,x)-G_{\sk}^{*}(z,x)\|&\lesssim   2^{-j\beta_Y}+(\log n)\cdot(\varepsilon_j^x)^{\beta_X}\\
       &\lesssim 2^{-2j}+(\log n)\cdot(\varepsilon_j^x)^{\alpha_X+\frac{\alpha_X}{\alpha_Y}}\\
       &\lesssim 2^{-2j}+(\log n)\cdot (\varepsilon_j^x)^{\alpha_X} \left(2^{\frac{Jd_Y}{2\alpha_X+d_X}} (\frac{n}{\log n})^{-\frac{1}{2\alpha_X+d_X}}\right)^{\frac{\alpha_X}{\alpha_Y}}\\
       &\lesssim 2^{-2j}+(\log n)\cdot (\varepsilon_j^x)^{\alpha_X} 2^{-J}\\
       &\lesssim 2^{-j}.
    \end{aligned}
\end{equation*}
Hence,
\begin{equation*}
\begin{aligned}
        &2^{\frac{j(D_Y-d_Y)}{2}}\cdot (E_A)=\sum_{k\in [K^*]\atop \|x^*-x_k^*\|\leq \tau_2+2\varepsilon_j^x} \int_{\mb B_{\mb R^{d_Y}}(\mathbf{0},\tau_1)}  \big(\psi^*(G_{[k]}^*(z,x))-\psi^*(G_{[k],x^*}^\dagger(z,x))\big) v_{\sk,x^*}^\dagger(z,x) \,\dd z\\
    &=\sum_{k\in [K^*]\atop \|x^*-x_k^*\|\leq \tau_2+2\varepsilon_j^x} \int_{\{z\in \mb B_{\mb R^{d_Y}}(\mathbf{0},\tau_1):\,\|G^*_\sk(z,x^*)-y_{\psi^*}\|\leq C_1\,2^{-j}\}} \big(\psi^*(G_{[k]}^*(z,x))-\psi^*(G_{[k],x^*}^\dagger(z,x))\big) v_{\sk,x^*}^\dagger(z,x) \,\dd z.
    \end{aligned}
\end{equation*}
Based on
 \begin{equation*} 
            \underset{z\in \mb B_{\mb R^{d_Y}}(\mathbf{0},\tau_1)\atop x\in \mb B_{\m M_X}(x^*,2\varepsilon_j^x)}{\sup} \| G_{[k],x^*}^{\dagger}(z,x)-G^*_{[k]}(z,x)\|   \lesssim 2^{-j\beta_Y}+ (\log n)\cdot(\varepsilon_j^x)^{\beta_X},
 \end{equation*}
 and  
 \begin{equation*}
   \underset{z\in \mb B_{\mb R^{d_Y}}(\mathbf{0},\tau_1)\atop x\in \mb B_{\m M_X}(x^*,2\varepsilon_j^x)}{\sup} \| v_{[k],x^*}^{\dagger}(z,x)-v^*_{[k]}(z,x)\|\lesssim 2^{-j\alpha_Y}+ (\log n) \cdot (\varepsilon_j^x)^{\alpha_X},
 \end{equation*}
 we can verify that
\begin{equation*}
\begin{aligned}
      &| (E_A)|\\
      &\lesssim 2^{\frac{j(d_Y-D_Y)}{2}}\sum_{k\in [K^*]\atop \|x^*-x_k^*\|\leq \tau_2+2\varepsilon_j^x} \int_{\{z\in \mb B_{\mb R^{d_Y}}(\mathbf{0},\tau_1):\,\|G^*_\sk(z,x^*)-y_{\psi^*}\|\leq C_1\,2^{-j}\}} \big|\psi^*(G_{[k]}^*(z,x))-\psi^*(G_{[k],x^*}^\dagger(z,x))\big|  \,\dd z\\
      &\lesssim 2^{\frac{d_Yj}{2}+j} \sum_{k\in [K^*]\atop \|x^*-x_k^*\|\leq \tau_2+2\varepsilon_j^x}\int_{\{z\in \mb B_{\mb R^{d_Y}}(\mathbf{0},\tau_1):\,\|G^*_\sk(z,x^*)-y_{\psi^*}\|\leq C_1\,2^{-j}\}} \big\| G_{[k]}^*(z,x))-G_{[k],x^*}^\dagger(z,x)\big\|\,\dd z\\
      &\lesssim 2^{-\frac{jd_Y}{2}}\cdot(2^{-j(\beta_Y-1)}+2^j\cdot \log n\cdot (\varepsilon_j^x)^{\beta_X}).
\end{aligned}
\end{equation*}
Let $\wt\beta_X=\alpha_X+\frac{\alpha_X}{\alpha_Y}$, using the Taylor's theorem for $\psi^*$, we have
\begin{equation*}
     \begin{aligned}
       & \psi^*(G_{[k]}^*(z,x)) =\sum_{l\in \mb N_{0}^{D_Y}\atop  0\leq |l|\leq \lfloor\wt \beta_X\rfloor} \frac{\psi^*{}^{(l)}(G_{[k]}^*(z,x^*))}{l!}(G_{[k]}^*(z,x)-G_{[k]}^*(z,x^*))^l\\
       &\qquad+\sum_{l\in \mb N_{0}^{D_Y}\atop|l|=\lfloor\wt \beta_X\rfloor+1} \frac{|l|}{l!} \int_{0}^{1}(1-t)^{\lfloor\wt \beta_X\rfloor} \psi^*{}^{(l)}(G_{[k]}^*(z,x^*)+t(G_{[k]}^*(z,x)-G_{[k]}^*(z,x^*))) \,\dd t \cdot (G_{[k]}^*(z,x)-G_{[k]}^*(z,x^*))^l
   \end{aligned}
 \end{equation*}
 and 
  \begin{equation*}
     \begin{aligned}
       & \psi^*(G_{[k],x^*}^\dagger(z,x)) =\sum_{l\in \mb N_{0}^{D_Y}\atop  0\leq |l|\leq \lfloor\wt \beta_X\rfloor} \frac{\psi^*{}^{(l)}(G_{[k],x^*}^\dagger(z,x^*))}{l!}(G_{[k],x^*}^\dagger(z,x)-G_{[k],x^*}^\dagger(z,x^*))^l\\
       &+\sum_{l\in \mb N_{0}^{D_Y}\atop|l|=\lfloor\wt \beta_X\rfloor+1} \frac{|l|}{l!} \int_{0}^{1}(1-t)^{\lfloor\wt \beta_X\rfloor} \psi^*{}^{(l)}(G_{[k],x^*}^\dagger(z,x^*)+t(G_{[k],x^*}^\dagger(z,x)-G_{[k],x^*}^\dagger(z,x^*))) \,\dd t\\
       &\qquad\qquad\cdot (G_{[k],x^*}^\dagger(z,x)-G_{[k],x^*}^\dagger(z,x^*))^l.
   \end{aligned}
 \end{equation*}
Then we can obtain
\begin{equation}\label{eqn:boundtermA}
\begin{aligned}
   & \bigg| (E_A)-\sum_{k\in [K^*]\atop \|x^*-x_k^*\|\leq \tau_2+2\varepsilon_j^x} \int_{\{z\in \mb B_{\mb R^{d_Y}}(\mathbf{0},\tau_1):\,\|G^*_\sk(z,x^*)-y_{\psi^*}\|\leq C\,2^{-j}\}} 2^{\frac{j(d_Y-D_Y)}{2}}\\
        & \Big(\sum_{l\in \mb N_{0}^{D_Y}\atop  0\leq |l|\leq \lfloor\wt\beta_X\rfloor} \frac{\psi^*{}^{(l)}(G_{[k]}^*(z,x^*))}{l!}(G_{[k]}^*(z,x)-G_{[k]}^*(z,x^*))^l- \frac{\psi^*{}^{(l)}(G_{[k],x^*}^\dagger(z,x^*))}{l!}(G_{[k],x^*}^\dagger(z,x)-G_{[k],x^*}^\dagger(z,x^*))^l\Big) \\
        &\qquad\qquad \qquad \cdot v^{\dagger}_{\sk,x^*}(z,x)\,\dd z\bigg|\\
        &\leq \sum_{k\in [K^*]\atop \|x^*-x_k^*\|\leq \tau_2+2\varepsilon_j^x}  \int_{\{z\in \mb B_{\mb R^{d_Y}}(\mathbf{0},\tau_1):\,\|G^*_\sk(z,x^*)-y_{\psi^*}\|\leq C\,2^{-j}\}} 2^{\frac{j(d_Y-D_Y)}{2}}\sum_{l\in \mb N_{0}^{D_Y}\atop|l|=\lfloor\wt\beta_X\rfloor+1} \frac{|l|}{l!} \int_{0}^{1}(1-t)^{\lfloor\wt\beta_X\rfloor} \\
&\quad\Big|\psi^*{}^{(l)}(G_{[k]}^*(z,x^*)+t(G_{[k]}^*(z,x)-G_{[k]}^*(z,x^*)))- \psi^*{}^{(l)}(G_{[k],x^*}^\dagger(z,x^*)+t(G_{[k],x^*}^\dagger(z,x)-G_{[k],x^*}^\dagger(z,x^*))) \Big|\,\dd t\\
&\underbrace{\qquad\cdot \big|(G_{[k]}^*(z,x)-G_{[k]}^*(z,x^*))^l\big|  \cdot|v_{\sk,x^*}^\dagger(z,x)\big|\,\dd z\qquad\qquad\qquad\qquad\qquad\qquad\qquad\qquad\qquad\qquad\qquad\qquad}_{(E_C)}\\
&+\sum_{k\in [K^*]\atop \|x^*-x_k^*\|\leq \tau_2+2\varepsilon_j^x}  \int_{\{z\in \mb B_{\mb R^{d_Y}}(\mathbf{0},\tau_1):\,\|G^*_\sk(z,x^*)-y_{\psi^*}\|\leq C\,2^{-j}\}} 2^{\frac{j(d_Y-D_Y)}{2}}\sum_{l\in \mb N_{0}^{D_Y}\atop|l|=\lfloor\wt\beta_X\rfloor+1} \frac{|l|}{l!} \cdot\Big|\int_{0}^{1}(1-t)^{\lfloor\wt\beta_X\rfloor} \\
&\quad \Big|\psi^*{}^{(l)}(G_{[k],x^*}^\dagger(z,x^*)+t(G_{[k],x^*}^\dagger(z,x)-G_{[k],x^*}^\dagger(z,x^*)))\,\dd t\Big|\\
&\underbrace{\qquad\cdot \big|(G_{[k]}^*(z,x)-G_{[k]}^*(z,x^*))^l-(G_{[k],x^*}^\dagger(z,x)-G_{[k],x^*}^\dagger(z,x^*))^l\big|  \cdot\big|v_{\sk,x^*}^\dagger(z,x)\big|\,\dd z.\qquad\qquad }_{(E_D)}\\
        \end{aligned}
\end{equation}
We first bound the term $(E_C)$. Notice that
 \begin{equation*}
     \begin{aligned}
          &\big|\psi^*{}^{(l)}(G_{[k]}^*(z,x^*)+t(G_{[k]}^*(z,x)-G_{[k]}^*(z,x^*)))-\psi^*{}^{(l)}(G_{[k],x^*}^\dagger(z,x^*)+t(G_{[k],x^*}^\dagger(z,x)-G_{[k],x^*}^\dagger(z,x^*)))\big| \\
          &\lesssim 2^{\frac{D_Yj}{2}} 2^{j(|l|+1)}\cdot(\|G_{[k],x^*}^\dagger(z,x)-G_{[k]}^*(z,x)\|+\|G_{[k],x^*}^\dagger(z,x^*)-G_{[k]}^*(z,x^*)\|)\\
          &\lesssim 2^{\frac{D_Yj}{2}} 2^{j(\lfloor\wt\beta_X\rfloor+2)} (2^{-j\beta_Y}+\log n\cdot(\varepsilon_j^x)^{\beta_X}) 
     \end{aligned}
 \end{equation*}
 and  
 \begin{equation*}
     \begin{aligned}
        \big| (G_{[k]}^*(z,x)-G_{[k]}^*(z,x^*))^l\big|\lesssim\left\{\begin{array}{cc}
          (\varepsilon_j^x)^{|l|}   & \beta_X\geq 1 \\
          (\varepsilon_j^x)^{\beta_X|l|}   & \beta_X<1.
        \end{array}  
        \right.
     \end{aligned}
 \end{equation*}
Using the conditions: $\beta_X\geq \alpha_X+\frac{\alpha_X}{\alpha_Y}$, $\beta_Y\geq \alpha_Y+1$, $\alpha_Y\geq \alpha_X$, and considering that for any $j\in \{0\}\cup [J]$ with $J=\lceil \frac{1}{2\alpha_Y+d_Y+d_X\frac{\alpha_Y}{\alpha_X}\cdot \log_2(\frac{n}{\log n})}\rceil$, it holds that
 \begin{equation*}
 \begin{aligned}
         2^{-j\alpha_Y}&\geq 2^{-J\alpha_Y}=(\frac{n}{\log n})^{-\frac{\alpha_Y}{2\alpha_Y+d_Y+d_X\frac{\alpha_Y}{\alpha_X}}}\\
         &=(\frac{n}{\log n})^{-\frac{\alpha_X}{2\alpha_X+d_X}}(\frac{n}{\log n})^{\frac{\alpha_X}{2\alpha_X+d_X}-\frac{\alpha_Y}{2\alpha_Y+d_Y+d_X\frac{\alpha_Y}{\alpha_X}}}\\
         &=(\frac{n}{\log n})^{-\frac{\alpha_X}{2\alpha_X+d_X}} 2^{J\frac{d_Y\alpha_X}{2\alpha_X+d_X}}\\
         &=(\varepsilon_J^x)^{\alpha_X}\geq (\varepsilon_j^x)^{\alpha_X}.
 \end{aligned}
 \end{equation*}
We can conclude, when $\beta_X\geq 1$,
 \begin{equation*}
     \begin{aligned}
&(E_C)=\sum_{k\in [K^*]\atop \|x^*-x_k^*\|\leq \tau_2+2\varepsilon_j^x}  \int_{\{z\in \mb B_{\mb R^{d_Y}}(\mathbf{0},\tau_1):\,\|G^*_\sk(z,x^*)-y_{\psi^*}\|\leq C\,2^{-j}\}} 2^{\frac{j(d_Y-D_Y)}{2}}\sum_{l\in \mb N_{0}^{D_Y}\atop|l|=\lfloor\wt\beta_X\rfloor+1} \frac{|l|}{l!} \int_{0}^{1}(1-t)^{\lfloor\wt\beta_X\rfloor} \\
&\quad\Big|\psi^*{}^{(l)}(G_{[k]}^*(z,x^*)+t(G_{[k]}^*(z,x)-G_{[k]}^*(z,x^*)))- \psi^*{}^{(l)}(G_{[k],x^*}^\dagger(z,x^*)+t(G_{[k],x^*}^\dagger(z,x)-G_{[k],x^*}^\dagger(z,x^*))) \Big|\,\dd t\\
&\qquad\cdot \big|(G_{[k]}^*(z,x)-G_{[k]}^*(z,x^*))^l\big|  \cdot|v_{\sk,x^*}^\dagger(z,x)\big|\,\dd z\\
         &\lesssim 2^{-\frac{jd_Y}{2}} 2^{j(\lfloor\wt\beta_X\rfloor+2)}  (2^{-j\beta_Y}+\log n\cdot(\varepsilon_j^x)^{\beta_X})(\varepsilon_j^x)^{\lfloor\wt\beta_X\rfloor+1}\\
         &=  2^{-\frac{jd_Y}{2}}(\varepsilon_j^x)^{\alpha_X}\cdot 2^{j(\lfloor\wt\beta_X\rfloor+2)} (2^{-j\beta_Y}+\log n\cdot(\varepsilon_j^x)^{\beta_X})(\varepsilon_j^x)^{\lfloor\wt\beta_X\rfloor+1-\alpha_X}\\
         &\lesssim   2^{-\frac{jd_Y}{2}}(\varepsilon_j^x)^{\alpha_X}\cdot \Big(2^{j(\lfloor\wt\beta_X\rfloor+2)} 2^{-j\beta_Y}(2^{-j\frac{\alpha_Y}{\alpha_X}})^{\lfloor\wt\beta_X\rfloor+1-\alpha_X}+  \log n\cdot(\varepsilon_j^x)^{\beta_X+\lfloor\wt\beta_X\rfloor+1-\alpha_X-\frac{\alpha_X}{\alpha_Y}(\lfloor\wt\beta_X\rfloor+2)}\Big)\\
         &= 2^{-\frac{jd_Y}{2}} (\varepsilon_j^x)^{\alpha_X}\cdot \Big(2^{-j(\frac{\alpha_Y}{\alpha_X}-1)(\lfloor\wt\beta_X\rfloor+1)} 2^{-j(\beta_Y-1-\alpha_Y)} +  \log n\cdot (\varepsilon_j^x)^{\beta_X-\alpha_X-\frac{\alpha_X}{\alpha_Y}+(\lfloor\wt\beta_X\rfloor+1)(1-\frac{\alpha_X}{\alpha_Y})}\Big)\\
         &\lesssim (\log n) \cdot 2^{-\frac{jd_Y}{2}}(\varepsilon_j^x)^{\alpha_X}.
     \end{aligned}
 \end{equation*}
When $\beta_X<1$, we have $\lfloor \wt \beta_X\rfloor=\lfloor \alpha_X+\frac{\alpha_X}{\alpha_Y}\rfloor\leq \lfloor\beta_X\rfloor=0$, and
 \begin{equation*}
     \begin{aligned}
(E_C)&\lesssim 2^{-\frac{jd_Y}{2}} 2^{2j}  (2^{-j\beta_Y}+\log n\cdot(\varepsilon_j^x)^{\beta_X})(\varepsilon_j^x)^{\beta_X}\\
         &=  2^{-\frac{jd_Y}{2}}(\varepsilon_j^x)^{\alpha_X}\cdot 2^{2j} (2^{-j\beta_Y}+\log n\cdot(\varepsilon_j^x)^{\beta_X})(\varepsilon_j^x)^{\beta_X-\alpha_X}\\
         &\lesssim   2^{-\frac{jd_Y}{2}}(\varepsilon_j^x)^{\alpha_X}\cdot \Big(2^{2j} 2^{-j\beta_Y}(2^{-j\frac{\alpha_Y}{\alpha_X}})^{\beta_X-\alpha_X}+  \log n\cdot(\varepsilon_j^x)^{\frac{2\alpha_X}{\alpha_Y}}2^{2j}\Big)\\
         &=  2^{-\frac{jd_Y}{2}}(\varepsilon_j^x)^{\alpha_X}\cdot \Big(2^{2j} 2^{-j\beta_Y}(2^{-j\frac{\alpha_Y}{\alpha_X}})^{\beta_X-\alpha_X}+  \log n\cdot(2^{-j\frac{\alpha_Y}{\alpha_X}})^{\frac{2\alpha_X}{\alpha_Y}}2^{2j}\Big)\\
         &\lesssim (\log n) \cdot 2^{-\frac{jd_Y}{2}}(\varepsilon_j^x)^{\alpha_X}.
     \end{aligned}
 \end{equation*}
Furthermore, for bounding the term $(E_D)$, notice that for any $x\in \mb B_{\m M_X}(x^*,2\varepsilon_j^x)$ and $z\in \mb B_{\mb R^{d_Y}}(\mathbf{0},\tau_1)$,  
 \begin{equation*}
     \|G_\sk(z,x)-G_\sk(z,x^*)\|\lesssim \|x-x^*\|^{1\wedge \beta_X}\lesssim (\varepsilon_j^x)^{1\wedge \beta_X},
 \end{equation*}
 and when $\beta_X\leq 1$, it holds that $\big\| G_{\sk,x^*}^\dagger(z,x)-G_{\sk,x^*}^\dagger(z,x^*)\|=0$; when $\beta_X>1$
 \begin{equation*}
     \begin{aligned}
         &\big\| G_{\sk,x^*}^\dagger(z,x)-G_{\sk,x^*}^\dagger(z,x^*)\|\\
         &=\Big\|\sum_{s=0}^{j} \sum_{\psi\in \wt\Psi_s^{d_Y}}\sum_{l\in \mb N_0^{D_X}\atop 1\leq |l|<\beta_X} \int_{\mb R^{d_Y}} \frac{1}{l!}G_{[k]}^{*(\mathbf{0},l)}(t,x^*) (x-x^*)^{l}\psi(t)\,\dd t\cdot \psi(z)\Big\| \\
         &=\Big\|\sum_{s=0}^{j} \sum_{\psi\in \wt\Psi_s^{d_Y}}\sum_{l\in \mb N_0^{D_X}\atop   |l|=1} \int_{\mb R^{d_Y}} \frac{1}{l!}G_{[k]}^{*(\mathbf{0},l)}(t,x^*) (x-x^*)^{l}\psi(t)\,\dd t\cdot \psi(z)\Big\| + o(\varepsilon_j^x)\\
         &\overset{(i)}{\lesssim} \varepsilon_j^x \sum_{s=0}^j 2^{-s(\beta_Y-\frac{\beta_Y}{\beta_X})}\lesssim \varepsilon_j^x \cdot(j \wedge \frac{1}{\beta_X-1}),
     \end{aligned}
 \end{equation*}
where $(i)$ uses that for any $l\in \mb N_0^{D_X}$ with  $|l|=1$, $G_{[k]}^{*(\mathbf{0},l)}(\cdot,x^*)\in \m H^{\beta_Y-\beta_Y/\beta_X}_{L,D_Y}(\mb R^{d_Y})$. Together with $\| G_{[k],x^*}^{\dagger}(z,x)-G^*_{[k]}(z,x)\|   \lesssim 2^{-j\beta_Y}+ (\log n)\cdot(\varepsilon_j^x)^{\beta_X}$, we can derive that, for any $l\in \mb N_{0}^{D_Y}$ with $|l|=\lfloor\wt\beta_X\rfloor+1$, and any $i\in [D_Y]$ with $l_i\geq 1$, 
 \begin{equation*}
     \begin{aligned}
     & \big|(G_{[k]i}^*(z,x)-G_{[k]i}^*(z,x^*))^{l_i}- (G_{[k]i}^\dagger(z,x)-G_{[k]i}^\dagger(z,x^*))^{l_i}\big|\\
     &=\big|(G_{[k]i}^*(z,x)-G_{[k]i}^*(z,x^*)- G_{[k]i}^\dagger(z,x)+G_{[k]i}^\dagger(z,x^*))\\
     &\qquad\cdot\sum_{i_1=1}^{l_i} (G_{[k]i}^*(z,x)-G_{[k]i}^*(z,x^*))^{l_i-i_1}(G_{[k]i}^\dagger(z,x)-G_{[k]i}^\dagger(z,x^*))^{i_1-1}\big|\\
     &\lesssim  \left\{\begin{array}{cc}
       (\log n)\cdot ((2^{-j\beta_Y}+(\varepsilon_j^x)^{\beta_X})\wedge \varepsilon_j^x) \big( (j\wedge \frac{1}{\beta_X-1})\cdot\varepsilon_j^x\big)^{l_i-1}   & \beta_X> 1 \\
      (\varepsilon_j^x)^{\beta_X}   & \beta_X\leq 1\\
     \end{array}
      \right.\\
 &\lesssim  \left\{\begin{array}{cc}
     (\log n)\cdot (j \wedge \frac{1}{\beta_X-1})^{l_i-1}\cdot (2^{-j\beta_Y}+(\varepsilon_j^x)^{\beta_X})^{\frac{\alpha_X}{\alpha_Y}}(\varepsilon_j^x)^{l_i-\frac{\alpha_X}{\alpha_Y}}  & \beta_X> 1 \\
      (\varepsilon_j^x)^{\beta_X}    & \beta_X\leq 1,\\
     \end{array}
      \right.\\
     \end{aligned}
 \end{equation*}
 where $G_{\sk i}^*(z,x)$ denote the $i$-th component of the $D_Y$-dimensional vector $G_\sk^*(z,x)$. Therefore, when $\beta_X>1$,
 \begin{equation*}
     \begin{aligned}
        &\big|(G_{[k]}^*(z,x)-G_{[k]}^*(z,x^*))^l- (G_{[k],x^*}^\dagger(z,x)-G_{[k],x^*}^\dagger(z,x^*))^l\big| \\
        &\lesssim (\log n)\cdot (j \wedge \frac{1}{\beta_X-1})^{\lfloor\wt\beta_X\rfloor} (2^{-j\beta_Y}+(\varepsilon_j^x)^{\beta_X})^{\frac{\alpha_X}{\alpha_Y}}(\varepsilon_j^x)^{\lfloor\wt\beta_X\rfloor+1-\frac{\alpha_X}{\alpha_Y}},
     \end{aligned}
 \end{equation*}
 and 
 \begin{equation*}
     \begin{aligned}
         &(E_D)=\sum_{k\in [K^*]\atop \|x^*-x_k^*\|\leq \tau_2+2\varepsilon_j^x}  \int_{\{z\in \mb B_{\mb R^{d_Y}}(\mathbf{0},\tau_1):\,\|G^*_\sk(z,x^*)-y_{\psi^*}\|\leq C\,2^{-j}\}} 2^{\frac{j(d_Y-D_Y)}{2}}\sum_{l\in \mb N_{0}^{D_Y}\atop|l|=\lfloor\wt\beta_X\rfloor+1} \frac{|l|}{l!} \cdot\Big|\int_{0}^{1}(1-t)^{\lfloor\wt\beta_X\rfloor} \\
&\quad \Big|\psi^*{}^{(l)}(G_{[k],x^*}^\dagger(z,x^*)+t(G_{[k],x^*}^\dagger(z,x)-G_{[k],x^*}^\dagger(z,x^*)))\,\dd t\Big|\\
&\qquad\cdot \big|(G_{[k]}^*(z,x)-G_{[k]}^*(z,x^*))^l-(G_{[k],x^*}^\dagger(z,x)-G_{[k],x^*}^\dagger(z,x^*))^l\big|  \cdot|v_{\sk,x^*}^\dagger(z,x)\big|\,\dd z\\
&\lesssim (\log n)\cdot (j\wedge \frac{1}{\beta_X-1})^{\lfloor\wt\beta_X\rfloor} \cdot 2^{-\frac{jd_Y}{2}}2^{j(\lfloor\wt\beta_X\rfloor+1)} (2^{-j\beta_Y}+(\varepsilon_j^x)^{\beta_X})^{\frac{\alpha_X}{\alpha_Y}}(\varepsilon_j^x)^{\lfloor\wt\beta_X\rfloor+1-\frac{\alpha_X}{\alpha_Y}}\\
&= (\log n) \cdot (j\wedge \frac{1}{\beta_X-1})^{\lfloor\wt\beta_X\rfloor} \cdot 2^{-\frac{jd_Y}{2}}(\varepsilon_j^x)^{\alpha_X}\Big(2^{j(\lfloor\wt\beta_X\rfloor+1)} (2^{-j\beta_Y}+(\varepsilon_j^x)^{\beta_X})^{\frac{\alpha_X}{\alpha_Y}}(\varepsilon_j^x)^{\lfloor\wt\beta_X\rfloor+1-\frac{\alpha_X}{\alpha_Y}-\alpha_X}\Big)\\
&\lesssim (\log n) \cdot (j \wedge \frac{1}{\beta_X-1})^{\lfloor\wt\beta_X\rfloor} \cdot 2^{-\frac{jd_Y}{2}} (\varepsilon_j^x)^{\alpha_X}\Big(2^{j(\lfloor\wt\beta_X\rfloor+1)} 2^{-j\beta_Y\frac{\alpha_X}{\alpha_Y}}2^{-j\frac{\alpha_Y}{\alpha_X}(\lfloor\wt\beta_X\rfloor+1-\frac{\alpha_X}{\alpha_Y}-\alpha_X)}\\
&\qquad\qquad\qquad+ (\varepsilon_j^x)^{-\frac{\alpha_X}{\alpha_Y}(\lfloor\wt\beta_X\rfloor+1)} (\varepsilon_j^x)^{\beta_X \frac{\alpha_X}{\alpha_Y}} (\varepsilon_j^x)^{\lfloor\wt\beta_X\rfloor+1-\frac{\alpha_X}{\alpha_Y}-\alpha_X}\Big)\\
&= (\log n) \cdot (j \wedge \frac{1}{\beta_X-1})^{\lfloor\wt\beta_X\rfloor} \cdot 2^{-\frac{jd_Y}{2}}(\varepsilon_j^x)^{\alpha_X}\Big(2^{-j\big((\frac{\alpha_Y}{\alpha_X}-1)(\lfloor\wt\beta_X\rfloor+1)+\beta_Y\frac{\alpha_X}{\alpha_Y}-1-\alpha_Y\big)}\\
&\qquad\qquad\qquad+ (\varepsilon_j^x)^{(\lfloor\wt\beta_X\rfloor+1)(1-\frac{\alpha_X}{\alpha_Y})+\beta_X\frac{\alpha_X}{\alpha_Y}-\frac{\alpha_X}{\alpha_Y}-\alpha_X}\Big)\lesssim \log n \cdot  2^{-\frac{jd_Y}{2}} (\varepsilon_j^x)^{\alpha_X},
     \end{aligned}
 \end{equation*}
 where the last inequality uses  that 
 \begin{equation*}
     \begin{aligned}
         (\frac{\alpha_Y}{\alpha_X}-1)(\lfloor\wt\beta_X\rfloor+1)+\beta_Y\frac{\alpha_X}{\alpha_Y}-1-\alpha_Y&\geq (\frac{\alpha_Y}{\alpha_X}-1) \wt\beta_X +\beta_Y\frac{\alpha_X}{\alpha_Y}-1-\alpha_Y\\
         &\geq   (\frac{\alpha_Y}{\alpha_X}-1) (\alpha_X+\frac{\alpha_X}{\alpha_Y}) +(\alpha_Y+1)\frac{\alpha_X}{\alpha_Y}-1-\alpha_Y=0,
     \end{aligned}
 \end{equation*}
 and
 \begin{equation*}
     \begin{aligned}
         (\lfloor\wt\beta_X\rfloor+1)(1-\frac{\alpha_X}{\alpha_Y})+\beta_X\frac{\alpha_X}{\alpha_Y}-\frac{\alpha_X}{\alpha_Y}-\alpha_X&\geq  \wt\beta_X(1-\frac{\alpha_X}{\alpha_Y})+\beta_X\frac{\alpha_X}{\alpha_Y}-\frac{\alpha_X}{\alpha_Y}-\alpha_X\\
           &\geq \wt\beta_X-\frac{\alpha_X}{\alpha_Y}-\alpha_X= 0,
     \end{aligned}
 \end{equation*}
 alongside the fact that $\lfloor\wt\beta_X\rfloor+1=\wt \beta_X$   only if $\wt \beta_X$ is an integer. 
Similarly, when $\beta_X\leq 1$, 
 \begin{equation*}
     \begin{aligned}
         (E_D)\lesssim   2^{-\frac{jd_Y}{2}+1} (\varepsilon_j^x)^{\beta_X}\lesssim   2^{-\frac{jd_Y}{2}} (\varepsilon_j^x)^{\alpha_X} 2^j (\varepsilon_j^x)^{\frac{\alpha_X}{\alpha_Y}}\lesssim  2^{-\frac{jd_Y}{2}} (\varepsilon_j^x)^{\alpha_X}.
     \end{aligned}
 \end{equation*}
By combining the bounds for terms $(E_C)$ and $(E_D)$, and using Equation~\eqref{eqn:boundtermA}, we can obtain that
 \begin{equation*}
     \begin{aligned}
        &  (E_A)=\sum_{k\in [K^*]\atop \|x^*-x_k^*\|\leq \tau_2+2\varepsilon_j^x} \int_{\{z\in \mb B_{\mb R^{d_Y}}(\mathbf{0},\tau_1):\,\|G^*_\sk(z,x^*)-y_{\psi^*}\|\leq C\,2^{-j}\}} 2^{\frac{j(d_Y-D_Y)}{2}}\\
        & \Big(\sum_{l\in \mb N_{0}^{D_Y}\atop  0\leq |l|\leq \lfloor\wt\beta_X\rfloor} \frac{\psi^*{}^{(l)}(G_{[k]}^*(z,x^*))}{l!}(G_{[k]}^*(z,x)-G_{[k]}^*(z,x^*))^l- \frac{\psi^*{}^{(l)}(G_{[k],x^*}^\dagger(z,x^*))}{l!}(G_{[k],x^*}^\dagger(z,x)-G_{[k],x^*}^\dagger(z,x^*))^l\Big) \\
        &\qquad\cdot v_{\sk,x^*}^\dagger(z,x) \,\dd z+\m O(\log n \cdot 2^{-\frac{jd_Y}{2}}(\varepsilon_j^x)^{\alpha_X}).
     \end{aligned}
 \end{equation*}
Given that for any $x\in \mb B_{\m M_X}(x^*,2\varepsilon_j^x)$,
 \begin{equation*}
     \begin{aligned}
       G_{[k]}^*(z,x)-G_{[k]}^*(z,x^*)=\sum_{s\in \mb N_{0}^{D_X}\atop 1\leq|s|\leq \lfloor\wt\beta_X\rfloor} \frac{G_{[k]}^*{}^{(\mathbf 0,s)}(z,x^*)}{s!}(x-x^*)^s+\m O((\varepsilon_j^x)^{\wt\beta_X})
     \end{aligned}
 \end{equation*}
 and  considering that $G_{[k],x^*}^\dagger(z,x)$ is polynomial in $x$, 
\begin{equation*}
     \begin{aligned}
       G_{[k],x^*}^\dagger(z,x)-G_{[k],x^*}^\dagger(z,x^*)=\sum_{s\in \mb N_{0}^{D_X}\atop 1\leq|s|\leq \lfloor\wt\beta_X\rfloor} \frac{G_{[k],x^*}^\dagger{}^{(\mathbf 0,s)}(z,x^*)}{s!}(x-x^*)^s,
     \end{aligned}
 \end{equation*}
 where  recall $G^{(\mathbf 0,s)}(z,x)$ denotes the partial derivative of $G(z,\cdot)$ of order $s$ evaluated at $x$. If $\wt\beta_X>1$, it holds
  for any $l\in \mb N_{0}^{D_Y}$ with $1\leq |l|\leq \lfloor\wt\beta_X\rfloor$ that,  
 \begin{equation*}
 \begin{aligned}
     \Big| (G_{[k]}^*(z,x)-G_{[k]}^*(z,x^*))^l-\big(\sum_{s\in \mb N_{0}^{D_X}\atop 1\leq|s|\leq \lfloor\wt\beta_X\rfloor} \frac{G_{[k]}^*{}^{(\mathbf 0,s)}(z,x^*)}{s!}(x-x^*)^s\big)^l\Big|\lesssim  (\varepsilon_j^x)^{\wt\beta_X+|l|-1}.
 \end{aligned}
 \end{equation*}
   Therefore, 
 \begin{equation*}
     \begin{aligned}
         &\sum_{k\in [K^*]\atop \|x^*-x_k^*\|\leq \tau_2+2\varepsilon_j^x} \int_{\{z\in \mb B_{\mb R^{d_Y}}(\mathbf{0},\tau_1):\,\|G^*_\sk(z,x^*)-y_{\psi^*}\|\leq C\,2^{-j}\}} 2^{\frac{j(d_Y-D_Y)}{2}}\\
        & \Big(\sum_{l\in \mb N_{0}^{D_Y}\atop  0\leq |l|\leq \lfloor\wt\beta_X\rfloor} \frac{\psi^*{}^{(l)}(G_{[k]}^*(z,x^*))}{l!}(G_{[k]}^*(z,x)-G_{[k]}^*(z,x^*))^l- \frac{\psi^*{}^{(l)}(G_{[k],x^*}^\dagger(z,x^*))}{l!}(G_{[k],x^*}^\dagger(z,x)-G_{[k],x^*}^\dagger(z,x^*))^l\Big) \\
        &\qquad\cdot v_{\sk,x^*}^\dagger(z,x) \,\dd z\\
        &=   \sum_{k\in [K^*]\atop \|x^*-x_k^*\|\leq \tau_2+2\varepsilon_j^x} \int_{\{z\in \mb B_{\mb R^{d_Y}}(\mathbf{0},\tau_1):\,\|G^*_\sk(z,x^*)-y_{\psi^*}\|\leq C\,2^{-j}\}} 2^{\frac{j(d_Y-D_Y)}{2}} \Big(\psi^*(G_{[k]}^*(z,x^*))-\psi^*(G_{[k],x^*}^\dagger(z,x^*))\\
        &+ \sum_{l\in \mb N_{0}^{D_Y}\atop  1\leq |l|\leq \lfloor\wt\beta_X\rfloor} \frac{\psi^*{}^{(l)}(G_{[k]}^*(z,x^*))}{l!}\big(\sum_{s\in \mb N_{0}^{D_X}\atop 1\leq|s|\leq\lfloor\wt\beta_X\rfloor} \frac{G_{[k]}^*{}^{(\mathbf 0,s)}(z,x^*)}{s!}(x-x^*)^s\big)^l\\
        &- \frac{\psi^*{}^{(l)}(G_{[k],x^*}^\dagger(z,x^*))}{l!}\big(\sum_{s\in \mb N_{0}^{D_X}\atop 1\leq|s|\leq \lfloor\wt\beta_X\rfloor} \frac{G_{[k],x^*}^\dagger{}^{(\mathbf 0,s)}(z,x^*)}{s!}(x-x^*)^s\big)^l\Big) v_{\sk,x^*}^\dagger(z,x) \,\dd z+\m O( 2^{-\frac{jd_Y}{2}}(\varepsilon_j^x)^{\alpha_X}),\\
     \end{aligned}
 \end{equation*}
 where  we have used the fact that 
 \begin{equation*}
    \begin{aligned}
          \sum_{l\in \mb N_{0}^{D_Y}\atop  1\leq |l|\leq \lfloor\wt\beta_X\rfloor} 2^{j|l|} (\varepsilon_j^x)^{\wt\beta_X+|l|-1} \lesssim \sum_{l\in \mb N_{0}^{D_Y}\atop  1\leq |l|\leq \lfloor\wt\beta_X\rfloor}   (\varepsilon_j^x)^{-\frac{\alpha_X}{\alpha_Y}|l|+\alpha_X+\frac{\alpha_X}{\alpha_Y}+|l|-1}
=   \sum_{l\in \mb N_{0}^{D_Y}\atop  1\leq |l|\leq \lfloor\wt\beta_X\rfloor}   (\varepsilon_j^x)^{\alpha_X}(\varepsilon_j^x)^{(1-\frac{\alpha_X}{\alpha_Y})(|l|-1)}
     \lesssim (\varepsilon_j^x)^{\alpha_X}.
    \end{aligned}
 \end{equation*}
Together with the fact that $  v_{[k],x^*}^\dagger(z,x)$ is polynomial in $x$ and
 \begin{equation*}
    \begin{aligned}
    v_{[k],x^*}^\dagger(z,x)=\sum_{s\in \mb N_{0}^{D_X}\atop 0\leq|s|\leq \lfloor\alpha_X\rfloor} \frac{v_{[k]}^\dagger{}^{(\mathbf 0,s)}(z,x^*)}{s!}(x-x^*)^s,
    \end{aligned}
 \end{equation*}
 we can then obtain 
 \begin{equation*}
     \begin{aligned}
         &  (E_A)= \sum_{k\in [K^*]\atop \|x^*-x_k^*\|\leq \tau_2+2\varepsilon_j^x} \int_{\{z\in \mb B_{\mb R^{d_Y}}(\mathbf{0},\tau_1):\,\|G^*_\sk(z,x^*)-y_{\psi^*}\|\leq C\,2^{-j}\}} 2^{\frac{j(d_Y-D_Y)}{2}} \Big(\psi^*(G_{[k]}^*(z,x^*))-\psi^*(G_{[k],x^*}^\dagger(z,x^*))\\
        &+ \sum_{l\in \mb N_{0}^{D_Y}\atop  1\leq |l|\leq \lfloor\wt\beta_X\rfloor} \frac{\psi^*{}^{(l)}(G_{[k]}^*(z,x^*))}{l!}\big(\sum_{s\in \mb N_{0}^{D_X}\atop 1\leq|s|\leq \lfloor\wt\beta_X\rfloor} \frac{G_{[k]}^*{}^{(\mathbf 0,s)}(z,x^*)}{s!}(x-x^*)^s\big)^l\\
        &- \frac{\psi^*{}^{(l)}(G_{[k],x^*}^\dagger(z,x^*))}{l!}\big(\sum_{s\in \mb N_{0}^{D_X}\atop 1\leq|s|\leq \lfloor\wt\beta_X\rfloor} \frac{G_{[k],x^*}^\dagger{}^{(\mathbf 0,s)}(z,x^*)}{s!}(x-x^*)^s\big)^l\Big) \sum_{s\in \mb N_{0}^{D_X}\atop 0\leq|s|\leq \lceil\alpha_X\rfloor} \frac{v_{[k]}^\dagger{}^{(\mathbf 0,s)}(z,x^*)}{s!}(x-x^*)^s\,\dd z\\
       & +\m O(\log n\cdot 2^{-\frac{jd_Y}{2}}(\varepsilon_j^x)^{\alpha_X}).\\
     \end{aligned}
 \end{equation*}
Also notice that we can rewrite 
 \begin{equation*}
     \begin{aligned}
         & \sum_{k\in [K^*]\atop \|x^*-x_k^*\|\leq \tau_2+2\varepsilon_j^x} \int_{\{z\in \mb B_{\mb R^{d_Y}}(\mathbf{0},\tau_1):\,\|G^*_\sk(z,x^*)-y_{\psi^*}\|\leq C\,2^{-j}\}} 2^{\frac{j(d_Y-D_Y)}{2}} \Big(\psi^*(G_{[k]}^*(z,x^*))-\psi^*(G_{[k],x^*}^\dagger(z,x^*))\\
        &+ \sum_{l\in \mb N_{0}^{D_Y}\atop  1\leq |l|\leq \lfloor\wt\beta_X\rfloor} \frac{\psi^*{}^{(l)}(G_{[k]}^*(z,x^*))}{l!}\big(\sum_{s\in \mb N_{0}^{D_X}\atop 1\leq|s|\leq \lfloor\wt\beta_X\rfloor} \frac{G_{[k]}^*{}^{(\mathbf 0,s)}(z,x^*)}{s!}(x-x^*)^s\big)^l\\
        &- \frac{\psi^*{}^{(l)}(G_{[k],x^*}^\dagger(z,x^*))}{l!}\big(\sum_{s\in \mb N_{0}^{D_X}\atop 1\leq|s|\leq\lfloor\wt\beta_X\rfloor} \frac{G_{[k],x^*}^\dagger{}^{(\mathbf 0,s)}(z,x^*)}{s!}(x-x^*)^s\big)^l\Big) \sum_{s\in \mb N_{0}^{D_X}\atop 0\leq|s|\leq \lceil\alpha_X\rfloor} \frac{v_{[k]}^\dagger{}^{(\mathbf 0,s)}(z,x^*)}{s!}(x-x^*)^s\,\dd z\\
        &=\sum_{s\in \mb N_0^{D_X}\atop 0\leq s\leq \lfloor\wt\beta_X\rfloor^2+\lfloor \alpha_X\rfloor} a_{\psi^*,x^*,s}(x-x^*)^s,
     \end{aligned}
 \end{equation*}
 where $|a_{\psi^*,x^*,s}|\leq C\, 2^{j\lfloor\wt\beta_X\rfloor}  (\log n)^{1+\lfloor\wt\beta_X\rfloor}\lesssim n$.  
 
 \medskip
 \noindent Then for term $ (E_B)$,  using the Taylor's theorem for $\psi^*(\cdot)$, $G^*_{\sk}(z,\cdot)$, and $v_{\sk}^*(z,\cdot)$ , we have
          \begin{equation*}
     \begin{aligned}
       & \psi^*(G_{[k]}^*(z,x)) =\sum_{l\in \mb N_{0}^{D_Y}\atop  0\leq |l|\leq \lfloor\wt \beta_X\rfloor} \frac{\psi^*{}^{(l)}(G_{[k]}^*(z,x^*))}{l!}(G_{[k]}^*(z,x)-G_{[k]}^*(z,x^*))^l+\m O(2^{\frac{jD_Y}{2}}((\varepsilon_j^x)^{\beta_X\wedge 1}2^{j})^{\lfloor\wt \beta_X\rfloor+1})\\
        & =\psi^*(G_{[k]}^*(z,x^*))+\sum_{l\in \mb N_{0}^{D_Y}\atop  1\leq |l|\leq \lfloor\wt \beta_X\rfloor} \frac{\psi^*{}^{(l)}(G_{[k]}^*(z,x^*))}{l!}\big(\sum_{s\in \mb N_{0}^{D_X}\atop 1\leq|s|\leq \lfloor\wt\beta_X\rfloor} \frac{G_{[k]}^*{}^{(\mathbf 0,s)}(z,x^*)}{s!}(x-x^*)^s\big)^l\\
        &+\m O(2^{\frac{jD_Y}{2}}((\varepsilon_j^x)^{\beta_X\wedge 1}2^{j})^{\lfloor\wt \beta_X\rfloor+1})+\m O\Big(2^{\frac{jD_Y}{2}} (\varepsilon_j^x)^{\wt \beta_X} 2^j\Big),
   \end{aligned}
 \end{equation*}
 \begin{equation*}
 \begin{aligned}
     \Big| v_{[k]}^*(z,x)-\sum_{s\in \mb N_{0}^{D_X}\atop 0\leq|s|\leq \lfloor \alpha_X\rfloor} \frac{v_{[k]}^*{}^{(\mathbf 0,s)}(z,x^*)}{s!}(x-x^*)^s\Big|\lesssim  (\varepsilon_j^x)^{\alpha_X} 
 \end{aligned}
 \end{equation*}
 and recall
 \begin{equation*}
 \begin{aligned}
    v_{[k],x^*}^\dagger(z,x)=\sum_{s\in \mb N_{0}^{D_X}\atop 0\leq|s|\leq \lfloor \alpha_X\rfloor} \frac{v_{[k]}^\dagger{}^{(\mathbf 0,s)}(z,x^*)}{s!}(x-x^*)^s.
 \end{aligned}
 \end{equation*}
 Combined with the fact that  $ \big|v^*_\sk(z,x)-v^\dagger_\sk(z,x)\big|\lesssim  2^{-j\alpha_Y}+\log n\cdot(\varepsilon_j^x)^{\alpha_X}$, and 
 \begin{equation*}
 \begin{aligned}
      &   (2^{-j\alpha_Y}+\log n\cdot(\varepsilon_j^x)^{\alpha_X}) ((\varepsilon_j^x)^{\beta_X\wedge 1}2^{j})^{\lfloor\wt \beta_X\rfloor+1}\\
      &\lesssim \log n\cdot(\varepsilon_j^x)^{\alpha_X}+2^{-j\alpha_Y}(\varepsilon_j^x2^{j})^{\alpha_X}+2^{-j\alpha_Y}(\varepsilon_j^x)^{\alpha_X+\frac{\alpha_X}{\alpha_Y}}2^{j}\\
       &\lesssim \log n\cdot(\varepsilon_j^x)^{\alpha_X},
 \end{aligned}
 \end{equation*}
  \begin{equation*}
 \begin{aligned}
      &   (2^{-j\alpha_Y}+\log n\cdot(\varepsilon_j^x)^{\alpha_X}) (\varepsilon_j^x)^{\wt \beta_X} 2^j\\
      &\lesssim  (\varepsilon_j^x)^{\wt \beta_X} (\varepsilon_j^x)^{-\frac{\alpha_X}{\alpha_Y}}\\
      &=(\varepsilon_j^x)^{\alpha_X}.
 \end{aligned}
 \end{equation*}
  We can get 
 \begin{equation*}
     \begin{aligned}
  (E_B)&=\sum_{k\in [K^*]\atop \|x^*-x_k^*\|\leq \tau_2+2\varepsilon_j^x} \int_{\mb B_{\mb R^{d_Y}}(\mathbf{0},\tau_1)} 2^{\frac{j(d_Y-D_Y)}{2}}\psi(G_{[k]}^*(z,x)) \big(v_\sk^*(z,x)-v_{\sk,x^*}^\dagger(z,x) \big)\,\dd z\\
 &=\sum_{k\in [{K^*}]\atop \|x^*-x_k^*\|\leq \tau_2+2\varepsilon_j^x}  \int_{\{z\in \mb B_{\mb R^{d_Y}}(\mathbf{0},\tau_1):\,\|G^*_\sk(z,x^*)-y_{\psi^*}\|\leq C\,2^{-j}\}}  2^{\frac{j(d_Y-D_Y)}{2}}\psi(G_{[k]}^*(z,x)) \big(v_\sk^*(z,x)-v_{\sk,x^*}^\dagger(z,x) \big)\,\dd z\\
 &=\sum_{k\in [{K^*}]\atop \|x^*-x_k^*\|\leq \tau_2+2\varepsilon_j^x}  \int_{\{z\in \mb B_{\mb R^{d_Y}}(\mathbf{0},\tau_1):\,\|G^*_\sk(z,x^*)-y_{\psi^*}\|\leq C\,2^{-j}\}}  2^{\frac{j(d_Y-D_Y)}{2}}\Big(\psi^*(G_{[k]}^*(z,x^*))\\
 & +\sum_{l\in \mb N_{0}^{D_Y}\atop  1\leq |l|\leq \lfloor\wt \beta_X\rfloor} \frac{\psi^*{}^{(l)}(G_{[k]}^*(z,x^*))}{l!}\big(\sum_{s\in \mb N_{0}^{D_X}\atop 1\leq|s|\leq \lfloor\wt\beta_X\rfloor} \frac{G_{[k]}^*{}^{(\mathbf 0,s)}(z,x^*)}{s!}(x-x^*)^s\big)^l \Big)\big(v_\sk^*(z,x)-v_{\sk,x^*}^\dagger(z,x) \big)\,\dd z\\
 &\qquad\qquad+\m O (\log n\cdot 2^{\frac{-jd_Y}{2}}(\varepsilon_j^x)^{\alpha_X})\\
& =\sum_{k\in [{K^*}]\atop \|x^*-x_k^*\|\leq \tau_2+2\varepsilon_j^x}  \int_{\{z\in \mb B_{\mb R^{d_Y}}(\mathbf{0},\tau_1):\,\|G^*_\sk(z,x^*)-y_{\psi^*}\|\leq C\,2^{-j}\}}  2^{\frac{j(d_Y-D_Y)}{2}}\Big(\psi^*(G_{[k]}^*(z,x^*))\\
 & +\sum_{l\in \mb N_{0}^{D_Y}\atop  1\leq |l|\leq \lfloor\wt \beta_X\rfloor} \frac{\psi^*{}^{(l)}(G_{[k]}^*(z,x^*))}{l!}\big(\sum_{s\in \mb N_{0}^{D_X}\atop 1\leq|s|\leq \lfloor\wt\beta_X\rfloor} \frac{G_{[k]}^*{}^{(\mathbf 0,s)}(z,x^*)}{s!}(x-x^*)^s\big)^l \Big)\\
 &\qquad\cdot\big(\sum_{s\in \mb N_{0}^{D_X}\atop 0\leq|s|\leq \lfloor\alpha_X\rfloor} \frac{v_{[k]}^*{}^{(\mathbf 0,s)}(z,x^*)}{s!}(x-x^*)^s-\sum_{s\in \mb N_{0}^{D_X}\atop 0\leq|s|\leq \lfloor\alpha_X\rfloor} \frac{v_{[k]}^\dagger{}^{(\mathbf 0,s)}(z,x^*)}{s!}(x-x^*)^s\big)\,\dd z+\m O(\log n\cdot 2^{-\frac{d_Yj}{2}}(\varepsilon_j^x)^{\alpha_X}).
     \end{aligned}
 \end{equation*}
Notice that we can write 
 \begin{equation*}
     \begin{aligned}
       &  \sum_{k\in [{K^*}]\atop \|x^*-x_k^*\|\leq \tau_2+2\varepsilon_j^x}  \int_{\{z\in \mb B_{\mb R^{d_Y}}(\mathbf{0},\tau_1):\,\|G^*_\sk(z,x^*)-y_{\psi^*}\|\leq C\,2^{-j}\}}  2^{\frac{j(d_Y-D_Y)}{2}}\Big(\psi^*(G_{[k]}^*(z,x^*))\\
 & +\sum_{l\in \mb N_{0}^{D_Y}\atop  1\leq |l|\leq \lfloor\wt \beta_X\rfloor} \frac{\psi^*{}^{(l)}(G_{[k]}^*(z,x^*))}{l!}\big(\sum_{s\in \mb N_{0}^{D_X}\atop 1\leq|s|\leq \lfloor\wt\beta_X\rfloor} \frac{G_{[k]}^*{}^{(\mathbf 0,s)}(z,x^*)}{s!}(x-x^*)^s\big)^l \Big)\\
 &\qquad\big(\sum_{s\in \mb N_{0}^{D_X}\atop 0\leq|s|\leq \lfloor\alpha_X\rfloor} \frac{v_{[k]}^*{}^{(\mathbf 0,s)}(z,x^*)}{s!}(x-x^*)^s-\sum_{s\in \mb N_{0}^{D_X}\atop 0\leq|s|\leq \lfloor\alpha_X\rfloor} \frac{v_{[k]}^\dagger{}^{(\mathbf 0,s)}(z,x^*)}{s!}(x-x^*)^s\big)\,\dd z\\
  &=\sum_{s\in \mb N_0^{D_X}\atop 0\leq s\leq \lfloor\wt\beta_X\rfloor^2+\lfloor \alpha_X\rfloor} a'_{\psi^*,x^*,s}(x-x^*)^s,
     \end{aligned}
 \end{equation*}
 where $|a'_{\psi^*,x^*,s}|\lesssim n$.    So by combining all pieces, we have for any $\psi^*\in \Psi_j^{D_Y}$ and  $x\in \mb B_{\m M_X}(x^*,2\varepsilon_j^x)$,
 \begin{equation*}
     \begin{aligned}
                &\sum_{k=1}^{K^*} \int_{\mb B_{\mb R^{d_Y}}(\mathbf{0},\tau_1)} 2^{\frac{j(d_Y-D_Y)}{2}}\psi^*(G_{[k]}^*(z,x)) v_\sk^*(z,x)\,\dd z\\
                &=\sum_{k=1}^{K^*} \int_{\mb B_{\mb R^{d_Y}}(\mathbf{0},\tau_1)} 2^{\frac{j(d_Y-D_Y)}{2}}\psi^*(G_{[k],x^*}^\dagger(z,x)) v_{\sk,x^*}^\dagger(z,x) \,\dd z+\sum_{s\in \mb N_0^{D_X}\atop 0\leq s\leq \lfloor\wt\beta_X\rfloor^2+\lfloor \alpha_X\rfloor} a^{*}_{\psi^*,x^*,s}(x-x^*)^s\\
                &+\m O( \log n\cdot 2^{-\frac{jd_Y}{2}}(\varepsilon_j^x)^{\alpha_X})),
             \end{aligned}
 \end{equation*}
    where $a^*_{\psi^*,x^*,s}=a_{\psi^*,x^*,s}+a'_{\psi^*,x^*,s}$. This completes the proof of the first statement.

    \medskip
  \noindent   For the second statement,  fix arbitrary $x^*\in \mb N_{\varepsilon_j^x}^x$, $\psi\in \Psi_j^{D_Y}\setminus \Psi_j^{D_Y}(x^*)$,   $x\in \mb B_{\m M_X}(x^*,\varepsilon_j^x)$, $x'\in \mb B_{\mb N^x_{\varepsilon_j^x}}(x,2\varepsilon_j^x)$, $k\in [K^*]$, and $z\in \mb B_{\mb R^{d_Y}}(\mathbf{0},\tau_1)$. There exists $z^*\in \mb N_{c2^{-j}}^z$ so that $\|z-z^*\|\leq c2^{-j}$ and when $c$ is small enough, it holds that 
    \begin{equation*}
        \|G_\sk^*(z,x)-G_\sk^*(z^*,x^*)\|< L\, c2^{-j}+L (\varepsilon_j^x)^{\beta_X\wedge 1}\leq  \frac{C}{2}\,2^{-j}.
    \end{equation*}
    Since for any $l\in \mb N_{0}^{D_Y}$ with $ |l|\leq \lfloor\wt \beta_X\rfloor$ 
    \begin{equation*}
         {\rm supp}(\psi^{*(l)})\cap \mb B_{\mb R^{D_Y}}(G_{\sk}^*(z^*,x^*),C2^{-j})=\emptyset,
    \end{equation*}
    we have $\psi^{*(l)}(G_\sk^*(z,x))=0$. Moreover, since  $\|x-x'\|\leq 2\varepsilon_j^x$, when $\|x'-x_k^*\|\leq \tau_2+2\varepsilon_j^x$ and $C$ is sufficiently large, we have
    \begin{equation*}
       \|G^{\dagger}_{\sk,x'}(z,x)-G_\sk^*(z,x)\|\leq C_1 \,(2^{-j\beta_Y}+ (\log n)\cdot(\varepsilon_j^x)^{\beta_X}\big)<\frac{C}{2}\,2^{-j}
    \end{equation*}
and 
\begin{equation*}
  \|G^{\dagger}_{\sk,x'}(z,x)-G_\sk^*(z^*,x^*)\|\leq    \|G_\sk^*(z,x)-G_\sk^*(z^*,x^*)\|+  \|G^{\dagger}_{\sk,x'}(z,x)-G_\sk^*(z,x)\|< C\,2^{-j},
    \end{equation*}
 and thus  $\psi^{*(l)}(G^{\dagger}_{\sk,x'}(z,x))=0$. Furthermore, since $\|x'-x^*\|\leq 3\varepsilon_j^x$, we have, when $C$  is sufficiently large,
 \begin{equation*}
       \|G^{\dagger}_{\sk,x'}(z,x')-G_\sk^*(z^*,x^*)\|\leq   \|G^{\dagger}_{\sk,x'}(z,x')-G_\sk^*(z,x')\|+\|G_\sk^*(z,x')-G_\sk^*(z^*,x^*)\|< C\,2^{-j},
 \end{equation*}
 and hence $\psi^{*(l)}(G^{\dagger}_{\sk,x'}(z,x'))=0$ and $\psi^{*(l)}(G^*_{\sk}(z,x'))=0$.
 So we can get
    \begin{equation*}
        \begin{aligned}
            &\sum_{k=1}^{K^*} \int_{\mb B_{\mb R^{d_Y}}(\mathbf{0},\tau_1)} 2^{\frac{j(d_Y-D_Y)}{2}}\psi^*(G_{[k]}^*(z,x)) v_\sk^*(z,x)\,\dd z=0,\\
            &\sum_{k=1}^{K^*} \int_{\mb B_{\mb R^{d_Y}}(\mathbf{0},\tau_1)} 2^{\frac{j(d_Y-D_Y)}{2}}\psi^*(G^{\dagger}_{\sk,x'}(z,x)) v^{\dagger}_{\sk,x'}(z,x)\,\dd z\\
            &=\sum_{k\in [K^*]\atop \|x'-x_k^*\|\leq \tau_2+2\varepsilon_j^x}  \int_{\mb B_{\mb R^{d_Y}}(\mathbf{0},\tau_1)} 2^{\frac{j(d_Y-D_Y)}{2}}\psi^*(G^{\dagger}_{\sk,x'}(z,x)) v^{\dagger}_{\sk,x'}(z,x)\,\dd z=0,
        \end{aligned}
    \end{equation*}
and
\begin{equation*}
    \begin{aligned}
    &\sum_{s\in \mb N_0^{D_X}\atop 0\leq s\leq \lfloor\wt\beta_X\rfloor^2+\lfloor \alpha_X\rfloor} a^*_{\psi^*,x',s}(x-x')^s\\
          &= \sum_{k\in [K^*]\atop \|x'-x_k^*\|\leq \tau_2+2\varepsilon_j^x} \int_{\{z\in \mb B_{\mb R^{d_Y}}(\mathbf{0},\tau_1):\,\|G^*_\sk(z,x')-y_{\psi^*}\|\leq C\,2^{-j}\}} 2^{\frac{j(d_Y-D_Y)}{2}} \Big(\psi^*(G_{[k]}^*(z,x'))-\psi^*(G_{[k],x'}^\dagger(z,x'))\\
        &+ \sum_{l\in \mb N_{0}^{D_Y}\atop  1\leq |l|\leq \lfloor\wt\beta_X\rfloor} \frac{\psi^*{}^{(l)}(G_{[k]}^*(z,x'))}{l!}\big(\sum_{s\in \mb N_{0}^{D_X}\atop 1\leq|s|\leq \lfloor\wt\beta_X\rfloor} \frac{G_{[k]}^*{}^{(\mathbf 0,s)}(z,x')}{s!}(x-x')^s\big)^l\\
        &- \frac{\psi^*{}^{(l)}(G_{[k],x'}^\dagger(z,x'))}{l!}\big(\sum_{s\in \mb N_{0}^{D_X}\atop 1\leq|s|\leq\lfloor\wt\beta_X\rfloor} \frac{G_{[k],x'}^\dagger{}^{(\mathbf 0,s)}(z,x')}{s!}(x-x')^s\big)^l\Big) \sum_{s\in \mb N_{0}^{D_X}\atop 0\leq|s|\leq \lceil\alpha_X\rfloor} \frac{v_{[k]}^\dagger{}^{(\mathbf 0,s)}(z,x')}{s!}(x-x')^s\,\dd z\\
          &+  \sum_{k\in [{K^*}]\atop \|x'-x_k^*\|\leq \tau_2+2\varepsilon_j^x}  \int_{\{z\in \mb B_{\mb R^{d_Y}}(\mathbf{0},\tau_1):\,\|G^*_\sk(z,x')-y_{\psi^*}\|\leq C\,2^{-j}\}}  2^{\frac{j(d_Y-D_Y)}{2}}\Big(\psi^*(G_{[k]}^*(z,x'))\\
 & +\sum_{l\in \mb N_{0}^{D_Y}\atop  1\leq |l|\leq \lfloor\wt \beta_X\rfloor} \frac{\psi^*{}^{(l)}(G_{[k]}^*(z,x'))}{l!}\big(\sum_{s\in \mb N_{0}^{D_X}\atop 1\leq|s|\leq \lfloor\wt\beta_X\rfloor} \frac{G_{[k]}^*{}^{(\mathbf 0,s)}(z,x')}{s!}(x-x')^s\big)^l \Big)\\
 &\qquad\big(\sum_{s\in \mb N_{0}^{D_X}\atop 0\leq|s|\leq \lfloor\alpha_X\rfloor} \frac{v_{[k]}^*{}^{(\mathbf 0,s)}(z,x')}{s!}(x-x')^s-\sum_{s\in \mb N_{0}^{D_X}\atop 0\leq|s|\leq \lfloor\alpha_X\rfloor} \frac{v_{[k]}^\dagger{}^{(\mathbf 0,s)}(z,x')}{s!}(x-x')^s\big)\,\dd z\\
&=0.
    \end{aligned}
\end{equation*}
    The proof is now complete.


\subsection{Proof of Lemma~\ref{lemma: smoothG}}\label{proof:{lemma: smoothG}}
 
Consider any $G=\sum_{j_1=0}^{J_1}\sum_{j_2=0}^{J_2}\sum_{\psi_1\in  \Psi_{j_1}^{d_Y} }\sum_{\psi_2\in {\Psi}_{j_2}^{D_X}}g_{\psi_1\psi_2}\psi_1(z)\psi_2(x)\in \wt{\m G}$,  then since $\beta_Y\geq 2$, we have
    \begin{equation*}
\begin{aligned}
       \| J_{G(\cdot,x)}(z)\|_{F}&=\left\|\sum_{j_1=0}^{J_1}\sum_{j_2=0}^{J_2}\sum_{\psi_1\in  \Psi_{j_1}^{d_Y} }\sum_{\psi_2\in {\Psi}_{j_2}^{D_X}}g_{\psi_1\psi_2}J_{\psi_1}(z)\psi_2(x)\right\|_{F}\\
       &\lesssim \sum_{j_1=0}^{J_1}\sum_{j_2=0}^{J_2}  2^{-((j_1\beta_Y)\vee (j_2\beta_X))} 2^{j_1}\\
       &\leq \sum_{j_1=0}^{J_1}\sum_{j_2=0}^{\lfloor\frac{j_1\beta_Y}{\beta_X}\rfloor}  2^{-j_1(\beta_Y-1)}+\sum_{j_1=0}^{J_1}\sum_{j_2=\lfloor\frac{j_1\beta_Y}{\beta_X}\rfloor+1}^{J_2}  2^{-j_2(\beta_X-\frac{\beta_X}{\beta_Y})}\\
    & =\m O(1).
\end{aligned}
\end{equation*}
For the second statement, define set $\m A_0=[1,\infty)$, and for any $j\in [J_1]$, define $\m A_j=[2^{-j},2^{-(j-1)})$, and $\m A_{J_1+1}=(0,2^{-J_1})$. Then $\cup_{j=0}^{J_1+1} \m A_j=(0,\infty)$.  If $\|z-z'\|\in \m A_0$, we have 
\begin{equation*}
      \| J_{G(\cdot,x)}(z)-J_{G(\cdot,x)}(z')\|_{F}\leq   \| J_{G(\cdot,x)}(z)\|_{F}+\|J_{G(\cdot,x)}(z')\|_{F}\leq L_2\leq L_2\|z-z'\|^{\beta-1}.
\end{equation*}
If $\|z-z'\|\in \m A_{j}$ with $j\in [J_1]$, we have 
\begin{equation*}
    \begin{aligned}
            \| J_{G(\cdot,x)}(z)-J_{G(\cdot,x)}(z')\|_{F}&=\Big\|\sum_{j_1=0}^{J_1}\sum_{j_2=0}^{J_2}\sum_{\psi_1\in  \Psi_{j_1}^{d_Y} }\sum_{\psi_2\in {\Psi}_{j_2}^{D_X}}g_{\psi_1\psi_2}(J_{\psi_1}(z)-J_{\psi_1}(z'))\psi_2(x)\Big\|_{F}\\
            &\lesssim \sum_{j_1=0}^{j}\sum_{j_2=0}^{J_2} 2^{-((j_1\beta_Y)\vee (j_2\beta_X))}   2^{2j_1}\|z-z'\| +\sum_{j_1=j+1}^{J}\sum_{j_2=0}^{J_2} 2^{-((j_1\beta_Y)\vee (j_2\beta_X))}   2^{j_1}\\
            &\lesssim \sum_{j_1=0}^{j}\sum_{j_2=0}^{J_2} \Big(2^{-((j_1\beta_Y)\vee (j_2\beta_X))}   2^{2j_1}2^{-j(2-\beta)}\|z-z'\|^{\beta-1} \Big)+2^{-j(\beta_Y-1)}\cdot j\\
            &\lesssim \|z-z'\|^{\beta-1} \cdot 2^{-j(2-\beta)} \cdot  \sum_{j_1=0}^j (1+j_1) 2^{-(\beta_Y-2)j_1}  +2^{-j(\beta-1)}\cdot 2^{-j(\beta_Y-\beta)} \cdot j\\
            &\lesssim \|z-z'\|^{\beta-1},
         \end{aligned}
\end{equation*}
where the last inequality uses $\beta<2\leq \beta_Y$. Similarly, if $\|z-z'\|\in \m A_{J_1+1}$, then 
\begin{equation*}
    \begin{aligned}
            \| J_{G(\cdot,x)}(z)-J_{G(\cdot,x)}(z')\|_{F}&=\Big\|\sum_{j_1=0}^{J_1}\sum_{j_2=0}^{J_2}\sum_{\psi_1\in  \Psi_{j_1}^{d_Y} }\sum_{\psi_2\in {\Psi}_{j_2}^{D_X}}g_{\psi_1\psi_2}(J_{\psi_1}(z)-J_{\psi_1}(z'))\psi_2(x)\Big\|_{F}\\
            &\lesssim \sum_{j_1=0}^{J_1}\sum_{j_2=0}^{J_2} 2^{-((j_1\beta_Y)\vee (j_2\beta_X))}   2^{2j_1}\|z-z'\|\\
            &\lesssim \sum_{j_1=0}^{J_1}\sum_{j_2=0}^{J_2} 2^{-((j_1\beta_Y)\vee (j_2\beta_X))}   2^{2j_1}2^{-J_1(2-\beta)}\|z-z'\|^{\beta-1} \\
            &\lesssim \|z-z'\|^{\beta-1}.
         \end{aligned}
\end{equation*}
The proof is complete.

\subsection{Proof of Lemma~\ref{lemma:coveringG}}\label{proof:lemma:coveringG}
 
    Consider $$G(z,x)=\sum_{j_1=0}^{J_1}\sum_{j_2=0}^{J_2}\sum_{\psi_1\in  \Psi_{j_1}^{d_Y} }\sum_{\psi_2\in {\Psi}_{j_2}^{D_X}}g_{\psi_1\psi_2}\psi_1(z)\psi_2(x)$$ and $$G'(z,x)=\sum_{j_1=0}^{J_1}\sum_{j_2=0}^{J_2}\sum_{\psi_1\in  \Psi_{j_1}^{d_Y} }\sum_{\psi_2\in {\Psi}_{j_2}^{D_X}}g_{\psi_1\psi_2}'\psi_1(z)\psi_2(x).$$ 
     Then 
    \begin{equation*}
        \begin{aligned}
           & \underset{x\in \mb R^{D_X}\atop z\in \mb R^{d_Y}}{\sup}\|G(z,x)-G'(z,x)\|\\
            &=\underset{x\in \mb R^{D_X}\atop z\in \mb R^{d_Y}}{\sup}\Big\|\sum_{j_1=0}^{J_1}\sum_{j_2=0}^{J_2}\sum_{\psi_1\in  \Psi_{j_1}^{d_Y} }\sum_{\psi_2\in {\Psi}_{j_2}^{D_X}}\big(g_{\psi_1\psi_2}-g_{\psi_1\psi_2}'\big)\psi_1(z)\psi_2(x)\Big\|\\
            &\leq \sum_{j_1=0}^{J_1}\sum_{j_2=0}^{J_2}\underset{\psi_1\in  \Psi_{j_1}^{d_Y} ,\psi_2\in {\Psi}_{j_2}^{D_X}}{\max}\big\|g_{\psi_1\psi_2}-g_{\psi_1\psi_2}'\big\| \cdot\underset{x\in \mb R^{D_X}\atop z\in \mb R^{d_Y}}{\sup} \sum_{\psi_1\in  \Psi_{j_1}^{d_Y} }\sum_{\psi_2\in {\Psi}_{j_2}^{D_X}}|\psi_1(z)\psi_2(x)|\\
            &\leq C\, \sum_{j_1=0}^{J_1}\sum_{j_2=0}^{J_2}\underset{\psi_1\in  \Psi_{j_1}^{d_Y} ,\psi_2\in {\Psi}_{j_2}^{D_X}}{\max}\big\|g_{\psi_1\psi_2}-g_{\psi_1\psi_2}'\big\| \cdot 2^{\frac{d_Yj_1+D_Xj_2}{2}}.
        \end{aligned}
    \end{equation*}
When $\frac{d_Y}{\beta_Y}+\frac{d_X}{\beta_X}\leq 2\gamma_1\leq 2$, there exists a constant $C_1$ so that 
    \begin{equation*}
        \begin{aligned}
            &\sum_{j_1=0}^{J_1}\sum_{j_2=0}^{J_2} 2^{\frac{d_Yj_1+d_Xj_2}{4\gamma_1}-\frac{j_1\beta_Y\vee j_2\beta_X}{2}}\\
            &\leq \frac{2^{d_X/4}}{2^{d_X/4}-1}\sum_{j_1=0}^{J_1} 2^{\frac{d_Xj_1\beta_Y}{4\beta_X\gamma_1}+\frac{d_Yj_1}{4\gamma_1}-\frac{j_1\beta_Y}{2}}+\frac{2^{d_Y/4}}{2^{d_Y/4}-1}\sum_{j_2=0}^{J_2} 2^{\frac{d_Yj_2\beta_X}{4\beta_Y\gamma_1}+\frac{d_Xj_2}{4\gamma_1}-\frac{j_2\beta_X}{2}}\\
        &\leq C_1( J_1+J_2).
                   \end{aligned}
    \end{equation*}
So if for any $j_1\in [J_1]$, $j_2\in [J_2]$, $\psi_1\in  \Psi_{j_1}^{d_Y}$, and $\psi_2\in {\Psi}_{j_2}^{D_X}$,
    \begin{equation*}
 \big\|g_{\psi_1\psi_2}-g_{\psi_1\psi_2}'\big\|\leq \frac{\varepsilon^{\frac{1}{\gamma_1}}}{CC_1(J_1+J_2)} 2^{\frac{d_Yj_1+d_Xj_2}{4\gamma_1}-\frac{j_1\beta_Y\vee j_2\beta_X}{2}} 2^{-\frac{d_Yj_1+D_Xj_2}{2}},
    \end{equation*}
    then 
    \begin{equation*}
        \underset{x\in \mb R^{D_X}\atop z\in \mb R^{d_Y}}{\sup}\|G(z,x)-G'(z,x)\|^{\gamma_1}\leq \varepsilon.
    \end{equation*}
    Therefore, we can get
    \begin{equation*}
    \begin{aligned}
      \mathbf N(\m G,d_{\infty}^{\gamma_1},\varepsilon)&\leq \prod_{j_1=0}^{J_1}\prod_{j_2=0}^{J_2}\prod_{\psi_1\in  \Psi_{j_1}^{d_Y}}\prod_{\psi_2\in {\Psi}_{j_2}^{D_X}} \m N([-L_1\, 2^{-\frac{d_Yj_1+D_Xj_2}{2}-((j_1\beta_Y)\vee (j_2\beta_X))},
        L_1\, 2^{-\frac{d_Yj_1+D_Xj_2}{2}-((j_1\beta_Y)\vee (j_2\beta_X))}]^{D_Y}\\
          &\qquad\qquad,\frac{\varepsilon^{\frac{1}{\gamma_1}}}{CC_1(J_1+J_2)} 2^{\frac{d_Yj_1+d_Xj_2}{4\gamma_1}-\frac{j_1\beta_Y\vee j_2\beta_X}{2}} 2^{-\frac{d_Yj_1+D_Xj_2}{2}},\|\cdot\|)\\
          &\leq \prod_{j_1=0}^{J_1}\prod_{j_2=0}^{J_2}\prod_{\psi_1\in  \Psi_{j_1}^{d_Y}}\prod_{\psi_2\in {\Psi}_{j_2}^{D_X}} \lceil\Big(\frac{6\sqrt{D_Y}L_1CC_1\, (J_1+J_2)2^{-\frac{d_Yj_1+d_Xj_2}{4\gamma_1}-\frac{(j_1\beta_Y)\vee (j_2\beta_X)}{2}}}{\varepsilon^{\frac{1}{\gamma_1}}}\Big)^{D_Y}\rceil \vee 1\\
           &\leq \prod_{j_1=0}^{J_1}\prod_{j_2=0}^{J_2}\prod_{\psi_1\in  \Psi_{j_1}^{d_Y}}\prod_{\psi_2\in {\Psi}_{j_2}^{D_X}} \Big(\frac{12\sqrt{D_Y}L_1CC_1\, (J_1+J_2)2^{-\frac{d_Yj_1+d_Xj_2}{4\gamma_1}-\frac{(j_1\beta_Y)\vee (j_2\beta_X)}{2}}}{\varepsilon^{\frac{1}{\gamma_1}}}\Big)^{D_Y}\vee 1.
         \end{aligned}
    \end{equation*}
   Moreover, for any $j\in [J_2]$, let $\m N_{2^{-j}}^x$ be the largest $2^{-j}$-packing set of $\m M_x$, then $|\m N_{2^{-j}}^x|\lesssim 2^{jd_X}$, and 
   \begin{equation*}
       |{\Psi}_{j}^{D_X}|\leq \sum_{x\in \m N_{2^{-j}}^x} \Big|\{\psi\in \ov\Psi_j^{D_X}:\, {\rm supp}(\psi)\cap \mb B_{\mb R^{D_X}}(x,2^{-j})\neq \emptyset\}\Big|\lesssim 2^{jd_X}.
   \end{equation*}
Hence there exists a constant $C_2$ so that for any $\gamma_1$ satisfying  $\frac{d_Y}{\beta_Y}+\frac{d_X}{\beta_X}\leq 2\gamma_1\leq 2$,  it holds that 
 \begin{equation*}
        \log   \mathbf N(\m G,d_{\infty}^{\gamma_1},\varepsilon)\leq C_2\sum_{j_1=0}^{J_1}\sum_{j_2=0}^{J_2} 2^{d_Yj_1+d_Xj_2}\log\left(\frac{C_2 (J_1+J_2)2^{-\frac{d_Yj_1+d_Xj_2}{4\gamma_1}-\frac{(j_1\beta_Y)\vee (j_2\beta_X)}{2}}}{\varepsilon^{\frac{1}{\gamma_1}}}\vee 1\right).
    \end{equation*}
    When $\frac{d_Y}{\beta_Y}+\frac{d_X}{\beta_X}>2\gamma_1$,  denote
    \begin{equation*}
        s_{j_1j_2}=\sqrt{\frac{2^{\frac{d_Yj_1+d_Xj_2}{2\gamma_1}-(j_1\beta_Y\vee j_2\beta_X)}}{2^{\frac{d_YJ_1+d_XJ_2}{2\gamma_1}-(J_1\beta_Y\vee J_2\beta_X)}}}.
    \end{equation*}
There exists constants $C_2,C_3$ so that for any $\gamma_1\in (0,\frac{d_Y}{2\beta_Y}+\frac{d_X}{2\beta_X})$,
    \begin{equation*}
        \begin{aligned}
            &S:=\sum_{j_1=0}^{J_1}\sum_{j_2=0}^{J_2} s_{j_1j_2}=\sqrt{\frac{1}{2^{\frac{d_YJ_1+d_XJ_2}{2\gamma_1}-(J_1\beta_Y\vee J_2\beta_X)}}}\cdot \sum_{j_1=0}^{J_1}\sum_{j_2=0}^{J_2} 2^{\frac{d_Yj_1+d_Xj_2}{4\gamma_1}-\frac{j_1\beta_Y\vee j_2\beta_X}{2}}\\
            &\leq\sqrt{\frac{1}{2^{\frac{d_YJ_1+d_XJ_2}{2\gamma_1}-(J_1\beta_Y\vee J_2\beta_X)}}}\cdot\Big( \frac{2^{d_X/4}}{2^{d_X/4}-1}\sum_{j_1=0}^{J_1} 2^{\frac{d_Xj_1\beta_Y}{4\beta_X\gamma_1}+\frac{d_Yj_1}{4\gamma_1}-\frac{j_1\beta_Y}{2}}+\frac{2^{d_Y/4}}{2^{d_Y/4}-1}\sum_{j_2=0}^{J_2} 2^{\frac{d_Yj_2\beta_X}{4\beta_Y\gamma_1}+\frac{d_Xj_2}{4\gamma_1}-\frac{j_2\beta_X}{2}}\Big)\\
           & \leq C_2\sqrt{\frac{1}{2^{\frac{d_YJ_1+d_XJ_2}{2\gamma_1}-(J_1\beta_Y\vee J_2\beta_X)}}}\cdot\Big( \frac{2^{\frac{(d_Y+d_X\beta_Y/\beta_X-2\beta_Y\gamma_1)(J_1+1)}{4\gamma_1}}-1}{2^{\frac{(d_Y+d_X\beta_Y/\beta_X-2\beta_Y\gamma_1)}{4\gamma_1}}-1}+\frac{2^{\frac{(d_X+d_Y\beta_X/\beta_Y-2\beta_X\gamma_1)(J_2+1)}{4\gamma_1}}-1}{2^{\frac{(d_X+d_Y\beta_X/\beta_Y-2\beta_X\gamma_1)}{4\gamma_1}}-1}\Big)\\
           &\leq C_3\Big((J_1\wedge \frac{2^{\frac{(d_Y+d_X\beta_Y/\beta_X-2\beta_Y\gamma_1)}{4\gamma_1}}}{2^{\frac{(d_Y+d_X\beta_Y/\beta_X-2\beta_Y\gamma_1)}{4\gamma_1}}-1}) +(J_2\wedge \frac{2^{\frac{(d_X+d_Y\beta_X/\beta_Y-2\beta_X\gamma_1)}{4\gamma_1}}}{2^{\frac{(d_X+d_Y\beta_X/\beta_Y-2\beta_X\gamma_1)}{4\gamma_1}}-1})\Big)\\
           &\leq C_3\Big( (J_1+J_2)  \wedge ( \frac{2^{\frac{(d_Y+d_X\beta_Y/\beta_X-2\beta_Y\gamma_1)}{4\gamma_1}}}{2^{\frac{(d_Y+d_X\beta_Y/\beta_X-2\beta_Y\gamma_1)}{4\gamma_1}}-1}+ \frac{2^{\frac{(d_X+d_Y\beta_X/\beta_Y-2\beta_X\gamma_1)}{4\gamma_1}}}{2^{\frac{(d_X+d_Y\beta_X/\beta_Y-2\beta_X\gamma_1)}{4\gamma_1}}-1})\Big).
                   \end{aligned}
    \end{equation*}
So if for any $j_1\in [J_1]$, $j_2\in [J_2]$, $\psi_1\in  \Psi_{j_1}^{d_Y}$, and $\psi_2\in {\Psi}_{j_2}^{D_X}$,
    \begin{equation*}
 \big\|g_{\psi_1\psi_2}-g_{\psi_1\psi_2}'\big\|\leq \frac{\varepsilon^{\frac{1}{\gamma_1}}s_{j_1j_2}}{CS} 2^{-\frac{d_Yj_1+D_Xj_2}{2}},
    \end{equation*}
    then 
    \begin{equation*}
        \underset{x\in \mb R^{D_X}\atop z\in \mb R^{d_Y}}{\sup}\|G(z,x)-G'(z,x)\|^{\gamma_1}\leq \varepsilon.
    \end{equation*}
     Therefore,  there exists a constant $C_4$ so that for any $\gamma_1\in (0,\frac{d_Y}{2\beta_Y}+\frac{d_X}{2\beta_X})$,
    \begin{equation*}
    \begin{aligned}
      & \log\mathbf N(\m G,d_{\infty}^{\gamma_1},\varepsilon)\leq \sum_{j_1=0}^{J_1}\sum_{j_2=0}^{J_2}\sum_{\psi_1\in  \Psi_{j_1}^{d_Y}}\sum_{\psi_2\in {\Psi}_{j_2}^{D_X}}\log \mathbf N\big( \\
       &\qquad\qquad\qquad [-L_1\, 2^{-\frac{d_Yj_1+D_Xj_2}{2}-((j_1\beta_Y)\vee (j_2\beta_X))},
        L_1\, 2^{-\frac{d_Yj_1+D_Xj_2}{2}-((j_1\beta_Y)\vee (j_2\beta_X))}]^{D_Y}
       , \frac{\varepsilon^{\frac{1}{\gamma_1}} s_{j_1j_2}}{CS} 2^{-\frac{d_Yj_1+D_Xj_2}{2}},\|\cdot\|\big)\\
          &\qquad\qquad\qquad\quad\leq C_4 \sum_{j_1=0}^{J_1}\sum_{j_2=0}^{J_2}  2^{d_Yj_1+d_Xj_2} \log\Big(\frac{C_4\, S2^{-((j_1\beta_Y)\vee (j_2\beta_X))}}{\varepsilon^{\frac{1}{\gamma_1}} s_{j_1j_2}}\vee 1\Big),
         \end{aligned}
    \end{equation*}
which completes the proof.

\section{Proof of Technical Details}
\subsection{Proof of Lemma~\ref{le:wavelet}}\label{proof:le:wavelet}
 
 Let $\zeta=(\lceil\alpha\rceil\vee \lceil\frac{d}{2}-\alpha\rceil)+1$ and let $\phi_{\mf M}$ and $\phi_{\mf F}$ be the Daubechies wavelet  and scaling function~\citep{daubechies1992ten,meyer1992wavelets} that are supported in a compact set $[-C,C]$,  have derivatives  up to order $\zeta$ and
\begin{equation*}
    \int_{\mathbb{R}} x^l \psi_{\mf M}(x) \,\dd x=0 \quad for \quad l=0, \ldots, \zeta.
\end{equation*}
Then by Proposition 1.51 of~\cite{Triebel2006},
\begin{equation*}
 \left\{
\begin{array}{ll}
   \psi_{\mf F} (x-k) & j=0, k\in \mb Z,  \\
    2^{(j-1)/2} \psi_{\mf M}(2^{j-1}x-k), &  j\in \mb N_{+}, k\in \mb Z,
\end{array}
\right.
\end{equation*}
is an orthonormal basis of $\m L^2(\mb R)$. Furthermore, by Proposition 1.53 of ~\cite{Triebel2006}, to obtain a basis of $\m L^2(\mb R^d)$ for an integer $d>1$, set 
\begin{equation*}
     \mf G = \{\mf F,\,\mf M\}^d\setminus \{(\mf F,\ldots,\mf F)\}.
\end{equation*}
Then for any multi-index $k\in \mb Z^d$, the level zero basis $\phi_k^{[d]}$ is obtained by translating the $d$-fold tensor product $\phi_{\mf F}^{\otimes d}$ by $k$ as $\phi_{k}^{[d]}(x) = \prod_{i=1}^d \phi_{\mf F}(x_i-k_i)$ for $x=(x_1,\ldots,x_d)\in\mb R^d$, and for any $j\geq 1$, the level $j$ basis $\big\{\psi_{ljk}^{[d]}:\, l\in[2^d-1]\big\}$ with translation $k$ is any ordering of the following $2^d-1$ functions,
\begin{align*}
    \phi_{gjk}^{[d]}(x)=2^{\frac{d(j-1)}{2}} \,\prod_{i=1}^d \phi_{g_i}\big(2^{j-1}x_i - k_i\big), \quad \forall g\in \mf G.
\end{align*}
This gives the orthornormal basis  
\begin{equation*}
\left\{
\begin{array}{ll}
   \phi_k^{[d]}(x), & j=0,l=0, k\in \mb Z^d,  \\
    \psi_{ljk}^{[d]}(x), &  j\in \mb N_{+},l\in [2^d-1], k\in \mb Z^d.
\end{array}
\right.
\end{equation*}
Denote $\ov\Psi^{d}_0=\{\phi_{k}^{[d]}(\cdot):\, k\in \mb Z^d\}$ as the set of level zero basis and $\ov\Psi^{d}_j=\{\psi_{ljk}^{[d]}(\cdot):\, l\in [2^d-1], k\in \mb Z^d\}$ as the set of level $j$ basis for $j\in \mb N_{+}$. Then use the fact that for any $s\in \mb N_0^d$ with $|s|\leq \alpha$,
\begin{equation*}
     \phi_{gjk}^{[d]}{}^{(s)}(x)=2^{\frac{d(j-1)}{2}} \,\prod_{i=1}^d 2^{(j-1)s_i}\phi_{g_i}^{(s_i)}\big(2^{j-1}x_i - k_i\big)\leq C_R 2^{\frac{dj}{2}+j|s|},
\end{equation*}
we can get the regularity condition. Moreover, by the compactness of the supports and smoothness of  $\phi_{\mf M}$  and $\phi_{\mf F}$, we have
\begin{equation*}
    {\rm supp}(\psi_{ljk}^{[d]}{}^{(s)})\subset \prod_{i=1}^d[\frac{-C+k_i}{2^{j-1}},\frac{C+k_i}{2^{j-1}}]=I_{\psi_{ljk}^{[d]}}.
\end{equation*}
\begin{equation*}
    {\rm supp}( \phi_k^{[d]}{}^{(s)})\subset \prod_{i=1}^d[-C+k_i,C+k_i]=I_{\phi_k^{[d]}}.
\end{equation*}
So for any $x\in \mb R^d$, $j\in \mb N$, and $l\in [2^d-1]$, there are only constant number of $k$ so that  $\psi_{ljk}^{[d]}(x)\neq 0$ ($j>0$) or $\phi_k^{[d]}(x)\neq 0$ ($j=0$). Hence $\sup_{x\in \mb R^d}\sum_{\psi\in \ov \Psi_j^d} \mathbf{1}(x\in I_{\psi})\leq C_L'$. Moreover, if $I_{\phi_k^{[d]}}\cap \mb B_{\mb R^d}(0,R)\neq \emptyset$, then $k\in [-C-R,C+R]^d$; if $I_{\psi_{ljk}^{[d]}}\cap \mb B_{\mb R^d}(0,R)\neq \emptyset$, then $k\in [2^{j-1}(-C-R),2^{j-1}(C+R)]^d$, so
$\big|\{\psi\in \ov\Psi_j^{d}\,: \,I_{\psi}\cap \mb B_{\mb R^d}(0,R)\neq \emptyset\}\big|\leq (2^d-1)(2^j(C+R)+1)^d\leq (2^d-1)(C+2)^d R^d 2^{jd}$; if $I_{\psi_{ljk}^{[d]}}\cap \mb B_{\mb R^d}(x,2^{-(j-1)})\neq \emptyset$, then for any $i\in [d]$, $k_i\in [2^{j-1}x_i-C,2^{j-1}x_i+C]$, which means $\big|\{\psi\in \ov\Psi_j^{d}\,: \,I_{\psi}\cap \mb B_{\mb R^d}(x,2^{-(j-1)})\neq \emptyset\}\big|\leq (2^d-1)(2C+1)^d$. For the third statement, since $f\in \m H_{r}^{\alpha_1}(\mb R^d)$, it holds for any $x,x_0\in \mb R^d$ that 
\begin{equation*}
   \Big| f(x)-\sum_{s\in \mb N_0^{d}\atop |s|<\alpha_1} \frac{  f^{(s)}(x_0)}{s!}(x-x_0)^{s}\Big|  \leq  r\|x-x_0\|^{\alpha_1}.
\end{equation*}
 Then for any $j\in \mb N$ and $\psi\in \ov\Psi_{j}^d$, we have 
 \begin{enumerate}
     \item If $j=0$, 
     \begin{equation*}
     \begin{aligned}
           &\int_{\mb R^{d}} f(x)\psi(x)\,\dd x=\int_{I_{\psi}} f(x)\psi(x)\,\dd x\leq \sqrt{\int_{I_{\psi}}\psi^2 (x)\,\dd x \int_{I_{\psi}} f^2(x)\,\dd x}\leq (2C)^{\frac{d}{2}} r.
     \end{aligned}
     \end{equation*}
 \item If $j>0$, then we have for any $l\in \mb N_0^d$ with $|l|<\alpha_1$, 
 \begin{equation*}
     \int_{\mb R^d} x^l \psi(x)\,\dd x=0
 \end{equation*}
 and thus for any $x_0\in I_{\psi}$, we have
 \begin{equation*}
     \begin{aligned}
        &\big|\int_{\mb R^d} f(x)\psi(x)\,\dd x\big|=\big|\int_{\mb R^d} (f(x)-f(x_0))\psi(x)\,\dd x\big|=\big|\int_{I_{\psi}} (f(x)-f(x_0))\psi(x)\,\dd x\big|\\
        &=\Big|\int_{\mb R^d} \sum_{s\in \mb N_0^{d}\atop 1\leq |s|<\alpha_1} \frac{  f^{(s)}(x_0)}{s!}(x-x_0)^{s}  \psi(x)\,\dd x+\int_{I_{\psi}} f(x)-\sum_{s\in \mb N_0^{d}\atop 1\leq |s|<\alpha_1} \frac{  f^{(s)}(x_0)}{s!}(x-x_0)^{s} \psi(x)\,\dd x\Big|\\
        &\leq \int_{I_{\psi}} \Big|f(x)-\sum_{s\in \mb N_0^{d}\atop 1\leq |s|<\alpha_1} \frac{  f^{(s)}(x_0)}{s!}(x-x_0)^{s} \Big|\cdot|\psi(x)|\,\dd x\\
        &\leq \int_{I_{\psi}}r \|x-x_0\|^{\alpha_1}\cdot|\psi(x)|\,\dd x\\
        &\leq r \cdot \underset{x\in I_{\psi}}{\sup}\|x-x_0\|^{\alpha_1}\cdot \underset{x\in I_{\psi}}{\sup}|\psi(x)| \cdot\int_{I_{\psi}} \,\dd x \\
        &\leq  r(\frac{4C}{2^j})^d (\frac{4C\sqrt{d}}{2^j})^{\alpha_1} C_R 2^{\frac{dj}{2}}\\
        &\lesssim 2^{-\frac{dj}{2}-j\alpha_1}.
     \end{aligned}
 \end{equation*}
 \end{enumerate}
 For the last statement.  When $j=0$, we have 
\begin{equation*}
   \Psi_0^{d}\subset \{\phi^{[d]}_k(x):\, k\in\mb Z \text{ and } k\in [-C_L-R',C_L+R']^d\}  
\end{equation*}
 Then we set 
 \begin{equation*}
     \ov{\ms I}_0=\{(\iota_1,0):\, \iota_1\in [0,1]^d \text{ and } (2\iota_1-1)\cdot(C_L+R')\in \mb{Z}^d\}
 \end{equation*}
and for any $\iota=(\iota_1,0)\in \ov{\ms I}_0$, we set 
\begin{equation*}
    \phi_{0\iota}(\cdot)=\phi^{[d]}_{(2\iota_1-1)\cdot(C_L+R')})(\cdot).
\end{equation*}
Let 
 \begin{equation*}
     {\ms I}_0=\{\iota\in  \ov{\ms I}_0\,:\,   \phi_{0\iota}(\cdot)\in  \Psi_0^{d}\},
 \end{equation*}
we have 
\begin{equation*}
      \Psi_0^{d}=\{\phi_{0\iota}(\cdot):\, \iota\in \ms I_0\subset [0,1]^{d+1}\},
\end{equation*}
and for any $\iota,\iota'\in  {\ms I}_0$ with $\iota\neq\iota'$, it holds that 
$\|\iota-\iota'\|\geq \frac{1}{2(C_L+R')}$. When $j>0$, we have 
\begin{equation*}
   \Psi_j^{d}\subset \{\psi^{[d]}_{ljk}(x):\, l\in [2^d-1],  k\in\mb Z \text{ and }k\in[-2^{j-1}C_L-R',2^{j-1}C_L+R']^d\}.
\end{equation*}
Then we set 
 \begin{equation*}
     \ov{\ms I}_j=\{(\iota_1,\iota_2):\, \iota_1\in [0,1]^d \text{ and } (2\iota_1-1)\cdot(2^{j-1}C_L+R')\in \mb{Z}^d, \iota_2\in [0,1]\text{ and }\iota_2(2^{d-1}-1)+1\in \mb Z\} , 
 \end{equation*}
 and for any $\iota=(\iota_1,\iota_2)\in \ov{\ms I}_j$, we set 
\begin{equation*}
    \phi_{j\iota}=\psi^{[d]}_{\iota_2(2^{d-1}-1)+1,j,(2\iota_1-1)\cdot(C_L+R')}(x).
\end{equation*}
Let 
 \begin{equation*}
     {\ms I}_j=\{\iota\in  \ov{\ms I}_j\,\,   \phi_{j\iota}\in  \Psi_j^{d}\},
 \end{equation*}
we have 
\begin{equation*}
      \Psi_j^{d}=\{\psi_{j\iota}(\cdot):\, \iota\in \ms I_j\subset [0,1]^{d+1}\},
\end{equation*}
and for any $\iota,\iota'\in  {\ms I}_j$ with $\iota\neq\iota'$, it holds that 
\begin{equation*}
    \|\iota-\iota'\|\geq \frac{1}{2^jC_L+2R'}\wedge \frac{1}{2^{d-1}-1}.
\end{equation*}
We can then get the desired result by combining all pieces.

\subsection{Proof of Lemma~\ref{le:poly}}\label{appendixE.2}
Without loss of generality, we may assume $\alpha_1\geq \alpha_2$.
  Given any $x\in \mb R^{d_1}$, and considering $f(x,\cdot)\in \m H^{\alpha_2}_r(\mb R^{d_2})$, it follows that for any $y_0,y\in \mb R^{d_2}$,
    \begin{equation*}
        \begin{aligned}
            &\Big| f(x,y)-\sum_{j_2\in \mb N_0^{d_2}\atop |j_2|<\alpha_2} \frac{  f^{(\mathbf{0},j_2)}(x,y_0)}{j_2!}(y-y_0)^{j_2}\Big| \\
            &=\left\{\begin{array}{ll}
               \bigg|\sum_{j_2\in \mb N_0^{d_2}\atop |j_2|=\lfloor\alpha_2\rfloor}\frac{\lfloor\alpha_2\rfloor}{j_2!} \int_0^1(1-t)^{ \lfloor\alpha_2\rfloor-1}  \big(  f^{(\mathbf 0,j_2)}(x,y_0+t(y-y_0))-  f^{(\mathbf 0,j_2)}(x,y_0)\big) \,\dd t \cdot (y-y_0)^{j_2}\bigg|,  & \, \alpha_2>1 \\
              |f(x,y)-f(x,y_0)|,   & \, \alpha_2\leq 1
            \end{array}
            \right.\\
            &=\m O(\|y-y_0\|^{\alpha_2}).
        \end{aligned}
    \end{equation*}
    Moreover, using $f\in   {\m H}^{\alpha_1,\alpha_2}_r(\mb R^{d_1},\mb R^{d_2})$, we have   for any $x,x_0\in \mb R^{d_1}$,
    \begin{equation*}
        \begin{aligned}
          &\Big|  \sum_{j_2\in \mb N_0^{d_2}\atop |j_2|<\alpha_2} \frac{  f^{(\mathbf{0},j_2)}(x,y_0)}{j_2!}(y-y_0)^{j_2}-\sum_{j_2\in \mb N_0^{d_2}\atop |j_2|<\alpha_2} \sum_{j_1\in \mb N_0^{d_1}\atop |j_1|+\frac{\alpha_1}{\alpha_2}|j_2|<\alpha_1}\frac{  f^{(j_1,j_2)}(x_0,y_0)}{j_1!j_2!}(x-x_0)^{j_1}(y-y_0)^{j_2}\Big|\\
          &=\Big|  \sum_{j_2\in \mb N_0^{d_2}\atop |j_2|<\alpha_2\cdot\frac{\alpha_1-1}{\alpha_1}}
          \sum_{j_1\in \mb N_0^{d_1}\atop |j_1|=\lfloor\alpha_1-\frac{\alpha_1}{\alpha_2}|j_2|\rfloor}\frac{\lfloor\alpha_1-\frac{\alpha_1}{\alpha_2}|j_2|\rfloor}{j_1!j_2!} \int_0^1(1-t)^{ \lfloor\alpha_1-\frac{\alpha_1}{\alpha_2}|j_2|\rfloor-1}  \big(  f^{(j_1,j_2)}(x_0+t(x-x_0),y_0)-  f^{(j_1,j_2)}(x_0,y_0)\big) \,\dd t \\
          &\qquad \cdot (x-x_0)^{j_1} (y-y_0)^{j_2}\Big|\\
          &+\Big|  \sum_{j_2\in \mb N_0^{d_2}\atop \alpha_2\cdot\frac{\alpha_1-1}{\alpha_1} \leq |j_2|<\alpha_2}
           \frac{1}{j_2!}   \big(  f^{(\mathbf{0},j_2)}(x,y_0)-  f^{(\mathbf{0},j_2)}(x_0,y_0)\big)(y-y_0)^{j_2}\Big|\\
          &=\m O\Big(  \sum_{j_2\in \mb N_0^{d_2}\atop |j_2|<\alpha_2} \sum_{j_1\in \mb N_0^{d_1}\atop |j_1|=\lfloor\alpha_1-\frac{\alpha_1}{\alpha_2}|j_2|\rfloor}\|x-x_0\|^{\alpha_1-\frac{\alpha_1}{\alpha_2}|j_2|}\|y-y_0\|^{|j_2|}\Big)=\m O(\|x-x_0\|^{\alpha_1}+\|y-y_0\|^{\alpha_2}),
             \end{aligned}
    \end{equation*}
    where the last inequality uses the Young's inequality for products.
    Therefore, we can get 
    \begin{equation*}
        \begin{aligned}
            & \Big| f(x,y)-\sum_{(j_1,j_2)\in \m J^{d_1,d_2}_{\alpha_1,\alpha_2}}  \frac{f^{(j_1,j_2)}(x_0,y_0)}{j_1!j_2!}(x-x_0)^{j_1}(y-y_0)^{j_1}\Big|\\
            &=\Big| f(x,y)-\sum_{j_2\in \mb N_0^{d_2}\atop |j_2|<\alpha_2} \sum_{j_1\in \mb N_0^{d_1}\atop |j_1|+\frac{\alpha_1}{\alpha_2}|j_2|<\alpha_1}\frac{f^{(j_1,j_2)}(x_0,y_0)}{j_2!}(x-x_0)^{j_1}(y-y_0)^{j_2}\Big|=\m O(\|x-x_0\|^{\alpha_1}+\|y-y_0\|^{\alpha_2}).
        \end{aligned}
    \end{equation*}

\subsection{Proof of Lemma~\ref{le:approwaveletsmooth}}\label{proof:le:approwaveletsmooth}
     For any $y\in \mb R^{d_2}$, the function $ f(\cdot,y)$ has the following wavelet expansion
\begin{equation*}
   f(\cdot,y)=\sum_{j=0}^{\infty}\sum_{\psi\in \ov\Psi_j^{d_1}}\psi(\cdot) f_{\psi}(y),\quad  f_{\psi}(y)=\int_{\mb R^{d_1}}  f(x,y)\psi(x)\,\dd x,
\end{equation*}
with $| f_{\psi}(y)|\leq C_W\,L\,2^{-\frac{d_1j}{2}-j\alpha_1}$ when $\psi\in \ov\Psi_j^{d_1}$. Then we have
\begin{equation*}
\begin{aligned}
       &\underset{x\in \mb R^{d_1}}{\sup}\Big|\sum_{j=J_1+1}^{\infty}\sum_{\psi\in \Psi_j}\psi(x)  f_{\psi}(y)\Big|\leq C_RC_L'C_W L\, \sum_{j=J_1+1}^{\infty} 2^{-j\alpha_1}\leq C_RC_L'C_W L\, 2^{-J_1\alpha_1}.
\end{aligned}
\end{equation*}
 Moreover, for any $j_1\in [J_1]$ and $\psi\in \ov\Psi_{j_1}^{d_1}$, it holds that
 \begin{equation*}
     \begin{aligned}
       &2^{\frac{d_1j_1}{2}}\int_{\mb R^{d_1}
       } |\psi(x)|\,\dd x\leq 2^{\frac{d_1j_1}{2}}\int_{I_{\psi}}\,\dd x\cdot \underset{x\in I_{\psi}}{\sup }|\psi(x)|\leq  (2C_L)^{d_1}C_R.
     \end{aligned}
 \end{equation*}
Furthermore, for any multi-index $\ell\in \mb N_0^{d_2}$ with $|\ell|<\alpha$, it holds that 
 \begin{equation*}
 \begin{aligned}
          2^{\frac{d_1j_1}{2}} f^{(\ell)}_{\psi}(y)&=   2^{\frac{d_1j_1}{2}}\left[\int_{\mb R^{d_1}}  f(x,\cdot)\psi(x)\,\dd x\right]^{(\ell)}(y)\\
          &= 2^{\frac{d_1j_1}{2}}\int_{\mb R^{d_1}} f^{(\mathbf{0}_{d_1},\ell)}(x,y)\psi(x)\,\dd x.
 \end{aligned}
 \end{equation*}
Therefore, there exists a constant $L_1=(2C_L)^{d_1}C_RL$ so that
 \begin{equation*}
     \begin{aligned}
         2^{\frac{d_1j_1}{2}} f_{\psi}(y)=2^{\frac{d_1j_1}{2}}\int_{\mb R^{d_1}}  f(x,y)\psi(x)\,\dd x\in \m H^{\alpha_2}_{L_1}(\mb R^{d_2}).
     \end{aligned}
 \end{equation*}
 For any $j_1\in \mb N$ and $\psi\in  \ov\Psi_{j_1}^{d_1}$, $ f_{\psi}(\cdot)$ has the following wavelet expansion
   \begin{equation*}
       f_{\psi}(y)=\sum_{j_2=0}^{\infty}\sum_{\phi\in \ov\Psi_{j_2}^{d_2}}\phi(y) f_{\psi,\phi},\quad  f_{\psi,\phi}=\int_{\mb R^{d_2}}  f_{\psi}(y)\phi(y)\,\dd y=\int_{\mb R^{d_2}} \int_{\mb R^{d_1}}   f(x,y)\psi(x)\phi(y)\,\dd x\dd y,
 \end{equation*}
with $|\wt f_{\psi,\phi}|\leq C_WL_12^{-\frac{d_1j_1+d_2j_2}{2}}2^{-j_2\alpha_2}$ for any $\psi\in \Psi_{j_1}^{d_1}$ and $\phi\in \Psi_{j_2}^{d_2}$.  Then let $$ f'(x,y)=\sum_{j_1=0}^{J_1}\sum_{j_2=0}^{J_2} \sum_{\psi\in \ov\Psi_{j_1}^{d_1}}\sum_{\psi\in \ov\Psi_{j_2}^{d_2}}   f_{\psi,\phi} \psi(x)\phi(y),$$
we have 
\begin{equation}\label{eqnapproxwavelet}
    \begin{aligned}
       &  |  f'(x,y)-  f(x,y)|\leq \Big|\sum_{j=J_1+1}^{\infty}\sum_{\psi\in \ov\Psi_j^{d_1}}\psi(x)  f_{\psi}(y)\Big|+\Big|\sum_{j_1=0}^{J_1}\sum_{\psi\in \ov\Psi_{j_1}^{d_1}}\sum_{j_2=J_2+1}^{\infty}\sum_{\phi\in \ov\Psi_{j_2}^{d_2}}  f_{\psi,\phi}\psi(x)\phi(y)\Big|\\
       &\leq C_RC_L'C_W L\, 2^{-J_1\alpha_1}+C_RC_L'C_WL_1\,\sum_{j_1=0}^{J_1}\sum_{\psi\in \Psi_{j_1}^{d_1}} 2^{-\frac{d_1j_1}{2}}2^{-J_2\alpha_2} |\psi(x)|\\
       &\leq C_RC_L'C_W L\, 2^{-J_1\alpha_1}+2^{d_1}C_R^3 C_L'{}^2C_W C_L^{d_1} L J_1 2^{-J_2\alpha_2}.
    \end{aligned}
\end{equation}

\subsection{Proof of Lemma~\ref{le:appro}}\label{appendixE.3}
Without loss of generality, we assume $U_1\subseteq \mb B_{\mb R^{d_1}}(\mathbf{0},1)$ and $U_2\subseteq \mb B_{\mb R^{d_2}}(\mathbf{0},1)$.   Then  consider a smooth transition function
 \begin{equation}\label{eqnrho}
      \rho(t)=\left\{
      \begin{array}{cc}
        0   & |t|\geq 2 \\
        1   & |t|\leq 1\\
       \frac{1}{1+ \exp(\frac{3-2t}{(t-1)(t-2)})}& 1<t<2\\
        \frac{1}{1+ \exp(\frac{2t+3}{(t+1)(2+t)})}& -2<t<-1.\\
      \end{array}
      \right.
  \end{equation}
Set $\wt f(x,y)=\ov f(x,y)\rho(\|x\|^2)\rho(\|y\|^2)$. We have $\wt f(x,y)\in   \ov{\m H}^{\alpha_1,\alpha_2}_{L'}(\mb R^{d_1},\mb R^{d_2})$, $\wt f|_{U_1\times U_2}=\ov f|_{U_1\times U_2}$ and the support of $\wt f$ is contained in $\mb B_{\mb R^{d_1}}(\mathbf{0},\sqrt{2})\times \mb B_{\mb R^{d_2}}(\mathbf{0},\sqrt{2})$.   Consider two wavelet basis $\{\ov\Psi_j^{d_1}\}_{j\geq 0}$ and $\{\ov\Psi_j^{d_2}\}_{j\geq 0}$ that both satisfy the properties in Lemma~\ref{le:wavelet} with smoothness $\alpha=\lceil\alpha_1\vee\alpha_2\rceil$ and constants $C_R,C_L,C_L',C_L^{\dagger},C_L^{\ddagger},C_W, C_I$. For any $j\in \mb N$, define 
    \begin{equation*}
      \Psi_j^{d_1}=\{\psi\in \ov \Psi_j^{d_1}:\, \text{supp}(\psi)\cap \mb B_{\mb R^{d_1}}(\mathbf{0},\sqrt{2})\neq \emptyset\},
\end{equation*}
and
    \begin{equation*}
      \Psi_j^{d_2}=\{\psi\in \ov \Psi_j^{d_2}:\, \text{supp}(\psi)\cap \mb B_{\mb R^{d_2}}(\mathbf{0},\sqrt{2})\neq \emptyset\}.
\end{equation*}
we have $|\Psi_j^{d_1}|\leq \sqrt{2}C_L^{\dagger}  2^{d_1j}$ and  $|\Psi_j^{d_2}|\leq \sqrt{2}C_L^{\dagger}  2^{d_2j}$. Set $$J_1=\lceil \frac{\log (2C_RC_L'C_WL')+ \log \frac{1}{\varepsilon}}{\alpha_1\log 2}\rceil$$ and  $$J_2=\lceil \frac{\log (2^{d_1+1} C_R^3 C_L'{}^2 C_W C_L^{d_1}  L'J_1)+\log \frac{1}{\varepsilon}}{\alpha_2\log 2} \rceil.$$ Define 
\begin{equation*}
    \begin{aligned}
        \wt f'(x,y)&=\sum_{j_1=0}^{J_1}\sum_{j_2=0}^{J_2} \sum_{\psi\in \Psi_{j_1}^{d_1}}\sum_{\psi\in \Psi_{j_2}^{d_2}} \wt f_{\psi,\phi} \psi(x)\phi(y)\\
        &= \sum_{j_1=0}^{J_1}\sum_{j_2=0}^{J_2} \sum_{\psi\in 
        \ov\Psi_{j_1}^{d_1}}\sum_{\psi\in \ov\Psi_{j_2}^{d_2}} \wt f_{\psi,\phi} \psi(x)\phi(y),\quad \wt f_{\psi,\phi}=\int_{\mb R^{d_2}}\int_{\mb R^{d_1}} \wt f(x,y)\psi(x)\phi(y)\,\dd x\dd y.\\
    \end{aligned}
\end{equation*}
 It holds that 
\begin{equation*}
    \begin{aligned}
    &  |\wt f'(x,y)-\wt f(x,y)|\leq   C_RC_L'C_W L'\, 2^{-J_1\alpha_1}+2^{d_1}C_R^3 C_L'{}^2C_W C_L^{d_1} L' J_1 2^{-J_2\alpha_2}\leq \varepsilon.     
    \end{aligned}
\end{equation*}
Now we show that $\wt f'(x,y)\in  {\m H}^{\alpha_1,\alpha_2}_{L_0J_1J_2}(\mb R^{d_1},\mb R^{d_2})$ for a constant $L_0$.  Notice that for any $\psi\in \Psi_{j_1}^{d_1}$ and $\phi\in  \Psi_{j_2}^{d_2}$, we have 
\begin{equation*}
    \begin{aligned}
       & |\wt f_{\psi,\phi}|=\big|\int_{\mb R^{d_2}} \int_{\mb R^{d_1}} \wt f(x,y)\psi(x)\phi(y)\,\dd x\dd y\big|\leq \int_{\mb R^{d_2}} |\wt f_{\psi}(y)|\cdot|\phi(y)|\,\dd y\\
        &\leq C_W L 2^{-\frac{d_1j_1}{2}-j_1\alpha_1}\int_{\mb R^{d_2}} |\phi(y)|\,\dd y\\
        &\leq  C_W L(2C_L)^{d_1}C_R\, 2^{-\frac{d_1j_1+d_2j_2}{2}-j_1\alpha_1}.
    \end{aligned}
\end{equation*}
Combined with $|\wt f_{\psi,\phi}|\leq C_WL_12^{-\frac{d_1j_1+d_2j_2}{2}}2^{-j_2\alpha_2}$,  we can estbalish that, for some constant $L_2$,
\begin{equation*}
     |\wt f_{\psi,\phi}|\leq L_2 2^{-\frac{d_1j_1+d_2j_2}{2}}2^{-((j_1\alpha_1)\vee (j_2\alpha_2))}.
\end{equation*}
Then, for any $(l_1,l_2)\in \m J^{d_1,d_2}_{\alpha_1,\alpha_2}=\{l_1\in \mb N_0^{d_1}, l_2\in \mb N_0^{d_2}:\,\frac{|l_1|}{\alpha_1}+\frac{|l_2|}{\alpha_2}<1\}$, we have 
\begin{equation*}
    \begin{aligned}
      &  |\wt f'{}^{(l_1,l_2)}(x,y)|=|\sum_{j_1=0}^{J_1}\sum_{j_2=0}^{J_2} \sum_{\psi\in \Psi_{j_1}^{d_1}}\sum_{\psi\in \Psi_{j_2}^{d_2}} \wt f_{\psi,\phi} \psi^{(l_1)}(x)\phi^{(l_2)}(y)|\\
      &\leq L_2 \sum_{j_1=0}^{J_1}\sum_{j_2=0}^{J_2} \sum_{\psi\in \Psi_{j_1}^{d_1}}\sum_{\phi\in \Psi_{j_2}^{d_2}}  2^{-((j_1\alpha_1)\vee (j_2\alpha_2))} 2^{-\frac{d_1j_1+d_2j_2}{2}} |\psi^{(l_1)}(x)\phi^{(l_2)}(y)|\\
      &\leq  L_2 \sum_{j_1=0}^{J_1}\sum_{j_2=0}^{J_2} 2^{-((j_1\alpha_1)\vee (j_2\alpha_2))} 2^{-\frac{d_1j_1+d_2j_2}{2}}    C_R^2 \,2^{j_1|l_1|+\frac{d_1j_1}{2}}2^{j_2|l_2|+\frac{d_2j_2}{2}} \sum_{\psi\in \Psi_{j_1}} \mathbf{1}(x\in I_{\psi})\sum_{\phi\in \Psi_{j_2}} \mathbf{1}(y\in I_{\phi})\\
      &\leq L_3 \sum_{j_1=0}^{J_1}\sum_{j_2=0}^{J_2} 2^{-((j_1\alpha_1)\vee (j_2\alpha_2))}  2^{|l_1|j_1+|l_2|j_2}.
    \end{aligned}
\end{equation*}
Notice that when $j_1\leq \frac{j_2\alpha_2}{\alpha_1}$, we have 
\begin{equation*}
    -j_2\alpha_2+|l_1|j_1+|l_2|j_2\leq  -j_2\alpha_2+|l_1|\frac{j_2\alpha_2}{\alpha_1}+|l_2|j_2=j_2(|l_1|\frac{\alpha_2}{\alpha_1}+|l_2|-\alpha_2)<0,
\end{equation*}
and when $j_1\geq \frac{j_2\alpha_2}{\alpha_1}$, we have 
\begin{equation*}
    -j_1\alpha_1+|l_1|j_1+|l_2|j_2\leq  -j_1\alpha_1+|l_1|j_1+|l_2|\frac{j_1\alpha_1}{\alpha_2}=j_1(|l_2|\frac{\alpha_1}{\alpha_2}+|l_1|-\alpha_1)<0
\end{equation*}
Therefore,
\begin{equation}\label{eqnwtf}
         |\wt f'{}^{(l_1,l_2)}(x,y)|\leq L_3 \sum_{j_1=0}^{J_1} 2^{j_1(|l_2|\frac{\alpha_1}{\alpha_2}+|l_1|-\alpha_1)}+L_3\sum_{j_2=0}^{J_2} 2^{j_2(|l_1|\frac{\alpha_2}{\alpha_1}+|l_2|-\alpha_2)}\leq L_4.
\end{equation}
Then consider $(l_1,l_2)\in \m J^{d_1,d_2}_{\alpha_1,\alpha_2}$ with $\frac{|l_1|}{\alpha_1}+\frac{|l_2|}{\alpha_2}+\frac{1}{\alpha_1\wedge \alpha_2}\geq 1$, we claim that
\begin{claim}\label{claim1}
 There exists a constant $L_4$ so that for any $x,x'\in \mb R^{d_1}$, $y,y'\in \mb R^{d_2}$,  $j_1\in [J_1]$ and $j_2\in [J_2]$,
\begin{enumerate}
    \item  for any $(l_1,l_2)\in \m J^{d_1,d_2}_{\alpha_1,\alpha_2}$ with $\frac{|l_1|}{\alpha_1}+\frac{|l_2|}{\alpha_2}+\frac{1}{\alpha_1}\geq 1$,
    \begin{equation*}
\begin{aligned}
        &|\sum_{\psi\in \Psi_{j_1}^{d_1}}\sum_{\phi\in \Psi_{j_2}^{d_2}} \wt f_{\psi,\phi} \psi^{(l_1)}(x)\phi^{(l_2)}(y)-\sum_{\psi\in \Psi_{j_1}^{d_1}}\sum_{\phi\in \Psi_{j_2}^{d_2}} \wt f_{\psi,\phi} \psi^{(l_1)}(x')\phi^{(l_2)}(y)|\\
        &\leq L_4\|x-x'\|^{\alpha_1-|l_1|-\frac{\alpha_1}{\alpha_2}|l_2|}.
\end{aligned}
\end{equation*}
 \item  for any $(l_1,l_2)\in \m J^{d_1,d_2}_{\alpha_1,\alpha_2}$ with $\frac{|l_1|}{\alpha_1}+\frac{|l_2|}{\alpha_2}+\frac{1}{\alpha_2}\geq 1$,
    \begin{equation*}
\begin{aligned}
        &|\sum_{\psi\in \Psi_{j_1}^{d_1}}\sum_{\phi\in \Psi_{j_2}^{d_2}} \wt f_{\psi,\phi} \psi^{(l_1)}(x)\phi^{(l_2)}(y)-\sum_{\psi\in \Psi_{j_1}^{d_1}}\sum_{\phi\in \Psi_{j_2}^{d_2}} \wt f_{\psi,\phi} \psi^{(l_1)}(x)\phi^{(l_2)}(y')|\\
        &\leq L_4\|y-y'\|^{\alpha_2-|l_2|-\frac{\alpha_2}{\alpha_1}|l_1|}.
\end{aligned}
\end{equation*}
\end{enumerate}

\end{claim}
\noindent Then given Claim~\ref{claim1}, we can derive that for any $(l_1,l_2)\in \m J^{d_1,d_2}_{\alpha_1,\alpha_2}$ with $\frac{|l_1|}{\alpha_1}+\frac{|l_2|}{\alpha_2}+\frac{1}{\alpha_1}\geq 1$
\begin{equation*}
    |\wt f'{}^{(l_1,l_2)}(x,y)-\wt f'{}^{(l_1,l_2)}(x',y)|\leq J_1J_2 L_4\|x-x'\|^{\alpha_1-|l_1|-\frac{\alpha_1}{\alpha_2}|l_2|} ,
\end{equation*}
and for any $(l_1,l_2)\in \m J^{d_1,d_2}_{\alpha_1,\alpha_2}$ with $\frac{|l_1|}{\alpha_1}+\frac{|l_2|}{\alpha_2}+\frac{1}{\alpha_2}\geq 1$,
\begin{equation*}
    |\wt f'{}^{(l_1,l_2)}(x,y)-\wt f'{}^{(l_1,l_2)}(x,y')|\leq J_1J_2 L_4\|y-y'\|^{\alpha_2-|l_2|-\frac{\alpha_2}{\alpha_1}|l_1|}.
\end{equation*}
Together with~\eqref{eqnwtf}, these results confirm that $\wt f'\in  {\m H}^{\alpha_1,\alpha_2}_{L_0J_1J_2}(\mb R^{d_1},\mb R^{d_2})$ with some constant $L_0$. Finally, by choosing $f=\wt f'$, we can get the the desired result. 
 
 \medskip
 \noindent We now present the proof of Claim~\ref{claim1}.  Consider an arbitrary pair  $j_1\in [J_1]$ and $j_2\in [J_2]$. Without loss of generality, we assume that  $j_1\leq j_2\frac{\alpha_2}{\alpha_1}$. The proof  for the case where $j_1 \geq j_2\frac{\alpha_2}{\alpha_1}$ follows a similar argument. For the first statement, consider an arbitrary $(l_1,l_2)\in \m J^{d_1,d_2}_{\alpha_1,\alpha_2}$ with $\frac{|l_1|}{\alpha_1}+\frac{|l_2|}{\alpha_2}+\frac{1}{\alpha_1}\geq 1$, then when $\|x-x'\|\geq 2^{-j_2\frac{\alpha_2}{\alpha_1}}$, there exists a constant $L_4$ so that the following inequality holds:
\begin{equation*}
\begin{aligned}
        &\Big|\sum_{\psi\in \Psi_{j_1}^{d_1}}\sum_{\phi\in \Psi_{j_2}^{d_2}}  \wt f_{\psi,\phi} \psi^{(l_1)}(x)\phi^{(l_2)}(y)-\sum_{\psi\in \Psi_{j_1}^{d_1}}\sum_{\phi\in \Psi_{j_2}^{d_2}}  \wt f_{\psi,\phi} \psi^{(l_1)}(x')\phi^{(l_2)}(y)\Big|  \\
        &\leq      \Big |\sum_{\psi\in \Psi_{j_1}^{d_1}}\sum_{\phi\in \Psi_{j_2}^{d_2}}  \wt f_{\psi,\phi} \psi^{(l_1)}(x)\phi^{(l_2)}(y)\Big|+\Big|\sum_{\psi\in \Psi_{j_1}^{d_1}}\sum_{\phi\in \Psi_{j_2}^{d_2}}  \wt f_{\psi,\phi} \psi^{(l_1)}(x')\phi^{(l_2)}(y)\Big|\\
        &\leq L_4\, 2^{-j_2\alpha_2+j_1|l_1|+j_2|l_2|}\\
        &=L_4\,2^{j_2(-\alpha_2+\frac{j_1}{j_2}|l_1|+|l_2|)}\\
        &\leq L_4\, 2^{-j_2(\alpha_2-\frac{\alpha_2}{\alpha_1}|l_1|-|l_2|)}\\
        &\leq L_4\,\|x-x'\|^{\alpha_1-|l_1|-\frac{\alpha_1}{\alpha_2}|l_2|}.
\end{aligned}
\end{equation*}
When $\|x-x'\|\leq 2^{-j_2\frac{\alpha_2}{\alpha_1}}$, we have 
\begin{equation*}
    \begin{aligned}
       & \Big|\sum_{\psi\in \Psi_{j_1}^{d_1}}\sum_{\phi\in \Psi_{j_2}^{d_2}}  \wt f_{\psi,\phi} \psi^{(l_1)}(x)\phi^{(l_2)}(y)-\sum_{\psi\in \Psi_{j_1}^{d_1}}\sum_{\phi\in \Psi_{j_2}^{d_2}}  \wt f_{\psi,\phi} \psi^{(l_1)}(x')\phi^{(l_2)}(y)\Big| \\
       &\leq  \sum_{\psi\in \Psi_{j_1}^{d_1}}\sum_{\phi\in \Psi_{j_2}^{d_2}} |\psi^{(l_1)}(x)-\psi^{(l_1)}(x')|\cdot| \wt f_{\psi,\phi}\phi^{(l_2)}(y)|\\
&\leq C_RC_L'L_2 \,2^{-j_2\alpha_2-\frac{j_1d_1}{2}+j_2|l_2|} \sum_{\psi\in \Psi_{j_1}^{d_1}} |\psi^{(l_1)}(x)-\psi^{(l_1)}(x')| \\
       & \leq  C_RC_L'L_2\,2^{-j_2\alpha_2-\frac{j_1d_1}{2}+j_2|l_2|} \sum_{\psi\in \Psi_{j_1}^{d_1}} |\psi^{(l_1)}(x)-\psi^{(l_1)}(x')|\cdot\big(\mathbf{1}(x\in I_{\psi})+\mathbf{1}(x'\in I_{\psi})\big)\\
      &\leq L_4  2^{-j_2\alpha_2+j_1(|l_1|+1)+j_2|l_2|}\|x-x'\|,
    \end{aligned}
\end{equation*}
where the  last inequality uses $\|\nabla \psi^{(l_1)}(x)\|\lesssim 2^{\frac{j_1d_1}{2}+|l_1|j_1+j_1}$.
Given that $\|x-x'\|\leq 2^{-j_2\frac{\alpha_2}{\alpha_1}}$, $\frac{|l_1|}{\alpha_1}+\frac{|l_2|}{\alpha_2}+\frac{1}{\alpha_1}\geq 1$ and $j_1\leq j_2\frac{\alpha_2}{\alpha_1}$, we deduce that
\begin{equation*}
    \begin{aligned}
         &2^{-j_2\alpha_2+j_1(|l_1|+1)+j_2|l_2|}\|x-x'\|=  2^{-j_2\alpha_2+j_1(|l_1|+1)+j_2|l_2|}\|x-x'\|^{\alpha_1-|l_1|-\frac{\alpha_1}{\alpha_2}|l_2|}\|x-x'\|^{1-\alpha_1+|l_1|+\frac{\alpha_1}{\alpha_2}|l_2|}\\
         &\leq  2^{-j_2\alpha_2+j_1(|l_1|+1)+j_2|l_2|-j_2\frac{\alpha_2}{\alpha_1}(1-\alpha_1+|l_1|+\frac{\alpha_1}{\alpha_2}|l_2|)}\|x-x'\|^{\alpha_1-|l_1|-\frac{\alpha_1}{\alpha_2}|l_2|}\\
         & = 2^{j_1(|l_1|+1)-j_2\frac{\alpha_2}{\alpha_1}(1+|l_1|)}\|x-x'\|^{\alpha_1-|l_1|-\frac{\alpha_1}{\alpha_2}|l_2|}\\
         &\leq \|x-x'\|^{\alpha_1-|l_1|-\frac{\alpha_1}{\alpha_2}|l_2|}.
    \end{aligned}
\end{equation*}
This completes the proof of the first statement in Claim~\ref{claim1}. Next, we prove the second statement. Consider an arbitrary $(l_1,l_2)\in \m J^{d_1,d_2}_{\alpha_1,\alpha_2}$ with $\frac{|l_1|}{\alpha_1}+\frac{|l_2|}{\alpha_2}+\frac{1}{\alpha_2}\geq 1$.  When $\|y-y'\|\geq 2^{-j_2} $, we have 
\begin{equation*}
\begin{aligned}
        &\Big|\sum_{\psi\in \Psi_{j_1}^{d_1}}\sum_{\phi\in \Psi_{j_2}^{d_2}}  \wt f_{\psi,\phi} \psi^{(l_1)}(x)\phi^{(l_2)}(y)-\sum_{\psi\in \Psi_{j_1}^{d_1}}\sum_{\phi\in \Psi_{j_2}^{d_2}}  \wt f_{\psi,\phi} \psi^{(l_1)}(x)\phi^{(l_2)}(y')\Big|  \\
        &\leq      \Big |\sum_{\psi\in \Psi_{j_1}^{d_1}}\sum_{\phi\in \Psi_{j_2}^{d_2}}  \wt f_{\psi,\phi} \psi^{(l_1)}(x)\phi^{(l_2)}(y)\Big|+\Big|\sum_{\psi\in \Psi_{j_1}^{d_1}}\sum_{\phi\in \Psi_{j_2}^{d_2}}  \wt f_{\psi,\phi} \psi^{(l_1)}(x)\phi^{(l_2)}(y')\Big|\\
          &\leq L_4\, 2^{-j_2(\alpha_2-\frac{\alpha_2}{\alpha_1}|l_1|-|l_2|)}\\
        &\leq L_4\,\|y-y'\|^{\alpha_2-\frac{\alpha_2}{\alpha_1}|l_1|-|l_2|}.
\end{aligned}
\end{equation*}
When $\|y-y'\|\leq 2^{-j_2}$, we have 
\begin{equation*}
    \begin{aligned}
       & \Big|\sum_{\psi\in \Psi_{j_1}^{d_1}}\sum_{\phi\in \Psi_{j_2}^{d_2}}  \wt f_{\psi,\phi} \psi^{(l_1)}(x)\phi^{(l_2)}(y)-\sum_{\psi\in \Psi_{j_1}^{d_1}}\sum_{\phi\in \Psi_{j_2}^{d_2}}  \wt f_{\psi,\phi} \psi^{(l_1)}(x)\phi^{(l_2)}(y')\Big| \\
       &\leq \sum_{\psi\in \Psi_{j_1}^{d_1}}\sum_{\phi\in \Psi_{j_2}^{d_2}} | \wt f_{\psi,\phi} \psi^{(l_1)}(x')|\cdot|\phi^{(l_2)}(y)-\phi^{(l_2)}(y')|\\
       & \leq C_RC_L'L_2  2^{-j_2\alpha_2-\frac{j_2d_2}{2}+j_1|l_1|}\sum_{\phi\in \Psi_{j_2}^{d_2}}|\phi^{(l_2)}(y)-\phi^{(l_2)}(y')|\\
       & \leq  C_RC_L'L_2 2^{-j_2\alpha_2-\frac{j_2d_2}{2}+j_1|l_1|}\sum_{\phi\in \Psi_{j_2}^{d_2}}|\phi^{(l_2)}(y)-\phi^{(l_2)}(y')|\cdot\big(\mathbf{1}(y\in I_{\phi})+\mathbf{1}(y'\in I_{\phi})\big)\\
       &\leq L_4  2^{-j_2\alpha_2+j_1|l_1|+j_2(|l_2|+1)}\|y-y'\|,
    \end{aligned}
\end{equation*}
where the  last inequality uses  $\|\nabla \phi^{(l_2)}(y)\|\lesssim 2^{\frac{j_2d_2}{2}+|l_2|j_2+j_2}$.
Then given that  $\|y-y'\|\leq 2^{-j_2}$, $\frac{|l_1|}{\alpha_1}+\frac{|l_2|}{\alpha_2}+\frac{1}{\alpha_2}\geq 1$ and $j_1\leq j_2\frac{\alpha_2}{\alpha_1}$, we obtain
\begin{equation*}
    \begin{aligned}
       & 2^{-j_2\alpha_2+j_1|l_1|+j_2(|l_2|+1)}\|y-y'\|\leq 2^{-j_2\alpha_2+j_1|l_1|+j_2(|l_2|+1)}\|y-y'\|^{\alpha_2-\frac{\alpha_2}{\alpha_1}|l_1|-|l_2|}\|y-y'\|^{1-\alpha_2+\frac{\alpha_2}{\alpha_1}|l_1|+|l_2|}\\
       &\leq 2^{j_1|l_1|-j_2\frac{\alpha_2}{\alpha_1}|l_1|}\|y-y'\|^{\alpha_2-\frac{\alpha_2}{\alpha_1}|l_1|-|l_2|}\\
        &\leq \|y-y'\|^{\alpha_2-\frac{\alpha_2}{\alpha_1}|l_1|-|l_2|}.\\
    \end{aligned}
\end{equation*}
This completes the proof.

\subsection{Proof of Lemma~\ref{le:defmanifold}}\label{proof:lemmamanifoldpro}
\subsubsection{$(3)\Rightarrow (2)$}\label{sec:3to1}
  Consider a small enough positive constant $\ov\tau_2\leq \frac{\ov\tau}{2}$ that will be specified later, and take an arbitrary  point $\omega_0=(x_0,y_0)\in \m M=\{(x,y):\, x\in \m M_X,y\in \m M_{Y|x}\}$. Let $V^*\in \mb R^{D_Y\times d_Y}$  be a matrix whose column forms an orthonormal basis of $T_{\m M_{Y|x_0}} y_0$ and let $V^*{}^\perp\in \mb R^{D_Y\times (D_Y-d_Y)}$ be the orthogonal complement of $V^*$. Consider $\ov F_{\omega_0}\in  \ov{\m H}^{\beta_Y,
 \beta_X}_{\ov L,D_Y-d_Y}(\mb  R^{D_Y},\mb R^{D_X})$ so that $\ov F_{\omega_0}|_{\mb B_{\mb R^{D_Y}}(y_0,\ov \tau)\times\mb B_{\m M_X}(x_0,\ov \tau)}=F_{\omega_0}$.
Define $\mathfrak{F}_{\omega_0}: \mb B_{\mb R^{d_Y}}(\mathbf{0},\ov\tau_2)\times  \mb B_{\mb R^{D_Y-d_Y}}(\mathbf{0},\frac{\ov\tau}{2})\times \mb B_{\mb R^{D_X}}(x_0,\ov\tau_2)\to \mb R^{D_Y-d_Y}$ as
$$\mathfrak{F}_{\omega_0}(z,s,x)=\ov F_{\omega_0}(V^*z+V^*{}^\perp s+y_0,x).$$ 

\medskip
\noindent\textbf{Step 1.}\emph{ We will first show that the equation system, $\mathfrak{F}_{\omega_0}(z,s,x)=\mathbf{0}$ admits a solution for $s$ for any given $z \in \mathbb{B}_{\mathbb{R}^{d_Y}}(\mathbf{0},  \overline{\tau}_2)$ and $x \in \mathbb{B}_{\mathbb{R}^{D_X}}(x_0, \ov{\tau}_2)$.}

\medskip
\noindent It is straightforward to verify that there exists a constant $\ov L_1$ so that $\mathfrak{F}_{\omega_0}\in  {\m H}^{\beta_Y,\beta_X}_{\ov L_1,D_Y-d_Y}( \mb B_{\mb R^{d_Y}}(\mathbf{0},\ov\tau_2)\times  \mb B_{\mb R^{D_Y-d_Y}}(\mathbf{0},\frac{\ov\tau}{2}), \mb B_{\mb R^{D_X}}(x_0,\ov\tau_2))$. Moreover, for any $(x,y)\in  \mb B_{\m M_X}(x_0,\ov\tau)\times \mb B_{\mb R^{D_Y}}(y_0,\ov\tau)$, it holds that  $J_{F_{\omega_0}(\cdot,x)}(y)J_{F_{\omega_0}(\cdot,x)}(y)^T\succeq  \ov\tau_1 I_{D_Y-d_Y}$ and thus
$$ J_{\ov F_{\omega_0}(\cdot,x_0)}(y_0)V^*{}^\perp(V^*{}^\perp)^TJ_{\ov F_{\omega_0}(\cdot,x_0)}(y_0)^T=J_{\ov F_{\omega_0}(\cdot,x_0)}(y_0)J_{\ov F_{\omega_0}(\cdot,x_0)}(y_0)^T\succeq  \ov\tau_1 I_{D_Y-d_Y}.$$
When  $\ov\tau,\ov\tau_2$ are small enough, there exists a constant $L_2$ so that for any $z,z'\in  \mb B_{\mb R^{d_Y}}(\mathbf{0},\ov\tau_2)$, $s,s'\in\mb B_{\mb R^{d_Y}}(\mathbf{0},\frac{\ov\tau}{2})$, 
 and $x,x'\in \mb B_{\mb R^{D_X}}(x_0,\ov\tau_2)$, the following conditions are satisfied:
 \begin{equation*}
     J_{F_{\omega_0}(\cdot,x)}(V^*z+V^*{}^\perp s+y_0)J_{F_{\omega_0}(\cdot,x)}(V^*z+V^*{}^\perp s+y_0)^T\succeq \frac{\ov\tau_1}{2} I_{D_Y-d_Y},
 \end{equation*}
 \begin{equation*}
    \| \mathfrak{F}_{\omega_0}(z,0,x)-\mathfrak{F}_{\omega_0}(z',0,x')\|\leq  L_2(\|z-z'\|+\|x-x'\|^{\beta_X\wedge 1}),
 \end{equation*}
\begin{equation*}
    J_{\mathfrak{F}_{\omega_0}(z,\cdot,x)}(s)J_{\mathfrak{F}_{\omega_0}(z,\cdot,x)}(s)^T= J_{\ov F_{\omega_0}(\cdot,x)}(V^*z+V^*{}^\perp s+y_0)V^*{}^\perp(V^*{}^\perp)^TJ_{\ov F_{\omega_0}(\cdot,x)}(V^*z+V^*{}^\perp s+y_0)^T\succeq \frac{\ov\tau_1}{2} I_{D_Y-d_Y},
\end{equation*}
and
\begin{equation*}
    \|\mathfrak{F}_{\omega_0}(z,s,x)-\mathfrak{F}_{\omega_0}(z,s',x)-  J_{\mathfrak{F}_{\omega_0}(z,\cdot,x)}(s')(s-s')\|\leq  L_2\|s-s'\|^2.
\end{equation*}
For any $z\in \mb B_{\mb R^{d_Y}}(\mathbf{0},\ov\tau_2)$ and $x\in \mb B_{\mb R^{D_X}}(x_0,\ov\tau_2)$, we construct a solution $s(z,x)$ to the equation system, $\mathfrak{F}_{\omega_0}(z,s,x)=\mathbf{0}$ in $s$ as follows: define $s_0(z,x)=\mathbf{0}$ and for $k=1,2,\cdots$, we recursively define $$s_k(z,x)=s_{k-1}(z,x)-( J_{\mathfrak{F}_{\omega_0}(z,\cdot,x)}(s_{k-1}(z,x)))^{-1}\mathfrak{F}_{\omega_0}(z,s_{k-1}(z,x),x).$$  Then define a sequence $b_k=\frac{\sqrt{\ov\tau_1}}{\sqrt{2} L_2}(\frac{4 L_2^2}{\ov\tau_1}\ov\tau_2)^{2^k}$. We can set $\ov\tau_2$ to be small enough so that $\sum_{k=0}^{\infty} b_k< \frac{\ov\tau}{2}\wedge \frac{\sqrt{\ov\tau_1}}{2\sqrt{2} L_2}$, and we can verify that   for any $k\in \mb N$,
\begin{equation*}
    \begin{aligned}
      &  \|s_{k+1}(z,x)-s_k(z,x)\|\leq b_k,\\
      &\|\mathfrak{F}_{\omega_0}(z,s_k(z,x),x)\|\leq \sqrt{\frac{\ov\tau_1}{2}}b_k.\\
    \end{aligned}
\end{equation*}
Hence $s(z,x)=\lim_{k\to\infty} s_k(z,x)$ exists, $\mathfrak{F}_{\omega_0}(z,s(z,x),x)=\mathbf 0$ and $\|s(z,x)\|< \ov\tau_3=\frac{\ov\tau}{2}\wedge \frac{\sqrt{\ov\tau_1}}{2\sqrt{2} L_2}$. 

\medskip
\noindent\textbf{Step 2.} \emph{Now we demonstrate that for any $z\in \mb B_{\mb R^{d_Y}}(\mathbf{0},\ov\tau_2)$ and $x\in \mb B_{\mb R^{D_X}}(x_0,\ov\tau_2)$, the equation $\mathfrak{F}_{\omega_0}(z,s,x)=0$ has a unique solution over $s\in \mb B_{\mb R^{D_Y-d_Y}}(\mathbf{0},\ov\tau_3)$.}

\medskip
\noindent Suppose there are two  solution $s,s'$ on $\mb B_{\mb R^{D_Y-d_Y}}(\mathbf{0},\ov\tau_3)$, then 
\begin{equation*}
      \sqrt{\frac{\ov\tau_1}{2}}\|s-s'\|\leq \frac{\|s-s'\|}{\|(J_{\mathfrak{F}_{\omega_0}(z,\cdot,x)}(s'))^{-1}\|_{\rm op}} \leq \|J_{\mathfrak{F}_{\omega_0}(z,\cdot,x)}(s')(s-s')\|\leq  L_2\|s-s'\|^2.
\end{equation*}
So we have
\begin{equation*}
    \|s-s'\|\geq \frac{\sqrt{\ov\tau_1}}{\sqrt{2} L_2},
\end{equation*}
which causes contradiction.  Then  we define  a function $\wt G_{\omega_0}: \mb B_{\mb R^{d_Y}}(\mathbf{0},\frac{\ov\tau_2}{2})\times \mb B_{\m M_x}(x_0,\frac{\ov\tau_2}{2})\to \mb R^{D_Y}$ as $\wt G_{\omega_0}(z,x)=V^*z+V^*{}^{\perp} s(z,x)+y_0$, where $s(z,x)$ is defined as the unique solution of $\mathfrak{F}_{\omega_0}(z,s,x)=0$ over $s\in B_{\mb R^{D_Y-d_Y}}(\mathbf{0},\ov\tau_3)$, and define
$\wt Q_{\omega_0}(y,x)=V^*{}^T(y-y_0)$.  

\medskip
\noindent\textbf{Step 3.} \emph{We will show that the pair $(\wt  G_{\omega_0}, \wt Q_{\omega_0})$ satisfies the conditions in Statement (2) of Lemma~\ref{le:defmanifold}.}

 \medskip
\noindent  Notice that for any $x\in  \mb B_{\m M_x}(x_0,\frac{\ov\tau_2}{2})$ and $y\in \mb B_{\m M_{Y|x}}(y_0,\frac{\ov\tau_2}{2}\wedge \ov\tau_3)$, we have $\|V^*{}^T(y-y_0)\|<\frac{\ov\tau_2}{2}$ and $F_{\omega_0}(y,x)=\mathfrak{F}_{\omega_0}(V^*{}^T(y-y_0),(V^*{}^\perp)^T(y-y_0),x)=\mathbf 0$. Therefore, for any $x\in \mb B_{\m M_x}(x_0,\frac{\ov\tau_2}{2})$, 
it holds that 
$$\mb B_{\m M_{Y|x}}(y_0,\frac{\ov\tau_2}{2}\wedge \ov\tau_3)\subset \wt G_{\omega_0}(\mb B_{\mb R^{d_Y}}(\mathbf{0},\frac{\ov\tau_2}{2}),x).$$  Furthermore,  for any  $x\in  \mb B_{\m M_x}(x_0,\frac{\ov\tau_2}{2})$ and $z\in \mb B_{\mb R^{d_Y}}(\mathbf{0},\frac{\ov\tau_2}{2})$,  it holds that
\begin{equation*}
  F_{\omega_0}(\wt G_{\omega_0}(z,x),x)=\mathbf{0}_{D_Y-d_Y} \text{ and }\|\wt G_{\omega_0}(z,x)-y_0\|\leq \ov\tau\,\Rightarrow\,  \wt G_{\omega_0}(z,x)\in \m M_{Y|x},
\end{equation*}
and 
\begin{equation*}
\wt Q_{\omega_0}(\wt G_{\omega_0}(z,x),x)=V^*{}^TV^*z=z.
\end{equation*}
 Now it only remains to show the smoothness of $\wt G_{\omega_0}$. For $a>1$, consider the smooth transition function
 \begin{equation} 
      \rho_a(t)=\left\{
      \begin{array}{cc}
        0   & |t|\geq a \\
        1   & |t|\leq 1\\
       \frac{1}{1+ \exp(\frac{(a+1)-2t}{(t-1)(t-a)})}& 1<t<a\\
        \frac{1}{1+ \exp(\frac{(a+1)+2t}{(t+1)(a+t)})}& -a<t<-1.\\
      \end{array}
      \right.
  \end{equation}
We define $\ov G_{\omega_0}:\mb R^{d_Y}\times \mb R^{D_X}\to \mb R^{D_Y}$ as 
\begin{equation*}
    \ov G_{\omega_0}(z,x)=\left\{
    \begin{array}{ll}
   (V^*z+V^*{}^\perp s(z,x)+y_0 )\rho_{\frac{9}{4}}(\frac{4\|z\|^2}{\ov\tau_2^2}) \rho_{\frac{9}{4}}(\frac{4\|x\|^2}{\ov\tau_2^2}),  &  \,z\in B_{\mb R^{d_Y}}(\mathbf{0},\frac{3\ov \tau_2}{4}), x\in \mb B_{\mb R^{D_X}}(x_0,\frac{3\ov \tau_2}{4})\\
       \mathbf{0},&\, o.w,
    \end{array} 
    \right.
\end{equation*}
Then it holds that $\wt G_{\omega_0}=\ov G_{\omega_0}|_{\mb B_{\mb R^{d_Y}}(\mathbf{0},\frac{\ov\tau_2}{2})\times \mb B_{\m M_X}(x_0,\frac{\ov\tau_2}{2})}$, and we will show that $\ov G_{\omega_0}\in \m H^{\beta_Y,\beta_X}_{L,D_Y}(\mb R^{d_Y},\mb R^{D_X})$.

\medskip
\noindent When $\beta_X> 1$,   by implicit function theorem (see for example, Theorem A.3 of~\cite{Eldering2013}), for any $z\in  \mb B_{\mb R^{d_Y}}(\mathbf{0},\ov\tau_2)$ 
 and $x\in \mb B_{\mb R^{D_X}}(x_0,\ov\tau_2)$
\begin{equation*}
    J_{s(\cdot,x)}(z)=-(J_{\mathfrak{F}_{\omega_0}(z,\cdot,x)}(s(z,x)))^{-1}J_{\mathfrak{F}_{\omega_0}(\cdot,s(z,x),x)}(z),
\end{equation*}
 and 
\begin{equation*}
   J_{s(z,\cdot)}(x)=-(J_{\mathfrak{F}_{\omega_0}(z,\cdot,x)}(s(z,x)))^{-1}J_{\mathfrak{F}_{\omega_0}(z,s(z,x),\cdot)}(x).
 \end{equation*}
Given that $J_{\mathfrak{F}_{\omega_0}(z,\cdot,x)}(s(z,x))J_{\mathfrak{F}_{\omega_0}(z,\cdot,x)}(s(z,x))^T\succeq \frac{\ov\tau_1}{2} I_{D_Y-d_Y}$, we can verify the following: for any multi-indices $j_1\in \mb N_0^{d_Y}$ and $j_2\in \mb N_0^{D_X}$, if for all $l_1\in \mb N_0^{D_Y},l_2\in \mb N_0^{D_X}$ satisfying $|l_1|+|l_2|\leq |j_1|+|j_2|$ and $|l_2|\leq |j_2|$,the partial derivatives $\mathfrak{F}_{\omega_0}{}^{(l_1,l_2)}((z,s),x)$ exist and are uniformly bounded in absolute value across $(z,s)\in  \mb B_{\mb R^{d_Y}}(\mathbf{0},\ov\tau_2)\times B_{\mb R^{D_Y-d_Y}}(\mathbf{0},\ov\tau_3)$ and $x\in \mb B_{\mb R^{D_X}}(x_0,\ov\tau_2)$, then the partial derivatives $s^{(j_1,j_2)}(z,x)$ exist and are uniformly bounded in absolute value for $z\in  \mb B_{\mb R^{d_Y}}(\mathbf{0},\ov\tau_2)$ and $x\in \mb B_{\mb R^{D_X}}(x_0,\ov\tau_2)$.  

\medskip
\noindent Therefore, note that $\mathfrak{F}_{\omega_0}\in  \ov{\m H}^{\beta_Y,\beta_X}_{\ov L_1,D_Y-d_Y}\big( \mb B_{\mb R^{d_Y}}(\mathbf{0},\ov\tau_2)\times  \mb B_{\mb R^{d_Y}}(\mathbf{0},\frac{\ov\tau}{2}), \mb B_{\mb R^{D_X}}(x_0,\ov\tau_2)\big)$ with $\beta_Y\geq \beta_X$ and $\beta_Y\geq 2$.  For any $(j_1,j_2)\in  \m J^{d_Y,D_X}_{\beta_Y,\beta_X}$ and $l_1\in \mb N_0^{D_Y},l_2\in \mb N_0^{D_X}$ satisfying $|l_1|+|l_2|\leq |j_1|+|j_2|$ and $|l_2|\leq |j_2|$, we have 
\begin{equation*}
\begin{aligned}
     & \frac{|l_1|}{\beta_Y}+\frac{|l_2|}{\beta_X}=\frac{|l_1|}{\beta_Y}+\frac{|l_2|}{\beta_Y}+|l_2|(\frac{1}{\beta_X}-\frac{1}{\beta_Y})\leq \frac{|j_1|+|j_2|}{\beta_Y}+|j_2|(\frac{1}{\beta_X}-\frac{1}{\beta_Y})=\frac{|j_1|}{\beta_Y}+\frac{|j_2|}{\beta_X}<1,
\end{aligned}
\end{equation*}
and thus $(l_1,l_2)\in \m J^{D_Y,D_X}_{\beta_Y,\beta_X}$ and  $\mathfrak{F}_{\omega_0}{}^{(l_1,l_2)}((z,s),x)$ are uniformly bounded in absolute values. Therefore, there exists a constant $L_3$ so that  for any $k\in [D_Y-d_Y]$, the $k$-th component $s_k(z,x)$ of $s(z,x)=(s_1(z,x),s_2(z,x),\cdots, s_{D_Y-d_Y}(z,x))$ satisfies
\begin{equation*}
    \sum_{(j_1,j_2)\in \m J^{d_Y,D_X}_{\beta_Y,\beta_X}} \sup_{(z,x)\in \mb B_{\mb R^{d_Y}}(\mathbf{0},\ov\tau_2) \times  \mb B_{\mb R^{D_X}}(x_0,\ov\tau_2)} |s_k^{(j_1,j_2)}(z,x)|\leq L_3.
\end{equation*}
Moreover, for any  $(j_1,j_2)\in  \m J^{d_Y,D_X}_{\beta_Y,\beta_X}$ with $\frac{|j_1|}{\beta_Y}+\frac{|j_2|}{\beta_X}+\frac{1}{\beta_Y}\geq 1$ and $l_1\in \mb N_0^{D_Y},l_2\in \mb N_0^{D_X}$ satisfying $|l_1|+|l_2|\leq |j_1|+|j_2|$ and $|l_2|\leq |j_2|$, 
\begin{enumerate}
    \item If $\frac{|l_1|}{\beta_Y}+\frac{|l_2|}{\beta_X}+\frac{1}{\beta_Y}\geq 1$, then for any $z,z'\in  \mb B_{\mb R^{d_Y}}(\mathbf{0},\ov\tau_2)$ and $x\in \mb B_{\mb R^{D_X}}(x_0,\ov\tau_2)$,
    \begin{equation*}
    \begin{aligned}
             & \|\mathfrak{F}_{\omega_0}{}^{(l_1,l_2)}((z,s(z,x)),x)-  \mathfrak{F}_{\omega_0}{}^{(l_1,l_2)}((z',s(z',x)),x)\|\\
             &\lesssim \|z-z'\|^{\beta_Y-|l_1|-\frac{\beta_Y}{\beta_X}|l_2|}+\|s(z,x)-s(z',x)\|^{\beta_Y-|l_1|-\frac{\beta_Y}{\beta_X}|l_2|}\\
             &\lesssim \|z-z'\|^{\beta_Y-|l_1|-\frac{\beta_Y}{\beta_X}|l_2|}\\
             &\lesssim  \|z-z'\|^{\beta_Y-|j_1|-\frac{\beta_Y}{\beta_X}|j_2|}.
    \end{aligned}
    \end{equation*}
     
    \item If $\frac{|l_1|}{\beta_Y}+\frac{|l_2|}{\beta_X}+\frac{1}{\beta_Y}<1$, then for any $z,z'\in  \mb B_{\mb R^{d_Y}}(\mathbf{0},\ov\tau_2)$ and $x\in \mb B_{\mb R^{D_X}}(x_0,\ov\tau_2)$,
    \begin{equation*}
    \begin{aligned}
             & \|\mathfrak{F}_{\omega_0}{}^{(l_1,l_2)}((z,s(z,x)),x)-  \mathfrak{F}_{\omega_0}{}^{(l_1,l_2)}((z',s(z',x)),x)\|\\
             &\lesssim \|z-z'\|+\|s(z,x)-s(z',x)\|\\
             &\lesssim \|z-z'\| \lesssim \|z-z'\|^{\beta_Y-|j_1|-\frac{\beta_Y}{\beta_X}|j_2|}.
    \end{aligned}
    \end{equation*}
\end{enumerate}
Therefore, there exists a constant $L_3$ so that  for any $k\in [D_Y-d_Y]$, the $k$-th component $s_k(z,x)$ of $s(z,x)$ satisfies
 \begin{equation*}
    \begin{aligned}
\sum_{(j_1,j_2)\in \m J^{d_Y,D_X}_{\beta_Y,\beta_X}\atop\frac{|j_1|}{\beta_Y}+\frac{|j_2|}{\beta_X}+\frac{1}{\beta_Y}\geq 1}\underset{z,z_0\in \mb B_{\mb R^{d_Y}}(\mathbf{0},\ov\tau_2),x\in \mb B_{\mb R^{D_X}}(x_0,\ov\tau_2)\atop z\neq z_0}{\sup}\frac{| s_k^{(j_1,j_2)}(z,x)- s_k^{(j_1,j_2)}(z_0,x)|}{\|z-z_0\|^{\beta_Y-|j_1|-\frac{\beta_Y}{\beta_X}|j_2|}}\leq L_3.
    \end{aligned}
\end{equation*}
Furthermore, for any  $(j_1,j_2)\in  \m J^{d_Y,D_X}_{\beta_Y,\beta_X}$ with $\frac{|j_1|}{\beta_Y}+\frac{|j_2|}{\beta_X}+\frac{1}{\beta_X}\geq 1$ and $l_1\in \mb N_0^{D_Y},l_2\in \mb N_0^{D_X}$ satisfying $|l_1|+|l_2|\leq |j_1|+|j_2|$ and $|l_2|\leq |j_2|$, 
\begin{enumerate}
    \item If $\frac{|l_1|}{\beta_Y}+\frac{|l_2|}{\beta_X}+\frac{1}{\beta_Y}\geq 1$, then for any $z\in  \mb B_{\mb R^{d_Y}}(\mathbf{0},\ov\tau_2)$ and $x,x'\in \mb B_{\mb R^{D_X}}(x_0,\ov\tau_2)$,
    \begin{equation*}
    \begin{aligned}
             & \|\mathfrak{F}_{\omega_0}{}^{(l_1,l_2)}((z,s(z,x)),x)-  \mathfrak{F}_{\omega_0}{}^{(l_1,l_2)}((z,s(z,x')),x')\|\\
             &\lesssim  \|s(z,x)-s(z,x')\|^{\beta_Y-|l_1|-\frac{\beta_Y}{\beta_X}|l_2|}+\|x-x'\|^{\beta_X-|l_2|-\frac{\beta_X}{\beta_Y}|l_1|}\\
             &\lesssim \|x-x'\|^{\beta_X-|l_2|-\frac{\beta_X}{\beta_Y}|l_1|}\\
             &\lesssim \|x-x'\|^{\beta_X-|j_2|-\frac{\beta_X}{\beta_Y}|j_1|}.
    \end{aligned}
    \end{equation*}
     
    \item If $\frac{|l_1|}{\beta_Y}+\frac{|l_2|}{\beta_X}+\frac{1}{\beta_Y}<1$ and $\frac{|l_1|}{\beta_Y}+\frac{|l_2|}{\beta_X}+\frac{1}{\beta_X}\geq 1$, then for any $z,z'\in  \mb B_{\mb R^{d_Y}}(\mathbf{0},\ov\tau_2)$ and $x\in \mb B_{\mb R^{D_X}}(x_0,\ov\tau_2)$,
    \begin{equation*}
    \begin{aligned}
             & \|\mathfrak{F}_{\omega_0}{}^{(l_1,l_2)}((z,s(z,x)),x)-  \mathfrak{F}_{\omega_0}{}^{(l_1,l_2)}((z,s(z,x')),x')\|\\
             &\lesssim \|s(z,x)-s(z,x')\|+\|x-x'\|^{\beta_X-|l_2|-\frac{\beta_X}{\beta_Y}|l_1|}\\
              &\lesssim \|x-x'\|^{\beta_X-|l_2|-\frac{\beta_X}{\beta_Y}|l_1|}\\
             &\lesssim \|x-x'\|^{\beta_X-|j_2|-\frac{\beta_X}{\beta_Y}|j_1|}.
    \end{aligned}
    \end{equation*}
      \item If  $\frac{|l_1|}{\beta_Y}+\frac{|l_2|}{\beta_X}+\frac{1}{\beta_X}<1$, then for any $z,z'\in  \mb B_{\mb R^{d_Y}}(\mathbf{0},\ov\tau_2)$ and $x\in \mb B_{\mb R^{D_X}}(x_0,\ov\tau_2)$,
    \begin{equation*}
    \begin{aligned}
             & \|\mathfrak{F}_{\omega_0}{}^{(l_1,l_2)}((z,s(z,x)),x)-  \mathfrak{F}_{\omega_0}{}^{(l_1,l_2)}((z,s(z,x')),x')\|\\
             &\lesssim \|s(z,x)-s(z,x')\|+\|x-x'\|\\
             &\lesssim \|x-x'\|^{\beta_X-|j_2|-\frac{\beta_X}{\beta_Y}|j_1|}.
    \end{aligned}
    \end{equation*}
\end{enumerate}
Therefore, there exists a constant $L_3$ so that  for any $k\in [D_Y-d_Y]$, the $k$-th component $s_k(z,x)$ of $s(z,x)$ satisfies
  \begin{equation*}
    \begin{aligned}
\sum_{(j_1,j_2)\in \m J^{d_Y,D_X}_{\beta_Y,\beta_X}\atop\frac{|j_1|}{\beta_Y}+\frac{|j_2|}{\beta_X}+\frac{1}{\beta_X}\geq 1}\underset{z\in \mb B_{\mb R^{d_Y}}(\mathbf{0},\ov\tau_2),x,x_0\in \mb B_{\mb R^{D_X}}(x_0,\ov\tau_2)\atop x\neq x_0}{\sup}\frac{| s_k^{(j_1,j_2)}(z,x)- s_k^{(j_1,j_2)}(z,x_0)|}{\|x-x_0\|^{\beta_X-|j_2|-\frac{\beta_X}{\beta_Y}|j_1|}}\leq L_3.
    \end{aligned}
\end{equation*}
 So by combining all pieces, we establish that for any $k\in [D_Y-d_Y]$
 \begin{equation*}
     \begin{aligned}
           &\sum_{(j_1,j_2)\in \m J^{d_Y,D_X}_{\beta_Y,\beta_X}} \sup_{(z,x)\in \mb B_{\mb R^{d_Y}}(\mathbf{0},\ov\tau_2) \times  \mb B_{\mb R^{D_X}}(x_0,\ov\tau_2)} |s_k^{(j_1,j_2)}(z,x)|\\
           &+\sum_{(j_1,j_2)\in \m J^{d_Y,D_X}_{\beta_Y,\beta_X}\atop\frac{|j_1|}{\beta_Y}+\frac{|j_2|}{\beta_X}+\frac{1}{\beta_Y}\geq 1}\underset{z,z_0\in \mb B_{\mb R^{d_Y}}(\mathbf{0},\ov\tau_2),x\in \mb B_{\mb R^{D_X}}(x_0,\ov\tau_2)\atop z\neq z_0}{\sup}\frac{| s_k^{(j_1,j_2)}(z,x)- s_k^{(j_1,j_2)}(z_0,x)|}{\|z-z_0\|^{\beta_Y-|j_1|-\frac{\beta_Y}{\beta_X}|j_2|}}\\
           &+\sum_{(j_1,j_2)\in \m J^{d_Y,D_X}_{\beta_Y,\beta_X}\atop\frac{|j_1|}{\beta_Y}+\frac{|j_2|}{\beta_X}+\frac{1}{\beta_X}\geq 1}\underset{z\in \mb B_{\mb R^{d_Y}}(\mathbf{0},\ov\tau_2),x,x_0\in \mb B_{\mb R^{D_X}}(x_0,\ov\tau_2)\atop x\neq x_0}{\sup}\frac{| s_k^{(j_1,j_2)}(z,x)- s_k^{(j_1,j_2)}(z,x_0)|}{\|x-x_0\|^{\beta_X-|j_2|-\frac{\beta_X}{\beta_Y}|j_1|}}\leq L_4.
     \end{aligned}
 \end{equation*}
Utilizing the fact that when $z\in B_{\mb R^{d_Y}}(\mathbf{0},\ov \tau_2)$ and $x\in \mb B_{\mb R^{D_X}}(x_0,\ov \tau_2)$, $\ov G_{\omega_0}(z,x)= (V^*z+V^*{}^\perp s(z,x)+y_0 )\rho_{\frac{9}{4}}(\frac{4\|z\|^2}{\ov\tau^2}) \rho_{\frac{9}{4}}(\frac{4\|x\|^2}{\ov\tau^2})$; and when $(z,x)\notin B_{\mb R^{d_Y}}(\mathbf{0},\frac{3\ov \tau_2}{4})\times B_{\mb R^{D_X}}(x_0,\frac{3\ov \tau_2}{4})$, $\ov G_{w_0}(z,x)=0$, we can obtain  $\ov G_{\omega_0}\in  {\m H}^{\beta_Y,\beta_X}_{L_5,D_Y}( \mb R^{d_Y},\mb R^{D_X})$ and thus $\wt G_{\omega_0}\in  {\m H}^{\beta_Y,\beta_X}_{L_5,D_Y}(\mb B_{\mb R^{d_Y}}(\mathbf{0},\frac{\ov\tau_2}{2}),\mb B_{\m M_X}(x_0,\frac{\ov\tau_2}{2}))$.  
 
 \medskip
 \noindent Then we consider the case  when $\beta_X\leq 1$, similar to the case for $\beta_X>1$, using  implicit function theorem, it is straightforward to show that for any  $x\in \mb R^{D_X}$,  $\ov G_{\omega_0}(\cdot,x)\in \m H^{\beta_Y}_{L_5,D_Y}( \mb R^{d_Y})$.  Next, we shall demonstrate that  for any $l\in \mb N_0^{d_Y}$ with $|l|<\beta_Y$,  and for any $z\in \mb R^{d_Y}$, $x,x'\in \mb R^{D_X}$, it holds that 
 \begin{equation*}
    \big\| \ov G^{(l,\mathbf{0})}(z,x)-\ov G^{(l,\mathbf{0})}(z,x')\big\|\leq L_5\,\|x-x'\|^{\beta_X-\frac{\beta_X}{\beta_Y}|l|}.
 \end{equation*}
To verify this result, it suffices to prove that for any $l\in \mb N_0^{d_Y}$ with $|l|<\beta_Y$, and any $z\in  \mb B_{\mb R^{d_Y}}(\mathbf{0},\ov\tau_2)$ and $x,x'\in \mb B_{\mb R^{D_X}}(x_0,\ov\tau_2)$, 
 \begin{equation*}
     \|s^{(l,\mathbf{0})}(z,x)-s^{(l,\mathbf{0})}(z,x')\|\leq C_1\,\|x-x'\|^{\beta_X-\frac{\beta_X}{\beta_Y}|l|}.
 \end{equation*}
To establish this, note that
 \begin{equation*}
 \begin{aligned}
    & \|\mathfrak{F}_{\omega_0}(z,s(z,x'),x)-\mathfrak{F}_{\omega_0}(z,s(z,x),x)\|\\
    &=\|\mathfrak{F}_{\omega_0}(z,s(z,x'),x)\|\\
    &\leq  \|\mathfrak{F}_{\omega_0}(z,s(z,x'),x')\|+\|\mathfrak{F}_{\omega_0}(z,s(z,x'),x')-\mathfrak{F}_{\omega_0}(z,s(z,x'),x)\|\\
    &=\|\mathfrak{F}_{\omega_0}(z,s(z,x'),x')-\mathfrak{F}_{\omega_0}(z,s(z,x'),x)\|\\
    &\leq C\,\|x-x'\|^{\beta_X},\\
    & \text{ and }J_{\mathfrak{F}_{\omega_0}(z,\cdot,x)}(s)J_{\mathfrak{F}_{\omega_0}(z,\cdot,x)}(s)^T\succeq \frac{\ov\tau_1}{2} I_{D_Y-d_Y},
     \end{aligned}
 \end{equation*}
 we can get $\|s(z,x)-s(z,x')\|\leq C_1\, \|x-x'\|^{\beta_X}$.   So for any $l\in \mb N_0^{d_Y}$ with $|l|<\beta_Y-1$,  it holds that 
 \begin{equation*}
     \begin{aligned}
     &  \big\|  s^{(l,\mathbf{0})}(z,x)- s^{(l,\mathbf{0})}(z,x')\big\|\\
       &\leq \|x-x'\|^{\beta_X-\frac{\beta_X}{\beta_Y}|l|}+\|s(z,x)-s(z,x')\|\\
       &\lesssim \|x-x'\|^{\beta_X-\frac{\beta_X}{\beta_Y}|l|}.
     \end{aligned}
 \end{equation*} 
For any  $l\in \mb N_0^{d_Y}$ with $|l|=\lfloor\beta_Y\rfloor$,  it holds that 
  \begin{equation*}
     \begin{aligned}
       &\big\|  s^{(l,\mathbf{0})}(z,x)- s^{(l,\mathbf{0})}(z,x')\big\|\\
       &\leq \|x-x'\|^{\beta_X-\frac{\beta_X}{\beta_Y}|l|}+\|s(z,x)-s(z,x')\|^{\beta_Y-|l|}\\
       &\lesssim \|x-x'\|^{\beta_X-\frac{\beta_X}{\beta_Y}|l|}+\|x-x'\|^{\beta_Y(\beta_X-\frac{\beta_X}{\beta_Y}|l|)}\\
       &\lesssim \|x-x'\|^{\beta_X-\frac{\beta_X}{\beta_Y}|l|}.
     \end{aligned}
 \end{equation*}
 We can then get the desired result by combining all pieces.

\subsubsection{$(2) \Rightarrow (1)$}\label{sec:1to2}

 {We first show that the conditions in (2) can imply that $\m M_{Y|x}$ has a reach that is uniformly lower bounded away from zero.} Suppose that there exists $x\in \m M_X$, so that the reach of $\m M_{Y|x}$ is smaller than $\tau$. Then by definition, there exists $y\in \mb R^{D_Y}$ and $y_1,y_2\in \m M_{Y|x}$, so  that $y_1\neq y_2$, $\|y-y_1\|=\|y-y_2\|<\tau$, 
$y-y_1\perp T_{y_1}\m M_{Y|x}$, and $y-y_2\perp T_{y_2}\m M_{Y|x}$.  Let $\omega=(x,y_1)$ and  consider the local parametrization $(\wt Q_{\omega}, \wt G_{\omega})$. It holds that $\wt G_{\omega}(\mathbf 0,x)=y_1$. Moreover, since $\|y_2-y_1\|\leq \|y-y_1\|+\|y-y_2\|<2\tau$, when $\tau\leq \frac{\wt\tau}{2}$,  it holds for $z_2=\wt Q_{\omega}(y_2,x)$ that 
\begin{equation*}
0<\|z_2\|=\|\wt Q_{\omega}(y_2,x)-\wt Q_{\omega}(y_1,x)\|\leq \wt L\|y_2-y_1\|<2\wt L\tau,
\end{equation*}
and  $\wt G_{\omega}(z_2,x)=y_2$. Furthermore, since $y-y_1\perp T_{y_1}\m M_{Y|x}$, let $V_{\omega}^{\perp}$ be a $D_Y$ by $(D_Y-d_Y)$  matrix whose columns form an orthornormal basis for the normal space of $T_{y_1}\m M_{Y|x}$, there exists a vector $s\in \mb R^{D_Y-d_Y}$ so that $y=y_1+V_{\omega}^{\perp} s$ and $\|s\|<\tau$.  Then by $y-y_2\perp T_{y_2}\m M_{Y|x}$, it holds that 
\begin{equation*}
   J_{\wt G_{\omega}(\cdot,x)}(z_2)^T (y_1+V_{\omega}^{\perp} s-y_2)=J_{\wt G_{\omega}(\cdot,x)}(z_2)^T(y-y_2)=\mathbf 0,
\end{equation*}
which implies that 
\begin{equation*}
    \|J_{\wt G_{\omega}(\cdot,x)}(z_2)^T (y_1-y_2)\|=\|J_{\wt G_{\omega}(\cdot,x)}(z_2)^TV_{\omega}^{\perp} s\|.
\end{equation*}
Then since $\wt G_{\omega}(\cdot,x)$ is $\beta_Y$-H\"{o}lder-smooth with $\beta_Y\geq 2$, we have 
\begin{equation*}
    \begin{aligned}
        &\|y_1-y_2+J_{\wt G_{\omega}(\cdot,x)}(z_2)z_2\|\\
        &=\|\wt G_{\omega}(\mathbf 0,x)-\wt G_{\omega}(z_2,x)-J_{\wt G_{\omega}(\cdot,x)}(\mathbf 0-z_2)\|\\
        &\leq  \wt L\sqrt{D_Yd_Y} \|z_2\|^2\\
        &< 2\wt L^2 \sqrt{D_Yd_Y}\tau\|z_2\|,
    \end{aligned}
\end{equation*}
and therefore,
\begin{equation*}
    \begin{aligned}
        & \|J_{\wt G_{\omega}(\cdot,x)}(z_2)^T (y_1-y_2)\|\\
         &>  \|J_{\wt G_{\omega}(\cdot,x)}(z_2)^TJ_{\wt G_{\omega}(\cdot,x)}(z_2)z_2\|-2\wt L^3 D_Yd_Y\tau\|z_2\|\\
         &>\sqrt{\lambda_{\min}(J_{\wt G_{\omega}(\cdot,x)}(z_2)^TJ_{\wt G_{\omega}(\cdot,x)}(z_2))}\|z_2\|-2\wt L^3 D_Yd_Y\tau\|z_2\|\\
         &\geq \left(\frac{1}{\wt L}-2\wt L^3  D_Yd_Y\tau\right)\|z_2\|,
    \end{aligned}
\end{equation*}
where the last inequality uses the $\wt L$-Lipschitzness of $\wt Q_{\omega}$.   Moreover, since $\|J_{\wt G_{\omega}(\cdot,x)}(\mathbf 0)^TV_{\omega}^{\perp}s\|=0$, we can obtain
\begin{equation*}
    \begin{aligned}
        &\|J_{\wt G_{\omega}(\cdot,x)}(z_2)^TV_{\omega}^{\perp} s\|=\|(J_{\wt G_{\omega}(\cdot,x)}(z_2)-J_{G_{\omega}(\cdot,x)}(\mathbf 0))^TV_{\omega}^{\perp} s\|\\
        &\leq \|J_{\wt G_{\omega}(\cdot,x)}(z_2)-J_{\wt G_{\omega}(\cdot,x)}(\mathbf 0)\|_{\rm F} \|s\|\\
        &\leq \wt Ld_Y\sqrt{D_Y}\|z_2\|\|s\|\\
        &< \wt L d_Y \sqrt{D_Y}\tau\|z_2\|.
    \end{aligned}
\end{equation*}
Therefore, we have 
\begin{equation*}
\begin{aligned}
        & \wt L d_Y \sqrt{D_Y}\tau\|z_2\|>\left(\frac{1}{\wt L}-2\wt L^3 D_Yd_Y\tau\right)\|z_2\|\\
         &\Rightarrow \tau> \Big(\wt L^2 d_Y\sqrt{D_Y}(1+2\wt L^2\sqrt{D_Y})\Big)^{-1}.
\end{aligned}
\end{equation*}
So by selecting $\tau=\frac{\wt \tau}{2}\wedge \big(\wt L^2 d_Y\sqrt{D_Y}(1+2\wt L^2\sqrt{D_Y})\big)^{-1}$, it holds for any $x\in \m M_X$ that the reach of $\m M_{Y|x}$ is lower bounded by $\tau$.

\medskip
 To complete our proof, it remains to show the smoothness of the inverse of the projection map onto the tangent space of the manifold. Notice that any tangent vector in $T_{\m M_{Y|x_0} y_0}$ can be uniquely represented by a $d_Y$-dimensional vector using an orthonormal basis of $T_{\m M_{Y|x_0} y_0}$. Therefore, by  
 selecting $V_{\omega_0}$ as an orthonormal basis of $T_{\m M_{Y|x_0} y_0}$ in Lemma~\ref{le:projection}, we can obtain the desired result.


\subsubsection {$(1)\Rightarrow (3)$}
 Take an arbitrary $\omega_0=(x_0,y_0)\in \m M$. Let $V_{\omega_0}\in \mb R^{D_Y\times d_Y}$  be  a matrix whose column forms an orthonormal basis of $T_{M_{Y|x_0}} y_0$ and  $V_{\omega_0}^\perp\in \mb R^{D_Y\times (D_Y-d_Y)}$ be the orthogonal complement of $V_{\omega_0}$. Given that the submanifold $\m M_{Y|x_0}$ has reach that is lower bounded by $\tau$, by Lemma 2 of~\cite{10.1214/18-AOS1685}, it holds with some constants $\tau_2,\tau_3>0$ so that $\mb B_{\m M_{Y|x}}(y_0,\tau_2)\subset \Phi_{\omega_0}(\mb B_{T_{\m M_{Y|x_0}} y_0}(0,\tau_1))\subset \mb B_{\m M_{Y|x}}(y_0,\tau_3)$, where $\Phi_{\omega_0}$ is defined as per Definition~\ref{def:manifold} in the main text. Now define 
\begin{equation*}
    F_{\omega_0}(y,x)=(V_{\omega_0}^\perp)^T(y-\Phi_{\omega_0}(V_{\omega_0}V_{\omega_0}^T(y-y_0),x)).
\end{equation*}
Then 
$J_{F_{\omega_0}(\cdot,x_0)}(y_0)J_{F_{\omega_0}(\cdot,x_0)}(y_0)^T=(V_{\omega_0}^\perp)^T(I_{D_Y}-V_{\omega_0}V_{\omega_0}^T)(I_{D_Y}-V_{\omega_0}V_{\omega_0}^T)V_{\omega_0}^\perp=I_{D_Y-d_Y}$. So there exist constants $0< \ov\tau<\tau_1\wedge \tau_2$ and $L_1$ so that $F_{\omega_0}\in  {\m H}^{\beta_Y,\beta_X}_{L_1,D_Y-d_Y}(\mb B_{\mb R^{D_Y}}(y_0,\ov\tau), \mb B_{\m M_X}(x_0,\ov\tau))$ and for any $(y,x)\in \mb B_{\mb R^{D_Y}}(y_0, \ov\tau)\times \mb B_{\m M_X}(x_0,\ov\tau)$,
$J_{F_{\omega_0}(\cdot,x)}(y)J_{F_{\omega_0}(\cdot,x)}(y)^T\succeq  \frac{1}{2} I_{D_Y-d_Y}$. 
Then we show that  for any $x\in\mb B_{\m M_X}(x_0,\ov \tau)$, $\mb B_{\m M_{Y|x}}(y_0,\ov\tau)=\{y\in\mb B_{\mb R^{D_Y}}(y_0,\ov\tau):\, F_{\omega_0}(y,x)=0\}$. Firstly, if  $y\in \mb B_{\m M_{Y|x} }(y_0,\ov \tau)$, then $F_{\omega_0}(y,x)=0$, which implies that 
$\mb B_{\m M_{Y|x}}(y_0,\ov\tau)\subset\{y\in\mb B_{\mb R^{D_Y}}(y_0,\ov\tau):\, F_{\omega_0}(y,x)=0\}$. Furthermore, if $y\in \mb B_{\mb R^{D_Y} }(y_0,\ov\tau)$ and $F_{\omega_0}(y,x)=\mathbf{0}$. Then define $y_1=\Phi_{\omega_0}(V_{\omega_0}V_{\omega_0}^T(y-y_0),x)\subset \m M_{Y|x}$. It holds that 
\begin{equation*}
    \|V_{\omega_0}V_{\omega_0}^T(y-y_0)\|\leq \|y-y_0\|\leq \ov \tau<\tau_1,
\end{equation*}
so
\begin{equation*}
   V_{\omega_0}V_{\omega_0}^T(y_1-y_0)={\rm Proj}_{T_{y_0}\m M_{Y|x_0}}(y_1-y_0)=V_{\omega_0}V_{\omega_0}^T(y-y_0) \,\Rightarrow \,V_{\omega_0}V_{\omega_0}^T(y-y_1)=\mathbf 0_{D_Y}\,\Rightarrow \,V_{\omega_0}^T(y-y_1)=\mathbf 0_{d_Y}.
\end{equation*}
Then combined with the fact that $F_{\omega_0}(y,x)=(V_{\omega_0}^\perp)^T(y-y_1)=\mathbf{0}_{D_Y-d_Y}$, we have $y=y_1\in  \m M_{Y|x}$. Therefore,
\begin{equation*}
    \mb B_{\m M_{Y|x}}(y_0,\ov \tau_2)=\{y\in  \mb B_{\mb R^{D_Y}}(y_0,\ov \tau_2):\, F_{\omega_0}(y,x)=\mathbf 0\},
\end{equation*}
this completes the proof.

\subsection{Proof of Lemma~\ref{le:projection}}\label{proof:le:projection}
 
Let $\tau_2\in(0,  \frac{\wt\tau}{2})$ be a sufficiently small positive constant, and take an arbitrary $\omega_0=(x_0,y_0)\in \m M$. Let $\ov G_{\omega_0}\in  \ov{\m H}^{\beta_Y,\beta_X}_{\wt L,D_Y}(\mb R^{d_Y},\mb R^{d_X})$ be a smooth extension of $\wt G_{\omega_0}$. For any $s\in \mb B_{\mb R^{d_Y}}(\mathbf{0},\tau_2)$ and $x\in \mb B_{\mb R^{D_X}}(x_0,\tau_2)$,  consider the following equation for $z\in \mb R^{d_Y}$:
 \begin{equation}\label{eqninverse}
     V_{\omega_0}^T(\ov G_{\omega_0}(z,x)-y_0)=s.
 \end{equation}
 Since $\wt Q_{\omega_0}(\cdot,x)$ is $\wt L$-Lipschitz, we have 
\begin{equation*}
J_{\ov G_{\omega_0}(\cdot,x_0)}(\mathbf 0)^T J_{\ov G_{\omega_0}(\cdot,x_0)}(\mathbf 0)\succeq \frac{1}{\wt L^2} I_{d_Y}.
\end{equation*}
Then let $\wt V_{\omega_0}\in \mb R^{D_Y\times d_Y}$ be an orthonormal matrix  with $\wt V_{\omega_0}\wt V_{\omega_0}^T=P_{\omega_0}$, since $$V_{\omega_0}^TP_{\omega_0}V_{\omega_0}=(V_{\omega_0}^T\wt V_{\omega_0})(V_{\omega_0}^T\wt V_{\omega_0})^T\succeq  \tau_0 I_{d_Y},$$ we have
\begin{equation*}
    \wt V_{\omega_0}^TV_{\omega_0}V_{\omega_0}^T \wt V_{\omega_0}=(V_{\omega_0}^T\wt V_{\omega_0})^T(V_{\omega_0}^T\wt V_{\omega_0})\succeq  \tau_0 I_{d_Y},
\end{equation*}
and
\begin{equation*}
J_{\ov G_{\omega_0}(\cdot,x_0)}(\mathbf 0)^T  V_{\omega_0}V_{\omega_0}^T J_{\ov G_{\omega_0}(\cdot,x_0)}(\mathbf 0)\succeq \frac{ \tau_0}{\wt L^2} I_{d_Y}.
\end{equation*} 
Then using the fact that  $\ov G_{\omega_0}\in  \ov{\m H}^{\beta_Y,\beta_X}_{\wt L,D_Y}(\mb R^{d_Y},\mb R^{d_X})$ with $\beta_Y\geq 2$ and $\beta_X>0$,  when $\wt\tau_1,\tau_2$ are small enough, we have for any $x\in \mb B_{\mb R^{D_X}}(x_0,\tau_2)$ and  $z\in \mb B_{\mb R^{d_Y}}(\mathbf{0},\wt\tau_1)$, 
\begin{equation*}
    J_{\ov G_{\omega_0}(\cdot,x)}(z)^T  V_{\omega_0}V_{\omega_0}^T   J_{\ov G_{\omega_0}(\cdot,x)}(z)\succeq \frac{\tau_0}{4\wt L^2} I_{d_Y}.
\end{equation*}
 So there exists a constant $L_1$ so that for any $x\in \mb B_{\mb R^{D_X}}(x_0,\tau_2)$  and $z,z'\in \mb B_{\mb R^{d_Y}}(\mathbf{0},\wt\tau_1)$, 
\begin{equation*}
      \|( V_{\omega_0}^TJ_{\ov G_{\omega_0}(\cdot,x)}(z))^{-1}\|_{\rm op}\leq L_1, 
\end{equation*}
and 
\begin{equation*}
    \|V_{\omega_0}^T\ov G_{\omega_0}(z,x)-V_{\omega_0}^T\ov G_{\omega_0}(z',x)-V_{\omega_0}^TJ_{\ov G_{\omega_0}(\cdot,x)}(z')(z-z')\|\leq L_1\|z-z'\|^2.
\end{equation*}
Then, by following a similar analysis to that outlined in the proof for  $(3)\Rightarrow (1)$ of Lemma~\ref{le:defmanifold} in Section~\ref{sec:3to1}, we can show that for sufficiently small $\tau_2$, there exists a function $\zeta: \mb B_{\mb R^{d_Y}}(\mathbf{0},\tau_2)\times  \mb B_{\mb R^{D_X}}(x_0,\tau_2)\to \mb B_{\mb R^{d_Y}}(\mathbf{0},\frac{\tau_1}{2}\wedge \frac{1}{2L_1^2})$, so that $\zeta(s,x)$ is the unique solution of $ V_{\omega_0}^T(\ov G_{\omega_0}(z,x)-y_0)=s$  over $z\in \mb B_{\mb R^{d_Y}}(\mathbf{0},\frac{\wt\tau_1}{2}\wedge \frac{1}{2 L_1^2})$. Then we can define $G_{\omega_0}: \mb B_{\mb R^{d_Y}}(\mathbf{0},\tau_2)\times \mb B_{\m M_X}(x_0,\tau_2) \to \mb R^{D_Y}$ as $G_{\omega_0}(z,x)=\ov G_{\omega_0}(\zeta(z,x),x)$.  Denote  $ Q_{\omega_0}(y)= V_{\omega_0}^T(y-y_0)$, for any $x\in  B_{\m M_X}(x_0,\frac{\tau_2}{2}\wedge \frac{\wt \tau_1}{4\wt L}\wedge \frac{1}{4\wt LL_1^2})$,
$y\in  \mb B_{\m M_{Y|x}}(y_0,\frac{\tau_2}{2}\wedge \frac{\wt \tau_1}{4\wt L}\wedge \frac{1}{4\wt LL_1^2})$ and $z\in \mb B_{\mb R^{d_Y}}(\mathbf{0},\frac{\tau_2}{2})$,  we have 
\begin{equation*}
\begin{aligned}
    & \|\wt Q_{\omega_0}(y,x)\|=\|\wt Q_{\omega_0}(y,x)-\wt Q_{\omega_0}(\wt G_{\omega_0}(\mathbf 0,x),x)\|\\
    &\leq L\|y-\wt G_{\omega_0}(\mathbf 0,x)\|\\
    &\leq L\|y-y_0\|+L\,\|\wt G_{\omega_0}(\mathbf 0,x_0)-\wt G_{\omega_0}(\mathbf 0,x)\|\\
    &< \frac{\wt \tau_1}{2}\wedge \frac{1}{2L_1^2},
\end{aligned}
\end{equation*}

\begin{equation*}
    G_{\omega_0}(Q_{\omega_0}(y),x)=  G_{\omega_0}(Q_{\omega_0}(\wt G_{\omega_0}(\wt Q_{\omega_0}(y,x),x)),x)=G_{\omega_0}(\wt Q_{\omega_0}(y,x),x)=y.
\end{equation*}
and 
\begin{equation*}
  Q_{\omega_0}(G_{\omega_0}(z,x))=z.
\end{equation*}
Therefore, for any $x\in  B_{\m M_X}(x_0,\frac{\tau_2}{2}\wedge \frac{\wt \tau_1}{4\wt L}\wedge \frac{1}{4\wt LL_1^2})$,  let $ U_{Y|x}= G_{\omega_0}(\mb B_{\mb R^{d_Y}}(\mathbf{0},\frac{\tau_2}{2}),x)$, it holds that (1) $Q_{\omega_0}$ is a diffeomorphism that maps $U_{Y|x}$ to  $\mb B_{\mb R^{d_Y}}(\mathbf{0},\frac{\tau_2}{2})$ with inverse $ G_{\omega_0}(\cdot,x)|_{\mb B_{\mb R^{d_Y}}(\mathbf{0},\frac{\tau_2}{2})}$. (2)  $\mb B_{\m M_{Y|x}}(y_0,\frac{\tau_2}{2}\wedge \frac{\wt \tau_1}{4L}\wedge \frac{1}{4LL_1^2})\subset  U_{Y|x}\subset \m M_{Y|x}$. 

\medskip
\noindent So it only remains to show the smoothness of $G_{\omega_0}$. 
By implicit function theorem, for any $ z\in B_{\mb R^{d_Y}}(\mathbf{0},\tau_2), x\in \mb B_{\mb R^{D_X}}(x_0,\tau_2)$, 
\begin{equation*}
    J_{\zeta(\cdot,x)}(z)=\left( V_{\omega_0}^TJ_{\ov G_{\omega_0}(\cdot,x)}(\zeta(z,x))\right)^{-1}.
\end{equation*}
\begin{equation*}
    J_{G_{\omega_0}(\cdot,x)}(z)=J_{\ov G_{\omega_0}(\cdot,x)}(\zeta(z,x))\left( V_{\omega_0}^TJ_{\ov G_{\omega_0}(\cdot,x)}(\zeta(z,x))\right)^{-1}.
\end{equation*}
 And  when $\beta_X>1$,
\begin{equation*}
    J_{\zeta(z,\cdot)}(x)=-\left( V_{\omega_0}^TJ_{\ov G_{\omega_0}(\cdot,x)}(\zeta(z,x))\right)^{-1}(V_{\omega_0}^T J_{\ov G_{\omega_0}(\zeta(z,x),\cdot)}(x)),
\end{equation*}
and
\begin{equation*}
    J_{ G_{\omega_0}(z,\cdot)}(x)=J_{\ov G_{\omega_0}(\cdot,x)}(\zeta(z,x))  J_{\zeta(z,\cdot)}(x)+J_{\ov G_{\omega_0}(\zeta(z,x),\cdot)}(x).
\end{equation*}
 Then  similar to the analysis  outlined in the proof  for  $(3)\Rightarrow (1)$ of Lemma~\ref{le:defmanifold}, using the fact that $\ov G_{\omega_0}\in  {\m H}^{\beta_Y,\beta_X}_{L,D_Y}(\mb R^{d_Y},\mb R^{d_X})$ with $\beta_Y\geq \beta_X$, and  $$ \|( V_{\omega_0}^TJ_{\ov G_{\omega_0}(\cdot,x)}(\zeta(z,x)))^{-1}\|_{\rm op}\leq L_1,$$  we can conclude that there exists a constant $L_2$ so that $G_{\omega_0}\in  {\m H}^{\beta_Y,\beta_X}_{L_2,D_Y}(\mb B_{\mb R^{d_Y}}(\mathbf{0},\tau_2), \mb B_{\m M_X}(x_0,\tau_2))$.

 \subsection{Proof of Lemma~\ref{le:defdistribution}}\label{proof:lemmadensitypro}
{We will begin by proving the first statement.} For $a>1$, consider the smooth transition function $\rho_a(\cdot)$ defined as
 \begin{equation*} 
      \rho_a(t)=\left\{
      \begin{array}{cc}
        0   & |t|\geq a \\
        1   & |t|\leq 1\\
       \frac{1}{1+ \exp(\frac{(a+1)-2t}{(t-1)(t-a)})}& 1<t<a\\
        \frac{1}{1+ \exp(\frac{(a+1)+2t}{(t+1)(a+t)})}& -a<t<-1.\\
      \end{array}
      \right.
  \end{equation*}
 Let $\{\omega_k=(x^*_k,y^*_k)\}_{k=1}^{K^*}\subset \m M$ be a $\frac{\tau}{\sqrt{2}}$-covering set of $\m M$. 
 For any $k\in [K^*]$, let $V_k$ be a matrix whose column forms an orthonormal basis of $T_{\m M_{Y|x_k^*}} y_k^*$, and denote $G_{[k]}(z,x)=\Phi_{\omega_k}(V_kz,x)$, $Q_{[k]}(y)=V_k^T(y-y_k^*)$, $\nu_{k}(z|x)= \nu_{\omega_k}(V_kz|x)$ and $U^{\omega_k}_{Y|x}=U_{\omega_k}\cap \m M_{Y|x}$.
 Then define the function
 \begin{equation*}
     u(y,x)=\sum_{k=1}^{K^*}  \frac{  \nu_k( Q_{[k]}(y)|x)\rho_{2}(\frac{2\|\omega_k-(x,y)\|^2}{\tau^2})\left({\rm det}\Big(J_{   G_{[k]}(\cdot,x)}\big( Q_{[k]}(y)\big)^TJ_{G_{[k]}(\cdot,x)}\big( Q_{[k]}(y)\big)\Big)\right)^{-\frac{1}{2}}}{\sum_{k=1}^{K^*} \rho_{2}(\frac{2\|\omega_k-(x,y)\|^2}{\tau^2})}.
 \end{equation*}
We will show that $u(\cdot,x)$ is the density function of $\mu^*_{Y|x}$ with respect to the volume measure of $\m M_{Y|x}$. In support of this objective, we present the following claim that will be proved later. 
\begin{claim}\label{claimdensity}
    For any $(x,y)\in \m M$ and $k\in [K^*]$, if $\|\omega_k-(x,y)\|<\tau$, then  $$u(y,x)=\nu_{k}(Q_{[k]}(y)|x)\left({\rm det}\Big(J_{  G_{[k]}(\cdot,x)}\big( Q_{[k]}(y)\big)^TJ_{G_{[k]}(\cdot,x)}\big( Q_{[k]}(y)\big)\Big)\right)^{-\frac{1}{2}}.$$ 
\end{claim}
\noindent Given the claim above, it follows that for any $x\in \m M_x$ and measurable function $f_1:\m M_{Y|x}\to \mb R$, 
 \begin{equation*}
 \begin{aligned}
      &\mb{E}_{\mu^*_{Y|x}}[f_1(Y)]\\
      &=\sum_{k=1}^{K^*} \mb{E}_{\mu^*_{Y|x}}\left[\frac{f_1(Y)\rho_{2}(\frac{2\|\omega_k-(x,Y)\|^2}{\tau^2})}{\sum_{k_1=1}^{K^*} \rho_{2}(\frac{2\|\omega_{k_1}-(x,Y)\|^2}{\tau^2})}\right]\\
      &\overset{(i)}{=}\sum_{k=1}^{K^*} \mb{E}_{\mu^*_{Y|x}|_{ U^{\omega_k}_{Y|x}}}\left[\frac{f_1(Y)\rho_{2}(\frac{2\|\omega_k-(x,Y)\|^2}{\tau^2})}{\sum_{k_1=1}^{K^*} \rho_{2}(\frac{2\|\omega_{k_1}-(x,Y)\|^2}{\tau^2})}\right]\\
      &=\sum_{k=1}^{K^*} \int_{\mb B_{\mb R^{d_Y}}(\mathbf{0},\tau_1)} \left[\frac{f_1(G_{[k]}(z,x))\rho_{2}(\frac{2\|\omega_k-(x,G_{[k]}(z,x))\|^2}{\tau^2})}{\sum_{k_1=1}^{K^*} \rho_{2}(\frac{2\|\omega_{k_1}-(x,G_{[k]}(z,x))\|^2}{\tau^2})}\right] \nu_k(z|x)\,\dd z\\
      &=\sum_{k=1}^{K^*} \int_{\mb B_{\mb R^{d_Y}}(\mathbf{0},\tau_1)} f_1(G_{[k]}(z,x))u(G_{[k]}(z,x),x)\cdot\sqrt{{\rm det}\Big(J_{  G_{[k]}(\cdot,x)}(z)^TJ_{G_{[k]}(\cdot,x)}(z)\Big)}\\
      &\qquad\qquad\cdot\frac{\rho_{2}(\frac{2\|\omega_k-(x,G_{[k]}(z,x))\|^2}{\tau^2})}{\sum_{k_1=1}^{K^*} \rho_{2}(\frac{2\|\omega_{k_1}-(x,G_{[k]}(z,x))\|^2}{\tau^2})}\,\dd z,\\
 \end{aligned}
 \end{equation*}
 where $(i)$ uses the fact that $\mb B_{\m M_{Y|x}}(y_k^*,\tau)\subset U^{\omega_k}_{Y|x}$. Therefore, $u(\cdot,x)$ is the density function of $\mu^*_{Y|x}$ with respect to the volume measure of $\m M_{Y|x}$. 
 
 \medskip
 \noindent {Now we will show the smoothness of $u$.}  Let $\ov G_{[k]}\in  {\m H}^{\beta_Y,
 \beta_X}_{L,D_Y}(\mb R^{d_Y}, \mb R^{D_X})$ be a smooth extension of $G_{[k]}$ and $\ov \nu_k\in   {\m H}^{\alpha_Y,
 \alpha_X}_{L,D_Y}(\mb R^{d_Y}, \mb R^{D_X})$ be a smooth extension of $\nu_k$. Then notice that for any $z\in \mb B_{\mb R^{d_Y}}(\mathbf{0},\tau_1)$,  and $x\in \mb B_{\m M_X}(x_0,\tau)$, 
 \begin{equation*}
     J_{\ov G_{[k]}(\cdot,x)}(z)^T   J_{\ov G_{[k]}(\cdot,x)}(z)\succeq \frac{1}{L^2}I_{d_Y}.
 \end{equation*}
When $\tau$ is small enough, it holds that for any $z\in \mb B_{\mb R^{d_Y}}(\mathbf{0},2\tau)$,  and $x\in \mb B_{\mb R^{D_X}}(x_0,2\tau)$,
  \begin{equation*}
     J_{\ov G_{[k]}(\cdot,x)}(z)^T   J_{\ov G_{[k]}(\cdot,x)}(z)\succeq \frac{1}{2L^2}I_{d_Y}.
 \end{equation*}
Then we define a function $s_k:\mb R^{D_Y}\times \mb R^{D_X}\to \mb R$ as 
\begin{equation*}
\begin{aligned}
  & s_k(y,x)=\\
  &\left\{
    \begin{array}{ll}
\ov \nu_k( Q_{[k]}(y)|x)\left({\rm det}\Big(J_{   \ov G_{[k]}(\cdot,x)}\big( Q_{[k]}(y)\big)^TJ_{\ov G_{[k]}(\cdot,x)}\big( Q_{[k]}(y)\big)\Big)\right)^{-\frac{1}{2}}\rho_{\frac{9}{4}}(\frac{\|\omega_k-(x,y)\|^2}{\tau^2}),  &\quad  (x,y)\in B_{\mb R^{D_Y+D_X}}(\omega_k,\frac{3\tau}{2})\\
       0,& \quad \text{otherwise},
    \end{array} 
    \right.
    \end{aligned}
\end{equation*}
and define $\ov u:\mb R^{D_Y}\times \mb R^{D_X}\to \mb R$ as
 \begin{equation*}
\begin{aligned}
  & \ov u(y,x)=\sum_{k=1}^{K^*}  \frac{  \rho_{2}(\frac{2\|\omega_k-(x,y)\|^2}{\tau^2}) s_k(y,x)}{\sum_{k=1}^{K^*} \rho_{2}(\frac{2\|\omega_k-(x,y)\|^2}{\tau^2})+\rho_2(2\sum_{k=1}^{K^*} \rho_{2}(\frac{2\|\omega_k-(x,y)\|^2}{\tau^2}))}.
    \end{aligned}
\end{equation*}
Then when $y\in \m M_{Y|x}$ and $x\in \m M_X$, since  $\{\omega_k=(x_k^*,y_k^*)\}_{k=1}^{K^*}\subset \m M$ is a $\frac{\tau}{\sqrt{2}}$-covering set of $\m M$, it holds that 
\begin{equation*}
    \sum_{k=1}^{K^*} \rho_{2}(\frac{2\|\omega_k-(x,y)\|^2}{\tau^2})\geq 1,
\end{equation*}
and thus 
\begin{equation*}
\begin{aligned}
      &\ov u(y,x)=\sum_{k=1}^{K^*}  \frac{  \rho_{2}(\frac{2\|\omega_k-(x,y)\|^2}{\tau^2}) s_k(y,x)}{\sum_{k=1}^{K^*} \rho_{2}(\frac{2\|\omega_k-(x,y)\|^2}{\tau^2})}\\
      &=\sum_{k=1}^{K^*}  \frac{  \rho_{2}(\frac{2\|\omega_k-(x,y)\|^2}{\tau^2})   \ov \nu_k( Q_{[k]}(y)|x) \left({\rm det}\Big(J_{   \ov G_{[k]}(\cdot,x)}\big( Q_{[k]}(y)\big)^TJ_{\ov G_{[k]}(\cdot,x)}\big( Q_{[k]}(y)\big)\Big)\right)^{-\frac{1}{2}} }{\sum_{k=1}^{K^*} \rho_{2}(\frac{2\|\omega_k-(x,y)\|^2}{\tau^2})}\\
     &= \sum_{k=1}^{K^*}  \frac{  \rho_{2}(\frac{2\|\omega_k-(x,y)\|^2}{\tau^2})    \nu_k( Q_{[k]}(y)|x) \left({\rm det}\Big(J_{    G_{[k]}(\cdot,x)}\big( Q_{[k]}(y)\big)^TJ_{ G_{[k]}(\cdot,x)}\big( Q_{[k]}(y)\big)\Big)\right)^{-\frac{1}{2}} }{\sum_{k=1}^{K^*} \rho_{2}(\frac{2\|\omega_k-(x,y)\|^2}{\tau^2})}\\
      &=u(y,x).
\end{aligned}
\end{equation*}
Moreover,  given that
\begin{equation*}
    \sum_{k=1}^{K^*} \rho_{2}(\frac{2\|\omega_k-(x,y)\|^2}{\tau^2})+\rho_2(2\sum_{k=1}^{K^*} \rho_{2}(\frac{2\|\omega_k-(x,y)\|^2}{\tau^2}))>\frac{1}{2},
\end{equation*}
in order to show that $\ov u \in  {\m H}^{\alpha_Y,\alpha_X}_{L_1}(\mb R^{D_Y},\mb R^{D_X})$ for some constant $L_1$, it suffices to 
 show that  each component of $J_{\ov G_{[k]}(\cdot,x)}(z)$, as a function  with input $(z,x)$, belongs to $ {\m H}^{\alpha_Y,\alpha_X}_{L_2}(\mb R^{d_Y},\mb R^{D_X})$ for a certain constant $L_2$. Then notice that   
$    \ov G_{[k]}\in  {\m H}^{\beta_Y,
 \beta_X}_{L,D_Y}(\mb R^{d_Y}, \mb R^{D_X})$ with $\beta_Y\geq \alpha_Y+1$ and $\beta_X\geq \alpha_X+\frac{\alpha_X}{\alpha_Y}$. For any $(j_1,j_2)\in \m J^{d_Y,D_X}_{\alpha_Y,\alpha_X}$, it holds that
 \begin{equation*}
 \begin{aligned}
       &\frac{|j_1|+1}{\beta_Y}+\frac{|j_2|}{\beta_X}\\
       &\leq \frac{|j_1|+1}{\alpha_Y+1}+\frac{|j_2|}{\alpha_X+\frac{\alpha_X}{\alpha_Y}}\\
       &=\frac{|j_1|+1+|j_2|\frac{\alpha_Y
       }{\alpha_X}}{\alpha_Y+1}\\
       &=\frac{\alpha_Y(\frac{|j_1|}{\alpha_Y}+\frac{|j_2|}{\alpha_X})+1}{\alpha+1}<1.
 \end{aligned}
 \end{equation*}
Hence, let $e_j\in \mb N_0^{d_Y}$ denote the multi-index with the $j$-th component being $1$ and all other components being $0$.  It holds for any $k\in [K^*]$ and $j\in [d_Y]$ that 
\begin{equation*}
  \sum_{(j_1,j_2)\in \m J^{d_Y,D_X}_{\alpha_Y,\alpha_X}}  \underset{(x,y)\in \mb R^{d_1}\times \mb R^{d_2}}{\sup}|G_{\sk}^{(e_j+j_1,j_2)}(z,x)|\leq L.
\end{equation*}
Furthermore, for any  $(j_1,j_2)\in \m J^{d_Y,D_X}_{\alpha_Y,\alpha_X}$ with $\frac{|j_1|}{\alpha_Y}+\frac{|j_2|}{\alpha_X}+\frac{1}{\alpha_Y}\geq 1$,   
\begin{enumerate}
    \item if $\frac{|j_1|+1}{\beta_Y}+\frac{|j_2|}{\beta_X}+\frac{1}{\beta_Y}<1$, then for any  $j\in [d_Y]$, $z,z_0\in \mb R^{d_z}$ with $z\neq z_0$ and $x\in \mb R^{D_X}$, 
    \begin{itemize}\normalfont
        \item if $\|z-z_0\|\geq 1$, then
        \begin{equation*}
            |G_{\sk}^{(e_j+j_1,j_2)}(z,x)- G_{\sk}^{(e_j+j_1,j_2)}(z_0,x)|\leq 2L\leq 2L\|z-z_0\|^{\alpha_Y-|j_1|-\frac{\alpha_Y}{\alpha_X}|j_2|}.
        \end{equation*}
        \item if $\|z-z_0\|<1$, then 
        \begin{equation*}
              |G_{\sk}^{(e_j+j_1,j_2)}(z,x)- G_{\sk}^{(e_j+j_1,j_2)}(z_0,x)|\leq \sqrt{d_Y}L\|z-z_0\|\leq \sqrt{d_Y}L\|z-z_0\|^{\alpha_Y-|j_1|-\frac{\alpha_Y}{\alpha_X}|j_2|}.
        \end{equation*}
    \end{itemize}
   \item if $\frac{|j_1|+1}{\beta_Y}+\frac{|j_2|}{\beta_X}+\frac{1}{\beta_Y}\geq 1$, then  since 
   \begin{equation*}
   \begin{aligned}
         &\beta_Y-(|j_1|+1)-\frac{\beta_Y}{\beta_X}|j_2|\\
         &=\beta_Y(1-\frac{|j_1|+1}{\beta_Y}-\frac{|j_2|}{\beta_X})\\
         &\geq (\alpha_Y+1)(1-\frac{|j_1|+1}{\beta_Y}-\frac{|j_2|}{\beta_X})\\
         &\geq (\alpha_Y+1)(1-\frac{|j_1|+1}{\alpha_Y+1}-\frac{|j_2|}{\alpha_X+\frac{\alpha_X}{\alpha_Y}})\\
       &=\alpha_Y-|j_1|-\frac{\alpha_Y}{\alpha_X}|j_2|,
   \end{aligned}
      \end{equation*}
   we have for any  $j\in [d_Y]$, $z,z_0\in \mb R^{d_z}$ with $z\neq z_0$ and $x\in \mb R^{D_X}$, 
     \begin{itemize}\normalfont
        \item if $\|z-z_0\|\geq 1$, then
        \begin{equation*}
            |G_{\sk}^{(e_j+j_1,j_2)}(z,x)- G_{\sk}^{(e_j+j_1,j_2)}(z_0,x)|\leq 2L\leq 2L\|z-z_0\|^{\alpha_Y-|j_1|-\frac{\alpha_Y}{\alpha_X}|j_2|}.
        \end{equation*}
        \item if $\|z-z_0\|<1$, then 
        \begin{equation*}
              |G_{\sk}^{(e_j+j_1,j_2)}(z,x)- G_{\sk}^{(e_j+j_1,j_2)}(z_0,x)|\leq L\|z-z_0\|^{\beta_Y-|j_1|-1-\frac{\beta_Y}{\beta_X}|j_2|}\leq L\|z-z_0\|^{\alpha_Y-|j_1|-\frac{\alpha_Y}{\alpha_X}|j_2|}.
        \end{equation*}
    \end{itemize}
    
\end{enumerate}
Therefore, there exists a constant $L'$ so that for any $j\in [d_Y]$,
\begin{equation*}
    \sum_{(j_1,j_2)\in \m J^{d_Y,D_X}_{\alpha_Y,\alpha_X}\atop\frac{|j_1|+1}{\alpha_Y}+\frac{|j_2|}{\alpha_X}\geq 1}\underset{z,z_0\in \mb R^{d_Y}, x\in \mb R^{D_X}\atop z\neq z_0}{\sup}\frac{| G_{\sk}^{(e_j+j_1,j_2)}(z,x)- G_{\sk}^{(e_j+j_1,j_2)}(z_0,x)|}{\|z-z_0\|^{\alpha_Y-|j_1|-\frac{\alpha_Y}{\alpha_X}|j_2|}}\leq L'.\\
    \end{equation*}
Similarly, using the fact that for any $j_1,j_2\in \m J^{d_Y,D_X}_{\alpha_Y,\alpha_X}$ with $\frac{|j_1|}{\alpha_Y}+\frac{|j_2|+1}{\alpha_X}\geq 1$,
\begin{equation*}
    \begin{aligned}
        \beta_X-|j_2|-\frac{\beta_X}{\beta_Y}(|j_1|+1)&=\beta_X(1-\frac{|j_2|}{\beta_X}-\frac{|j_1|+1}{\beta_Y})\\
        &\geq (\alpha_X+\frac{\alpha_X}{\alpha_Y})\cdot(1-\frac{|j_2|}{\alpha_X+\frac{\alpha_X}{\alpha_Y}}-\frac{|j_1|+1}{\alpha_Y+1})\\
        &=\alpha_X-|j_2|-\frac{\alpha_X}{\alpha_Y}|j_1|.
    \end{aligned}
\end{equation*}
 We can also show that for any $j\in [d_Y]$,
 \begin{equation*}
    \sum_{(j_1,j_2)\in \m J^{d_Y,D_X}_{\alpha_Y,\alpha_X}\atop\frac{|j_1|}{\alpha_Y}+\frac{|j_2|+1}{\alpha_X}\geq 1}\underset{z\in \mb R^{d_Y}, x,x_0\in \mb R^{D_X}\atop x\neq x_0}{\sup}\frac{| G_{\sk}^{(e_j+j_1,j_2)}(z,x)- G_{\sk}^{(e_j+j_1,j_2)}(z,x_0)|}{\|x-x_0\|^{\alpha_X-|j_2|-\frac{\alpha_X}{\alpha_Y}|j_1|}}\leq L'.\\
 \end{equation*}
By combining all pieces, we can obtain that there exists a constant $L_1$ so that
\begin{equation*}
  J_{\ov G_{[k]}(\cdot,x)}(\cdot)\in   {\m H}^{\alpha_Y,
 \alpha_X}_{L_1,D_Yd_Y}(\mb R^{d_Y},\mb R^{D_X}),
\end{equation*}
which implies that there exists a constant $L_2$ so that $\ov u \in  {\m H}^{\alpha_Y,\alpha_X}_{L_2}(\mb R^{D_Y},\mb R^{D_X})$. 

\medskip
\noindent{Now it remains to show Claim~\ref{claimdensity}.}  For any $(x,y)\in \m M$, if there exists $k_1\neq k_2$ so that $\|\omega_{k_1}-(x,y)\|<\tau$ and $\|\omega_{k_2}-(x,y)\|<\tau$, then by change of variable formula, we have 
\begin{equation*}
    \nu_{k_1}(Q_{[k_1]}(y)|x)= \nu_{k_2}(Q_{[k_2]}(y)|x)\sqrt{\big({\rm det}(J_{Q_{[k_2]}(G_{[k_1]}(\cdot,x))}(Q_{[k_1]}(y))^TJ_{Q_{[k_2]}(G_{[k_1]}(\cdot,x))}(Q_{[k_1]}(y)))\big)},
\end{equation*}
and 
\begin{equation*}
    \begin{aligned}
        &\nu_{k_1}(Q_{[k_1]}(y)|x)\left({\rm det}\Big(J_{  G_{[k_1]}(\cdot,x)}\big( Q_{[k_1]}(y)\big)^TJ_{G_{[k_1]}(\cdot,x)}\big( Q_{[k_1]}(y)\big)\Big)\right)^{-\frac{1}{2}}\\
        &= \nu_{k_2}(Q_{[k_2]}(y)|x)\\
        &\cdot\sqrt{\frac{{\rm det}\Big(J_{Q_{[k_2]}(G_{[k_1]}(\cdot,x))}(Q_{[k_1]}(y))^TJ_{Q_{[k_2]}(G_{[k_1]}(\cdot,x))}(Q_{[k_1]}(y))\Big)}{{\rm det}\Big(J_{  G_{[k_1]}(\cdot,x)}\big( Q_{[k_1]}(y)\big)^TJ_{G_{[k_1]}(\cdot,x)}\big( Q_{[k_1]}(y)\big)\Big)}}\\
        &= \nu_{k_2}(Q_{[k_2]}(y)|x)\\
        &\cdot\sqrt{\frac{{\rm det}\Big(J_{  G_{[k_1]}(\cdot,x)}\big( Q_{[k_1]}(y)\big)^TJ_{Q_{[k_2]}}(y)^TJ_{Q_{[k_2]}}(y)J_{  G_{[k_1]}(\cdot,x)}\big( Q_{[k_1]}(y)\big)\Big)}{{\rm det}\Big(J_{  G_{[k_1]}(\cdot,x)}\big( Q_{[k_1]}(y)\big)^TJ_{G_{[k_1]}(\cdot,x)}\big( Q_{[k_1]}(y)\big)\Big)}}.\\
    \end{aligned}
\end{equation*}
 Then using the fact that for any $y'\in \mb B_{\m M_{Y|x}}(y_{k_1}^*,\tau)\cap \mb B_{\m M_{Y|x}}(y_{k_2}^*,\tau)$
 \begin{equation}
     G_{[k_1]}(Q_{[k_1]}(y'),x)=G_{[k_2]}(Q_{[k_2]}(G_{[k_1]}(Q_{[k_1]}(y'),x)),x),
 \end{equation}
 we can get
 \begin{equation*}
     J_{G_{[k_1]}(\cdot,x)}(Q_{[k_1]}(y))=J_{G_{[k_2]}(\cdot,x)}(Q_{[k_2]}(y))J_{Q_{[k_2]}}(y)J_{G_{[k_1]}(\cdot,x)}(Q_{[k_1]}(y)).
 \end{equation*}
 So we can write
\begin{equation*}
    \begin{aligned}
     & {\rm det}\Big(J_{  G_{[k_1]}(\cdot,x)}\big( Q_{[k_1]}(y)\big)^TJ_{  G_{[k_1]}(\cdot,x)}\big( Q_{[k_1]}(y)\big)\Big)\\
    &  = {\rm det}\Big(J_{G_{[k_1]}(\cdot,x)}(Q_{[k_1]}(y))^TJ_{Q_{[k_2]}}(y)^TJ_{G_{[k_2]}(\cdot,x)}(Q_{[k_2]}(y))^T J_{G_{[k_2]}(\cdot,x)}(Q_{[k_2]}(y))J_{Q_{[k_2]}}(y)J_{G_{[k_1]}(\cdot,x)}(Q_{[k_1]}(y))\Big)\\
    & = {\rm det}\Big(J_{G_{[k_1]}(\cdot,x)}(Q_{[k_1]}(y))^TJ_{Q_{[k_2]}}(y)^T\Big){\rm det }\Big(J_{G_{[k_2]}(\cdot,x)}(Q_{[k_2]}(y))^T J_{G_{[k_2]}(\cdot,x)}(Q_{[k_2]}(y))\Big)\\
    &\qquad \qquad\cdot{\rm det}\Big(J_{Q_{[k_2]}}(y)J_{G_{[k_1]}(\cdot,x)}(Q_{[k_1]}(y))\Big)\\
    &={\rm det}\Big(J_{  G_{[k_1]}(\cdot,x)}\big( Q_{[k_1]}(y)\big)^TJ_{Q_{[k_2]}}(y)^TJ_{Q_{[k_2]}}(y)J_{  G_{[k_1]}(\cdot,x)}\big( Q_{[k_1]}(y)\big)\Big)\\
    &\qquad\qquad\cdot {\rm det }\Big(J_{G_{[k_2]}(\cdot,x)}(Q_{[k_2]}(y))^T J_{G_{[k_2]}(\cdot,x)}(Q_{[k_2]}(y))\Big).
    \end{aligned}
\end{equation*}
 Therefore, we have 
 \begin{equation*}
 \begin{aligned}
         &\nu_{k_1}(Q_{[k_1]}(y)|x)\left({\rm det}\Big(J_{  G_{[k_1]}(\cdot,x)}\big( Q_{[k_1]}(y)\big)^TJ_{G_{[k_1]}(\cdot,x)}\big( Q_{[k_1]}(y)\big)\Big)\right)^{-\frac{1}{2}}\\
         &=  \nu_{k_2}(Q_{[k_2]}(y)|x)\left({\rm det}\Big(J_{  G_{[k_2]}(\cdot,x)}\big( Q_{[k_2]}(y)\big)^TJ_{G_{[k_2]}(\cdot,x)}\big( Q_{[k_2]}(y)\big)\Big)\right)^{-\frac{1}{2}}.\\
 \end{aligned}
 \end{equation*}
That implies that for any $(x,y)\in \m M$ and $k\in [K^*]$, if $\|\omega_k-(x,y)\|<\tau$, then  $$u(y,x)=\nu_{k}(Q_{[k]}(y)|x)\left({\rm det}\Big(J_{  G_{[k]}(\cdot,x)}\big( Q_{[k]}(y)\big)^TJ_{G_{[k]}(\cdot,x)}\big( Q_{[k]}(y)\big)\Big)\right)^{-\frac{1}{2}},$$  
which concludes the proof of Claim~\ref{claimdensity}. The proof of  the first statement in Lemma~\ref{le:defdistribution} is now concluded.

\medskip
\noindent {Then we show the second statement in Lemma~\ref{le:defdistribution}. }  For any $\omega_0=(x_0,y_0)\in \m M$, we can express $\wt v_{\omega_0}$ as 
 $$\wt v_{\omega_0}(z,x)=u(\wt G_{\omega_0}(z,x)|x)\cdot \sqrt{{\rm det}(J_{\wt G_{\omega_0}(\cdot,x)}(z)^TJ_{\wt G_{\omega_0}(\cdot,x)}(z))}.$$ 
 Let $\ov G_{\omega_0}\in  \ov{\m H}^{\beta_Y,
 \beta_X}_{L,D_Y}(\mb R^{d_Y}, \mb R^{D_X})$ be a smooth extension of $\wt G_{\omega_0}$ and $\ov u\in   {\m H}^{\alpha_Y,
 \alpha_X}_{L,D_Y}(\mb R^{D_Y}, \mb R^{D_X})$ be a smooth extension of $u$, using $\beta_Y\geq 2\vee (\alpha_Y+1)$, $\beta_X\geq \alpha_X+\frac{\alpha_X}{\alpha_Y}$, $\alpha_Y\geq \alpha_X$, we have $\ov u(\ov G_{\omega_0}(z,x)|x)\in  {\m H}^{\alpha_Y,
 \alpha_X}_{L_1,D_Y}(\mb R^{d_Y}, \mb R^{D_X})$ for a constant $L_1$. Then notice that for any $z\in \mb B_{\mb R^{d_Y}}(\mathbf{0},\tau_1)$,  and $x\in \mb B_{\m M_X}(x_0,\tau)$, 
 \begin{equation*}
     J_{\ov G_{\omega_0}(\cdot,x)}(z)^T   J_{\ov G_{\omega_0}(\cdot,x)}(z)\succeq \frac{1}{L^2}I_{d_Y}.
 \end{equation*}
So there exist  a constant $\tau_2$  so that  for any $z\in \mb B_{\mb R^{d_Y}}(\mathbf{0},\tau_1+\tau_2)$,  and $x\in \mb B_{\mb R^{D_X}}(x_0,\tau+\tau_2)$,
  \begin{equation*}
     J_{\ov G_{\omega_0}(\cdot,x)}(z)^T   J_{\ov G_{\omega_0}(\cdot,x)}(z)\succeq \frac{1}{2L^2}I_{d_Y}.
 \end{equation*}
Therefore, consider the smooth transition function $\rho_a(\cdot)$, we define $\ov v_{\omega_0}:\mb R^{d_Y}\times \mb R^{D_X}\to \mb R$ as 
 
\begin{equation*}
\begin{aligned}
   & \ov v_{\omega_0}(z,x)= \ov u(\ov G_{\omega_0}(z,x)|x)  \sqrt{{\rm det}(J_{\ov G_{\omega_0}(\cdot,x)}(z)^TJ_{\ov G_{\omega_0}(\cdot,x)}(z))}    \rho_{(1+\frac{\tau_2}{2\tau_1})^2}(\frac{\|z\|^2}{\tau_1^2}) \rho_{(1+\frac{\tau_2}{2\tau})^2}(\frac{\|x-x_0\|^2}{\tau^2})
  \end{aligned}
\end{equation*}
 By applying the same argument as in the proof of statement 1, we can establish that 
  $ J_{\ov G_{\omega_0}(\cdot,x)}(z)\in   {\m H}^{\alpha_Y,
 \alpha_X}_{L_1,D_Yd_Y}(\mb R^{d_Y},  \mb R^{D_X})$ and therefore $ \ov v_{\omega_0}\in  {\m H}^{\alpha_Y,\alpha_X}_{L_2}(\mb R^{d_Y},\mb R^{D_X})$, for some constants $L_1,L_2$.  Additionally, $\ov\nu_{\omega_0}(z,x)=\wt \nu_{\omega_0}(z,x)$  holds for any $z\in \mb B_{\mb R^{d_Y}}(\mathbf{0},\tau_1)$ and $x\in B_{\m M_X}(x_0,\tau)$. Consequently, $\wt v_{\omega_0}\in  {\m H}^{\alpha_Y,\alpha_X}_{L_2}(\mb B_{\mb R^{d_Y}}(\mathbf{0},\tau_1),\mb B_{\m M_X}(x_0,\tau))$.

 \subsection{Proof of Lemma~\ref{legenerative}}\label{prooflegenerative}

 For any $\omega^*=(x^*,y^*)$, let $V_{\omega^*}$ be an arbitrary orthonormal basis of $T_{\m M_{Y|x^*}} y^*$, and denote $G_{\omega^*}(z,x)=\Phi_{\omega^*}(V_{\omega^*}z,x)$, $Q_{\omega^*}(y)=V_{\omega^*}^T(y-y^*)$ and $U^{\omega^*}_{Y|x}=U_{\omega^*}\cap \m M_{Y|x}$. Then since $\m M_{Y|x}$ has a reach no smaller than $\tau$, by  by Lemma 2 of~\cite{10.1214/18-AOS1685}, it holds that 
 $\mb B_{\m M_{Y|x}}(y^*,\frac{7\tau_1}{8}\wedge \frac{7\tau}{16})\subset U^{\omega^*}_{Y|x}$. Moreover, according to Lemma~\ref{le:defdistribution}, the density of the push forward measure $[Q_{\omega^*}(\cdot)]_{\#}(\mu^*_{Y|x}|_{U^{\omega^*}_{Y|x}})$, denoted as  $  v_{\omega^*}(z|x)$, satisfies that $v_{\omega^*}(z,|,x)\in  {\m H}^{\alpha_Y,\alpha_X}_{L_1}(\mb B_{\mb R^{d_Y}}(\mathbf{0},\tau_1),B_{\m M_X}(x_0,\tau) )$. Based on the  aforementioned facts, we can first show that
\begin{claim}\label{claim4}
 For any $(x^*,y^*)\in \m M$ and $x\in \mb B_{\m M_X}(x^*,\tau)$, it holds for any $r\leq \frac{7\tau_1}{8}\wedge \frac{7\tau}{16}$  and any measurable function $g:\m M_{Y|x}\to \mb R$ that 
         \begin{equation*}
             \mb{E}_{ \mu^*_{Y|x}}[g(Y)\cdot\mathbf{1}(Y\in \mb B_{\m M_{Y|x}}(y^*,r))]=\int_{\mb B_{\mb R^{d_Y}}(\mathbf{0},\tau_1)} g(G_{\omega^*}(z,x))\mathbf{1}(G_{\omega^*}(z,x)\in \mb B_{\m M_{Y|x}}(y^*,r))) v_{\omega^*}(z|x)\,\dd z.
         \end{equation*}
\end{claim}
 \noindent Indeed, denote ${\rm vol}_{\m M_{Y|x}}$ as the volume measure of $\m M_{Y|x}$, and  let $u(y\,|\,x)$ be the density of $\mu^*_{Y|x}$ with respect to ${\rm vol}_{\m M_{Y|x}}$. We can obtain that
    \begin{equation*}
    \begin{aligned}
                  & \mb{E}_{\mu^*_{Y|x}}[g(Y)\cdot\mathbf{1}(Y\in \mb B_{\m M_{Y|x}}(y^*,r))]=\int g(y)\cdot\mathbf{1}(y\in \mb B_{\m M_{Y|x}}(y^*,r)) u(y\,|\,x)\,\dd  {\rm vol}_{\m M_{Y|x}}(y)\\
                   &=\int_{U^{\omega^*}_{Y|x}} g(y)\cdot\mathbf{1}(y\in \mb B_{\m M_{Y|x}}(y^*,r)) u(y\,|\,x)\,\dd  {\rm vol}_{\m M_{Y|x}}(y)\\
                   &=\int_{\mb B_{\mb R^{d_Y}}(\mathbf{0},\tau_1)} g(G_{\omega^*}(z,x))\mathbf{1}(G_{\omega^*}(z,x)\in \mb B_{\m M_{Y|x}}(y^*,r)) u(G_{\omega^*}(z,x)|x)\sqrt{{\rm det}(J_{G_{\omega^*}(\cdot,x)}(z)^TJ_{G_{\omega^*}(\cdot,x)}(z))}\,\dd z\\
                   &=\int_{\mb B_{\mb R^{d_Y}}(\mathbf{0},\tau_1)} g(G_{\omega^*}(z,x))\mathbf{1}(G_{\omega^*}(z,x)\in \mb B_{\m M_{Y|x}}(y^*,r))) v_{\omega^*}(z|x)\,\dd z,
    \end{aligned}
   \end{equation*}
   which proves Claim~\ref{claim4}.
 Then let $\{\omega_k^*=(x_k^*,y_k^*)\}_{k=1}^{K^*}\subset \m M$ be a $\tau_2$-covering set of $\m M$, and consider a smooth transition function $\rho:\mb R\to [0,1]$  that satisfies $\rho(t)=1$ when $t\in[0,1]$ and $\rho(t)=0$ when $t\in [2,\infty)$ (for example, the function defined in~\eqref{eqnrho}). For $k\in [K^*]$, define
         \begin{equation*}
             \wt\rho_\sk(x,y)=\rho(\frac{\|x-x_k^*\|^2}{\tau_2^2})\rho(\frac{\|y-y_k^*\|^2}{\tau_2^2}).
         \end{equation*}
and 
\begin{equation*}
    \rho_\sk(x,y)=\frac{\wt\rho_\sk(x,y)}{\kappa\Big(\sum_{k'=1}^{K^*} \wt\rho_{[k']}(x,y)\Big)},
\end{equation*}
with
\begin{equation*}
    \kappa(t)=t(1-\rho(2t))+\frac{\rho(2t)}{2}.
\end{equation*}
We can verify that $\kappa(t)\geq 1/2$ holds for any $t>0$ and $\kappa(t)=t$ if $t\geq 1$. Consequently, $\rho_\sk$ is a smooth function defined over the entire space of $\mb R^{D_X}\times \mb R^{D_Y}$. Additionally, for any $(x,y)\in \m M$,   there exists $k'\in [K^*]$ so that $\|(x_{k'}^*,y_{k'}^*)-(x,y)\|\leq \tau_2$. Consequently, $\sum_{k=1}^{K^*}\wt\rho_\sk(x,y)\geq \wt \rho_{[k']}(x,y)\geq 1$. Therefore,  when $(x,y)\in \m M$, it holds that $ \rho_\sk(x,y)=\wt\rho_\sk(x,y)/\sum_{k'=1}^{K^*} \wt\rho_{[k']}(x,y)$  and $\sum_{k=1}^{K^*} \rho_\sk(x,y)=1$.
Furthermore, given that for any $k\in [K^*]$,  $G_{\omega^*_k}(z,x)\in \m H^{\beta_Y,\beta_X}_L(\mb B_{\mb R^{d_Y}}(\mathbf{0},\tau_1), B_{\m M_X}(x_0,\tau) )$ and $v_{\omega_k^*}(z,|,x)\in {\m H}^{\alpha_Y,\alpha_X}_{L_1}(\mb B_{\mb R^{d_Y}}(\mathbf{0},\tau_1),B_{\m M_X}(x_0,\tau) )$, there exist $G_\sk^*\in \m H^{\beta_Y,\beta_X}_L( \mb R^{d_Y} , \mb R^{D_X} )$ and  $ \wt\nu_\sk^*\in {\m H}^{\alpha_Y,\alpha_X}_{L_1}(\mb R^{d_Y} ,\mb R^{D_X})$ such that for any $z\in \mb B_{\mb R^{d_Y}}(\mathbf{0},\tau_1)$ and $x\in B_{\m M_X}(x_0,\tau)$, $G_\sk^*(z,x)=G_{\omega_k^*}(z,x)$ and $ \wt\nu_\sk^*(z|x)=\nu_{\omega_k^*}(z|x)$.  Then based on $\sqrt{2}\tau_2\leq \frac{\sqrt{2}}{4}(\tau\wedge \tau_1)<\frac{7\tau_1}{8}\wedge \frac{7\tau}{16}$ and Claim~\ref{claim4}, we have 
   \begin{equation*}
       \begin{aligned}
            \mb{E}_{\mu^*_{Y|x}}[g(Y)]&= \sum_{k=1}^{K^*} \mb{E}_{ \mu^*_{Y|x}}[g(Y)\rho_\sk(x,Y)]\\
            &=\sum_{k=1}^{K^*} \mb{E}_{\mu^*_{Y|x}}[g(Y)  \rho_\sk(x,Y)\cdot\mathbf{1}(Y\in \mb B_{\m M_{Y|x}}(y_k,2\tau_2))]\\
            &=\sum_{k=1}^{K^*} \int_{\mb B_{\mb R^{d_Y}}(\mathbf{0},\tau_1)} g(G^*_{[k]}(z,x)) \rho_\sk(x,G^*_{[k]}(z,x)) \mathbf{1}(G^*_{[k]}(z,x)\in \mb B_{\m M_{Y|x}}(y_k,\sqrt{2}\tau_2))  \wt v_{[k]}^*(z|x)\,\dd z\\
            &=\sum_{k=1}^{K^*} \int_{\mb B_{\mb R^{d_Y}}(\mathbf{0},\tau_1)} g(G^*_{[k]}(z,x)) \rho_\sk(x,G^*_{[k]}(z,x))  \wt v_{[k]}^*(z|x)\,\dd z.\\
       \end{aligned}
   \end{equation*}
Then let $v_\sk^*(z,x)= \rho_\sk(x,G^*_{[k]}(z,x))\wt v_{[k]}^*(z|x)$, we can get the desired result.
 \subsection{Proof of Lemma~\ref{lemmalowerboundsmooth}}\label{proof:lemmalowerboundsmooth}
   The proof  uses a similar argument as in the proof of Lemma~\ref{le:appro} (see Appendix~\ref{appendixE.3}). For any $(l_1,l_2)\in \m J^{d_1,d_2}_{\alpha_1,\alpha_2}=\{l_1\in \mb N_0^{d_1}, l_2\in \mb N_0^{d_2}:\,\frac{|l_1|}{\alpha_1}+\frac{|l_2|}{\alpha_2}<1\}$, since $\phi_1$,$\phi_2$ are smooth compactly supported, we have 
\begin{equation*}
    \begin{aligned}
      &  |f^{(l_1,l_2)}(x,y)|=\left|\frac{m_1^{|l_1|}m_2^{|l_2|}}{(m_1)^{\alpha_1}}\sum_{\xi_1\in [m_1]^{d_1}}\sum_{\xi_2\in [m_2]^{d_2}}  \omega_{\xi_1,\xi_2}\phi_1^{(l_1)}(m_1x-\xi_1)\phi_2^{(l_2)}(m_2y-\xi_2)\right|\\
      &\leq L \frac{m_1^{|l_1|}m_2^{|l_2|}}{(m_1)^{\alpha_1}}\leq L_1 \frac{m_1^{|l_1|}m_1^{\alpha_1|l_2|/\alpha_2}}{(m_1)^{\alpha_1}}\leq L_1.
    \end{aligned}
\end{equation*}
Furthermore, by employing a similar approach to that used in the proof of Claim~\ref{claim1}, and considering the relationship \( |\omega_{\xi_1,\xi_2}| \lesssim m_1^{-\alpha_1} \asymp m_2^{-\alpha_2} \), we can demonstrate that for any \(x, x' \in \mathbb{R}^{d_1}\), \(y \in \mathbb{R}^{d_2}\), and any \((l_1, l_2) \in \mathcal{J}^{d_1, d_2}_{\alpha_1, \alpha_2}\) with \( \frac{|l_1|}{\alpha_1} + \frac{|l_2|}{\alpha_2} + \frac{1}{\alpha_1} \geq 1\), it holds that  
    \begin{equation*}
\begin{aligned}
        &\Big|\frac{m_1^{|l_1|}m_2^{|l_2|}}{(m_1)^{\alpha_1}}\sum_{\xi_1\in [m_1]^{d_1}}\sum_{\xi_2\in [m_2]^{d_2}}  \omega_{\xi_1,\xi_2}\phi_1^{(l_1)}(m_1x-\xi_1)\phi_2^{(l_2)}(m_2y-\xi_2)\\
        &\quad-\frac{m_1^{|l_1|}m_2^{|l_2|}}{(m_1)^{\alpha_1}}\sum_{\xi_1\in [m_1]^{d_1}}\sum_{\xi_2\in [m_2]^{d_2}}  \omega_{\xi_1,\xi_2}\phi_1^{(l_1)}(m_1x'-\xi_1)\phi_2^{(l_2)}(m_2y-\xi_2)\Big|\\
        &\leq L_1\|x-x'\|^{\alpha_1-|l_1|-\frac{\alpha_1}{\alpha_2}|l_2|}.
\end{aligned}
\end{equation*}
Moreover, for any \((l_1, l_2) \in \mathcal{J}^{d_1, d_2}_{\alpha_1, \alpha_2}\) with \( \frac{|l_1|}{\alpha_1} + \frac{|l_2|}{\alpha_2} + \frac{1}{\alpha_2} \geq 1\), and for any \(x \in \mathbb{R}^{d_1}\), \(y, y' \in \mathbb{R}^{d_2}\),
    \begin{equation*}
\begin{aligned}
         &\Big|\frac{m_1^{|l_1|}m_2^{|l_2|}}{(m_1)^{\alpha_1}}\sum_{\xi_1\in [m_1]^{d_1}}\sum_{\xi_2\in [m_2]^{d_2}}  \omega_{\xi_1,\xi_2}\phi_1^{(l_1)}(m_1x-\xi_1)\phi_2^{(l_2)}(m_2y-\xi_2)\\
        &\quad-\frac{m_1^{|l_1|}m_2^{|l_2|}}{(m_1)^{\alpha_1}}\sum_{\xi_1\in [m_1]^{d_1}}\sum_{\xi_2\in [m_2]^{d_2}}  \omega_{\xi_1,\xi_2}\phi_1^{(l_1)}(m_1x-\xi_1)\phi_2^{(l_2)}(m_2y'-\xi_2)\Big|\\
        &\leq L_1\|y-y'\|^{\alpha_2-|l_2|-\frac{\alpha_2}{\alpha_1}|l_1|}.
\end{aligned}
\end{equation*}
We can then conclude that there exists a constant $L_1$ so that $f\in  {\m H}^{\alpha_1,\alpha_2}_{L_1}(\mb R^{d_1},\mb R^{d_2})$.

  \subsection{Proof of Theorem~\ref{theoremjointregression}}\label{proof:theoremjointregression}
 
Denote the loss function
\begin{equation*}
    \begin{aligned}
       \ell(x,y,S)=\sum_{\lambda\in \Lambda} S(\lambda,x)^2-2\psi_{\lambda}(y) S(\lambda,x).
    \end{aligned}
\end{equation*}
Then we have 
\begin{equation}\label{eqn:1}
   \begin{aligned}
           \wh S&=  \underset{S\in \m S}{\arg\min}\frac{1}{n} \sum_{i=1}^n\sum_{\lambda\in \Lambda}(S(\lambda,X_i)-\psi_{\lambda}(Y_i))^2\\
           &= \underset{S\in \m S}{\arg\min}\frac{1}{n} \sum_{i=1}^n\sum_{\lambda\in \Lambda}S(\lambda,X_i)^2-2\psi_{\lambda}(Y_i)S(\lambda,X_i) \\
           &=\underset{S\in \m S}{\arg\min}\frac{1}{n} \sum_{i=1}^n \ell(X_i,Y_i,S).
   \end{aligned}
\end{equation}
 Denote $\mu^*=\mu^*_X\mu^*_{Y|X}$ as the joint distribution of $(X,Y)$. We have 
\begin{equation}\label{eqn:2}
    \begin{aligned}
        &\mb{E}_{\mu^*}[\ell(X,Y,S)]\\
        &=\mb{E}_{\mu^*}\big[\sum_{\lambda\in \Lambda}\big(S(\lambda,X)^2-2\psi_{\lambda}(Y)S(\lambda,X) \big) \big]\\
        &=\mb{E}_{\mu^*_X}\big[\sum_{\lambda\in \Lambda} S(\lambda,X)^2\big]-2\cdot\mb{E}_{\mu^*_X}\Big[\sum_{\lambda\in \Lambda} \mb{E}_{\mu^*_{Y|X}}\big[\psi_{\lambda}(Y)\big]S(\lambda,X)\Big]\\
        &\quad+ \mb{E}_{\mu^*_X}\Big[ \sum_{\lambda\in \Lambda}\Big(\mb{E}_{\mu^*_{Y|X}}\big[\psi_{\lambda}(Y)\big]\Big)^2\Big]-\mb{E}_{\mu^*_X}\Big[ \sum_{\lambda\in \Lambda}\Big(\mb{E}_{\mu^*_{Y|X}}\big[\psi_{\lambda}(Y)\big]\Big)^2\Big]\\
&=\mb{E}_{\mu^*_X}\Big[\sum_{\lambda\in \Lambda}  \Big(S(\lambda,X)-\mb{E}_{\mu^*_{Y|X}}\big[\psi_{\lambda}(Y)\big]\Big)^2 \Big]-\mb{E}_{\mu^*_X}\Big[ \sum_{\lambda\in \Lambda}\Big(\mb{E}_{\mu^*_{Y|X}}\big[\psi_{\lambda}(Y)\big]\Big)^2\Big].\\
    \end{aligned}
\end{equation}
Furthermore,  for $\rho=\max\{\sup_{(x,y)\in \m M}\sup_{S\in \m S}|\ell(x,y,S)|,1\}$, it holds that 
\begin{equation*}
\begin{aligned}
     \rho&\leq \max\{\sup_{(x,y)\in \m M}\sup_{S\in \m S}\sum_{\lambda\in \Lambda}\big(S(\lambda,x)^2+2|\psi_{\lambda}(y)S(\psi,x)|\big),1\}\leq  \max\{2C,1\}.
\end{aligned}
\end{equation*}
Then let
\begin{equation*}
\begin{aligned}
    S^*&\in {\arg\min}_{S\in \m S}\mb{E}_{\mu^*} [\ell(X,Y,S)]= {\arg\min}_{S\in \m S}\mb{E}_{\mu^*_X}\Big[\sum_{\lambda\in \Lambda}  \Big(S(\lambda,X)-\mb{E}_{\mu^*_{Y|X}}\big[\psi_{\lambda}(Y)\big]\Big)^2 \Big].
    \end{aligned}
\end{equation*}
Consider the function class
\begin{equation*}
\ms G^*=\{g(x,y)=\ell(x,y,S)-\ell(x,y,S^*):\, S\in \m S\}
\end{equation*}
and the star hull 
\begin{equation*}
   \ov{\ms G}^*=\{g(x,y)=a(\ell(x,y,S)-\ell(x,y,S^*)):\, a\in [0,1], S\in \m S\}.
\end{equation*}
Define the local Rademacher complexity
\begin{equation*}
 \m R_n(\ov{\ms G}^*,r)=\mathbb{E}_{\mu^{*,\otimes n}} \left[\underset{g \in \overline{\ms G}^* \atop  \mb{E}_{\mu^*}[g^2] \leq r^2}{\sup} \left|\frac{1}{n}\sum_{i=1}^n g(X_i,Y_i)-\mathbb{E}_{\mu^*}\left[g(x,y)\right]\right| \right],
 \end{equation*}
 where we uses the notation $\mb{E}_{\mu^*}[g^2]=\mb{E}_{\mu^*}[g(X,Y)^2]$ for simplicity.  We claim that the critical radius associated with $\ov{\ms G}^*$ is $\delta_{n}=c_1  \sqrt{\frac{W_n(\log n+\log T_n)}{n}}$ for a large enough $c_1$. This implies that 
 \begin{equation}\label{crir}
 \begin{aligned}
  \overline{R}_n(\overline{\ms G}^*,\delta_{n}) \leq \delta_{n}^2. \\
 \end{aligned}
 \end{equation}
 The Claim~\eqref{crir} will be proved later. Then define
\begin{equation*}
    M_{n}(S)=\frac{1}{n}\sum_{i=1}^n\ell(X_i,Y_i,S) \text{ and }    M^*(S)= \mb{E}_{\mu^*}[\ell(X,Y,S)].
\end{equation*}
Utilizing the uniform law (see for example, Theorem 14.20 of~\cite{wainwright2019high}) in conjunction with the aforementioned Claim~\eqref{crir}, we can get that,  there exists a constant $C_1$  so that it holds with probability larger than $1-n^{-c}$ that

  \begin{equation}\label{eqn:3}
  \begin{aligned}
  & \forall S\in \m S,\quad\frac{\left|M_{n}(S)-M_{n}(S^*)-M^*(S)+M^*(S^*)\right|}{\delta_{n}+\sqrt{\mb{E}_{\mu^*}[(\ell(x,y,S)-\ell(x,y,S^*))^2]}}  \leq  C_1 \delta_{n}.
 \end{aligned}
  \end{equation}

\noindent By the assumption that for any $S,S'\in \m S$,
\begin{equation*}
           \mb{E}_{\mu^*}\Big[\big(  \ell(X,Y,S)-  \ell(X,Y,S')\big)\Big]\leq C\, \mb{E}_{\mu^*_X}\Big[\sum_{\lambda\in \Lambda} \big(S(\lambda,X)-S'(\lambda,X)\big)^2 \Big],
       \end{equation*}
 we can get
\begin{equation*}
\begin{aligned}
     &\mb{E}_{\mu^*}[(\ell(X,Y,S)-\ell(X,Y,S^*))^2]\leq  C\, \mb{E}_{\mu^*_X}\Big[\sum_{\lambda\in \Lambda} \big(\wh S(\lambda,X)-S^*(\lambda,X)\big)^2\Big]\\
     &\leq 2C\, \mb{E}_{\mu^*_X}\Big[\sum_{\lambda\in \Lambda} \big(\wh S(\lambda,X)-\mb{E}_{\mu^*_{Y|X}}\big[\psi_{\lambda}(Y)\big]\big)^2 \Big] + 2C\, \mb{E}_{\mu^*_X}\Big[\sum_{\lambda\in \Lambda} \big( S^*(\lambda,X)-\mb{E}_{\mu^*_{Y|X}}\big[\psi_{\lambda}(Y)\big]\big)^2\Big].
\end{aligned}
\end{equation*}
Then, combined with~\eqref{eqn:1},~\eqref{eqn:2},~\eqref{eqn:3}, we can get
\begin{equation*}
\begin{aligned}
      &  \mb{E}_{\mu^*_X}\Big[\sum_{\lambda\in \Lambda} \big(\wh S(\lambda,X)-\mb{E}_{\mu^*_{Y|X}}\big[\psi_{\lambda}(Y)\big]\big)^2 \Big]-\underset{S\in \m S}{\min}\, \mb{E}_{\mu^*_X}\Big[\sum_{\lambda\in \Lambda} \big( S(\lambda,X)-\mb{E}_{\mu^*_{Y|X}}\big[\psi_{\lambda}(Y)\big]\big)^2\Big] \\
          & =  \mb{E}_{\mu^*_X}\Big[\sum_{\lambda\in \Lambda} \big(\wh S(\lambda,X)-\mb{E}_{\mu^*_{Y|X}}\big[\psi_{\lambda}(Y)\big]\big)^2 \Big]- \mb{E}_{\mu^*_X}\Big[\sum_{\lambda\in \Lambda} \big( S^*(\lambda,X)-\mb{E}_{\mu^*_{Y|X}}\big[\psi_{\lambda}(Y)\big]\big)^2\Big]\\
      &=  M^*(\wh S)-M^*(S^*)\leq M^*(\wh S)-M^*(S^*)+M_{n}(S^*)-M_{n}( S)\\
       &\leq C_1 \delta_{n}^2+C_1\sqrt{2C}\delta_{n}\cdot\sqrt{  \mb{E}_{\mu^*_X}\Big[\sum_{\lambda\in \Lambda} \big(\wh S(\lambda,X)-\mb{E}_{\mu^*_{Y|X}}\big[\psi_{\lambda}(Y)\big]\big)^2 \Big]
+ \underset{S\in \m S}{\min}\, \mb{E}_{\mu^*_X}\Big[\sum_{\lambda\in \Lambda} \big( S(\lambda,X)-\mb{E}_{\mu^*_{Y|X}}\big[\psi_{\lambda}(Y)\big]\big)^2\Big]}.
\end{aligned}
\end{equation*}
So by combining all pieces, we can get that it holds with probability at least $1-n^{-c}$ that
\begin{equation*}
\begin{aligned}
&\mb{E}_{\mu^*_X}\Big[\sum_{\lambda\in \Lambda} \big(\wh S(\lambda,X)-\mb{E}_{\mu^*_{Y|X}}\big[\psi_{\lambda}(Y)\big]\big)^2 \Big] \\
& \leq C_2\Big( \delta_{n}^2+ \underset{S\in \m S}{\min}\, \mb{E}_{\mu^*_X}\Big[\sum_{\lambda\in \Lambda} \big( S(\lambda,X)-\mb{E}_{\mu^*_{Y|X}}\big[\psi_{\lambda}(Y)\big]\big)^2\Big]\Big)\Big)\\
&\leq C_2\Big( \frac{W_n(\log n+\log T_n)}{n}+ \underset{S\in \m S}{\min}\, \mb{E}_{\mu^*_X}\Big[\sum_{\lambda\in \Lambda} \big( S(\lambda,X)-\mb{E}_{\mu^*_{Y|X}}\big[\psi_{\lambda}(Y)\big]\big)^2\Big]\Big)\Big).
\end{aligned}
\end{equation*}
 Now it only remains to show Claim~\eqref{crir}. Using standard symmetrization, we can get for any $r>0$,
\begin{equation*}
 \m R_n(\ov{\ms G}^*,r)=\mathbb{E}_{\mu^{*,\otimes n}} \left[\underset{g \in \overline{\ms G}^* \atop\mb{E}_{\mu^*}[g^2] \leq r^2}{\sup} \left|\frac{1}{n}\sum_{i=1}^n g(X_i,Y_i)-\mathbb{E}_{\mu^*}\left[g(x,y)\right]\right| \right]\leq \mathbb{E}_{\mu^{*,\otimes n}}\mathbb{E}_\epsilon \left[\underset{g \in  \overline{\ms G}^* \atop \mb{E}_{\mu^*}[g^2] \leq r^2}{\sup}\left|\frac{2}{n}\sum_{i=1}^n\epsilon_i g(X_i,Y_i)\right|\right],
 \end{equation*}
  where \(\left\{\epsilon_{i}\right\}_{i=1}^{n}\) are \(n\) i.i.d. copies from Rademacher distribution, i.e. \(\mb P\left(\epsilon_{i}=1\right)=\mb P\left(\epsilon_{i}=-1\right)=\frac{1}{2}\).
  
 \quad\\
Define $d_n^g(g , g' )=\sqrt{\frac{1}{n}\sum_{i=1}^n(g(X_i,Y_i)-g'(X_i,Y_i))^2}$,
then 
\begin{equation*}
r_{n}:=\underset{g  ,g' \in \overline{\ms G}^* \atop \mb{E}_{\mu^*}[g^2],\mb{E}_{\mu^*}[(g')^2]\leq r^2}{\max} d_n(g, g')\leq 2\rho,
\end{equation*}
and by equation (3.84) of~\cite{wainwright2019high}, there exists a constant $C_3$ such that,
\begin{equation*}
\begin{aligned}
\mathbb{E}_{\mu^{*,\otimes n}}[r_{n}^2]&\leq  \mathbb{E}_{\mu^{*,\otimes n}}\left[\underset{g_ \in\overline{\ms G}^* \atop \mb{E}_{\mu^*}[g^2] \leq r^2}{\sup} \frac{4}{n} \sum_{i=1}^n g^2 (X_i,Y_i)\right]\\
&\leq  \mathbb{E}_{\mu^{*,\otimes n}}\left[\underset{g \in\overline{\ms G}^* \atop \mb{E}_{\mu^*}[g^2] \leq r^2}{\sup} \frac{8}{n}\sum_{i=1}^n \left(g (X_i,Y_i)-\mathbb{E}_{\mu^*}[g(x,y)]\right)^2 \right]+8r^2\\
&\leq C_3(r^2+\rho\m {R}_n(\overline{\ms G}^*,r)).
\end{aligned}
\end{equation*}
Moreover, for any  $g \in \ms G^*$ and $a \in (0,1]$,  there exists an integer $\kappa\in \mathbb{N}$, such that $\kappa\frac{\varepsilon}{2\rho}<a\leq (\kappa+1)\frac{\varepsilon}{2\rho} 
$ and $d_n^g((\kappa+1)\frac{\varepsilon}{2\rho}g ,ag)\leq \frac{\varepsilon}{2\rho}\rho=\frac{\varepsilon}{2}$. Therefore it follows that the $\varepsilon$-covering number of $\overline{\ms G}^*$ with respect to $d_n^g$ satisfies that, $\mathbf{N}(\overline{\ms G}^*,d_n^g,\varepsilon)\leq \mathbf{N}(\ms G^*,d_n^g,\frac{\varepsilon}{2})\cdot\frac{2\rho}{\varepsilon}$. Therefore, we can obtain  for any $0<\varepsilon\leq r_n$
\begin{equation*}
    \begin{aligned}
        \log \mathbf{N}(\overline{\ms G}^* ,d_n^g,\varepsilon)&\leq \log \mathbf{N}(\ms G^* ,d_n^g,\frac{\varepsilon}{2})+ \log 
\frac{2\rho}{\varepsilon} =\log \mathbf{N}(\m S ,d_n,\frac{\varepsilon}{2})+ \log 
\frac{2\rho}{\varepsilon} \\
&\leq W_n\log \frac{2T_n}{\varepsilon} + \log 
\frac{2\rho}{\varepsilon}\leq  W_n\log \frac{4T_n\rho}{\varepsilon^2}.
    \end{aligned}
\end{equation*}
Then, by Dudley entropy integral bound~\cite{wainwright2019high, vershynin2018high}, we have 
\begin{equation*}
\begin{aligned}
        &\overline{R}_n(\overline{\ms G}^*,r)\leq \frac{C_4}{\sqrt{n}}\mb{E}_{\mu^{*,\otimes n}}\left[\int_{0}^{r_{n}}\sqrt{ W_n\log \frac{4T_n\rho}{\varepsilon^2}}\,\dd \varepsilon\right]\\
&=\frac{C_4}{\sqrt{n}}\mb{E}_{\mu^{*,\otimes n}}\left[r_{n}\int_{0}^{1}\sqrt{ W_n\log \frac{4T_n\rho}{\varepsilon^2r_n^2} }\,\dd \varepsilon\right]\\
&\leq \frac{C_4}{\sqrt{n}}\mb{E}_{\mu^{*,\otimes n}}\left[r_{n}\int_{0}^{1}\sqrt{ W_n\log \frac{T_n}{\varepsilon^2\rho} }\,\dd \varepsilon\right]+\frac{\sqrt{2}C_4}{\sqrt{n}}\mb{E}_{\mu^{*,\otimes n}}\left[r_{n}\int_{0}^{1}\sqrt{ W_n\log \frac{2\rho}{r_n} }\,\dd \varepsilon\right]\\
&\leq C_4\Big(\log(T_n)+\int_0^1 \sqrt{2\log  \frac{1}{\varepsilon}}\,\dd \varepsilon\Big)\sqrt{\frac{W_n}{n}}\mb{E}_{\mu^{*,\otimes n}}[r_{n}]+2C_4\rho\sqrt{\frac{W_n}{n}} \mb{E}_{\mu^{*,\otimes n}}\left[\sqrt{(\frac{r_n}{2\rho})^2-(\frac{r_n}{2\rho})^2\log((\frac{r_n}{2\rho})^2)}\right]\\
&\leq C_5\sqrt{\frac{W_n}{n}} \sqrt{-  \mathbb{E}_{\mu^{*,\otimes n}}\left[r_{n}^2\right] \log \mathbb{E}_{\mu^{*,\otimes n}}\left(\frac{r_{n}}{2 \rho}\right)^2 + \mathbb{E}_{\mu^{*,\otimes n}}[r_{n}^2]\cdot\log (T_n)},
\end{aligned}
\end{equation*}
where the last inequality uses that $\sqrt{- y\log y+ y}$ is concave and non-decreasing when $y=(\frac{r_{n}}{2\rho})^2\leq 1$. Then  by $\mathbb{E}_{\mu^{*,\otimes n}}[r_{n}^2]\leq C_3(r^2+\rho\overline{R}_n(r,\overline{\ms G}^*))$,  there exists some constant $C_6$ so that
\begin{equation*}
\begin{aligned}
    \overline{R}_n(\overline{\ms G}^*,r) &\leq C_6\, \sqrt{\frac{W_n}{n}} (r^2+\rho\overline{R}_n(\overline{\ms G}^*,r))^{\frac{1}{2}}\sqrt{\log \frac{1}{r}+\log T_n}.
\end{aligned}
 \end{equation*}
 Choose $\delta_{n}=c_1 \sqrt{\frac{W_n(\log n+\log T_n)}{n}}$ with $c_1>1$. If $\overline{R}_n(\delta_{n},\overline{\ms G}^*)> \delta_{n}^2$,  then 
 \begin{equation*}
 \begin{aligned}
 \overline{R}_n(\overline{\ms G}^*, \delta_{n}) \leq C_6\sqrt{\frac{W_n}{n}} \sqrt{2(1+C)} \overline{R}_n(\overline{\ms G}^*, \delta_{n})^{\frac{1}{2} }\sqrt{\log n+\log T_n}
  \end{aligned}
 \end{equation*}
 which means
 \begin{equation*}
 \begin{aligned}
  \overline{R}_n( \overline{\ms G}^*,\delta_n)  \leq 2(1+C)C_6^2 \frac{W_n}{n} (\log n+\log T_n) \leq \frac{2(1+C)C_6^2}{c_1^2}\delta_{n}^2. \\
 \end{aligned}
 \end{equation*}
So $ \overline{R}_n(\delta_{n},\overline{\ms G}^*)  \leq \delta_{n}^2$ holds if $c_1>\sqrt{2(1+C)C_6^2}\vee 1$. This completes the proof.

\subsection{Proof of Theorem~\ref{th:mainifold}}\label{prooftheorem1}

The proof of the lower bound is derived directly from the proof of the lower bound in Theorem~\ref{th2:lower} as detailed in Appendix~\ref{proofth2lower}. Specifically, consider the construction of the submanifolds  described in Appendix~\ref{proofth2lower}. For any \(j, k \in [H]\) with \(j \neq k\), it is established that:
\begin{equation*}
    \begin{aligned}
       & \underset{x\in \m M_X}{\sup} \mb H(\m M^{(j)}_{Y|x},\m M^{(k)}_{Y|x})=\underset{x\in \m M_X}{\sup} \mb H(\m M^{\omega^{(j)}}_{Y|x},\m M^{\omega^{(k)}}_{Y|x})\\
       &\geq  \underset{x\in \m M_X}{\sup}\underset{z\in \mb B_{\mb R^{d_Y}}(\mathbf{0},1)}{\sup}\|g_{\omega^{(j)}}(z,x)-g_{\omega^{(k)}}(z,x)\|\\
       &\gtrsim \frac{1}{m_1^{\beta_Y}}\asymp n^{\frac{1}{\frac{d_Y}{\beta_Y}+\frac{d_X}{\beta_X}}}.
    \end{aligned}
\end{equation*}
The desired result then follows in a manner similar to that outlined in Appendix~\ref{proofth2lower}, utilizing Fano's lemma. Now we will show the upper bound. We begin with the construction of the estimator.
  Consider the data points $\big\{(X_i,Y_i)\big\}_{i=1}^n$, for each $k\in [n]$, we define the local polynomial estimator $(\wh V_k, (\wh a_{j_1j_2k})_{(j_1,j_2)\in \m J_{\beta_Y,\beta_X}^{d_Y,D_X}})$ at $(X_k,Y_k)$ to be any element of 
    \begin{equation*}
    \begin{aligned}
             \underset{V\in \mb O(D_Y,d_Y)\atop \underset{(j_1,j_2)\in \m J_{\beta_Y,\beta_X}^{d_Y,D_X}}{\sup} |a_{j_1,j_2}|\leq L_1}  {\arg\min}\frac{1}{n}&\sum_{i=1}^n \|Y_i-\sum_{(j_1,j_2)\in \m J_{\beta_Y,\beta_X}^{d_Y,D_X}} \frac{a_{j_1j_2}}{j_1!j_2!}(V^T(Y_i-Y_k))^{j_1}(X_i-X_k)^{j_2}\|^2\\
             &\cdot \textbf{1}(Y_i\in \mb B_{\mb R^{D_Y}}(Y_k,h_1))\textbf{1}(X_i\in \mb B_{\mb R^{D_X}}(X_k,h_2)),
    \end{aligned}
    \end{equation*}
    where $h_1=b_1\,(\frac{\log n}{n})^{\frac{1}{d_Y+\frac{d_X\beta_Y}{\beta_X}}}$, $h_2=b_2\,(\frac{\log n}{n})^{\frac{1}{d_X+\frac{d_Y\beta_X}{\beta_Y}}}$ and $b_1$, $b_2$  are large enough constants. Then for any $x\in \m M_X$, consider the estimator $\wh {\m M}_{Y|x}$ of $\m M_{Y|x}$ defined as
\begin{equation*}
    \wh {\m M}_{Y|x}=\bigcup_{i\in[n]\atop \|X_i-x\|\leq h_2} \Big\{y=\sum_{(j_1,j_2)\in \m J_{\beta_Y,\beta_X}^{d_Y,d_X}}   \frac{1}{j_1!j_2!} \wh a_{j_1j_2k}z^{j_1}(x-X_i)^{j_2}\,:\, z\in \mb B_{\mb R^{d_Y}}(\mathbf{0},h_1)\Big\}.
\end{equation*}
  We will show  $\wh {\m M}_{Y|x}$  can achieve the upper bound in Theorem~\ref{th:mainifold}.  Let $V_{k}\in \mb R^{D_Y\times d_Y}$ be a matrix whose column forms an orthonormal basis of $T_{\m M_{Y|X_k}} Y_k$. Consider the function $G^*_\sk(z,x)=\Phi_{(X_k,Y_k)}(V_{k}^*z,x)$, where $\Phi_{(X_k,Y_k)}$ is the one defined in Definition~\ref{def:manifold} of the main text. It holds with a constant $L$ that  $G^*_\sk(z,x)\in  {\m H}^{\beta_Y,\beta_X}_{L,D_Y}(\mb B_{\mb R^{d_Y}}(\mathbf{0},\tau_1),\mb B_{\m M_X}(X_k,\tau))$. 
Moreover,  notice that $\m M_X$ is a $\beta_X$-smooth manifold, let $W_k^*$  be a matrix whose column forms an orthornormal basis of $T_{\m M_X} X_k$ and define $g_\sk(s)=\phi_{X_k}(W_k^*s)$, where $\phi_{X_k}$ is the one defined in Definition~\ref{def:manifold1}.  Denote $\wt G_\sk(z,s)=G_{\sk}^*(z,g_\sk(s))$, it holds that 
\begin{equation*}
    \begin{aligned}
          \| G_{\sk}^*(z,g_\sk(s))-\sum_{(j_1,j_2)\in \m J_{d_Y,d_X}^{\beta_Y,\beta_X}}\frac{1}{j_1!j_2!}\wt G_\sk{}^{(j_1,j_2)}(\mathbf{0}_{d_Y},\mathbf{0}_{d_X})z^{j_1}s^{j_2}\|\lesssim 
      \|z\|^{\beta_Y}+\|s\|^{\beta_X}.
    \end{aligned}
\end{equation*}
Denote 
\begin{equation*}
    \wh f_\sk(z,x)=\sum_{(j_1,j_2)\in \m J_{\beta_Y,\beta_X}^{d_Y,d_X}}   \frac{1}{j_1!j_2!} \wh a_{j_1j_2k}(\wh V_k^T(G_{\sk}^*(z,x)-Y_k))^{j_1}(x-X_k)^{j_2},
    \end{equation*}
and $\wt f_\sk(z,s)=\wh f_\sk(z,g_\sk(s))$.
Then 
 \begin{equation*}
      \| \wh f_{\sk}(z,g_\sk(s))-\sum_{(j_1,j_2)\in \m J_{d_Y,d_X}^{\beta_Y,\beta_X}}\frac{1}{j_1!j_2!}\wt f_{\sk}^{(j_1,j_2)}(\mathbf{0}_{d_Y},\mathbf{0}_{d_X})z^{j_1}s^{j_2}\|\lesssim \|z\|^{\beta_Y}+\|s\|^{\beta_X}.
    \end{equation*}
Therefore, denote $Z_{ik}=V_k^*{}^T(Y_i-Y_k)$ and $S_{ik}=W_k^*{}^T(X_i-X_k)$, we have 
    \begin{equation*}
        \begin{aligned}
         &\frac{1}{n}\sum_{i=1}^n \Big\|\sum_{(j_1,j_2)\in \m J_{d_Y,d_X}^{\beta_Y,\beta_X}}\frac{1}{j_1!j_2!}\wt G_{\sk}{}^{(j_1,j_2)}(\mathbf{0}_{d_Y},\mathbf{0}_{d_X})Z_{ik}^{j_1}S_{ik}^{j_2}-\sum_{(j_1,j_2)\in \m J_{d_Y,D_X}^{\beta_Y,\beta_X}}\frac{1}{j_1!j_2!}\wt f_{\sk}^{(j_1,j_2)}(\mathbf{0}_{d_Y},\mathbf{0}_{d_X})Z_{ik}^{j_1}S_{ik}^{j_2}\Big\|^2\\
         &\qquad\cdot \mathbf{1}(Y_i\in \mb B_{\mb R^{D_Y}}(Y_k,h_1))\textbf{1}(X_i\in \mb B_{\mb R^{D_X}}(X_k,h_2))\\
            &\leq \frac{3}{n}\sum_{i=1}^n \|\sum_{(j_1,j_2)\in \m J_{d_Y,d_X}^{\beta_Y,\beta_X}}\frac{1}{j_1!j_2!}\wt G_{\sk}{}^{(j_1,j_2)}(\mathbf{0}_{d_Y},\mathbf{0}_{d_X})Z_{ik}^{j_1}S_{ik}^{j_2}- G^*_\sk(Z_{ik},g(S_{ik}))\|^2\\
           & \qquad\cdot \mathbf{1}(Z_{ik}\in \mb B_{\mb R^{d_Y}}(\mathbf{0},h_1))\mathbf{1}(S_{ik}\in \mb B_{\mb R^{d_X}}(\mathbf{0},h_2))\\
            &+\frac{3}{n}\sum_{i=1}^n \|Y_i-\sum_{(j_1,j_2)\in \m J_{\beta_Y,\beta_X}^{d_Y,D_X}} \wh a_{j_1j_2k}(\wh V_k^T(Y_i-Y_k))^{j_1}(X_i-X_k)^{j_2}\|^2\cdot \textbf{1}(Y_i\in \mb B_{\mb R^{D_Y}}(Y_k,h_1))\mathbf{1}(X_i\in \mb B_{\mb R^{D_X}}(X_k,h_2))\\
             &+ \frac{3}{n}\sum_{i=1}^n \|\wh f_\sk(Z_{ik},g(S_{ik}))-\sum_{(j_1,j_2)\in \m J_{d_Y,D_X}^{\beta_Y,\beta_X}}\frac{1}{j_1!j_2!}\wt f_{\sk}^{(j_1,j_2)}(\mathbf{0}_{d_Y},\mathbf{0}_{d_X})Z_{ik}^{j_1}S_{ik}^{j_2}\|^2\\
             &\qquad\cdot \textbf{1}(Z_{ik}\in \mb B_{\mb R^{d_Y}}(\mathbf{0},h_1))\textbf{1}(S_{ik}\in \mb B_{\mb R^{d_X}}(\mathbf{0},h_2))\\
&\leq \frac{3}{n}\sum_{i=1}^n \|G_{\sk}^*(V_{k}^*{}^T(Y_i-Y_k),X_i)-\sum_{(j_1,j_2)\in \m J_{\beta_Y,\beta_X}^{d_Y,D_X}} \frac{1}{j_1!j_2!}G_{\sk}^*{}^{(j_1,j_2)}(\mathbf{0}_{d_Y},X_k)(V_k^*{}^T(Y_i-Y_k))^{j_1}(X_i-X_k)^{j_2}\|^2\\
&\quad\qquad\cdot \textbf{1}(Y_i\in \mb B_{\mb R^{D_Y}}(Y_k,h_1))\textbf{1}(X_i\in \mb B_{\mb R^{D_X}}(X_k,h_2))\\
&\quad+ C\,(h_1^{2\beta_Y}+h_2^{2\beta_X})\cdot \frac{1}{n}\sum_{i=1}^n \textbf{1}(Z_{ik}\in \mb B_{\mb R^{d_Y}}(\mathbf{0},h_1))\textbf{1}(S_{ik}\in \mb B_{\mb R^{d_X}}(\mathbf{0},h_2))\\
&\lesssim  (h_1^{2\beta_Y}+h_2^{2\beta_X})\cdot \frac{1}{n}\sum_{i=1}^n \textbf{1}(Z_{ik}\in \mb B_{\mb R^{d_Y}}(\mathbf{0},h_1))\textbf{1}(S_{ik}\in \mb B_{\mb R^{d_X}}(\mathbf{0},h_2)).
        \end{aligned}
    \end{equation*}
Building on the analysis presented in~\cite{10.1214/23-AOS2291,10.1214/18-AOS1685}, we can derive the following lemma, whose proof is given in Section~\ref{prooflemma13}.
    \begin{lemma}\label{lemmamanifold1}
  For any positive constant $c$, there exists a constant $C$ so that  it holds with probability at least $1-n^{-c}$ that for any $k\in [n]$,
    \begin{equation*}
        \begin{aligned}
           & \mb{E}_{\mu^*}\bigg[\Big\|\sum_{(j_1,j_2)\in \m J_{d_Y,d_X}^{\beta_Y,\beta_X}}\frac{1}{j_1!j_2!}\wt G_{\sk}{}^{(j_1,j_2)}(\mathbf{0}_{d_Y},\mathbf{0}_{d_X})(V_k^*{}^T(Y-Y_k))^{j_1}(W_k^*{}^T(X-X_k))^{j_2}\\
           &-\sum_{(j_1,j_2)\in \m J_{d_Y,D_X}^{\beta_Y,\beta_X}}\frac{1}{j_1!j_2!}\wt f_{\sk}^{(j_1,j_2)}(\mathbf{0}_{d_Y},\mathbf{0}_{d_X})(V_k^*{}^T(Y-Y_k))^{j_1}(W_k^*{}^T(X-X_k))^{j_2}\Big\|^2 \\
           &\qquad\qquad\cdot \mathbf{1}(Y\in \mb B_{\mb R^{D_Y}}(Y_k,h_1))\mathbf{1}(X\in \mb B_{\mb R^{D_X}}(X_k,h_2))\bigg]\\
           &\leq C\,\sqrt{\frac{\log n}{n}} h_1^{\frac{d_Y}{2}} h_2^{\frac{d_X}{2}}   \cdot\Big[\sum_{(j_1,j_2)\in \m J_{d_Y,d_X}^{\beta_Y,\beta_X}}\frac{1}{j_1!j_2!}\Big\|\wt G_{\sk}{}^{(j_1,j_2)}(\mathbf{0}_{d_Y},\mathbf{0}_{d_X})-\wt f_{\sk}^{(j_1,j_2)}(\mathbf{0}_{d_Y},\mathbf{0}_{d_X})\Big\|h_1^{j_1}h_2^{j_2}\Big]^2\\
           &+C\,(h_1^{2\beta_Y}+h_2^{2\beta_X}) (h_1^{d_Y}h_2^{d_X}+\frac{\log n}{n}).\\
        \end{aligned}
    \end{equation*}
    \end{lemma}

    \medskip
   \noindent On the other hand, notice that there exists a small enough constant $c$ so that for any $z\in \mb B_{\mb R^{d_Y}}(\mathbf{0},c\,h_1)$ and $s\in \mb B_{\mb R^{d_X}}(\mathbf{0},c\, h_2)$,  it holds that 
    \begin{equation*}
       \| g_\sk(s)-X_k\|=\|g_\sk(s)-g_\sk(\mathbf 0)\|\leq h_2,
    \end{equation*}
    \begin{equation*}
        \|G_\sk^*(z,g_\sk(s))-Y_k\|=   \|G_\sk^*(z,g_\sk(s))-G_\sk^*(\mathbf 0,g_\sk(\mathbf 0))\|\leq \frac{h_1+h_2^{\beta_X\wedge 1}}{2}\leq h_1.
    \end{equation*}
Therefore, we can obtain the following lower bound 
    \begin{equation*}
        \begin{aligned}
              & \mb{E}_{\mu^*}\bigg[\Big\|\sum_{(j_1,j_2)\in \m J_{d_Y,d_X}^{\beta_Y,\beta_X}}\frac{1}{j_1!j_2!}\wt G_{\sk}{}^{(j_1,j_2)}(\mathbf{0}_{d_Y},\mathbf{0}_{d_X})(V_k^*{}^T(Y-Y_k))^{j_1}(W_k^*{}^T(X-X_k))^{j_2}\\
           &-\sum_{(j_1,j_2)\in \m J_{d_Y,D_X}^{\beta_Y,\beta_X}}\frac{1}{j_1!j_2!}\wt f_{\sk}^{(j_1,j_2)}(\mathbf{0}_{d_Y},\mathbf{0}_{d_X})(V_k^*{}^T(Y-Y_k))^{j_1}(W_k^*{}^T(X-X_k))^{j_2}\Big\|^2 \\
            &\qquad\qquad\cdot \mathbf{1}(Y\in \mb B_{\mb R^{D_Y}}(Y_k,h_1))\mathbf{1}(X\in \mb B_{\mb R^{D_X}}(X_k,h_2))\bigg]\\
           &\geq    
           \mb{E}_{\mu^*}\bigg[\Big\|\sum_{(j_1,j_2)\in \m J_{d_Y,d_X}^{\beta_Y,\beta_X}}\frac{1}{j_1!j_2!}\wt G_{\sk}{}^{(j_1,j_2)}(\mathbf{0}_{d_Y},\mathbf{0}_{d_X})(V_k^*{}^T(Y-Y_k))^{j_1}(W_k^*{}^T(X-X_k))^{j_2}\\
           &-\sum_{(j_1,j_2)\in \m J_{d_Y,D_X}^{\beta_Y,\beta_X}}\frac{1}{j_1!j_2!}\wt f_{\sk}^{(j_1,j_2)}(\mathbf{0}_{d_Y},\mathbf{0}_{d_X})(V_k^*{}^T(Y-Y_k))^{j_1}(W_k^*{}^T(X-X_k))^{j_2}\Big\|^2 \\
            &\underbrace{\qquad\qquad\cdot \mathbf{1}(V_k^*{}^T(Y-Y_k)\in \mb B_{\mb R^{d_Y}}(\mathbf{0},c\,h_1))\mathbf{1}(W_k^*{}^T(X-X_k)\in \mb B_{\mb R^{d_X}}(\mathbf{0},c\,h_2))\bigg]}_{(E_A)}.
        \end{aligned}
    \end{equation*}
    The term  $(E_A)$ can be further lower bounded by 
    \begin{equation*}
        \begin{aligned}
         &(E_A)= \int_{z\in \mb B_{\mb R^{d_Y}}(\mathbf{0},c\,h_1) }\int_{s\in\mb B_{\mb R^{d_X}}(\mathbf{0},c\,h_2)} \Big\|\sum_{(j_1,j_2)\in \m J_{d_Y,d_X}^{\beta_Y,\beta_X}}\frac{1}{j_1!j_2!}\Big(\wt G_{\sk}{}^{(j_1,j_2)}(\mathbf{0}_{d_Y},\mathbf{0}_{d_X})-\wt f_{\sk}^{(j_1,j_2)}(\mathbf{0}_{d_Y},\mathbf{0}_{d_X})\Big)z^{j_1}s^{j_2}\Big\|^2 \\
            & \cdot \mu^*_X(g_\sk(s))\mu^*_{Y|g_\sk(s)}(\wt G_\sk(z,s)) \sqrt{{\rm det}\Big(J_{\wt G_\sk(\cdot,s)}(z)^T J_{\wt G_\sk(\cdot,s)}(z)\Big)} \sqrt{{\rm det}\Big(J_{g_\sk}(s)^T J_{g_\sk}(s)\Big)}\,\dd z\dd s\\
           &\gtrsim h_1^{d_Y}h_2^{d_X} \\
           &\cdot\int_{ \mb B_{\mb R^{d_Y}}(\mathbf{0},1) }\int_{\mb B_{\mb R^{d_X}}(\mathbf{0},1)} \Big\|\sum_{(j_1,j_2)\in \m J_{d_Y,d_X}^{\beta_Y,\beta_X}}\frac{1}{j_1!j_2!}\Big(\wt G_{\sk}{}^{(j_1,j_2)}(\mathbf{0}_{d_Y},\mathbf{0}_{d_X})-\wt f_{\sk}^{(j_1,j_2)}(\mathbf{0}_{d_Y},\mathbf{0}_{d_X})\Big)h_1^{j_1}h_2^{j_2}z^{j_1}s^{j_2}\Big\|^2   \dd z\dd s\\
           &\gtrsim  h_1^{d_Y}h_2^{d_X} \Big(\sum_{(j_1,j_2)\in \m J_{d_Y,d_X}^{\beta_Y,\beta_X}}\frac{1}{j_1!j_2!}\Big\|\wt G_{\sk}{}^{(j_1,j_2)}(\mathbf{0}_{d_Y},\mathbf{0}_{d_X})-\wt f_{\sk}^{(j_1,j_2)}(\mathbf{0}_{d_Y},\mathbf{0}_{d_X})\Big\| h_1^{j_1}h_2^{j_2}  \Big)^2,
        \end{aligned}
    \end{equation*}
where the last inequality uses the fact that for any $d$-variate polynomial $\mathcal{S}(y)=\sum_{j\in \mb N_0^d,\,|j|\leq k} a_{j} y^{j}$, $y\in\mb R^d$, there exists some positive constant $C(d,k)$ only depending on $(d,k)$ such that 
\begin{equation*}
\int_{\mb B_1^d} \mathcal{S}^2(y) \, \dd y \geq C(d,k) \sum_{j\in \mb N_0^d,\,|j|\leq k} a_{j}^2.
  \end{equation*}
Therefore, combined with Lemma~\ref{lemmamanifold1}, when  $h_1=b_1\,(\frac{\log n}{n})^{\frac{1}{d_Y+\frac{d_X\beta_Y}{\beta_X}}}$, $h_2=b_2\,(\frac{\log n}{n})^{\frac{1}{d_X+\frac{d_Y\beta_X}{\beta_Y}}}$  with sufficiently large $b_1,b_2$, we can obtain that 
    \begin{equation}\label{eqndiffGF}
        \sum_{(j_1,j_2)\in \m J_{d_Y,d_X}^{\beta_Y,\beta_X}}\frac{1}{j_1!j_2!}\Big\|\wt G_{\sk}{}^{(j_1,j_2)}(\mathbf{0}_{d_Y},\mathbf{0}_{d_X})-\wt f_{\sk}^{(j_1,j_2)}(\mathbf{0}_{d_Y},\mathbf{0}_{d_X})\Big\| h_1^{j_1}h_2^{j_2} \lesssim (\frac{\log n}{n})^{\frac{1}{\frac{d_Y}{\beta_Y}+\frac{d_X}{\beta_X}}}.
    \end{equation}
In order to show $\wh {\m M}_{Y|x}$ satisfies the desired result, we will also use the following lemma, whose proof is provided in Appendix~\ref{prooflemma14}.
\begin{lemma}\label{lemmamanifold2}
    It holds with probability at least $1-n^{-1}$ that for any $(x,y)\in \m M$, there exists $i\in [n]$ so that $\|y-Y_i\|< h_1$ and $\|x-X_i\|< h_2$. 
\end{lemma}
\noindent Using Lemma~\ref{lemmamanifold2} and inequality~\eqref{eqndiffGF}, for any $x\in \m M_X$ and $y\in \m M_{Y|x}$, there exists $k\in [n]$ so that $\|y-Y_k\|\leq h_1$, $\|x-X_k\|\leq h_2$, and
\begin{equation*}
    \begin{aligned}
    y&=\Phi_{(X_k,Y_k)}(V_k^*V_k^*{}^T(y-Y_k),x)=G_\sk^*(V_k^*{}^T(y-Y_k),x)=\wt G_\sk(V_k^*{}^T(y-Y_k), W_k^*{}^T(x-X_k))\\
    &=\sum_{(j_1,j_2)\in \m J_{d_Y,d_X}^{\beta_Y,\beta_X}}\frac{1}{j_1!j_2!}\wt G_{\sk}{}^{(j_1,j_2)}(\mathbf{0}_{d_Y},\mathbf{0}_{d_X})(V_k^*{}^T(y-Y_k))^{j_1}(W_k^*{}^T(x-X_k))^{j_2}+\m O\big( h_1^{\beta_X}+h_2^{\beta_Y}\big)\\
    &=\sum_{(j_1,j_2)\in \m J_{d_Y,d_X}^{\beta_Y,\beta_X}}\frac{1}{j_1!j_2!}\wt f_{\sk}{}^{(j_1,j_2)}(\mathbf{0}_{d_Y},\mathbf{0}_{d_X})(V_k^*{}^T(y-Y_k))^{j_1}(W_k^*{}^T(x-X_k))^{j_2}+\m O\big(h_1^{\beta_X}+h_2^{\beta_Y}\big)\\
    &=\wh f_\sk(V_k^*{}^T(y-Y_k),x)+\m O\big(h_1^{\beta_X}+h_2^{\beta_Y}\big)\\
    &=\sum_{(j_1,j_2)\in \m J_{\beta_Y,\beta_X}^{d_Y,d_X}}   \frac{1}{j_1!j_2!} \wh a_{j_1j_2k}(\wh V_k^T(y-Y_k))^{j_1}(x-X_k)^{j_2}+\m O\big(h_1^{\beta_X}+h_2^{\beta_Y}\big).
      \end{aligned}
\end{equation*}
    Moreover, we have $\sum_{(j_1,j_2)\in \m J_{\beta_Y,\beta_X}^{d_Y,d_X}}   \frac{1}{j_1!j_2!} \wh a_{j_1j_2k}(\wh V_k^T(y-Y_k))^{j_1}(x-X_k)^{j_2}\in \wh{\m M}_{Y|x}$,
    Thus 
    \begin{equation*}
        \underset{y\in \m M_{Y|x}}{\sup}\underset{y'\in \wh{\m M}_{Y|x}}{\inf}\|y-y'\|\lesssim (\frac{\log n}{n})^{\frac{1}{\frac{d_Y}{\beta_Y}+\frac{d_X}{\beta_X}}}.
    \end{equation*}
On the other side, for a fixed $x\in \m M_X$, consider any $k\in [n]$ with $\|X_k-x\|\leq h_2$ and $z\in \mb B_{\mb R^{d_Y}}(\mathbf{0},h_1)$. Then 
    \begin{equation*}
        \begin{aligned}
            &\sum_{(j_1,j_2)\in \m J_{\beta_Y,\beta_X}^{d_Y,d_X}}   \frac{1}{j_1!j_2!} \wh a_{j_1j_2k}z^{j_1}(x-X_k)^{j_2}=\wh f_\sk(z,g_\sk(W_k^*{}^T(x-X_k)))\\
            &=\sum_{(j_1,j_2)\in \m J_{d_Y,d_X}^{\beta_Y,\beta_X}}\frac{1}{j_1!j_2!}\wt f_{\sk}{}^{(j_1,j_2)}(\mathbf{0}_{d_Y},\mathbf{0}_{d_X})z^{j_1}(W_k^*{}^T(x-X_k))^{j_2}+\m O\big(h_1^{\beta_X}+h_2^{\beta_Y}\big)\\
            &=\sum_{(j_1,j_2)\in \m J_{d_Y,d_X}^{\beta_Y,\beta_X}}\frac{1}{j_1!j_2!}\wt G_{\sk}{}^{(j_1,j_2)}(\mathbf{0}_{d_Y},\mathbf{0}_{d_X})z^{j_1}(W_k^*{}^T(x-X_k))^{j_2}+\m O\big(h_1^{\beta_X}+h_2^{\beta_Y}\big)\\
            &=\wt G_\sk(z,W_k^*{}^T(x-X_k))+\m O\big(h_1^{\beta_X}+h_2^{\beta_Y}\big)\\
            &=G_\sk^*(z,x)+\m O\big(h_1^{\beta_X}+h_2^{\beta_Y}\big).\\
        \end{aligned}
    \end{equation*}
 Then since $G_\sk^*(z,x)\in \m M_{Y|x}$,  we have 
     \begin{equation*}
       \underset{y'\in \wh{\m M}_{Y|x}}{\sup} \underset{y\in \m M_{Y|x}}{\inf}\|y-y'\|\lesssim (\frac{\log n}{n})^{\frac{1}{\frac{d_Y}{\beta_Y}+\frac{d_X}{\beta_X}}}.
    \end{equation*}
 Therefore,  it holds with probability at least $1-2\,n^{-1}$ that $$\underset{x\in \m M_X}{\sup} \mb H(\wh {\m M}_{Y|x},\m M_{Y|x})\lesssim (\frac{\log n}{n})^{\frac{1}{\frac{d_X}{\beta_X}+\frac{d_Y}{\beta_Y}}},$$ which can lead to $$\mb{E}_{\mu^{*,\otimes n}}\big[\underset{x\in \m M_X}{\sup} \mb H(\wh {\m M}_{Y|x},\m M_{Y|x})\big]\lesssim  (\frac{\log n}{n})^{\frac{1}{\frac{d_X}{\beta_X}+\frac{d_Y}{\beta_Y}}}.$$

\subsubsection{Proof of Lemma~\ref{lemmamanifold1}}\label{prooflemma13}
The proof directly follows~\cite{10.1214/23-AOS2291}, we include it here for completeness.  Since there exists a constant $C_0$ so that $\|\wt{f}_\sk^{(j_1,j_2)}(\mathbf{0}_{d_Y},\mathbf{0}_{d_X})\|_2\leq C_0$ holds for any possible $k,j_1,j_2$. For any fixed $k\in [n]$ and $\widetilde\delta>0$, let 
\begin{equation*}
\begin{aligned}
&\bar{\mathcal{T}}(\widetilde{\delta})=\Big\{T=\{T_{j_1,j_2}\}_{(j_1,j_2)\in \m J^{\beta_Y,\beta_X}_{d_Y,D_X}}\in [-C_0,\,C_0]^{D\times |\m J^{\beta_Y,\beta_X}_{d_Y,D_X}|}\, :\, \\
&\qquad\qquad\sum_{(j_1,j_2)\in \m J^{\beta_Y,\beta_X}_{d_Y,D_X}}\frac{1}{j_1!j_2!}\, \big\|T_{j_1,j_2}-\wt G_\sk^{(j_1,j_2)}(\mathbf{0}_{d_Y},\mathbf{0}_{d_X})\big\|_2\, h_1^{|j_1|}h_2^{|j_2|}\leq \widetilde{\delta}\Big\}.
\end{aligned}
\end{equation*}
We also define the following supreme of an empirical process indexed by $T\in \bar{\mathcal{T}}(\widetilde{\delta})$,
\begin{equation*}
 \begin{aligned}
&\quad Z_n(\widetilde{\delta})=\\
 &\underset{T\in \bar{\mathcal{T}}(\widetilde{\delta}) }{\sup}\Bigg| \, \mathbb{E}_{\mu^{\ast}} \bigg[\Big\|\sum_{(j_1,j_2)\in \m J^{\beta_Y,\beta_X}_{d_Y,D_X}} \frac{1}{j_1!j_2!}\,\big(\wt G_\sk^{(j_1,j_2)}(\mathbf{0}_{d_Y},\mathbf{0}_{d_X}) -T_{j_1,j_2}\big)\,(V_k^*{}^T(Y-Y_k))^{j_1}(W_k^*{}^T(X-X_k))^{j_2} \Big\|_2^2\\
 &\cdot    \mathbf{1}(Y\in \mb B_{\mb R^{D_Y}}(Y_k,h_1))\textbf{1}(X\in \mb B_{\mb R^{D_X}}(X_k,h_2))\bigg]\\
& -n^{-1}\sum_{i\in [n]\atop i\neq k}\bigg[\Big\|\sum_{(j_1,j_2)\in \m J^{\beta_Y,\beta_X}_{d_Y,D_X}} \frac{1}{j_1!j_2!}\,\big(\wt G_\sk^{(j_1,j_2)}(\mathbf{0}_{d_Y},\mathbf{0}_{d_X}) -T_{j_1,j_2}\big)\,(V_k^*{}^T(Y_i-Y_k))^{j_1}(W_k^*{}^T(X_i-X_k))^{j_2} \Big\|_2^2\\
 &\cdot  \mathbf{1}(Y_i\in \mb B_{\mb R^{D_Y}}(Y_k,h_1))\textbf{1}(X_i\in \mb B_{\mb R^{D_X}}(X_k,h_2))\bigg]\Bigg|,
  \end{aligned}
\end{equation*}
and $R_n(\widetilde{\delta})=\mathbb{E}_{\mu^\ast{}^{\otimes n}}\big[ Z_n(\widetilde{\delta}) \big]$. We will first prove a concentration inequality for a fixed radius $\widetilde{\delta}>0$, and then using the peeling technique to allow the radius to be random, which leads to the desired result.

To apply the Talagrand concentration inequality (see, for example, Theorem 3.27 of~\cite{wainwright2019high}) for bounding the difference $|Z_n(\widetilde{\delta}) - R_n(\widetilde \delta)|$ for a fixed $\widetilde{\delta}>0$, we notice that each additive component in the second empirical sum above has second moment uniformly bounded by
\begin{equation*}
 \begin{aligned}
& \mathbb{E}_{\mu^{\ast}} \bigg[\underset{T\in \bar{\mathcal{T}}(\widetilde{\delta}) }{\sup}\Big( \Big\|\sum_{(j_1,j_2)\in \m J^{\beta_Y,\beta_X}_{d_Y,D_X}} \frac{1}{j_1!j_2!}\,\big(\wt G_\sk^{(j_1,j_2)}(\mathbf{0}_{d_Y},\mathbf{0}_{d_X}) -T_{j_1,j_2}\big)\,(V_k^*{}^T(Y-Y_k))^{j_1}(W_k^*{}^T(X-X_k))^{j_2} \Big\|_2^4\\
 & \qquad\qquad \cdot \mathbf{1}(Y\in \mb B_{\mb R^{D_Y}}(Y_k,h_1))\textbf{1}(X\in \mb B_{\mb R^{D_X}}(X_k,h_2))\Big)\bigg]\\
&\leq \underset{z\in \mb B_{\mb R^{d_Y}}(\mathbf{0},h_1), s\in \mb B_{\mb R^{d_X}}(\mathbf{0},h_2)\atop T\in \bar{\mathcal{T}}(\widetilde{\delta})} {\sup} \Big\|\sum_{(j_1,j_2)\in \m J^{\beta_Y,\beta_X}_{d_Y,D_X}} \frac{1}{j_1!j_2!}\,\big(\wt G_\sk^{(j_1,j_2)}(\mathbf{0}_{d_Y},\mathbf{0}_{d_X}) -T_{j_1,j_2}\big)\,z^{j_1}s^{j_2}  \Big\|_2^4 \\
&\cdot \mathbb{P}_{\mu^{\ast}}\big(Y\in \mb B_{\mb R^{D_Y}}(Y_k,h_1),\,X\in \mb B_{\mb R^{D_X}}(X_k,h_2)\big)\\
&\leq C \underset{T\in \bar{\mathcal{T}}(\widetilde{\delta}) }{\sup}\bigg(\sum_{(j_1,j_2)\in \m J^{\beta_Y,\beta_X}_{d_Y,D_X}} \frac{1}{j_1!j_2!}\,\big(\wt G_\sk^{(j_1,j_2)}(\mathbf{0}_{d_Y},\mathbf{0}_{d_X}) -T_{j_1,j_2}\big)\,h_1^{|j_1|}h_2^{|j_2|}\bigg)^4\cdot  h_1^{d_Y}h_2^{d_X}\\
&\leq C\, \wt \delta^4 h_1^{d_Y}h_2^{d_X}.\\
  \end{aligned}
\end{equation*}
Moreover, each additive component can be almost surely bounded by 
\begin{equation*}
 \begin{aligned}
& \quad \underset{z\in \mb B_{\mb R^{d_Y}}(\mathbf{0},h_1), s\in \mb B_{\mb R^{d_X}}(\mathbf{0},h_2)\atop T\in \bar{\mathcal{T}}(\widetilde{\delta})} {\sup}\Big\|\sum_{(j_1,j_2)\in \m J^{\beta_Y,\beta_X}_{d_Y,D_X}} \frac{1}{j_1!j_2!}\,\big(\wt G_\sk^{(j_1,j_2)}(\mathbf{0}_{d_Y},\mathbf{0}_{d_X}) -T_{j_1,j_2}\big)\,z^{j_1}s^{j_2} \Big\|_2^2\\
&\leq C\, \underset{T\in \bar{\mathcal{T}}(\widetilde{\delta}) }{\sup} \bigg(\sum_{(j_1,j_2)\in \m J^{\beta_Y,\beta_X}_{d_Y,D_X}} \frac{1}{j_1!j_2!}\,\big\|\wt G_\sk^{(j_1,j_2)}(\mathbf{0}_{d_Y},\mathbf{0}_{d_X}) -T_{j_1,j_2}\big\|\,h_1^{|j_1|}h_2^{|j_2|}\bigg)^2
\leq C\, \widetilde \delta^2.
   \end{aligned}
\end{equation*}
Based on these two bounds, we can apply the Talagrand concentration inequality to obtain that for any $s\geq 0$,
\begin{equation}~\label{Talagrand}
\mathbb{P} \big(Z_n(\widetilde{\delta}) \geq R_n(\widetilde{\delta}) +s^2\big )\leq 2 \exp\left(-\frac{c\,ns^4}{s^2 \,\widetilde{\delta}^2 +  \wt \delta^4 h_1^{d_Y}h_2^{d_X}}\right).
   \end{equation}
It remains to bound the expectation $R_n(\widetilde{\delta})$ via the symmetrization technique and chaining. 
By a standard symmetrization, we can get 
\begin{equation*}
\begin{aligned}
&R_n(\widetilde{\delta})\leq \frac{2}{\sqrt{n}}\, \mathbb{E}\Bigg[\underset{T\in \bar{\mathcal{T}}(\widetilde{\delta})}{\sup}\\
& \Bigg|\frac{1}{\sqrt{n}} \sum_{i\in [n]\atop i\neq k} \varepsilon_i\bigg[\Big\|\sum_{(j_1,j_2)\in \m J^{\beta_Y,\beta_X}_{d_Y,D_X}} \frac{1}{j_1!j_2!}\,\big(\wt G_\sk^{(j_1,j_2)}(\mathbf{0}_{d_Y},\mathbf{0}_{d_X}) -T_{j_1,j_2}\big)\,(V_k^*{}^T(Y_i-Y_k))^{j_1}(W_k^*{}^T(X_i-X_k))^{j_2} \Big\|_2^2\\
 &\cdot  \mathbf{1}(Y_i\in \mb B_{\mb R^{D_Y}}(Y_k,h_1))\textbf{1}(X_i\in \mb B_{\mb R^{D_X}}(X_k,h_2))\bigg]\Bigg|\Bigg],
\end{aligned}
\end{equation*}
where $\{\varepsilon_i\}_{i=1}^n$ are $n$ i.i.d.~copies from the Rademacher distribution, i.e.~$\mb P(\varepsilon_i = 1) = \mb P(\varepsilon_i= -1) = 0.5$. 
Since given $\{X_i,Y_i\}_{i\in [n],i\neq k}$, the stochastic process inside the supreme is a sub-Gaussian process with intrinsic metric 
 \begin{equation*}
 \begin{aligned}
& \quad d_n^2(T,\,\widetilde{T})\\
&=\frac{1}{n} \sum_{i\in [n]\atop i\neq k}  \bigg(\Big\|\sum_{(j_1,j_2)\in \m J^{\beta_Y,\beta_X}_{d_Y,D_X}} \frac{1}{j_1!j_2!}\,\big(\wt G_\sk^{(j_1,j_2)}(\mathbf{0}_{d_Y},\mathbf{0}_{d_X}) - T_{j_1,j_2}\big)\,(V_k^*{}^T(Y_i-Y_k))^{j_1}(W_k^*{}^T(X_i-X_k))^{j_2} \Big\|_2^2 \\
&\qquad\qquad  - \Big\|\sum_{(j_1,j_2)\in \m J^{\beta_Y,\beta_X}_{d_Y,D_X}} \frac{1}{j_1!j_2!}\,\big(\wt G_\sk^{(j_1,j_2)}(\mathbf{0}_{d_Y},\mathbf{0}_{d_X}) -\wt T_{j_1,j_2}\big)\,(V_k^*{}^T(Y_i-Y_k))^{j_1}(W_k^*{}^T(X_i-X_k))^{j_2} \Big\|_2^2\bigg)^2 \\
&\qquad\qquad \cdot \mathbf{1}(Y_i\in \mb B_{\mb R^{D_Y}}(Y_k,h_1))\textbf{1}(X_i\in \mb B_{\mb R^{D_X}}(X_k,h_2))\\
 &\leq C\,\widetilde{\delta}^4\, \frac{1}{n} \sum_{i\in [n]\atop i\neq k}  \mathbf{1}(Y_i\in \mb B_{\mb R^{D_Y}}(Y_k,h_1))\textbf{1}(X_i\in \mb B_{\mb R^{D_X}}(X_k,h_2)),
 \end{aligned}
\end{equation*}
for any $T,\widetilde{T}\in\bar{\mathcal{T}}(\widetilde \delta)$, where the last step uses the definition of $\bar{\mathcal{T}}(\widetilde \delta)$. So we have
\begin{equation*}
\mathbb{E}_{\mu^\ast}\Big[ \underset{T,\widetilde{T}\in \bar{\mathcal{T}}(\delta)}{\sup} d_n^2(T,\widetilde{T})\Big]\leq C\, \widetilde{\delta}^4 \cdot h_1^{d_Y}h_2^{d_X}\quad \mbox{and}\quad
d_n(T,\widetilde{T}) \leq C \widetilde{\delta}\sum_{(j_1,j_2)\in \m J^{\beta_Y,\beta_X}_{d_Y,D_X}}\frac{1}{j_1!j_2!}\, \|T_{(j_1,j_2)}-\widetilde{T}_{(j_1,j_2)}\|_{2}\, h_1^{|j_1|}h_2^{|j_2|}.
 \end{equation*}
 Lastly, let $\m K_n(\delta)= \underset{T,\widetilde{T}\in \bar{\mathcal{T}}(\delta)}{\sup} d_n(T,\widetilde{T})$, by applying the standard chaining via Dudley's inequality, we can get 
\begin{equation}\label{dudley}
 \begin{aligned}
 R_n(\widetilde{\delta})&\leq C\, \frac{1}{\sqrt{n}}\, \mathbb{E}_{\mu^\ast} \Big[ \int_{0}^{ \m K_n(\widetilde\delta)}\sqrt{\log \frac{C_1 \widetilde{\delta}}{u}}\,\dd u\Big]\\
 &= C\, \frac{1}{\sqrt{n}}\, \mathbb{E}_{\mu^\ast} \Big[ \m K_n(\widetilde\delta)\cdot\int_{0}^{1}\sqrt{\log \frac{C_1 \widetilde{\delta}}{u\cdot \m K_n(\widetilde\delta)}}\,\dd u\Big]\\
 &= C\, \frac{1}{\sqrt{n}}\, \mathbb{E}_{\mu^\ast} \Big[   \m K_n(\widetilde\delta)\cdot \mathbf{1}( \m K_n(\widetilde\delta)\leq \wt \delta^2 h_1^{d_Y/2}h_2^{d_X/2})\int_{0}^{1}\sqrt{\log \frac{C_1 \widetilde{\delta}}{u \cdot\m K_n(\widetilde\delta)}}\,\dd u\Big]\\
 &+C\, \frac{1}{\sqrt{n}}\, \mathbb{E}_{\mu^\ast} \Big[   \m K_n(\widetilde\delta)\cdot \mathbf{1}( \m K_n(\widetilde\delta)> \wt \delta^2 h_1^{d_Y/2}h_2^{d_X/2})\int_{0}^{1}\sqrt{\log \frac{C_1 \widetilde{\delta}}{u \cdot\m K_n(\widetilde\delta)}}\,\dd u\Big]\\
 &\leq C_1\, h_1^{\frac{d_Y}{2}} h_2^{\frac{d_X}{2}} \cdot \sqrt{\frac{-\log (\widetilde\delta h_1 h_2)}{n}}\cdot \widetilde{\delta}^2,
    \end{aligned}
\end{equation}
where we have used the fact that the $u$-covering entropy of $\bar{\mathcal{T}}(\widetilde \delta)$ relative to metric $d_n$ is at most $C_2\log\frac{C_1\widetilde\delta}{u}$ for $u\in(0,C_1\widetilde\delta)$. By combining this with inequality~\eqref{Talagrand}, we obtain that for all $t\geq 1$,
\begin{align}\label{fix_delta_tala}
    \mathbb{P} \Big(Z_n(\widetilde{\delta}) \geq  C\, t^2\, h_1^{\frac{d_Y}{2}} h_2^{\frac{d_X}{2}} \cdot  \sqrt{\frac{-\log (\widetilde\delta h_1 h_2)}{n}}\wt\delta^2\Big )\leq 2 \exp\Big(-c\,t^2\, \log (n/\widetilde \delta)\Big).
\end{align}
 Finally, we apply the peeling technique to extend the above high probability bound on $Z_n(\widetilde\delta)$ to the random radius $\widetilde\delta=\sum_{(j_1,j_2)\in \m J^{\beta_Y,\beta_X}_{d_Y,D_X}} \frac{1}{j_1!j_2!}\,\big\|\wt G_\sk^{(j_1,j_2)}(\mathbf{0}_{d_Y},\mathbf{0}_{d_X}) -T_{j_1,j_2}\big\|\,h_1^{|j_1|}h_2^{|j_2|}$. Specifically, we first set the basic level $\bar{\delta} =h_1^{\beta_Y}+h_2^{\beta_X}$, and for $s=1,\cdots, S$ with $S\leq \log \frac{C}{\bar{\delta}}$, define sets
\begin{equation*}
\begin{aligned}
&\widetilde{\mathcal{T}}_0=\Big\{T=\{T_{(j_1,j_2)}\}_{(j_1,j_2)\in \m J^{\beta_Y,\beta_X}_{d_Y,D_X}}\in [-C_0(\log n)^2,\,C_0(\log n)^2]^{D\times |\m J^{\beta_Y,\beta_X}_{d_Y,D_X}|}\, :\, \\
&\qquad\qquad\sum_{(j_1,j_2)\in \m J^{\beta_Y,\beta_X}_{d_Y,D_X}}\frac{1}{j_1!j_2!}\, \big\|T_{j_1,j_2}-\wt G_\sk^{(j_1,j_2)}(\mathbf{0}_{d_Y},\mathbf{0}_{d_X})\big\|_2\, h_1^{|j_1|}h_2^{|j_2|}\leq \ov{\delta}\Big\};\\
&\widetilde{\mathcal{T}}_s=\Big\{T=\{T_{(j_1,j_2)}\}_{(j_1,j_2)\in \m J^{\beta_Y,\beta_X}_{d_Y,D_X}}\in [-C_0(\log n)^2,\,C_0(\log n)^2]^{D\times |\m J^{\beta_Y,\beta_X}_{d_Y,D_X}|}\, :\, \\
&\qquad\qquad 2^{s-1}\bar{\delta} \leq \sum_{(j_1,j_2)\in \m J^{\beta_Y,\beta_X}_{d_Y,D_X}}\frac{1}{j_1!j_2!}\, \big\|T_{j_1,j_2}-\wt G_\sk^{(j_1,j_2)}(\mathbf{0}_{d_Y},\mathbf{0}_{d_X})\big\|_2h_1^{|j_1|}h_2^{|j_2|}  \leq 2^s\bar{\delta}\Big\}.\\
\end{aligned}
\end{equation*}
 By applying inequality~\eqref{fix_delta_tala} to $\widetilde\delta = 2^s\bar\delta$ for $s\in[S]$ with sufficiently large constant $t>0$, we obtain that 
 \begin{equation*}
\mathbb{P}\left(Z_n(\bar{\delta})\geq C\, h_1^{\frac{d_Y}{2}}h_1^{\frac{d_X}{2}}\sqrt{\frac{\log n}{n}}\bar{\delta}^2\right) + \sum_{s=1}^S \mathbb{P}\left(Z_n(2^s \bar{\delta})\geq C\, h_1^{\frac{d_Y}{2}}h_1^{\frac{d_X}{2}}\sqrt{\frac{\log n}{n}}\, 4^s \bar{\delta}^2\right)\leq n^{-(c+1)}.
\end{equation*}
Note that for any $T\in \widetilde{\mathcal{T}}_s$ and any $s\in\{0\}\cup [S]$, the event $Z_n(2^s \bar{\delta})\leq C\,b_2^{\frac{d}{2}} \,\frac{\log n}{n}\, 4^s \bar{\delta}^2$ implies
 \begin{align*}
 & \underset{T\in\widetilde{\mathcal{T}}_s }{\sup}\Bigg| \, \mathbb{E}_{\mu^{\ast}} \bigg[\Big\|\sum_{(j_1,j_2)\in \m J^{\beta_Y,\beta_X}_{d_Y,D_X}} \frac{1}{j_1!j_2!}\,\big(\wt G_\sk^{(j_1,j_2)}(\mathbf{0}_{d_Y},\mathbf{0}_{d_X}) -T_{j_1,j_2}\big)\,(V_k^*{}^T(Y-Y_k))^{j_1}(W_k^*{}^T(X-X_k))^{j_2} \Big\|_2^2\\
 &\cdot   \mathbf{1}(Y\in \mb B_{\mb R^{D_Y}}(Y_k,h_1))\textbf{1}(X\in \mb B_{\mb R^{D_X}}(X_k,h_2))\bigg]\\
& -n^{-1}\sum_{i\in [n]\atop i\neq k}\bigg[\Big\|\sum_{(j_1,j_2)\in \m J^{\beta_Y,\beta_X}_{d_Y,D_X}} \frac{1}{j_1!j_2!}\,\big(\wt G_\sk^{(j_1,j_2)}(\mathbf{0}_{d_Y},\mathbf{0}_{d_X}) -T_{j_1,j_2}\big)\,(V_k^*{}^T(Y_i-Y_k))^{j_1}(W_k^*{}^T(X_i-X_k))^{j_2} \Big\|_2^2\\
 &\cdot  \mathbf{1}(Y_i\in \mb B_{\mb R^{D_Y}}(Y_k,h_1))\textbf{1}(X_i\in \mb B_{\mb R^{D_X}}(X_k,h_2))\bigg]\Bigg|\\
&\leq  c_1\, h_1^{\frac{d_Y}{2}}h_1^{\frac{d_X}{2}}\sqrt{\frac{\log n}{n}}\, \Bigg\{\bar\delta^2 + \bigg(\sum_{(j_1,j_2)\in \m J^{\beta_Y,\beta_X}_{d_Y,D_X}}\frac{1}{j_1!j_2!}\, \big\|T_{j_1,j_2}-\wt G_\sk^{(j_1,j_2)}(\mathbf{0}_{d_Y},\mathbf{0}_{d_X})\big\|_2h_1^{|j_1|}h_2^{|j_2|}\bigg)^2\Bigg\}.
 \end{align*}
Furthermore, 
 \begin{equation*}
     \begin{aligned}
      &   n^{-1}\sum_{i\in [n]\atop i\neq k}\bigg[\Big\|\sum_{(j_1,j_2)\in \m J^{\beta_Y,\beta_X}_{d_Y,D_X}} \frac{1}{j_1!j_2!}\,\big(\wt G_\sk^{(j_1,j_2)}(\mathbf{0}_{d_Y},\mathbf{0}_{d_X}) -\wt f_\sk^{(j_1,j_2)}(\mathbf{0}_{d_Y},\mathbf{0}_{d_X})\big)\,(V_k^*{}^T(Y_i-Y_k))^{j_1}(W_k^*{}^T(X_i-X_k))^{j_2} \Big\|_2^2\\
 &\cdot  \mathbf{1}(Y_i\in \mb B_{\mb R^{D_Y}}(Y_k,h_1))\textbf{1}(X_i\in \mb B_{\mb R^{D_X}}(X_k,h_2))\\
 &\lesssim \big(h_1^{2\beta_Y}+h_2^{2\beta_X}\big)\cdot \frac{1}{n}\sum_{i\in [n]\atop i\neq k} \textbf{1}(Z_{ik}\in \mb B_{\mb R^{d_Y}}(\mathbf{0},h_1))\textbf{1}(S_{ik}\in \mb B_{\mb R^{d_X}}(\mathbf{0},h_2)).
     \end{aligned}
 \end{equation*}
Then since 
 \begin{equation*}
 \begin{aligned}
         & \mb{E}_{\mu^*}[ \big(\textbf{1}(Z_{ik}\in \mb B_{\mb R^{d_Y}}(\mathbf{0},h_1))\textbf{1}(S_{ik}\in \mb B_{\mb R^{d_X}}(\mathbf{0},h_2))\big)^2]=\mathbb{P}_{\mu^{\ast}}\big(Y\in \mb B_{\mb R^{D_Y}}(Y_k,h_1),\,X\in \mb B_{\mb R^{D_X}}(X_k,h_2)\big)\leq C\, h_1^{d_Y}h_2^{d_X}.
 \end{aligned}
 \end{equation*}
By Bernstein's inequality, it holds with probability at least $1-n^{-c-1}$ that
 \begin{equation*}
     \begin{aligned}
         & \frac{1}{n}\sum_{i\in [n]\atop i\neq k} \textbf{1}(Z_{ik}\in \mb B_{\mb R^{d_Y}}(\mathbf{0},h_1))\textbf{1}(S_{ik}\in \mb B_{\mb R^{d_X}}(\mathbf{0},h_2))\\
         &\leq \frac{1}{n-1}\sum_{i\in [n]\atop i\neq k} \textbf{1}(Z_{ik}\in \mb B_{\mb R^{d_Y}}(\mathbf{0},h_1))\textbf{1}(S_{ik}\in \mb B_{\mb R^{d_X}}(\mathbf{0},h_2))\\
         &\leq \Big|\frac{1}{n-1}\sum_{i=1}^n \textbf{1}(Z_{ik}\in \mb B_{\mb R^{d_Y}}(\mathbf{0},h_1))\textbf{1}(S_{ik}\in \mb B_{\mb R^{d_X}}(\mathbf{0},h_2))-\mathbb{P}_{\mu^{\ast}}\big(Y\in \mb B_{\mb R^{D_Y}}(Y_k,h_1),\,X\in \mb B_{\mb R^{D_X}}(X_k,h_2)\big)\Big|\\
         &\qquad +\mathbb{P}_{\mu^{\ast}}\big(Y\in \mb B_{\mb R^{D_Y}}(Y_k,h_1),\,X\in \mb B_{\mb R^{D_X}}(X_k,h_2)\big).
     \end{aligned}
 \end{equation*}
Then using the fact that $\m M_{X}$ and $\m M_{Y|X}$ are smooth submaifolds with reach bounded away from zero, and $f_X$, $f_{Y|X}$ are uniformly lower bounded by a constant, using Lemma B.7 of~\cite{10.1214/18-AOS1685}, we can get 
 \begin{equation*}
    \frac{1}{n}\sum_{i\in [n]\atop i\neq k} \textbf{1}(Z_{ik}\in \mb B_{\mb R^{d_Y}}(\mathbf{0},h_1))\textbf{1}(S_{ik}\in \mb B_{\mb R^{d_X}}(\mathbf{0},h_2))\lesssim \frac{\log n}{n}+ h_1^{d_Y}h_2^{d_X}.
 \end{equation*}
 So by combining all pieces, we can get that it holds with probability at least $1-n^{-c-1}$ that 
\begin{align*}
 &  \, \mathbb{E}_{\mu^{\ast}} \bigg[\Big\|\sum_{(j_1,j_2)\in \m J^{\beta_Y,\beta_X}_{d_Y,D_X}} \frac{1}{j_1!j_2!}\,\big(\wt G_\sk^{(j_1,j_2)}(\mathbf{0}_{d_Y},\mathbf{0}_{d_X}) -\wt f_\sk^{(j_1,j_2)}(\mathbf{0}_{d_Y},\mathbf{0}_{d_X})\big)\,(V_k^*{}^T(Y-Y_k))^{j_1}(W_k^*{}^T(X-X_k))^{j_2} \Big\|_2^2 \\
 &\cdot   \mathbf{1}(Y\in \mb B_{\mb R^{D_Y}}(Y_k,h_1))\textbf{1}(X\in \mb B_{\mb R^{D_X}}(X_k,h_2))\bigg]\\
&\lesssim h_1^{\frac{d_Y}{2}}h_1^{\frac{d_X}{2}}\sqrt{\frac{\log n}{n}}\,  \bigg(\sum_{(j_1,j_2)\in \m J^{\beta_Y,\beta_X}_{d_Y,D_X}}\frac{1}{j_1!j_2!}\, \big\|\wt f_\sk^{(j_1,j_2)}(\mathbf{0}_{d_Y},\mathbf{0}_{d_X})-\wt G_\sk^{(j_1,j_2)}(\mathbf{0}_{d_Y},\mathbf{0}_{d_X})\big\|_2h_1^{|j_1|}h_2^{|j_2|}\bigg)^2\\
\quad &+\big(h_1^{2\beta_Y}+h_2^{2\beta_X}\big)\cdot \big(\frac{\log n}{n}+ h_1^{d_Y}h_2^{d_X}\big).
 \end{align*}
Then the claimed result is a consequence of a simple union bound over $k\in [n]$.

\subsubsection{Proof of Lemma~\ref{lemmamanifold2}}\label{prooflemma14}
Recall $h_1=b_1\,(\frac{\log n}{n})^{\frac{1}{d_Y+\frac{d_X\beta_Y}{\beta_X}}}$ and $h_2=b_2\,(\frac{\log n}{n})^{\frac{1}{d_X+\frac{d_Y\beta_X}{\beta_Y}}}$, since $\beta_Y\geq \beta_X$ and $\beta_Y\geq 2$, we have  $h_2^{1\wedge \beta_X}\leq \frac{b_2}{b_1}h_1$. Then when $\frac{b_2}{b_1}$ is small enough,  it holds for some positive constants $C,C_1$ that,  
\begin{equation*}
    \begin{aligned}
   \forall (x^*,y^*)\in \m M,\qquad     &\mb P_{\mu^*}(\|y-y^*\|<h_1/2,\|x-x^*\|< h_2/2)\\
        &\geq \mb P_{\mu^*_X}(\|x-x^*\|<h_2/2)\cdot \underset{x\in \mb B_{\m M_X}(x^*,h_2/2)}{\inf}\mb P_{\mu^*_{Y|x}}(\|y-y^*\|< h_1/2)\\
        &\overset{(i)}{\geq} C\, h_2^{d_X}\cdot  \underset{x\in \mb B_{\m M_X}(x^*,h_2/4)}{\inf}\mb P_{\mu^*_{Y|x}}(\|y-\Phi_{(x^*,y^*)}(\mathbf{0},x)\|< \frac{h_1}{4})\\
        &\geq C_1\, h_2^{d_X} h_1^{d_Y},\\
    \end{aligned}
\end{equation*}
 where $(i)$ uses the fact that $\mb P_{\mu^*_X}(\|x-x^*\|<h_2/2)\gtrsim h_2^{d_X}$ and $\|y^*-\Phi_{(x^*,y^*)}(\mathbf 0,x)\|=\|\Phi_{(x^*,y^*)}(\mathbf{0},x^*)-\Phi_{(x^*,y^*)}(\mathbf{0},x)\|\leq L\, \|x^*-x\|^{\beta_X\wedge 1}< \frac{h_1}{4}$ when $\frac{b_2}{b_1}$ is sufficiently small.  Furthermore, by Bernstein's inequality, there exists a constant $C_2$  so that for any $t>0$, it holds with probability at least $1-\exp(-t)$ that 
 \begin{equation*}
     \begin{aligned}
        & \frac{1}{n}\sum_{i=1}^n \mathbf{1}(\|Y_i-y_0\|< h_1/2,\|X_i-x_0\|< h_2/2)-\mb P_{\mu^*}(\|y-y_0\|< h_1/2,\|x-x_0\|< h_2/2)\\
         &\geq -\sqrt{\frac{t}{n}}\sqrt{\mb P_{\mu^*}\Big(\|Y-y_0\|< h_1/2,\|X-x_0\|< h_2/2\Big)}-\frac{t}{3n}\\
         &\geq -\frac{t}{3n}-C_2 \sqrt{\frac{t}{n}}h_2^{d_X/2} h_1^{d_Y/2}.
     \end{aligned}
 \end{equation*}
 Consider $\varepsilon_1=c_1 h_1$ and $\varepsilon_2=c_1 h_2$ with $c_1=(\frac{b_1}{2b_1+2Lb_2})^{\frac{1}{\beta_X\wedge 1}}$.  Let $N^x_{\varepsilon_2}$ be the largest $\varepsilon_2$-packing of $\m M_X$, then by Lemma B.7 of~\cite{10.1214/18-AOS1685}, it holds that $|N^x_{\varepsilon_2}|\lesssim \varepsilon_2^{-d_X}$. Moreover, for each $x\in N^x_{\varepsilon_2}$, let $N^y_{\varepsilon_1}(x)$ be the largest $\varepsilon_1$-packing of $\m M_{Y|x}$, then  $|N^y_{\varepsilon_1}(x)|\lesssim \varepsilon_1^{-d_Y}$. So for any $(x^*,y^*)\in \m M_X$, there exists $x_0\in N^x_{\varepsilon_2}$ so that $\|x^*-x_0\|\leq \varepsilon_2$. Moreover, there exists $y_0\in N^y_{\varepsilon_1}(x_0)$ so that $\|y_0-\Phi_{(x^*,y^*)}(\mathbf 0,x)\|\leq \varepsilon_1$ and thus $\|y_0-y^*\|\leq \varepsilon_1+\|\Phi_{(x^*,y^*)}(\mathbf 0,x)-\Phi_{(x^*,y^*)}(\mathbf 0,x^*)\|\leq \varepsilon_1+L\,\varepsilon_2^{\beta_X\wedge 1}$.   By a union argument over $\{(x,y)\,:\,x\in N^x_{\varepsilon_2}, y\in N^y_{\varepsilon_1}(x)\}$, there exists a constant $C_3$ so that it holds with probability at least $1-n^{-1}$ that for any $x_0\in N^x_{\varepsilon_2}$ and $y_0\in N^y_{\varepsilon_1}(x_0)$,
\begin{equation*}
    \begin{aligned}
      & \frac{1}{n}\sum_{i=1}^n \mathbf{1}(\|Y_i-y_0\|<h_1/2,\|X_i-x_0\|< h_2/2)\\
        &= \mb P_{\mu^*}(\|y-y_0\|<h_1/4,\|x-x_0\|< h_2/4)\\
        &+\frac{1}{n}\sum_{i=1}^n \mathbf{1}(\|Y_i-y_0\|< h_1/2,\|X_i-x_0\|< h_2/2)-\mb P_{\mu^*}(\|y-y_0\|< h_1/2,\|x-x_0\|< h_2/2) \\
        &\geq  C_1\, h_2^{d_X} h_1^{d_Y}- C_3\frac{\log n}{3n}- C_2\sqrt{\frac{C_3\log n}{n}}h_2^{d_X/2} h_1^{d_Y/2}.
    \end{aligned}
\end{equation*}
When $b_1,b_2$ are sufficiently large, we have $C_1\, h_2^{d_X} h_1^{d_Y}- C_3\frac{\log n}{3n}- C_2\sqrt{\frac{C_3\log n}{n}}h_2^{d_X/2} h_1^{d_Y/2}>0$,
which means for any $x_0\in N^x_{\varepsilon_2}$ and $y_0\in N^y_{\varepsilon_1}(x_0)$, there exists $i\in [n]$ so that $\|Y_i-y_0\|< h_1/2$ and $\|X_i-x_0\|<h_2/2$. Then, combined with the fact that for any $(x^*,y^*)\in \m M$, there exists  $x_0\in N^x_{\varepsilon_2}$ and $y_0\in N^y_{\varepsilon_1}(x_0)$ so that  $\|x^*-x_0\|\leq \varepsilon_2< \frac{h_2}{2}$ and $\|y^*-y_0\|\leq \varepsilon_1+L\,\varepsilon_2^{\beta_X\wedge 1}\leq \frac{h_1}{2}$, we can get the desired result.

\end{document}